\documentclass[a4paper, twoside, titlepage]{book}
\pdfoutput=1

\usepackage[latin1]{inputenc}
\usepackage[T1]{fontenc}
\usepackage{lmodern, hyphenat}
\usepackage[ngerman, english]{babel}
\usepackage{fancyhdr, titletoc, enumitem, mathtools, microtype}
\usepackage[newparttoc]{titlesec}

\usepackage{xcolor,framed}
\usepackage{makeidx}
\usepackage{array, graphicx, pgf, tikz}
\usetikzlibrary{arrows,calc,decorations.markings}
\usepackage[footskip=40pt,width=420pt,height=665pt,top=3.15cm,headheight=12.6pt]{geometry}
\usepackage[Lenny]{fncychap}
\usepackage{etoolbox}
\makeatletter
\patchcmd{\chapter}{\if@openright\cleardoublepage\else\clearpage\fi}{}{}{}
\makeatother
\usepackage{amsmath, amssymb, amsthm, stmaryrd, cancel}

\usepackage[colorlinks]{hyperref}

\usepackage[acronym,toc,style=super3col]{glossaries}

\newglossary{symbols}{sym}{sbl}{List of Symbols}
\newglossary{foundations}{fou}{fbl}{List of Foundations}
\newglossary{results}{res}{rbl}{List of Results}
\makeglossaries

\makeindex

\begin{document}

\mathtoolsset{centercolon}

\def\mathllap{\mathpalette\mathllapinternal}
\def\mathllapinternal#1#2{%
\llap{$\mathsurround=0pt#1{#2}$}
}
\def\clap#1{\hbox to 0pt{\hss#1\hss}}
\def\mathclap{\mathpalette\mathclapinternal}
\def\mathclapinternal#1#2{%
\clap{$\mathsurround=0pt#1{#2}$}%
}
\def\mathrlap{\mathpalette\mathrlapinternal}
\def\mathrlapinternal#1#2{%
\rlap{$\mathsurround=0pt#1{#2}$}
}

\def\itemMath#1{\raisebox{-\abovedisplayshortskip}
{\parbox{1.0\linewidth}{\begin{equation}\begin{split}#1\end{split}\end{equation}}}}

\def\itemMathtwo#1#2{\raisebox{-#1pt}
{\parbox{1.0\linewidth}{\begin{equation}\begin{split}#2\end{split}\end{equation}}}}

\renewcommand{\AA}{\mathbb{A}}
\newcommand{\CC}{\mathbb{C}}
\newcommand{\DD}{\mathbb{D}}
\newcommand{\EE}{\mathbb{E}}
\newcommand{\FF}{\mathbb{F}}
\newcommand{\HH}{\mathbb{H}}
\newcommand{\KK}{\mathbb{K}}
\newcommand{\LL}{\mathbb{L}}
\newcommand{\MM}{\mathbb{M}}
\newcommand{\NN}{\mathbb{N}}
\newcommand{\OO}{\mathbb{O}}
\newcommand{\PP}{\mathbb{P}}
\newcommand{\QQ}{\mathbb{Q}}
\newcommand{\RR}{\mathbb{R}}
\newcommand{\ZZ}{\mathbb{Z}}

\newcommand{\ssi}{\scriptsize}
\newcommand{\draws}[1]{\draw[postaction=decorate,decoration={markings,mark=at position #1 with {\arrow[scale=1.25]{>}}}]}
\newcommand{\drawd}[1]{\draw[double,postaction=decorate,decoration={markings,mark=at position #1 with {\arrow[scale=0.75]{>}}}]}

\newcommand{\AAA}{\mathcal{A}}
\newcommand{\BBB}{\mathcal{B}}
\newcommand{\CCC}{\mathcal{C}}
\newcommand{\EEE}{\mathcal{E}}
\newcommand{\FFF}{\mathcal{F}}
\newcommand{\GGG}{\mathcal{G}}
\newcommand{\HHH}{\mathcal{H}}
\newcommand{\LLL}{\mathcal{L}}
\newcommand{\MMM}{\mathcal{M}}
\newcommand{\OOO}{\mathcal{O}}
\newcommand{\PPP}{\mathcal{P}}
\newcommand{\QQQ}{\mathcal{Q}}
\newcommand{\RRR}{\mathcal{R}}
\newcommand{\SSS}{\mathcal{S}}
\newcommand{\TTT}{\mathcal{T}}
\newcommand{\UUU}{\mathcal{U}}
\newcommand{\VVV}{\mathcal{V}}
\newcommand{\XXX}{\mathcal{X}}
\newcommand{\ZZZ}{\mathcal{Z}}

\newcommand{\nt}{\trianglelefteq}
\newcommand{\R}{\Rightarrow}
\newcommand{\sm}{\setminus}

\newcommand{\Char}{\mathrm{Char}\,}
\newcommand{\grad}{\mathrm{grad}\,}

\newcommand{\p}{\varphi}
\newcommand{\Bild}{\mathrm{Bild}\,}
\newcommand{\Kern}{\mathrm{Ker}\,}
\newcommand{\Rang}{\mathrm{Rang}\,}

\newcommand{\Emb}{\mathrm{Emb}}
\newcommand{\Fix}{\mathrm{Fix}}
\newcommand{\Quot}{\mathrm{Quot}}

\newcommand{\IB}{Integrit\"atsbereich}
\newcommand{\EK}{Erweiterungskörper}
\newcommand{\HIR}{Hauptidealring}

\newcommand{\id}{\mathrm{id}}
\newcommand{\Aut}{\mathrm{Aut}}
\newcommand{\End}{\mathrm{End}}
\newcommand{\Hom}{\mathrm{Hom}}

\newcommand{\Alt}{\mathrm{Alt}}
\newcommand{\Sym}{\mathrm{Sym}}
\newcommand{\Trans}{\mathrm{Trans}}

\newcommand{\GL}{\mathit{\Gamma}L}

\newcommand{\sign}{\mathrm{sign}}

\newcommand{\Lin}{\mathrm{Lin}}

\newcommand{\ggT}{\mathrm{ggT}}
\newcommand{\kgV}{\mathrm{kgV}}
\newcommand{\lcm}{\mathrm{lcm}}

\makeatletter
\renewcommand*{\thepart}{\@Roman\c@part}
\titleformat{\part}[display]{\center\normalfont\huge\bfseries}{\partname\ \thepart}{20pt}{\Huge} 
\titlecontents{part}[0pt]{\vspace*{7pt}}{\large\bfseries Part\ \thecontentslabel\ \ }{\large\bfseries }{}[]
\titleformat{\chapter}{\normalfont\Large\bfseries}{Chapter \thechapter}{1em}{}
\titleformat{\section}{\normalfont\large\bfseries}{§ \thesection}{1em}{}
\titleformat{\subsection}{\normalfont\normalsize\bfseries}{§ \thesubsubsection}{1em}{}
\titlespacing*{\chapter}{0pt}{0pt}{8pt}
\makeatother

\newtheoremstyle{diplom}{15pt}{15pt}{}{}{\bf }{ }{1.5em}{}
\newtheoremstyle{beweis}{15pt}{15pt}{}{}{\bfseries\itshape}{}{\newline}{}
\newtheoremstyle{beweiszwei}{15pt}{15pt}{}{}{\bfseries\itshape}{}{0pt}{}

\swapnumbers

\theoremstyle{diplom}
\newtheorem{satz}{Theorem}[chapter]
\newtheorem{de}[satz]{Definition}
\newtheorem{no}[satz]{Notation}
\newtheorem{kor}[satz]{Corollary}
\newtheorem{bem}[satz]{Remark}
\newtheorem{lemma}[satz]{Lemma}
\newtheorem{prop}[satz]{Proposition}
\newtheorem{bsp}[satz]{Example}
\newtheorem{bspe}[satz]{Examples}

\theoremstyle{beweis}
\newtheorem*{bew}{Proof}

\theoremstyle{beweiszwei}
\newtheorem*{bewzwei}{Proof}

\newenvironment{sonst}[1]{\noindent\textbf{\stepcounter{satz}\thesatz\ #1}\hspace*{1,5em}}{}

\makeatletter
\g@addto@macro\th@diplom{\renewcommand*{\p@enumi}{\thesatz~}}
\makeatother

\newcommand*\alt{}
\let\alt\thesatz
\renewcommand\thesatz{(\alt)}

\newenvironment{fshaded}[1][]{%
\def\FrameCommand{\fcolorbox{framecolor}{shadecolor}}%
\MakeFramed{\hsize\textwidth\advance\hsize-\width\FrameRestore}}%
{\endMakeFramed}






\definecolor{darkblue}{rgb}{0,0.1,0.4}
\hypersetup{breaklinks=true, linkcolor=darkblue, menucolor=darkblue, urlcolor=darkblue}

\thispagestyle{empty} 	
\begin{center}
\begin{flushright}
\vspace*{-1cm}
\noindent\rule[5pt]{\textwidth}{1pt}
\includegraphics[width=0.30\textwidth]{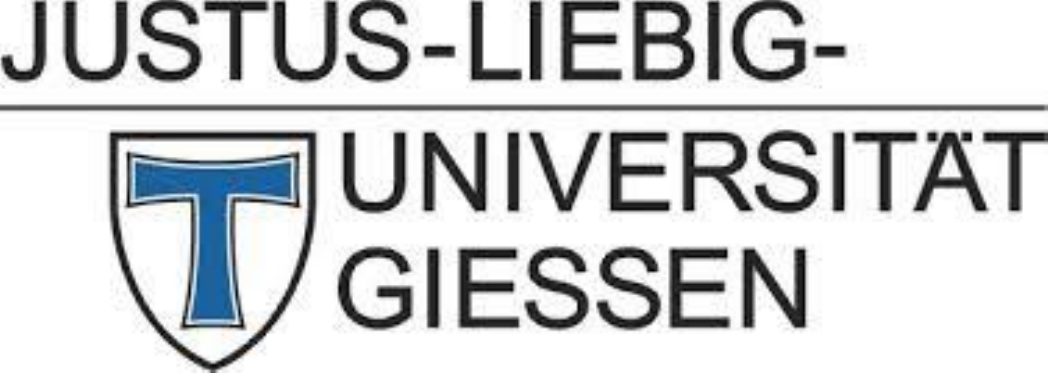}
\noindent\rule[8pt]{\textwidth}{1pt}
\end{flushright}
\ \vspace{2.0cm}\\
\textbf{\huge{Integrability of}\vspace{10pt}\\
\huge{Moufang Foundations}}
\\\vspace{2cm}\Large{\textbf{Dissertation}}
\\\vspace{1cm}\large{zur Erlangung des akademischen Grades doctor rerum naturalium}
\\\vspace{1cm}\normalsize{vorgelegt von}
\\\vspace{1cm}
\textbf{\large{Herrn Dipl.-Math. Sebastian Weiß}}\\
geb. am 18.01.1985 in Trostberg
\vspace{1cm}\\
\normalsize{am}
\vspace{1cm}\\
\large{Mathematischen Institut \\ der Justus-Liebig-Universität Gießen}\\
\vspace{1cm}

September 2013
\ \vspace{3.0cm}\\
\noindent\rule{\textwidth}{1pt}
Betreuer: Prof. Dr. Bernhard Mühlherr \vspace*{-1cm}
\noindent\rule[7pt]{\textwidth}{1pt}
\end{center}

\clearpage

\thispagestyle{empty} 

\ \\\vspace*{18cm}

\noindent \rule[-5pt]{\textwidth}{1pt}\\ \ \\
\textit{Wer fragt, ist ein Narr für fünf Minuten. Wer nicht fragt, bleibt ein Narr für immer.}
\begin{flushright}
Aus China
\end{flushright}
\noindent \rule[7pt]{\textwidth}{1pt}

\clearpage

\thispagestyle{empty} 	

\begin{center} \  \\\vspace{1.5cm}
\noindent \underline{\hspace{420pt}}\\ \ \\ \ \\
\textbf{\begin{LARGE}Integrability of Moufang Foundations\end{LARGE}}\\\vspace{0.5cm}
\textbf{\begin{large}A Contribution to the Classification of Twin Buildings\end{large}}\\\vspace{1cm}
Dedicated to My Siblings
\\ \  \ \\
\noindent \underline{\hspace{420pt}}
\ \vspace{1cm}\\
\hspace*{0.75cm}\includegraphics[width=0.90\textwidth]{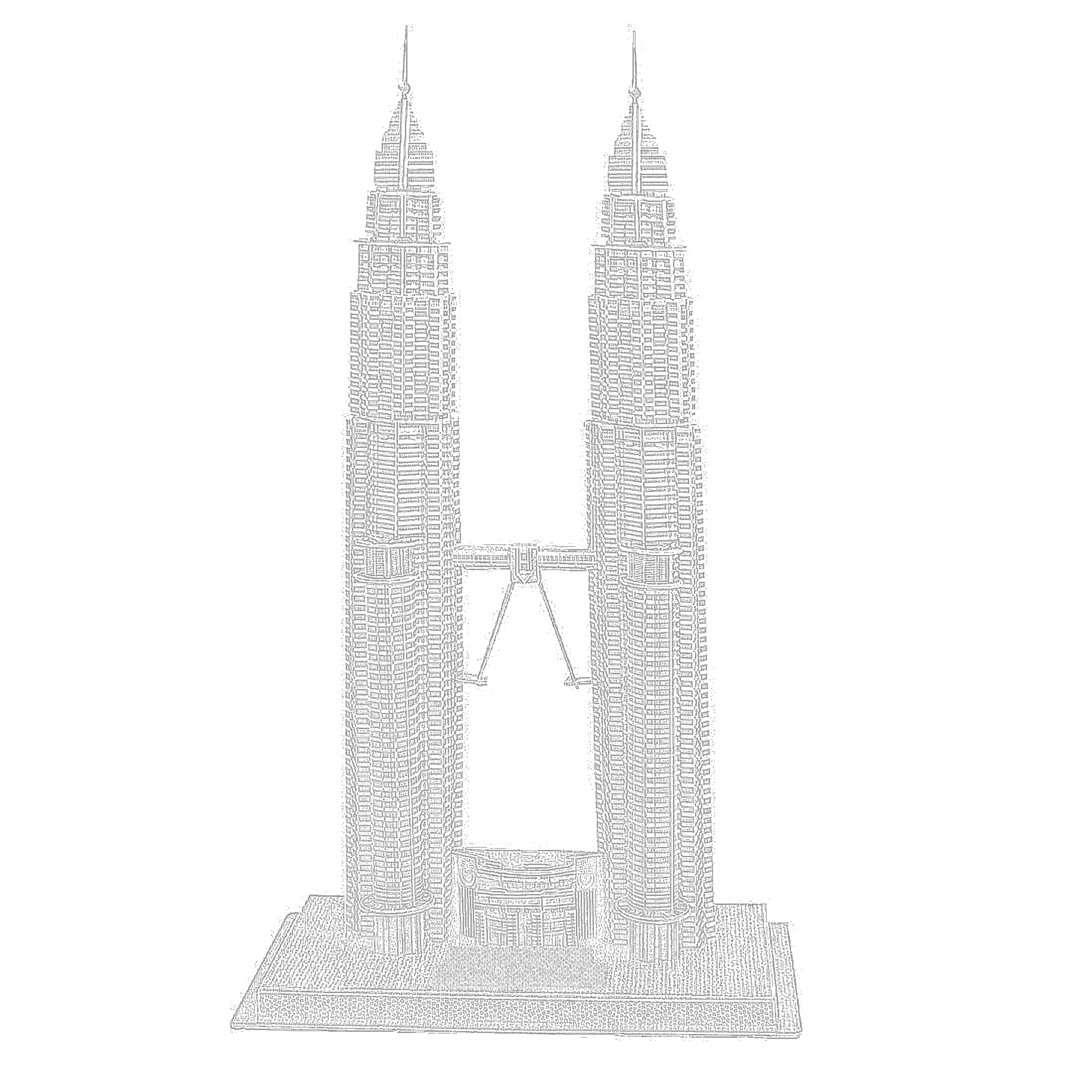}
\end{center}

\clearpage

\thispagestyle{empty}

\ \\\vspace*{18cm}

\noindent \rule[-5pt]{\textwidth}{1pt}\\ \ \\
\textit{Die Zukunft gehört denjenigen, die an die Wahrhaftigkeit ihrer Träume glauben.}
\begin{flushright}
Eleanor Roosevelt
\end{flushright}
\noindent \rule[7pt]{\textwidth}{1pt}

\cleardoublepage

\pagestyle{fancy}

\renewcommand{\headrulewidth}{0.4pt}
\renewcommand{\footrulewidth}{0.4pt}
\renewcommand{\footruleskip}{20pt}

\fancypagestyle{plain}{\lhead[\nouppercase{\leftmark}]{}
\rhead[]{\nouppercase{\rightmark}}
\cfoot{-\ \thepage\ -}}

\pagestyle{fancy}
\pagenumbering{Roman}\cfoot{-\ \thepage\ -} 	
\lhead[Contents]{}
\rhead[]{Contents}
\tableofcontents
\newpage

\renewcommand{\partmark}[1]{%
\markboth{\partname
\ \thepart\ #1}{}}
\renewcommand{\chaptermark}[1]{%
\markright{\chaptername
\ \thechapter\ #1}}
\renewcommand{\sectionmark}[1]{{}}
\lhead[\leftmark]{}
\rhead[]{\rightmark}
\pagenumbering{arabic}\cfoot{-\ \thepage\ -}

\addtocontents{toc}{\vfill}
\addtocontents{toc}{\noindent\protect\mbox{}\protect\hrulefill\par}
\part*{Introduction}\markboth{Introduction}{Introduction}
\addtocontents{toc}{\noindent\protect\mbox{}\protect\hrulefill\par}

\section*{Historical and Theoretical Context}

\noindent The description below closely follows those given in \cite{M} and \cite{AB}.

\subsection*{Twin Buildings}

Buildings have been introduced by J. Tits in order to study semi-simple algebraic groups from a geometrical point of view. One of the most important results in the theory of buildings is the classification of irreducible spherical buildings of rank at least 3 in \cite{Ti74}. Meanwhile, there is a simplified proof in \cite{TW} which makes use of the classification of Moufang polygons.

About 25 years ago, M. Ronan and J. Tits defined a new class of buildings, which generalize spherical buildings in a natural way, namely the class of twin buildings. The motivation
of their definition is provided by the theory of Kac-Moody groups, and we refer to \cite{Ti} for further general information about twin buildings.

The sense in which twin buildings generalize spherical buildings is the following: Given a building of spherical type, there is a natural \textit{opposition relation} on the set of its chambers. This relation restricts the structure of spherical buildings essentially. The classification of irreducible spherical buildings of rank at least 3 mentioned above is in fact based on this opposition relation. The idea of a twin building is to introduce a symmetric relation between the chambers of two different buildings of the same type which has properties similar to the opposition relation of a spherical building. Thus a twin building is a triple consisting of two buildings of the same type and an opposition relation between the chambers of the two ``halves'' of the twin building.

\subsection*{The Classification Program for 2-Spherical Twin Buildings}

In view of the classification of spherical buildings, it is natural to ask whether it is possible to classify higher rank twin buildings. A large part of \cite{Ti} deals with this question. As a first observation, it turns out that such a classification seems only to be feasible under the additional assumption that the entries in the corresponding Coxeter matrices are all finite. We call these buildings \textit{2-spherical}. The classification program described in \cite{Ti} is based on the conjecture that there is a bijective correspondence between twin buildings of type $M$ and certain Moufang foundations of type $M$ for each 2-spherical Coxeter diagram of type $M$.

Foundations have been introduced by M. Ronan and J. Tits in \cite{RT} in order to describe chamber systems which are candidates for being the local structure of a building. Roughly
speaking, foundations can be seen as amalgams of rank 2 buildings which are glued along certain rank 1 residues. Given a chamber $c$ of a building $\BBB$ of type $M$, the union $E_2(c)$ of the rank 2 residues which contain this chamber constitutes a foundation of type $M$, the \textit{foundation of $\mathit{\BBB}$ at $\mathit{c}$}. Thus the term ``local structure'' above has to be understood as a kind of 2-neighbourhood of a given chamber of a building. 

It is a (not completely trivial) fact that if two chambers are contained in the same half of a twin building, the foundations at these chamber are isomorphic. Moreover, if one knows the isomorphism class of the foundation of one half of a twin building, then the isomorphism class of the foundations of the other half is uniquely determined. Conversely, a generalization of Tits' extension theorem by B. Mühlherr and M. Ronan in \cite{MR} states that a twin building is uniquely determined by the foundation of one of its halves in almost all cases, cf. (5.10), (*5.11), (*9.11) and (*9.12) of \cite{AB} for a summary. Thus the foundation at a chamber of a twin building is a classifying invariant of the corresponding twin building if the following condition is satisfied:
\begin{enumerate}[leftmargin=30pt]
\item[(CO)]
No rank 2 residue is isomorphic to one of the buildings which are associated with the groups $B_2(2)$, $G_2(2)$, $G_2(3)$ and $^{2}{F_4}(2)$.
\end{enumerate}
This condition guarantees that for every chamber $c\in \BBB_\epsilon$ ($\epsilon\in \{\pm\}$), the set $c^o$ of chambers opposite $c$ is a gallery-connected subset of $\BBB_{-\epsilon}$.

In view of what has been mentioned so far, the classification of 2-spherical twin buildings reduces to the classification of all foundations which can be realized as the local structure of a twin building. We call such a foundation \textit{integrable}. In order to determine the integrable foundations, one proceeds in two steps.

\newpage

\subsection*{Step 1: Exclude Non-Integrable Foundations}

It is proved in \cite{Ti} that an integrable foundation is Moufang, which means that the rank 2 buildings in the foundation are Moufang, i.e., they are Moufang polygons, and that the glueings are compatible with the Moufang structures induced on the rank 1 residues. Thus a first necessary condition for the integrability of a foundation is that it is Moufang.

As a consequence, the classification of Moufang polygons in \cite{TW} and the solution of the isomorphism problem for Moufang sets are essential to work out which Moufang polygons fit together in order to form a foundation. Moreover, one can reduce the list of possibly integrable foundations by considering certain automorphisms of the twin building, the so-called \textit{Hua automorphisms}, which are closely related to the double $\mu$-maps of the appearing Moufang sets.

\subsection*{Step 2: Existence / Integrability Proof}

Finally, one has to prove that each of the remaining candidates is in fact integrable, i.e., realized by a twin building, which is then unique up to isomorphism. In \cite{M} and his Habilitationsschrift \cite{MHab}, B. Mühlherr developed techniques which produce certain twin buildings as fixed point structures in twin buildings coming from Kac-Moody groups. He, H. Petersson and R. Weiss actually prepare a book which provides further well-founded background.

\section*{Goals and Main Results}

The present thesis contributes to establish complete lists of integrable foundations for certain types of diagrams. We closely follow the approach for the classification of spherical buildings in \cite{TW}. However, we have to refine the techniques, since in general, foundations don't only depend on the diagram and the defining field. For example, there may be several non-isomorphic foundations of type $\tilde{A}_n$ with respect to a given skew-field $\AA$: Automorphisms of $\AA$ are involved as well, which represents the fact that there are several possibilities for glueing Moufang polygons along a rank 1 residue.

The main question is how to parametrize sequences of Moufang polygons with respect to the usual commutator relations in order to make the glueings visible. The crucial subtlety is the following: Each Moufang polygon is parametrized twice, once for each direction in which the underlying root group sequence can be read. As a consequence, we obtain glueings between directed Moufang polygons, and it's a difference whether we look at $\id_\AA:\AA\to \AA$ or $\id_\AA^o:\AA\to \AA^o$, where $\AA^o$ is the opposite with respect to $\AA$: The former is an isomorphism, while the latter is an anti-isomorphism of skew-fields.

As mentioned above, excluding non-integrable foundations is closely related to the investigation of Moufang sets and their isomorphisms. Therefore, a large part deals with the introduction of underlying parameter systems and, in the sequel, with the solution of the isomorphism problem for Moufang sets. Many questions have already been answered, cf. \cite{K}, but we need to refine and extend the existing results for our purposes and translate their proofs into our setup. 

\subsection*{Simply Laced Foundations}

The main result of this thesis is the complete classification of simply laced twin buildings via their foundations. Of course, the basic requirement for a foundation to be integrable is that it is Moufang: Its glueings are Jordan isomorphisms, i.e., they preserve the Jordan product $xyx$.

A powerful tool is Hua's theorem, cf. \cite{H} for a reference, which answers the isomorphism problem for Moufang sets of skew-fields: Each Jordan isomorphism is in fact an iso- or anti-isomorphism of skew-fields. However, the class of parameter systems for Moufang triangles additionally includes octonion division algebras, which cause a lot of trouble due to the lack of associativity. A byproduct is the existence of Jordan isomorphisms which are neither iso- nor anti-isomorphisms of alternative rings. The most sophisticated part is the handling of the exceptional cases where octonions occur.

We give an overview of the restriction process and point out the main ideas.\bigskip

\noindent The following observations yield the first restriction of possibilities: 
\begin{enumerate}[label=(\arabic*)]
\item Each Moufang triangle is defined over the same alternative division ring $\AA$.
\item An integrable foundation of type $A_3$ is necessarily defined over a skew-field, and the glueing is necessarily an isomorphism of skew-fields.
\end{enumerate}

\noindent Thus the crucial step is the classification of integrable foundations of type $\tilde{A}_2$ since these are the smallest ones which allow weird ``non-standard'' things to appear. The theory of affine and Bruhat-Tits buildings and the theory of composition algebras which are complete with respect to a discrete valuation enable us to get further restrictions:
\begin{enumerate}
\item[(3)] Given an octonion division algebra $\OO$, there is only one twin building of type $\tilde{A}_2$ with respect to $\OO$.
\item[(4)] An integrable foundation of type $\tilde{A}_2$ whose glueings are anti-isomorphisms is necessarily defined over a quaternion division algebra, and given a quaternion division algebra $\HH$, there is only one such ``positive'' twin building of type $\tilde{A}_2$ with respect to $\HH$.
\end{enumerate}

\noindent A closer look at the group of Jordan automorphisms of octonion division algebras completes the classification of integrable foundations which are defined over octonions:

\begin{enumerate}
\item[(5)] There are no integrable foundations over octonions such that the corresponding graph is a tetrahedron. In particular, up to isomorphism, the only integrable foundations with respect to an octonion division algebra $\OO$ are $\AAA_2(\OO)=\TTT(\OO)$ and $\tilde{\AAA}_2(\OO)$. 
\end{enumerate}

\noindent Finally, in connection with (4), the following observation heavily restricts the list of integrable foundations over non-commutative skew-fields which are not quaternion division algebras:
\begin{enumerate}
\item[(6)] An integrable foundation of type $D_4$ is necessarily defined over a field.
\end{enumerate}

\noindent Kac-Moody theory provides the integrability proofs as the corresponding Coxeter diagram is a tree. The remaining integrability proofs rely on techniques developed by B. Mühlherr. 

\subsection*{Jordan Automorphisms of Alternative Division Rings}

In view of Hua's theorem
$$\Aut_J(\DD)=\Aut(\DD)\cup \Aut^o(\DD)$$
for any skew-field $\DD$, its group $\Aut_J(\DD)$ of Jordan automorphisms, its subgroup $\Aut(\DD)$ of automorphisms and its set $\Aut^o(\DD)$ of anti-automorphisms, the question arises whether it is possible to get a similar result for octonion division algebras.

In the proof that integrable tetrahedron-foundations over an octonion division algebra $\OO$ do not exist, we define a subset $\Gamma\subseteq \Aut_J(\OO)$ which turns out to not contain the standard involution $\sigma_s$. The elements of $\Gamma$ are automorphisms of $\OO$ multiplied with one of the ``exceptional'' Jordan automorphisms as defined in \cite{TW}, which fix a quaternion subalgebra $\HH$ pointwise and which act on the orthogonal complement of $\HH$ as conjugation.

The fact that $\Gamma$ is a subgroup of $\Aut_J(\OO)$ can be deduced from the knowledge about the automorphism group of the corresponding Moufang triangle $\TTT(\OO)$. This subgroup $\Gamma$ corresponds to the subgroup $\Aut(\DD)$ in Hua's theorem, i.e., we obtain
$$\Aut_J(\OO)=\langle \sigma_s,\Gamma\rangle=\Gamma\cup \sigma_s\Gamma\ .$$
The strategy for the proof is as follows:
\begin{enumerate}[label=(\arabic*)]
\item Jordan automorphisms restricted to subfields are monomorphisms of rings, i.e., the image of a subfield is again a subfield.
\item As an immediate consequence, Jordan automorphisms of octonions are norm similarities.
\item The results of \cite{S} allow us to restrict to isometries which fix a quaternion subalgebra pointwise.
\item Hua's theorem and the Skolem-Noether theorem allow us to show that any Jordan automorphism is indeed a product in $\langle \sigma_s,\Gamma\rangle$.
\end{enumerate}

\subsection*{443-Foundations}

The second result in connection with the classification of twin buildings is the completion of step 1 for 443-foundations whose diagram is a triangle and whose Moufang polygons are two quadrangles and one triangle. Although we only deal with a single diagram in this case, there is a rich variety of integrable 443-foundations as there are six families of Moufang quadrangles which often fit together in this configuration. Nevertheless, quadrangles of type $E_n$, of type $F_4$ and of indifferent type don't appear since their Moufang sets are not of linear type, i.e., they aren't projective lines.

The same holds for Moufang sets of pseudo-quadratic form and involutory type, but the second panel of the corresponding \textit{unitary} quadrangle is of linear type so that there is exactly one possibility for the orientation of the quadrangles. The solution of the isomorphism problem for the appearing Moufang sets and the knowledge about the automorphism group of a unitary quadrangle allow us to show the following: 
\begin{enumerate}[label=(\arabic*)]
\item The appearing pseudo-quadratic spaces are defined over a quaternion division algebra $\HH$ or over a separable quadratic extension $\EE$.
\item In the former case, there is exactly one possibly integrable 443-foundation with respect to such a pseudo-quadratic space $\Xi$.
\item In the latter case, the isomorphism class of a possibly integrable 443-foundations additionally depends on an automorphism $\gamma\in \Aut(\EE)$.
\item The appearing involutory sets are defined over a quaternion division algebra $\HH$, and there is exactly one possibly integrable 443-foundation with respect to such an involutory set $\Xi$.
\end{enumerate}
Finally, quadrangles of quadratic form type are the most flexible ones since there are Moufang sets which are both of quadratic form type and of linear type so that they can glued together in any orientation. Furthermore, there is one point where we need to restrict to proper quadratic spaces as parametrizing structures to exclude characteristic 2 phenomenons in order to obtain a satisfying description. 

In contrast to the classification of integrable simply laced foundations however, we omit step 2 in the classification program as the proofs require different kinds of techniques, established by B. Mühlherr, H. Petersson and R. Weiss. As before, there are two possibilities how to show the integrability of a given foundation: Either the universal cover is isomorphic to a canonical foundation, which is a foundation such that each glueing is the identity map and thus integrable if it the corresponding diagram is a tree, or the foundation can be obtained as a fixed point structure via a Tits index. The former method applies to 443-foundations with quadrangles of quadratic form type, while the latter applies to 443-foundations with unitary quadrangles.

\subsection*{Jordan Isomorphisms of Pseudo-Quadratic Spaces}

Hua's theorem is essential for the classification of integrable simply laced foundations. In the same spirit, the solution of the isomorphism problem for the appearing Moufang sets is essential for the classification of integrable 443-foundations. As mentioned above, R. Knop handles a lot of cases in his PhD thesis \cite{K}. However, he only deals with commutative Moufang sets. Thus we need to establish the corresponding results for Moufang sets of pseudo-quadratic form type.

We obtain that Jordan isomorphism between two Moufang sets of pseudo-quadratic form type are induced by isomorphisms of the corresponding pseudo-quadratic spaces in almost all cases, i.e., whenever the dimension is at least 3 or the involved involutory set is proper. As a consequence, exceptions necessarily involve pseudo-quadratic spaces of small dimensions which are defined over a quaternion division algebra or over a separable quadratic extension. Luckily, these exceptional cases don't occur in the classification of integrable 443-foundations. As a consequence, both the quadrangles are defined over the same pseudo-quadratic space $\Xi$.
\section*{Outlook and Open Problems}

\subsection*{Jordan Isomorphisms}

In the theory of Moufang sets, the $\mu$-maps and the Hua maps play a central role as they carry a lot of information. As a consequence, Jordan isomorphisms -- which are additive isomorphisms preserving the Hua maps -- are closely related to isomorphisms of Moufang sets. In fact, each isomorphism of Moufang sets is a Jordan isomorphism since the Hua maps can be expressed in terms of sums and the permutation $\tau$.

In this context, the following question naturally arises: Is each Jordan isomorphism an isomorphism of Moufang sets? Of course, the Hua maps of sharply 2-transitive Moufang sets are trivial. Therefore, the question has to be answered negatively for these ``non-proper'' Moufang sets. But experts in the area such as R. Weiss and T. De Medts are optimistic that both the definitions are equivalent if we restrict to proper Moufang sets.

\subsection*{The Classification Program}

The main conjecture in connection with the classification program is the following, cf. page 5 in \cite{MHab}:
\begin{itemize}
\item[] A Moufang foundation of 2-spherical type is integrable if and only if each of its rank 3 residues is integrable.
\end{itemize}
In his Habilitationsschrift \cite{MHab}, B. Mühlherr indicates that one could prove the conjecture under the additional assumption that all rank 3 residues are spherical, which is of course a severe restriction. However, there isn't any written proof yet.

Once one has proved the conjecture, the classification program reduces to the classification of integrable Moufang foundations of rank 3. Most of them can be handled with the methods established in \cite{MHab} and \cite{M}. However, there are some exceptions, the most complicated of which are foundations of type $\tilde{C}_2$, $\tilde{A}_2$ and $443$-foundations. The $\tilde{A}_2$- and the $443$-case are solved in the present thesis, while there are (unpublished) partial results for the $\tilde{C}_2$-case by T. De Medts, B. Mühlherr, H. Van Maldeghem and R. Weiss.

\subsection*{The Classification of Simply Laced Twin Buildings}

Although the classification of integrable simply laced foundations is complete, we don't make any statement whether two given foundations in our list are isomorphic. By taking classifying invariants into account and introducing suitable parameters, one could create a list with pairwise non-isomorphic foundations.

If the underlying Coxeter diagram $\GGG_F$ is a tree, the foundations $\FFF$ depends only on the defining field. Circles in the diagram cause an additional dependence on ``twists", i.e., on automorphisms of the defining field $\AA$. More precisely:
\begin{itemize}
\item If $\AA$ is a field, an integrable foundation $\FFF$ is uniquely determined by $\GGG_F$ and a homomorphism $\p: \Pi_1(\GGG_F)\to \Aut(\AA)/\mathrm{Inn}(\AA)\cong \Aut(\AA)$, where $\Pi_1(\GGG_F)$ is the fundamental group of $\GGG_F$.
\item If $\AA$ is a skew-field distinct from a quaternion division algebra and $\FFF$ is an integrable foundation of type $\tilde{A}_n$, the foundation is uniquely determined by $n$ and an element of $\Aut(\AA)/\mathrm{Inn}(\AA)$.
\item If the defining field is a quaternion division algebra, then a similar result as in the field case holds.
\end{itemize}
Moreover, the integrability proofs might be improved at some points as soon as the applied theory is developed properly by B. Mühlherr, H. Petersson and R. Weiss.

\subsection*{Finite Moufang Foundations}

The introduced terminology and the methods of \cite{MHab} can be used to show that each locally finite twin building of 2-spherical type is the fixed point building of a Galois action in the sense of B. Rémy, which means that it is of algebraic origin.

\section*{Acknowledgments}

I would like to express my gratitude towards my primary advisor Bernhard Mühlherr for drawing my attention to the interesting field of Moufang foundations and their Moufang sets: It was a pleasure to contribute to the classification program for twin buildings. Many fruitful discussions showed me the right direction, i.e., the assertions that one should be able to prove and the way how to achieve it. His intuition is tremendous.

Furthermore, I would like to thank Richard Weiss, who raised the question for the generalization of Hua's theorem and who laid the foundation for this thesis with his wonderful and detailed work on Moufang polygons, spherical and affine buildings as well as their classification. Many little questions could be answered with the aid of \cite{TW}, and if not, he always had the right idea where to look for the solution.

Ralf Köhl enhanced the work group with many nice people and gave me the opportunity to take part in a research project about compact subgroups of Kac-Moody groups. His enthusiasm and his dedication are impressive and build the basis for our prospering group.

Thanks go also to the following people: Tom De Medts who provided me a pleasant stay in Ghent during which we worked on the isomorphism problem for Moufang sets of pseudo-quadratic form type. And besides my family finally, there are so many friends who enriched my life with time, conversations and activities that were an essential contrast to the ``abstract nonsense'' called mathematics.

The most important factor which made the last five years a wonderful time of my life is the following: Bernhard and Ralf both endowed me with maximal flexibility in my academic work. In particular, this allowed me and two old friends of mine, Steffen Presse and Joram Gornowitz, to realize our ambitious cinema project.

\newpage

\addtocontents{toc}{\noindent\protect\mbox{}\protect\hrulefill\par}
\part{Preliminaries}
\addtocontents{toc}{\noindent\protect\mbox{}\protect\hrulefill\par}

\noindent As we will have to reconstruct structure out of some given identities and as we will make use of several structures at the same time, it is important to have exact definitions and notations and to make precise statements, e.g., it is important to know whether we talk about an isomorphism as an isomorphism of algebras or as an isomorphism of vector spaces, and we avoid ``identification''. 

We start with giving the definitions of elementary structures such as vector spaces, algebras and graphs, then we proceed by introducing Coxeter matrices and Coxeter diagrams before we turn to twin buildings themselves.\\
 
\chapter{Notations}

We fix some notations.

\begin{de}[\textbf{Vector Spaces}]\ 
\begin{itemize} 
\item A \textit{(right) vector space}\index{vector space} is a pair $(V,\KK)$ consisting of a commutative group $V$ and a skew-field $\KK$ together with a \textit{scalar multiplication} $\cdot: V\times \KK\to V$\index{scalar multiplication} satisfying
\begin{align*}
\forall\ v\in V:\  v\cdot 1_\KK=v\ , && \forall\ v\in V,\ s,t\in \KK:\ (v\cdot s)\cdot t=v\cdot (st)
\end{align*}
and
\begin{align*}
\forall\ v,w\in V,\ s,t\in \KK:&& (v+w)\cdot s=v\cdot s+w\cdot s\ , && v\cdot (s+t)=v\cdot s+v\cdot t\ .
\end{align*}
\item If $(V,\KK)$ is a vector space, we say that \textit{$\mathit{V}$ is a $\mathit{\KK}$-vector space} or that \textit{$\mathit{V}$ is a vector space over $\mathit{\KK}$}.
\item Two vector spaces $(V,\KK)$ and $(\tilde{V},\tilde{\KK})$ are \textit{isomorphic}\index{isomorphism!of vector spaces} if there is a pair $(\p,\phi)$ of isomorphisms $\p:V\to \tilde{V}$ and $\phi:\KK\to \tilde{\KK}$ of groups and skew-fields, resp., satisfying
$$\forall\ s\in \KK,\ v\in V:\qquad \p(v\cdot s)=\p(v)\cdot \phi(s)\ .$$
\item If $(\p,\phi):(V,\KK)\to (\tilde{V},\tilde{\KK})$ is an isomorphism of vector spaces, we say that $\p$ is a \textit{$\mathit{\phi}$-isomorphism}.
\item Let $(V,\KK)$ be a vector space. 
\begin{enumerate}[label=$\circ$]
\item An automorphism $(\p,\phi)$ of $(V,\KK)$ with $\phi=\id_\KK$ is \textit{linear}\index{automorphism!linear}.

\newglossaryentry{SLAV}{type=symbols,name={\ensuremath{\Gamma L(V,\KK)}},description=group of (semi-linear) automorphisms of the vector space ${(V,\KK)}$, sort=vector space}
\item We denote the \textit{group of (semi-linear) automorphisms of $\mathit{(V,\KK)}$} by \gls{SLAV}.

\newglossaryentry{LAV}{type=symbols,name={\ensuremath{GL(V,\KK)}},description=group of linear automorphisms of the vector space ${(V,\KK)}$, sort=vector space}
\item We denote the \textit{group of linear automorphisms of $\mathit{(V,\KK)}$} by \gls{LAV}.
\end{enumerate}
\end{itemize}
\end{de}

\begin{de}[\textbf{Algebras}]\ 
\begin{itemize}
\item An \textit{algebra}\index{algebra} is a pair $(A,\KK)$ such that $A$ is a vector space over a field $\KK$ together with a $\KK$-bilinear map $\cdot: A\times A\to A$.
\item If $(A,\KK)$ is an algebra, we say that \textit{$\mathit{A}$ is a $\mathit{\KK}$-algebra} or that \textit{$\mathit{\AA}$ is an algebra over $\mathit{\KK}$}.
\item An algebra $(A,\KK)$ is \textit{associative}\index{algebra!associative} if the map $\cdot:A\times A\to A$ is associative.
\item Two algebras $(A,\KK)$ and $(\tilde{A},\tilde{\KK})$ are \textit{isomorphic}\index{isomorphism!of algebras} if there is an isomorphism $(\p,\phi)$ of vector spaces satisfying
$$\forall\ x,y\in A:\qquad \p(x\cdot y)=\p(x)\cdot \p(y)\ .$$
\item If $(\p,\phi):(A,\KK)\to (\tilde{A},\tilde{\KK})$ is an isomorphism of algebras, we say that $\p$ is a \textit{$\mathit{\phi}$-isomorphism}.
\item Let $(A,\KK)$ be an algebra. 
\begin{enumerate}[label=$\circ$]
\item An automorphism $(\p,\phi)$ of $(A,\KK)$ with $\phi=\id_\KK$ is \textit{linear}\index{automorphism!linear}.

\newglossaryentry{SLAA}{type=symbols,name={\ensuremath{\Aut(A,\KK)}},description=group of (semi-linear) automorphisms of the algebra ${(A,\KK)}$, sort=algebra}
\item We denote the \textit{group of (semi-linear) automorphisms of $\mathit{(A,\KK)}$} by \gls{SLAA}.

\newglossaryentry{LAA}{type=symbols,name={\ensuremath{\Aut_\KK(A,\KK)}},description=group of linear automorphisms of the algebra ${(A,\KK)}$, sort=algebra}
\item We denote the \textit{group of linear automorphisms of $\mathit{(A,\KK)}$} by \gls{LAA}.
\end{enumerate}
\end{itemize}\end{de}

\newpage

\begin{de}[\textbf{Graphs}]\ 
\begin{itemize}

\newglossaryentry{V(G)}{type=symbols,name={\ensuremath{V(\GGG)}},description=set of vertices of the graph $\GGG$, sort=graph}
\newglossaryentry{E(G)}{type=symbols,name={\ensuremath{E(\GGG)}},description=set of edges of the graph $\GGG$, sort=graph}
\newglossaryentry{A(G)}{type=symbols,name={\ensuremath{A(\GGG)}},description=set of directed edges of the graph $\GGG$, sort=graph}
\newglossaryentry{G(G)}{type=symbols,name={\ensuremath{G(\GGG)}},description=set of rank 3 subdiagrams of the graph $\GGG$, sort=graph}
\newglossaryentry{B(v)}{type=symbols,name={\ensuremath{B_1(v)}},description=set of neighbours of the vertex $v$, sort=graph}
\item A \textit{graph}\index{graph} is a pair $\GGG=(V,E)=\big( \gls{V(G)},\gls{E(G)}     \big)$ consisting of a set of \textit{vertices} \index{vertex} $V$ and a set of \textit{edges}\index{edge}
$$E\subseteq \binom{V}{2}:=\{ X\subseteq V\mid |X|=2\}\ .$$
\item Given two graphs $\GGG=(V,E)$, $\tilde{\GGG}=(\tilde{V},\tilde{E})$, a \textit{morphism $\mathit{\p}:\GGG\to \tilde{\GGG}$ of graphs}\index{morphism!of graphs} is a map $\p:V\to \tilde{V}$ such that
$$\forall\ v,w\in V:\qquad \{v,w\}\in E\ \Rightarrow\ \{ \p(v),\p(w)\}\in \tilde{E}\ .$$
\item Given a graph $\GGG=(V,E)$, we set
$$ \gls{A(G)}:=\{ (i,j)\in V^2 \mid \{i,j\}\in E\}\ ,$$
which is the set of directed edges, and
$$ \gls{G(G)}:=\{(i,j,k)\in V^3 \mid (i,j)\neq(k,j)\in A(\GGG) \}\ .$$
\item Given a graph $\GGG=(V,E)$ and a vertex $v\in V$, the \textit{set of neighbours of $\mathit{v}$}\index{set of neighbours} is
$$\gls{B(v)}:=\{ w\in V \mid \{v,w\}\in E\}\ .$$
\item Given a graph $\GGG=(V,E)$, a \textit{cover of $\mathit{\GGG}$}\index{cover} is a pair $(\tilde{\GGG},\p)$ consisting of a graph $\tilde{\GGG}=(\tilde{V},\tilde{E})$ and an epimorphism $\varphi: \tilde{V}\to V$ of graphs such that for each $v\in \tilde{V}$, the map 
$$\p_{|B_1(v)}: B_1(v)\to B_1\big( \p(v)\big)$$
is a bijection.
\end{itemize}
\end{de}

\begin{bsp}
Given the graph
\begin{center}\begin{tikzpicture}[scale=0.7,>=stealth,thick]
\begin{scope}
\coordinate (1) at (0,0);                    
\coordinate (2) at (3,0);
\coordinate (3) at (60:3);
\draw (1)--(2)--(3)--cycle;
\node () at ($(3)+(0,0.5)$) {\ssi${\bar{3}}$};
\node () at ($(1)+(0,-0.5)$) {\ssi${\bar{1}}$};
\node () at ($(1)+(3,-0.5)$) {\ssi${\bar{2}}$};
\node () at (4,0) {,};
\node () at (-1,0) {};
\foreach \i in {1,...,3} {\fill (\i) circle (3pt);}
\end{scope}\end{tikzpicture}\end{center}
the graphs
\begin{center}\begin{tikzpicture}[scale=0.7,>=stealth,thick]
\begin{scope}
\coordinate (1) at (0,0);                    
\coordinate (2) at (3,0);
\coordinate (3) at ($(2)+(60:3)$);
\coordinate (4) at ($(3)+(120:3)$);
\coordinate (5) at ($(4)+(-3,0)$);
\coordinate (6) at ($(5)+(240:3)$);
\draw (1)--(2)--(3)--(4)--(5)--(6)--cycle;
\node () at ($(3)+(0.5,0)$) {\ssi${{3}}$};
\node () at ($(1)+(0,-0.5)$) {\ssi${{1}}$};
\node () at ($(2)+(0,-0.5)$) {\ssi${{2}}$};
\node () at ($(4)+(0,0.5)$) {\ssi${{4}}$};
\node () at ($(5)+(0,0.5)$) {\ssi${{5}}$};
\node () at ($(6)+(-0.5,0)$) {\ssi${{6}}$};
\foreach \i in {1,...,6} {\fill (\i) circle (3pt);}
\end{scope}\end{tikzpicture}\end{center}
and
\begin{center}\begin{tikzpicture}[scale=0.7,>=stealth,thick]
\begin{scope}
\coordinate (1) at (0,0);                    
\coordinate (2) at (3,0);
\coordinate (3) at (6,0);
\coordinate (4) at (9,0);
\coordinate (5) at (12,0);
\coordinate (6) at (14,0);
\coordinate (7) at (-2,0);
\draw (1)--(2)--(3)--(4)--(5);
\draw[dashed] (5)--(6);
\draw[dashed] (7)--(1);
\node () at ($(3)+(0,-0.5)$) {\ssi${0}$};
\node () at ($(1)+(0,-0.5)$) {\ssi${{-2}}$};
\node () at ($(2)+(0,-0.5)$) {\ssi${{-1}}$};
\node () at ($(4)+(0,-0.5)$) {\ssi${{1}}$};
\node () at ($(5)+(0,-0.5)$) {\ssi${{2}}$};
\foreach \i in {1,...,5} {\fill (\i) circle (3pt);}
\end{scope}\end{tikzpicture}\end{center}
are covers, where $$\p: \ZZ\to \ZZ/3\ZZ,\ z\mapsto \bar{z}$$
is the natural homomorphism in both cases.
\end{bsp}

\newpage

\begin{bem}
We only deal with Coxeter matrices such that $m_{ij}\neq \infty$ for all $i,j\in I$.
\end{bem}

\begin{de}[\textbf{Coxeter Matrices}]\ 
 \begin{itemize}
\item A \textit{(2-spherical) Coxeter matrix over (an index set) $\mathit{I}$}\index{Coxeter!matrix (2-spherical)} is a map $M:I\times I\to \NN^*$ such that
\begin{align*}
\forall\ i\in I:\ m_{ii}=1\ , && \forall\ i\neq j\in I:\ m_{ji}=m_{ij}>1\ ,
\end{align*}
where $m_{ij}:=M(i,j)$ for all $i,j\in I$.
\item Given two Coxeter matrices $M$ over $I$ and $\tilde{M}$ over $\tilde{I}$, a \textit{morphism $\mathit{\p:M\to \tilde{M}}$ of Coxeter matrices}\index{morphism!of Coxeter matrices} is a map $\p:I\to \tilde{I}$ such that
$$\forall\ i,j\in I:\qquad \tilde{m}_{\p(i)\p(j)}=m_{ij}\ .$$

\newglossaryentry{M_J}{type=symbols,name={\ensuremath{M_J}},description=submatrix with entries in $J$, sort=Coxeter}
\item Given a Coxeter matrix $M$ over $I$ and a subset $J\subseteq I$, we set $\gls{M_J}:= M_{|J\times J}$.
\end{itemize}
\end{de}

\begin{de}[\textbf{Coxeter Diagrams}]\ 
\begin{itemize}
\item A \textit{(2-spherical) Coxeter diagram}\index{Coxeter!diagram (2-spherical)} is a pair $(\GGG,\nu)$ consisting of a graph $\GGG$ and a map $\nu:E(\GGG)\to \NN_{\geq 3}$.
\item Given two Coxeter diagrams $(\GGG,\nu)$ and $(\tilde{\GGG},\tilde{\nu})$, a \textit{morphism $\mathit{\p:(\GGG,\nu)\to (\tilde{\GGG},\tilde{\nu})}$ of Coxeter diagrams}\index{morphism!of Coxeter diagrams} is a morphism $\varphi:\GGG\to \tilde{\GGG}$ of graphs such that
$$\forall\ \{i,j\}\in E(\GGG):\qquad \tilde{\nu}(\{ \p(i),\p(j)\})=\nu(\{i,j\})\ .$$
\end{itemize}
\end{de}

\begin{bem}\ 
\begin{enumerate}[label=(\alph*)]
\item Given a Coxeter diagram $(\GGG,\nu)$, we indicate an edge such that $\nu(\{i,j\})=3$ by a single edge, an edge $\{i,j\}$ such that $\nu(\{i,j\})=4$ by a double edge, an edge $\{i,j\}$ such that $\nu(\{i,j\})=6$ by a triple edge and an edge $\{i,j\}$ such that $\nu(\{i,j\})=8$ by a quadruple edge.

\newglossaryentry{Pi_M}{type=symbols,name={\ensuremath{\Pi_M}},description=Coxeter diagram with respect to the Coxeter matrix $M$, sort=Coxeter}
\newglossaryentry{V(M)}{type=symbols,name={\ensuremath{V(\GGG)}},description=set of vertices of Coxeter matrix $M$, sort=Coxeter}
\newglossaryentry{E(M)}{type=symbols,name={\ensuremath{E(\GGG)}},description=set of edges of the Coxeter matrix $M$, sort=Coxeter}
\newglossaryentry{A(M)}{type=symbols,name={\ensuremath{A(\GGG)}},description=set of directed edges of the Coxeter matrix $M$, sort=Coxeter}
\newglossaryentry{G(M)}{type=symbols,name={\ensuremath{G(\GGG)}},description=set of rank 3 subdiagrams of the Coxeter matrix $M$, sort=Coxeter}
\item Given a Coxeter matrix $M$ over $I$, the corresponding \textit{Coxeter diagram}\index{Coxeter!diagram} is $\gls{Pi_M}:=\big(\mathcal{G}_M, \nu_M\big)$ with  $V(\GGG_M):=I$ and
\begin{align*}
\forall\ i,j\in I:\ \{i,j\}\in E(\GGG_M)\ :\Leftrightarrow\ m_{ij}\geq 3\ , &&\forall\ \{i,j\}\in E(\GGG_M):\ \nu_M(\{i,j\}):=m_{ij}\ .
\end{align*}
We set
\begin{align*}
\gls{V(M)}:=V(\GGG_M)\ , && \gls{E(M)}:=E(\GGG_M)\ , && \gls{A(M)}:=A(\GGG_M)\ , && \gls{G(M)}:=G(\GGG_M)\ .
\end{align*}
\item Let $CM$ be the set of Coxeter matrices and let $CD$ be the set of Coxeter diagrams. Then the map
$$\Pi:CM\to CD,\ M\mapsto \Pi_M$$
is a bijection such that $M\cong \tilde{M}\ \Leftrightarrow\ \Pi_M\cong \Pi_{\tilde{M}}$.
\end{enumerate}
\end{bem}

\begin{de}[\textbf{Coxeter Systems}] Let $M$ be a Coxeter matrix over $I$.

\newglossaryentry{W_M}{type=symbols,name={\ensuremath{W_M}},description=Coxeter group with respect to the Coxeter matrix $M$, sort=Coxeter}
\newglossaryentry{W_Mr}{type=symbols,name={\ensuremath{(W_M,r)}},description=Coxeter system with respect to the Coxeter matrix $M$, sort=Coxeter}
\begin{itemize}
\item The \textit{Coxeter group of type $\mathit{M}$}\index{Coxeter!group} is the group $$\gls{W_M}:=\langle \{r_i\mid i\in I\} \mid \{(r_ir_j)^{m_{ij}}=1\mid i,j\in J\}\rangle\ .$$ 
\item The \textit{Coxeter system of type $\mathit{M}$}\index{Coxeter!system} is the pair \gls{W_Mr}, where $$r:\mathrm{Mon}(I)\to W_M,\ f\mapsto r_f$$
is the unique extension of the map $r:I\to W_M,\ i\mapsto r_i$ to a homomorphism from the free monoid $\mathrm{Mon}(I)$ on $I$ to $W_M$.
\end{itemize}
\end{de}

\chapter{Twin Buildings}

As the main issue of this thesis is the classification of twin buildings, we give a rough overview of the main concepts and the main results which allow us to pass from the whole building to its local structure without loss of information. By theorem \ref{454} and \ref{455} below, each irreducible residue of rank 2 of an irreducible twin building of rank at least 3 is a Moufang polygon. As a consequence, we may exploit the classification of Moufang polygons in \cite{TW}.

\begin{de}[\textbf{Chamber Systems}]\ 
\begin{itemize}
\item A \textit{chamber system over an index set $\mathit{I}$}\index{chamber!system} is a set $\Delta$ (whose elements are called \textit{chambers})\index{chamber} together with an equivalence relation $\sim_i$ on $\Delta$ (called \textit{$\mathit{i}$-equivalence}\index{$i$-equivalence}) for each $i\in I$.
\item Given a chamber system $\Delta$ over $I$ and $i\in I$, an \textit{$\mathit{i}$-panel}\index{$i$-panel}\index{panel!of type $i$} is an $i$-equivalence class, and a \textit{panel}\index{panel} is an $i$-panel for some $i\in I$.
\item Given a chamber system $\Delta$ over $I$ and $i\in I$, two distinct chambers $x,y\in \Delta$ such that $x\sim_i y$ are called \textit{$\mathit{i}$-adjacent}\index{$i$-adjacent}, and they are \textit{adjacent}\index{adjacent} if they are $i$-adjacent for some $i\in I$.

\newglossaryentry{xJy}{type=symbols,name={\ensuremath{x\sim_J y}},description=the chambers $x$ and $y$ are connected by a $J$-gallery, sort=Building}
\item Given a chamber system $\Delta$ over $I$, chambers $x,y\in \Delta$ and $J\subseteq I$, a \textit{$\mathit{J}$-gallery of length $\mathit{k}$ from $\mathit{x}$ to $\mathit{y}$}\index{$J$-gallery} is a sequence $\gamma=(x_0,\ldots,x_k)\subseteq \Delta^{k+1}$ for some $k\in \NN$ such that
\begin{align*}
x_0=x\ , && x_k=y\ , &&  \forall\ j\in \{1,\ldots,k\}\ \exists\ i_j\in J:\ x_{j-1}\sim_{i_j} x_j\ \wedge\ x_{j-1}\neq x_j\ ,
\end{align*}
a \textit{gallery from $\mathit{x}$ to $\mathit{y}$}\index{gallery} is an $I$-gallery from $x$ to $y$, and we write \gls{xJy} if there is a $J$-gallery from $x$ to $y$.

\newglossaryentry{d(x,y)}{type=symbols,name={\ensuremath{\mathrm{dist}(x,y)}},description=distance of the chambers $x$ and $y$, sort=Building}
\item Given a chamber system $\Delta$ and chambers $x,y\in \Delta$, the \textit{distance}\index{distance} \gls{d(x,y)} \textit{from $\mathit{x}$ to $\mathit{y}$} is the length of a shortest gallery from $x$ to $y$ if there is one and $\infty$ otherwise.
\item Given a chamber system $\Delta$ over $I$, chambers $x,y\in \Delta$ and a gallery $\gamma=(x_0,\ldots,x_k)$ from $x$ to $y$, the \textit{type of $\mathit{\gamma}$}\index{type of a gallery} is the word $i_1\cdots i_k\in \mathrm{Mon}(I)$.

\newglossaryentry{D(x)}{type=symbols,name={\ensuremath{\Delta_J(x)}},description=$J$-residue of the chamber $x$, sort=Building}
\item Given a chamber system $\Delta$ over $I$, a chamber $x\in \Delta$ and $J\subseteq I$, the \textit{$\mathit{J}$-residue of $\mathit{x}$}\index{$J$-residue!of chamber system} is
$$\gls{D(x)}:=\{ y\in \Delta\mid x\sim_J y \}\ .$$
A \textit{residue}\index{residue!of a chamber system} is a $J$-residue $\Delta_J(x)$ for some chamber $x\in \Delta$ and some $J\subset I$.
\end{itemize}
\end{de}

\begin{de}[\textbf{Buildings}]
Let $M$ be a Coxeter diagram over $I$ and let $(W_M,r)$ be the corresponding Coxeter system. A \textit{building of type $\mathit{M}$}\index{building} is pair $\BBB=(\Delta,\delta)$, where $\Delta$ is a chamber system over $I$ endowed with a function $\delta:\Delta\times \Delta\to W_M$ such that the following holds:
\begin{enumerate}[label=(B\arabic*)]
\item Each panel contains at least two chambers.
\item For each reduced word $f\in \mathrm{Mon}(I)$ and for each ordered pair $(x,y)$ of chambers, we have $\delta(x,y)=r_f$ if and only if there is a gallery of type $f$ from $x$ to $y$. \end{enumerate}
\end{de}

\begin{bem}
Cf. definition (39.10) of \cite{TW} for the definition of a reduced word.
\end{bem}

\begin{de}[\textbf{Standard Thin Buildings}]
Let $M$ be a Coxeter Matrix. Then the building $\Sigma(M):=(W_M,\delta_{W_M})$ with
$$\delta_{W_M}:W_M\times W_M\to W_M,\ (w_1,w_2)\mapsto w_1^{-1}w_2$$
is the \textit{standard thin building of type $\mathit{M}$}\index{building!standard thin}\index{standard thin building}.
\end{de}

\begin{bem}
Cf. example (5.7) of \cite{AB} that $\Sigma(M)$ is a building of type $M$.
\end{bem}

\begin{de}[\textbf{Apartments}]
Let $\BBB=(\Delta,\delta)$ be a building of type $M$ and let $X\subseteq W_M$.
\begin{itemize}
\item An \textit{isometry from $\mathit{X}$ to $\mathit{\BBB}$}\index{isometry} is a map $\pi:X\to \Delta$ such that
$$\forall\ x,y\in X:\qquad \delta(x^\pi,y^\pi)=x^{-1}y\ .$$
\item An \textit{apartment of $\mathit{\BBB}$}\index{apartment} is the image $\Sigma$ of some isometry $\pi:W_M\to \Delta$.
\end{itemize}
\end{de}

\begin{satz}[\textbf{$\boldsymbol{J}$-Residues}]
Let $M$ be a Coxeter matrix over $I$, let $\BBB=(\Delta,\delta)$ be a building of type $M$, let $x\in \Delta$ and let $J\subseteq I$. Then the following holds:

\newglossaryentry{B(x)}{type=symbols,name={\ensuremath{\BBB_J(x)}},description=the $J$-residue of the chamber $x$, sort=Building}
\begin{enumerate}[label=(\alph*)]
\item The \textit{$\mathit{J}$-residue}\index{$J$-residue!of a building}\index{residue!of a building} $$\BBB_J(x):=\big(\Delta_J(x),\delta_{|\Delta_J(x)\times \Delta_J(x)}\big)$$
is a building of type $M_J$.
\item \label{438} If $\Sigma$ is an apartment of $\BBB$ such that $\Sigma\cap \BBB_J(x)\neq \emptyset$, then $\Sigma_J:=\BBB_J(x)\cap \Sigma$
is an apartment of $\BBB_J(x)$.
\item If $\Sigma_J$ is an apartment of $\BBB_J(x)$, then we have $\Sigma_J=\BBB_J(x)\cap \Sigma$ for some apartment $\Sigma$ of $\BBB$.
\end{enumerate}
\end{satz}

\begin{bew}
This results from (39.52) of \cite{TW}.\qed
\end{bew}

\begin{de}[\textbf{Roots}]
Let $\BBB$ be a building of type $M$, let $\Sigma$ be an apartment of $\BBB$ and let $c$ be a chamber of $\Sigma$.
\begin{itemize}
\item A \textit{root of $\mathit{\Sigma}$}\index{root!simple} is a subset $\alpha\subset \Sigma$ such that
$$\alpha=\{ w\in \Sigma \mid \mathrm{dist}(w,x)<\mathrm{dist}(w,y) \}$$
for some ordered pair $(x,y)$ of adjacent chambers. We denote the set of roots of $\Sigma$ by $\Phi(\BBB,\Sigma)$.
\item A \textit{root of $\mathit{\BBB}$} is a root of some apartment $\Sigma\subseteq \BBB$.
\item Given $i\in I$, the \textit{simple root $\alpha_i$ with respect to $(\Sigma,c)$} is the root $$\alpha_i:=\{ w\in \Sigma \mid \mathrm{dist}(w,c)<\mathrm{dist}(w,c_i) \}\ , $$
where $c_i$ is the unique chamber of $\Sigma$ which is $i$-adjacent to $c$. We write $\Phi(\BBB,\Sigma,c)$ instead of $\Phi(\BBB,\Sigma)$ if we additionally take the simple roots with respect to $(\Sigma,c)$ into account.
\end{itemize}
\end{de}

\begin{de}[\textbf{Standard Root Systems}]
Let $M$ be a Coxeter Matrix.
\begin{itemize}
\item The set $\Phi(M):=\Phi\big(\Sigma(M),\Sigma(M),1_{W_M}\big)$ of roots of $\Sigma(M)$ is the \textit{standard root system of type $\mathit{M}$}\index{root system!standard}\index{standard root system}.
\item Given $\alpha,\beta\in \Phi(M)$, the pair $\{\alpha,\beta\}$ is \textit{prenilpotent}\index{prenilpotent} if we have
$$\alpha\cap \beta\neq \emptyset \neq (-\alpha)\cap(-\beta)\ .$$
In this case, we set
\begin{align*}
[\alpha,\beta]:=\{ \gamma\in \Phi(M) \mid \alpha\cap \beta\subseteq \gamma,\ (-\alpha)\cap(-\beta)\subseteq -\gamma \}\ , && (\alpha,\beta):=[\alpha,\beta]\sm \{\alpha,\beta\}\ .
\end{align*}
\end{itemize}
\end{de}

\begin{bem}
Let $M$ be a Coxeter matrix over $I$, let $\Phi(M)$ be the standard root system of type $M$ and let $\alpha\in \Phi(M)$. Then we have $$\alpha=v\alpha_i=\{ v\cdot w \mid w\in \alpha_i\}$$ for some $i\in I$ and some $v\in W_M$, cf. proposition (5.81) of \cite{AB}.
\end{bem}

\begin{de}[\textbf{Twin Buildings}]
Let $M$ be a Coxeter diagram over $I$. A \textit{twin building of type $\mathit{M}$}\index{twin building} is a triple $\BBB=(\BBB_+,\BBB_-,\delta^*)$, where each \textit{half}\index{half of a twin building} $\BBB_\epsilon=(\Delta_\epsilon,\delta_\epsilon)$ with $\epsilon\in \{\pm\}$ is a building of type $M$ and
$$\delta^*:(\Delta_+\times \Delta_-)\cup (\Delta_-\times \Delta_+)\to W_M$$
is a \textit{codistance}\index{codistance}, i.e., given $\epsilon\in \{\pm\}$, $x\in \Delta_\epsilon$, $y\in \Delta_{-\epsilon}$ and $w:=\delta^*(x,y)$, the following holds:
\begin{enumerate}[label=(C\arabic*)]
\item We have $\delta^*(y,x)=w^{-1}$.
\item Given $z\in \Delta_{-\epsilon}$, $i\in I$ such that $\delta_{-\epsilon}(y,z)=r_i$ and $l(wr_i)=l(w)-1$, we have $\delta^*(x,z)=wr_i\ .$
\item Given $i\in I$, there exists a chamber $z\in \Delta_{-\epsilon}$ such that $\delta_{-\epsilon}(y,z)=r_i$ and $\delta^*(x,z)=wr_i$.
\end{enumerate}
Here $l:W\to \NN^*$ is the length function with respect to the set $\{r_i \mid i\in I\}$ of generators.
\end{de}

\begin{de}[\textbf{Opposite Chambers}]
Let $M$ be a Coxeter diagram over $I$, let $\BBB$ be a twin building of type $M$, let $J\subseteq I$ and let $\epsilon\in \{\pm\}$.
\begin{itemize}
\item Two chambers $x\in \BBB_\epsilon$ and $y\in \BBB_{-\epsilon}$ such that $\delta^*(x,y)=1$ are called \textit{opposite}\index{opposite!chambers}\index{chamber!opposite}. We set
\newglossaryentry{Op}{type=symbols,name={\ensuremath{\OOO_B}},description=set of opposite chambers, sort=Building}
$$\gls{Op}:=\{ (x,y)\in \BBB_+\times \BBB_- \mid \delta^*(x,y)=1 \}\ .$$
\item Two residues $\RRR_+\subseteq \BBB_+$ and $\RRR_-\subseteq \BBB_-$ such that $$\RRR_+\times \RRR_{-\epsilon}\cap \OOO_B\neq \emptyset$$ are called \textit{opposite}\index{residue!opposite}.
\end{itemize}
\end{de}

\begin{lemma}\label{435}
Let $\BBB$ be a twin building, let $\epsilon\in \{\pm\}$ and let $x\in \BBB_\epsilon$. Then there exists a chamber $y\in \BBB_{-\epsilon}$ such that $\delta^*(x,y)=1$.
\end{lemma}

\begin{bew}
This results from corollary (5.141) of \cite{AB}.\qed
\end{bew}

\begin{satz}[\textbf{$\boldsymbol{J}$-Residues}]\label{436}
Let $M$ be a Coxeter diagram over $I$, let $\BBB$ be a twin building of type $M$, let $J\subseteq I$, let $(x,y)\in \OOO_B$ and let $\BBB_J(x):=(\BBB_+)_J(x)$, $\BBB_J(y):=(\BBB_-)_J(y)$. Then the \textit{$\mathit{J}$-residue}\index{$J$-residue!of a twin building}\index{residue!of a twin building}
\newglossaryentry{B(xy)}{type=symbols,name={\ensuremath{\BBB_J(x,y)}},description=$J$-residue of the pair ${(x,y)}$ of opposite chambers, sort=Building}
$$\gls{B(xy)}:=\big( \BBB_J(x),\BBB_J(y),\delta^*_{| (\BBB_J(x)\times \BBB_J(y))\cup (\BBB_J(y)\times \BBB_J(y))} \big)$$
is a twin building of type $M_J$.
\end{satz}

\begin{bew}
By lemma (5.148) of \cite{AB}, we have
$$\delta(\bar{x},\bar{y})\in W_{M_J}\delta^*(x,y) W_{M_j}=W_{M_J}\cdot 1\cdot W_{M_J}=W_{M_J}$$
for all $\bar{x}\in \BBB_J(x)$, $\bar{y}\in \BBB_J(y)$.\qed
\end{bew}

\begin{kor}\label{437}
Let $M$ be a Coxeter diagram over $I$, let $\BBB$ be a twin building of type $M$, let $\epsilon\in \{\pm\}$, let $x\in \BBB_\epsilon$ and let $J\subseteq I$. Then $(\BBB_\epsilon)_J(x)$ is the half of a twin building.
\end{kor}

\begin{bew}
By lemma \ref{435}, there is a chamber $y\in \BBB_{-\epsilon}$ such that $\delta^*(x,y)=1$, thus $(\BBB_\epsilon)_J(x)$ is the half of a twin building by theorem \ref{436}.\qed
\end{bew}

\begin{no}
Let $M$ be a Coxeter matrix over $I$ and let $\BBB$ be a building of type $M$. Given a chamber $c\in \BBB$, we define
$$E_2(c):=\{ \BBB_{\{i,j\}}(c) \mid i\neq j\in I \}\ .$$
\end{no}

\begin{de}[\textbf{Twin Apartments}]
Let $\BBB$ be a twin building. 
\begin{itemize}
\item A \textit{twin apartment of $\mathit{\BBB}$}\index{twin apartment} is a pair $\Sigma=(\Sigma_+,\Sigma_i)$ of apartments $\Sigma_\epsilon$ of $\BBB_\epsilon$ such that each chamber of $\Sigma_+\cup \Sigma_-$ is opposite precisely one other chamber $\mathrm{op}_\Sigma(c) \in \Sigma_+\cup \Sigma_-$.
\item Given a twin apartment $\Sigma$, the map 
$$\mathrm{op}_\Sigma:\Sigma_+\cup \Sigma_-\to \Sigma_+\cup \Sigma_-,\ c\mapsto \mathrm{op}_\Sigma(c)$$
is the \textit{opposition involution with respect to $\Sigma$}\index{opposition involution}.
\end{itemize}
\end{de}

\begin{de}[\textbf{Twin Roots}]
Let $\BBB$ be a twin building and let $\Sigma$ be a twin apartment of $\BBB$.
\begin{itemize}
\item A \textit{twin root of $\mathit{\Sigma}$}\index{twin root} is a pair $\alpha=(\alpha_+,\alpha_-)$ of roots $\alpha_\epsilon$ of $\Sigma_\epsilon$ such that
$$\mathrm{op}_\Sigma(\alpha)=-\alpha\ .$$
\item We denote the set of twin roots of $\Sigma$ by $\Phi(\BBB,\Sigma)$.
\end{itemize}
\end{de}

\begin{bem}\label{493}
Let $\BBB$ be a twin building. Given a twin apartment $\Sigma$ of $\BBB$ and a root $\alpha_+$ of $\Sigma_+$, then $\alpha_-:=-\mathrm{op}_\Sigma(\alpha_+)$ is the unique root of $\Sigma_-$ such that $\alpha:=(\alpha_+,\alpha_-)$ is a twin root of $\Sigma$. As a consequence, the map 
$$f:\Phi(\BBB,\Sigma)\to \Phi(\BBB_+,\Sigma_+),\ (\alpha_+,\alpha_-)\to \alpha_+$$
is a bijection.
\end{bem}

\begin{de}[\textbf{Isometries and Automorphisms}]
Let $\BBB=(\BBB_+,\BBB_-,\delta^*)$ be a twin building of type $M$ and let $\tilde{\BBB}=(\tilde{\BBB}_+,\tilde{\BBB}_-,\tilde{\delta}^*)$ be a twin building of type $\tilde{M}$.
\begin{itemize}
\item An \textit{isometry of twin buildings}\index{isometry!of twin buildings} is a triple $\phi=(\sigma,\phi_+,\phi_-)$ consisting of an isomorphism $\sigma:M\to \tilde{M}$ of Coxeter diagrams and maps $\phi_\epsilon:\BBB_\epsilon\to \BBB_\epsilon$ such that
\begin{align*}
\forall\ c_\epsilon,d_\epsilon\in \BBB_\epsilon: \qquad \tilde{\delta}_\epsilon\big(\phi_\epsilon(c_\epsilon),\phi_\epsilon(d_\epsilon)\big)=\sigma\big(\delta_\epsilon(c_\epsilon,d_\epsilon)\big)
\end{align*}
and 
\begin{align*}
\forall\ c_\epsilon\in \BBB_\epsilon,c_{-\epsilon}\in \BBB_{-\epsilon}: \qquad \tilde{\delta}^*\big(\phi_\epsilon(c_\epsilon),\phi_{-\epsilon}(c_{-\epsilon})\big)=\sigma\big(\delta^*(c_\epsilon,c_{-\epsilon})\big)
\end{align*}
\item An \textit{isomorphism of twin buildings}\index{isomorphism!of twin buildings} is a surjective isometry $\phi:\BBB\to \tilde{\BBB}$. 
\item An \textit{automorphism of $\mathit{\BBB}$}\index{automorphism!of twin buildings} is an isomorphism $\phi:\BBB\to \BBB$. We denote the group of automorphisms of $\BBB$ by $\Aut(\BBB)$.
\item An automorphism $\phi\in \Aut(\BBB)$ is \textit{special} if we have $\sigma=\id_M$. We denote the group of special automorphisms of $\BBB$ by $\Aut_0(\BBB)$.
\end{itemize}
\end{de}

\begin{de}[\textbf{Strongly Transitive Actions}]
Let $\BBB$ be a twin building and let $G$ be a group.
\begin{itemize}
\item An \textit{action of $\mathit{G}$ on $\mathit{\BBB}$}\index{action} is a homomorphism 
$$\p:G\to \Aut(\BBB)\ .$$
\item An action $\p:G\to \Aut(\BBB)$ is \textit{strongly transitive}\index{action!strongly transitive} if the action is transitive on the set
$$\{ (\Sigma,c) \mid \Sigma\ \textrm{is a twin apartment of $\BBB$},\  c\in \OOO_\Sigma \}\ .$$
\end{itemize}
\end{de}

\begin{satz}\label{497}
Let $\BBB$ be a Moufang twin building. Then $\Aut_0(\BBB)$ acts strongly transitively on $\BBB$.
\end{satz}

\begin{bew}
This results from proposition (8.19) of \cite{AB}.\qed
\end{bew}

\newglossaryentry{Extension}{type=results,name={{Extension Theorem}},description={},sort=res}

\begin{satz}[\textbf{\gls{Extension}}]\label{498}
Let $\BBB$, $\tilde{\BBB}$ be thick (2-spherical) twin buildings of the same type. Assume that $\BBB$ and $\tilde{\BBB}$ satisfy condition (CO). Given $c\in \OOO_\BBB$ and $\tilde{c}\in \OOO_{\tilde{\BBB}}$ and a surjective isometry
$$\phi: E_2(c_+)\cup\{c_-\}\to E_2(\tilde{c}_+)\cup \{\tilde{c}_-\}\ ,$$
there is a unique extension of $\phi$ to an isomorphism $\phi:\BBB\to \tilde{\BBB}$ of twin buildings. (A building is \textit{thick}\index{thick} if each panel contains at least three chambers.)
\end{satz}

\begin{bew}
Cf. theorem (5.213) of \cite{AB}. Notice that our buildings are 2-spherical by definition.\qed
\end{bew}

\begin{de}[\textbf{Root Groups and the Moufang Property}] Let $\BBB$ be a twin building of rank at least 2, where the \textit{rank of a building of type $\mathit{M}$}\index{rank} is just $|V(M)|$.
\begin{itemize}
\item Given a twin root $\alpha$, the corresponding \textit{root group}\index{root group} is
\newglossaryentry{Ua}{type=symbols,name={\ensuremath{U_\alpha}},description=root group with respect to the root $\alpha$, sort=Building}
$$\gls{Ua}:=\{ g\in \Aut(\BBB)\mid g\ \textrm{acts trivially on each panel of}\ \alpha^o\}\ ,$$
where $\alpha^o$ is the set of all panels of $\BBB$ which contain at least two chambers in $\alpha$.
\item The building $\BBB$ is \textit{Moufang}\index{Moufang property} if it is thick and if for each root $\alpha$ of $\BBB$, the root group $U_\alpha$ acts transitively on the set of twin apartments containing $\alpha$.
\item The building $\BBB$ is \textit{strictly Moufang}\index{Moufang property!strict} if the actions of the root groups are simply transitive.
\end{itemize}
\end{de}

\begin{bem}
A Moufang twin building whose Coxeter diagram has no isolated nodes is strictly Moufang, cf. p. 455 of \cite{AB}.
\end{bem}

\begin{satz}\label{454} Every thick, irreducible twin building of rank at least 3 that satisfies condition (CO) is Moufang.
\end{satz}

\begin{bew}
This is theorem (8.27) of \cite{AB}. Notice that our buildings are 2-spherical by definition.\qed
\end{bew} 

\begin{satz}\label{455} Let $\BBB$ be a Moufang twin building. Then every spherical residue of $\BBB$ is Moufang.
\end{satz}

\begin{bew}
This is proposition (8.21) of \cite{AB}.
\end{bew}

\begin{kor}\label{494}
Let $\BBB$ be a thick, irreducible twin building of rank at least 3. Then each residue of rank 2 is also Moufang. In particular, the irreducible residues of rank 2 are Moufang polygons.
\end{kor}

\begin{bew}
This results from remark (8.30)(a) of \cite{AB}. Notice that our buildings are 2-spherical by definition.\qed
\end{bew}

\begin{de}[\textbf{RGD System}]
Let $M$ be a Coxeter matrix, let $\Phi:=\Phi(M)$, let $G$ be a group, let $(U_\alpha)_{\alpha\in \Phi}$ be a family of non-trivial subgroups of $G$ and let
$$T:=\bigcap_{\mathclap{\alpha\in \Phi}} N_G(U_\alpha)\ .$$
Then the pair $\big(G,(U_\alpha)_{\alpha\in \Phi}\big)$ is an \textit{RGD system of type $\mathit{M}$}\index{RGD system} if the following holds:
\begin{enumerate}[label=(RGD\arabic*), leftmargin=40pt]
\item We have
$$[U_\alpha,U_\beta]\leq U_{(\alpha,\beta)}$$
for all $\alpha\neq \beta$ such that $\{\alpha,\beta\}$ is prenilpotent.
\item Given $i\in I$, there is a function $\mu:U_{\alpha_i}^*\to G$ such that we have
\begin{align*}
\forall\ u\in U_{\alpha_i}^*:\ \mu(u)\in U_{-\alpha_i} u U_{-\alpha_i}\ , && \forall\ u\in U_{\alpha_i}^*,\ \alpha\in \Phi:\ \mu(u)U_{\alpha_i}\mu(u)^{-1}=U_{r_i\alpha_i}\ .
\end{align*}
\item Given $i\in I$, we have
$$U_{-\alpha_i}\not\leq U_+\ ,$$
where 
\begin{align*}
U_{\epsilon}:=\langle U_\alpha \mid \alpha\in \Phi_\epsilon\rangle\ , && \Phi_+:=\{ \alpha\in \Phi \mid 1\in \alpha\}\ , && \Phi_-:=\{\alpha\in \Phi \mid 1\notin \alpha\}\ .
\end{align*}
\item We have $$G=T\langle U_\alpha \mid \alpha\in \Phi\rangle\ .$$
\end{enumerate}
\end{de}

\begin{satz} The following holds:
\begin{enumerate}[label=(\alph*)]
\item Each RGD system $\big(G,(U_\alpha)_{\alpha\in \Phi}\big)$ gives rise to a twin building $\BBB\big(G,(U_\alpha)_{\alpha\in \Phi}\big)$.
\item Let $\BBB$ be a twin building of type $M$, let $\Sigma$ be a twin apartment of $\BBB$, let $c\in \OOO_\Sigma$, let $G:=\Aut(\BBB)$, let $\Phi:=\Phi(\BBB,\Sigma,c)$ and let $(U_\alpha)_{\alpha\in \Phi}$ be the family of root groups with respect to $(\Sigma,c)$. Then the pair $\big(G,(U_\alpha)_{\alpha\in \Phi}\big)$ is an RGD system of type $M$, and $\BBB$ is uniquely determined by this RGD system, i.e., we have $$\BBB\cong \BBB\big(G,(U_\alpha)_{\alpha\in \Phi}\big)\ .$$
\end{enumerate}
\end{satz}

\begin{bewzwei}\ 
\begin{enumerate}[label=(\alph*)]
\item This is theorem (8.81) of \cite{AB}.
\item The first statement is example (8.47)(a) of \cite{AB} while the second one results from theorem (8.9) of \cite{AB}.
\end{enumerate}
\qed
\end{bewzwei}

\begin{prop}\label{495}
Let $M$ be a Coxeter matrix over $I$, let $\BBB$ be a strictly Moufang twin building of type $M$, let $\Sigma$ be a twin apartment of $\BBB$ and let $c\in \OOO_\Sigma$. Let $J\subseteq I$ be such that $M_J$ has no isolated nodes, let $\BBB_J:=\BBB_J(c)$ and let $\alpha_J$ be a root of $\Sigma_J=\Sigma\cap \BBB_J$. Then there is a unique root $\alpha$ of $\Sigma$ such that $\alpha_J=\alpha\cap \Sigma_J$, and the restriction map
$$\rho: U_\alpha\to U_{\alpha_J}$$
is an isomorphism of groups.
\end{prop}

\begin{bew}
For spherical buildings this is proposition (7.32) of \cite{AB}, and the arguments given in its proof go through in the twin case without much change. \qed
\end{bew}

\newpage

\begin{satz}\label{453}
Let $M$ be an irreducible Coxeter matrix over $I$ such that $|I|\geq 3$, let $\BBB$ be a thick Moufang twin building of type $M$, let $\Sigma$ be a twin apartment of $\BBB$, let $c\in \OOO_\Sigma$ and let $\Phi:=\Phi(\BBB,\Sigma,c)$. Let $(i,j)\in A(M)$, $\BBB_{ij}:=\BBB_{\{i,j\}}(c)$ and $n:=m_{ij}$. Then the following holds:

\begin{enumerate}[label=(\alph*)]
\item The residue $\BBB_{ij}$ is a Moufang $n$-gon.
\item The intersection $\Sigma_{ij}:=\Sigma\cap \BBB_{ij}$ is an apartment of $\BBB_{ij}$, and the roots $\alpha_i\cap \BBB_{ij}$ and $\alpha_j\cap \BBB_{ij}$ form a root basis of $\BBB_{ij}$.
\item Let
$$( \bar{\omega}_1=\alpha_i\cap \BBB_{ij},\bar{\omega}_2,\ldots,\bar{\omega}_{n-1}, \bar{\omega}_n=\alpha_j\cap \BBB_{ij} )$$
be the root sequence of $\BBB_{ij}$ from $\alpha_i\cap \BBB_{ij}$ to $\alpha_j\cap \BBB_{ij}$. Then there are exactly $n$ roots $\omega_1=\alpha_i,\omega_2,\ldots,\omega_n=\alpha_j$ of $\Sigma$ such that
$$\forall\ 1\leq i\leq n:\qquad \bar{\omega}_i=\omega_i\cap \BBB_{ij}\ .$$
\item \label{499} For $i=1,\ldots,n$ let $U_i:=U_{\omega_i}$, let $U_{[\alpha_i,\alpha_j]}:=U_1\cdots U_n$ and let
$$\Theta_{(i,j)}:=\big(U_{[\alpha_i,\alpha_j]},U_1,\ldots,U_n\big)\ .$$
Then $\Theta_{(i,j)}$ is isomorphic to the root group sequence of $\BBB_{ij}$ from $\alpha_i \cap \BBB_{ij}$ to $\alpha_j\cap \BBB_{ij}$.
\end{enumerate}
\end{satz}

\begin{bewzwei}\ 
\begin{enumerate}[label=(\alph*)]
\item This is corollary \ref{494}.
\item The first assertion results from theorem \ref{438}, and by definition, the roots $\alpha_i\cap \BBB_{ij}$ and $\alpha_j\cap \BBB_{ij}$ are simple roots of $\Sigma_{ij}$.
\item This results from proposition \ref{495}.
\item This results from proposition \ref{495}.
\end{enumerate}\qed
\end{bewzwei} 

\begin{de}[\textbf{Double $\boldsymbol{\mu}$-Maps}]
\newglossaryentry{muab}{type=symbols,name={\ensuremath{h_{a,b}}},description=double $\mu$-map with respect to $a$ and $b$, sort=Building}
Let $M$ be a Coxeter matrix over $I$, let $\BBB$ be a Moufang twin building of type $M$, let $\Sigma$ be a twin apartment of $\BBB$, let $c\in \OOO_\Sigma$ and let $\Phi:=\Phi(\BBB,\Sigma,c)$. Let $i\in I$, let $a,b\in U_{\alpha_i}^*$ and let $\mu(a),\mu(b)$ be as in (RGD2). Then the map
$$\gls{muab}:\langle U_\alpha \mid \alpha\in \Phi \rangle \to \langle U_\alpha \mid \alpha\in \Phi \rangle, u\mapsto \mu(a)^{-1}\mu(b) u\mu(b)^{-1}\mu(a)$$
is the \textit{double $\mathit{\mu}$-map with respect to $\mathit{a,b}$}\index{double $\mu$-map}\index{$\mu$-map}.
\end{de}

\begin{bem}\ 
\begin{enumerate}[label=(\alph*)]
\item We have $h_{a,b}\in T$ for all $a,b\in U_{\alpha_i}^*$.
\item \label{496} By (3) in section (7.8.2) of \cite{AB}, the $\mu$-maps in (RGD2) are uniquely determined.
\end{enumerate}
\end{bem}

\begin{satz}\label{452}
Let $M$ be an irreducible Coxeter matrix over $I$ such that $|I|\geq 3$, let $\BBB$ be a thick Moufang twin building of type $M$, let $\Sigma$ be a twin apartment of $\BBB$, let $c\in \OOO_\Sigma$ and let $\Phi:=\Phi(\BBB,\Sigma,c)$. Let $(i,j)\in A(M)$, let $\BBB_{ij}:=\BBB_{\{i,j\}}(c)$, let $\Sigma_{ij}:=\Sigma\cap \BBB_{ij}$ and let $\Phi_{ij}:=\Phi(\BBB_{ij},\Sigma_{ij},c)$. Then we have
$$\forall\ a\in U_{\alpha_i}^*:\qquad \mu^{\BBB_{ij}}\big(\rho(a)\big)=\rho\big( \mu^\BBB(a) \big)\ ,$$
where $\rho: U_{[\alpha_i,\alpha_j]}\to \Aut(\BBB_{ij})$ is the restriction homomorphism. In particular, we have
$$\forall\ a,b\in U_{\alpha_i}^*:\qquad h^{\BBB_{ij}}_{\rho(a),\rho(b)}=\rho\circ h^\BBB_{a,b}\circ \rho^{-1}\ .$$
\end{satz}

\begin{bew}
Given a root $\bar{\alpha}\in \Phi_{ij}$, let ${\alpha}\in \Phi$ be the unique root of $\Sigma$ such that $\bar{\alpha}={\alpha}\cap \Sigma_{ij}$. By proposition \ref{495}, the map $\rho:U_{{\alpha}}\to U_{\bar{\alpha}} $ is an isomorphism of groups for each root $\bar{\alpha}\in \Phi_{ij}$. Given $a\in U_{\alpha_i}^*$, we have
$$\rho\big(\mu^\BBB(a)\big)\in \rho( U_{-\alpha_i} aU_{-\alpha_i} )=U_{-\bar{\alpha}_i} \rho(a)U_{-\bar{\alpha}_i}$$
and
\begin{align*}\forall\ \alpha\in \Phi_{ij}:\qquad \rho\big(\mu^\BBB(a)\big) U_{\bar{\alpha}}\rho\big(\mu^\BBB(a)\big)^{-1}=\rho\big( \mu^\BBB(a)U_{{\alpha}}\mu^\BBB(a)^{-1}\big)=\rho(U_{r_i{\alpha}})=U_{r_i\bar{\alpha}}\ .
\end{align*}
and thus
$$\mu^{\BBB_{ij}}\big(\rho(a)\big)=\rho\big( \mu^\BBB(a) \big)$$
by remark \ref{496}. 
\qed
\end{bew}

\newpage

\addtocontents{toc}{\noindent\protect\mbox{}\protect\hrulefill\par}
\part{Parameter Systems}
\addtocontents{toc}{\noindent\protect\mbox{}\protect\hrulefill\par}

\noindent In this part, we introduce the parameter systems which arise in the description of Moufang triangles and the six families of Moufang quadrangles. Moreover, we collect the basic results which will be needed for the classification of twin buildings, and, closely related, the solution of the isomorphism problem for Moufang sets.

For a detailed reference on these subjects, see \cite{TW}. Concerning alternative rings, we additionally refer to \cite{Sch}.\\
 
\chapter{Alternative Rings}

Alternative division rings are the parametrizing structures for Moufang triangles, the building bricks for simply laced twin buildings.

\section{Basic Definitions and Basic Properties}

\begin{de} An \textit{alternative ring}\index{alternative ring} is a triple $(\AA,+,\cdot)$ such that the following holds:
\begin{enumerate}[label=(A\arabic*)]
\item The pair $(\AA,+)$ is a commutative group.
\item The multiplication $\cdot:\AA\times \AA\to \AA$ is biadditive.

\newglossaryentry{ass}{type=symbols,name={\ensuremath{{[x,y,z]}}},description=associator of the elements ${x,y,z}$, sort=alternative ring}
\item The multiplication $\cdot:\AA\times \AA\to\AA$ is \textit{alternative}, i.e., it satisfies
$$\forall\ x,y\in \AA:\qquad [x,x,y]=0_\AA=[y,x,x]\ .$$
where $\gls{ass}:=(xy)z-x(yz)$ is the \textit{associator of $\mathit{x,y,z\in\AA}$}\index{associator}
\item There is an identity element $1_\AA$.
\end{enumerate}
\end{de}

\begin{lemma}\label{148}
An alternative ring $\AA$ is \textit{flexible}\index{flexible}, i.e., given $x,y\in \AA$, we have
$$[x,y,x]=0_\AA\ .$$
\end{lemma}

\begin{bew}
Cf. page 27 of \cite{Sch}.\qed
\end{bew}

\newglossaryentry{Moufang}{type=results,name={{Moufang Identities}},description={},sort=res}

\begin{lemma}[\textbf{\gls{Moufang}}] An alternative ring $\AA$ satisfies the \textit{Moufang identities}, i.e., given $x,y,z\in \AA$, we have
\begin{align*} 
(xyx)z=x\big(y(xz)\big)\ , && z(xyx)=\big((zx)y\big)x\ , && (xy)(zx)=x(yz)x\ .
\end{align*}
\end{lemma}

\begin{bew}
Cf. page 28 of \cite{Sch}.\qed
\end{bew}

\begin{de} An alternative ring $\AA$ is an \textit{alternative division ring}\index{alternative division ring} if the maps
\begin{align*}
\rho_w:\AA\to\AA,\ x\mapsto xw\ , && \lambda_w:\AA\to \AA,\ x\mapsto wx
\end{align*}
are bijective for each $w\in \AA^*$.
\end{de}

\begin{bem}
Let $\AA$ be an alternative division ring. Given $x\in \AA^*$, there are unique elements $x^{-l},x^{-r}\in \AA^*$ such that
\begin{align*}
x^{-l}\cdot x=1_\AA=x\cdot x^{-r}\ .
\end{align*}
By lemma \ref{148}, we have
\begin{align*}
\lambda_{x^{-l}}(xx^{-l})=x^{-l}\cdot xx^{-l}=x^{-l}x \cdot x^{-l}=1_\AA\cdot x^{-l}=x^{-l}=\lambda_{x^{-l}}(1_\AA)\ , && xx^{-l}=1_\AA\end{align*}
and therefore $x^{-1}:= x^{-l}=x^{-r}$.
\end{bem}

\begin{lemma}\label{124} An alternative division ring $\AA$ has the \textit{inverse properties}\index{inverse properties}, i.e., given $x,y\in \AA^*$, we have
\begin{align*} 
x^{-1}(xy)=y\ , && (yx)x^{-1}=y\ , && (xy)^{-1}=y^{-1}x^{-1}\ .
\end{align*}
\end{lemma}

\begin{bew}
This results from the Moufang identities.\qed
\end{bew}

\begin{de}\label{355}
Let $(\AA,+,\cdot)$ be an alternative ring. Then $(\AA,+,\circ)$ with the multiplication
$$\circ:\AA\times \AA\to \AA,\ x\circ y:=y\cdot x$$
is the \textit{opposite alternative ring}\index{alternative ring!opposite}\index{opposite!alternative ring}.
\end{de}

\begin{de}\label{467} Let $\AA,\tilde{\AA}$ be alternative rings. 
\begin{itemize}
\item An \textit{(anti-)isomorphism of alternative rings}\index{isomorphism!of alternative rings}\index{anti-isomorphism!of alternative rings} is an additive isomorphism $\gamma:\AA\to \tilde{\AA}$ such that
\begin{align*}
\forall\ x,y\in \AA:\qquad \gamma(xy)=\gamma(x)\gamma(y)\qquad \big(\ \gamma(xy)=\gamma(y)\gamma(x)\ \big)\ .
\end{align*} 
\item A \textit{Jordan homomorphism} is an additive monomorphism $\gamma:\AA\to \tilde{\AA}$ such that 
\begin{align*}
\gamma(1_{\AA})=1_{\tilde{\AA}}\ , && \forall\ x,y\in \AA:\ \gamma(xyx)=\gamma(x)\gamma(y)\gamma(x)\ .
\end{align*}
\end{itemize}
\end{de}

\begin{no} Let $\AA$ be an alternative division ring.
\begin{itemize}
\newglossaryentry{auta}{type=symbols,name={\ensuremath{\Aut(\AA)}},description=group of automorphisms of the alternative division ring $\AA$, sort=alternative ring}
\item We denote the \textit{group of automorphisms of $\mathit{\AA}$} by \gls{auta}. 
\newglossaryentry{aauta}{type=symbols,name={\ensuremath{\Aut^o(\AA)}},description=set of anti-automorphisms of the alternative division ring $\AA$, sort=alternative ring}
\item We denote the \textit{set of anti-automorphisms of $\mathit{\AA}$} by \gls{aauta}.
\newglossaryentry{leftw}{type=symbols,name={\ensuremath{\lambda_w}},description=left multiplication with the element $w$, sort=alternative ring}
\newglossaryentry{rightw}{type=symbols,name={\ensuremath{\rho_w}},description=right multiplication with the element $w$, sort=alternative ring}
\newglossaryentry{conjw}{type=symbols,name={\ensuremath{\gamma_w}},description=conjugation with the element $w$, sort=alternative ring}
\item Given $w\in \AA$, we set
\begin{align*}
\gls{leftw}:\AA\to\AA,\ x\mapsto wx\ , &&\gls{rightw}:\AA\to \AA,\ x\mapsto xw\ , && \gls{conjw}:\AA\to \AA,\ x\mapsto w^{-1}xw\ .
\end{align*}
Notice that the conjugation map $\gamma_w$ is well-defined by (9.23)(ii) of \cite{TW} with $c:=a$.
\newglossaryentry{ao}{type=symbols,name={\ensuremath{\AA^o}},description=opposite of the alternative division ring $\AA$, sort=alternative ring}
\item We denote the \textit{opposite alternative division ring} by \gls{ao}.
\newglossaryentry{jora}{type=symbols,name={\ensuremath{\Aut_J(\AA)}},description=group of Jordan automorphisms of the alternative division ring $\AA$, sort=alternative ring}
\item We denote the \textit{group of Jordan automorphisms of $\AA$} by \gls{jora}.
\end{itemize}
\end{no}

\begin{de}
\newglossaryentry{ZA}{type=symbols,name={\ensuremath{Z(\AA)}},description=center of the alternative ring $\AA$, sort=alternative ring}
\newglossaryentry{comm}{type=symbols,name={\ensuremath{[x,y]}},description=commutator of the elements ${x,y}$, sort=alternative ring}
The \textit{center}\index{center} of an alternative ring $\AA$ is
$$\gls{ZA}:=\{ x\in \AA \mid [x,\AA,\AA]=[x,\AA]=0_\AA\}\ ,$$
where $\gls{comm}:=xy-yx$ is the \textit{commutator of $\mathit{x,y\in\AA}$}\index{commutator}.
\end{de}

\begin{lemma} 
Let $\AA$ be an alternative division ring. Then $\KK:=Z(\AA)$ is a field and $\AA$ is an algebra over $\KK$.
\end{lemma}

\begin{bew}
This results from (9.18) and (9.23) of \cite{TW}.\qed
\end{bew}

\section{Octonion Division Algebras}
The Bruck-Kleinfeld theorem states that a non-associative alternative division ring is an octonion division algebra. First of all we give the exact definition of such an algebra before we collect some basic concepts, including the doubling process.

\begin{bem}
The construction here is taken from \cite{TW}.
\end{bem}

\begin{de} \newglossaryentry{SQE}{type=symbols,name={\ensuremath{\EE/\KK}},description=separable quadratic extension, sort=alternative ring}
Let \gls{SQE} be a separable quadratic extension and let $\sigma$ be the non-trivial element of $\mathrm{Gal}(\EE/\KK)$. 
\begin{itemize}
\item Given $x\in \EE$, we write $\bar{x}:=\sigma(x)$.
\item We denote the norm map and the trace map of $\EE/\KK$ by $N$ and $T$, respectively.

\newglossaryentry{Quat}{type=symbols,name={\ensuremath{(\EE/\KK,\beta)}},description=quaternion division algebra with respect to $\EE/\KK$ and $\beta\in \KK^*$, sort=alternative ring}
\item Given $\beta\in \KK^*$, we set
\begin{align*}\gls{Quat}:=\left\{ \begin{pmatrix} x & \beta \bar{y} \\y & \bar{x} \end{pmatrix} \mid x,y\in \EE\right\}\subseteq M_2(\EE)\ , && e:=\begin{pmatrix} 0_\EE & \beta \\ 1_\EE & 0_\EE \end{pmatrix} \in (\EE/\KK,\beta)\ .\end{align*}
\end{itemize}
\end{de}

\begin{lemma}
Let $\EE/\KK$ be a separable quadratic extension and let $\beta \in \KK^*\sm N(\EE)$. Then the following holds:
\begin{enumerate}[label=(\alph*)]
\item The set $$\HH:=(\EE/\KK,\beta)\subseteq M_2(\EE)$$ is an associative division ring.
\item We have
$$\HH=1_\EE\cdot \EE+ e\cdot \EE\ .$$
\item The map
$$\sigma_s:\HH\to \HH,\ x+e\cdot y\mapsto \bar{x}-e\cdot y$$
is an involution of $\HH$ extending $\sigma$.
\item The maps
\begin{align*}
N:\HH\to \KK,\ x+e\cdot y\mapsto N(x)-\beta\cdot N(y)\ , && T:\HH\to \KK,\ x+e\cdot y\mapsto T(x)
\end{align*}
are extensions of $N$ and $T$.
\end{enumerate}
\end{lemma}

\begin{bew}
Cf. (9.2), (9.3) and (9.4) of \cite{TW}.\qed
\end{bew}

\begin{de} \newglossaryentry{SI}{type=symbols,name={\ensuremath{\sigma_s}},description=standard involution, sort=alternative ring}
A \textit{quaternion division algebra}\index{quaternion division algebra} is an algebra $\HH$ isomorphic to $(\EE/\KK,\beta)$ for some separable quadratic extension $\EE/\KK$ and some $\beta\in \KK^*\sm N(\EE)$. The map \gls{SI} is the \textit{standard involution of $\mathit{\HH}$}\index{standard involution}.
\end{de}

\begin{de}
\newglossaryentry{Oct}{type=symbols,name={\ensuremath{(\HH,\beta)}},description=octonion algebra with respect to the quaternion division algebra $\HH$ and $\beta$, sort=alternative ring}
Let $\HH$ be a quaternion division algebra. 
\begin{itemize}
\item Given $x\in \HH$, we write $\bar{x}:=\sigma_s(x)$.
\item Given $\beta\in \KK^*$, we set
\begin{align*}\gls{Oct}:=\left\{ \begin{pmatrix} x & \beta \bar{y} \\y & \bar{x} \end{pmatrix} \mid x,y\in \HH\right\}\subseteq M_2(\HH)\ , && e:=\begin{pmatrix} 0_\HH & \beta \\ 1_\HH & 0_\HH \end{pmatrix} \in (\HH,\beta)\ .\end{align*}
\item We define a multiplication on $(\HH,\beta)$ by
\begin{align*}
\begin{pmatrix} x & \beta \bar{y} \\ y & \bar{x}\end{pmatrix}\cdot \begin{pmatrix} u & \beta \bar{v} \\ v & \bar{u}\end{pmatrix}:=\begin{pmatrix} xu+\beta v\bar{y} & \beta(\bar{v} x+\bar{y}\bar{u}) \\ \bar{x} v+uy & \bar{u}\bar{x}+\beta y\bar{v}\end{pmatrix}\ ,
\end{align*}
which is non-associative.
\end{itemize}
\end{de}

\begin{lemma}
Let $\HH$ be a quaternion division algebra and let $\beta \in \KK^*\sm N(\HH)$. Then the following holds:
\begin{enumerate}[label=(\alph*)]
\item With the ordinary matrix addition and the above multiplication, the set $$\OO:=(\HH,\beta)\subseteq M_2(\HH)$$ is an alternative division ring.
\item We have
$$\OO= 1_\HH\cdot \HH+ e\cdot \HH\ .$$
\item The map
$$\sigma_s:\OO\to \OO,\ x+e\cdot y\mapsto \bar{x}-e\cdot y$$
is an involution of $\OO$ extending $\sigma_s$.
\item The maps
\begin{align*}
N:\OO\to \KK,\ x+e\cdot y\mapsto N(x)-\beta\cdot N(y)\ , && T:\OO\to \KK,\ x+e\cdot y\mapsto T(x)
\end{align*}
are extensions of $N$ and $T$.
\end{enumerate}
\end{lemma}

\begin{bew}
Cf. (9.8) of \cite{TW}.\qed
\end{bew}

\begin{de} 
An \textit{octonion division algebra}\index{octonion division algebra} is an algebra $\OO$ isomorphic to $(\HH,\beta)$ for some quaternion division algebra $\HH$ and some $\beta\in \KK^*\sm N(\HH)$. The map \gls{SI} is the \textit{standard involution of $\mathit{\OO}$}\index{standard involution}.
\end{de}

\begin{bem}
In the following, we list the basic properties of an octonion division algebra $\OO$ which will be needed in the sequel. Since each subalgebra of $\OO$ is a division algebra by (20.8) of \cite{TW}, we will omit the term ``division'' whenever we deal with an octonion division algebra and its subalgebras.
\end{bem}

\begin{no} Throughout the rest of this paragraph, $\OO$ denotes an octonion division algebra.
\end{no}

\begin{bem}
The norm $N:\OO\to \KK$ is a quadratic form with associated bilinear form
\begin{align*}
\langle \cdot,\cdot\rangle:\OO\times \OO\to \KK,\ (x,y)\mapsto x\bar{y}+y\bar{x}=T(x\bar{y})\ .
\end{align*}
\end{bem}

\begin{de}
Let $V$ be a vector space over $\KK$. A quadratic form $q:V\to \KK$ is \textit{non-defective}\index{non-defective}\index{quadratic form (anisotropic)!(non-)defective} if the associated bilinear form
$$f_q:V\to \KK,\ (x,y)\mapsto q(x+y)-q(x)-q(y)$$
is non-degenerate, i.e., we have $V^\bot=\{0_V\}$, cf. definition \ref{460}. Otherwise, it is \textit{defective}\index{defective}.
\end{de}

\begin{lemma}\label{114} There exists an element $x\in \OO$ such that $\bar{x}\neq x$.
\end{lemma}

\begin{bew}
Cf. (20.15) of \cite{TW}.\qed
\end{bew}

\begin{kor}\label{126}
We have
$$\langle \cdot,\cdot \rangle\not\equiv 0_\KK\ .$$
As a consequence, $\langle\cdot,\cdot\rangle$ is non-degenerate, hence $N$ is non-defective.
\end{kor}

\begin{bewzwei}\ 
\begin{itemize}
\item $\Char\ \OO\neq 2$: In this case, we have
$$\langle 1_\OO,1_\OO\rangle=2\cdot N(1_\OO)=2\cdot 1_\OO\neq 0_\OO\ .$$
\item $\Char \OO=2$: By lemma \ref{114}, there is an element $x\in \OO$ such that $x\neq\bar{x}$. We obtain
$$\langle x,1_\OO\rangle =x+\bar{x}\neq 0_\OO\ .$$
\end{itemize}
Now the map $\langle\cdot,\cdot\rangle=\bar{T}$ is non-degenerate by (20.16) of \cite{TW}, cf. definition (20.12) of \cite{TW}.\qed
\end{bewzwei}

\begin{lemma}
Let $x,y\in \OO$. Then the following holds:
\begin{enumerate}[label=(\alph*)]
\item \label{111} There is an associative subalgebra $\AA$ containing both $x$ and $y$.
\item \label{115} If we have $\bar{x}\neq x$, then $\AA$ can be chosen to be a quaternion subalgebra $\HH$.
\item \label{93} There exists a quaternion subalgebra $\HH$ containing $x$.
\end{enumerate}
\end{lemma}

\begin{bew}
Parts (a) and (b) result from the proof of (20.22) in \cite{TW}. Part (c) is (20.23) of \cite{TW}.\qed
\end{bew}

\newglossaryentry{Doubling}{type=results,name={{Doubling Process}},description={},sort=res}

\begin{lemma}[\textbf{\gls{Doubling}}]\label{116}
Let $\AA$ be a subalgebra such that $\AA^\bot \nsubseteq \AA$, let $e\in \AA^\bot\sm \AA$ and let $u:=-N(e)$. Then the following holds:
\begin{enumerate}[label=(\alph*)]
\item The set $\tilde{\AA}:=\AA+e\cdot \AA$ is a subalgebra.
\item We have
\begin{align*}
\bar{e}=-e\ , && e\cdot \AA\subseteq \AA^\bot\ , && \AA\cap e\cdot \AA=\{0_\OO\}\ .
\end{align*}
\item \label{250} Given $x,y\in \AA$, we have
\begin{align*}
(e\cdot x)(e\cdot y)=u(y\bar{x})\ , && (e\cdot x)y=e\cdot yx\ , && x(e\cdot y)=e\cdot \bar{x}y\ .
\end{align*}
\end{enumerate}
\end{lemma}
\begin{bew}
This is (20.17) of \cite{TW}.\qed
\end{bew}

\begin{bem}
Let $\AA$ be an alternative division ring. By definition, we have $$Z(\AA)=\{ x\in \AA \mid [x,\AA]=[x,\AA,\AA]=0_\AA \}\ .$$
\end{bem}

\begin{lemma}\label{123}
Let $\AA$ be an alternative division ring. Then we have
$$Z(\AA)=\{ x\in \AA \mid  [x,\AA]=0_\AA \}\ .$$
\end{lemma}

\begin{bew}
The assertion is clearly true if $\AA$ is a skew-field, thus we may suppose $\OO:=\AA$ to be an octonion division algebra. By proposition (1.9.2) of \cite{S}, we have
$$\{ x\in \OO \mid [x,\OO,\OO]=0_\OO \}\subseteq \{ x\in \OO \mid [x,\OO]=0_\OO\}$$
and therefore
\begin{align*} Z(\OO)&=\{ x\in \OO\mid [x,\OO]=[x,\OO,\OO]=0_\OO\} \\
&\subseteq \{x\in \OO \mid [x,\OO,\OO]=0_\OO\}\subseteq \{x\in \OO\mid [x,\OO]=[x,\OO,\OO]=0_\OO\}=Z(\OO)\ .\end{align*}\qed
\end{bew}

\section{The Bruck-Kleinfeld Theorem}

On the one hand, we will need the minimum equation in §\ref{130}, and on the other hand, we will need the classification of alternative division rings which are quadratic over a subfield of its center in §\ref{147}. The main steps in the proof of the Bruck-Kleinfeld theorem in \cite{TW} involve those algebras and the corresponding classification result, thus we mention them at this point and dedicate a short paragraph to this fundamental theorem.

\begin{de} Let $\AA$ be an alternative division ring, let $\KK:=Z(\AA)$ and let $\FF$ be a subfield of $\KK$. Then $\AA$ is \textit{quadratic over $\mathit{\FF}$}\index{quadratic over a subfield of the center} if there are maps $T=T^\AA_\FF,N=N^\AA_\FF:\AA\to \FF$ such that
\begin{align}
\forall\ a\in \AA:\ a^2-T(a)a+N(a)=0_\AA\ , && \forall\ a\in \FF:\ T(a)=2a,\  N(a)=a^2\ .\label{129}
\end{align}
The maps $T$ and $N$ are the \textit{trace}\index{trace} and the \textit{norm}\index{norm}, respectively.
\end{de}

\begin{bem}
Trace and norm are uniquely determined by the \textit{minimum equation} \eqref{129}.
\end{bem}

\begin{prop}\label{125}
A non-associative alternative division ring $\AA$ is quadratic over its center. In particular, an octonion division algebra $\OO$ is quadratic over its center.
\end{prop}

\begin{bew}
This is theorem (20.2) of \cite{TW}.\qed
\end{bew}

\begin{prop}\label{121}
Let $\AA$ be an alternative division ring which is quadratic over some subfield $\FF$ of its center $\KK:=Z(\AA)$, let $T$ and $N$ be the trace and the norm, respectively, and let 
$$\sigma:\AA\to \AA,\ x\mapsto T(x)-x\ .$$
Then exactly one of the following holds:
\begin{enumerate}[label=(\roman*)]
\item $\AA=\KK,\ \Char \KK=2,\ \KK^2\subseteq \FF\neq \KK$ and $\sigma=\id_\AA$.
\item $\AA=\KK=\FF$ and $\sigma=\id_\AA$.
\item $\AA=\KK$, $\KK/\FF$ is a separable quadratic extension and $\langle \sigma\rangle=\mathrm{Gal}(\KK/\FF)$.
\item $\AA$ is a quaternion division algebra over $\KK$, $\FF=\KK$ and $\sigma=\sigma_s$.
\item $\AA$ is an octonion division algebra over $\KK$, $\FF=\KK$ and $\sigma=\sigma_s$.
\end{enumerate}
In each case, we have $$N(x)=xx^\sigma=x^\sigma x$$ for each $x\in \AA$.
\end{prop}

\begin{bew}
This is theorem (20.3) of \cite{TW}.\qed
\end{bew}

\newglossaryentry{Bruck}{type=results,name={{Bruck-Kleinfeld-Theorem}},description={},sort=res}

\begin{satz}[\textbf{\gls{Bruck}}]\label{490}
A non-associative alternative division ring is an octonion division algebra.
\end{satz}

\begin{bew} This is a consequence of proposition \ref{125} and proposition \ref{121}.\qed
\end{bew}

\section{Discrete Valuations and Composition Algebras}
Given a Bruhat-Tits building, the defining field for the building at infinity is complete with respect to a discrete valuation. In particular, we will have to deal with octonions and thus composition algebras which are complete with respect to a discrete valuation.

\begin{de}
Let $\AA$ be an alternative division ring. A \textit{discrete valuation of $\mathit{\AA}$}\index{discrete valuation} is a map $\nu:\AA^*\to \ZZ$ such that
\begin{align*}
\nu(xy)=\nu(x)+\nu(y)\ , &&\nu(x+y)\geq \min\{ \nu(x),\nu(y)\}
\end{align*}
for all $x,y\in\AA^*$. A \textit{$\mathit{\nu}$-uniformizer}\index{$\nu$-uniformizer}\index{uniformizer} is an element $\pi\in \AA^*$ such that $\langle \nu(\pi)\rangle=\nu(\AA^*)$.
\end{de}

\begin{lemma}
Let $\AA$ be an alternative division ring with discrete valuation $\nu$. Then the map
$$\delta_\nu:\AA\times \AA\to \RR,\ (x,y)\mapsto\begin{cases} 2^{-\nu(x-y)} & ,\ x\neq y \\
0 & ,\ x=y 
\end{cases}$$
is a metric.
\end{lemma}

\begin{bew}
This is lemma (9.18) of \cite{W}.\qed
\end{bew}

\begin{de}
Let $\AA$ be an alternative division ring with discrete valuation $\nu$. Then $\AA$ is \textit{complete with respect to $\mathit{\nu}$}\index{complete w.r.t. a discrete valuation} if $(\AA,\delta_\nu)$ is a complete metric space.
\end{de}

\begin{de}
A \textit{composition algebra over a field $\mathit{\KK}$}\index{composition algebra} is a unital algebra $\AA$ over $\KK$ together with a non-defective quadratic form $N:\AA\to \KK$ which \textit{permits composition}, i.e., we have
$$\forall\ x,y\in \AA:\qquad N(xy)=N(x)N(y)\ .$$
\end{de}

\begin{lemma} An octonion division algebra $\OO$ is a composition algebra over $\KK:=Z(\OO)$.
\end{lemma}

\begin{bew}
The norm $N$ is non-defective by corollary \ref{126} and  multiplicative by (9.9)(iii) of \cite{TW}.\qed
\end{bew}

\begin{lemma}
Let $\AA$ be an alternative division ring with discrete valuation $\nu$. Then the following holds:
\begin{enumerate}[label=(\alph*)]
\item\label{149} The algebra $\AA$ is complete with respect to $\nu$ if and only if the center $Z(\AA)$ is complete with respect to ${\nu}_{|Z(\AA)}$.
\item \label{53} If $\AA$ is a composition algebra which is complete with respect to $\nu$, we have
\begin{align*}
\forall\ x\in\AA: \qquad \nu(x)=\frac{\nu\big(N(x)\big)}{2}\left(=\frac{\nu\big(N(-x)\big)}{2}=\nu(-x)\right)\ .
\end{align*}
\end{enumerate}
\end{lemma}

\begin{bewzwei}\ 
\begin{enumerate}[label=(\alph*)]
\item This is proposition (23.14) of \cite{W}.
\item This results from proposition 1 of \cite{P}.\qed
\end{enumerate}
\end{bewzwei}

\begin{de} 
Let $\AA$ be a composition division algebra over $\KK$ which is complete with respect to a discrete valuation $\nu$ and such that its residue field $\bar{\AA}$ is a composition algebra over the residue field $\bar{\KK}$. Then $\AA$ is \textit{unramified}\index{unramified}\index{composition algebra!(un)ramified} if we have $\nu(\AA)=\nu(\KK)$, and \textit{ramified}\index{ramified} otherwise.
\end{de}

\chapter{Quadratic Spaces}

Quadrangles of quadratic form type are parametrized by quadratic spaces.

\section{Basic Definitions and Basic Properties}

\begin{de}\ 
\begin{itemize}
\item An \textit{(anisotropic) quadratic space}\index{quadratic space (anisotropic)} is a triple $(L_0,\KK,q)$ such that $\KK$ is a field, $L_0$ is a right vector space over $\KK$ and $q$ is an \textit{(anisotropic) quadratic form on $\mathit{L_0}$}\index{quadratic form (anisotropic)}, i.e., $q:L_0\to \KK$ is a map such that the following holds:
\begin{enumerate}[leftmargin=25pt,label=(Q\arabic*)]
\item $\forall\ a\in L_0,\ t\in \KK:\ q(at)=q(a)t^2$.

\newglossaryentry{Bil}{type=symbols,name={\ensuremath{f_q}},description=bilinear form with respect to the quadratic form $q$, sort=quadratic space}
\item The map
$$\gls{Bil}:L_0\times L_0\to \KK,\ (a,b)\mapsto q(a+b)-q(a)-q(b)$$
is bilinear.
\item $\forall\ a\in L_0:\ q(a)=0_\KK\ \Leftrightarrow\ a=0_{L_0}$.
\end{enumerate}
\item A quadratic space $(L_0,\KK,q)$ is \textit{proper}\index{proper!quadratic space}\index{quadratic space (anisotropic)!proper} if we have $f_q\not\equiv 0_\KK$.
\item Two quadratic spaces $(L_0,\KK,q)$ and $(\tilde{L}_0,\tilde{\KK},\tilde{q})$ are \textit{isomorphic}\index{isomorphism!of quadratic spaces} if there is an isomorphism
\begin{align*}
\Phi=(\p,\phi):(L_0,\KK)\to(\tilde{L}_0,\tilde{\KK})
\end{align*}
of vector spaces such that $$ \tilde{q}\circ \p=\phi\circ q\ .$$
\item A quadratic space $(L_0,\KK,q)$ is \textit{unital}\index{unital!quadratic space}\index{quadratic space (anisotropic)!unital} if there is a \textit{basepoint}\index{basepoint} $\epsilon\in L_0$ such that $q(\epsilon)=1_\KK$.
\item Given a quadratic space $(L_0,\KK,q)$ with basepoint $\epsilon$, we set
\begin{align*}
T:L_0\to \KK,\ x\mapsto f_q(\epsilon,x)=f_q(x,\epsilon)\ , && \sigma:L_0\to L_0,\ x\mapsto \bar{x}:=\epsilon\cdot T(x)-x\ .
\end{align*}
Given $a\in L_0^*$, we set
\begin{align*}
\pi_a:L_0\to L_0,\ v\mapsto v-a\cdot f_q(a,v)/q(a)\ , && h_a:L_0\to L_0,\ v\mapsto \pi_a\pi_\epsilon(v)\cdot q(a)\ .
\end{align*}
\end{itemize}
\end{de}

\begin{lemma}\label{476}
Let $(L_0,\KK,q)$ and $(\tilde{L}_0,\tilde{\KK},\tilde{q})$ be quadratic spaces with basepoints $\epsilon$ and $\tilde{\epsilon}$, respectively, and let $(\p,\phi):(L_0,\KK,q)\to (\tilde{L}_0,\tilde{\KK},\tilde{q})$ be an isomorphism of quadratic spaces such that $\p(\epsilon)=\tilde{\epsilon}$. Then we have
\begin{align*}
\tilde{T}\circ \p=\phi\circ T\ , && \tilde{\sigma}\circ \p=\p\circ \sigma\ .
\end{align*}
\end{lemma}

\begin{bew}
Given $x,y\in L_0$, we have
\begin{align*}
f_{\tilde{q}}\big(\p(x),\p(y)\big)&=\tilde{q}\big( \p(x)+\p(y)\big)-\tilde{q}\big(\p(x)\big)-\tilde{q}\big(\p(y)\big) \\
&=\phi\big( q(x+y)-q(x)-q(y)\big)=\phi\big( f_q(x,y)\big)\ .
\end{align*}
In particular, we have $$\tilde{T}\big(\p(x)\big)=f_{\tilde{q}}\big(\p(x),\tilde{\epsilon}\big)=\phi\big(f_q(x,\epsilon)\big)=\phi\big(T(x)\big)$$
and thus
$$\p(x)^{\tilde{\sigma}}=\tilde{\epsilon}\cdot \tilde{T}\big(\p(x)\big)-\p(x)=\tilde{\epsilon}\cdot \phi\big(T(x)\big)-\p(x)=\p\big(\epsilon\cdot T(x)-x\big)=\p(x^\sigma)$$
for each $x\in L_0$.\qed
\end{bew}

\begin{no}
Throughout this paragraph, $(L_0,\KK,q)$ is a quadratic space with basepoint $\epsilon$.
\end{no}

\begin{lemma}\label{330}
The maps $T:L_0\to \KK$ and $\sigma:L_0\to L_0$ are $\KK$-linear.
\end{lemma}

\begin{bew}
Given $x\in L_0$ and $s\in \KK$, we have
$$T(x\cdot s)=f_q(\epsilon, x\cdot s)=f_q(\epsilon, x)s=T(x)s$$
and thus
$$(x\cdot s)^\sigma=\epsilon\cdot T(x\cdot s)-x\cdot s=\epsilon\cdot T(x)s-x\cdot s=\big(\epsilon\cdot T(x)-x\big)\cdot s=x^\sigma\cdot s\ .$$
Given $x,y\in L_0$, we have
$$T(x+y)=f_q(\epsilon, x+y)=f_q(\epsilon,x)+f_q(\epsilon,y)=T(x)+T(y)$$
and thus
$$(x+y)^\sigma=\epsilon\cdot T(x+y)-(x+y)=\epsilon\cdot T(x)-x+\epsilon\cdot T(y)-y=x^\sigma+y^\sigma\ .$$
\qed
\end{bew}

\begin{lemma}\label{329}
Given $x\in L_0$, $a\in L_0^*$, we have
$$h_a(x)=a\cdot f_q(a,x^\sigma)-x^\sigma\cdot q(a)\ .$$
\end{lemma}

\begin{bew}
Given $x\in L_0$, $a\in L_0^*$, we have
\begin{align*}
h_a(x)&=\pi_a\pi_\epsilon(x)\cdot q(a)=\pi_a\big(x-\epsilon\cdot f_q(\epsilon,x)\big)\cdot q(a) \\
&=-\pi_a(x^\sigma)\cdot q(a)=-\big(x^\sigma-a\cdot f_q(a,x^\sigma)/q(a)\big)q(a)=a\cdot f_q(a,x^\sigma)-x^\sigma\cdot q(a)\ .
\end{align*}\qed
\end{bew}

\begin{kor}\label{327}
Given $a\in L_0^*$, the corresponding \textit{Hua map}\index{Hua map} $h_a$ is $\KK$-linear.
\end{kor}

\begin{bew}
Let $a\in L_0^*$, $x,y\in L_0$ and $s\in \KK$. By lemma \ref{329} and lemma \ref{330}, we have
\begin{align*}
h_a(x\cdot s)=a\cdot f_q(a,x^\sigma\cdot s )-x^\sigma\cdot sq(a)=\big(a\cdot f_q(a,x^\sigma)-x^\sigma\cdot q(a)\big)\cdot s=h_a(x)\cdot s
\end{align*}
and
\begin{align*}
h_a(x+y)&=a\cdot f_q\big(a,(x+y)^\sigma\big)-(x+y)^\sigma\cdot q(a)\\
&=a\cdot \big(f_q(a,x^\sigma)+f_q(a,y^\sigma)\big)-x^\sigma\cdot q(a)-y^\sigma\cdot q(a)=h_a(x)+h_a(y)\ .
\end{align*}\qed
\end{bew}

\begin{lemma}\label{328} Given $x\in L_0$, $a\in L_0^*$ and $s\in \KK$, we have
\begin{align*} h_{a\cdot s}(x)=h_a(x\cdot s^2)=h_a(x)\cdot s^2\ , && h_{a\cdot s}=s^2\cdot h_a\ .\end{align*}
\end{lemma}

\begin{bew}
Let $x\in L_0$, $a\in L_0^*$ and $s\in \KK$. By lemma \ref{329}, lemma \ref{330} and corollary \ref{327}, we have
\begin{align*}
h_{a\cdot s}(x)&=a\cdot sf_q(a\cdot s,x^\sigma)-x^\sigma\cdot q(a\cdot s) =a\cdot f_q(a,x^\sigma)s^2-x^\sigma\cdot q(a)s^2\\
&=a\cdot f_q\big(a,(xs^2)^\sigma\big)-(xs^2)^\sigma\cdot q(a)=h_a(x\cdot s^2)=h_a(x)\cdot s^2\ .
\end{align*}
\qed
\end{bew}

\begin{lemma}\label{326}
Given $x,y\in L_0$, we have $$f_q(x^\sigma,y)=f_q(y^\sigma,x)=f_q(x,y^\sigma)\ .$$
\end{lemma}

\begin{bew}
Given $x,y\in L_0$, we have
\begin{align*}
f_q(x^\sigma,y)=f_q\big(\epsilon\cdot T(x)-x,y\big)=-f_q(x,y)+f_q(\epsilon,y)T(x)=-f_q(x,y)+T(x)T(y)\ ,
\end{align*}
which is symmetric in $x$ and $y$.\qed
\end{bew}

\begin{bem}\ 
\begin{enumerate}[label=(\alph*)]
\item Given $x\in L_0$, we have
$$f_q(x,x)=q(2x)-q(x)-q(x)=4q(x)-2q(x)=2q(x)\ .$$
\item \label{325} We have $$\epsilon^\sigma=\epsilon\cdot f_q(\epsilon,\epsilon)-\epsilon=\epsilon\cdot 2-\epsilon=\epsilon\ .$$
\end{enumerate}
\end{bem}

\begin{kor}
Given $x\in L_0$, we have $$T(x^\sigma)=T(x)\ .$$
\end{kor}

\begin{bew}
Let $x\in L_0$. By remark \ref{325} and lemma \ref{326}, we have
$$T(x)=f_q(\epsilon,x^\sigma)=f_q({\epsilon}^\sigma,x^\sigma)=f_q\big(\epsilon,(x^{\sigma})^\sigma\big)=T(x^\sigma)\ .$$
\qed
\end{bew}

\begin{kor}
We have $\sigma^2=\id_{L_0}$.
\end{kor}

\begin{bew}
Given $x\in L_0$, we have
$$(x^\sigma)^\sigma=\epsilon\cdot T(x^\sigma)-x^\sigma=\epsilon\cdot T(x)-\big(\epsilon\cdot T(x)-x\big)=x\ .$$\qed
\end{bew}

\section{Small Dimensions}

Quadratic spaces of  small dimension are in fact quadratic spaces corresponding to fields.

\begin{satz}\label{408} Let $(L_0,\KK,q)$ be a quadratic space with basepoint $\epsilon$ such that $\dim_\KK L_0\leq 2$. Then there is a unique multiplication $\ast:L_0\times L_0\to L_0$ such that the following holds:
\begin{enumerate}[label=(\roman*)]
\item The triple $\tilde{\FF}:=(L_0,+,\ast)$ is a field.
\item The subspace $\tilde{\KK}:=\langle \epsilon\rangle_\KK$ is a subfield of $\tilde{\FF}$.
\item The map $\phi:\KK\to \tilde{\KK},\ s\mapsto \epsilon\cdot s $ is an isomorphism of fields.
\item Given $x\in L_0$, we have $$\phi\big(q(x)\big)=x\ast x^{\sigma}\ .$$
\end{enumerate}
\end{satz}

\begin{bew}
Condition (iii) forces
$$\forall\ s,t\in \KK:\qquad (\epsilon\cdot s)\ast (\epsilon\cdot t)=\phi(s)\ast \phi(t)=\phi(st)=\epsilon \cdot st\ ,$$
which makes $\tilde{\KK}$ really into a field, hence (ii) and (iii) hold.
\begin{itemize}
\item If we have $\dim_{\KK} L_0=1$, then $\tilde{\FF}=\tilde{\KK}$ is a field. By lemma \ref{330} and remark \ref{325}, we have
$$\forall\ s\in \KK:\qquad (\epsilon\cdot s)^\sigma=\epsilon^\sigma\cdot s =\epsilon \cdot s$$
and thus
$$\forall\ s\in \KK:\qquad \phi\big(q(\epsilon\cdot s)\big)=\phi(s^2)=\epsilon\cdot s^2=(\epsilon\cdot s)\ast (\epsilon\cdot s)=(\epsilon\cdot s)\ast (\epsilon\cdot s)^\sigma\ .$$
\item If we have $\dim_\KK L_0=2$, there is an element $\tilde{x}\in L_0\sm \tilde{\KK}$. Conditions (i) and (ii) force $\tilde{x}\ast \epsilon:=\tilde{x}=:\epsilon\ast \tilde{x}$, and condition (iv) forces
\begin{align*}
\tilde{x}\cdot T(\tilde{x})-\tilde{x}\ast \tilde{x}=\tilde{x}\ast(\epsilon\cdot T(\tilde{x})-\tilde{x})=\tilde{x}\ast \tilde{x}^\sigma=\phi\big(q(\tilde{x})\big)=\epsilon\cdot q(\tilde{x})
\end{align*}
and thus $\tilde{x}\ast \tilde{x}:=\tilde{x}\cdot T(\tilde{x})-\epsilon\cdot q(\tilde{x})$. Given $s,t,\tilde{s},\tilde{t}\in \KK$, we set
$$(\epsilon\cdot s+\tilde{x}\cdot t)\ast (\epsilon\cdot \tilde{s}+\tilde{x}\cdot \tilde{t}):=\epsilon\cdot (s\tilde{s}-q(\tilde{x})t\tilde{t})+\tilde{x}\cdot \big(s\tilde{t}+\tilde{s}t+T(\tilde{x})t\tilde{t}\big)\ .$$
Let $f:=y^2-yT(\tilde{x})+q(\tilde{x})\in \KK[y]$. Given $s\in \KK$, we have
$$s^2-sT(\tilde{x})+q(\tilde{x})=q\big(\epsilon\cdot (s)\big)-sf_q(\epsilon,\tilde{x})+q(\tilde{x})=q\big(\epsilon\cdot (s)\big)+q(\tilde{x})-f_q(\epsilon\cdot s,\tilde{x})=q(\epsilon\cdot s+\tilde{x})\neq 0_\KK\ .$$
Therefore, the polynomial $f\in \KK[y]$ is irreducible. Let $x$ be an element in the algebraic closure of $\KK$ such that $f({x})=0_{\KK}$. Then $\FF:=\KK(x)$ is a field with $\dim_\KK \FF=\deg(f)=2$ and multiplication given by
$$x^2=x T(\tilde{x})-q(\tilde{x})\ .$$
Therefore, the map 
$$\Phi:\FF\to \tilde{\FF},\ s+xt\mapsto \epsilon\cdot s+\tilde{x}\cdot t$$
is an isomorphism of fields, hence (i) holds. Given $s,t\in \KK$, we finally have
\begin{align*}
\phi\big(q(\epsilon\cdot s+\tilde{x}\cdot t)\big)&=\epsilon\cdot \big(s^2+q(\tilde{x})t^2+f_q(\epsilon\cdot s,\tilde{x}\cdot t)\big)=\epsilon\cdot\big(s^2+f_q(\epsilon,\tilde{x})st+q(\tilde{x})t^2\big) \\
&=\epsilon\cdot(s^2+T(\tilde{x})st+q(\tilde{x})t^2)+\tilde{x}\cdot(-st+st+T(\tilde{x})t^2-T(\tilde{x})t^2\big) \\
&=\big(\epsilon\cdot s+\tilde{x}\cdot t\big)\ast\big(\epsilon\cdot s+\epsilon\cdot T(\tilde{x}\cdot t)-\tilde{x}\cdot t\big) \\
&=(\epsilon\cdot s+\tilde{x}\cdot t)\ast (\epsilon\cdot s+\tilde{x}\cdot t)^\sigma\ ,
\end{align*}
hence (iv) holds.
\end{itemize}
\qed
\end{bew}

\begin{de} \newglossaryentry{field}{type=symbols,name={\ensuremath{\FF(L_0,\KK,q)}},description=field associated with the quadratic space ${(L_0,\KK,q)}$ with $\dim_\KK L_0\leq 2$, sort=quadratic space}
Given a quadratic space $(L_0,\KK,q)$ with basepoint $\epsilon$ and $\dim_\KK L_0\leq 2$, we set
$$\gls{field}:=\tilde{\FF}\ ,$$
where $\tilde{\FF}$ is as in theorem \ref{408}.
\end{de}

\begin{kor}\label{432}
Let $(L_0,\KK,q)$ be a quadratic space with basepoint $\epsilon$ and $\dim_\KK L_0\leq 2$, let $\tilde{\FF}:=\FF(L_0,\KK,q)$, let $\tilde{\KK}:=\langle \epsilon\rangle_\KK$ and let $\phi:\KK\to \tilde{\KK},\ s\mapsto \epsilon\cdot s$. Then $\tilde{\FF}$ is quadratic over $\tilde{\KK}$, we have $N^{\tilde{\FF}}_{\tilde{\KK}}=\phi\circ q$, the map $(\id_{L_0},\phi):(L_0,\KK,q)\to (\tilde{\FF},\tilde{\KK},N^{\tilde{\FF}}_{\tilde{\KK}})$ is an isomorphism of quadratic spaces, and exactly one of the following holds:
\begin{enumerate}[label=(\roman*)]
\item $\Char \tilde{\FF}=2$, $\tilde{\FF}^2\subseteq \tilde{\KK}\neq \tilde{\FF}$ and $\sigma=\id_{\tilde{\FF}}$, which means that $\tilde{\FF}/\tilde{\KK}$ is inseparable.
\item $\tilde{\KK}=\tilde{\FF}$ and $\sigma=\id_{\tilde{\FF}}$.
\item $\tilde{\FF}/\tilde{\KK}$ is a separable quadratic extension and $\langle \sigma\rangle=\mathrm{Gal}(\tilde{\FF}/\tilde{\KK})$.
\end{enumerate}
\end{kor}

\newpage

\begin{bew}
By construction, we have
$$\forall\ x\in \tilde{\FF}:\qquad x\ast x-x\ast\big(\epsilon\cdot T(x)\big)+\epsilon\cdot q(x)=0_{\tilde{\FF}}\ ,$$
thus $\tilde{\FF}$ is quadratic over $\tilde{\KK}=\langle \epsilon\rangle_\KK$ with $N^{\tilde{\FF}}_{\tilde{\KK}}=\phi\circ q$, $T^{\tilde{\FF}}_{\tilde{\KK}}=\phi\circ T$ and
$$\forall\ x\in L_0,\ s\in \KK:\qquad \id_{L_0}(x\cdot s)=x\cdot s=x\ast(\epsilon\cdot s)=\id_{L_0}(x)\ast \phi(s)$$ (which shows that $(\id_{L_0},\phi)$ is an isomorphism of quadratic spaces), and we have
$$\forall\ x\in \tilde{\FF}:\qquad \sigma(x)=\epsilon\cdot T(x)-x=T^{\tilde{\FF}}_{\tilde{\KK}}(x)-x\ ,$$
which is just the map $\sigma$ in proposition \ref{121} (which we apply). We have $\dim_{\tilde{\KK}}\tilde{\FF}\leq 2$, thus $(\tilde{\FF},\tilde{\KK},\sigma)$ is neither of type (iv) nor of type (v).\qed
\end{bew}

\begin{lemma} Let $\AA$ be an alternative division ring which is quadratic over a subfield $\FF$ of its center. Then $(\AA,\FF,N)$ with $N:=N^\AA_\FF$ is a quadratic space.
\end{lemma}

\begin{bew}
By proposition \ref{121}, we have $N(x)=xx^\sigma$ for each $x\in \AA$.
\begin{enumerate}[label=(Q\arabic*),leftmargin=27pt]
\item Given $s\in \FF$, we have
$$N(x\cdot s)=N(x)N(s)=N(x)s^2\ .$$
\item Given $x,y\in \AA$, we have
\begin{align*}
f_q(x,y)&=N(x+y)-N(x)-N(y)\\
&=(x+y)(x+y)^\sigma-xx^\sigma-yy^\sigma=xy^\sigma+yx^\sigma=xy^\sigma+(xy^\sigma)^\sigma=T(xy^\sigma)\ ,
\end{align*}
which is $\KK$-linear in $x$ and $y$ by lemma \ref{330}.
\item Given $x\in \AA$, we have
$$N(x)\neq 0_\AA\ \Leftrightarrow\ x\in U_\AA=\AA^*\ .$$
\end{enumerate}\qed
\end{bew}

\begin{de}
A quadratic space $(L_0,\KK,q)$ with basepoint $\epsilon$ is \textit{(linear) of type (m)} if we have $(L_0,\KK,q)\cong(\AA,\FF,N^\AA_\FF)$ as in (m) of proposition \ref{121}, i.e., the quadratic space $(L_0,\KK,q)$ is of type
\begin{enumerate}[label=(\roman*)]
\item if $\FF:=L_0$ is a field with $\Char \FF=2,\ \FF^2\subseteq \KK\neq \FF$, $\sigma=\id_\FF$ and $q=N^\FF_\KK$,
\item if we have $L_0=\KK$, $\sigma=\id_\KK$ and $q=N^\KK_\KK$,
\item if $\EE:=L_0$ is a field, $\EE/\KK$ is a separable quadratic extension, $\langle \sigma\rangle=\mathrm{Gal}(\EE/\KK)$ and $q=N^\EE_\KK$,
\item if $\HH:=L_0$ is a quaternion division algebra over $\KK$, $\sigma=\sigma_s$ and $q=N^\HH_\KK$,
\item if $\OO:=L_0$ is an octonion division algebra over $\KK$, $\sigma=\sigma_s$ and $q=N^\OO_\KK$.
\end{enumerate}\index{quadratic space (anisotropic)!(linear) of type (m)}
\end{de}

\begin{bem}
These are exactly the quadratic spaces such that the corresponding Moufang set of quadratic form type is isomorphic to a Moufang set of linear type, cf. theorem \ref{379} and lemma \ref{456}. This explains the terminology.
\end{bem}

\section{Clifford Algebras and the Clifford Invariant}

\begin{de}\ 
\begin{itemize}
\item A \textit{central simple algebra}\index{central simple}\index{algebra!central simple} is a finite-dimensional associative algebra $(A,\KK)$ which is simple as a ring and such that $Z(A)\cong\KK$ as fields.
\item Two central simple algebras $(A,\KK)$ and $(\tilde{A},\tilde{\KK})$ are \textit{isomorphic}\index{isomorphism!of algebras} if there is an isomorphism  $(\gamma,\phi):(A,\KK)\to (\tilde{A},\tilde{\KK})$ of vector spaces such that $\gamma$ is an isomorphism of rings.

\newglossaryentry{Clifa}{type=symbols,name={\ensuremath{C(q)}},description=Clifford algebra with respect to the quadratic form $q$, sort=quadratic space}
\item Let $(L_0,\KK,q)$ be a quadratic space, let $T(L_0)$ be the tensor algebra of $L_0$ and let
$$I(q):=\langle u\otimes u-1_\KK\cdot q(u) \mid u\in L_0 \rangle\ .$$
Then $\gls{Clifa}:=T(L_0)/I(q)$ is the \textit{Clifford algebra of $\mathit{q}$}\index{Clifford algebra}.
\end{itemize}
\end{de}

\begin{lemma}
Let $(L_0,\KK,q)$ be a non-defective quadratic space such that $\dim_\KK L_0$ is even. Then $\big(C(q),\KK\big)$ is a central simple algebra. 
\end{lemma}

\begin{bew}
This results from proposition (11.6) of \cite{EKM}.\qed
\end{bew}

\newglossaryentry{Wedderburn}{type=results,name={{Wedderburn}},description={},sort=res}

\begin{satz}[\textbf{\gls{Wedderburn}}]\label{401}
Given a central simple algebra $(A,\KK)$, there are a unique natural number $n\in \NN^*$ and an associative division algebra $(\DD,\KK)$ such that $(A,\KK)\cong \big(M_n(\DD),\KK\big)$ as algebras. The algebra $(\DD,\KK)$ is unique up to isomorphism of algebras.
\end{satz}

\begin{bew}
Cf. theorem (1.1) of \cite{KMRT}.\qed
\end{bew}

\begin{de}\
\begin{itemize}
\item Given a central simple algebra $(A,\KK)$, we set $\mathrm{S}(A,\KK):=[(\DD,\KK)]$, where $(\DD,\KK)$ is a division algebra as in theorem \ref{401} and $[(\DD,\KK)]$ denotes its isomorphism class.

\newglossaryentry{Clifi}{type=symbols,name={\ensuremath{\mathrm{Clif}(q)}},description=Clifford invariant of the quadratic form $q$, sort=quadratic space}
\item Let $(L_0,\KK,q)$ be a non-defective quadratic space such that $\dim_\KK L_0$ is even. Then $$\gls{Clifi}:=\mathrm{S}\big(C(q),\KK\big)$$ is the \textit{Clifford invariant of $\mathit{q}$}\index{Clifford invariant}.
\end{itemize}
\end{de}

\begin{lemma}\label{374}
Given two isomorphic central simple algebras $(A,\KK)$ and $(\tilde{A},\tilde{\KK})$, we have
$$\mathrm{S}(A,\KK)=S(\tilde{A},\tilde{\KK})\ .$$\end{lemma}

\begin{bew}
Let $n,\tilde{n}$ and $(\DD,\KK),(\tilde{\DD},\tilde{\KK})$ be as in theorem \ref{401}. Then we have
$$\big(M_n(\DD),\KK\big)\cong (A,\KK)\cong (\tilde{A},\tilde{\KK})\cong \big(M_{\tilde{n}}(\tilde{\DD}),\tilde{\KK}\big)$$
as algebras and thus $(\DD,\KK)\cong (\tilde{\DD},\tilde{\KK})$ by theorem \ref{401}.\qed
\end{bew}

\begin{lemma}\label{376}
Let $(L_0,\KK,q)$, $(\tilde{L}_0,\tilde{\KK},\tilde{q})$ be isomorphic quadratic spaces. Then we have $$\big(C(q),\KK\big)\cong \big(C(\tilde{q}),\tilde{\KK}\big)$$
as algebras. In particular, we have $\mathrm{Clif}(q)=\mathrm{Clif}(\tilde{q})$ if the dimensions are even and at least one (and thus both) quadratic spaces are non-defective.
\end{lemma}

\begin{bew}
This results from (12.23) of \cite{TW}. In particular, we may apply lemma \ref{374}.\qed
\end{bew}

\newpage

\section{Norm Splittings}

\begin{de}
Let $(L_0,\KK,q)$ be a quadratic space. A \textit{norm splitting of $\mathit{q}$}\index{norm splitting} is a triple $(\EE,\cdot,\{v_1,\ldots,v_d\})$ such that the following holds:
\begin{enumerate}[label=(N\arabic*),leftmargin=27pt]
\item $\EE/\KK$ is a separable quadratic extension,
\item $\cdot:L_0\times \EE\to L_0$ is a scalar multiplication extending the scalar multiplication $\cdot:L_0\times \KK\to L_0$,
\item $\{v_1,\ldots,v_d\}$ is an $\EE$-Basis of $L_0$ with
$$\forall\ t_1,\ldots,t_d\in \EE:\qquad q\big( \sum_{i=1}^d v_it_i\big)=\sum_{i=1}^d s_iN(t_i)\ ,$$
where $s_i=q(v_i)$ for each $i\in\{1,\ldots,n\}$ and $N=N^\EE_\KK$.
\end{enumerate}
The elements $s_1,\ldots,s_d\in \KK$ are called the \textit{constants of the norm splitting}\index{constants of a norm splitting}.
\end{de}

\begin{lemma}\label{375}
Let $\OO$ be an octonion division algebra with center $\KK:=Z(\OO)$ and let $\EE$ be a subfield such that $\EE/\KK$ is a separable quadratic extension (which exists by lemma \ref{114} and (20.19) of \cite{TW}). Then there are $v_1,\ldots,v_4\in \OO$ such that $(\EE,\cdot,\{v_1,\ldots,v_4\})$ is a norm splitting of $(\OO,\KK,N)$, satisfying $s_1\cdots s_4\in N(\EE)$. 
\end{lemma}

\begin{bew}
Let $v_1:=1_\OO$. By (20.20) of \cite{TW}, there is an element $v_2\in \EE^\bot\sm\EE$, and $\HH_2:=\EE+v_2\cdot \EE$ is a quaternion division algebra. By (20.21) of \cite{TW}, there is an element $v_3\in \HH_2^\bot\sm \HH_2$. Finally, let $v_4:=v_2v_3\in \HH_2^\bot \sm \HH_2$. Then $\{v_1,\ldots,v_4\}$ is an $\EE$-Basis of $\OO$. By construction and (20.20) of \cite{TW} again,
$$\HH_i:=\EE+v_i\cdot \EE$$
is a quaternion division algebra for each $i\in \{2,3,4\}$, satisfying $v_j\in \HH_i^\bot \sm \HH_i$ for all $i\neq j\in \{2,3,4\}$. By lemma \ref{116}, we have
\begin{align*}
N\big(\sum_{i=1}^4 t_iv_i\big)&= (t_1v_1+t_2v_2+t_3v_3+t_4v_4)\cdot (\bar{v}_1\bar{t_1}+\bar{v}_2\bar{t}_2+\bar{v}_3\bar{t}_3+\bar{v}_4\bar{t}_4) \\
&=\sum_{i=1}^4 t_iv_i\bar{v_i}\bar{t}_i+\sum_{i\neq j} \big( (t_iv_i)(\bar{v}_j\bar{t}_j)+(t_jv_j)(\bar{v}_i\bar{t_i})\big)\\
&=\sum_{i=1}^4 N(v_i)N(t_i)+\sum_{i\neq j} \big(\bar{v}_j(\bar{v}_i\bar{t}_i\bar{t_j})+(v_j\bar{t}_j)(\bar{v_i}\bar{t}_i)\big) \\
&=\sum_{i=1}^4 N(v_i)N(t_i)+\sum_{i\neq j} \big(\bar{v}_j(\bar{v}_i\bar{t}_i\bar{t_j})+v_j(\bar{v}_i\bar{t}_i\bar{t}_j)\big) \\
&=\sum_{i=1}^4 N(v_i)N(t_i)+\sum_{i\neq j} \big( (\bar{v}_j+v_j)(\bar{v}_i\bar{t}_i\bar{t}_j)\big) \\
&=\sum_{i=1}^4 N(v_i)N(t_i)+\sum_{i\neq j}\big( (-v_j+v_j)(\bar{v}_i\bar{t}_i\bar{t}_j)\big)=\sum_{i=1}^4 N(v_i)N(t_i)\ .
\end{align*}
In particular, we have $s_i=N(v_i)$ for each $i\in\{1,\ldots,4\}$, and thus, by lemma \ref{116} again, 
$$s_1\cdots s_4=N(v_1)\cdots N(v_4)=N(v_1\cdots v_4)=N(v_2v_3v_2v_3)=N(-v_2^2v_3^2)\in N(\KK)\subseteq N(\EE)\ .$$
\qed
\end{bew}

\begin{de}
A quadratic space $(L_0,\KK,q)$\index{quadratic space (anisotropic)!of type $E_n$} is 
\begin{itemize}
\item \textit{of type $\mathit{E_6}$} if $\dim_\KK L_0=6$ and $q$ has a norm splitting,
\item \textit{of type $\mathit{E_7}$} if $\dim_\KK L_0=8$ and $q$ has a norm splitting $(\EE,\cdot,\{v_1,\ldots,v_4\})$ such that $$s_1\cdots s_4\notin N(\EE)\ ,$$
\item \textit{of type $\mathit{E_8}$} if $\dim_\KK L_0=12$ and $q$ has a norm splitting $(\EE,\cdot,\{v_1,\ldots,v_6\})$ such that $$-s_1\cdots s_6\in N(\EE)\ .$$
\end{itemize}
\end{de}

\begin{bem}\label{382} As we only deal with anisotropic quadratic spaces and since each anisotropic space having a norm splitting is automatically non-defective by remark (12.12) of \cite{TW} (and thus proper), we may reformulate remark (12.30) of \cite{TW} as follows:
\end{bem}

\begin{lemma}
Let $(L_0,\KK,q)$ be a quadratic space having a norm splitting $(\EE,\cdot,\{v_1,\ldots,v_d\})$ with constants $s_1,\ldots,s_d$ and let
$$\gamma:=(-1)^{[d/2]}s_1\cdots s_d\ .$$
Then the following holds:
\begin{enumerate}[label=(\alph*)]
\item\label{377} We have $C(q)\cong M_{2^d}(\KK)$ if $\gamma\in N(\EE)$ and thus $\mathrm{Clif}(q)=[(\KK,\KK)]$.
\item\label{378} We have $C(q)\cong M_{2^{d-1}}(\HH)$ if $\gamma\notin N(\EE)$ and thus $\mathrm{Clif}(q)=[(\HH,\KK)]$, where $\HH=(\EE/\KK, \gamma)$.
\end{enumerate}
\end{lemma}

\begin{bew}
This is remark (12.30) of \cite{TW}.\qed
\end{bew}

\begin{kor} Let $\OO$ be an octonion division algebra with norm $N$ and center $\KK:=Z(\OO)$ and let $(\tilde{L}_0,\tilde{\KK},\tilde{q})$ be a quadratic space of type $E_7$. Then the following holds:
\begin{enumerate}[label=(\alph*)]
\item  We have $\mathrm{Clif}(N)=[(\KK,\KK)]$.
\item  We have $\mathrm{Clif}(\tilde{q})=[(\tilde{\HH},\tilde{\KK})]$ for some quaternion division algebra $\tilde{\HH}$ with center $\tilde{\KK}$.
\item \label{381} We have $(\OO,\KK,N)\not\cong (\tilde{L}_0,\tilde{\KK},\tilde{q})$ as quadratic spaces.
\end{enumerate}
\end{kor}

\begin{bewzwei}\ 
\begin{enumerate}[label=(\alph*)]
\item By lemma \ref{375}, $N$ has a norm splitting $(\EE,\cdot,\{v_1,\ldots,v_4\})$ such that $$\gamma=s_1\cdots s_4\in N(\EE)\ ,$$ hence $\mathrm{Clif}(N)=[(\KK,\KK)]$ by lemma \ref{377}.
\item By definition, $q$ has a norm splitting $(\tilde{\EE},\cdot,\{v_1,\ldots,v_4\})$ such that $$\gamma=s_1\cdots s_4\notin N(\tilde{\EE})\ ,$$ hence we may apply lemma \ref{378}.
\item Since $\KK$ is a field and $\tilde{\HH}$ is non-commutative, we have $(\KK,\KK)\not\cong (\tilde{\HH},\tilde{\KK})$ as algebras and thus $$\mathrm{Clif}(N)\neq \mathrm{Clif}(q)\ .$$
Now the quadratic spaces can't be isomorphic by lemma \ref{376}.
\end{enumerate}\qed
\end{bewzwei}

\section[Quadratic Spaces of Type \texorpdfstring{$F_4$}{F4}]{Quadratic Spaces of Type \texorpdfstring{$\boldsymbol{F_4}$}{F4}}

\begin{de}\label{460}
Let $(L_0,\KK,q)$ be a quadratic space. \newglossaryentry{bot}{type=symbols,name={\ensuremath{W^\bot}},description=orthogonal complement of the subspace $W$, sort=quadratic space}
\begin{itemize}
\item Given a subset $W\subseteq L_0$, we set
$$\gls{bot}:=\{ v\in L_0 \mid f_q(v,W)=0_\KK\}\ .$$

 \newglossaryentry{Def}{type=symbols,name={\ensuremath{\mathrm{Def(q)}}},description=defect of the quadratic form $q$, sort=quadratic space}
\item The subspace $\gls{Def}:=L_0^\bot$ is the \textit{defect of $\mathit{q}$}\index{defect}.
\item The quadratic space $(L_0,\KK,q)$ is \textit{non-defective}\index{non-defective}\index{quadratic space (anisotropic)!(non-)defective} if $\mathrm{Def}(q)=\{0_{L_0}\}$, and \textit{defective}\index{defective} otherwise.
\end{itemize}
\end{de}

\begin{bem}
Let $(L_0,\KK,q)$ be a non-proper quadratic space, i.e., we have $f_q\equiv 0_\KK$. Then we  have $\mathrm{Def}(q)=L_0$, and $(L_0,\KK,q)$ is defective. In particular, a quadratic space is proper if it is non-defective.
\end{bem}

\begin{de}
Let $(L_0,\KK,q)$ be a quadratic space and let $R_0:=\mathrm{Def}(q)$. Then $(L_0,\KK,q)$ is a quadratic space \textit{of type $\mathit{F_4}$}\index{quadratic space (anisotropic)!of type $F_4$} if we have $\Char \KK=2$ and the following holds:
\begin{itemize}
\item $q(R_0)/q(\rho)$ is a subfield of $\KK$ for some $\rho\in R_0^*$.
\item For some complement $S_0$ of $R_0$ in $L_0$, the restriction of $q$ to $S_0$ has a norm splitting $(\EE,\cdot,\{v_1,v_2\})$ with constants $s_1,s_2$ such that $s_1s_2\in q(R_0)/q(\rho)$.
\end{itemize}
\end{de}

\begin{bem}
By (14.2) of \cite{TW}, the field $q(R_0)/q(\rho)$ is independent of the choice of $\rho\in R_0^*$. In particular, we have $$\FF:=g(R_0)/q(\epsilon)=q(R_0)$$ if $(L_0,\KK,q)$ is a quadratic space with basepoint $\epsilon\in R_0^*$.
\end{bem}

\begin{no}
Let $(L_0,\KK,q)$ be a quadratic space of type $F_4$ and let $\rho\in R_0^*$. We set $$\FF:=q(R_0)/q(\rho)\ .$$
\end{no}

\begin{lemma}\label{383}
A quadratic space $(L_0,\KK,q)$ of type $F_4$ is proper.
\end{lemma}

\begin{bew}
By definition, there is a complement $S_0$ of $R_0$ in $L_0$ such that the restriction of $q$ to $S_0$ has a norm splitting $(\EE,\cdot,\{v_1,v_2\})$ with constants $s_1,s_2$ such that $s_1s_2\in \FF$. Moreover, $q$ is anisotropic, hence $(S_0,\KK,q_{|S_0})$ is non-defective by (12.12) of \cite{TW}. In particular, $(S_0,\KK,q_{|S_0})$ and thus $(L_0,\KK,q)$ is proper.\qed
\end{bew}

\begin{lemma}\label{385}
Let $(L_0,\KK,q)$ be a quadratic space of type $F_4$ and let $(\tilde{\AA},\tilde{\FF},\tilde{N})$ be a quadratic space of type (m). Then we have $$(L_0,\KK,q)\not\cong (\tilde{\AA},\tilde{\FF},\tilde{N})$$ as quadratic spaces.
\end{lemma}

\begin{bew}
Assume $(L_0,\KK,q)\cong (\tilde{\AA},\tilde{\FF},\tilde{N})$. Since $(L_0,\KK,q)$ is proper by lemma \ref{383}, $(\tilde{\FF},\tilde{\AA},\tilde{N})$ is proper as well and thus (m)$\notin\{$(i),(ii)$\}$ since we have $\Char \KK=2$. But a quadratic space of type (m)$\in\{$(iii),\ldots,(v)$\}$ is non-defective by (20.15) of \cite{TW} and the proof of corollary \ref{126}, while a quadratic space of type $F_4$ is defective by definition since $\FF=q(R_0)/q(\rho)$ is a field. \qed
\end{bew}

\newpage

\chapter{Involutory Sets}

Quadrangles of purely involutory type are parametrized by proper involutory sets.

\section{Basic Definitions}
\begin{de}\label{302}\  \newglossaryentry{trace}{type=symbols,name={\ensuremath{\KK_\sigma}},description=the set of traces with respect to the involution $\sigma$, sort=involutory set}
\newglossaryentry{fix}{type=symbols,name={\ensuremath{\mathrm{Fix}(\sigma)}},description=the set of fixed points with respect to the involution $\sigma$, sort=involutory set}
\begin{itemize}
\item An \textit{involutory set}\index{involutory set} is a triple $(\KK,\KK_0,\sigma)$, where $\KK$ is a skew-field, $\sigma$ is an involution of $\KK$ and $\KK_0$ is an additive subgroup of $\KK$ such that
\begin{align*}
1_\KK\in \KK_0\ , && \{ a+a^\sigma\mid a\in \KK\}=:\gls{trace}\subseteq \KK_0\subseteq \gls{fix}\ , && \forall\ a\in \KK:\ a^\sigma \KK_0 a\subseteq \KK_0\ .
\end{align*}
\item An involutory set $(\KK,\KK_0,\sigma)$ is \textit{proper}\index{proper!involutory set}\index{involutory set!proper} if we have $\sigma\neq \id_\KK$ and $\langle \KK_0\rangle=\KK$ as rings.
\item Two involutory sets $(\KK,\KK_0,\sigma)$ and $(\tilde{\KK},\tilde{\KK}_0,\tilde{\sigma})$ are \textit{isomorphic}\index{isomorphism!of involutory sets} if there is an isomorphism $\phi:\KK\to\tilde{\KK}$ of skew-fields such that
\begin{align*}
\phi(\KK_0)=\tilde{\KK}_0\ , && \phi\circ\sigma=\tilde{\sigma}\circ\phi\ .
\end{align*}
\item Let $(\KK,\KK_0,\sigma),\ (\tilde{\KK},\tilde{\KK}_0,\tilde{\sigma})$ be two involutory sets. A \textit{Jordan homomorphism}\index{Jordan homomorphism} is an additive monomorphism $\gamma:\KK_0\to\tilde{\KK}_0$ such that
\begin{align*}
\gamma(1_{\KK})=1_{\tilde{\KK}}\ , && \forall\ x,y\in \KK_0:\ \gamma(xyx)=\gamma(x)\gamma(y)\gamma(x)\ .
\end{align*}
\end{itemize}
\end{de}

\begin{lemma}\label{157}
Let $(\KK,\KK_0,\sigma)$ be an involutory set. If $\KK_0$ is commutative, then $(\KK,\KK_0,\sigma)$ is non-proper. In particular, $\KK$ is non-commutative if $(\KK,\KK_0,\sigma)$ is proper.
\end{lemma}

\begin{bew}
Suppose that we have $\langle \KK_0\rangle=\KK$ and $\sigma\neq\id_\KK$. Then $\KK$ is commutative and  $\KK_0=\Fix(\sigma)\subsetneq \KK$ is a field, cf. remark (11.3) of \cite{TW}. But then we have
$$\KK=\langle \KK_0\rangle =\KK_0\neq \KK\qquad \lightning\ .$$
\qed
\end{bew}

\section{Jordan Isomorphisms of Involutory Sets}

The following result is essential for the classification of Jordan isomorphisms of pseudo-quadratic spaces on the one hand, and on the other hand, it is essential for the classification of a certain class of 443-foundations.

Parts of the solution can be deduced from paragraph 4.10 of \cite{K}. However, there are still some details which have to be worked out. But for reasons of brevity, we don't try to complete the proof at this point and suppose the result to be true.

\begin{satz}\label{175}
Let $(\KK,\KK_0,\sigma)$ be a proper involutory set, let $(\tilde{\KK},\tilde{\KK}_0,\tilde{\sigma})$ be an involutory set and let $\gamma:\KK_0\to \tilde{\KK}_0$ be a Jordan isomorphism such that $\gamma(1_\KK)=1_{\tilde{\KK}}$. Then $\gamma$ is induced by an isomorphism $$\phi:(\KK,\KK_0,\sigma)\to (\tilde{\KK},\tilde{\KK}_0,\tilde{\sigma})\ .$$ In particular, $(\tilde{\KK},\tilde{\KK}_0,\tilde{\sigma})$ is proper as well.
\end{satz}

\section{Involutory Sets of Quadratic Type}

Quadrangles parametrized by non-proper involutory sets can equally described by quadratic spaces, cf. chapter 38 of \cite{TW}. Although they don't appear explicitly in the description of the six families, they occur as substructures of proper pseudo-quadratic spaces, the parametrizing structures for quadrangles of purely pseudo-quadratic form type.

\begin{de} Let $(\AA,\FF,\sigma)$ be an involutory set (with $\AA$ possibly an alternative division ring) and $\KK:=Z(\AA)$. Then the involutory set is \textit{quadratic of type}\index{involutory set!of quadratic type}
\begin{enumerate}[label=(\roman*)]
\item if we have $\AA=\KK,\ \Char \KK=2,\ \KK^2\subseteq \FF\neq \KK$ and $\sigma=\id_\AA$,
\item if we have $\AA=\KK=\FF$ and $\sigma=\id_\AA$,
\item if we have $\AA=\KK$, $\KK/\FF$ is a separable quadratic extension and $\langle \sigma\rangle=\mathrm{Gal}(\KK/\FF)$,
\item if $\AA$ is a quaternion division algebra over $\KK$, $\FF=\KK$ and $\sigma=\sigma_s$,
\item if $\AA$ is an octonion division algebra over $\KK$, $\FF=\KK$ and $\sigma=\sigma_s$,
\end{enumerate}
where $\sigma_s$ denotes the standard involution, cf. proposition \ref{121}. We denote the corresponding norm by $N$ and the corresponding trace by $T$.
\end{de}

\begin{bem}
The following lemma gives a criterion when an involutory set is of quadratic type. We will need it for the classification of a certain class of 443-foundations.
\end{bem}

\begin{lemma}
Let $(\KK,\KK_0,\sigma)$ be an involutory set, let $\FF$ be a subfield of the center $Z(\KK)$ and let $\gamma\in \Aut(\KK,+)$ such that
\begin{align*}
\gamma(1_\KK)=1_\KK\ , && \forall\ x\in \KK:\ \gamma(x) x\in \FF\ .
\end{align*}
Then the following holds:
\begin{enumerate}[label=(\alph*)]
\item $\KK$ is quadratic over $\FF$.
\item\label{160} If $\KK$ is non-commutative, then $\KK$ is a quaternion division algebra and we have
\begin{align*}
\FF=Z(\KK)\ , && \gamma= \sigma_s\ .
\end{align*}
\end{enumerate}
\end{lemma}

\begin{bewzwei}\ 
\begin{enumerate}[label=(\alph*)]
\item Given $x\in \KK$, we have
\begin{align*}
\gamma(x)+x&=1_\KK+\gamma(x)+x+\gamma(x) x-1_\KK-\gamma(x) x \\
&=\big(1_\KK+\gamma(x)\big)(1_\KK+x)-1_\KK-\gamma(x) x \\
&=\gamma(1_\KK+x)(1_\KK+x)-1_\KK-\gamma(x) x \in \FF\ .
\end{align*}
If we set
$$m_x:=t^2-\big(\gamma(x) +x\big)\cdot t+ \gamma(x) x\in \FF[t]\ ,$$
we have
$$m_x(x)=x^2-\gamma(x) x-x^2+\gamma(x) x=0_\KK\ ,$$
thus $\KK$ is quadratic over $\FF$.
\item If $\KK$ is non-commutative, proposition \ref{121} implies that $\KK$ is a quaternion division algebra with $\FF=Z(\KK)$. Moreover, we have
$$\forall\ x\in \KK:\qquad \sigma_s(x)=T(x)-x=\gamma(x)+x-x=\gamma(x)$$
and therefore
$$\gamma=\sigma_s\ .$$
\end{enumerate}\qed
\end{bewzwei}

\newpage

\begin{bem}
The following results will be needed for the solution of the isomorphism problem for Moufang sets of pseudo-quadratic form type in §\ref{314}.
\end{bem}

\begin{lemma}\label{463}
Let $\DD$ be a skew-field and let $\EE$ be a maximal subfield of $\DD$. Then we have
$$C_\DD(\EE)=\EE\ .$$
\end{lemma}

\begin{bew}
This results from corollary (4.9) in chapter 8 of \cite{WSch}.\qed
\end{bew}

\begin{no}
Let $\HH$ be a quaternion division algebra with center $\KK$. Given $x\in \HH\sm \KK$, let $\EE_x:=\langle 1_\HH, x\rangle_\KK$ be the quadratic subfield of $\HH$ generated by $1_\HH$ and $x$, cf. (20.9) of \cite{TW}.
\end{no}

\begin{kor}\label{231}
Let $\HH$ be a quaternion division algebra and let $x,y\in \HH\setminus Z(\HH)$. Then we have
$$xy=yx\ \R\ \EE_x=\EE_y\ .$$
\end{kor}

\begin{bew}
By lemma \ref{463}, we have $$\EE_x=\langle 1,x\rangle_{Z(\HH)}\subseteq C_\HH(\EE_y)=\EE_y\ .$$
\qed
\end{bew}

\begin{lemma}\label{252}
For $i=1,2$, let $(\AA_i,\FF_i,\sigma_i)$ be quadratic of type (iii), (iv) or (v) with corresponding norms and traces $N_i$ and $T_i$, respectively, and let $\phi:\AA_1\to \AA_2$ be an isomorphism of alternative rings such that $\phi(\FF_1)=\FF_2$. Then we have
\begin{align*}\forall\ x\in \AA: && \phi\big(N_1(x)\big)=N_2\big(\phi(x)\big)\ , && \phi\big(T_1(x)\big)=T_2\big(\phi(x)\big)\ .\end{align*}
\end{lemma}

\begin{bew}
Given $x\in\AA_i$, we have
\begin{align*} x^2-T_1(x)x+N_1(x)=0_{\AA_1}\ , &&\phi(x)^2-T_2\big(\phi(x)\big)\phi(x)+N_2\big(\phi(x)\big)=0_{\AA_2}
\end{align*} and hence
\begin{align*}
0_{\AA_1}&=\phi^{-1}\big(\phi(x)^2-T_2(\phi(x))\phi(x)+N_2(\phi(x))\big)\\
&=\phi^{-1}\phi\big(x^2-\phi^{-1}(T_2(\phi(x)))x+\phi^{-1}(N_2(\phi(x)))\big)=x^2-\phi^{-1}\big(T_2(\phi(x))\big)x+\phi^{-1}\big(N_2(\phi(x))\big)\ .
\end{align*}
As the maps $T_1$ and $N_1$ are uniquely determined by the minimum equation, we obtain
\begin{align*}
N_1(x)&=\phi^{-1}\big(N_2(\phi(x))\big)\ , & \phi\big(N_1(x)\big)&=N_2\big(\phi(x)\big)\ , \\
T_1(x)&=\phi^{-1}\big(T_2(\phi(x))\big)\ , & \phi\big(T_1(x)\big)&=T_2\big(\phi(x)\big)\ .
\end{align*}
\qed
\end{bew}

\begin{kor}\label{251}
For $i=1,2$, let $(\AA_i,\FF_i,\sigma_i)$ be quadratic of type (iii), (iv) or (v) with corresponding norms and traces $N_i$ and $T_i$, respectively, and let $\phi:\AA_1\to \AA_2$ be an isomorphism of alternative rings such that $\phi(\FF_1)=\FF_2$. Then we have
$$\phi\circ \sigma_1 = \sigma_2\circ \phi\ .$$
\end{kor}

\begin{bew}
Given $x\in \AA_1$, we have
\begin{align*}
\phi\sigma_1(x)=\phi(\bar{x})=\phi\big(N_1(x)\cdot x^{-1}\big)=\phi\big(N_1(x)\big)\cdot \phi(x^{-1})=N_2\big(\phi(x)\big)\cdot \phi(x)^{-1}=\overline{\phi(x)}=\sigma_2\phi(x)\ .
\end{align*}\qed
\end{bew}

\begin{lemma}\label{253}
Let $\HH$ be a quaternion division algebra, let $\EE$ be a separable quadratic subfield and let $y\in \HH\sm \EE$. Then we have
$$y\in \EE^\bot\ \Leftrightarrow\ \forall\ x\in \EE:\ xy=y\bar{x}\ .$$
\end{lemma}

\begin{bewzwei}\
\begin{itemize}
\item[``$\Rightarrow$''] This holds by lemma \ref{250}.
\item[``$\Leftarrow$''] Let $e\in \EE^\bot$ and $y_1,y_2\in \EE$ such that $y=y_1+e y_2$. Given $x\in \EE\sm Z(\HH)$, we have
$$y_1x+ey_2\bar{x}=y_1x+e\bar{x}y_2=x(y_1+e\cdot y_2)=xy=y\bar{x}=y_1\bar{x}+ey_2\bar{x}\ ,$$
hence
\begin{align*}
y_1x=y_1\bar{x}\ , && y_1(x-\bar{x})=0_\HH\ , && y_1=0_\HH\ , && y=ey_2\in \EE^\bot\ .
\end{align*}
Notice that we use lemma \ref{116} several times.
\end{itemize}
\qed
\end{bewzwei}

\begin{bem}
The following results give a description for some extensions of isomorphisms between subfields of two given composition algebras.
\end{bem}

\begin{lemma}\label{225}
For $i=1,2$, let $\EE_i/\KK_i$ be a separable quadratic extension and let $t_i\in \EE_i\sm \KK_i$. Let $\phi:\KK_1\to\KK_2$ be an isomorphism of fields. Then the map
$$\tilde{\phi}: \EE_1\to \EE_2,\ x+t_1\cdot y\mapsto \phi(x)+t_2\cdot \phi(y)$$
is an isomorphism of fields if and only if we have
\begin{align*} \phi\big(N_1(t_1)\big)=N_2(t_2)\ , && \phi\big(T_1(t_1)\big)=T_2(t_2)\ .\end{align*}
\end{lemma}

\begin{bewzwei}\ 
\begin{itemize}
\item[``$\Rightarrow$''] This holds by lemma \ref{252}.
\item[``$\Leftarrow$''] This is a direct calculation using the minimum equation.
\end{itemize}\qed
\end{bewzwei}

\begin{lemma}\label{257}
For $i=1,2$, let $\HH_i=(\EE_i/\KK_i,\beta_i)$ be a quaternion division algebra and let $t_i\in \EE_i^\bot$. Let $\phi:\EE_1\to\EE_2$ be an isomorphism of fields. Then the map
$$\tilde{\phi}: \HH_1\to \HH_2,\ x+t_1\cdot y\mapsto \phi(x)+t_2\cdot \phi(y)$$
is an isomorphism of skew-fields if and only if we have
\begin{align*} \phi\big(N_1(t_1)\big)=N_2(t_2)\ .\end{align*}
\end{lemma}

\begin{bewzwei}\ 
\begin{itemize}
\item[``$\Rightarrow$''] This holds by lemma \ref{252}.
\item[``$\Leftarrow$''] This is a direct calculation using the minimum equation and lemma \ref{253}.
\end{itemize}\qed
\end{bewzwei}

\begin{bem}
We finally recall the list of quadratic spaces corresponding to involutory sets of quadratic type.
\end{bem}

\begin{de}\label{462}
A quadratic space $(L_0,\KK,q)$ is \textit{of type}\index{quadratic space (anisotropic)!of type (m)}
\begin{enumerate}[label=(\roman*)]
\item if $\FF:=L_0$ is a field with $\Char \FF=2,\ \FF^2\subseteq \KK\neq \FF$, $\sigma=\id_\FF$ and $q=N^\FF_\KK$,
\item if we have $L_0=\KK$, $\sigma=\id_\KK$ and $q=N^\KK_\KK$,
\item if $\EE:=L_0$ is a field, $\EE/\KK$ is a separable quadratic extension, $\langle \sigma\rangle=\mathrm{Gal}(\EE/\KK)$ and $q=N^\EE_\KK$,
\item if $\HH:=L_0$ is a quaternion division algebra over $\KK$, $\sigma=\sigma_s$ and $q=N^\HH_\KK$,
\item if $\OO:=L_0$ is an octonion division algebra over $\KK$, $\sigma=\sigma_s$ and $q=N^\OO_\KK$.
\end{enumerate}
\end{de}

\begin{bem}\ 
\begin{enumerate}[label=(\alph*)]
\item By corollary \ref{432}, a quadratic space $(L_0,\KK,q)$ with $\dim_\KK L_0\leq 2$ is of type (i)-(iii).
\item We will see that the Moufang sets of a quadratic space $(\AA,\FF,N^\AA_\FF)$ of type (m) and of the corresponding alternative division ring $\AA$ coincide, cf. lemma \ref{456}.
\end{enumerate}
\end{bem}

\chapter{Indifferent Sets}

Quadrangles of purely indifferent type are parametrized by proper indifferent sets.

\begin{de}\ 
\begin{itemize}
\item An \textit{indifferent set}\index{indifferent set} is a triple $(\KK,\KK_0,\LL_0)$, where $\KK$ is a field with $\Char \KK=2$ and $\KK_0$ and $\LL_0$ are additive subgroups of $\KK$ containing $1_\KK$ such that
\begin{align*}
\KK_0^2\LL_0\subseteq \LL_0\ , && \LL_0\KK_0\subseteq \KK_0\ , && \langle \KK_0\rangle=\KK\ \textrm{as rings}\ .
\end{align*}
\item An indifferent set is \textit{proper}\index{proper!indifferent set}\index{indifferent set!proper} if we have $\KK_0\neq \KK$ and $\LL_0\neq \LL:=\langle \LL_0\rangle$.
\item Two indifferent sets $(\KK,\KK_0,\LL_0)$ and $(\tilde{\KK},\tilde{\KK}_0,\tilde{\LL}_0)$ are \textit{isomorphic}\index{isomorphism!of indifferent sets} if there is an isomorphism $\phi:\KK\to\tilde{\KK}$ of fields such that
\begin{align*}
\phi(\KK_0)=\tilde{\KK}_0\ , && \phi(\LL_0)=\tilde{\LL}_0\ .
\end{align*}
\end{itemize}
\end{de}

\begin{lemma}
Let $(\KK,\KK_0,\LL_0)$ be an indifferent set. Then $(\LL,\LL_0,\KK_0^2)$ is an indifferent set.
\end{lemma}

\begin{bew}
This is (10.2) of \cite{TW}.\qed
\end{bew}

\begin{de}
The \textit{opposite of an indifferent set $\mathit{(\KK,\KK_0,\LL_0)}$}\index{indifferent set!opposite}\index{opposite!indifferent set} is the indifferent set $(\LL,\LL_0,\KK_0^2)$.
\end{de}

\begin{lemma}\label{395}
Given a proper indifferent set $(\KK,\KK_0,\LL_0)$, its opposite $(\LL,\LL_0,\KK_0^2)$ is proper.
\end{lemma}

\begin{bew}
Since $(\KK,\KK_0,\LL_0)$ is proper, we have $\LL_0\neq \LL$. We have to show $\KK_0^2\neq \langle \KK_0^2\rangle$. By remark (10.8) of \cite{TW}, the opposite of $(\LL,\LL_0,\KK_0^2)$ is $(\KK^2,\KK_0^2,\LL_0^2)$, and we have
$$(\KK^2,\KK_0^2,\LL_0^2)\cong (\KK,\KK_0,\LL_0)$$
as indifferent sets. In particular, we have $\langle \KK_0^2\rangle=\KK^2$. Moreover, we have $\KK_0^2\neq \KK^2$ since the map $\mathrm{Fr}:\KK\to \KK,\ x\mapsto x^2$ is injective and $\KK_0\neq \KK$ by assumption. We finally obtain
$$\KK_0^2\neq \KK^2=\langle \KK_0^2\rangle\ .$$
\qed
\end{bew}

\chapter{Pseudo-Quadratic Spaces}

Quadrangles of purely pseudo-quadratic form type are parametrized by proper pseudo-quadratic spaces.

\section{Basic Definitions and Basic Properties}
First of all we give the basic definitions and introduce the Moufang set of pseudo-quadratic form type, more precisely, its associated group, corresponding to a pseudo-quadratic space.

\begin{de}\ 
\begin{itemize}
\item An \textit{(anisotropic) right (resp. left) pseudo-quadratic space}\index{pseudo-quadratic space (anisotropic)} is a quintuple $\Xi=(\KK,\KK_0,\sigma,L_0,q)$ such that $(\KK,\KK_0,\sigma)$ is an involutory set, $L_0$ is a right (resp. left) vector space over $\KK$ and $q$ is an \textit{(anisotropic) pseudo-quadratic form on $\mathit{L_0}$ with respect to $\mathit{\sigma}$}\index{pseudo-quadratic form (anisotropic)}, i.e., there is a skew-hermitian form $f$ on $L_0$ such that the following holds:
\begin{enumerate}[leftmargin=25pt, label=(P\arabic*)]
\item\label{183} $\forall\ a,b\in L_0:\ q(a+b)\equiv q(a)+q(b)+f(a,b)\mod \KK_0$,
\item\label{315}$\forall\ a\in L_0,\ t\in \KK:\ q(at)\equiv t^\sigma q(a)t\mod \KK_0\ (\textrm{resp.}\ q(ta)\equiv tq(a)t^\sigma \mod \KK_0)$,
\item $q(a)\equiv 0_\KK\mod \KK_0\ \Leftrightarrow\ a=0_{L_0}$.
\end{enumerate}
\item A pseudo-quadratic space $(\KK,\KK_0,\sigma,L_0,q)$ is \textit{proper}\index{proper!pseudo-quadratic space}\index{pseudo-quadratic space (anisotropic)!proper} if we have $\sigma\neq \id_\KK$, $L_0\neq \{0\}$ and if the associated skew-hermitian form $f$ is non-degenerate.
\item Two pseudo-quadratic space $\Xi$ and $\tilde{\Xi}$ are \textit{isomorphic}\index{isomorphism!of pseudo-quadratic spaces} if there is an isomorphism
\begin{align*}
\Phi=(\p,\phi): (L_0,\KK)\to (\tilde{L}_0,\tilde{\KK})
\end{align*}
of vector spaces such that $\phi:(\KK,\KK_0,\sigma)\to (\tilde{\KK},\tilde{\KK}_0,\tilde{\sigma})$ is an isomorphism of involutory sets and such that
$$\phi \circ q\equiv \tilde{q}\circ \p \mod \tilde{\KK}_0\ .$$
\end{itemize}
\end{de}

\begin{lemma}\label{177}
The skew-hermitian form $f$ is uniquely determined by \ref{183}  and satisfies
\begin{align*}
\forall\ a\in {L}_0:\qquad {f}(a,a)={q}(a)-{q}(a)^{{\sigma}}\ .
\end{align*}
\end{lemma}

\begin{bew}
This is (11.19) of \cite{TW}. Notice that we have $\KK_0\neq \KK$ since $\Xi$ is proper.\qed
\end{bew}

\begin{kor}\label{303} Given an isomorphism $\Phi=(\p,\phi):\Xi\to\tilde{\Xi}$ of pseudo-quadratic spaces, we have
$$\forall\ a,b\in L_0:\qquad \tilde{f}\big(\p(a),\p(b)\big)=\phi\big(f(a,b)\big)\ .$$
\end{kor}

\begin{bew}
We have to show that 
\begin{align*}
\forall\ a,b\in L_0:\qquad \phi^{-1}\big(\tilde{f}(\p(a),\p(b))\big)\equiv q(a+b)-q(a)-q(b) \mod \KK_0\ .
\end{align*}
Given $a,b\in L_0$, we have
\begin{align*} \phi^{-1}\big(\tilde{f}(\p(a),\p(b))\big)&\in \phi^{-1}\big( \tilde{q}(\p(a+b))-\tilde{q}(\p(a))-\tilde{q}(\p(b))+\tilde{\KK}_0 \big) \\
&=\phi^{-1}\big( \phi(q(a+b))-\phi(q(a))-\phi(q(b))+\tilde{\KK}_0\big) \\
&=q(a+b)-q(a)-q(b)+\KK_0\ .
\end{align*}\qed
\end{bew}

\begin{bem}\label{320} Let $\Xi$ be a pseudo-quadratic space and let $a\in L_0$.
\begin{itemize}
\item Assume $\Char \KK\neq 2$. Then we have
\begin{align*}
4q(a)\equiv q(2a)\equiv 2q(a)+f(a,a) \mod \KK_0\ , && q(a)\equiv \frac{f(a,a)}{2} \mod \KK_0\ .
\end{align*}
\item Assume $\Char \KK=2$. Then we have $f(a,a)=q(a)+q(a)^\sigma \in \KK_0$.
\end{itemize}
\end{bem}

\begin{de}\label{190} \newglossaryentry{T}{type=symbols,name={\ensuremath{T(\Xi)}},description=associated group with the Moufang set with respect to the pseudo-quadratic space $\Xi$, sort=pseudoquadratic space}
Given a pseudo-quadratic space $\Xi$, we set
$$T:=\gls{T}:=\{ (a,t)\in L_0\times \KK \mid q(a)-t\in \KK_0\}\ .$$
\end{de}

\begin{no}
Throughout the rest of this chapter, let $\Xi$ be a proper pseudo-quadratic space and let $T$ be the corresponding set as in definition \ref{190}.
\end{no}

\begin{lemma}\label{182} 
Given $(a,t)\in T,\ k\in \KK$, we have
$$(a,t+k)\in T\ \Leftrightarrow\ k\in\KK_0\ .$$
\end{lemma}

\begin{bew}
Given $(a,t)\in T$ and $k\in \KK$, we have
\[ (a,t+k)\in T\ \Leftrightarrow\ q(a)-t-k\in \KK_0\ \Leftrightarrow\ k\in q(a)-t+\KK_0=\KK_0\ . \]
\qed
\end{bew}

\begin{kor}\label{187}
Given $(a,t)\in T$, we have $ f(a,a)=t-t^\sigma$.
\end{kor}

\begin{bew} 
Let $(a,t)\in T$. By lemma \ref{182}, there is an element $k\in\KK_0$ such that $t=q(a)+k$, hence
\[t-t^\sigma=q(a)+k-q(a)^\sigma-k^\sigma=q(a)-q(a)^\sigma=f(a,a)\]
by lemma \ref{177}.\qed
\end{bew}

\begin{kor}\label{321}\label{184} Let
$$ \cdot:T\times T\to T,\  (a,t)\cdot (b,v):=\big(a+b,t+v+f(b,a)\big)\ .$$
Then $(T,\cdot)$ is a group with $Z(T)=\{ (0_{L_0},t) \mid t\in \KK_0\}\cong \KK_0$ and $(a,t)^{-1}=(-a,-t^\sigma)$ for each $(a,t)\in T$.
\end{kor}

\begin{bew}
This results from (11.24) and (38.10) of \cite{TW}. Notice that $\Xi$ is proper.\qed
\end{bew}

\begin{bem}
In the following, we don't distinguish between $Z(T)$ and $\KK_0$, i.e.,  we consider $\KK_0$ to be a subset of $T$ via the above identification.
\end{bem}

\begin{lemma}\label{164}
Given $(a,t),(b,v)\in T$, we have
$$(a,t)\cdot (b,v)\in \KK_0\ \Leftrightarrow\ a=-b\ .$$
\end{lemma}

\begin{bew}
Given $(a,t),(b,v)\in T$, we have
\[
(a,t)\cdot (b,v)\in \KK_0\ \Leftrightarrow\ \big(a+b,t+v+f(b,a)\big)\in \KK_0\ \Leftrightarrow\ a+b=0_{L_0}\ \Leftrightarrow\ a=-b\ .
\]\qed
\end{bew}

\section{Jordan Isomorphisms}
We collect some first results about isomorphisms preserving the Moufang set structure.

\begin{no}
Throughout this paragraph, let $\tilde{\Xi}$ be an additional pseudo-quadratic space, let $\tilde{T}$ be the corresponding group as in definition \ref{190} and let $\gamma:T\to \tilde{T}$ be an isomorphism of groups.
\end{no}

\begin{bem} As $T$ is non-abelian and $\gamma$ is an isomorphism of groups, we have $\tilde{\sigma}\neq \id_{\tilde{\KK}}$, $\tilde{L}_0\neq \{0\}$, and the associated skew-hermitian form is not identically zero, i.e., $\tilde{\Xi}$ is pre-proper. Moreover, it is proper if we have $\Char \tilde{\KK}\neq 2$, and, if $\Char \tilde{\KK}=2$, we may replace $\tilde{\Xi}$ by a proper pseudo-quadratic space, cf. definition (35.5) of \cite{TW}.

Thus we get a satisfying solution of the isomorphism problem for Moufang sets of pseudo-quadratic form type if we restrict to the case of two proper pseudo-quadratic spaces.
\end{bem}

\begin{kor}\label{168}
We have $\gamma(\KK_0)=\tilde{\KK}_0.$
\end{kor}

\begin{bew}
We have $\gamma(\KK_0)=\gamma\big(Z(T)\big)=Z(\tilde{T})=\tilde{\KK}_0$.
\qed
\end{bew}

\begin{lemma}
\label{167} Let $\p_1:T \to \tilde{L}_0,\ \p_2:T\to \tilde{\KK}$ defined by
$$\gamma(a,t)=\big(\p_1(a,t),\p_2(a,t)\big)\ .$$
Then we have
$$\forall\ (a,t)\in T:\qquad \p_1(a,t)=\p_1(a)\ .$$
Moreover, the $\p_1:L_0\to\tilde{L}_0$ is an isomorphism of groups.
\end{lemma}

\begin{bewzwei}\ 
\begin{itemize}
\item Given $(a,t),(a,u)\in T$, we have
$$(a,t)\cdot (a,u)^{-1}=(a,t)\cdot (-a,-u^{\sigma})\in \KK_0=Z(T)$$
by lemma \ref{164}; therefore,
\begin{align*}\big(\p_1(a,t),\p_2(a,t)\big)\cdot \big(-\p_1(a,u),-\p_2(a,u)^{\tilde{\sigma}}\big)&=\gamma_2(a,t)\cdot \gamma_2(a,u)^{-1}\\
&=\gamma_2\big((a,t)\cdot (a,u)^{-1}\big)\in Z(\tilde{T})=\tilde{\KK}_0 \end{align*}
and thus $\p_1(a,t)=\p_1(a,u)$ by lemma \ref{164} again.
\item We have $\big(a,q(a)\big)\in T$ for each $a\in L_0$, hence $\p_1$ is well-defined.
\item As the multiplication in $T$ is additive in the first component, $\p_1$ is additive.
\item Let $(a,t)\in T$ such that $a\in \Kern \p_1$. Then we have
$$\gamma(a,t)=\big(0_{L_0},\p_2(a,t)\big)\in \tilde{\KK}_0\subseteq Z(\tilde{T})\ ,$$
hence $(a,t)\in Z(T)=\KK_0$ by corollary \ref{184} and thus $a=0_{L_0}$.
\item We have $\big(a,\tilde{q}(a)\big)\in \tilde{T}$ for each $a\in \tilde{L}_0$. As $\gamma$ is surjective, $\p_1$ is surjective as well.
\end{itemize}\qed
\end{bewzwei}

\begin{de}\label{354}
A \textit{Jordan isomorphism}\index{Jordan isomorphism} is an isomorphism of groups $\gamma:T\to\tilde{T}$ with $\gamma(0,1_\KK)=(0,1_{\tilde{\KK}})$ such that the maps $\p_1:T \to \tilde{L}_0,\ \p_2:T\to \tilde{\KK}$ defined by
$$\gamma(a,t)=\big(\p_1(a,t)=\p_1(a),\p_2(a,t)\big)$$
satisfy 
\begin{align}
&\p_1\big(b t^\sigma-at^{-1}f(a,b)t^\sigma\big)=\p_1(b)\p_2(a,t)^{\tilde{\sigma}}-\p_1(a)\p_2(a,t)^{-1}\tilde{f}\big(\p_1(a),\p_1(b)\big)\p_2(a,t)^{\tilde{\sigma}}\ , \label{169} \\
&\p_2\big(b t^\sigma-at^{-1}f(a,b)t^\sigma,t vt^\sigma\big)=\p_2(a,t)\p_2(b,v)\p_2(a,t)^{\tilde{\sigma}} \label{170}
\end{align}
for all $(a,t),(b,v)\in T$.
\end{de}

\begin{de}\label{eq:hua} \newglossaryentry{huap}{type=symbols,name={\ensuremath{h_{(a,t)}}},description=Hua-map with respect to ${(a,t)\in T}$, sort=pseudoquadratic space}
Given $(a,t)\in T^*$, the \textit{Hua-map with respect to $\mathit{(a,t)}$}\index{Hua map} is
$$\gls{huap}: T\to T,\ (b,v)\mapsto \big(b t^\sigma-at^{-1}f(a,b)t^\sigma,t vt^\sigma\big)\ .$$
\end{de}

\begin{bem}
A Jordan isomorphism could equally be defined as an isomorphism of groups ${\gamma}:T\to \tilde{T}$ satisfying $\gamma(0,1_\KK)=(0,1_{\tilde{\KK}})$ and preserving the Hua-maps, i.e, given $(a,t),(b,v)\in T^*$, we have
$$ {\gamma}\big( h_{(a,t)}(b,v)\big)=\tilde{h}_{{\gamma}(a,t)} \big({\gamma}(b,v)\big)\ .$$
\end{bem}

\begin{lemma}\label{292}
We have $h_{(a,t)}\in \Aut(T)$ for each $(a,t)\in T^*$.
\end{lemma}

\begin{bew}
This is theorem 2 of \cite{DW}.\qed
\end{bew}

\begin{no}
Since the first component in $h_{(a,t)}(b,v)$ is independent of $v$, we may restrict $h_{(a,t)}$ to the first component, i.e., given $b\in L_0$, we set
$$h_{(a,t)}(b):=b t^\sigma-at^{-1}f(a,b)t^\sigma\ .$$
\end{no}

\begin{lemma} Let $\gamma:T\to\tilde{T}$ be a Jordan isomorphism. Then the following holds:
\begin{enumerate}[label=(\alph*)]
\item \label{172} Given $b\in L_0,\ (0_{L_0},t)\in\KK_0$, we have  $\p_1(b\cdot t)=\p_1(b)\cdot \p_2(0_{L_0},t)$.
\item \label{171} Given $(a,t)\in T,\ (0_{L_0},v)\in \KK_0$, we have $\p_2(0_{L_0},t vt^\sigma)=\p_2(a,t) \cdot \p_2(0_{L_0},v)\cdot \p_2(a,t)^{\tilde{\sigma}}$.
\item \label{204} The map $\phi:\KK_0\to\tilde{\KK}_0,\ t\mapsto \p_2(0_{L_0},t)$ is a Jordan isomorphism as in definition \ref{302}.
\item \label{221} 
Given $(a,t)\in T$ and $s\in \KK_0$, we have $\p_2(a,t+s)=\p_2(a,t)+\phi(s)$.
\end{enumerate}
\end{lemma}

\begin{bewzwei}\ 
\begin{enumerate}[label=(\alph*)]
\item This results from identity \eqref{169} with $a=0_{L_0}$. Notice that we have 
\begin{align*} \KK_0\subseteq \Fix(\sigma)\ ,&& \gamma(\KK_0)=\tilde{\KK}_0\subseteq \Fix(\tilde{\sigma})\ .\end{align*}
\item This results from identity \eqref{170} with $b=0_{L_0}$.
\item By corollary \ref{168}, $\phi$ is an isomorphism of groups. Given $(0_{L_0},t),(0_{L_0},v)\in \KK_0$, we have
$$\p_2(0,tvt)=\p_2(0,t v t^\sigma)=\p_2(0,t)\cdot \p_2(0,v)\cdot \p_2(0,t)^{\tilde{\sigma}}=\p_2(0,t)\cdot \p_2(0,v)\cdot \p_2(0,t)\ .$$
\item Given $(a,t)\in T$ and $s\in \KK_0$, we have
\begin{align*}
  \big(\p_1(a),\p_2(a,t+s)\big)&=\gamma(a,t+s)=\gamma(a,t)\cdot \gamma(0,s)\\
  &=\big(\p_1(a),\p_2(a,t)\big)\cdot \big(0,\phi(s)\big)=\big(\p_1(a),\p_2(a,t)+\phi(s)\big)\ .
  \end{align*}
\end{enumerate}\qed
\end{bewzwei}

\newpage

\addtocontents{toc}{\protect\newpage}
\addtocontents{toc}{\noindent\protect\mbox{}\protect\hrulefill\par}
\part{Jordan Isomorphisms of Pseudo-Quadratic Spaces}
\addtocontents{toc}{\noindent\protect\mbox{}\protect\hrulefill\par}

\noindent For the classification of 443-foundations in part \ref{201}, we will need a part of the solution of the isomorphism problem for Moufang sets of pseudo-quadratic form type. We state a version of the result at this point and break up the proof into several steps. The main point is that under certain conditions, a Jordan isomorphism of pseudo-quadratic spaces is induced by an isomorphism of the corresponding pseudo-quadratic spaces, and we handle the possible cases one by one.

Then we have a closer look at the exceptions, which only occur in small dimensions, before we finally show that each of the appearing maps really induces a Jordan isomorphism.\\

\chapter{A Partial Result}

\begin{satz}\label{283}
Let $\Xi$ and $\tilde{\Xi}$ be proper pseudo-quadratic spaces, let $\gamma:T\to\tilde{T}$ be a Jordan isomorphism and suppose that one of the following holds:
\begin{enumerate}[label=(\roman*)]
\item The involutory set $(\KK,\KK_0,\sigma)$ is proper.
\item The involutory set $(\KK,\KK_0,\sigma)$ is quadratic of type (iii) or (iv) and $\dim_{\KK} L_0\geq 3$.
\item\label{254} The involutory set $(\KK,\KK_0,\sigma)$ is quadratic of type (iii) or (iv), $\dim_\KK L_0\leq 2$ and $\tilde{\KK}\cong \KK\not\cong \FF_4$.
\end{enumerate}
Then $\gamma$ is induced by an isomorphism $\Phi:\Xi\to\tilde{\Xi}$ of pseudo-quadratic spaces. 
\end{satz}

\begin{bew}
This results from theorem \ref{174}, theorem \ref{178}, theorem \ref{271} and theorem \ref{270}.\qed
\end{bew}

\begin{no}\label{276}
Throughout this part, let $\Xi$ and $\tilde{\Xi}$ be proper pseudo-quadratic spaces, let $\gamma:T\to \tilde{T}$ be a Jordan isomorphism and let $\phi:\KK_0\to\tilde{\KK}_0,\ t\mapsto \p_2(0_{L_0},t)$.
\end{no}

\chapter{The Involutory Set Is Proper}

The first case is that of a proper involutory set $(\KK,\KK_0,\sigma)$. The solution of the isomorphism problem for proper involutory sets yields an isomorphism $\phi:\KK\to \tilde{\KK}$. From identity \eqref{169} and the fact that $\KK_0$ generates $\KK$ as a ring, we deduce that $(\p_1,\phi):(L_0,\KK)\to (\tilde{L}_0,\tilde{\KK})$ is an isomorphism of vector spaces.

Finally, we show that the second component of $\gamma$ is induced by $\phi$, using identity \eqref{170} and the fact that the dimension of $\langle\KK_0\rangle_{Z(\KK)}$ over $Z(\KK)$ is at least 2.

\begin{no} Throughout this chapter, let $(\KK,\KK_0,\sigma)$ be proper.
\end{no}

\begin{lemma}\label{166}
Let $(\hat{\KK},\hat{\KK}_0,\hat{\sigma})$ be a proper involutory set and let $s,t\in \hat{\KK}$ be such that
$$\forall\ u\in \hat{\KK}_0:\qquad sus^{\hat{\sigma}}=tut^{\hat{\sigma}}\ .$$
Then we have $t^{-1}s\in Z(\hat{\KK})$.
\end{lemma}

\begin{bew}
First of all we notice that
$$ss^{\hat{\sigma}}=s\cdot 1_{\hat{\KK}}\cdot s^{\hat{\sigma}}=t\cdot 1_{\hat{\KK}}\cdot t^{\hat{\sigma}}=tt^{\hat{\sigma}}$$
and thus
$$\forall\ u\in \hat{\KK}_0:\qquad sus^{-1}=sus^{{\hat{\sigma}}}(ss^{{\hat{\sigma}}})^{-1}=tut^{{\hat{\sigma}}}(tt^{{\hat{\sigma}}})^{-1}=tut^{-1}\ .$$
It follows that
$$\forall\ u\in \langle \hat{\KK}_0\rangle=\hat{\KK}:\qquad t^{-1}su(t^{-1}s)=u\ ,$$
hence $t^{-1}s\in Z(\hat{\KK})$. \qed
\end{bew}

\newpage

\begin{bem}\label{173}
 As $(\KK,\KK_0,\sigma)$ is proper, the map $\phi:\KK_0\to\tilde{\KK}_0$ as defined in notation \ref{276} is induced by an isomorphism $\phi:(\KK,\KK_0,\sigma)\to (\tilde{\KK},\tilde{\KK}_0,\tilde{\sigma})$ of involutory sets, and $(\tilde{\KK},\tilde{\KK}_0,\tilde{\sigma})$ is proper, cf. theorem \ref{175}.
 \end{bem}
 
 \begin{kor}\label{180}
Let $(a,t)\in T$. Then there is an element $\lambda\in Z(\tilde{\KK})$ such that
\begin{align*} \p_2(a,t)=\lambda\cdot \phi(t)\ .\end{align*}
\end{kor}

\begin{bew}
By lemma \ref{171}, we have
\begin{align*} \phi(t)\phi(v)\phi(t)^{\tilde{\sigma}}=\phi(tv t^\sigma )=\p_2(0,t vt^\sigma)=\p_2(a,t) \p_2(0,v) \p_2(a,t)^{\tilde{\sigma}}=\p_2(a,t) \phi(v)\p_2(a,t)^{\tilde{\sigma}}\end{align*}
for all $(a,t)\in T$, $v\in {\KK}_0$ and thus
\begin{align*} \phi(t)\cdot \tilde{v}\cdot \phi(t)^{\tilde{\sigma}}=\p_2(a,t)\cdot \tilde{v}\cdot \p_2(a,t)^{\tilde{\sigma}}\end{align*}
for all $(a,t)\in T$, $\tilde{v}\in \tilde{\KK}_0$. Now the assertion results from lemma \ref{166}. \qed
\end{bew}

\begin{lemma}\label{188} We have $$\dim_{Z(\KK)} \langle \KK_0 \rangle_{Z(\KK)}\geq 2\ .$$
\end{lemma}

\begin{bew}
As $(\KK,\KK_0,\sigma)$ is proper, there is an element $x\in \KK_0\setminus Z(\KK)$ by lemma \ref{157}.
\qed
\end{bew}

\begin{kor}\label{189}
Let $(a,s)\in T$. Then there is an element $t\in \KK$ such that 
\begin{align*}
(a,t)\in T\ ,&&  s \notin \langle t^\sigma\rangle_{Z(\KK)}\ .
\end{align*}
\end{kor}

\begin{bew}
We have $(a,-s^\sigma)\in T$. By lemma \ref{188}, there is an element $u\in \KK_0$ such that $s\notin \langle u\rangle_{Z(\KK)}$, thus
$$s\notin \langle -s+u\rangle_{Z(\KK)}= \langle(-s+u)^{\sigma^2}\rangle_{Z(\KK)}=\langle (-s^\sigma+u)^\sigma\rangle_{Z(\KK)}\ ,$$
and by lemma \ref{182}, we have $(a, -s^\sigma+u)\in T$.\qed
\end{bew}

\begin{lemma}\label{179}
Given  $(a,s)\in T$, we have
$$\p_2(a,s)=\phi(s)\ .$$
\end{lemma}

\begin{bew} Let $(a,s)\in T$ and let $t\in \KK$ be as in corollary \ref{189}. Notice that $s+t^\sigma\in \KK_0$. By corollary \ref{180}, there are  $\lambda_1,\lambda_2\in Z(\tilde{\KK})$ such that
\begin{align*}
\p_2(a,s)=\lambda_1 \cdot {\phi}(s)\ , && \p_2(-a,t)=\lambda_2\cdot {\phi}(t)\ .
\end{align*}
Observing corollary \ref{187}, we obtain
\begin{align*}
\big(0,{\phi}(s)+\phi(t^\sigma)\big)&=\gamma(0,s+t^\sigma)=\gamma\big(0,s+t-f(a,a)\big) \\
&=\gamma(a,s)\cdot \gamma(-a,t)=\big(\p_1(a),\p_2(a,s)\big)\cdot \big(\p_1(-a),\p_2(-a,t)\big) \\
&=\big(0,\p_2(a,s)+\p_2(-a,t)-\tilde{f}(\p_1(-a),\p_1(-a))\big)\\
&=\big(0,\p_2(a,s)+\p_2(-a,t)^{\tilde{\sigma}}\big)=\big(0, \lambda_1\cdot {\phi}(s)+\lambda_2\cdot {\phi}(t^\sigma)\big)
\end{align*}
and thus $\lambda_1=1_{\tilde{\KK}}$ by the linear independence of ${\phi}(s)$ and ${\phi}(t^\sigma)$, cf. remark \ref{173}.\qed
\end{bew}

\begin{satz}\label{174} If $(\KK,\KK_0,\sigma)$ is proper, the map $\Phi:\Xi\to\tilde{\Xi}$ defined by
\begin{align*}
\Phi:=(\p_1,\phi): (L_0,\KK)\to(\tilde{L},\tilde{\KK}),\ (a,t)\mapsto \big(\p_1(a), \phi(t)\big)
\end{align*}
is an isomorphism of pseudo-quadratic spaces satisfying
$$\forall\ (a,t)\in T:\qquad \Phi(a,t)=\gamma(a,t)\ .$$
\end{satz}

\begin{bewzwei}\ 
\begin{itemize}
\item By remark \ref{173}, the map $\phi:(\KK,\KK_0,\sigma)\to (\tilde{\KK},\tilde{\KK}_0,\tilde{\sigma})$ is an isomorphism of involutory sets.
\item By lemma \ref{167}, the map $\p_1:L_0\to \tilde{L}_0$ is an isomorphism of groups. Given $a\in L_0$ and $t_1,\ldots,t_n\in \KK_0$, we have
\begin{align*} \p_1\big(a\cdot t_1\cdots t_n\big)&=\p_1(a)\cdot \phi(t_1)\cdots \phi(t_n)=\p_1(a)\cdot \phi( t_1\cdots t_n)  \end{align*}
by lemma \ref{172}, thus $(\p_1,\phi)$ is an isomorphism of vector spaces as we have $\langle \KK_0\rangle =\KK$.
\item Given $(a,t)\in T$, we have
\begin{align*}
\tilde{q}\big(\p_1(a)\big)&\in \p_2(a,t)+\tilde{\KK}_0=\phi(t)+\tilde{\KK}_0= \phi\big(t+\KK_0\big)=\phi\big(q(a)+\KK_0\big)=\phi\big(q(a)\big)+\tilde{\KK}_0
\end{align*}
by lemma \ref{179}.
\item Given $(a,t)\in T$, we have
$$\Phi(a,t)=\big(\p_1(a),\phi(t)\big)=\big(\p_1(a),\p_2(a,t)\big)=\gamma(a,t)\ .$$
\end{itemize}\qed
\end{bewzwei}

\begin{bem}\label{191}
If $(\KK,\KK_0,\sigma)$ is non-proper, then $(\AA,\FF,\sigma):=(\KK,\KK_0,\sigma)$ is quadratic of type (iii) or (iv), cf. (38.14) of \cite{TW}. It would be nice to have the same result in this case. However, this is false in general, thus we will need some additional assumptions.
\end{bem}

\chapter{The Involutory Set Is of Quadratic Type}
The second case is that of an involutory set which is of quadratic type and $\dim_\KK L_0\geq 3$. The crucial step is to show that $\p_1$ maps orthogonal vectors to orthogonal vectors and hence subspaces to subspaces. In particular, we may apply the fundamental theorem of projective geometry to get an isomorphism $(\p_1,\phi):(L_0,\KK)\to (\tilde{L}_0,\KK_0)$ of vector spaces. 

In order to prove this, we introduce a technical condition concerning orthogonality which is preserved by the Hua-maps. However, this condition is only enough to handle separable elements which results in a separate treatment of inseparable elements.

Finally, we show that the second component of $\gamma$ is induced by $\phi$, using identity \eqref{169} and the fact that the dimension of $L_0$ over $\KK$ is at least 3, which ensures the existence of enough orthogonal vectors.

\begin{no}
Throughout this chapter, we suppose $(\AA,\FF,\sigma):=(\KK,\KK_0,\sigma)$ \big(and therefore $(\tilde{\AA},\tilde{\FF},\tilde{\sigma}):=(\tilde{\KK},\tilde{\KK}_0,\tilde{\sigma})$\big) to be quadratic of type (iii) or (iv).
\end{no}

\begin{de}\
\begin{itemize}
\item An element $a\in L_0$ with $f(a,a)=0_\AA$ is called \textit{inseparable}\index{inseparable}.
\item Otherwise, it is called \textit{separable}\index{separable}.
\end{itemize}
\end{de}

\newpage

\begin{bem}\ 
\begin{enumerate}[label=(\alph*)]
\item If $a\in L_0$ is separable, we have
\begin{align*}
q(a)-q(a)^\sigma=f(a,a)\neq 0_\AA\ , && q(a)^\sigma\neq q(a)\ .
\end{align*}
\item \label{195} If $a\in L_0^*$ is inseparable, we have
\begin{align*} q(a)^\sigma=q(a)\ , && q(a)\in \Fix(\sigma)\sm \FF\ ,
\end{align*}
which implies $\Char \AA=2$. Moreover, $(\AA,\FF,\sigma)$ is quadratic of type (iv) in this case.
\end{enumerate}
\end{bem}

\begin{lemma}
An element $a\in L_0^*$ is inseparable iff we have $(a,t)^2=0_T$ for each $(a,t)\in T$.
\end{lemma}

\begin{bew} Notice that we have $\Char \AA=2$ by remark \ref{195}. Given $(a,t)\in T$, we have
$$0_T=(a,t)^2=\big(a+a,t+t+f(a,a)\big)=\big(0_{L_0},f(a,a)\big)\ \Leftrightarrow\ f(a,a)=0_\AA\ .$$\qed
\end{bew}

\begin{kor}\label{197}
An element $a\in L_0^*$ is inseparable if and only if $\p_1(a)\in \tilde{L}_0^*$ is inseparable.
\end{kor}

\begin{bew}
Given $(a,t)\in T$, we have
$$(a,t)^2=0_T\ \Leftrightarrow\ \gamma(a,t)^2=0_{\tilde{T}}\ .$$
\qed
\end{bew}

\begin{no}\ \newglossaryentry{QSx}{type=symbols,name={\ensuremath{\EE_x}},description=quadratic subfield containing $x$, sort=pseudoquadratic space}
\newglossaryentry{QSa}{type=symbols,name={\ensuremath{\EE_a}},description=quadratic subfield containing $q(a)$, sort=pseudoquadratic space}
\newglossaryentry{Sa}{type=symbols,name={\ensuremath{R_a}},description=$2$-dimensional $\FF$-subspace with respect to $a$, sort=pseudoquadratic space}
\newglossaryentry{Xa}{type=symbols,name={\ensuremath{X_a}},description=admissible elements with respect to $a$, sort=pseudoquadratic space}
\begin{itemize}
\item We set
\begin{equation} g:L_0\times L_0\to \FF,\  (a,b)\mapsto f(b,a) - q(a+b) + q(a) + q(b)\ .\label{eq:defg}\end{equation}
 \item Given $x\in \AA\sm \FF$, let $\gls{QSx}$ be the quadratic subfield of $\AA$ generated by $1_\AA$ and $x$.
	\item Given $a \in L_0^*$, we set $\gls{QSa}:=\EE_{q(a)}$. It is a separable quadratic subfield iff $a$ is separable.
	\item 	Given $a \in L_0^*$, we set
	\[ \gls{Sa} := \langle a, a\cdot q(a) \rangle_\FF=a\cdot \EE_a \ , \]
	so that $R_a$ is a $2$-dimensional $\FF$-subspace of $L_0$.
	\item Suppose $(\AA,\FF,\sigma)$ to be quadratic of type (iv). Given a separable element $a\in L_0$, let $e_a$ be an element of $\AA$ orthogonal to $\EE_a$ (with respect to the standard trace of $\AA$).
	For each $x \in \AA$, let $\alpha_a(x)$ and $\beta_a(x)$ be the unique elements of $\EE_a$ such that
	$$x = \alpha_a(x) + e_a \beta_a(x)\ .\footnote{We assume $e_{a'} = e_{a}$ for each $a'\in R_a$, hence $\alpha_{a'}=\alpha_a$ and $\beta_{a'}=\beta_a$ for each $a' \in R_{a}$}$$
	\item Given $a\in L_0^*$, we set
	$$\gls{Xa}:=\{ t\in \KK \mid (a,t)\in T \}=\{ s+q(a) \mid s\in \FF \}\subseteq \EE_a\ .$$
		\end{itemize}
\end{no}

\begin{bem}\label{235}
Given $a\in L_0^*$, we have $X_a\subseteq \EE_a\sm \FF$ which implies
$$\forall\ t\in X_a:\qquad \EE_t=\EE_a\ .$$
In particular, we have $\EE_a=\EE_t$ for each $(a,t)\in T$ such that $a\neq 0_{L_0}$.
\end{bem}

\begin{lemma}\label{306}
Let $a\in L_0^*$ be such that $N(t)=T(t)^2$ for each $t\in X_a$. Then we have $|\FF|=2$.
\end{lemma}

\begin{bew}
Let $t\in X_a$. Given $s\in \FF$, we have
\begin{align*}
0_\AA&=T(t)^2-N(t)=T(t+s)^2-N(t)=N(t+s)-N(t)=sT(t)+N(s)=s\big(T(t)+s\big)
\end{align*}
and therefore
$$\forall\ s\in \FF^*:\qquad s=-T(t)\ ,$$
which implies $|\FF|=2$.\qed
\end{bew}

\begin{lemma}\label{210}
Let $a\in L_0^*$ and $|\FF|\geq 3$. Then we have
$$\{ t^\sigma-t^{-1}f(a,a)t^\sigma \mid t\in X_a\}\not\subseteq \FF \ .$$
\end{lemma}

\begin{bew}
Notice that we have
\begin{align}\{ t^\sigma-t^{-1}f(a,a)t^\sigma \mid t\in X_a\}=\{ \big(1_\AA-t^{-1}(t-t^\sigma)\big)t^\sigma \mid t\in X_a \}=\{ t^{-1}(t^\sigma)^2 \mid t\in X_a\}\ .\label{272}\end{align}
Assume $\{ t^{-1}(t^\sigma)^2 \mid t\in X_a\}\subseteq \FF$. Then we have
\begin{align*}
\forall\ t\in X_a:&& (t^\sigma)^3=t^{-1}(t^\sigma)^2\cdot N(t)\in \FF\ , && t^3\in \FF\ .
\end{align*}
Therefore, we have
\begin{align}
\big(T(t)^2-N(t)\big)\cdot t-T(t)N(t)=\big(T(t)t-N(t)\big)\cdot t=t^3\in \FF \label{222}
\end{align}
and thus $T(t)^2=N(t)$ for each $t\in X_a$. Now lemma \ref{306} yields $|\FF|=2$.\qed
\end{bew}

\begin{bem}\ 
\begin{enumerate}[label=(\alph*)]
\item \label{205}
The map $$\phi:\FF\to\tilde{\FF},\ t\mapsto \p_2(0,t)$$ is a Jordan isomorphism by lemma \ref{204}, hence an isomorphism of fields by Hua's theorem, cf. theorem \ref{386}. Therefore, the map $\p_1:L_0\to \tilde{L}_0$ is an isomorphism of vector spaces over $\FF$ by lemma \ref{172}.
\item \label{211} Given $(a,t)\in T$, we have $\big(\p_1(a),\p_2(a,t)\big)\in \tilde{T}$, hence $\tilde{\EE}_{\p_2(a,t)}= \tilde{\EE}_{\p_1(a)}$ by remark \ref{235}.
\end{enumerate}
\end{bem}

\begin{kor}\label{212}
Given $a\in L_0^*$, we have 
$$\p_1(R_a)=\tilde{R}_{\p_1(a)}\ .$$
\end{kor}

\begin{bewzwei}\ 
\begin{itemize}
\item Assume $|\FF|\geq 3$. By lemma \ref{210}, we have
$$\langle 1_\AA, \{ t^\sigma-t^{-1}f(a,a)t^\sigma \mid t\in X_a\} \rangle_\FF=\EE_{a}\ ,$$
hence
$$\p_1(R_a)=\p_1(a\cdot  \EE_{a})\subseteq \p_1(a)\cdot \tilde{\EE}_{\p_1(a)}=\tilde{R}_{\p_1(a)}$$
by remark \ref{205}, identity \eqref{169} with $(b,v):=(a,t)$, corollary \ref{187} and remark \ref{211}.
\item It remains to consider the case $|\FF|=2$ which implies $\AA=\EE_a\cong\FF_4\cong\tilde{\EE}_{\p_1(a)}=\tilde{\AA}$. Given $a\in L_0^*$ and $t\in X_a$, we have
\begin{align*}
N(t)=1_\AA=T(t)=f(a,a)\ , 
\end{align*}
hence
$$(at,t)=(at, t^\sigma tt)\in T\ .$$
Substituting $(a,t)$ by $(at,t)$ and $(b,v)$ by $(a,t)$ in identity \eqref{169}  yields
\begin{align*}
\p_1(a)&=\p_1\big(a\cdot (t^\sigma+t)\big)=\p_1\big(a\cdot t^\sigma+(at)\cdot t^{-1}f(at,a)t^\sigma\big)\\
&=\p_1(a)\cdot \p_2(at,t)^{\tilde{\sigma}}+\p_1(at)\cdot \p_2(at,t)^{-1}\tilde{f}\big(\p_1(at),\p_1(a)\big)\p_2(at,t)^{\tilde{\sigma}}\ .
\end{align*}
Because of $\p_2(at,t)\notin \tilde{\FF}$ we have
$$\p_2(at,t)^{-1}\tilde{f}\big(\p_1(at),\p_1(a)\big)\p_2(at,t)^{\tilde{\sigma}}\neq 0_{\tilde{\AA}}$$
and thus
\begin{align*} \p_1(at)=\big(\p_1(a)+\p_1(a)\cdot \p_2(at,t)^\sigma\big)\cdot \big( \p_2(at,t)^{-1}\tilde{f}(\p_1(at),\p_1(a))\p_2(at,t)^{\tilde{\sigma}}\big)^{-1}\in \tilde{R}_{\p_1(a)}\ .
\end{align*}
\end{itemize}
\qed
\end{bewzwei}

\begin{bem}
The following lemma is due to Tom De Medts.
\end{bem}
\begin{lemma}\label{le:one}
	Let $x \in \AA$ with $x^\sigma \neq x$.
	Then the set
	\[ S := \left\{ (s + x)^{-1} (s + x)^\sigma \mid s \in \FF \right\} \]
	cannot be completely contained in a one-dimensional $\FF$-subspace of $\AA$.
\end{lemma}

\begin{bew}
  Suppose that $|\FF|=2$. Then we have $\AA\cong \FF_4$, hence $$S=\{ x, x+1_\AA \}\ ,$$ and the assertion is true. So assume $|\FF|\geq 3$ and $S \subseteq y \cdot \FF$ for some $y \in \AA^*$.
	Then for each $s \in \FF$, there is an element $t_s \in \FF$ such that $(s+x)^\sigma = (s+x) y t_s$.
	In particular, $x^\sigma = x y t_0$, and hence we also have $(s+x)^\sigma = s + x^\sigma = s + x y t_0$.
	Therefore, $(s+x) y t_s = s + x y t_0$ for all $s \in \FF$.
	Multiplying on the right by $s^{-1} y^{-1}$ yields
	\begin{equation}\label{eq:one}
		s^{-1} \bigl( (s+x) t_s - x t_0 \bigr) = y^{-1}
	\end{equation}
	for all $s \in \FF^*$.
	It follows that
	\[ r^{-1} \bigl( (r+x) t_r - x t_0 \bigr) = s^{-1} \bigl( (s+x) t_s - x t_0 \bigr) \]
	for all $r,s \in \FF^*$.
	This can be rearranged to get
	\[ rs ( t_s - t_r ) = x ( s t_r - s t_0 - r t_s + r t_0 ) \]
	for all $r,s \in \FF^*$.
	Since $x \not\in \FF$, this can only happen if both sides are zero.
	Hence $t_s = t_r$ for all $r,s \in \FF^*$, and substituting this in the right hand side gives
	$(s-r)(t_r - t_0) = 0_\AA$ for all $r,s\in \FF^*$ and hence $t_r = t_0$ for all $r \in \FF$, where we use the fact that $\FF^*$ has at least two elements. Substituting this in equation \eqref{eq:one} yields $t_0 = y^{-1}$,
	but then $x^\sigma = x y t_0 = x$.
	This contradiction finishes the proof.\qed
\end{bew}

\begin{bem}\label{300}
Let $d \in L_0^*$ be separable. Then for all $s \in \FF$ and $x \in \AA$, we have
	\[ \bigl( s + q(d) \bigr)^{-1} x \bigl( s + q(d) \bigr)^\sigma
		= \alpha_d(x) \bigl( s + q(d) \bigr)^{-1} \bigl( s + q(d) \bigr)^\sigma + e_d \beta_d(x) \ , \]
	and by lemma \ref{le:one}, the set
	$$S := \left\{ \bigl( s + q(d) \bigr)^{-1} \bigl( s + q(d) \bigr)^\sigma \mid s \in \FF \right\}$$
	cannot be completely contained in a one-dimensional $\FF$-subspace of $\AA$.
\end{bem}

\begin{lemma}\label{196}
Let $\AA$ be a quaternion division algebra and let $a,b\in L_0$ be separable. Then there are $a'\in R_{a}$ and $b'\in R_{b}$ such that $a'+b'$ and $a'\cdot q(a')+b'$ are separable.
\end{lemma}

\begin{bewzwei}\ 
\begin{itemize}
\item Let $\Char \AA\neq 2$. By remark \ref{195}, it is enough to choose $a':=a$ and $s\in \FF$ such that $$b':=b\cdot s\notin \{-a,-a\cdot q(a)\}\ .$$
\item Let $\Char \AA=2$. Given $x\in L$, we have
$$f(x,x)=q(x)+q(x)^\sigma\in \FF$$
and thus
\begin{align*}
F_x(s):=f(x+b\cdot s,x+b\cdot s)&=s^2\cdot f(b,b)+s\cdot \big(f(x,b)+f(b,x)\big)+f(x,x) \\
&=s^2\cdot f(b,b)+s\cdot \big(f(x,b)+f(x,b)^\sigma\big)+f(x,x)\in \FF[s]\ .
\end{align*}
Since $\AA$ is non-commutative, we have $|\FF|=\infty$. As a consequence, there is an element $s\in \FF$ not contained in the set of zeroes of $F_a(s)$ and $F_{a\cdot q(a)}(s)$. Then $a':=a$ and $b'=b\cdot s$ satisfy the required conditions.
\end{itemize}\qed
\end{bewzwei}

\begin{no}
	Given $a,b\in L_0^*, c\in L_0$ and $r,s,t\in \FF$, we set
	\[ M_{(a,s),(b,t),r}(c) := h_{(a,s+q(a)) \cdot (b,t+q(b)) \cdot (0,r)}(c) - h_{(a,s+q(a))}(c) - h_{(b,t+q(b))}(c)\ . \]
\end{no}

\begin{bem} Given $a,b,c\in L_0^*$, $r,s,t\in \FF$ and $r':=r+s+t+g(a,b)$, we have
$$\big(a,s+q(a)\big)\cdot (b,t+q(b)\big)\cdot (0,r)=\big(a+b,r'+q(a+b)\big)\ .$$
Now it follows from equation \eqref{eq:defg} and definition \ref{eq:hua} that we have
	\begin{multline*}
		M_{(a,s),(b,t),r}(c)  = c \cdot \bigl( r - f(a,b) \bigr) \\
		+ a \cdot \Bigl( \bigl( s + q(a) \bigr)^{-1} f(a,c) \bigl( s + q(a) \bigr)^\sigma \\
			\shoveright{- \bigl( r' + q(a+b) \bigr)^{-1} f(a+b,c) \bigl( r' + q(a+b) \bigr)^\sigma \Bigr)} \\
		+ b \cdot \Bigl( \bigl( t + q(b) \bigr)^{-1} f(b,c) \bigl( t + q(b) \bigr)^\sigma \\
			- \bigl( r' + q(a+b) \bigr)^{-1} f(a+b,c) \bigl( r' + q(a+b) \bigr)^\sigma \Bigr)\ .
	\end{multline*}
\end{bem}

\begin{bem}
The following lemma shows that orthogonality can be encoded in a condition that is preserved by Jordan isomorphisms. Notice that we need three pairwise orthogonal vectors. This is the point where the assumption about the dimension will come in.
\end{bem}

\begin{lemma}
	\label{203} Let $x,y,z \in L_0^*$. Suppose that $f(x,y) = f(x,z) = f(y,z) = 0_\AA$. Then we have
	\begin{equation}\label{eq:mucond}
		\forall\ r, s, t \in \FF,\ a \in R_{a'}, b \in R_{b'}, c \in R_{c'}:\qquad	M_{(a,s),(b,t),r}(c) = c \cdot r 
	\end{equation}
	for each permutation $(a', b', c')$ of $(x,y,z)$.
	\end{lemma}

\begin{bew}	
	For each permutation $(a', b', c')$ of $(x,y,z)$ and for all
	$a \in R_{a'}$, $b \in R_{b'}$ and $c \in R_{c'}$, we have
	$$f(a,b) = f(a,c) = f(b,c) = 0_\AA$$ as well,
	and hence $M_{(a,s),(b,t),r}(c) = c \cdot r$ for all $r,s,t \in \FF$, which shows \eqref{eq:mucond}.\qed
\end{bew}

\begin{bem}\label{213}
Let $a\in L_0^*$. Then we have $$\forall\ t\in X_a:\qquad \p_2(a,t)\in \tilde{q}\big(\p_1(a)\big)+\tilde{\FF}=\tilde{X}_{\p_1(a)}\ .$$
Since $\gamma$ is surjective, it follows that
$$\{ \p_2(a,t) \mid t\in X_a\}=\tilde{X}_{\p_1(a)}\ .$$
\end{bem}

\begin{kor}\label{215}
Let $x,y,z \in L_0^*$. Suppose that $$f(x,y) = f(x,z) = f(y,z) = 0_\AA\ .$$ Then \eqref{eq:mucond} holds for each permutation $(a', b', c')$ of $\big(\p_1(x),\p_1(y),\p_1(z)\big)$.
\end{kor}

\begin{bew}
By lemma \ref{203}, \eqref{eq:mucond} holds for each permutation $(a', b', c')$ of $(x,y,z)$. Since $\gamma$ preserves the Hua-maps and $\p_1:L_0\to \tilde{L}_0$ is an isomorphism of vector spaces over $\FF$ by remark \ref{205}, we have
\begin{align*}
	\tilde{M}_{(\p_1(a),\p_2(a,s)),(\p_1(b),\p_2(b,t)),\phi(r)}(\p_1(c))=\p_1\big(M_{(a,s),(b,t),r}(c)\big) = \p_1(c \cdot r)=\p_1(c)\cdot \phi(r)
\end{align*}
for all $a\in R_{a'}, b\in R_{b'}, c\in R_{c'}$, $r, s, t \in \FF$ and for each permutation $(a',b',c')$ of $(x,y,z)$. Corollary \ref{212} and remark \ref{213} yield
\begin{align*}
	\forall\ a\in  \tilde{R}_{a'}, b\in \tilde{R}_{b'}, c\in \tilde{R}_{c'},\ r,s,t\in \tilde{\FF}:\qquad \tilde{M}_{(a,s),(b,t),r}(c)=c\cdot r
\end{align*}
for each permutation $(a',b',c')$ of $\big(\p_1(x),\p_1(y),\p_1(z)\big)$.\qed
\end{bew}

\begin{bem}
The following lemma is essentially due to Tom De Medts. It shows that we can reconstruct the orthogonality from the above condition if we suppose all appearing elements to be separable. In this situation, we know that we have a one-dimensional $\FF$-subspace on the right side of \eqref{eq:mucond}, but on the left side we have terms which are not contained in a one-dimensional $\FF$-subspace if we vary the coefficients. As a consequence, the occurring scalar products necessarily vanish.

Afterwards we will have a closer look at inseparable elements. In order to obtain the same result in this case, we need to establish a connection between the corresponding quadratic extensions. For this purpose, it is convenient to use an inseparable element in \eqref{eq:mucond} twice so that we can deduce more information about the occurring terms.
\end{bem}

\begin{lemma} \label{202} \label{206} Let $x,y,z \in L_0^*$. Suppose that \eqref{eq:mucond} holds for each permutation $(a', b', c')$ of $(x,y,z)$. If $x,y,z$ are separable, we have
	$$f(x,y)=f(x,z)=f(y,z)=0_\AA\ .$$
\end{lemma}

\begin{bew}
	Assume that \eqref{eq:mucond} holds for each permutation $(a', b', c')$ of $(x,y,z)$.
	In particular, for all $a \in R_{a'}$, $b \in R_{b'}$ and $c \in R_{c'}$,
	the set
	\[ \left\{ M_{(a,s),(b,t),r}(c) \mid r,s,t \in \FF \right\} \]
	is a one-dimensional $\FF$-subspace of $L_0$.
	\begin{itemize}
	\item Suppose that $\AA$ is a quaternion division algebra. By lemma \ref{196}, there are $a\in R_{a'},b\in R_{b'}$ such that $a+b$ and $a\cdot q(a)+b$ are   separable. By varying $s$ over $\FF$ in $M_{(a,s),(b,t),r}(c)$ but keeping $r':=r+s+t+g(a,b)$ and $t$ invariant (by the right choice for $r$), remark \ref{300} yields\footnote{Notice that we need $a,b$ and $a+b$ to be separable so that we can apply Lemma \ref{le:one} at this point. This follows from the assumption that $x,y,z$ are separable and from the choice of $a,b$. As we will replace $a$ by $a\cdot q(a)$, we additionally need $a\cdot q(a)+b$ to be separable. }  $\alpha_a\big(f(a,c)\big) = 0_\AA$.
	Similarly, $$ \alpha_b\big(f(b,c)\big) = 0_\AA=\alpha_{a+b}\big(f(a+b,c)\big)\ .$$
	Now the expression for $M_{(a,s),(b,t),r}(c)$ can be simplified\footnote{Notice that we have $e_a\beta_a\big(f(a,c)\big)=\alpha_a\big(f(a,c)\big)+e_a\beta_a\big(f(a,c)\big)=f(a,c)$ etc.} to
	\[ M_{(a,s),(b,t),r}(c) = c \cdot \big( r - f(a,b) \big) - a \cdot f(b,c) - b \cdot f(a,c) \]
	for all $c\in R_{c'}$, $r,s,t \in \FF$. By assumption, we have
	\begin{equation}\label{eq:hh}
		c \cdot f(a,b) + a \cdot f(b,c) + b \cdot f(a,c) = 0_\AA
	\end{equation}
for each $c \in R_{c'}$.

	Now suppose that $f(a,b) \neq 0_\AA$; we will derive a contradiction.
	If we interchange $a$ and $b$ in equation \eqref{eq:hh}, we get $f(a,b) = f(b,a)$, and if we replace $a$
	by $a\cdot q(a) \in R_{a'}$, then $f\big(a\cdot q(a), b\big) = f\big(b, a\cdot q(a)\big)$ and hence
	\begin{equation}\label{eq:pi}
		q(a)^\sigma f(a,b) = f(a,b) q(a)\ .
	\end{equation}
	Now let $t$ be an arbitrary element of $\EE_c$.
	On the one hand, we can multiply equation \eqref{eq:hh} by $t$; on the other hand, we can
	replace $c$ by $ct \in R_c$.
	Comparing these two resulting equations yields $t f(a,b) = f(a,b) t$ for all $t \in \EE_c$, which implies
	$f(a,b) \in C_\AA(\EE_c) = \EE_c$. 	If we now replace $a$ by $a\cdot q(a)\in R_a=a\cdot \EE_a$, we get $\EE_a f(a,b) \subseteq \EE_c$,
	which implies $\EE_a = \EE_c$ and hence $f(a,b) \in \EE_a$.
	But then $q(a) f(a,b) = f(a,b) q(a)$, and comparing this with equation \eqref{eq:pi} yields
	$q(a) = q(a)^\sigma$, which gives us the required contradiction and hence
	$$f(a',b')=f(a,b)=0_\AA\ .$$

	Permuting $a'$, $b'$ and $c'$ now yields $$f(x,y) = f(x,z) = f(y,z) = 0_\AA\ .$$
	\item If $\AA$ is commutative, we immediately obtain
	$$f(a,c)=0_\AA=f(b,c)$$
	by the same arguments, followed by
	$$f(a+b,c)=f(a,c)+f(b,c)=0_\AA$$
	and finally $f(a,b)=0_\AA$.
	\end{itemize}  \qed
\end{bew}

\begin{bem}
Lemma \ref{202} is enough to handle the cases where $\Char \AA\neq2$ or where $(\AA,\FF,\sigma)$ is quadratic of type (iii) since there are no inseparable elements in this situation.  But with some technical effort, we can handle the remaining case as well.
\end{bem}

\begin{no} 
Until proposition \ref{198}, we suppose $(\AA,\FF,\sigma)$ to be quadratic of type (iv) with $\Char \AA=2$.
\end{no}

\begin{bem}\label{220}
Let $x\in L_0$ be inseparable and let $y\in L_0$ be such that $f(x,y)=f(x,y)^\sigma$. Then we have
$$f(x+y,x+y)=0_\AA\ \Leftrightarrow\ f(x,x)+f(x,y)+f(y,x)+f(y,y)=0_\AA\ \Leftrightarrow\ f(y,y)=0_\AA\ .$$
\end{bem}

\begin{lemma}\label{214}
Let $x,y,z\in L_0^*$ be such that \eqref{eq:mucond} holds for each permutation $(a',b',c')$ of $(x,y,z)$. If $x=y$ is inseparable and $f(x,z)\neq 0_\AA$, we have $f(x,z)\in\EE_x=\EE_{x+z}$, and $z$ is inseparable.
\end{lemma}

\begin{bew}
Since $x$ is inseparable and $f(x,z)\neq 0_\AA$, we have $z\notin \langle x\rangle_\AA$. Putting $a'=c'=x,\ b'=z$ and $a=c$ in \eqref{eq:mucond} and comparing the coefficients of $a=c$ yield
$$f(a,b)=\big(r'+q(a+b)\big)^{-1} f(b,a)\big(r'+q(a+b)\big)^\sigma$$
for all $a\in R_x, b\in R_z$, and putting $a'=z,\ b'=c'=x$ and $b=c$ in \eqref{eq:mucond} and comparing the coefficients of $b=c$ yield
\begin{align} f(a,b)=\big(r'+q(a+b)\big)^{-1} f(a,b)\big(r'+q(a+b)\big)^\sigma \label{301} \end{align}
for all $a\in R_z, b\in R_x$, hence
$$\forall\ a\in R_x, b\in R_z:\qquad f(a,b)=f(b,a)=f(a,b)^\sigma\ .$$
In particular, we have
\begin{align*} f(x,z)=f(z,x)\ , && q(x)f(x,z)=f(x,z)q(x)\ ,&& f(x,z)\in C_\AA\big(q(x)\big)=\EE_x\end{align*}
Putting $a'=z,\ b'=c'=x,\ a=z,\ b=x$ and $c=x\cdot q(x)$ in \eqref{eq:mucond} and comparing the coefficients of $x$ yield
\begin{align*}
q(x)f(z,x)=q(x+z)^{-1}f(z,x)q(x)q(x+z)^\sigma=q(x+z)^{-1}q(x)f(z,x)q(x+z)^\sigma\ .
\end{align*}
Using equation \eqref{301} with $a=z$ and $b=c=x$ and multiplying by $q(x)$ yield
\begin{align*}
q(x)f(z,x)=q(x)q(x+z)^{-1}f(z,x)q(x+z)^\sigma\ .
\end{align*}
As a consequence, we have
\begin{align*} q(x)\in C_\AA\big(q(x+z)\big)=\EE_{x+z}\ ,&& \EE_x=\EE_{x+z}\ .\end{align*}
Since $x$ is inseparable, $x+z$ has to be inseparable as well. Finally, $z$ is inseparable by remark \ref{220}.
\qed
\end{bew}

\begin{kor}\label{216}
Let $x\in L_0^*$ be inseparable and let $y\in x^\bot$ be separable. Then we have $$\p_1(x)\in \p_1(y)^\bot\ .$$
\end{kor}

\begin{bew}
By corollary \ref{215}, identity \eqref{eq:mucond} holds for each permutation $(a',b',c')$ of $\big(\p_1(x),\p_1(x),\p_1(y)\big)$. By corollary \ref{197}, $\p_1(x)$ is inseparable and $\p_1(y)$ is separable, thus lemma \ref{214} yields
$$\big(\p_1(x),\p_1(y)\big)=0_{\tilde{\AA}}\ .$$\qed
\end{bew}

\begin{lemma}\label{217}
Let $x\in L_0^*$ and $y\in x^\bot\sm\langle x\rangle_\AA$ both be inseparable. Then we have $$\p_1(x)\in \p_1(y)^\bot\ .$$
\end{lemma}

\begin{bew}
By corollary \ref{215}, identity \eqref{eq:mucond} holds for each permutation $(a',b',c')$ of $\big(\p_1(x),\p_1(x),\p_1(y)\big)$ and $\big(\p_1(x),\p_1(y),\p_1(y)\big)$, respectively. Suppose that $\tilde{f}\big(\p_1(x),\p_1(y)\big)\neq 0_{\tilde{\AA}}$. By corollary \ref{197}, $\p_1(x)$ and $\p_1(y)$ are inseparable, hence lemma \ref{214} yields
$$\tilde{f}\big(\p_1(x),\p_1(y)\big)\in \tilde{\EE}_{\p_1(x)}=\tilde{\EE}_{\p_1(x+y)}=\tilde{\EE}_{\p_1(y)}\ .$$
By identity \eqref{169} with $b:=y$ and $(a,t):=\big(x,q(x)\big)$, we have
$$\p_1\big(y\cdot q(x)^\sigma\big)=\p_1(y)\p_2\big(x,q(x)\big)^{\tilde{\sigma}}-\p_1(x)\cdot \p_2\big(x,q(x)\big)^{-1}\tilde{f}\big(\p_1(x),\p_1(y)\big)\p_2\big(x,q(x)\big)^{\tilde{\sigma}} $$
with
\begin{align*} \p_2\big(x,q(x)\big)^{\tilde{\sigma}}\in \tilde{\EE}_{\p_1(x)}=\tilde{\EE}_{\p_1(y)}\ , &&  \p_2\big(x,q(x)\big)^{-1}\tilde{f}\big(\p_1(x),\p_1(y)\big)\p_2\big(x,q(x)\big)^{\tilde{\sigma}}\in \tilde{\EE}_{\p_1(x)}\ .\end{align*}
By corollary \ref{212}, there are elements $s,t\in \AA$ such that
\begin{align*} \p_1(y\cdot s)=\p_1(y)\p_2\big(x,q(x)\big)^{\tilde{\sigma}}\ , && \p_1(x\cdot t)=\p_1(x)\cdot \p_2\big(x,q(x)\big)^{-1}\tilde{f}\big(\p_1(x),\p_1(y)\big)\p_2\big(x,q(x)\big)^{\tilde{\sigma}}\ ,
\end{align*}
which yields
\begin{align*} y\cdot q(x)^\sigma+y\cdot s+ x\cdot t=0_{L_0}\ , && s=q(x)^\sigma,\ t=0_\AA
\end{align*}
and finally
\begin{align*} 0_{L_0}=\p_1(x\cdot t)=\p_1(x)\cdot \p_2\big(x,q(x)\big)^{-1}\tilde{f}\big(\p_1(x),\p_1(y)\big)\p_2\big(x,q(x)\big)^{\tilde{\sigma}}\ , && \tilde{f}\big(\p_1(x),\p_1(y)\big)=0_{\tilde{\AA}}\ .
\end{align*}\qed
\end{bew}

\begin{prop}\label{198}
Let $x\in L_0^*$ be inseparable and let $y\in x^\bot$. If $\dim_\AA L_0\geq 3$, we have $$\p_1(x)\in \p_1(y)^\bot\ .$$
\end{prop}

\begin{bew}
If $y$ is separable, we may apply lemma \ref{216}, and if $x$ and $y$ are linearly independent over $\AA$, we may apply lemma \ref{217} so that we may assume $y\in \langle x\rangle_\AA$. By assumption, there is an element $z\in x^\bot\sm \langle x\rangle_\AA$. As a consequence, we have
$$y+z\in x^{\bot}\sm \langle x\rangle_\AA\ .$$
Now corollary \ref{216} and lemma \ref{217}, respectively,  yield
$$\tilde{f}\big(\p_1(x),\p_1(y)\big)=\tilde{f}\big(\p_1(x),\p_1(y+z)\big)+\tilde{f}\big(\p_1(x),\p_1(z)\big)=0_{\tilde{\AA}}+0_{\tilde{\AA}}=0_{\tilde{\AA}}\ .$$\qed
\end{bew}

\begin{lemma}\label{199}
	Assume $\dim_\AA L_0 \geq 3$ and let $a,b \in L_0^*$.
	Then $f(a,b) = 0_\AA$ if and only if there is an element $c \in L_0^*$ such that
	$$f(a,b) = f(a,c) = f(b,c) = 0_\AA\ .$$
\end{lemma}

\begin{bew}
If $a$ or $b$ is inseparable, we may choose $c:=a$ or $c:=b$, respectively. So assume that $a$ and $b$ are separable and $f(a,b) = 0_\AA$.
	Since $\dim_\AA L_0 \geq 3$ and $\dim_\AA \langle a,b \rangle_\AA \leq 2$, there is some element $d \in X \setminus \langle a,b \rangle_\AA$.
	Then
	\[ c := d - a \cdot f(a,a)^{-1} f(a,d) - b \cdot f(b,b)^{-1} f(b,d) \neq 0_\AA \]
	satisfies the required conditions.\qed
\end{bew}

\begin{bem}
We return to the general case.
\end{bem}

\begin{prop}\label{200} Let $x\in L_0$ be separable and let $y\in x^\bot$. If $\dim_\AA L_0\geq 3$, we have $$\p_1(x)\in \p_1(y)^\bot\ .$$
\end{prop}

\begin{bewzwei}\ 
\begin{itemize}
\item If $y$ is inseparable, we may apply proposition  \ref{198}.
\item If $y$ is separable, lemma \ref{199} yields an element $z\in L_0^*$ such that
$$f(x,z)=0_\AA=f(y,z)\ ,$$
thus \eqref{eq:mucond} holds for each permutation $(a',b',c')$ of $\big(\p_1(x),\p_1(y),\p_1(z)\big)$ by corollary \ref{215}. 
\begin{enumerate}[label=$\circ$]
\item If $z$ is inseparable, proposition \ref{198} yields
$$\tilde{f}\big(\p_1(x),\p_1(z)\big)=0_{\tilde{\AA}}=\tilde{f}\big(\p_1(y),\p_1(z)\big)\ ,$$
followed by
$$\tilde{f}\big(\p_1(x+y),\p_1(z)\big)=0_{\tilde{\AA}}$$
and finally $\tilde{f}\big(\p_1(x),\p_1(y)\big)=0_{\tilde{\AA}}$ by \eqref{eq:mucond}.
\item If $z$ is separable, then $\p_1(x),\p_1(y),\p_1(z)$ are separable by corollary \ref{197}, thus we may apply lemma \ref{202} to obtain
$$\tilde{f}\big(\p_1(x),\p_1(y)\big)=0_{\tilde{\AA}}\ .$$
\end{enumerate}
\end{itemize}\qed
\end{bewzwei}

\begin{kor}\label{207}
If we have $\dim_\AA L_0\geq 3$, the map $\p_1:L_0\to \tilde{L}_0$ is an isomorphism of vector spaces.
\end{kor}

\begin{bew}
By proposition \ref{198} and proposition \ref{200}, we have 
$$\forall\ a\in L_0:\qquad \p_1(a^\bot)=\p_1(a)^\bot\ .$$
Since $f$ is non-degenerate, this implies that we have
$$\forall\ a\in L_0:\qquad \p_1(\langle a\rangle_\AA)=\p_1(a^{\bot\bot})=\p_1(a)^{\bot\bot}=\langle \p_1(a)\rangle_{\tilde{\AA}}\ .$$
Now the assertion results from the fundamental theorem of projective geometry.\qed
\end{bew}

\begin{no}\label{208}
Let $\phi:\AA\to\tilde{\AA}$ be the isomorphism of skew-fields defined by
$$\forall\ a\in L_0,\ t\in \AA:\qquad \p_1(a\cdot t)=\p_1(a)\cdot \phi(t)\ .$$
\end{no}

\begin{bem}\label{232}
Notice that $\phi$ is an extension of the isomorphism $\phi:\FF\to \tilde{\FF}$ of fields.
\end{bem}

\begin{bem}
We state the theorem using the general notation.
\end{bem}

\begin{satz}\label{178}
Let $(\KK,\KK_0,\sigma)$ be quadratic of type (iii) or (iv) and let $\dim_\KK L_0\geq 3$. Then the map $\Phi:\Xi\to\tilde{\Xi}$ defined by
\begin{align*}
\Phi:=(\p_1,\phi): (L_0,\KK)\to(\tilde{L}_0,\tilde{\KK}),\ (a,t)\mapsto \big(\p_1(a), \phi(t)\big)
\end{align*}
is an isomorphism of pseudo-quadratic spaces satisfying
$$\forall\ (a,t)\in T:\qquad \Phi(a,t)=\gamma(a,t)\ .$$
\end{satz}

\begin{bewzwei}\ 
\begin{itemize}
\item By corollary \ref{207} and notation \ref{208}, the map $(\p_1,\phi):(L_0,\KK)\to (\tilde{L}_0,\tilde{\KK})$ is an isomorphism of vector spaces. 
\item By remark \ref{232}, we have $\phi(\KK_0)=\tilde{\KK}_0$.
\item By corollary \ref{251}, we have
$$\forall\ x\in \KK:\qquad \phi(x^\sigma)=\phi(x)^{\tilde{\sigma}}\ .$$
\item Let $(a,t)\in T$ and $0_{L_0}\neq b\in a^\bot$. Then we have
$$\p_1(b)\cdot \phi(t)^{\tilde{\sigma}}=\p_1(b)\cdot \phi(t^\sigma) =\p_1(b\cdot t^\sigma)=\p_1(b)\cdot \p_2(a,t)^{\tilde{\sigma}}$$
by identity \eqref{169}, proposition \ref{198} and proposition \ref{200}, thus
$$\p_2(a,t)=\phi(t)$$
and therefore
\begin{align*}
\tilde{q}\big(\p_1(a)\big)&\in \p_2(a,t)+\tilde{\KK}_0=\phi(t)+\phi(\KK_0)= \phi\big(t+\KK_0\big)=\phi\big(q(a)+\KK_0\big)=\phi\big(q(a)\big)+\tilde{\KK}_0
\end{align*}
as well as
$$\Phi(a,t)=\big(\p_1(a),\phi(t)\big)=\big(\p_1(a),\p_2(a,t)\big)=\gamma(a,t)\ .$$
\end{itemize}\qed
\end{bewzwei}

\chapter{Small Dimensions I}\label{314}
First we refine some results of the previous chapter, without assuming additional assumptions. The first step is to show that the map $\p_1: R_a\to R_{\p_1(a)}$ is an isomorphism of vector spaces for each $a\in L_0$. We manage to do this by proving that the map $\p_2:X_a\to X_{\p_1(a)},\ t\mapsto \p_2(a,t)$ is induced by an isomorphism between the associated separable extensions.

Once we have done this, it is easy to prove theorem \ref{254} if the involutory sets are quadratic of type (iii). Afterwards we will need some more considerations to handle the case of involutory sets which are quadratic of  type (iv), cf. chapter \ref{268}.

\begin{no}
Throughout this chapter, the involutory sets $(\AA,\FF,\sigma)$ and $(\tilde{\AA},\tilde{\FF},\tilde{\sigma})$ are still quadratic of type (iii) or (iv).
\end{no}

\begin{lemma}\label{224}
Let $a\in L_0$. Then the following holds:
\begin{enumerate}[label=(\alph*)]
\item \label{273} $\forall\ t\in X_a:\ \phi\big(N(t)\big)=\tilde{N}\big(\p_2(a,t)\big)$.
\item $\forall\ t\in X_a:\ \phi\big(T(t)\big)=\tilde{T}\big(\p_2(a,t)\big)$.
\end{enumerate}
\end{lemma}

\begin{bewzwei}\ 
\begin{enumerate}[label=(\alph*)]
\item By lemma \ref{171} with $(b,v):=(0_{L_0},1_\AA)$, we have
$$\phi(tt^\sigma)=\p_2(0_{L_0},t\cdot 1_\AA\cdot t^\sigma)=\p_2(a,t)\cdot \p_2(0_{L_0},1_\AA)\cdot \p_2(a,t)^{\tilde{\sigma}}=\tilde{N}\big(\p_2(a,t)\big)\ .$$
\item By part (a) and lemma \ref{221}, we have
\begin{align*}
\phi\big(N(t)\big)+\phi\big(T(t)\big)+1_{\tilde{A}}&=\phi\big((t+1_\AA)(t+1_\AA)^\sigma\big)=\phi\big(N(t+1_\AA)\big) \\
&=\tilde{N}\big(\p_2(a,t+1_\AA)\big)=\tilde{N}\big(\p_2(a,t)+\phi(1_\AA)\big)\\
&=\tilde{N}\big(\p_2(a,t)\big)+\tilde{T}\big(\p_2(a,t)\big)+1_{\tilde{\AA}}\ ,
\end{align*}
hence
$$\phi\big(T(t)\big)=\tilde{T}\big(\p_2(a,t)\big)$$
by part (a) again.\qed
\end{enumerate}
\end{bewzwei}

\begin{lemma}\label{226}
Let $a\in L_0^*$ and $t\in X_a$. Then the map
$$\phi_{(a,t)}:\EE_a\to \tilde{\EE}_{\p_1(a)},\ x+t\cdot y \mapsto \phi(x)+\p_2(a,t)\cdot \phi(y)\qquad (x,y\in \FF)$$
is an isomorphism of fields such that
$$\forall\ u\in X_a:\qquad  \p_2(a,u)=\phi_{(a,t)}(u)\ .$$
\end{lemma}

\begin{bew}
By lemma \ref{224}, we have
\begin{align*} \phi\big(N(t)\big)=\tilde{N}\big(\p_2(a,t)\big)\ , && \phi\big(T(t)\big)=\tilde{T}\big(\p_2(a,t)\big)\end{align*}
so that we may apply lemma \ref{225}. By lemma \ref{221}, we have
$$\forall\ s\in \FF:\qquad \p_2(a,t+s)=\p_2(a,t)+\phi(s)=\phi_{(a,t)}(t)+\phi_{(a,t)}(s)=\phi_{(a,t)}(t+s)\ .$$
\qed
\end{bew}

\begin{kor}\label{304}
Let $a\in L_0^*$. Given $t,u\in X_a$, we have $$\phi_a:=\phi_{(a,t)}=\phi_{(a,u)}\ .$$
\end{kor}

\begin{bew}
We have $\phi_{(a,u)}(u)=\p_2(a,u)=\phi_{(a,t)}(u)$ and $u\notin \FF$.\qed
\end{bew}

\begin{bem} Let $a\in L_0^*$.
\begin{enumerate}[label=(\alph*)]
\item\label{236} Notice that $\phi_a$ is an extension of the isomorphism $\phi:\FF\to \tilde{\FF}$ of fields.
\item\label{234} We have $\phi_a(\EE_a)=\tilde{\EE}_{\p_1(a)}$.
\item \label{223}
Given $t\in \AA$, we have
\begin{align*}
t^{-1} (t^{\sigma})^2&=\frac{(t^\sigma)^3}{N(t)}\ \ \overset{\mathclap{\eqref{222}}}{=}\ \ \frac{\big(T(t)^2-N(t)\big)\cdot t^\sigma-T(t)N(t)}{N(t)} \\
&=\frac{\big(T(t)^2-N(t)\big) \big(T(t)-t\big)-T(t)N(t)}{N(t)}
=\frac{T(t)^3-2T(t)N(t)}{N(t)}+\frac{N(t)-T(t)^2}{N(t)}\cdot t\ .
\end{align*}
\item \label{307} Assume $|\FF|\geq 3$. By lemma \ref{306}, there is an element $t\in X_a$ such that $T(t)^2\neq N(t)$.
\end{enumerate}
\end{bem}

\begin{lemma}\label{227}
Assume $|\FF|\geq 3$. Let $a\in L_0^*$ and $t\in X_a$ be as in remark \ref{307}, i.e., we have $N(t)\neq T(t)^2$ and thus $\phi\big(N(t)\big)\neq \phi\big(T(t)\big)^2$. Then we have
$$\p_1(a\cdot t)=\p_1(a)\cdot \p_2(a,t)\ .$$
\end{lemma}

\begin{bew}
By remark \ref{223}, equation \eqref{272} and identity \eqref{169} with $b=a$, we have
\begin{align*}
&\ \qquad \p_1(a)\cdot \frac{\phi\big(T(t)\big)^3-2\phi\big(T(t)\big)\phi\big(N(t)\big)}{\phi\big(N(t)\big)}+\p_1(a\cdot t)\cdot \frac{\phi\big(N(t)\big)-\phi\big(T(t)\big)^2}{\phi\big(N(t)\big)} \\
\overset{\mathclap{\textrm{\ref{205}}}}{=}&\ \qquad \p_1\left( a\cdot \left( \frac{T(t)^3-2T(t)N(t)}{N(t)}+\frac{N(t)-T(t)^2}{N(t)}\cdot t \right) \right) \\
\overset{\mathclap{\textrm{\ref{223}}}}{\underset{\mathclap{\eqref{272},\eqref{169}}}{=}}&\ \qquad \p_1(a)\cdot\left( \frac{\tilde{T}\big(\p_2(a,t)\big)^3-2\tilde{T}\big(\p_2(a,t)\big)\tilde{N}\big(\p_2(a,t)\big)}{\tilde{N}\big(\p_2(a,t)\big)}+\frac{\tilde{N}\big(\p_2(a,t)\big)-\tilde{T}\big(\p_2(a,t)\big)^2}{\tilde{N}\big(\p_2(a,t)\big)}\cdot \p_2(a,t)\right) \\
\overset{\mathclap{\textrm{\ref{224}}}}{=}&\ \qquad \p_1(a)\cdot \frac{\phi\big(T(t)\big)^3-2\phi\big(T(t)\big)\phi\big(N(t)\big)}{\phi\big(N(t)\big)}+\p_1(a)\cdot \p_2(a,t)\cdot \frac{\phi\big(N(t)\big)-\phi\big(T(t)\big)^2}{\phi\big(N(t)\big)}
\end{align*}
Notice that we may cancel scalars by assumption.\qed
\end{bew}

\begin{bem}\label{308}
By corollary \ref{212}, we have $\p_1(R_a)= \tilde{R}_{\p_1(a)}$ even in the case $|\FF|=2$. Notice that we have $$\AA=\EE_a\cong \FF_4\cong  \tilde{\EE}_{\p_1(a)}=\tilde{\AA}$$ in this situation which shows that we have $\p_1(\langle a \rangle_\AA)=\langle \p_1(a)\rangle_{\tilde{\AA}}$. As a consequence, the map $\phi_a:\AA\to \tilde{\AA}$ defined by
$$\forall\ t\in \AA:\qquad \p_1(a\cdot t)=\p_1(a)\cdot \phi_a(t)$$
is an isomorphism of skew-fields. However, it doesn't necessarily satisfy $\phi_a(t)=\p_2(a,t)$ for each $t\in X_a$. We will discuss this case later and go on with assuming $|\FF|\geq 3$.
\end{bem}

\begin{prop}\label{228}
Given $a\in L_0^*$, the map 
$$(\p_1,\phi_{a}): (R_a,\EE_a)\to (\tilde{R}_{\p_1(a)},\tilde{\EE}_{\p_1(a)})$$
is an isomorphism of vector spaces such that $\phi_a(t)=\p_2(a,t)$ for each $t\in X_a$.
\end{prop}

\begin{bew}
Let $t\in X_a$ be as in remark \ref{307}. By remark \ref{205}, lemma \ref{227} and lemma \ref{226}, we have
$$\forall\ x,y\in \FF:\qquad \p_1\big(a\cdot (x+ty)\big)=\p_1(a)\cdot \big(\phi(x)+\p_2(a,t)\phi(y)\big)=\p_1(a)\cdot \phi_{a}(x+ty)\ .$$
The second assertion results from lemma \ref{226}.\qed
\end{bew}

\begin{lemma}\label{259}
Let $(\AA,\FF,\sigma)$ and $(\tilde{\AA},\tilde{\FF},\tilde{\sigma})$ both be quadratic of type (iii). Then we have
$$\forall\ a,b\in L_0^*:\qquad \phi:=\phi_a=\phi_b\ .$$
\end{lemma}

\begin{bew}
By assumption, we have $\AA=\EE_a=\EE_b$ and $\tilde{\AA}=\tilde{\EE}_{\p_1(a)}=\tilde{\EE}_{\p_1(b)}$ for all $a,b\in L_0^*$.
\begin{itemize}
\item Let $b=a\cdot s$ for some $s\in \AA^*$. Given $t\in \AA$, we have
\begin{align*}
\p_1(a)\cdot \phi_a(t)&=\p_1(a\cdot t)=\p_1(as \cdot s^{-1}t)\\
&=\p_1(b)\cdot \phi_b(s^{-1}t)=\p_1(b)\cdot \phi_b(s^{-1})\cdot \phi_b(t)=\p_1(a)\cdot \phi_b(t)\ .\end{align*}
\item Let $b\notin \langle a\rangle_\AA=R_a$, hence $\p_1(b)\notin \tilde{R}_{\p_1(a)}=\langle \p_1(a)\rangle_{\tilde{\AA}}$ by proposition \ref{228}. Given $t\in \AA$, we have
\begin{align*}
\p_1(a)\cdot \phi_a(t)+\p_1(b)\cdot \phi_b(t)&=\p_1(a\cdot t)+\p_1(b\cdot t)=\p_1\big((a+b)\cdot t\big)\\
&=\p_1(a+b)\cdot \phi_{a+b}(t)=\p_1(a)\cdot \phi_{a+b}(t)+\p_1(b)\cdot \phi_{a+b}(t)\ ,
\end{align*}
thus $$\phi_a(t)=\phi_{a+b}(t)=\phi_b(t)$$
\end{itemize}\qed
\end{bew}

\begin{bem}\ 
\begin{enumerate}[label=(\alph*)] 
\item \label{305} The proof shows that we have $$\forall\ t\in \EE_a:\qquad \phi_a=\phi_{at}\ ,$$ independent of the type of $(\AA,\FF,\sigma)$.
\item \label{233} Notice that $\phi$ is an extension of the (Jordan) isomorphism $\phi:\FF\to \tilde{\FF}$.
\item We state the theorem using the general notation.
\end{enumerate}
\end{bem}

\newpage

\begin{satz}\label{271}
Let $(\KK,\KK_0,\sigma)$ and $(\tilde{\KK},\tilde{\KK}_0,\tilde{\sigma})$ both be quadratic of type (iii) such that $|\KK_0|\geq 3 $. Then the map $\Phi:\Xi\to\tilde{\Xi}$ defined by
\begin{align*}
\Phi:=(\p_1,\phi): (L_0,\KK)\to(\tilde{L},\tilde{\KK}),\ (a,t)\mapsto \big(\p_1(a), \phi(t)\big)
\end{align*}
is an isomorphism of pseudo-quadratic spaces satisfying
$$\forall\ (a,t)\in T:\qquad \Phi(a,t)=\gamma(a,t)\ .$$
In particular, we have $\KK\cong \tilde{\KK}$.
\end{satz} 

\begin{bewzwei}\ 
\begin{itemize}
\item By proposition \ref{228} and lemma \ref{259}, the map $(\p_1,\phi):(L_0,\KK)\to(\tilde{L}_0,\tilde{\KK})$ is an isomorphism of vector spaces.
\item By remark \ref{233}, we have $\phi(\KK_0)=\tilde{\KK}_0$.
\item By corollary \ref{251}, we have $\phi\circ \sigma=\tilde{\sigma}\circ \phi$.
\item Let $a\in L_0^*$ and $t\in X_a$. By proposition \ref{228}, we have $$\p_2(a,t)=\phi_a(t)={\phi}(t)\ ,$$ hence
\begin{align*}
\tilde{q}\big(\p_1(a)\big)&\in \p_2(a,t)+\tilde{\KK}_0={\phi}(t)+\phi(\KK_0)= {\phi}\big(t+\KK_0\big)={\phi}\big(q(a)+\KK_0\big)={\phi}\big(q(a)\big)+\tilde{\KK}_0\ .
\end{align*}
\item Given $(a,t)\in T$, we have
$$\Phi(a,t)=\big(\p_1(a),{\phi}(t)\big)=\big(\p_1(a),\p_2(a,t)\big)=\gamma(a,t)\ .$$\qed
\end{itemize}
\end{bewzwei}

\chapter{Small Dimensions II} \label{268}

Now we prove theorem \ref{254} for the case that $(\AA,\FF,\sigma)$ is quadratic of type (iv). As in the previous paragraph, we exploit the identities \eqref{169} and \eqref{170} to show that $\phi_a:\EE_a\to \tilde{\EE}_{\p_1(a)}$ is induced by an isomorphism $\phi_{a,e}:\AA\to \bar{\AA}$ of skew-fields, where $a\in L_0^*$ is a separable element and $e\in \EE_a^\bot$. 

Notice that we restrict to a separable element $a\in L_0^*$ because of the helpful decomposition $\AA=\EE_a\oplus e\EE_a$ with $e\EE_a=\EE_a^\bot$. If we assume $\p_1(a\cdot e)=\p_1(a)\cdot \tilde{e}$ for some $\tilde{e}\in \tilde{\AA}$, the isomorphisms $\phi_a$ and $\phi_{ae}$ of fields induce a map
$$\phi_{a,e}:\AA\to\tilde{\AA},\ s+et\mapsto \phi_a(s)+\tilde{e}\phi_{ae}(t)\qquad (s,t\in\EE_a)$$
which turns out to be an isomorphism of skew-fields. As a consequence, the two isomorphisms $\p_1:R_a\to \tilde{R}_{\p_1(a)}$ and $\p_1:R_{ae}\to \tilde{R}_{\p_1(ae)}$ of vector spaces over $\EE_a$ are induced by an isomorphism $\p_1:\langle a\rangle_\AA\to \langle \p_1(a)\rangle_{\tilde{\AA}}$ of vector spaces over $\AA$. 

\begin{no}
Throughout this chapter, let $(\AA,\FF,\sigma)$ be quadratic of type (iv), and let $(\tilde{\AA},\tilde{\FF},\tilde{\sigma})$ be quadratic of type (iii) or (iv). Until corollary \ref{267}, let $a\in L_0^*$ be separable and let $e\in \EE_a^\bot$.
\end{no}

\begin{lemma}\label{274}
Given $t\in X_a\subseteq \EE_a$, we have
$$ae\cdot t=ae\cdot t^\sigma-a\cdot t^{-1}f(a,ae)t^\sigma\ .$$
\end{lemma}

\newpage

\begin{bew}
We have
\begin{align*}
ae\cdot t^\sigma-a\cdot t^{-1}f(a,ae)t^\sigma&=ae\cdot t^\sigma-a\cdot t^{-1}f(a,a)et^\sigma \\
&=ae\cdot t^\sigma -a\cdot  t^{-1}(t-t^\sigma)et^\sigma =ae\cdot t^\sigma- ae\cdot  (t^\sigma-t)t^{-\sigma}t^\sigma=ae\cdot t\ .
\end{align*}\qed
\end{bew}

\begin{lemma}\label{230}
We have $$X_{ae}=-N(e)\cdot X_a\ .$$
\end{lemma}

\begin{bew}
We have
$$q(ae)\equiv e^\sigma q(a)e=N(e)q(a)^\sigma \mod \FF\ ,$$
hence
\begin{align*} -N(e)\cdot X_a&=-N(e)\cdot \{ s+q(a) \mid s\in \FF \}=-N(e)\cdot \{ s-q(a)^\sigma \mid s\in \FF \} \\
&=\{ -N(e) s+N(e)q(a)^\sigma \mid s\in \FF \}= \{ s+q(ae) \mid s\in \FF \}=X_{ae}\ .
\end{align*}
\qed
\end{bew}

\begin{lemma}\label{261}
We have $$\EE_{ae}=\EE_a\ .$$
\end{lemma}

\begin{bew}
We have
$$q(ae)\in N(e)q(a)^\sigma+\FF\subseteq \EE_a\ .$$\qed
\end{bew}

\begin{lemma}
The isomorphisms $\phi_{ae}$ and $\phi_a$ of fields as in remark \ref{304} satisfy
$$\phi_{ae}(\EE_a)=\phi_a(\EE_a)\ .$$
\end{lemma}

\begin{bew}
Let $t\in X_a\subseteq \EE_a=\EE_{ae}$. Then we have $-N(e)t\in X_{ae}$ by lemma \ref{230}, hence
\begin{align*}
(a,t)\in T\ , && \big(ae,-N(e)t\big)\in T\ .
\end{align*}
Notice that we have
$$\phi_a(t)\phi_a(t)^{\tilde{\sigma}}=\tilde{N}\big(\phi_a(t)\big)=\tilde{N}\big(\p_2(a,t)\big)=\phi\big(N(t)\big)$$
by lemma \ref{273}. Remarks \ref{236},  \ref{305} and identity \eqref{170} with $(b,v):=\big(ae,-N(e)t\big)$ yield
\begin{align*}
\phi_{ae}(t)\phi_a(t)\phi_a(t)^{\tilde{\sigma}}\phi\big(-N(e)\big)&=\phi_{ae}(t) \phi\big(N(t)\big)\phi \big(-N(e)\big) \\
&=\phi_{ae}(t)\phi_{ae}\big(-N(e)\big)\phi_{ae}\big(N(t)\big)=\phi_{ae}\big(t(-N(e))N(t)\big) \\
&=\phi_{aet}\big(t\cdot (-N(e)t)\cdot t^\sigma\big)=\p_2\big(ae\cdot t,t\cdot (-N(e)t)\cdot t^\sigma\big) \\
&=\p_2\big(ae\cdot t^\sigma-a\cdot t^{-1}f(a,ae)t^\sigma,t\cdot (-N(e)t)\cdot t^\sigma\big) \\
&=\p_2(a,t)\p_2\big(ae,-N(e)t\big)\p_2(a,t)^{\tilde{\sigma}}=\phi_a(t) \phi_{ae}\big(-N(e) t\big)\phi_a(t)^{\tilde{\sigma}} \\
&=\phi_a(t)\phi_{ae}(t)\phi_{ae}\big(-N(e)\big)\phi_a(t)^{\tilde{\sigma}}=\phi_a(t)\phi_{ae}(t)\phi_a(t)^{\tilde{\sigma}}\phi\big(-N(e)\big)\ ,
\end{align*}
which implies $\phi_{ae}(t)\phi_a(t)=\phi_a(t)\phi_{ae}(t)$ and thus
\begin{align*}
 \phi_{ae}(\EE_a)=\phi_{ae}(\EE_{ae})=\tilde{\EE}_{\p_1(ae)}=\tilde{\EE}_{\phi_{ae}(t)}=\tilde{\EE}_{\phi_a(t)}=\tilde{\EE}_{\p_1(a)}=\phi_a(\EE_a)\ ,
\end{align*}
cf. lemma \ref{261}, remark \ref{234}, lemma \ref{230}, remark \ref{211} and corollary \ref{231}.\qed
\end{bew}

\begin{kor}\label{237}
We have $$\phi_{ae}\in \{\phi_a,\tilde{\sigma}\circ \phi_a\}\ .$$
\end{kor}

\begin{bew}
The field $\tilde{\EE}:=\phi_{ae}(\EE_a)=\phi_a(\EE_a)$ is a quadratic separable extension with $\Aut(\tilde{\EE}:\tilde{\FF})=\langle \tilde{\sigma}\rangle$, and by remark \ref{236}, we have
\begin{align*}
\phi_{ae}^{\vphantom{-1}}\circ \phi_a^{-1}\in \Aut(\tilde{\EE}:\tilde{\FF})\ .
\end{align*}
\qed
\end{bew}

\begin{bem}
The following lemma provides the assumptions which are necessary to show that we can extend the isomorphism $\phi_a:\EE_a\to \tilde{\EE}_{\p_1(a)}$ of fields to an isomorphism $\phi:\AA\to \tilde{\AA}$ of skew-fields which induces $\phi_{ae}$. As a consequence, the map $\p_1: \langle a\rangle_\AA\to \langle\p_1(a)\rangle_{\tilde{\AA}}$ is an isomorphism of vector spaces.
\end{bem}

\begin{lemma}\label{256}
Suppose that we have $\p_1(ae)=\p_1(a)\tilde{e}$ for some $\tilde{e}\in \tilde{\AA}$ and let $t\in X_a$. Then the following holds:
\begin{enumerate}[label=(\alph*)]
\item \label{258} We have $\phi:=\phi_{ae}=\phi_a$.
\item We have
\begin{align*}
\tilde{f}\big(\p_1(ae),\p_1(ae)\big)=-\phi\big(N(e)\big)\tilde{f}\big(\p_1(a),\p_1(a)\big)\ .
\end{align*}
\item We have
$$\tilde{f}\big(\p_1(a),\p_1(a)\big)\tilde{e}\\
=-\tilde{e}\tilde{f}\big(\p_1(a),\p_1(a)\big)\ .$$
\item We have $\tilde{N}(\tilde{e})=\phi\big(N(e)\big)$.
\item We have $\tilde{e}\in \tilde{\EE}_{\p_1(a)}^\bot$.
\end{enumerate}
\end{lemma}

\begin{bewzwei}\ 
\begin{enumerate}[label=(\alph*)]
\item By corollary \ref{197}, $\p_1(a)$ is separable, hence
$$\tilde{f}\big(\p_1(a),\p_1(ae)\big)=\tilde{f}\big(\p_1(a),\p_1(a)\big)\tilde{e}\neq 0_{\tilde{\AA}}\ .$$
By identity \eqref{169} with $(b,v):=\big(ae,-N(e)t\big)$ and lemma \ref{274}, we have
\begin{align}
\p_1(ae)\phi_{ae}(t)&=\p_1(ae\cdot t)=\p_1\big(ae\cdot t^\sigma-a\cdot t^{-1}f(a,ae)t^\sigma\big) \nonumber \\
&=\p_1(ae)\cdot \p_2(a,t)^{\tilde{\sigma}}-\p_1(a)\cdot \p_2(a,t)^{-1}\tilde{f}\big(\p_1(a),\p_1(ae)\big)\p_2(a,t)^{\tilde{\sigma}} \label{238} \\
&=\p_1(ae)\cdot \phi_a(t)^{\tilde{\sigma}}-\p_1(a)\cdot \underbrace{\phi_a(t)^{-1}\tilde{f}\big(\p_1(a),\p_1(ae)\big)\phi_a(t)^{\tilde{\sigma}}}_{\neq 0}\ , \nonumber
\end{align}
hence
$$\phi_{ae}(t)\neq \phi_a(t)^{\tilde{\sigma}}$$
for each $t\in X_a$ and thus $\phi_{ae}\neq \tilde{\sigma}\circ \phi_a$. Now corollary \ref{237} yields
$$\phi_{ae}=\phi_a\ .$$
\item By lemma \ref{230}, we have $-N(e)t\in X_{ae}$, hence
\begin{align*}
\tilde{f}\big(\p_1(ae),\p_1(ae)\big)&=\phi_{ae}\big(-N(e)t\big)-\phi_{ae}\big(-N(e)t\big)^{\tilde{\sigma}}\\
&=-\phi\big(N(e)\big)\big(\phi_a(t)-\phi_a(t)^{\tilde{\sigma}}\big)=-\phi\big(N(e)\big)\tilde{f}\big(\p_1(a),\p_1(a)\big)\ ,
\end{align*}
cf. remark \ref{213} and proposition \ref{228}.
\item By (a), we have
\begin{align}\begin{aligned}
\p_1(a)\tilde{e}\cdot \phi(t)&=\p_1(ae)\cdot \phi(t)\overset{\mathclap{\eqref{238}}}{=}\p_1(ae)\cdot \phi(t)^{\tilde{\sigma}}-\p_1(a)\cdot {\phi(t)^{-1}\tilde{f}\big(\p_1(a),\p_1(ae)\big)\phi(t)^{\tilde{\sigma}}} \\
&=\p_1(ae)\cdot \phi(t)^{\tilde{\sigma}}-\p_1(a)\cdot \phi(t)^{-1}\tilde{f}\big(\p_1(a),\p_1(a)\big)\tilde{e}\phi(t)^{\tilde{\sigma}} \\
&=\p_1(a)\tilde{e}\cdot \phi(t)^{\tilde{\sigma}}-\p_1(a)\cdot \phi(t)^{-1}\big(\phi(t)-\phi(t)^{\tilde{\sigma}}\big)\tilde{e}\phi(t)^{\tilde{\sigma}} \\
&=\p_1(a)\cdot \phi(t)^{-1}\phi(t)^{\tilde{\sigma}}\tilde{e}\phi(t)^{\tilde{\sigma}}\ ,
\end{aligned} \label{255}\end{align}
hence
\begin{align*}
\tilde{e}\phi(t)=\phi(t)^{-1}\phi(t)^{\tilde{\sigma}}\tilde{e}\phi(t)^{\tilde{\sigma}}\ , && \phi(t)\tilde{e}\phi(t)=\phi(t)^{\tilde{\sigma}}\tilde{e}\phi(t)^{\tilde{\sigma}}
\end{align*}
Since $t\in X_a$ is arbitrary, we may replace $t$ by $t+1_\AA$ to obtain
\begin{align*} \tilde{e}\phi(t)+\phi(t)\tilde{e}=\tilde{e}\phi(t)^{\tilde{\sigma}}+\phi(t)^{\tilde{\sigma}}\tilde{e}\ , && \tilde{e}\big(\phi(t)-\phi(t)^{\tilde{\sigma}}\big)=-\big( \phi(t)-\phi(t)^{\tilde{\sigma}}\big)\tilde{e}\ .\end{align*}
\item We have
\begin{align*}
-\tilde{e}\tilde{N}(\tilde{e})^{-1}\phi\big(N(e)\big)\tilde{f}\big(\p_1(a),\p_1(a)\big)&\overset{(\textrm{b})}{=}\tilde{e}^{-\tilde{\sigma}}\tilde{f}\big(\p_1(ae),\p_1(ae)\big)=\tilde{f}(\p_1(ae)\tilde{e}^{-1},\p_1(ae)\big)\\
&\vphantom{\overset{\textrm{(c)}}{=}}=\tilde{f}\big(\p_1(a),\p_1(ae)\big)=\tilde{f}\big(\p_1(a),\p_1(a)\big)\tilde{e}\\
&\overset{\textrm{(c)}}{=}-\tilde{e}\tilde{f}\big(\p_1(a),\p_1(a)\big)\ ,
\end{align*}
hence
\begin{align*}
\tilde{N}(\tilde{e})^{-1}\phi\big(N(e)\big)=1_{\tilde{\AA}}\ , && \tilde{N}(\tilde{e})=\phi\big(N(e)\big)\ .
\end{align*}
\item We have
$$\p_1(a)\cdot \tilde{e}\phi(t)\overset{\eqref{255}}{=}\p_1(a)\cdot \tilde{e}\phi(t)^{\tilde{\sigma}}-\p_1(a)\cdot \phi(t)^{-1}\tilde{f}\big(\p_1(a),\p_1(a)\big)\tilde{e}\phi(t)^{\tilde{\sigma}}\ ,$$
hence
\begin{align*}\p_1(a)\cdot \tilde{e}\tilde{f}\big(\p_1(a),\p_1(a)\big)&=-\p_1(a)\cdot \phi(t)^{-1}\tilde{f}\big(\p_1(a),\p_1(a)\big)\tilde{e}\phi(t)^{\tilde{\sigma}} \\
&\overset{\mathclap{\textrm{(c)}}}{=}\p_1(a)\cdot  \phi(t)^{-1}\tilde{e}\phi(t)^{\tilde{\sigma}}\tilde{f}\big(\p_1(a),\p_1(a)\big)
\end{align*}
which yields
\begin{align*}
\tilde{e}=\phi(t)^{-1}\tilde{e}\phi(t)^{\tilde{\sigma}}\ , && \phi(t)\tilde{e}=\tilde{e}\phi(t)^{\tilde{\sigma}}\ .
\end{align*}
Because of $\phi(t)\in \tilde{\EE}_{\p_1(a)}\sm \tilde{\FF}$, we have $\tilde{e}\in \tilde{\EE}_{\p_1(a)}^\bot$ by lemma \ref{253}.
\end{enumerate}\qed
\end{bewzwei}

\begin{bem}\label{262}
Suppose that we have $\p_1(ae)=\p_1(a)\tilde{e}$ for some $\tilde{e}\in \tilde{\AA}$. Since we have $\p_1(R_a)=\tilde{R}_{\p_1(a)}$ and $a\cdot e\notin R_a$, the skew-field $\tilde{\AA}$ is necessarily a quaternion division algebra since $(\tilde{\AA},\tilde{\FF},\tilde{\sigma})$ is quadratic of type (iii) or (iv). Moreover, we have $\p_1(\langle a\rangle_\AA)=\langle \p_1(a)\rangle_{\tilde{\AA}}$.
\end{bem}

\begin{lemma}
Suppose that we have $\p_1(ae)=\p_1(a)\tilde{e}$ for some $\tilde{e}\in \tilde{\AA}$. Then the map $$\phi_{a,e}:\AA\to \tilde{\AA},\ x+ey\mapsto \phi(x)+\tilde{e}\phi(y)\qquad (x,y\in \EE_a)$$
is an isomorphism of skew-fields such that
$$\forall\ t\in X_a:\qquad  \p_2(a,t)=\phi_{a,e}(t)\ .$$
\end{lemma}

\begin{bew}
By lemma \ref{256}, we have
\begin{align*} \phi\big(N(e)\big)=\tilde{N}(\tilde{e})\ , && \tilde{e}\in \tilde{\EE}_{\p_1(a)}^\bot\end{align*}
so that we may apply lemma \ref{257}. By lemma \ref{226}, we have
$$\forall\ t\in X_a\subseteq \EE_a:\qquad \p_2(a,t)=\phi_a(t)=\phi_{a,e}(t)\ .$$
\qed
\end{bew}

\begin{prop}\label{260}
Suppose that we have $\p_1(ae)=\p_1(a)\tilde{e}$ for some $\tilde{e}\in \tilde{\AA}$. Then the map
$$(\p_1,\phi_{a,e}): (\langle a\rangle_\AA,\AA)\to (\langle \p_1(a)\rangle_{\tilde{\AA}},\tilde{\AA})$$
is an isomorphism of vector spaces.
\end{prop}

\begin{bew}
By lemma \ref{261} and lemma \ref{258}, we have
$$\forall\ x,y\in \EE_a:\qquad \p_1\big(a\cdot (x+ey)\big)=\p_1(a)\cdot \big(\phi(x)+\tilde{e}\phi(y)\big)=\p_1(a)\cdot \phi_{a,e}(x+ey)\ .$$
\qed
\end{bew}

\begin{kor}\label{267}
Suppose that we have $\p_1(ae)=\p_1(a)\tilde{e}$ for some $\tilde{e}\in \tilde{\AA}$ and let $f\in \EE_a^\bot$. Then we have $\p_1(af)=\p_1(a)\tilde{f}$ for some $\tilde{f}\in \tilde{\AA}$ and $\phi_a:=\phi_{a,e}=\phi_{a,f}$.
\end{kor}

\begin{bew}
By remark \ref{262}, the first assertion holds, hence $\phi_{a,f}$ is well-defined. Given $x\in \AA$, we have
$$ \p_1(a)\cdot \phi_{a,f}(x)=\p_1(a\cdot x)=\p_1(a)\cdot \phi_{a,e}(x)\ .$$
\qed
\end{bew}

\begin{lemma}\label{263}
Suppose that we have $\p_1(\langle a\rangle_\AA)=\langle \p_1(a)\rangle_{\tilde{\AA}}$ for each separable element $a\in L_0^*$ and let $a,b\in L_0^*$ be separable. Then we have
$$\phi:=\phi_a=\phi_b\ .$$
\end{lemma}

\begin{bewzwei}\ 
\begin{itemize}
\item Let $b=a\cdot s$ for some $s\in \AA^*$. Given $t\in \AA$, we have
\begin{align*}
\p_1(a)\cdot \phi_a(t)&=\p_1(a\cdot t)=\p_1(as \cdot s^{-1}t)\\
&=\p_1(b)\cdot \phi_b(s^{-1}t)=\p_1(b)\cdot \phi_b(s^{-1})\cdot \phi_b(t)=\p_1(a)\cdot \phi_b(t)\ .\end{align*}
\item Let $b\notin \langle a\rangle_\AA$, hence $\p_1(b)\notin \langle \p_1(a)\rangle_{\tilde{\AA}}$. By lemma \ref{196}, there are $a'\in R_{a}$ and $b'\in R_{b}$ such that $a'+b'$ is separable. Given $t\in \AA$, we have
\begin{align*}
\p_1(a')\cdot \phi_{a}(t)+\p_1(b')\cdot \phi_{b}(t)&=\p_1(a')\cdot \phi_{a'}(t)+\p_1(b')\cdot \phi_{b'}(t)=\p_1(a'\cdot t)+\p_1(b'\cdot t) \\
&=\p_1\big((a'+b')\cdot t\big)=\p_1(a'+b')\cdot \phi_{a'+b'}(t)\\
&=\p_1(a')\cdot \phi_{a'+b'}(t)+\p_1(b')\cdot \phi_{a'+b'}(t)\ ,
\end{align*}
thus
$$\phi_a(t)=\phi_{a'+b'}(t)=\phi_b(t)\ .$$
\end{itemize}\qed
\end{bewzwei}

\begin{prop}\label{265}
Suppose that we have $\p_1(\langle a\rangle_\AA)=\langle \p_1(a)\rangle_{\tilde{\AA}}$ for each separable element $a\in L_0^*$ and that $L_0$ has a basis $\{a_i \mid i\in I\}$ of separable elements. Then the map
$$(\p_1,\phi):(L_0,\AA)\to (\tilde{L}_0,\tilde{\AA})$$
is an isomorphism of vector spaces.
\end{prop}

\begin{bew}
Let $a=\sum_{i\in I} a_i \lambda_i \in L_0$ and $t\in \AA$. By lemma \ref{263}, we have
\begin{align*}
\p_1(a\cdot t)&=\p_1 \Big(\big(\sum_{i\in I} a_i \lambda_i\cdot  t\big)\Big)=\sum_{i\in I} \p_1(a_i\lambda_i\cdot   t) \\
&=\sum_{i\in I} \p_1(a_i\lambda_i )\cdot  \phi_{a_i \lambda_i}(t)=\sum_{i\in I} \p_1(a_i\lambda_i )\cdot  \phi(t)=\Big(\sum_{i\in I} \p_1(a_i\lambda_i)\Big)\cdot \phi(t)=\p_1(a)\cdot \phi(t)\ .
\end{align*}
\qed
\end{bew}

\begin{bem}\label{264}
Let $\dim_\AA L_0=2$ and $a\in L_0^*$. Then $\langle a\rangle_\AA$ is of the form $L_{a,b}=\langle R_a,R_b\rangle_\FF$ for some element $b\in L_0$, and for each $c\in L_{a,b}$, we have $R_c\subseteq L_{a,b}$. Conversely, if we have a subspace $L_{a,b}=\langle R_a,R_b\rangle_\FF$ such that $R_c\subseteq L_{a,b}$ for each $c\in L_{a,b}$, we have $L_{a,b}=\langle a\rangle_{\AA_i}$ for one of the three pairwise non-isomorphic quaternion division algebras $\AA=:\AA_1,\AA_2,\AA_3$ mentioned in \cite{DM}, cf. proposition (3.2) in \cite{DM}. 
\end{bem}

\begin{lemma}\label{266}
Let $\dim_\AA L_0=2$, let $a\in L_0^*$ be separable and let $\tilde{\AA}$ be a quaternion division algebra. Then we have
$$\p_1(\langle a\rangle_\AA)=\langle \p_1(a)\rangle_{\tilde{\AA}_i}$$
for one of the three pairwise non-isomorphic quaternion division algebras $\tilde{\AA}=:\tilde{\AA}_1,\tilde{\AA}_2,\tilde{\AA}_3$ mentioned in \cite{DM}, and the map $(\p_1,\phi_a):(\langle a\rangle_\AA, \AA)\to (\langle \p_1(a)\rangle_{\tilde{\AA}_i},\tilde{\AA}_i)$ is an isomorphism of vector spaces. In particular, we have $\tilde{\AA}_i=\tilde{\AA}$ if we assume $\tilde{\AA}\cong \AA$.
\end{lemma}

\begin{bew}
By remark \ref{264}, the subspace $\langle a\rangle_\AA$ is of the form $L_{a,b}=\langle R_a,R_b\rangle_\FF$ for some element $b\in L_0$, and for each $c\in L_{a,b}$, we have $R_c\subseteq L_{a,b}$. As $\p_1:L_0\to \tilde{L}_0$ is $\FF$-linear, we have
$$\p_1(L_{a,b})=\p_1(\langle R_a,R_b\rangle_\FF)=\langle \tilde{R}_{\p_1(a)},\tilde{R}_{\p_1(b)}\rangle_{\tilde{\FF}}=\tilde{L}_{\p_1(a),\p_1(b)}$$
by corollary \ref{212}. Now let $\tilde{c}\in \tilde{L}_{\p_1(a),\p_1(b)}$. Because of $c:=\p_1^{-1}(\tilde{c})\in L_{a,b}$ we have $R_c\subseteq L_{a,b}$, hence
$$\tilde{R}_{\tilde{c}}=\p_1(R_c)\subseteq \p_1(L_{a,b})=\tilde{L}_{\p_1(a),\p_1(b)}\ ,$$
which shows that the first assertion holds. The second assertion holds by proposition \ref{260}. In particular, we have $\tilde{\AA}_i\cong \AA$. As $\AA_1,\AA_2,\AA_3$ are pairwise non-isomorphic, we must have $\tilde{\AA}_i=\tilde{\AA}$ if we assume $\tilde{\AA}\cong \AA$.\qed
\end{bew}

\begin{bem}
Once again we switch to the general notation.
\end{bem}

\begin{satz}\label{270}
Let $(\KK,\KK_0,\sigma)$ be quadratic of type (iv), let $\dim_\KK L_0\leq 2$ and $\tilde{\KK}\cong \KK$. Then the map $\Phi:\Xi\to\tilde{\Xi}$ defined by
\begin{align*}
\Phi:=(\p_1,\phi): (L_0,\KK)\to(\tilde{L},\tilde{\KK}),\ (a,t)\mapsto \big(\p_1(a), \phi(t)\big)
\end{align*}
is an isomorphism of pseudo-quadratic spaces satisfying
$$\forall\ (a,t)\in T:\qquad \Phi(a,t)=\gamma(a,t)\ .$$
\end{satz} 

\begin{bewzwei}\ 
\begin{itemize}
\item In the case $\dim_\KK L_0=1$, each $a\in L_0^*$ is separable, and by assumption, we have 
$$\p_1(\langle a\rangle_\KK)=\p_1(L_0)=\tilde{L}_0=\langle \p_1(a) \rangle_{\tilde{\KK}}$$
so that we may apply proposition \ref{265}. Now assume $\dim_\KK L_0=2$. By theorem (6.3) in chapter 7 of \cite{WSch}, $L_0$ has an orthogonal basis, and lemma \ref{266} implies that we have $\p_1(\langle a\rangle_\KK)=\langle \p_1(a)\rangle_{\tilde{\KK}}$ for each separable element $a\in L_0^*$ so that we may apply proposition \ref{265}. In both cases, the map $(\p_1,\phi):(L_0,\KK)\to (\tilde{L}_0,\tilde{\KK})$ is an isomorphism of vector spaces.
\item By remark \ref{233}, we have $\phi(\KK_0)=\tilde{\KK}_0$.
\item By corollary \ref{251}, we have $\phi\circ \sigma=\tilde{\sigma}\circ \phi$.
\item Let $a\in L_0^*$ and let $t\in X_a$ be as in remark \ref{307}. By proposition \ref{228}, we have
$$  \phi(t)=\phi_a(t)=\p_2(a,t)\ ,$$
hence
\begin{align*}
\tilde{q}\big(\p_1(a)\big)&\in \p_2(a,t)+\tilde{\KK}_0={\phi}(t)+\phi(\KK_0)= {\phi}\big(t+\KK_0\big)={\phi}\big(q(a)+\KK_0\big)={\phi}\big(q(a)\big)+\tilde{\KK}_0
\end{align*}
\item Given $(a,t)\in T$, we have\footnote{Although we obtained $\phi$ only by taking separable elements into account, $\phi$ of course extends the isomorphism associated with inseparable elements as well since $\p_1:L_0\to \tilde{L}_0$ is an extension of $\p_1:R_a\to \tilde{R}_{\p_1(a)}$.}
\begin{align*} \Phi(a,t)=\big(\p_1(a),{\phi}(t)\big)=\big(\p_1(a),\p_2(a,t)\big)=\gamma(a,t)\ .
\end{align*}
\end{itemize}
\end{bewzwei}

\chapter{Exceptonial Isomorphisms I}\label{269}

Now we drop the condition $\tilde{\AA}\cong \AA$, i.e., we assume $\tilde{\AA}\not \cong \AA$. By theorem \ref{178}, this can only occur in small dimensions. As a consequence, there aren't many possibilities for those exceptional isomorphisms which are, of course, not induced by an isomorphism of pseudo-quadratic spaces.

\begin{lemma}\label{275}
Suppose that $\tilde{\AA}\not\cong \AA$. Then exactly one of the following holds:
\begin{enumerate}[label=(\roman*)]
\item The involutory sets $(\AA,\FF,\sigma)$ and $(\tilde{\AA},\tilde{\FF},\tilde{\sigma})$ both are quadratic of type (iv) and we have
$$\dim_\AA L_0=2=\dim_{\tilde{\AA}} \tilde{L}_0\ .$$
\item The involutory sets $(\AA,\FF,\sigma)$ and $(\tilde{\AA},\tilde{\FF},\tilde{\sigma})$ are quadratic of type (iv) and (iii), respectively, and we we have
\begin{align*} \dim_\AA L_0=1\ , &&  \dim_{\tilde{\AA}} \tilde{L}_0=2\ .\end{align*}
\end{enumerate}
\end{lemma}

\begin{bew}
By theorem \ref{271} and remark \ref{308}, we have $\AA\cong \tilde{\AA}$ if $(\AA,\FF,\sigma)$ and $(\tilde{\AA},\tilde{\FF},\tilde{\sigma})$ both are quadratic of type (iii). Thus we may suppose $(\AA,\FF,\sigma)$ to be quadratic of type (iv). By theorem \ref{178}, we have $\dim_\AA L_0\leq 2$ and $\dim_{\tilde{\AA}} \tilde{L}_0\leq 2$. Notice that we have
$$\dim_\FF L_0=\dim_{\tilde{\FF}} \tilde{L}_0\ .$$
\begin{enumerate}[label=(\roman*)]
\item In the case $\dim_\AA L_0=1=\dim_{\tilde{\AA}} \tilde{L}_0$, we may apply proposition \ref{265} so that there is an isomorphism $\phi:\AA\to\tilde{\AA}$ of skew-fields$\qquad\lightning$.
\item As we have $\dim_\FF L_0\geq 4$ and $\dim_{\tilde{\FF}} \tilde{L}_0\leq 4$, the assertion follows immediately.
\end{enumerate}\qed
\end{bew}

\newpage

\begin{no}
Throughout this chapter, let $(\AA,\FF,\sigma)$ be quadratic of type (iv).
\end{no}

\begin{bem}
By lemma \ref{275}, there are two cases left. First of all, we deal with the case where both the involutory sets are quadratic of type (iv). We do this by using lemma \ref{266} in a suitable way. The appearing isomorphisms turn out to be almost isomorphisms of pseudo-quadratic spaces, but modified by switching the parametrizing space.
\end{bem}

\begin{no}\label{279}
Let $(\AA,\FF,\sigma)$ be quadratic of type (iv) and suppose that $\dim_\AA L_0=2$. By \cite{DM}, there are exactly three pseudo-quadratic spaces 
\begin{align*}
(\AA,\FF,\sigma,L_0,q)=(\AA_1,\FF,\sigma,L_0,q_1)=\Xi_1\ ,&& (\AA_2,\FF,\sigma,L_0,q_2)=\Xi_2\ , &&(\AA_3,\FF,\sigma,L_0,q_3)=\Xi_3
\end{align*}
with pairwise non-isomorphic quaternion division algebras $\AA_1,\AA_2,\AA_3$ which define the group $T$. When we switch between the parametrizing pseudo-quadratic spaces, we indicate this by the map
$$\id_T^i:T\to T,\ (a,t)\mapsto (a,t)\ ,$$
i.e., after applying $\id_T^i$, we consider $T$ to be defined by $\Xi_i$.
\end{no}

\begin{prop}\label{280}
Let $\tilde{\AA}$ be quadratic of type (iv) and suppose that $\dim_\AA L_0=2$. Then there are an $i\in \{1,2,3\}$ and an isomorphism $\Phi:\Xi \to \tilde{\Xi}_i$ of pseudo-quadratic spaces such that $\gamma$ is induced by $(\id_{\tilde{T}}^i)^{-1}\circ \Phi$.
\end{prop}

\begin{bew}
Let $a\in L_0^*$ be separable. By lemma $\ref{266}$, there is an $i\in \{1,2,3\}$ such that
$$\p_1(\langle a\rangle_\AA)=\langle \p_1(a)\rangle_{\tilde{\AA}_i}\ ,$$
and the map $(\p_1,\phi_a):(\langle a\rangle_\AA, \AA)\to (\langle \p_1(a)\rangle_{\tilde{\AA}_i},\tilde{\AA}_i)$ is an isomorphism of vector spaces. In particular, we have $\tilde{\AA}_i\cong \AA$, thus $i\in \{1,2,3\}$ is independent of the choice of $a$. Now the Jordan isomorphism
$$\id_{\tilde{T}}^i\circ \gamma: T\to{\tilde{T}}$$
is induced by an isomorphism $\Phi:\Xi\to \tilde{\Xi}_i$ of pseudo-quadratic spaces by theorem \ref{270}.\qed
\end{bew}

\begin{bem}
Now we consider the last case, where a quaternion division algebra ``splits'' into two separable quadratic extensions.
\end{bem}

\begin{lemma}\label{277}
Let $a\in L_0^*$ be separable and let $e\in \EE_a^\bot$. Suppose that $\p_1(ae)\notin \langle\p_1(a)\rangle_{\tilde{\AA}}$. Then we have
\begin{align*}
\tilde{f}\big(\p_1(ae),\p_1(a)\big)=0_{\tilde{\AA}}\ , && \phi_{ae}=\tilde{\sigma}\circ \phi_a\ .
\end{align*}
\end{lemma}

\begin{bew}
By equation \eqref{238}, we have
\begin{align*} \p_1(ae)\cdot \phi_{ae}(t)&=\p_1(ae)\cdot \p_2(a,t)^{\tilde{\sigma}}-\p_1(a)\cdot \p_2(a,t)^{-1}\tilde{f}\big(\p_1(a),\p_1(ae)\big)\p_2(a,t)^{\tilde{\sigma}}
\end{align*}
for each $t\in X_a$, hence
\begin{align*}
\p_2(a,t)^{-1}\tilde{f}\big(\p_1(a),\p_1(ae)\big)\p_2(a,t)^{\tilde{\sigma}}=0_{\tilde{\AA}}\ , && \forall\ t\in X_a:\ \phi_{ae}(t)=\p_2(a,t)^{\tilde{\sigma}}=\phi_a(t)^{\tilde{\sigma}}\ .
\end{align*}
\qed
\end{bew}

\begin{no}
Throughout the rest of this chapter, we suppose $\tilde{\AA}$ to be quadratic of type (iii).
\end{no}

\begin{bem}
By lemma \ref{275}, we have $\dim_\AA L_0=1$ and $\dim_{\tilde{\AA}} \tilde{L}_0=2$. As a consequence, each element $a\in L_0^*$ is separable. Moreover, we have $$\p_1(ae)\notin \tilde{R}_{\p_1(a)}=\langle \p_1(a)\rangle_{\tilde{\AA}}$$ for each element $e\in \EE_a^\bot$. As a consequence, lemma \ref{277} applies.
\end{bem}

\begin{no}\ 
\begin{itemize}
\item Throughout the rest of this chapter, let $a\in L_0^*$ (which is separable) and let $e\in \EE_a^{\bot}$.
\item We set
\begin{align*} \tilde{a}:=\p_1(a)\ , && \tilde{b}:=\p_1(ae)\in \tilde{a}^\bot\ , && \phi:=\phi_a\ .
\end{align*}
\end{itemize}
\end{no}

\begin{bem}
Given $x\in \AA$, we have
$$q(a\cdot x)\equiv x^\sigma q(a) x \mod \FF\ ,$$
cf. definition \ref{315}, hence $$\big(a\cdot x, x^\sigma q(a) x \big)\in T\ .$$
\end{bem}

\begin{lemma}\label{278}
Given $x=s+et\in \AA$, we have
\begin{align*}
\gamma\big( a\cdot x,x^\sigma q(a)x \big)=\big( \tilde{a}\cdot \phi(s)+\tilde{b}\cdot  \phi(t)^{\tilde{\sigma}}, \phi\big( N(x)q(a)\big)\big)\ .
\end{align*}
\end{lemma}

\begin{bewzwei}\
\begin{itemize}
\item We have
$$\p_1(a\cdot x)=\p_1(a\cdot s+ae\cdot t)=\p_1(a)\cdot \phi(s)+\p_1(ae)\cdot \phi(t)^{\tilde{\sigma}}=\tilde{a}\cdot \phi(s)+\tilde{b}\cdot  \phi(t)^{\tilde{\sigma}}$$
by lemma \ref{277}.
\item By Proposition \ref{228}, we have
$$\p_1\big(ax\cdot x^\sigma q(x)x\big)=\p_1(ax)\cdot \p_2(ax,x^\sigma q(a) x)\ .$$
On the other hand, we have
\begin{align*}
\p_1\big(a\cdot xx^\sigma q(a)x\big)&=\p_1\big(a\cdot N(x)q(a)(s+et)\big)\\
&=\p_1\big(as \cdot N(x)q(a)+aet \cdot N(x)q(a)^\sigma \big) \\
&=\p_1(as)\cdot \phi\big(N(x)q(a)\big)+\p_1(aet)\cdot \phi\big( N(x)q(a)^\sigma\big)^{\tilde{\sigma}} \\
&=\p_1(as+aet)\cdot \phi\big(N(x)q(a)\big)\\
&=\p_1(a x)\cdot \phi\big(N(x)q(a)\big)\ .
\end{align*}
\end{itemize}\qed
\end{bewzwei}

\begin{prop}\label{281}
Given $x=s+et\in \AA$ and $u\in \FF$, we have
\begin{align*}
\gamma( a\cdot x,x^\sigma q(a)x+u)= \big(\tilde{a}\cdot\phi(s)+\tilde{b}\cdot \phi(t)^{\tilde{\sigma}}, \phi\big(N(x)q(a)+u\big)\big)\ .
\end{align*}
\end{prop}

\begin{bew}
This results from lemma \ref{278} and lemma \ref{221}.\qed
\end{bew}

\begin{bem}
This describes $\gamma$ completely.
\end{bem}

\begin{satz}\label{284}
Suppose that $\tilde{\KK}\not\cong \KK$. Then exactly one of the following holds:
\begin{enumerate}[label=(\roman*)]
\item The involutory sets $(\KK,\KK_0,\sigma)$ and $(\tilde{\KK},\tilde{\KK}_0,\tilde{\sigma})$ both are quadratic of type (iv), we have $\dim_\KK L_0=2=\dim_{\tilde{\KK}} \tilde{L}_0$ and there are an $i\in \{2,3\}$ and an isomorphism  $\Phi:\Xi \to \tilde{\Xi}_i$ of pseudo-quadratic spaces such that $\gamma$ is induced by $(\id_{\tilde{T}}^i)^{-1}\circ \Phi$, where $\id_{\tilde{T}}^i$ and $\tilde{\Xi}=:\tilde{\Xi}_1,\tilde{\Xi}_2,\tilde{\Xi}_3$ are as in notation \ref{279}.
\item The involutory sets $(\KK,\KK_0,\sigma)$ and $(\tilde{\KK},\tilde{\KK}_0,\tilde{\sigma})$ are quadratic of type (iv) and (iii), respectively, we have $\dim_\KK L_0=1,\ \dim_{\tilde{\KK}} \tilde{L}_0=2$ and $\gamma$ can be described by
$$\forall\ x=s+et\in \KK,\ u\in \KK_0:\  \gamma\big( ax,x^\sigma q(a)x+u\big)= \big(\p_1(a)\phi(s)+\p_1(ae)\phi(t)^{\tilde{\sigma}}, \phi\big(N(x)q(a)+u\big)\big)\ ,$$
where $a\in L_0^*$ is arbitrary, $\phi=\phi_a$, $e\in \EE_a^\bot$ and $\p_1(ae)\in \p_1(a)^\bot$.
\end{enumerate}
\end{satz}

\begin{bew}
This results from lemma \ref{275}, proposition \ref{280}, lemma \ref{277} and proposition \ref{281}. Notice that we have $i\neq 1$ in (i) as we have $\KK\not\cong \tilde{\KK}$.\qed
\end{bew}

\chapter{The Reverse Direction}

Now we consider the reverse direction, i.e., we prove that each map as above is a Jordan isomorphism. 

\begin{bem} Notice that we don't assume $\gamma$ to be a Jordan isomorphism any longer.\end{bem}

\begin{satz}\label{288}
Let $\Xi$ and $\tilde{\Xi}$ be proper pseudo-quadratic spaces and let $\gamma:T\to\tilde{T}$ be a map that is induced by an isomorphism $\Phi=(\p,\phi):\Xi\to \tilde{\Xi}$ of pseudo-quadratic spaces. Then $\gamma$ is a Jordan isomorphism. 
\end{satz}

\begin{bewzwei}\ 
\begin{itemize}
\item We have
$$\gamma(0_{L_0},1_\KK)=\big(\p(0_{L_0}),\phi(1_\KK)\big)=(0_{\tilde{L}_0},1_{\tilde{\KK}})\ .$$
\item Given $(a,t),(b,v)\in T$, we have
\begin{align*}
\gamma\big((a,t)\cdot (b,v)\big)&=\gamma\big( a+b,t+v+f(b,a) \big)=\big(\p(a+b),\phi(t+v+f(b,a))\big)\\
&=\big(\p(a)+\p(b),\phi(t)+\phi(v)+\phi(f(a,b))\big) \\
&=\big(\p(a)+\p(b), \phi(t)+\phi(v)+\tilde{f}\big(\p(a),\p(b)\big) \big)\\
&=\big(\p(a),\phi(t)\big)\cdot \big(\p(b),\phi(v)\big)=\gamma\big((a,t)\big)\cdot \gamma\big((b,v)\big)\ .
\end{align*}
\item Let $(a,t),(b,v)\in T$. By corollary \ref{303}, we have
\begin{align*}
\gamma\big(h_{(a,t)}(b,v)\big)&=\gamma\big(b\cdot t^\sigma-a\cdot t^{-1}f(a,b)t^\sigma,tvt^\sigma\big) \\
&=\big(\p(b\cdot t^\sigma-a\cdot t^{-1}f(a,b)t^\sigma), \phi(tvt^\sigma) \big) \\
&=\big( \p(b)\cdot \phi(t^\sigma)-\p(a)\cdot  \phi(t)^{-1}\phi\big(f(a,b)\big)\phi(t^\sigma), \phi(t)\phi(v)\phi(t^\sigma)\big) \\
&=\big(\p(b)\cdot \phi(t)^{\tilde{\sigma}}-\p(a)\cdot \phi(t^{-1})\tilde{f}\big(\p(a),\p(b)\big)\phi(t)^{\tilde{\sigma}},\phi(t)\phi(v)\phi(t)^{\tilde{\sigma}}\big) \\
&=\tilde{h}_{(\p(a),\phi(t))}\big((\p(b),\phi(v))\big)=\tilde{h}_{\gamma(a,t)}\big(\gamma(b,v)\big)\ .
\end{align*}
\end{itemize}\qed
\end{bewzwei}

\begin{bem}
This shows that each map as in proposition \ref{280} is a Jordan isomorphism since $\id_T^i$ is a Jordan isomorphism for each $i\in\{1,2,3\}$.
\end{bem}

\begin{no}\ 
\begin{itemize}
\item Throughout the rest of this chapter, let $\Xi$ and $\tilde{\Xi}$ be proper pseudo-quadratic spaces, let $(\AA,\FF,\sigma)$ and $(\tilde{\AA},\tilde{\FF},\tilde{\sigma})$ be quadratic of type (iv) and (iii), respectively, let
\begin{align*}
\dim_\AA L_0=1\ , && \dim_{\tilde{\AA}} \tilde{L}_0=2\ ,
\end{align*}
let $a\in L_0^*$ and $e\in \EE_a^\bot$, and let $\phi:\EE_a\to \tilde{\AA}$ be an isomorphism of fields.
\item Moreover, let $\gamma:T\to \tilde{T}$ be a map such that
$$\forall\ x=s+et\in \AA,\ u\in \FF:\qquad \gamma( ax,x^\sigma q(a)x+u)=\big(\tilde{a}\phi(s)+\tilde{b}\phi(t)^{\tilde{\sigma}}, \phi\big(N(x)q(a)+u\big)\big)$$
for some $\tilde{a},\tilde{b}\in \tilde{L}_0$ such that $\tilde{f}(\tilde{a},\tilde{b})=0_{\tilde{\AA}}$.
\item We set $b:=ae$ and
\begin{align*}
A&:=\{ ( a\cdot x, t) \mid x\in \EE_a, t\in X_{ax} \}\leq T\ , & B&:=\{ ( b\cdot x, t) \mid x\in \EE_a, t\in X_{bx} \}\leq T\ , \\
\tilde{A}&:=\{ ( \tilde{a}\cdot x, t) \mid x\in \tilde{\AA}, t\in \tilde{X}_{\tilde{a}x} \}\leq \tilde{T}\ , & \tilde{B}&:=\{ ( \tilde{b}\cdot x, t) \mid x\in \tilde{\AA}, t\in \tilde{X}_{\tilde{b}x} \}\leq \tilde{T}\ .
\end{align*}
\end{itemize}
\end{no}

\begin{lemma}\label{285}
Given $x,y\in \EE_a$, we have
$$f(a\cdot x,b\cdot y)=f(b\cdot y,a\cdot x)\ .$$
\end{lemma}

\begin{bew}
Given $x,y\in \EE_a$, we have
\begin{align*}
f(a\cdot x,b\cdot y)&=f(a\cdot x,ae\cdot y)=x^\sigma f(a,a) ey=ey xf(a,a)^\sigma \\
&=-y^\sigma ef(a,a)x=(ey)^\sigma f(a,a) x=f(ae\cdot y,a\cdot x)=f(b\cdot y,a\cdot x)\ .
\end{align*}
\qed
\end{bew}

\begin{lemma}
Given $(a,t),(b,v)\in T$, we have
$$(a,t)\in C_T\big((b,v)\big)\ \Leftrightarrow\ f(a,b)=f(b,a)\ .$$
\end{lemma}

\begin{bew}
We have
\begin{align*}
(a,t)(b,v)=(b,v)(a,t)\ &\Leftrightarrow\ \big(a+b,t+v+f(b,a)\big)=\big(b+a,v+t+f(a,b)\big) \\
&\Leftrightarrow\ f(b,a)=f(a,b)\ .
\end{align*}
\qed
\end{bew}

\begin{kor}\label{286}
We have
\begin{align*}
A\subseteq C_T(B)\ , && \tilde{A}\subseteq C_{\tilde{T}}(\tilde{B})\ .
\end{align*}
\end{kor}

\begin{bew}
The first assertion results from lemma \ref{285}, and given $x,y\in \tilde{\AA}$, we have
$$f(\tilde{a}\cdot x,\tilde{b}\cdot y)=0_{\tilde{\AA}}=f(\tilde{b}\cdot y,\tilde{a}\cdot x)\ .$$\qed
\end{bew}

\begin{lemma}\label{287}
We have
\begin{align*}
T=AB\ , && \tilde{T}=\tilde{A}\tilde{B}\ .
\end{align*}
\end{lemma}

\begin{bew}
Let $(a\cdot x+b\cdot y,t)\in T$ where $x,y\in \EE_a$. Then we have
$$t\equiv q(a\cdot x+b\cdot y)\equiv q(a\cdot x)+q(b\cdot y)+f(a\cdot x,b\cdot y) \mod \FF\ ,$$
hence
$$\tilde{t}:=t-q(a\cdot x)-q(b\cdot y)-f(a\cdot x,b\cdot y)\in \FF\ .$$
Observe that we have $f(a\cdot x,b\cdot y)=f(b\cdot y,a\cdot x)$ by lemma \ref{285}, thus
\begin{align*}
(a\cdot x+b\cdot y,t)&=(a\cdot x+b\cdot y,q(a\cdot x)+q(b\cdot y)+f(a\cdot x,b\cdot y)+\tilde{t})\\
&=(a\cdot x+b\cdot y,q(a\cdot x)+q(b\cdot y)+f(b\cdot y,a\cdot x)+\tilde{t})\\
&=\big(a\cdot x,q(a\cdot x)\big)\cdot \big(b\cdot y,q(b\cdot y)+\tilde{t}\big)\in AB\ .
\end{align*}
The second assertion follows analogously.
\qed
\end{bew}

\begin{lemma}\label{289}
The maps
\begin{align*}
\Phi_1&: ( \EE_a,\FF,\sigma, R_a, q) \to ( \tilde{\AA},\tilde{\FF}, \tilde{\sigma}, \langle \tilde{a} \rangle_{\tilde{\AA}},\tilde{q}),\ (ax, t)  \mapsto \big(\tilde{a}\phi(x), \phi(t)\big)\ , \\
\Phi_2&: ( \EE_a,\FF,\sigma, R_b, q) \to ( \tilde{\AA},\tilde{\FF}, \tilde{\sigma}, \langle \tilde{b} \rangle_{\tilde{\AA}},\tilde{q}),\ (bx, t) \mapsto \big(\tilde{b}\phi(x)^{\tilde{\sigma}}, \phi(t)^{\tilde{\sigma}}\big)
\end{align*}
are isomorphisms of pseudo-quadratic spaces inducing $\gamma_A:=\gamma_{|A}:A\to \tilde{A}$ and $\gamma_B:=\gamma_{|B}:B\to \tilde{B}$, respectively.
\end{lemma}

\begin{bew} We consider $\Phi_2$, the first assertion follows analogously.
\begin{itemize}
\item Given $x=et\in \EE_a^\bot$ and $u\in \FF$, we have
\begin{align*}
\gamma_B\big( a\cdot x, x^\sigma q(a)x +u \big)&=\big(\tilde{b}\cdot \phi(t)^{\tilde{\sigma}}, \phi(N(x)q(a)+u)\big)\\
&=\big(\tilde{b}\cdot \phi(t)^{\tilde{\sigma}}, \phi\big((x^\sigma xq(a)^{\tilde{\sigma}}+u\big)^{\tilde{\sigma}}\big)\big) \\
&=\big(\tilde{b}\cdot \phi(t)^{\tilde{\sigma}}, \phi(x^\sigma q(a)x+u)^{\tilde{\sigma}}\big)=\Phi_2( a\cdot x,x^\sigma q(a)x+u)\ .
\end{align*}
\item Given $x\in \EE_a$, we have
$$\p_1(b\cdot x)=\tilde{b}\cdot \phi(x)^{\tilde{\sigma}}\ ,$$
thus $({\p_1}_{|R_b},\tilde{\sigma}\circ \phi): (R_b,\EE_a)\to (\langle \tilde{b}\rangle_{\tilde{\AA}},\tilde{\AA})$ is an isomorphism of vector spaces.
\item The map $\tilde{\sigma}\circ \phi:\EE_a\to \tilde{\AA}$ is clearly an isomorphism of involutory sets, cf. corollary \ref{251}.
\item Given $x=et\in \EE_a^\bot$, we have $\big(\p_1(a\cdot x), \tilde{q}(\p_1(a\cdot x))\big)\in \tilde{B}$ and $\big(\p_1(a\cdot x),\phi\big(N(x)q(a)\big)\big)\in \tilde{B}$, thus
\begin{align*}
\tilde{q}\big(\p_1(x)\big)\in \phi\big(N(x)q(a)\big)+\tilde{\FF}&=\phi\big(x^\sigma q(a)^\sigma x+\FF\big)=\phi\big((x^\sigma q(a) x)^\sigma+\FF\big)\\
&=\phi\big(q(a\cdot x)+\FF\big)^{\tilde{\sigma}}=\phi\big(q(a\cdot x)\big)^{\tilde{\sigma}}+\tilde{\FF}\ .
\end{align*}
\end{itemize}
\qed
\end{bew}

\begin{prop}\label{293}
The map $\gamma:T\to \tilde{T}$ is an isomorphism of groups.
\end{prop}

\begin{bew}
By corollary \ref{286}, it is enough to consider the following two cases:
\begin{itemize}
\item By lemma \ref{289} and theorem \ref{288}, $\gamma_A: A\to \tilde{A}$ and $\gamma_B:B\to \tilde{B}$ are Jordan isomorphisms. In particular, they are isomorphisms of groups.
\item Given $x\in \EE_a,\ y=et\in \EE_a^\bot$ and $u,v\in \FF$, we have
\begin{align*}
&\gamma\big((a\cdot x, x^\sigma q(a) x+u)\cdot(a\cdot y,y^\sigma q(a)y+v)\big)\\
&\qquad\qquad =\gamma\big(a\cdot (x+y), x^\sigma q(a) x+y^\sigma q(a)y+f(ax,ay)+u+v\big)  \\
&\qquad\qquad =\gamma\big(a\cdot (x+y), x^\sigma q(a) x+y^\sigma q(a)y+x^\sigma f(a,a)y +u+v\big)  \\
&\qquad\qquad =\gamma\big(a\cdot (x+y), x^\sigma q(a) x+y^\sigma q(a)y+ x^\sigma q(a)y- x^\sigma q(a)^\sigma y+u+v\big)\\
&\qquad\qquad =\gamma\big(a\cdot (x+y), x^\sigma q(a) x+y^\sigma q(a)y+ x^\sigma q(a)y+ y^\sigma q(a) x+u+v\big)\\
&\qquad\qquad =\gamma\big(a\cdot (x+y), (x+y)^\sigma q(a)(x+y)+u+v\big)\\
&\qquad\qquad=\big(\tilde{a}\cdot \phi(x)+\tilde{b}\cdot \phi(t)^{\tilde{\sigma}}, \phi(N(x+y)q(a)+u+v)\big) \\
&\qquad\qquad =\big(\tilde{a}\cdot \phi(x)+\tilde{b}\cdot \phi(t)^{\tilde{\sigma}},\phi(N(x)q(a)+N(y)q(a)+u+v)\big) \\
&\qquad\qquad =\big(\tilde{a}\cdot \phi(x),\phi(N(x)q(a)+u)\big)\cdot \big(\tilde{b}\cdot \phi(t)^{\tilde{\sigma}},\phi(N(y)q(a)+v)\big) \\
&\hspace*{5cm} \qquad =\gamma\big((a\cdot x, x^\sigma q(a)x+u)\big)\cdot \gamma\big((a\cdot y, y^\sigma a y+v)\big)\ .
\end{align*}
\end{itemize}
\qed
\end{bew}

\begin{lemma}\label{291}
Given $x\in \AA,\ t:=x^\sigma q(a)x+s\in X_{ax}$ and $\tilde{t}:=N(x)q(a)+s\in \EE_a$, we have
$$N(\tilde{t})=N(t)\ .$$
\end{lemma}

\begin{bew}
We have
\begin{align*}
N\big(\tilde{t}\big)-N(t)&=N\big(N(x)q(a)+s\big)-N\big(x^\sigma q(a)x+s\big)\\
&=N\big(N(x)q(a)\big)+N(s)+T\big(N(x)q(a)s\big)-N\big(x^\sigma q(a) x\big)-N(s)-T\big(x^\sigma q(a) xs\big) \\
&=N(x)^2N\big(q(a)\big)-N(x)N\big(q(a)\big)N(x)+s\big( T(N(x)q(a))-T(x^\sigma q(a) x)\big) \\
&=s\big( x^\sigma x (q(a)+q(a)^\sigma)- x^\sigma(q(a)+q(a)^\sigma)x \big)=sx^\sigma \big( xT(q(a))-T(q(a))x\big)=0_\AA\ .
\end{align*}
\qed
\end{bew}

\begin{lemma}\label{297}
Given $x\in \AA,\ t:=x^\sigma q(a) x+s\in X_{ax}$ and $\tilde{t}:=N(x)q(a)+s\in \EE_a$, we have
$$xt^{-1}t^{\sigma}x^{-1}=\tilde{t}^{-1}\tilde{t}^\sigma\in \EE_a\ .$$
\end{lemma}

\begin{bew}
By lemma \ref{291}, we have 
\begin{align*}
xt^{-1}t^{\sigma}x^{-1}=(xt^\sigma x^{-1})^2\cdot N(t)^{-1}=(\tilde{t}^\sigma)^2\cdot N(\tilde{t})^{-1}=\tilde{t}^{-1}\tilde{t}^\sigma\in \EE_a\ .
\end{align*}\qed
\end{bew}

\begin{kor}\label{290}
Given $x\in \AA$ and $t:=x^\sigma q(a) x+s\in X_{ax}$, we have
$$(x^{-\sigma}tt^{-\sigma}x^\sigma)q(a)(xt^{-1}t^\sigma x^{-1})=q(a)\ .$$
\end{kor}

\begin{bew}
First of all observe that we have
\begin{align*}
x^{-\sigma}tt^{-\sigma}x^\sigma=N(x)^{-1}x t^{-\sigma} t x^{-1}N(x)=(xt^{-1}t^\sigma x^{-1})^{-1}\ .
\end{align*}
Now we have $xt^{-1}t^{\sigma}x^{-1}\in \EE_a=C_\AA(\EE_a)$ by lemma \ref{297}.\qed
\end{bew}

\begin{lemma}\label{298}
Given $x,y\in \AA,\ t:=x^\sigma q(a)x+s\in X_{ax}$ and $u\in \FF$, we have
$$\gamma\big(h_{(ax,t)}\big(a\cdot y,y^\sigma q(a)y+u\big)\big)=\big(\p_1(a\cdot z),\phi\big(N(t)N(y)q(a)+N(t)u\big)\big)\ ,$$
where $z=xt^{-1}t^\sigma x^{-1}yt^\sigma$.
\end{lemma}

\begin{bew}
The first component of $h_{(ax,t)}\big(a\cdot y,y^\sigma q(a)y+u\big)$ is
\begin{align*}
ayt^\sigma-axt^{-1}f(ax,ay)t^\sigma&=ayt^\sigma-axt^{-1}f(ax,ax)x^{-1}yt^\sigma=ayt^\sigma-axt^{-1}(t-t^\sigma)x^{-1}yt^\sigma \\
&=ayt^\sigma-axt^{-1}tx^{-1}yt^\sigma+axt^{-1}t^\sigma x^{-1}yt^\sigma=axt^{-1}t^\sigma x^{-1}yt^\sigma\ ,
\end{align*}
and the second component of $h_{(ax,t)}\big(ay,y^\sigma q(a)y+u\big)$ is
\begin{align*}
t(y^\sigma q(a) y+u)t^\sigma=ty^\sigma q(a) yt^\sigma+N(t)u\ .
\end{align*}
By corollary \ref{290}, we have
\begin{align*}
z^\sigma q(a)z&=\big(xt^{-1}t^\sigma x^{-1}yt^\sigma\big)^\sigma q(a) \big(xt^{-1}t^\sigma x^{-1}yt^\sigma\big) \\
&=ty^\sigma (x^{-\sigma}tt^{-\sigma}x^\sigma)q(a)(xt^{-1}t^\sigma x^{-1})yt^\sigma=ty^\sigma q(a)yt^\sigma\ .
\end{align*}
Therefore, we have
\begin{align*}
\gamma\big(h_{(ax,t)}\big(a\cdot y,y^\sigma q(a)y+u\big)\big)&=\gamma\big(( a\cdot z,z^\sigma q(a)z+N(t)u)\big)=\big(\p_1(a\cdot z),\phi\big(N(z)q(a)+N(t)u\big)\big) \\
&=\big(\p_1(a\cdot z),\phi\big(N(t)N(y)q(a)+N(t)u\big)\big)
\end{align*}\qed
\end{bew}

\begin{prop}\label{294}
The map $\gamma$ preserves the second component of the Hua-maps.
\end{prop}

\begin{bew}
Let $x,y\in \AA,\ t:=x^\sigma q(a)x+s\in X_{ax},\ \tilde{t}:=N(x)q(a)+s\in \EE_a$ and $u\in \FF$. The second component of
\begin{align*}
\tilde{h}_{\gamma (ax,t)}\big(\gamma((a\cdot y,y^\sigma q(a)y+u))\big)&=\tilde{h}_{(\p_1(ax),\phi(\tilde{t}))}\big( \p_1(a\cdot y), \phi(N(y)q(a)+u)\big)
\end{align*}
is
\begin{align*}
\phi(\tilde{t})\phi\big(N(y)q(a)+u\big)\phi(\tilde{t})^{\tilde{\sigma}}&=\phi\big( N(\tilde{t})\cdot(N(y)q(a)+u)\big) \\
&=\phi\big(N(t)\cdot(N(y)q(a)+u)\big)=\phi\big( N(t)N(y)q(a)+N(t)u\big)
\end{align*}
by lemma \ref{291}. Now the assertion results from Lemma \ref{298}.
\qed
\end{bew}





\begin{prop}\label{296}
Let $x=\lambda+e\mu\in \AA,\ t:=x^{\sigma}q(a)x+s\in A_{ax}$ and $\tilde{t}:=N(x)q(a)+s\in \EE_a$. Then the following holds:
\begin{enumerate}[label=(\alph*)]
\item Given $(a\cdot y,u)\in A$, we have
$$\gamma\big(h_{(ax,t)}(a\cdot y,u)\big)=\tilde{h}_{\gamma(ax,t)}\big(\gamma(a\cdot y,u)\big)\ .$$
\item Given $(b\cdot y,u)\in B$, we have
$$\gamma\big(h_{(ax,t)}(b\cdot y,u)\big)=\tilde{h}_{\gamma(ax,t)}\big(\gamma(b\cdot y,u)\big)\ .$$
\end{enumerate}
\end{prop}

\newpage

\begin{bew}
By proposition \ref{294}, it remains to check the first component.
\begin{enumerate}[label=(\alph*)]
\item By lemma \ref{298} and lemma \ref{297}, the first component of $h_{(ax,t)}(a\cdot y,u)$
is
\begin{align*}
a\cdot \tilde{t}^{-1}\tilde{t}^\sigma y t^\sigma=a\cdot y\tilde{t}^{-1}\tilde{t}^\sigma\cdot \big( N(\lambda)q(a)^\sigma+N(e)N(\mu)q(a)+f(a\cdot \lambda,ae\cdot \mu)^\sigma+s\big)\ .
\end{align*}
Therefore, the first component of $\gamma\big(h_{(ax,t)}(a\cdot y,u)\big)$ is
\begin{align*}
\tilde{a}\cdot \phi\big(y\tilde{t}^{-1}\tilde{t}^\sigma\cdot \big( N(\lambda)q(a)^\sigma+N(e)N(\mu)q(a)+s\big)\big)+\tilde{b}\cdot \phi\big(y^\sigma\tilde{t}^{-\sigma}\tilde{t}f(a,a)\lambda\mu\big)^{\tilde{\sigma}}\ .
\end{align*}
On the other hand, the first component of $$\tilde{h}_{(\tilde{a}\phi(\lambda)+\tilde{b}\phi(\mu)^{\tilde{\sigma}},\phi(\tilde{t}))}\big( \tilde{a}\cdot \phi(y),\p_2(a\cdot y,u)\big)$$
is
\begin{align*}
&\tilde{a}\phi(y)\phi(\tilde{t})^{\tilde{\sigma}}-\big(\tilde{a}\phi(\lambda)+\tilde{b}\phi(\mu)^{\tilde{\sigma}}\big)\cdot\phi(\tilde{t})^{-1}\tilde{f}\big(\tilde{a}\phi(\lambda)+\tilde{b}\phi(\mu)^{\tilde{\sigma}},\tilde{a}\phi(y)\big)\phi(\tilde{t})^{\tilde{\sigma}} \\
&\quad= \tilde{a}\phi(y \tilde{t}^{-1}\tilde{t}^\sigma)\phi(\tilde{t})-\tilde{a}\phi(\lambda\tilde{t}^{-1}\lambda^{\sigma}\tilde{f}(\tilde{a},\tilde{a})y\tilde{t}^\sigma)
-\tilde{b}\phi\big(\mu^\sigma\tilde{t}^{-1}\lambda^\sigma \phi^{-1}\big(\tilde{f}(\tilde{a},\tilde{a})\big)y \tilde{t}^\sigma\big) \\
&\quad =\tilde{a}\phi(y \tilde{t}^{-1}\tilde{t}^\sigma)\phi(\tilde{t})-\tilde{a}\phi(y\tilde{t}^{-1}\tilde{t}^\sigma)\phi\big(N(\lambda)(\phi(q(a))-\phi(q(a))^\sigma)\big)
-\tilde{b}\phi\big(y\tilde{t}^{-1}\tilde{t}^\sigma f(a,a) \lambda^\sigma \mu^\sigma \big)  \\
&\quad =\tilde{a}\phi\big( y\tilde{t}^{-1}\tilde{t}^\sigma\cdot \big(N(x)q(a)+s-N(\lambda)q(a)+N(\lambda)q(a)^\sigma\big)\big)
-\tilde{b}\phi\big(y^\sigma\tilde{t}^{-\sigma}\tilde{t} f(a,a)^\sigma \lambda \mu\big)^{\tilde{\sigma}} \\
&\quad =\tilde{a}\phi\big( y\tilde{t}^{-1}\tilde{t}^\sigma\cdot\big(N(e)N(\mu)q(a)+N(\lambda)q(a)^\sigma+s\big)\big)
+\tilde{b}\phi\big(y^\sigma\tilde{t}^{-\sigma}\tilde{t} f(a,a) \lambda \mu\big)^{\tilde{\sigma}}\ .
\end{align*}
\item This follows analogously.
\end{enumerate}
\qed
\end{bew}

\begin{satz}\label{299}
Let $\Xi$ and $\tilde{\Xi}$ be proper pseudo-quadratic spaces, let $(\AA,\FF,\sigma)$ and $(\tilde{\AA},\tilde{\FF},\tilde{\sigma})$ be quadratic of type (iv) and (iii), respectively, let
\begin{align*}
\dim_\AA L_0=1\ , && \dim_{\tilde{\AA}} \tilde{L}_0=2\ ,
\end{align*}
let $a\in L_0^*$ and $e\in \EE_a^\bot$, and let $\phi:\EE_a\to \tilde{\AA}$ be an isomorphism of fields. Moreover, let $\gamma:T\to \tilde{T}$ be a map such that
$$\forall\ x=s+et\in \AA,\ u\in \FF:\qquad \gamma\big( a\cdot x,x^\sigma q(a)x+u\big)=\big(\tilde{a}\cdot \phi(s)+\tilde{b}\cdot \phi(t)^{\tilde{\sigma}}, \phi\big(N(x)q(a)+u\big)\big)$$
for some $\tilde{a},\tilde{b}\in \tilde{L}_0$ such that $\tilde{f}(\tilde{a},\tilde{b})=0_{\tilde{\AA}}$. Then $\gamma$ is a Jordan isomorphism.
\end{satz}

\begin{bewzwei}\ 
\begin{itemize}
\item Putting $x:=0_\AA$ and $u:=1_{\AA}$ yields
$$\gamma(0_{L_0},1_\AA)=\big(0_{\tilde{L}_0},\phi(1_\AA)\big)=(0_{\tilde{L}_0},1_{\tilde{\AA}})\ .$$
\item By proposition \ref{293}, $\gamma$ is an isomorphism of groups.
\item Let $x,y\in T$. By lemma \ref{287}, there are $\hat{a}\in A$ and $\hat{b}\in B$ such that $y=\hat{a}\cdot \hat{b}$. By lemma \ref{292} and proposition \ref{296}, we have
\begin{align*}
\gamma\big(h_{x}(y)\big)&=\gamma\big( h_{x}(\hat{a}\cdot \hat{b})\big)=\gamma\big( h_{x}(\hat{a})\cdot h_{x}(\hat{b})\big)=\gamma\big(h_{x}(\hat{a})\big)\cdot \gamma\big(h_{x}(\hat{b})\big)\\
&=\tilde{h}_{\gamma(x)}\big(\gamma(\hat{a})\big)\cdot \tilde{h}_{\gamma(x)}\big(\gamma(\hat{b})\big) = \tilde{h}_{\gamma(x)}\big(\gamma(\hat{a})\cdot\gamma(\hat{b})\big)= \tilde{h}_{\gamma(x)}\big(\gamma(\hat{a}\cdot \hat{b})\big)=\tilde{h}_{\gamma(x)}\big(\gamma(y)\big)
\end{align*}
\end{itemize}\qed
\end{bewzwei}

\chapter{Exceptional Isomorphisms II}
We show that a Moufang set $T$ defined by a one-dimensional pseudo-quadratic space over a quaternion division algebra could equally defined by a 2-dimensional pseudo-quadratic space over a separable quadratic extension and vice versa.

\begin{lemma}\label{319}
Let $(\AA, \FF,\sigma)$ be quadratic of type (iii) and assume $\dim_{L_0} \AA=2$. Let $f$ be a skew-hermitian form on $L_0$ and let $\{a,b\}$ be an orthogonal basis of $L_0$. Then $f$ is anisotropic if and only if we have $$f(b,b)f(a,a)^{-1}\notin -N(\AA)\ .$$
\end{lemma}

\begin{bew}
Given $s,t\in \AA$, we have
\begin{align*}
f(as+bt,as+bt)=0_\AA\ &\Leftrightarrow\ N(s)f(a,a)+N(t)f(b,b)=0_\AA \\
&\Leftrightarrow\ f(b,b)f(a,a)^{-1}=-N(st^{-1})\in -N(\AA)\ .
\end{align*}\qed
\end{bew}

\begin{lemma}
Let $\Xi$ be a pseudo-quadratic space such that dim $\dim_\AA L_0=1$ and such that $(\AA,\FF,\sigma)$ is quadratic of type (iv). Let $a\in L_0^*$, $e\in \EE_a^\bot$, $\beta:=-N(e)$ and $b:=ae$. Then the following holds:
\begin{enumerate}[label=(\alph*)]
\item  The set $\{a,b\}$ is an $\EE_a$-basis of $L_0$.
\item The skew-hermitian form $\tilde{f}:L_0\times L_0\to \EE_a$ defined by
\begin{align*}
\tilde{f}(a,a):=f(a,a)\ , && \tilde{f}(b,b):=N(e)f(a,a)\ , && \tilde{f}(a,b):=0_\AA
\end{align*}
is anisotropic.
\item There is a pseudo-quadratic form $\tilde{q}$ on $L_0$ with respect to $\FF$, $\sigma$ and $\tilde{f}$ such that $$\tilde{\Xi}:=(\EE_a,\FF,\sigma,L_0,\tilde{q})$$
is a pseudo-quadratic space, satisfying $\tilde{q}(a)\equiv q(a)\mod \FF$ and $\tilde{q}(b)\equiv N(e)q(a)\mod \FF$.
\item The map $\gamma: T\to \tilde{T}$ defined by
$$\forall\ x=s+et\in \AA,\ u\in \FF:\qquad \big(a\cdot x,x^\sigma q(a)x+u\big)\mapsto \big(a\cdot s+b\cdot t^{\sigma}, N(x)q(a)+u\big)$$
is a Jordan isomorphism.
\end{enumerate}
\end{lemma}

\begin{bewzwei}\ 
\begin{enumerate}[label=(\alph*)]
\item We have $L_0=a\cdot (\EE_a+e\EE_a)=a\cdot \EE_a+b\cdot \EE_a$.
\item Since we have $\AA\cong(\EE_a/\FF,\beta)$, we have $\tilde{f}(b,b)\tilde{f}(a,a)^{-1}=N(e)=-\beta\notin -N(\EE_a)$.
\item This results from (11.28) and (11.30) of \cite{TW}, respectively, and remark \ref{320}. Notice that $\tilde{f}$ is trace-valued with
\begin{align*}
\tilde{f}(a,a)=f(a,a)=q(a)+q(a)^\sigma\ , && \tilde{f}(b,b)=N(e)f(a,a)=N(e)q(a)+\big(N(e)q(a)\big)^\sigma
\end{align*} if we have $\Char \AA=2$, so that we choose $\tilde{\beta}_a:=q(a),\ \tilde{\beta}_b:=N(e)q(a)$ in this case.
\item By theorem \ref{299}, it suffices to show that $\gamma$ is well-defined. Given $s,t\in \EE_a$, we have
\begin{align*}
\tilde{q}(a\cdot s+b\cdot t^\sigma)&\equiv s^\sigma\tilde{q}(a)s+t^\sigma\tilde{q}(b)t+\tilde{f}(a\cdot s,b\cdot t^\sigma) \\
&\equiv N(s)q(a)+N(e)N(t)q(a)=N(x)q(a)\equiv N(x)q(a)+u \mod \FF\ . 
\end{align*}
\end{enumerate}\qed
\end{bewzwei}

\newpage

\begin{lemma}
Let $\Xi$ be a pseudo-quadratic space such that dim $\dim_\AA L_0=2$ and such that $(\AA,\FF,\sigma)$ is quadratic of type (iii). Let $\{ a,b\}$ be an orthogonal basis of $L_0$ (which exists by theorem (6.3) in chapter 7 of \cite{WSch}) and let $$\beta:=-f(b,b)f(a,a)^{-1}\ .$$ Then the following holds:
\begin{enumerate}[label=(\alph*)]
\item  We have $\beta\in \FF\sm N(\AA)$. In particular, $\tilde{\HH}:=(\AA/\FF,\beta)$ is a quaternion division algebra.
\item Let $e\in \AA^\bot$ such that $N(e)=-\beta$. We extend the scalar multiplication on $L_0$ to $\tilde{\HH}$ by $ae:=b$. Then the skew-hermitian form $\tilde{f}:L_0\times L_0\to \tilde{\HH}$ defined by $\tilde{f}(a,a):=f(a,a)$ is anisotropic, satisfying 
\begin{align*}
\tilde{f}(b,b)=-N(e)f(a,a)\ , && \tilde{f}(a,b)=-ef(a,b)\ .
\end{align*}
\item There is a pseudo-quadratic form $\tilde{q}$ on $L_0$ with respect to $\FF$, $\sigma$ and $\tilde{f}$ such that $$\tilde{\Xi}:=(\tilde{\HH},\FF,\sigma,L_0,\tilde{q})$$
is a pseudo-quadratic space, satisfying $q(a)\equiv\tilde{q}(a)\mod \FF$ and $q(b)\equiv N(e)\tilde{q}(a)\mod \FF$.
\item The map $\gamma: \tilde{T}\to T$ defined by
$$\forall\ x=s+et\in \tilde{\HH},\ u\in \FF:\qquad \big(a\cdot x,x^\sigma \tilde{q}(a)x+u\big)\mapsto \big(a\cdot s+b\cdot t^{\sigma}, N(x)\tilde{q}(a)+u\big)$$
is a Jordan isomorphism.
\end{enumerate}
\end{lemma}

\begin{bew}
 Notice that we have $\Fix(\sigma)=\FF$ since $(\AA,\FF,\sigma)$ is quadratic of type (iii).
\begin{enumerate}[label=(\alph*)]
\item We have 
\begin{align*} \beta^\sigma=-f(b,b)^\sigma f(a,a)^{-\sigma}=-f(b,b)f(a,a)^{-1}=\beta\ , &&\beta \in \Fix(\sigma)=\FF\ .
\end{align*}
By lemma \ref{319}, we have $\beta=-f(b,b)f(a,a)^{-1}\notin N(\AA)$.
\item We have
\begin{align*}
\tilde{f}(a,a)=f(a,a)_\AA\neq 0_{\AA}\ , && \tilde{f}(b,b)=\tilde{f}(ae,ae)=e^\sigma f(a,a)e=-N(e)f(a,a)\ .
\end{align*}
\item This results from (11.28) and (11.30) of \cite{TW}, respectively, and remark \ref{320}.
Notice that
\begin{align*}
q(b)+q(b)^\sigma=N(e)\big(q(a)+q(a)^\sigma\big)\ , && q(b)+N(e)q(a)=\big(q(b)+N(e)q(a)\big)^\sigma\in \Fix(\sigma)=\FF
\end{align*}
and that $\tilde{f}$ is trace-valued with $\tilde{f}(a,a)=f(a,a)=q(a)+q(a)^\sigma$ if we have $\Char \AA= 2$, so that we choose $\tilde{\beta}_a:=q(a)$ in this case to obtain
\begin{align*} q(b)&=f(b,b)+q(b)^\sigma=N(e)f(a,a)+q(b)^\sigma\\
&=N(e)q(a)+N(e)q(a)^\sigma+q(b)^\sigma\equiv N(e)q(a)\equiv N(e)\tilde{q}(a) \mod \FF\ .\end{align*}
\item By theorem \ref{299}, it suffices to show that $\gamma$ is well-defined. Given $s,t\in \AA$, we have
\begin{align*}
q(a\cdot s+b\cdot {t}^\sigma)&\equiv s^\sigma {q}(a)s+t^\sigma{q}(b)t+{f}(a\cdot s,b\cdot t^\sigma) \\
&\equiv N(s)\tilde{q}(a)+N(e)N(t)\tilde{q}(a)=N(x)\tilde{q}(a)\equiv N(x)\tilde{q}(a)+u \mod \FF\ .
\end{align*}
\end{enumerate}\qed
\end{bew}

\chapter[The Field \texorpdfstring{$\FF_4$}{F4}]{The Field \texorpdfstring{$\boldsymbol{\FF_4}$}{F4}}
We return to the case $\AA\cong \FF_4$ which we excluded in chapter \ref{314}.

\begin{no}
Throughout this chapter, let $a\in L_0^*$ and $\AA\cong \FF_4$, which implies $\tilde{\AA}\cong \FF_4$ and thus $\FF\cong \FF_2 \cong \tilde{\FF}$, and let $\gamma:T\to \tilde{T}$ be a Jordan isomorphism.
\end{no}

\begin{bem}
By remark \ref{308}, the map $(\p_1,\phi_a):(\langle a\rangle_\AA,\AA)\to (\langle \p_1(a)\rangle_{\tilde{\AA}},\tilde{\AA})$ is an isomorphism of vector spaces, where $\phi_a:\AA\to\tilde{\AA}$ is defined by
$$\forall\ t\in \AA:\qquad \p_1(a\cdot t)=\p_1(a)\cdot \phi_a(t)\ .$$
Moreover, we have $\tilde{X}_{\p_1(a)}=\{ \p_2(a,t),\p_2(a,t)^{\tilde{\sigma}}\}$ and thus either $\phi_a(t)=\p_2(a,t)$ for each $t\in X_a$ or $\phi_a(t)=\p_2(a,t)^{\tilde{\sigma}}$ for each $t\in X_a$.
\end{bem}

\begin{lemma}
Assume $\dim_\AA L_0\geq 2$. Then we have
$$\forall\ t\in X_a:\qquad \phi_a(t)=\p_2(a,t)\ .$$
\end{lemma}

\begin{bew} Notice that $(\AA,\FF,\sigma)$ is quadratic of type (iii), hence there are no inseparable elements. Let $b\in a^\bot$ and $t\in X_a$. By identity \eqref{169}, we have
\begin{align*}
\p_1(a)\cdot \phi_a(t^\sigma)&=\p_1(a\cdot t^\sigma)=\p_1(a\cdot t^\sigma-b\cdot t^{-1}f(a,b)t^\sigma)\\
&=\p_1(a)\cdot \p_2(a,t)^{\tilde{\sigma}}-\p_1(b)\cdot \p_2(a,t)^{-1}\tilde{f}\big(\p_1(a),\p_1(b)\big)\p_2(a,t)^{\tilde{\sigma}}\ .
\end{align*}
The linear independence of $\p_1(a)$ and $\p_1(b)$ yields $\phi_a(t)^{\tilde{\sigma}}=\phi_a(t^\sigma)=\p_2(a,t)^{\tilde{\sigma}}$.\qed
\end{bew}

\begin{bem}\label{309}
Now we can go on as in chapter \ref{314} and we obtain that $\gamma$ is induced by an isomorphism $\Phi:\Xi\to\tilde{\Xi}$ of pseudo-quadratic spaces. In the case $\dim_\AA L_0=1$ however, each isomorphism of groups turns out to be a Jordan isomorphism. At this point, we drop the assumption that $\gamma$ is a Jordan isomorphism.
\end{bem}

\begin{lemma}\label{357}
Assume $\dim_\AA L_0=1$. Then we have $h_a=\id_T$ for each $a\in L_0^*$.
\end{lemma}

\begin{bew}
Let $a\in L_0^*$. Given $x,y\in \AA$ and $s\in X_{ax},\ t\in X_{ay}$, we have
\begin{align*}
h_{(ax,s)}(a\cdot y,t)&=\big(a\cdot ys^\sigma-a\cdot x s^{-1}f(ax,ay)s^\sigma, sts^\sigma  \big)\\
&=\big(a\cdot ys^\sigma-a\cdot yN(x)f(a,a)s, N(s)t\big)=\big(a\cdot y(s^\sigma+s),t\big)=(a\cdot y,t)\ .
\end{align*}\qed
\end{bew}

\begin{kor}\label{311}
Assume $\dim_\AA L_0=1$ and let $\gamma:T\to \tilde{T}$ be an isomorphism of groups. Then $\gamma$ is a Jordan isomorphism.
\end{kor}

\begin{bew}
Because of $|Z(\tilde{T})|=2$, we have $\gamma(0_{L_0},1_\AA)=(0_{\tilde{L}_0},1_{\tilde{\AA}})$, and because of $|\tilde{T}|=|T|=8$, we have $\dim_{\tilde{\AA}} \tilde{L}_0=1$. By lemma \ref{357}, we have
\begin{align*}
\forall\ a\in L_0^*:&& h_a=\id_T\ , && \tilde{h}_{\gamma(a)}=\id_{\tilde{T}}
\end{align*}
and thus
$$\forall\ a\in L_0^*,\ x\in L_0:\qquad \gamma\big(h_a(x)\big)=\gamma(x)=\tilde{h}_{\gamma(a)}\big(\gamma(x)\big)\ .$$
\qed
\end{bew}

\newpage

\begin{bem}
As a consequence, it is convenient to determine the isomorphism class of the group $T$ and hence the structure of  $\Aut(T)$.
\end{bem}

\begin{lemma}
We have $T\cong Q_8$, where $Q_8$ denotes the quaternion group.
\end{lemma}

\begin{bew}
We have $|T|=8$, and given $(a,t)\in T\sm Z(T)$, we have $$(a,t)^2=\big(a+a,t+t+f(a,a)\big)=(0_{L_0},1_\AA)\neq (0_{L_0},0_\AA)\ .$$\qed
\end{bew}

\begin{bem}\label{313}
The outer automorphisms of $T$ are represented by the restrictions of isomorphisms of pseudo-quadratic spaces, i.e., by the group
$$\langle \gamma_\sigma, \rho_s \mid \gamma_\sigma: T\to \tilde{T},\ (a\cdot x,t)\mapsto (a\cdot x^\sigma, t^\sigma),\ \rho_s:T\to \tilde{T},\ (a\cdot x,t)\mapsto (a\cdot xs,t)\rangle\cong \Sigma_3\ ,$$
where $s\in \AA\sm \FF$. The three non-trivial inner automorphisms yield some exceptional Jordan automorphisms:
\begin{align*}
\gamma_a&: T\to T: &(a,t)&\mapsto (a,t)\ , && (as,t)\mapsto (as,t^\sigma)\ , &&(as^\sigma,t) \mapsto (as^\sigma,t^\sigma)\ , \\
\gamma_{as}&:T\to T: &(a,t)&\mapsto (a,t^\sigma)\ , &&(as,t)\mapsto (as,t)\ , &&(as^\sigma,t)\mapsto (as^\sigma, t^\sigma)\ , \\
\gamma_{as^\sigma}&:T\to T:&(a,t)&\mapsto (a,t^\sigma)\ , &&(as,t)\mapsto (as,t^\sigma)\ , &&(as^\sigma,t)\mapsto (as^\sigma, t)\ .
\end{align*}
\end{bem}

\begin{lemma}\label{310}
Suppose that $\dim_\AA L_0=1$ and let $\gamma:T\to \tilde{T}$ be a Jordan isomorphism. Then there are an isomorphism $\Phi:\Xi\to \tilde{\Xi}$ of pseudo-quadratic spaces and an inner automorphism $\tilde{\gamma}\in \Aut(\tilde{T})$ such that $\gamma$ is induced by $\tilde{\gamma}\circ \Phi$.
\end{lemma}

\begin{bew}
By remark \ref{308}, the map $$(\p_1,\phi_a):(\langle a\rangle_\AA,\AA)\to (\langle \p_1(a)\rangle_{\tilde{\AA}},\tilde{\AA})$$ is an isomorphism of vector spaces. Let $\Phi:\Xi\to\tilde{\Xi}$ be defined by
$$\Phi:=(\p_1,\phi_a): (L_0,\AA)\to (\tilde{L}_0,\tilde{\AA}),\ (a,t)\mapsto \big(\p_1(a),\phi_a(t)\big)\ .$$
Then the map $\tilde{\gamma}:=\gamma\circ \Phi_{|\tilde{T}}^{-1}:\tilde{T}\to \tilde{T}$ is an automorphism of $\tilde{T}$ such that $\tilde{\p}_1=\id_{\tilde{L}_0}$, hence $\tilde{\gamma}$ is inner.\qed
\end{bew}

\begin{satz}\label{312}
Assume $\KK\cong \FF_4$. A map $\gamma:T\to\tilde{T}$ is a Jordan isomorphism if and only if one of the following holds:
\begin{enumerate}[label=(\roman*)]
\item There is an isomorphism $\Phi:\Xi\to \tilde{\Xi}$ that induces $\gamma$.
\item We have $\dim_\KK L=1$ and there are an isomorphism $\Phi:\Xi\to \tilde{\Xi}$ of pseudo-quadratic spaces and a non-trivial inner automorphism $\tilde{\gamma}\in \Aut(\tilde{T})$ such that $\gamma$ is induced by $\tilde{\gamma}\circ \Phi$.
\end{enumerate}
\end{satz}

\begin{bewzwei}\ 
\begin{itemize}
\item[``$\Rightarrow$''] This results from remark \ref{309} and lemma \ref{310}.
\item[``$\Leftarrow$''] This result from theorem \ref{288} and corollary \ref{311}.
\end{itemize}\qed
\end{bewzwei}

\chapter{Conclusion}

The complete description of the Jordan isomorphisms between the Moufang sets of two proper pseudo-quadratic spaces is as follows:

\newglossaryentry{JPQ}{type=results,name={{Jordan Isomorphisms of Moufang Sets of Pseudo-Quadratic Form Type}},description={},sort=res}

\begin{satz}[\textbf{\gls{JPQ}}]\label{353}
Let $\Xi$ and $\tilde{\Xi}$ be proper pseudo-quadratic spaces. A map $\gamma:T\to\tilde{T}$ is a Jordan isomorphism if and only if one of the following holds:
\begin{enumerate}[label=(\roman*)]
\item There is an isomorphism $\Phi:\Xi\to \tilde{\Xi}$ of pseudo-quadratic spaces that induces $\gamma$.
\item The involutory sets $(\KK,\KK_0,\sigma)$ and $(\tilde{\KK},\tilde{\KK}_0,\tilde{\sigma})$ both are quadratic of type (iv), we have
\begin{align*}
\KK\not\cong \tilde{\KK}\ , && \dim_\KK L_0=2=\dim_{\tilde{\KK}} \tilde{L}_0
\end{align*}
and there are an $i\in \{2,3\}$ and an isomorphism  $\Phi:\Xi \to \tilde{\Xi}_i$ of pseudo-quadratic spaces such that $\gamma$ is induced by $(\id_{\tilde{T}}^i)^{-1}\circ \Phi$, where $\id_{\tilde{T}}^i$ and $\tilde{\Xi}=:\tilde{\Xi}_1,\tilde{\Xi}_2,\tilde{\Xi}_3$ are as in notation \ref{279}.
\item The involutory sets $(\KK,\KK_0,\sigma)$ and $(\tilde{\KK},\tilde{\KK}_0,\tilde{\sigma})$ are quadratic of type (iv) and (iii), respectively, we have $\dim_\KK L_0=1,\ \dim_{\tilde{\KK}} \tilde{L}_0=2$ and $\gamma$ can be described by
$$\forall\ x=s+et\in \KK,\ u\in \KK_0:\qquad \gamma( ax,x^\sigma q(a)x+u)= \big(\tilde{a}\phi(s)+\tilde{b}\phi(t)^{\tilde{\sigma}}, \phi\big(N(x)q(a)+u\big)\big)\ ,$$
where $a\in L_0^*$ is arbitrary, $\phi:\EE_a\to \tilde{\KK}$ is an isomorphism of fields, $e\in \EE_a^\bot$, $\tilde{a}\in \tilde{L}_0$ and $\tilde{b}\in \tilde{a}^\bot$.
\item We have $\KK\cong \FF_4\cong \tilde{\KK}$, $\dim_\KK L=1$ and there are an isomorphism $\Phi:\Xi\to \tilde{\Xi}$ of pseudo-quadratic spaces and a non-trivial inner automorphism $\tilde{\gamma}\in \Aut(\tilde{T})$ such that $\gamma$ is induced by $\tilde{\gamma}\circ \Phi$.
\end{enumerate}
\end{satz}

\begin{bewzwei}\ 
\begin{itemize}
\item[``$\Rightarrow$''] If (i) or (iv) holds, we have $\KK\cong \tilde{\KK}$, thus neither (ii) nor (iii) holds. If (iv) holds, (i) can't hold since $\tilde{\gamma}$ is not induced by an isomorphism of pseudo-quadratic spaces by remark \ref{313}. Suppose that neither (i) nor (iv) holds. By theorem \ref{283}, $(\KK,\KK_0,\sigma)$ and $(\tilde{\KK},\tilde{\KK}_0,\tilde{\sigma})$ are non-proper, hence they are of quadratic type by remark \ref{191}, and by theorem \ref{312}, we have $\KK\not\cong \FF_4$. But then we have $\KK\not\cong \tilde{\KK}$ by theorem \ref{283}.  Finally, either (ii) or (iii) holds by theorem \ref{284}.
\item[``$\Leftarrow$''] This results from theorem \ref{288}, theorem \ref{299} and theorem \ref{312}.
\end{itemize}\qed
\end{bewzwei}

\newpage

\addtocontents{toc}{\noindent\protect\mbox{}\protect\hrulefill\par}
\part{Simply Laced Foundations}\label{418}
\addtocontents{toc}{\noindent\protect\mbox{}\protect\hrulefill\par}

\noindent The subject of this part is the classification of simply laced twin buildings via foundations, which are amalgams of parametrized Moufang triangles. Given a simply laced twin building, we obtain a foundation by taking the set of rank 2 residues, which are Moufang triangles, and by parametrizing the corresponding root group sequences. These parametrizations make the glueings visible, and they turn out to be Jordan isomorphisms.

By looking at foundations of rank 3, we can deduce more information about the appearing glueings, e.g., the glueing of a foundation of type $A_3$ is an isomorphism of skew-fields, which is quite restrictive, e.g., concerning foundations involving octonions or residues of type $D_4$.\\

\chapter{Parametrized Moufang Triangles}

\newglossaryentry{tri}{type=foundations,name={\ensuremath{\TTT(\AA)}},description=parametrized standard Moufang triangle with respect to the alternative division ring $\AA$, sort=Moufang polygon}
\newglossaryentry{trio}{type=foundations,name={\ensuremath{\TTT^o(\AA)}},description=parametrized opposite Moufang triangle with respect to the alternative division ring $\AA$, sort=Moufang polygon}

\begin{de} Let $\AA$ be an alternative division ring. 
\begin{itemize}
\item The root group sequence
$$\gls{tri}:=\big( U_{[1,3]},x_1(\AA),x_2(\AA),x_3(\AA)\big)$$
with commutator relations
\begin{align*}
\forall\ s,t\in \AA:\qquad [x_1(s),x_3(t)]&:=x_2(st)
\end{align*}
is the \textit{parametrized standard triangle with respect to $\mathit{\AA}$}\index{Moufang triangle!parametrized standard}.
\item The root group sequence
$$\gls{trio}:=\big(U_{[1,3]},x_1(\AA),x_2(\AA),x_3(\AA)\big)$$
with commutator relations
\begin{align*}
\forall\ s,t\in \AA:\qquad [x_1(s),x_3(t)]&:=x_2(-st)
\end{align*}
is the \textit{parametrized opposite triangle with respect to $\mathit{\AA}$}\index{Moufang triangle!parametrized opposite}.
\end{itemize}
\end{de}

\begin{bem}
For reasons of brevity, we will write $$\TTT^{(o)}(\AA)=\big( x_1(\AA),\ldots, x_3(\AA)\big)$$ instead of 
$\TTT^{(o)}(\AA)=\big( U_{[1,3]}, x_1(\AA),\ldots, x_3(\AA)\big)$.
\end{bem}

\begin{lemma}\label{360}
Given an alternative division ring $\AA$, we have
$$\TTT^o(\AA)\cong \TTT(\AA)\ .$$
\end{lemma}

\begin{bew}
Let 
\begin{align*}
\TTT^o(\AA)=\big(x_1(\AA),x_2(\AA), x_3(\AA)\big)\ , && \TTT(\AA)=\big(\tilde{x}_1(\AA), \tilde{x}_2(\AA), \tilde{x}_3(\AA)\big)\ .
\end{align*}
Then $\alpha=(\alpha^1,\alpha^2,\alpha^3)$ with
\begin{align*}
\alpha_1:x_1(\AA)\to\tilde{x}_1(\AA),\ &x_1(t)\mapsto \tilde{x}_1(t)\ , \\
\alpha_2:x_2(\AA)\to\tilde{x}_2(\AA),\ &x_2(t)\mapsto \tilde{x}_2(-t)\ , \\
\alpha_3:x_3(\AA)\to\tilde{x}_3(\AA),\ &x_3(t)\mapsto \tilde{x}_3(t)
\end{align*}
preserves the commutator relations: Given $s,t\in \AA$, we have
$$\alpha\big([x_1(s),x_3(t)]\big)=\alpha\big(x_2(-st)\big)=\tilde{x}_2(st)=[\tilde{x}_1(s),\tilde{x}_3(t)]=[\alpha\big(x_1(s)\big),\alpha\big(x_3(t)\big)]\ .$$\qed
\end{bew}

\begin{lemma}\label{359}
Let $\TTT(\AA)=\big(x_1(\AA), x_2(\AA), x_3(\AA)\big)$ be a parametrized standard triangle. Then we have
$$\big(x_3(\AA^o), x_2(\AA^o), x_1(\AA^o)\big)=\TTT^o(\AA^o)\ .$$
\end{lemma}

\newpage 

\begin{bew}
Given $s,t\in \AA^o$, we have
\begin{align*}
[x_3(s),x_1(t)]=[x_1(t),x_3(s)]^{-1}=x_2(ts)^{-1}=x_2(-s\circ t)\ .
\end{align*}
\qed
\end{bew}

\begin{no} In the following, a \textit{(parametrized) Moufang triangle}\index{Moufang triangle} always denotes a parametrized standard Moufang triangle.
\end{no}

\begin{de} Let $\TTT(\AA)$, $\TTT(\tilde{\AA})$ be parametrized Moufang triangles. 
\begin{itemize}
\item An \textit{isomorphism}\index{isomorphism!of parametrized Moufang triangles} $\alpha:\TTT(\AA)\to \TTT(\tilde{\AA})$ is a triple $(\alpha_1,\alpha_2,\alpha_3)$ such that $\alpha_1,\alpha_2,\alpha_3:\AA\to \tilde{\AA}$ are isomorphisms of additive groups satisfying
$$\forall\ s,t\in\AA:\qquad \alpha_2(st)=\alpha_1(s)\alpha_3(t)\ .$$
\item A \textit{reparametrization for $\mathit{\TTT(\AA)}$}\index{reparametrization!for a Moufang triangle} is an ordered set $\alpha=(\tilde{\AA},\alpha_1,\alpha_2,\alpha_3)$ such that $\tilde{\AA}$ is an alternative division ring and $\alpha_1,\alpha_2,\alpha_3:\tilde{\AA}\to \AA$ are isomorphisms of additive groups satisfying
$$\forall\ s,t\in\tilde{\AA}:\qquad \alpha_2(st)=\alpha_1(s)\alpha_3(t)\ .$$
\end{itemize}
\end{de}

\begin{lemma}
Let $\TTT(\AA)$ be a parametrized Moufang triangle and let $\alpha=(\tilde{\AA},\alpha_1,\alpha_2,\alpha_3)$ be a reparametrization for $\TTT(\AA)$. Then we have
\begin{align*} \big(\tilde{x}_1(\tilde{\AA}),\tilde{x}_2(\tilde{\AA}),\tilde{x}_3(\tilde{\AA})\big)=\TTT(\tilde{\AA})\ ,&& \forall\ i=1,2,3:\ \tilde{x}_i:=x_i\circ \alpha_i\ .\end{align*}
\end{lemma}

\begin{bew}
Given $s,t\in \tilde{\AA}$, we have
$$[\tilde{x}_1(s),\tilde{x}_3(t)]=[x_1\big(\alpha_1(s)\big),x_3\big(\alpha_3(t)\big)]=x_2\big(\alpha_1(s)\alpha_3(t)\big)=x_2\big(\alpha_2(st)\big)=\tilde{x}_2(st)\ .$$
\qed
\end{bew}

\begin{lemma}\label{362}
Let $\TTT(\AA)$ be a parametrized Moufang triangle and let $a,b\in\AA^*$. Then there are reparametrizations $\alpha=(\AA,\alpha_1,\alpha_2,\alpha_3)$ and $\beta=(\AA,\beta_1,\beta_2,\beta_3)$ for $\TTT(\AA)$ such that
\begin{align*}
x_1\big(\alpha_1(1_\AA)\big)=x_1(1_\AA)\ , && x_3\big(\alpha_3(1_\AA)\big)=x_3(a)\ , && x_1\big(\beta_1(1_\AA)\big)=x_1(b)\ , && x_3\big(\beta_3(1_\AA)\big)=x_3(1_\AA)\ .
\end{align*}
\end{lemma}

\begin{bew}
Set $\alpha:=(\AA,\phi,\rho_a\phi,\rho_a \phi)$ with $\phi$ as in (20.25) of \cite{TW}. For the second statement, apply the first result to $\TTT(\AA^o)$.\qed
\end{bew}

\begin{lemma}\label{361}
Let $\TTT(\AA)=\big(x_1(\AA), x_2(\AA), x_3(\AA)\big)$ be a parametrized Moufang triangle. Then the action of the \textit{Hua automorphism} $h_1(s):=\mu\big(x_1(1_\AA)\big)^{-1}\mu\big(x_1(s)\big)$ on $x_1(\AA)\times x_3(\AA)$ corresponds to the map
$$(t,u)\mapsto (sts, s^{-1}u)\ ,$$
and the action of the \textit{Hua automorphism}\index{Hua automorphism! of a Moufang triangle} $h_3(s):=\mu\big(x_3(1_\AA)\big)^{-1}\mu\big(x_3(s)\big)$ on $x_1(\AA)\times x_3(\AA)$ corresponds to the map
$$(t,u)\mapsto (ts^{-1}, sus )\ .$$
\end{lemma}

\begin{bew}
This is (33.10) of \cite{TW}.
\qed
\end{bew}

\newpage

\chapter{Foundations}\label{491}

In this chapter, we introduce the objects we will mainly deal with and which turn out to be a classifying invariant for simply laced twin buildings.

\section{Definition}

The definition given here combines ideas and concepts of B. Mühlherr and R. Weiss. In this part, we only consider simply laced foundations, cf. part \ref{492} for a general definition.

\begin{de}\ 
\begin{itemize}
\item Let $M$ be a \textit{simply laced} Coxeter matrix\index{simply laced}, i.e., we have $m_{ij}\in\{2,3\}$ for all $i,j\in I$. A \textit{foundation of type $\mathit{M}$}\index{foundation} is a set
$$\FFF:=\{\TTT(\AA_{(i,j)}) , \gamma_{(i,j,k)} \mid (i,j)\in A(M),(i,j,k)\in G(M) \}$$
such that:
\begin{enumerate}[label=(F\arabic*),leftmargin=25pt]
\item Given $(i,j)\in A(M)$, then $\TTT(\AA_{(i,j)})$ is a Moufang triangle over $\AA_{(i,j)}$.
\item Given $(i,j)\in A(M)$, we have $\AA_{(i,j)}=\AA_{(j,i)}^o$.
\item Given $(i,j,k)\in G(M)$, then $\gamma_{(i,j,k)}:\AA_{(i,j)}\to \AA_{(j,k)}$ is an isomorphism of additive groups satisfying
\begin{align*} \gamma_{(i,j,k)}(1)=1\ , && \gamma_{(i,j,k)}=\id^o\circ \gamma_{(k,j,i)}^{-1}\circ \id^o\ .\end{align*}
\item Given $(i,j,k),(i,j,l),(l,j,k)\in G(M)$, we have
$$\gamma_{(i,j,k)}=\gamma_{(l,j,k)}\circ \id^o\circ \gamma_{(i,j,l)}\ .$$
\end{enumerate}
\item Given a foundation $\FFF$, we denote the corresponding Coxeter Matrix by $F$.
\item A foundation $\FFF$ is a \textit{Moufang foundation}\index{Moufang foundation}\index{foundation!Moufang} if each \textit{glueing} $\gamma:=\gamma_{(i,j,k)}$ is a Jordan isomorphism, i.e., we have
\begin{align*}
\forall\ s,t\in\AA_{(i,j)}:\qquad\gamma(sts)=\gamma(s)\gamma(t)\gamma(s)\ .
\end{align*}
\end{itemize}
\end{de}

\begin{de}  \newglossaryentry{Fres}{type=symbols,name={\ensuremath{\FFF_J}},description=$J$-residue of the foundation $\FFF$, sort=foundations}
Let $\FFF$ be a foundation over $I=V(F)$ and let $J\subseteq I$. The \textit{$\mathit{J}$-residue of $\mathit{\FFF}$}\index{$J$-residue!of a foundation}\index{residue!of a foundation} is the foundation
$$\gls{Fres}:=\{ \TTT(\AA_{(i,j)}),\gamma_{(i,j,k)} \mid (i,j)\in J^2\cap A(F), (i,j,k)\in J^3\cap G(F)\}\ .$$
\end{de}

\begin{bem}
Since a foundation is, in fact, an amalgam of Moufang triangles, an isomorphism of foundations is a system of isomorphism of Moufang triangles preserving the glueings.
\end{bem}

\begin{de}
Let $\FFF,\tilde{\FFF}$ be foundations.
\begin{itemize}
\item An \textit{isomorphism}\index{isomorphism!of foundations} $\alpha:\FFF\to \tilde{\FFF}$ is a system $\alpha=\{\pi, \alpha_{(i,j)} \mid (i,j)\in A(F)\}$ of isomorphisms
\begin{align*}
\pi:F\to \tilde{F}\ , && \alpha_{(i,j)}=(\alpha_{(i,j)}^i,\alpha_{(i,j)}^{ij},\alpha_{(i,j)}^j): \TTT(\AA_{(i,j)})\to \TTT(\tilde{\AA}_{(\pi(i),\pi(j))})
\end{align*}
such that
\begin{align*}
\forall\ (i,j,k)\in G(F): \qquad \tilde{\gamma}_{(\pi(i),\pi(j),\pi(k))}\circ \alpha_{(i,j)}^j=\alpha_{(j,k)}^j\circ \gamma_{(i,j,k)}
\end{align*}
and $\alpha_{(i,j)}=\alpha_{(j,i)}^o$ for each $(i,j)\in A(F)$.
\item An isomorphism $\alpha:\FFF\to\tilde{\FFF}$ is \textit{special}\index{isomorphism!special} if $F=\tilde{F}$ and $\pi=\id_F$.
\item An \textit{automorphism of $\mathit{\FFF}$}\index{automorphism!of a foundation} is an isomorphism $\alpha:\FFF\to\FFF$.
\end{itemize}
\end{de}

\newpage

\section{Visualizing Foundations}

Given a foundation $\FFF$ of type $M$, we can extend the corresponding Coxeter diagram $\Pi_M$ in such a way that it contains all the information of the given foundation $\FFF$:
\begin{itemize}
\item Given an edge $\{i,j\}\in E(M)$, we label it by either $\TTT(\AA_{(i,j)})$ or $\TTT(\AA_{(j,i)})$ and add an arrow to indicate in which direction we have the given standard root group sequence.
\item Given $(i,j,k)\in G(M)$, we choose either $(i,j,k)$ or $(k,j,i)$, we add a directed arc from $\{i,j\}$ to $\{j,k\}$ or vice versa, and label it by $\gamma_{(i,j,k)}$, resp. $\gamma_{(k,j,i)}$.
\end{itemize}
The remaining information can be deduced from the given ones. Notice that the constructed diagram is not uniquely determined by $\FFF$ as there is a choice in the directions.

\begin{bsp}
An arbitrary foundation of type $A_3$ is given by
\begin{center}\begin{tikzpicture}[scale=0.7,>=stealth,thick]
\begin{scope}
\coordinate (1) at (0,0);                    
\coordinate (2) at (3,0);
\coordinate (3) at (6,0);
\draws{0.55} (1)--(2);
\draws{0.55} (2)--(3);
\node () at (0,-0.5) {\ssi$1$};
\node () at (3,-0.5) {\ssi$2$};
\node () at (6,-0.5) {\ssi$3$};
\node () at (3,0.4) {\ssi$\gamma_{(1,2,3)}$};
\node () at (1.5,-0.5) {\ssi$\TTT(\AA_{(1,2)})$};
\node () at (4.5,-0.5) {\ssi$\TTT(\AA_{(2,3)})$};
\node () at (8,0) {,};
\node () at (-2,0) {};
\node () at (3,1.5) {};
\draw[->] (2,0) arc[radius=1,start angle=180, end angle=0];
\foreach \i in {1,...,3} {\fill (\i) circle (3pt);}
\end{scope}
\end{tikzpicture}\end{center}
a concrete example is
\begin{center}\begin{tikzpicture}[scale=0.7,>=stealth,thick]
\begin{scope}
\coordinate (1) at (0,0);                    
\coordinate (2) at (3,0);
\coordinate (3) at (6,0);
\draws{0.55} (1)--(2);
\draws{0.55} (2)--(3);
\node () at (0,-0.5) {\ssi$1$};
\node () at (3,-0.5) {\ssi$2$};
\node () at (6,-0.5) {\ssi$3$};
\node () at (3,0.4) {\ssi$\id_\AA^o$};
\node () at (1.5,-0.5) {\ssi$\TTT(\AA)$};
\node () at (4.5,-0.5) {\ssi$\TTT(\AA^o)$};
\node () at (8,0) {,};
\node () at (-2,0) {};
\draw[->] (2,0) arc[radius=1,start angle=180, end angle=0];
\foreach \i in {1,...,3} {\fill (\i) circle (3pt);}
\end{scope}
\end{tikzpicture}\end{center}
where $\AA$ is an arbitrary alternative division ring. We will see that an integrable foundation of type $A_3$ is isomorphic to the foundation
\begin{center}\begin{tikzpicture}[scale=0.7,>=stealth,thick]
\begin{scope}
\coordinate (1) at (0,0);                    
\coordinate (2) at (3,0);
\coordinate (3) at (6,0);
\draws{0.55} (1)--(2);
\draws{0.55} (2)--(3);
\node () at (0,-0.5) {\ssi$1$};
\node () at (3,-0.5) {\ssi$2$};
\node () at (6,-0.5) {\ssi$3$};
\node () at (3,0.4) {\ssi$\id_\DD$};
\node () at (1.5,-0.5) {\ssi$\TTT(\DD)$};
\node () at (4.5,-0.5) {\ssi$\TTT(\DD)$};
\node () at (8,0) {};
\node () at (-2,0) {};
\node () at (-0.5,1) {$\AAA_3(\DD)$};
\draw[->] (2,0) arc[radius=1,start angle=180, end angle=0];
\foreach \i in {1,...,3} {\fill (\i) circle (3pt);}
\end{scope}
\end{tikzpicture}\end{center}
for some skew-field $\DD$, i.e., the previous example is not integrable.
\end{bsp}

\begin{bem}\ 
\begin{enumerate}[label=(\alph*)]
\item Concerning a given problem, we sometimes don't need the whole visualization to get a feeling for the crucial step in the solution. In this case, we restrict to a diagram with the essential information, e.g., we just indicate whether some glueings are iso- or anti-isomorphisms of skew-fields.
\item Notice that in the above construction, the resulting diagram possibly carries redundant information: Given $(i,j,k),(i,j,l)\in G(M)$ (and thus $(l,j,k)\in G(M)$), we have 
$$\gamma_{(i,j,k)}=\gamma_{(l,j,k)}\circ \id^o\circ \gamma_{(i,j,l)}\ ,$$
which means that
\begin{center}\begin{tikzpicture}[scale=0.7,>=stealth,thick]
\begin{scope}[xshift=-4cm]
\coordinate (1) at (0,0);                    
\coordinate (2) at (-30:3);
\coordinate (3) at (-150:3);
\coordinate (4) at (90:3);
\draw (1)--(2);
\draw (1)--(3);
\draw (1)--(4);
\node () at (270:0.6) {\ssi$\gamma_{1}$};
\node () at (30:0.6) {\ssi$\gamma_{2}$};
\node () at (150:0.6) {\ssi$\gamma_{3}$};
\draw[->] (-150:1) arc[radius=1,start angle=-150, end angle=-30];
\draw[->] (-30:1) arc[radius=1,start angle=-30, end angle=90];
\draw[->] (210:1) arc[radius=1,start angle=210, end angle=90];
\foreach \i in {1,...,4} {\fill (\i) circle (3pt);}
\end{scope}
\node () at (0,0) {and};
\begin{scope}[xshift=4cm]
\coordinate (1) at (0,0);                    
\coordinate (2) at (-30:3);
\coordinate (3) at (-150:3);
\coordinate (4) at (90:3);
\draw (1)--(2);
\draw (1)--(3);
\draw (1)--(4);
\node () at (270:0.6) {\ssi$\gamma_{1}$};
\node () at (30:0.6) {\ssi$\gamma_{2}$};
\draw[->] (-150:1) arc[radius=1,start angle=-150, end angle=-30];
\draw[->] (-30:1) arc[radius=1,start angle=-30, end angle=90];
\foreach \i in {1,...,4} {\fill (\i) circle (3pt);}
\end{scope}
\end{tikzpicture}\end{center}
carry the same information, where $\gamma_1=\gamma_{(i,j,l)}$, $\gamma_2=\gamma_{(l,j,k)}$ and $\gamma_3=\gamma_{(i,j,k)}$.
\end{enumerate}
\end{bem}

\section{Root Group Systems}

The fact that a root group systems is a classifying invariant of the corresponding twin building is a fundamental result in twin building theory.
\newglossaryentry{rgs}{type=symbols,name={\ensuremath{U_{(i,j)}}},description=root group sequence from $\alpha_i$ to $\alpha_j$, sort=foundation}
\newglossaryentry{rgsys}{type=symbols,name={\ensuremath{\UUU(\BBB,M,\Sigma,c)}},description=root group system of $\BBB$ based at ${(\Sigma,c)}$, sort=foundation}

\begin{de}
Let $\BBB$ be a simply laced twin building of type $M$, let $\Sigma$ be a twin apartment of $\BBB$ and let $c\in \OOO_\Sigma$.  
\begin{itemize}
\item Given $(i,j)\in A(M)$, let $\alpha_i,\alpha_j$ be the simple roots with respect to $(\Sigma,c)$ and let $\Theta_{(i,j)}$ be as in theorem \ref{499}. Then
$$\gls{rgs}:=( U_{[i,j]}, U_{(i,j)}^i, U_{(i,j)}^{ij}, U_{(i,j)}^j ):=\Theta_{(i,j)}$$
denotes the root group sequence of $\BBB$ from $\alpha_i$ to $\alpha_j$, which is isomorphic to the root group sequence of $\BBB_{ij}$ from $\alpha_i \cap \BBB_{ij}$ to $\alpha_j\cap \BBB_{ij}$.
\item The resulting set
$$\gls{rgsys}:=\{ U_{(i,j)} \mid (i,j)\in A(M)\}$$
is the \textit{root group system of $\mathit{\BBB}$ based at $\mathit{(\Sigma,c)}$}\index{root group system}.
\end{itemize}
\end{de}

\begin{lemma}\label{59}
Given $(i,j,k)\in G(F)$, we have $U_{(i,j)}^j=U_{(j,k)}^j$.
\end{lemma}

\begin{bew}
This holds by definition.\qed
\end{bew}

\begin{de}
Let $\UUU:=\UUU(\BBB,M,\Sigma,c)$ and $\tilde{\UUU}:=\UUU(\tilde{\BBB},\tilde{M},\tilde{\Sigma},\tilde{c})$ be root group systems.\begin{itemize}
\item An \textit{isomorphism}\index{isomorphism!of root group systems} $\alpha:\UUU\to\tilde{\UUU}$ is a system
$$\alpha=\{ \pi,\alpha_{(i,j)} \mid (i,j)\in A(M)\} $$
of isomorphisms
\begin{align*}
\pi:M\to \tilde{M}\ , && \alpha_{(i,j)}:U_{(i,j)}\to \tilde{U}_{(\pi(i),\pi(j))}
\end{align*}
such that
\begin{align*}
\forall\ (i,j,k)\in G(M):\ {\alpha_{(i,j)}}_{|U_{(i,j)}^j}={\alpha_{(j,k)}}_{|U_{(j,k)}^j}\ , && \forall\ (i,j)\in A(M):\ \alpha_{(i,j)}=\alpha_{(j,i)}^o\ .
\end{align*}
\item An isomorphism $\alpha:\UUU\to\tilde{\UUU}$ is \textit{special}\index{isomorphism!special} if $M=\tilde{M}$ and $\pi=\id_M$.
\item An \textit{automorphism of $\mathit{\UUU}$}\index{automorphism!of a root group system} is an isomorphism $\alpha:\UUU\to\UUU$.
\end{itemize}
\end{de}

\begin{satz}\label{64} 
Two root group systems $\UUU(\BBB,M,\Sigma,c)$ and $\UUU(\BBB,M,\tilde{\Sigma},\tilde{c})$ of of a twin building $\BBB$ are specially isomorphic.
\end{satz}

\begin{bew} This is a consequence of theorem \ref{497}. \qed
\end{bew}

\begin{satz}\label{137}
Let $\UUU:=\UUU(\BBB,M,\Sigma,c)$ be a root group system of a twin building $\BBB$. Then the isomorphism class of $\UUU$ is a classifying invariant of the isomorphism class of $\BBB$.
\end{satz}

\begin{bew}
This is a consequence of the extension theorem \ref{498}.\qed
\end{bew}

\newpage

\section{Foundations and Root Group Systems}

Given a root group system, there is a natural way to attach a foundation to it.

\begin{de}\label{60}
Let $\UUU(\BBB,M,\Sigma,c)$ be a root group system. 
\begin{itemize}
\item Given $(i,j)\in A(M)$, there is an alternative division ring $\AA_{(i,j)}$ such that $U_{(i,j)}\cong \TTT(\AA_{(i,j)})$. In particular, there is a system of parametrizations 
\begin{align*}
x_{(i,j)}^*:\AA_{(i,j)}\to U_{(i,j)}^*\ ,\ t\mapsto x_{(i,j)}^*(t)\ , && *\in\{i,ij,j\}
\end{align*}
extending to the defining relations for $\TTT(\AA_{(i,j)})$, i.e., we have
$$\forall\ s,t\in\AA_{(i,j)}:\qquad [x_{(i,j)}^i(s),x_{(i,j)}^j(t)]=x_{(i,j)}^{ij}(st)\ .$$
By lemma \ref{359} and lemma \ref{360}, such a parametrization yields an opposite system of parametrizations
\begin{align*}
x_{(j,i)}^j&:\AA_{(i,j)}^{o}\to U_{(j,i)}^j,\ t\mapsto x_{(i,j)}^j\big(\id^o(t)\big)\ , \\
x_{(j,i)}^{ji}&:\AA_{(i,j)}^{o}\to U_{(j,i)}^{ji},\ t\mapsto x_{(i,j)}^{ij}\big(\id^o(-t)\big)\ , \\
x_{(j,i)}^i&:\AA_{(i,j)}^{o}\to U_{(j,i)}^i,\ t\mapsto x_{(i,j)}^i\big(\id^o(t)\big)\ .
\end{align*}
The resulting set $\Lambda:=\{ \TTT(\AA_{(i,j)}) \mid (i,j)\in A(M)\}$ is a \textit{parameter system for $\mathit{\UUU}$}\index{parameter system}.
\item Given $(i,j,k)\in G(M)$ and parametrizations $\TTT(\AA_{(i,j)})$ and $\TTT(\AA_{(j,k)})$, we define the glueing  $\gamma_{(i,j,k)}:\AA_{(i,j)}\to \AA_{(j,k)}$ by
$$x_{(i,j)}^j(t)=x_{(j,k)}^j\big(\gamma_{(i,j,k)}(t)\big)$$
which is justified by lemma \ref{59}. Then $\gamma_{(i,j,k)}$ is an isomorphism of additive groups satisfying $\gamma_{(i,j,k)}=\id^o\circ\gamma_{(k,j,i)}^{-1}\circ\id^o$. By lemma \ref{362}, we may adjust all the parametrizations such that $$\forall\ (i,j,k)\in G(F):\qquad \gamma_{(i,j,k)}(1)=1\ .$$
Notice that for this purpose we need the following fact: The adjustment of a glueing $\gamma_{(i,j,k)}$ can be realized by a reparametrization for $\TTT(\AA_{(i,j)})$ which fixes $x_{(i,j)}^i(1)$. Thus we can make sure that we don't alter glueings which have already been adjusted before.
\end{itemize}
\end{de}

\newglossaryentry{fl}{type=symbols,name={\ensuremath{\FFF(U,\Lambda)}},description=foundation with respect to the root group system $\UUU$ and the parametrization $\Lambda$, sort=foundation}

\begin{lemma} Given a root group system $\UUU:=\UUU(\BBB,M,\Sigma,c)$, a parameter system $\Lambda$ as in definition \ref{60} induces a foundation $$\gls{fl}:=\{ \TTT(\AA_{(i,j)}),\gamma_{(i,j,k)}\mid (i,j)\in A(M), (i,j,k)\in G(M)\}\ .$$
\end{lemma}

\begin{bew}
We emphasize that the glueings in definition \ref{60} are identifications with respect to directed edges. Given $(i,j,k),(i,j,l),(l,j,k)\in G(M)$ and $t\in \AA_{(i,j)}$, we have
\begin{align*} 
x_{(j,k)}^j\big(\gamma_{(i,j,k)}(t)\big)&=x_{(i,j)}^j(t)=x_{(j,l)}^j\big(\gamma_{(i,j,l)}(t)\big)\\
&=x_{(l,j)}^j\big(\id^o\circ \gamma_{(i,j,l)}(t)\big)=x_{(j,k)}^j\big(\gamma_{(l,j,k)}\circ\id^o\circ \gamma_{(i,j,l)}(t)\big)\end{align*}
and thus $\gamma_{(i,j,k)}=\gamma_{(l,j,k)}\circ \id^o\circ \gamma_{(i,j,l)}$.\qed
\end{bew}

\begin{de} A foundation $\FFF$ is \textit{integrable}\index{integrable}\index{foundation!integrable} if it is the foundation of a twin building $\BBB$, i.e., if there are a root group system $\UUU:=\UUU(\BBB,M,\Sigma,c)$ and a parameter system $\Lambda$ for $\UUU$ such that $$\FFF=\FFF(\UUU,\Lambda)\ .$$
\end{de}

\begin{bem}\ 
\begin{enumerate}[label=(\alph*)]
\item We will see in lemma \ref{33} that an integrable foundation is necessarily a Moufang foundation.
\item The next step is to show that the foundation attached to a root group system is unique up to isomorphism. Moreover, we want to prove that the building corresponding to an integrable foundation is unique up to isomorphism.
\end{enumerate}
\end{bem}

\begin{prop}\label{71}
Let $\UUU:=\UUU(\BBB,M,\Sigma,c)$ and $\tilde{\UUU}:=\UUU(\tilde{\BBB},\tilde{M},\tilde{\Sigma},\tilde{c})$ be root group systems and let $\Lambda$ and $\tilde{\Lambda}$ be parameter systems for $\UUU$ and $\tilde{\UUU}$, respectively. Then the following holds:
\begin{enumerate}[label=(\alph*)]
\item An isomorphism $\tilde{\alpha}:\FFF(\UUU,\Lambda)\to \FFF(\tilde{\UUU},\tilde{\Lambda})$ induces an isomorphism $\alpha:\UUU\to \tilde{\UUU}$.
\item An isomorphism $\alpha:\UUU\to \tilde{\UUU}$ induces an isomorphism $\tilde{\alpha}: \FFF(\UUU,\Lambda)\to \FFF(\tilde{\UUU},\tilde{\Lambda})$.
\end{enumerate}
\end{prop}

\begin{bew}
Each isomorphism $$\alpha_{(i,j)}:U_{(i,j)}\to \tilde{U}_{(\pi(i),\pi(j))}$$ induces an isomorphism
$$\tilde{\alpha}_{(i,j)}:\TTT(\AA_{(i,j)})\to \TTT(\tilde{\AA}_{(\pi(i),\pi(j))})$$
and vice versa. Given $(i,j)\in A(M)$, we have
$$ \alpha_{(i,j)}=\alpha_{(j,i)}^o\ \Leftrightarrow\ \tilde{\alpha}_{(i,j)}=\tilde{\alpha}_{(j,i)}^o\ .$$
\begin{enumerate}[label=(\alph*)]
\item We show that
$$\alpha:=\{ \pi,\alpha_{(i,j)} \mid (i,j)\in A(M) \}:\UUU\to\tilde{\UUU}$$
is an isomorphism. 

Given $(i,j,k)\in G(M)$ and $t\in \AA_{(i,j)}$, we have
\begin{align*}
\big(x_{(i,j)}^j(t)\big)^{\alpha_{(i,j)}}&=\tilde{x}_{(\pi(i),\pi(j))}^{\pi(j)}\big(\tilde{\alpha}_{(i,j)}^j(t)\big)=\tilde{x}_{(\pi(j),\pi(k))}^{\pi(j)}\big(\tilde{\gamma}_{(\pi(i),\pi(j),\pi(k))}\circ \tilde{\alpha}_{i,j}^j(t)\big) \\
&=\tilde{x}_{(\pi(j),\pi(k))}^{\pi(j)}\big(\tilde{\alpha}_{(j,k)}^j\circ \gamma_{(i,j,k)}(t)\big)=\big(x_{(j,k)}^j(\gamma_{(i,j,k)}(t))\big)^{\alpha_{(j,k)}}
\end{align*}
and therefore
$${\alpha_{(i,j)}}_{|U_{(i,j)}^j}={\alpha_{(j,k)}}_{|U_{(j,k)}^j}\ .$$
\item We show that
$$\tilde{\alpha}:=\{ \pi,\tilde{\alpha}_{(i,j)} \mid (i,j)\in A(M) \}:\FFF(\UUU,\Lambda)\to \FFF(\tilde{\UUU},\tilde{\Lambda})$$
is an isomorphism.

Given $(i,j,k)\in G(M)$ and $t\in \AA_{(i,j)}$, we have
$$ U_{(i,j)}^j\ni x_{(i,j)}^j (t)=x_{(j,k)}^j\big(\gamma_{(i,j,k)}(t)\big)\in U_{(j,k)}^j\ .$$
As we have \begin{align}{\alpha_{(i,j)}}_{|U_{(i,j)}^j}={\alpha_{(j,k)}}_{|U_{(j,k)}^j}\ , \label{136}\end{align}
it follows that
\begin{align*}
\tilde{x}_{(\pi(j),\pi(k))}^{\pi(j)}\big(\tilde{\gamma}_{(\pi(i),\pi(j),\pi(k))}\circ \tilde{\alpha}_{i,j}^j(t)\big)&=\tilde{x}_{(\pi(i),\pi(j))}^{\pi(j)}\big(\tilde{\alpha}_{(i,j)}^j(t)\big) \\
&\overset{\mathclap{\eqref{136}}}{=}\tilde{x}_{(\pi(j),\pi(k))}^{\pi(j)}\big(\tilde{\alpha}_{(j,k)}^j\circ \gamma_{(i,j,k)}(t)\big)
\end{align*}
and therefore
$$\tilde{\gamma}_{(\pi(i),\pi(j),\pi(k))}\circ \tilde{\alpha}_{(i,j)}^j=\tilde{\alpha}_{(j,k)}^j\circ \gamma_{(i,j,k)}\ .$$
\end{enumerate}\qed
\end{bew}

\section{Reparametrizations and Isomorphisms}

The concept of reparametrizations is quite similar to that of isomorphisms. However, we deal with a single foundation and produce (in fact, all the) foundations which are isomorphic to a given one. Moreover, this concept allows us to complete the proof that a foundation is a classifying invariant of the corresponding twin building.

\begin{de}\label{142}
Let $\FFF$ be a foundation.
\begin{itemize}
\item A system of reparametrizations
$$\alpha:=\{ \alpha_{(i,j)} \mid (i,j)\in A(F) \}$$
satisfying $\alpha_{(i,j)}=\alpha_{(j,i)}^o$ for each $(i,j)\in A(F)$ and
$$\gamma_{(i,j,k)}\circ \alpha_{(i,j)}^j(1)=\alpha_{(j,k)}^j(1)$$
for each $(i,j,k)\in G(F)$ is a \textit{reparametrization for $\mathit{\FFF}$}.

\newglossaryentry{fr}{type=symbols,name={\ensuremath{\FFF_\alpha}},description=foundation with respect to the reparametrization $\alpha$, sort=foundation}
\item Given a reparametrization $\alpha$ for $\FFF$, we set
$$\gls{fr}:=\{ \TTT(\tilde{\AA}_{(i,j)}),\tilde{\gamma}_{(i,j,k)}\mid (i,j)\in A(F),(i,j,k)\in G(F)\}$$
with
$$\tilde{\gamma}_{(i,j,k)}:=(\alpha_{(j,k)}^j)^{-1}\circ \gamma_{(i,j,k)}\circ \alpha_{(i,j)}^j$$
for each $(i,j,k)\in G(F)$.
\end{itemize}
\end{de}

\begin{bsp} Given the foundation
\begin{center}
\begin{tikzpicture}[scale=0.7,thick,>=stealth]
\node() at (-4,0) {\ssi{$\FFF$}};
\draws{0.55} (-3,0)--(0,0);
\draws{0.55} (0,0)--(3,0);
\fill (-3,0) circle (3pt);
\fill (0,0) circle (3pt);
\fill (3,0) circle (3pt);
\draw[->](-1,0)arc[radius=1,start angle=180, end angle=0];
\node()at (-3.0,-0.5) {\ssi$1$};
\node()at (0.0,-0.5) {\ssi$2$};
\node()at (3.0,-0.5) {\ssi$3$};
\node()at (-1.5,-0.5) {\ssi$\TTT(\AA)$};
\node()at (1.5,-0.5) {\ssi$\TTT(\AA)$};
\node()at (0,0.5) {\ssi{$\gamma$}};
\node() at (4,0) {};
\end{tikzpicture}
\end{center}
with $\gamma\in \Aut(\AA)$ and $\alpha:=\{ \alpha_{(1,2)}:=(\AA,\id_\AA,\id_\AA,\id_\AA), \alpha_{(2,3)}:=(\AA,\gamma,\gamma,\gamma)\}$, we have
\begin{center}
\begin{tikzpicture}[scale=0.7,thick,>=stealth]
\node() at (-4,0) {\ssi{$\FFF_\alpha$}};
\draws{0.55} (-3,0)--(0,0);
\draws{0.55} (0,0)--(3,0);
\fill (-3,0) circle (3pt);
\fill (0,0) circle (3pt);
\fill (3,0) circle (3pt);
\draw[->](-1,0)arc[radius=1,start angle=180, end angle=0];
\node()at (-3.0,-0.5) {\ssi$1$};
\node()at (0.0,-0.5) {\ssi$2$};
\node()at (3.0,-0.5) {\ssi$3$};
\node()at (-1.5,-0.5) {\ssi$\TTT(\AA)$};
\node()at (1.5,-0.5) {\ssi$\TTT(\AA)$};
\node()at (0,0.5) {\ssi{$\id_\AA$}};
\node() at (4,0) {.};
\end{tikzpicture}
\end{center}
\end{bsp} 

\vfill

\begin{lemma}\label{84}
Let $\UUU:=\UUU(\BBB,M,\Sigma,c)$ be a root group system, let $\FFF:=\FFF(\UUU,\Lambda)$ for some parameter system $\Lambda$ for $\UUU$, let 
$\alpha$ be a reparametrization for $\FFF$ and let $\tilde{\Lambda}$ be the parameter system induced by $\alpha$. Then we have $\tilde{\FFF}:=\FFF(\UUU,\tilde{\Lambda})=\FFF_\alpha$.
\end{lemma}

\begin{bew}
We have
\begin{align*}
\tilde{x}_{(j,k)}^j\big(\tilde{\gamma}_{(i,j,k)}(t)\big)&=\tilde{x}_{(i,j)}^j(t)=x_{(i,j)}^j\big(\alpha_{(i,j)}^j(t)\big) \\
&=x_{(j,k)}^j\big(\gamma_{(i,j,k)}\circ \alpha_{(i,j)}^j(t)\big)=\tilde{x}_{(j,k)}^j\big((\alpha_{(j,k)}^j)^{-1}\circ \gamma_{(i,j,k)}\circ \alpha_{(i,j)}^j(t)\big)
\end{align*}
for each $t\in \tilde{\AA}_{(i,j)}$.\qed
\end{bew}

\begin{kor}\label{82}
Let $\UUU:=\UUU(\BBB,M,\Sigma,c)$ be a root group system, let $\FFF:=\FFF(\UUU,\Lambda)$ for some parameter system $\Lambda$ for $\UUU$ and let 
$$\alpha=\{\pi,\alpha_{(i,j)}\mid (i,j)\in A(F)\}:\FFF\to \tilde{\FFF}$$
be an isomorphism. Then $\tilde{\FFF}$ is integrable.
\end{kor}

\newpage

\begin{bew}
Take $\big(\tilde{\AA}_{(i,j)}:=\tilde{\AA}_{(\pi(i),\pi(j))}, (\alpha_{(i,j)}^i)^{-1},(\alpha_{(i,j)}^{ij})^{-1},(\alpha_{(i,j)}^j)^{-1}\big)$ as reparametrization for $\TTT(\AA_{(i,j)})$, then replace $i\in I$ by $\pi(i)\in \tilde{I}$. The resulting parameter system $\tilde{\Lambda}$
satisfies
$$\FFF(\UUU,\tilde{\Lambda})=\FFF_\alpha=\tilde{\FFF}\ .$$\qed \end{bew}

\begin{satz}
The isomorphism class of an integrable foundations $\FFF=\FFF(\UUU,\Lambda)$ is a classifying invariant of the isomorphism class of the corresponding building.
\end{satz}

\begin{bew}
This results from corollary \ref{82}, proposition \ref{71} and theorem \ref{137}.\qed
\end{bew}

\begin{bem}
The following theorem shows that the concept of reparametrization is useful if we want to determine all the foundations isomorphic to a given foundation $\FFF$.
\end{bem}

\begin{satz}
Let $\FFF,\tilde{\FFF}$ be foundations with $F=\tilde{F}$. Then the following holds:
\begin{enumerate}[label=(\alph*)]
\item \label{141} Let $\tilde{\alpha}=\{ \tilde{\alpha}_{(i,j)}\mid (i,j)\in A(F)\}:\FFF\to \tilde{\FFF}$ be a special isomorphism. Then there is a reparametrization $\alpha$ of $\FFF$ such that $\FFF_{\alpha}=\tilde{\FFF}$.
\item Let $\alpha=\{ \alpha_{(i,j)} \mid (i,j)\in A(F)$ be a reparametrization for $\FFF$ such that $\FFF_\alpha=\tilde{\FFF}$. Then there is a special isomorphism $\tilde{\alpha}:\FFF\to \tilde{\FFF}$.
\end{enumerate}
\end{satz}

\begin{bewzwei}\ 
\begin{enumerate}[label=(\alph*)]
\item If we take $\alpha:=\{ \alpha_{(i,j)} \mid (i,j)\in A(F)\}$ with
$$\alpha_{(i,j)}:=\{ \tilde{\AA}_{(i,j)}, (\tilde{\alpha}_{(i,j)}^i)^{-1},(\tilde{\alpha}_{(i,j)}^{ij})^{-1},(\tilde{\alpha}_{(i,j)}^j)^{-1} \}$$
as reparametrization for $\FFF$, then $\FFF_\alpha=\tilde{\FFF}$.
\item We have $$\alpha_{(i,j)}=(\tilde{\AA}_{(i,j)}, \alpha_{(i,j)}^i,\alpha_{(i,j)}^{ij},\alpha_{(i,j)}^j)$$
for each $(i,j)\in A(F)$, thus $\tilde{\alpha}:=\{ \id_F,\tilde{\alpha}_{(i,j)} \mid (i,j)\in A(F)\}:\FFF\to \tilde{\FFF}$ with
$$\tilde{\alpha}_{(i,j)}:=\big((\alpha_{(i,j)}^i)^{-1},(\alpha_{(i,j)}^{ij})^{-1},(\alpha_{(i,j)}^j)^{-1}\big)$$
is an isomorphism.
\end{enumerate}\qed
\end{bewzwei}

\begin{bem}\label{85}\ 
\begin{enumerate}[label=(\alph*)]
\item Let $\FFF$ and $\tilde{\FFF}$ be foundations and let $$\alpha=\{\pi,\alpha_{(i,j)}\mid (i,j)\in A(F)\}:\FFF\to \tilde{\FFF}$$ be an isomorphism. As we may replace $i\in V(F)$ by $\pi(i)\in V(\tilde{F})$, we may consider $\alpha$ as special. Thus it suffices to determine all foundations which are specially isomorphic to $\FFF$. The remaining foundations isomorphic to $\FFF$ are obtained by relabelings of the vertex set.
\item \label{143} The theorem is useful if we want to show that two given foundations $\FFF$ and $\tilde{\FFF}$ with isomorphic residues $\RRR$ and $\tilde{\RRR}$ are isomorphic. In this case we may replace $\RRR$ by $\tilde{\RRR}$, observing that there is a relabeling of the corresponding vertices involved.
\end{enumerate}
\end{bem}

\section{Glueings}

At this point we collect some first results about glueings of integrable foundations.

\begin{no} Throughout the rest of this part, $\FFF$ denotes an integrable foundation.
\end{no}

\begin{satz}
Let $(i,j)\in A(F)$. Then a Hua automorphism of $\TTT(\AA_{(i,j)})$ is induced by an automorphism of $\FFF$.
\end{satz}

\begin{bew}
This results from theorem \ref{452}.\qed
\end{bew}

\begin{lemma}
Let $(i,j,k)\in G(F)$ and let $\AA:=\AA_{(i,j)},\ \tilde{\AA}:=\AA_{(j,k)}$. Then the following holds:
\begin{enumerate}[label=(\alph*)]
\item \label{33} The glueing $\gamma:=\gamma_{(i,j,k)}$ is a Jordan isomorphism. In particular, $\FFF$ is a Moufang foundation.
\item \label{138} If $\gamma_{(i,j,k)}$ is an isomorphism of alternative rings, then $\TTT(\tilde{\AA})\cong \TTT(\AA)$.
\item \label{139} If $\gamma_{(i,j,k)}$ is an anti-isomorphism of alternative rings, then $\TTT(\tilde{\AA})\cong \TTT(\AA^{o})$.
\end{enumerate}
\end{lemma}

\begin{bewzwei}\ 
\begin{enumerate}[label=(\alph*)]
\item If we set
\begin{align*}
h(t):=\mu\big(x_{(i,j)}^j(1_\AA)\big)^{-1}\mu\big(x_{(i,j)}^j(t)\big),\ t\in \AA\ , &&\tilde{h}(\tilde{t}):=\mu\big(x_{(j,k)}^j({1_{\tilde{\AA}}})\big)^{-1}\mu\big(x_{(j,k)}^j(\tilde{t})\big),\  \tilde{t}\in \tilde{\AA}\ , 
\end{align*}
we have
$$\tilde{h}\big(\gamma(t)\big)=\mu\big(x_{(j,k)}^j(1_{\tilde{\AA}})\big)^{-1}\mu\big(x_{(j,k)}^j(\gamma(t))\big)=\mu\big(x_{(i,j)}^j(1_\AA)\big)^{-1}\mu\big(x_{(i,j)}^j(t)\big)=h(t)$$
for each $t\in \AA$. Moreover, we have
\begin{align*}
x_{(i,j)}^j(t)^{h(s)}=x_{(i,j)}^j(sts)\ , && x_{(j,k)}^j(\tilde{t})^{\tilde{h}(\tilde{s})}=x_{(j,k)}^j(\tilde{s}\tilde{t}\tilde{s})
\end{align*}
for all $s,t\in \AA,\ \tilde{s},\tilde{t}\in\tilde{\AA}$, cf. lemma \ref{361}. Combining these two facts yields
\begin{align*}
x_{(j,k)}^j\big(\gamma(s)\gamma(t)\gamma(s)\big)&=x_{(j,k)}^j\big(\gamma(t)\big)^{\tilde{h}(\gamma(s))}=x_{(i,j)}^j(t)^{h(s)}=x_{(i,j)}^j(sts)=x_{(j,k)}^j\big(\gamma(sts)\big)
\end{align*}
and therefore
$$\gamma(sts)=\gamma(s)\gamma(t)\gamma(s)$$
for all $s,t\in \AA$, thus $\gamma$ is a Jordan isomorphism. 
\item If $\gamma$ is an isomorphism, then $(\AA,\gamma,\gamma,\gamma)$ is a reparametrization for $\TTT(\tilde{\AA})$, thus $\TTT(\tilde{\AA})\cong \TTT(\AA)$.
\item If $\gamma$ is an anti-isomorphism, then $(\AA^o,\gamma^o,\gamma^o,\gamma^o)$ is a reparametrization for $\TTT(\tilde{\AA})$, thus $\TTT(\tilde{\AA})\cong \TTT(\AA^o)$.
\end{enumerate}\qed
\end{bewzwei}

\begin{de}\ 
\begin{itemize}
\item A glueing is \textit{negative}\index{negative}\index{glueing!negative} if it is an isomorphism of alternative rings.
\item A glueing is \textit{positive}\index{positive}\index{glueing!positive} if it is an anti-isomorphism of alternative rings.
\item A glueing is \textit{exceptional}\index{exceptional}\index{glueing!exceptional} if it is neither positive nor negative.
\item A foundation is \textit{negative}\index{foundation!negative} if each glueing is negative.
\item A foundation is \textit{positive}\index{foundation!positive} if each glueing is positive.
\item A foundation is \textit{mixed}\index{mixed}\index{foundation!mixed} if there are both positive and negative glueings.
\end{itemize}
\end{de}

\newpage
\begin{prop}\label{22} A foundation $$\FFF=\{\TTT(\AA_{(1,2)}),\TTT(\AA_{(2,3)}),\gamma:=\gamma_{(1,2,3)}\}$$ of type $A_3$ is negative. Moreover, $\AA$ (and thus $\tilde{\AA}$) is associative.
\end{prop}

\begin{bew} 
If we set $\AA:=\AA_{(1,2)}, \tilde{\AA}:=\AA_{(2,3)}$, 
\begin{align*}
x_1:=x_{(1,2)}^1\ , && x_2:=x_{(1,2)}^2\ , && \tilde{x}_2:=x_{(2,3)}^2\ , && \tilde{x}_3:=\tilde{x}_{(2,3)}^3
\end{align*}
and
\begin{align*}
h_1(t)&:=\mu\big(x_1(1_\AA)\big)^{-1}\mu\big(x_1(t^{-1})\big)\in U_{(1,2)}^1,\qquad t\in \AA\ ,\\
h_3(\tilde{t})&:=\mu\big(\tilde{x}_3(1_{\tilde{\AA}})\big)^{-1}\mu\big(\tilde{x}_3(\tilde{t}^{-1})\big)\in U_{(2,3)}^3,\qquad \tilde{t}\in\tilde{\AA}\ ,
\end{align*}
lemma \ref{361} yields
\begin{align*}
x_2(s)^{h_1(t)}=x_2(t\cdot s)\ , && \tilde{x}_2(\tilde{s})^{h_3(\tilde{t})}=\tilde{x}_2(\tilde{s}\ast \tilde{t})
\end{align*}
for all $s,t \in \AA,\ \tilde{s},\tilde{t}\in\tilde{\AA}$. Combining these two facts and observing $[U_{(1,2)}^1,U_{(2,3)}^3]=1$ yield
\begin{align*}
\tilde{x}_2\big(\gamma(t)\ast \gamma(s)\big)&=\tilde{x}_2\big(\gamma(t)\big)^{h_3(\gamma(s))}=x_2(t)^{h_3(\gamma(s))}=x_2(1_\AA)^{h_1(t)h_3(\gamma(s))}\\
&=\tilde{x}_2(1_{\tilde{\AA}})^{h_3(\gamma(s))h_1(t)}=\bar{x}_2\big(\gamma(s)\big)^{h_1(t)}=x_2(s)^{h_1(t)}=x_2(t\cdot s)=\tilde{x}_2\big(\gamma(t\cdot s)\big)
\end{align*}
and thus $\gamma(t\cdot s)=\gamma(t)\ast\gamma(s)$ for all $s,t\in \AA$. As a consequence, we obtain
\begin{align*}
x_2(s\cdot t)=\tilde{x}_2\big(\gamma(s\cdot t)\big)=\tilde{x}_2\big(\gamma(s)\ast\gamma(t)\big)=\tilde{x}_2\big(\gamma(s)\big)^{h_3(\gamma(t))}=x_2(s)^{h_3(\gamma(t))}
\end{align*}
for all $s,t\in\AA$. This implies
\begin{align*} 
x_2\big((s\cdot t)\cdot u\big)&=x_2(s\cdot t)^{h_3(\gamma(u))}=x_2(t)^{h_1(s)h_3(\gamma(u))}\\
&=x_2(t)^{h_3(\gamma(u))h_1(s)}=x_2(t\cdot u)^{h_1(s)}=x_2\big(s\cdot (t\cdot u)\big)\end{align*}
and therefore $(s\cdot t)\cdot u=s\cdot (t\cdot u)$ for all $s,t,u\in \AA$. 
\qed
\end{bew}

\begin{kor}\label{144}
Let $\FFF$ be a positive foundation. Then $\GGG_F$ is a complete graph.
\end{kor}

\newglossaryentry{Hua}{type=results,name={{Hua's Theorem}},description={},sort=res}

\begin{satz}[\textbf{\gls{Hua}}]\label{386}
Let $\AA$ and $\tilde{\AA}$ be alternative division rings and let $\gamma:\AA\to\tilde{\AA}$ be a Jordan homomorphism. If $\tilde{\AA}$ is associative, then $\gamma:\AA\to \gamma(\AA)$ is an iso- or anti-isomorphism of alternative rings. In particular, $\AA$ and $\gamma(\AA)$ are also skew-fields, and, if $\AA$ or $\gamma(\AA)$ is a field, the map $\gamma$ is an isomorphism of fields.
\end{satz}

\begin{bewzwei}\ 
\begin{enumerate}[label=(\roman*)]
\item We show: Given $s,t\in \AA$, we have $\gamma(st)=\gamma(s)\gamma(t)$ or $\gamma(st)=\gamma(t)\gamma(s)$.

The assertion is clearly true for $s=0_\AA$ or $t=0_\AA$, so assume $s\neq0_\AA\neq t$. As we have
$$\gamma(u^2)=\gamma(u\cdot 1_\AA\cdot u)=\gamma(u)\gamma(1_\AA)\gamma(u)=\gamma(u)\cdot 1_{\tilde{\AA}}\cdot\gamma(u)=\gamma(u)^2$$
for each $u\in \AA$, it follows that
\begin{align*}
\gamma\big((s+t)^2\big)&=\gamma(s^2+st+ts+t^2)=\gamma(s)^2+\gamma(st)+\gamma(ts)+\gamma(t)^2\ ,\\
\gamma\big((s+t)^2\big)&=\gamma(s+t)^2=\big(\gamma(s)+\gamma(t)\big)^2=\gamma(s)^2+\gamma(s)\gamma(t)+\gamma(t)\gamma(s)+\gamma(t)^2
\end{align*}
and thus
\begin{align}
\gamma(st)+\gamma(ts)=\gamma(s)\gamma(t)+\gamma(t)\gamma(s)\ . \label{31}
\end{align}
\newpage \noindent On the other hand, we have
\begin{align*}
s\big(t(st)^{-1}t\big)s&=s([t(t^{-1}s^{-1})]t)s=[s(s^{-1}t)]s=ts
\end{align*}
by the Moufang identities and the inverse properties and, by lemma \ref{127}, therefore,
\begin{align}
\gamma(ts)&=\gamma\big(s(t(st)^{-1}t)s\big)=\gamma(s)\big(\gamma(t)\gamma(st)^{-1}\gamma(t)\big)\gamma(s)\ . \label{32}
\end{align}
Observing the associativity of $\tilde{\AA}$, we obtain
\begin{align*}
&\ [1_{\tilde{\AA}}-\gamma(s)\gamma(t)\gamma(st)^{-1}][\gamma(st)-\gamma(t)\gamma(s)] \\
=&\ \gamma(st)-\gamma(s)\gamma(t)-\gamma(t)\gamma(s)+\gamma(s)\gamma(t)\gamma(st)^{-1}\gamma(t)\gamma(s) \\
\overset{\mathclap{\eqref{32}}}{=}&\ \gamma(st)-\gamma(s)\gamma(t)-\gamma(t)\gamma(s)+\gamma(ts)\overset{\mathclap{\eqref{31}}}{=}0_{\tilde{\AA}}\ .
\end{align*}
Since $\tilde{\AA}$ has no zero divisors, it finally follows that
$$1_{\tilde{\AA}}-\gamma(s)\gamma(t)\gamma(st)^{-1}=0_{\tilde{\AA}}\qquad \vee\qquad \gamma(st)-\gamma(t)\gamma(s)=0_{\tilde{\AA}}\ .$$
\item Given $s\in \AA$, let
\begin{align*}
N_s:=\{t\in \AA\mid \gamma(st)=\gamma(s)\gamma(t)\}\ , && P_s:=\{t\in \AA\mid \gamma(st)=\gamma(t)\gamma(s)\}\ .
\end{align*}
By step (i), the subgroups $N_s$ and $P_s$ of $(\AA,+)$ satisfy $\AA=N_s\cup P_s$. As no group is the union of two proper subgroups, we obtain
$$N_s=\AA\qquad\vee\qquad P_s=\AA\ .$$
\item If we set 
\begin{align*}
N:=\{s\in \AA\mid N_s=\AA\}\ , && P:=\{s\in \AA\mid P_s=\AA\}\ ,
\end{align*}
step (ii) shows that $N$ and $P$ are subgroups of $(\AA,+)$ satisfying $\AA=N\cup P$, hence
$$N=\AA\qquad\vee\qquad P=\AA\ .$$
\end{enumerate}\qed
\end{bewzwei}

\begin{bem}
If we suppose $\AA$ instead of $\tilde{\AA}$ to be associative and if $\gamma(\AA)$ is an alternative division ring, the theorem remains true as we may apply it to $\gamma^{-1}:\gamma(\AA)\to \AA$.
\end{bem}

\begin{kor}\label{34}
Let $\FFF$ be a foundation such that there exists an edge $(a,b)\in A(F)$ with $\AA:=\AA_{(a,b)}$ associative. Then we have
\begin{align*} \forall\ (i,j)\in A(F):&&  \AA_{(i,j)}\cong \AA\ \vee\ \AA_{(i,j)}\cong \AA^o\ , && \TTT(\AA_{(i,j)})\cong \TTT(\AA)\ \vee\ \TTT(\AA_{(i,j)})\cong \TTT(\AA^o)\ .
\end{align*}
In particular, the alternative division ring $\AA_{(i,j)}$ is associative for each $(i,j)\in A(F)$. Moreover, the assumption is satisfied if $\FFF$ has a residue of type $A_3$. 
\end{kor}

\begin{bew}
This results from lemma \ref{33}, Hua's theorem and lemma \ref{138}, (c), using an easy induction. The final assertion results from proposition \ref{22}. \qed
\end{bew}

\begin{kor}\label{50}
Let $\FFF$ be a foundation such that there exists an edge $(a,b)\in A(F)$ with $\AA:=\AA_{(a,b)}$ non-associative, i.e., $\AA$ is an octonion division algebra. Then $\GGG_F$ is a complete graph.
\end{kor}

\begin{bew}
By proposition \ref{22}, $\FFF$ satisfies the assumption of corollary \ref{34} if it has a residue of type $A_3$.\qed
\end{bew}

\newpage

\chapter{Integrability of Certain Foundations}\label{150}
In this chapter, we prove the integrability of certain foundations which are not negative or whose defining field is an octonion division algebra. Afterwards we will prove that these foundations are the only integrable foundations of this type, all of them defined over an octonion or quaternion division algebra. For this purpose, we need the theory of fixed point structures developed in \cite{MHab}.

The idea is as follows: Given an integrable foundation $\FFF$ and an automorphism $\alpha$ of the corresponding twin building $\BBB$, the foundation $\tilde{\FFF}$ corresponding to the twin building $\tilde{\BBB}:=\Fix(\alpha)$ can be constructed out of $\FFF$ and is, in fact, the fixed point structure of the automorphism induced on the group generated by all the root groups. Conversely, all such fixed point structures arise in this way, i.e., if we can realize a given foundation $\tilde{\FFF}$ as the fixed point structure of an automorphism $\alpha$ of an integrable foundation $\FFF$, then $\tilde{\FFF}$ itself is integrable.

We apply the so called theory of Tits indices without introducing it in detail since this would involve a lot of technical considerations. We just follow the recipes, hopefully in a natural and intuitive way, except for the canonical triangle over an octonion division algebra, where we just indicate the main idea.

But first of all we start with some general integrability criterions. In combination with Kac-Moody theory, they allow us to handle all the foundations that are defined over skew-fields distinct from a quaternion division algebra (which includes fields).

\section{Integrability Criterions}

The first one is a straight forward result as each residue of a twin building is again a twin building, the second one is taken from \cite{M}, and the third relies on Kac-Moody theory.

\begin{satz}
Let $M$ be a Coxeter matrix over $I=V(M)$, let $\FFF$ be an integrable foundation of type $M$ and let $J\subseteq I$ such that $|J|\geq 2$. Then $\FFF_J$ is integrable.
\end{satz}

\begin{bew}
By assumption, we have $\FFF=\FFF(\UUU,\Lambda)$ for some root group system $\UUU=\UUU(\BBB,M,\Sigma,c)$ and some parameter system $\Lambda$. By theorem \ref{436}, $\BBB_J:=\BBB_J(c)$ is a twin building, and by theorem \ref{438}, $\Sigma_J:=\Sigma_J(c)$ is a twin apartment of $\BBB_J$ containing $c$. Since we have
$$\forall\ i,j\in A(M_J):\qquad \Sigma\cap \BBB_{\{i,j\}}(c)=\Sigma_J\cap (\BBB_J)_{\{i,j\}}(c)\ ,$$
we have
$$\UUU_J:=\UUU(\BBB_J,M_J,\Sigma_J,c)=\{ U_{(i,j)} \mid i,j\in A(M_J)\}$$
and thus $\FFF_J=\FFF(\UUU_J,\Lambda_J)$, where $\Lambda_J$ is the parameter system induced by $\Lambda$.\qed
\end{bew}

\newglossaryentry{fc}{type=symbols,name={\ensuremath{\FFF(\tilde{F},\p)}},description=foundation with respect to the cover ${(\tilde{F},\p)}$, sort=foundation}

\begin{de} Let $\FFF$ be a foundation.
\begin{itemize}
\item Let $(\tilde{F},\p)$ be a cover of $F$. Then the foundation
$$\gls{fc}:=\{ \TTT(\tilde{\AA}_{(i,j)}), \tilde{\gamma}_{(i,j,k)} \mid (i,j)\in A(\tilde{F}), (i,j,k)\in G(\tilde{F})\}$$
with
\begin{align*}
\forall\ (i,j)\in A(\tilde{F}):\ \tilde{\AA}_{(i,j)}=\AA_{(\p(i),\p(j))}\ , && \forall\ (i,j,k)\in G(\tilde{F}):\ \tilde{\gamma}_{(i,j,k)}=\gamma_{(\p(i),\p(j),\p(k))}
\end{align*}
is the \textit{cover corresponding to $\mathit(\tilde{F},\p)$}\index{cover!of a foundation}.
\item A foundation $\tilde{\FFF}$ is a \textit{cover of $\mathit{\FFF}$} if there is a cover $(\tilde{F},\p)$ of $F$ such that $$\tilde{\FFF}\cong \FFF(\tilde{F},\p)\ .$$
\end{itemize}
\end{de}

\begin{bsp} Given the foundation
\begin{center}
\begin{tikzpicture}[scale=0.7,>=stealth,thick]
\coordinate (1) at (0,0);                    
\coordinate (2) at (3,0);
\coordinate (3) at (60:3);
\draws{0.55} (1)--(2);
\draws{0.55} (2)--(3);
\draws{0.55} (3)--(1);
\fill (1) circle (3pt);
\fill (3,0) circle (3pt);
\fill (60:3) circle (3pt);
\node()at (0.55,0.25) {\ssi{$\gamma_1$}};
\node()at (2.45,0.25) {\ssi$\gamma_2$};
\node()at (1.5,1.9) {\ssi$\gamma_3$};
\node()at (0.0,-0.5) {\ssi$\bar{1}$};
\node()at (3,-0.5) {\ssi$\bar{2}$};
\node()at (1.5,3) {\ssi$\bar{3}$};
\node()at (1.5,-0.5) {\ssi{$\TTT(\AA_{(1,2)})$}};
\node()at (3.3,1.5) {\ssi{$\TTT(\AA_{(2,3)})$}};
\node()at (-0.3,1.5) {\ssi{$\TTT(\AA_{(3,1)})$}};
\node()at (5.0,1.5) {\normalsize{,}};
\node()at (-2.0,1.5) {};
\draw[<-](60:2)arc[radius=1,start angle=240, end angle=300];
\draw[<-] (1,0) arc[radius=1,start angle=0, end angle=60];
\draw[->] (2,0) arc[radius=1,start angle=180, end angle=120];
\end{tikzpicture}\end{center}
a cover is given by
\begin{center}
\begin{tikzpicture}[scale=0.7,thick,>=stealth]
\draws{0.55} (-3,0)--(0,0);
\draws{0.55} (0,0)--(3,0);
\draws{0.55} (3,0)--(6,0);
\draws{0.55} (6,0)--(9,0);
\draw[dashed] (-5,0)--(-3,0);
\draw[dashed] (9,0)--(11,0);
\fill (-3,0) circle (3pt);
\fill (0,0) circle (3pt);
\fill (3,0) circle (3pt);
\fill (6,0) circle (3pt);
\fill (9,0) circle (3pt);
\node()at (-3.0,-0.5) {\ssi$0$};
\node()at (0.0,-0.5) {\ssi$1$};
\node()at (3.0,-0.5) {\ssi$2$};
\node()at (6.0,-0.5) {\ssi$3$};
\node()at (9.0,-0.5) {\ssi$4$};
\node()at (-3,0.5) {\ssi{$\gamma_3$}};
\node()at (0,0.5) {\ssi{$\gamma_1$}};
\node()at (3,0.5) {\ssi{$\gamma_2$}};
\node()at (6,0.5) {\ssi{$\gamma_3$}};
\node()at (9,0.5) {\ssi{$\gamma_1$}};
\node()at (-1.5,-0.5) {\ssi{$\TTT(\AA_{(3,1)})$}};
\node()at (1.5,-0.5) {\ssi{$\TTT(\AA_{(1,2)})$}};
\node()at (4.5,-0.5) {\ssi{$\TTT(\AA_{(2,3)})$}};
\node()at (7.5,-0.5) {\ssi{$\TTT(\AA_{(3,1)})$}};
\node()at (12,0) {\normalsize{,}};
\node()at (-6,0) {};
\draw[->](-4,0)arc[radius=1,start angle=180, end angle=0];
\draw[->](-1,0)arc[radius=1,start angle=180, end angle=0];
\draw[->](2,0)arc[radius=1,start angle=180, end angle=0];
\draw[->](5,0)arc[radius=1,start angle=180, end angle=0];
\draw[->](8,0)arc[radius=1,start angle=180, end angle=0];
\end{tikzpicture}\end{center}
where $\p:\ZZ\to \ZZ/3\ZZ,\ z\mapsto \bar{z}$ is the natural homomorphism.
\end{bsp}

\begin{satz}\label{145}
Let $\FFF$ be a foundation and let $\tilde{\FFF}$ be a cover of $\FFF$. Then $\FFF$ is integrable if $\tilde{\FFF}$ is integrable.
\end{satz}

\begin{bew}
This is a consequence of theorem C in \cite{M}.\qed
\end{bew}

\begin{de}
A foundation $\FFF$ such that $$\forall\ (i,j,k)\in G(F):\qquad \gamma_{(i,j,k)}=\id$$ is a \textit{canonical}\index{canonical}\index{foundation!canonical} foundation.
\end{de}

\begin{lemma}\label{446}
Let $\FFF$ be a negative foundation such that $\GGG_F$ is a tree. Then $\FFF$ is isomorphic to the corresponding canonical foundation.
\end{lemma}

\begin{bew}
Since $\GGG_F$ is a tree, it suffices to show the following: Given $(i,j,k)\in G(F)$, there is a reparametrization $$\alpha_{(j,k)}=(\AA_{(i,j)},\alpha_{(j,k)}^j,\alpha_{(j,k)}^{jk},\alpha_{(j,k)}^k)$$ such that $(\alpha_{(j,k)}^j)^{-1}\circ \gamma_{(i,j,k)}=\id_{\AA_{(i,j)}}$. This holds for $\alpha_{(j,k)}:=(\AA_{(i,j)},\gamma_{(i,j,k)},\gamma_{(i,j,k)},\gamma_{(i,j,k)})$.
\qed
\end{bew}

\begin{satz}\label{146}
Let $\FFF$ be a canonical foundation such that one of the following holds:
\begin{enumerate}[label=(\alph*)]
\item The defining field $\AA$ (cf. definition \ref{501}) is a field, and $F$ is a tree.
\item The defining field is a non-commutative skew-field, and $F$ is a string, a ray or a chain.
\end{enumerate}
Then $\FFF$ is integrable.
\end{satz}

\begin{bewzwei}\ 
\begin{enumerate}[label=(\alph*)]
\item Kac-Moody theory provides the existence of an integrable foundation $\tilde{\FFF}$ over $\AA$ such that $\tilde{F}=F$, cf. \cite{Ti}. Since $\tilde{\FFF}$ is negative by proposition \ref{22}, we have $\tilde{\FFF}\cong \FFF$ by lemma \ref{446}. Therefore, $\FFF$ is integrable by corollary \ref{82}.
\item The corresponding twin building is the limit of a sequence of twin buildings $A_n(\AA)$.
\end{enumerate}\qed
\end{bewzwei}

\section{The Canonical Triangle over an Octonion Division Algebra}

\newglossaryentry{TO}{type=foundations,name={\ensuremath{\tilde{\AAA}_2(\OO)}},description=canonical triangle with respect to the octonion division algebra $\OO$, sort=simply laced}

\begin{de} Let $\OO$ be an octonion division algebra. The foundation
$$\gls{TO}:=\{ \TTT(\AA_{(1,2)}):=\TTT(\AA_{(2,3)}):=\TTT(\AA_{(3,1)}):=\TTT(\OO), \gamma_{(1,2,3)}:=\gamma_{(2,3,1)}:=\gamma_{(3,1,2)}:=\id_\OO\}$$
is the \textit{canonical triangle over $\mathit{\OO}$}.
\end{de}

\begin{satz}\label{445}
The foundation $\tilde{A}_2(\OO)$ is integrable.
\end{satz}

\begin{bew}[\textbf{Sketch}]
Let $\KK:=Z(\OO)$ and let $\EE\subseteq \OO$ be such that $\EE/\KK$ is a quadratic separable extension. The foundation $\tilde{A}_2(\OO)$ is the fixed point structure of the Tits index
\begin{center}
\begin{tikzpicture}[scale=0.35,>=stealth,thick]
\coordinate (1) at (0,0);                    
\coordinate (2) at (0,3);
\coordinate (3) at (0,6);
\coordinate (4) at (-30:3);
\coordinate (5) at (-30:6);
\coordinate (6) at (-150:3);
\coordinate (7) at (-150:6);
\draw (1)--(2)--(3);
\draw (1)--(4)--(5);
\draw (1)--(6)--(7);
\foreach \i in {1,2,3,4,5,6,7} {\fill (\i) circle (3pt);}
\foreach \i in {3,5,7} {\draw (\i) circle (7pt);}
\node () at (0,-3) {$\EE/\KK$};
\node () at (6,0) {.};
\node () at (-6,0) {};
\node () at (1,6) {\ssi{$0$}};
\node () at (1,3) {\ssi{$3$}};
\node () at (0,-1) {\ssi{$6$}};
\node () at ($(4)-(0,1)$) {\ssi{$5$}};
\node () at ($(5)-(0,1)$) {\ssi{$2$}};
\node () at ($(6)-(0,1)$) {\ssi{$4$}};
\node () at ($(7)-(0,1)$) {\ssi{$1$}};
\end{tikzpicture}\end{center}
Let $\BBB$ be the twin building associated with the diagram, let $\RRR:=\RRR_{\{3,4,5,6\}}$ and let $\tau$ be the triality associated with $\RRR$. Given $k\in\{0,1,2\}$, let $\bar{k}:=\{0,\ldots,6\}\sm\{k\}$ and  $\tilde{k}:=\{0,1,2\}\sm\{k\}$. By similar arguments as in \cite{MGeo}, there are automorphisms $\p_k\in \Aut(\RRR_{\bar{k}})$ such that
\begin{align*}
\Fix(\p_k)\cong \TTT(\OO)\ , && \p_k(R_{\tilde{k}})=R_{\tilde{k}}\ , && {\p_k}_{|\RRR}\circ \tau =\tau\circ {\p_k}_{|\RRR}\ ,
\end{align*}
and $\p_0$ can be extended to an automorphism $\p\in \Aut(\BBB)$ such that $\p_{|\RRR_{\bar{k}}}=\p_k$. By a generalization of the arguments in chapter 3 of \cite{MPhD}, the fixed point set of $\p$ is a twin building of type $\tilde{A}_2$ whose foundation is $\tilde{A}_2(\OO)$.\qed
\end{bew}

\section{Positive and Mixed Foundations over Quaternions}\label{154}

\begin{no}
Throughout this paragraph, $\HH:=(\EE/\KK,\beta)$ is a quaternion division algebra with standard involution $\sigma_s$. 
\end{no}

\begin{no} Let $B:=(\begin{smallmatrix} 0 & \beta \\ 1 & 0\end{smallmatrix})\in GL_2(\EE)$ and let $$\sigma:M_2(\EE)\to M_2(\EE),\ X\mapsto B\bar{X}B^{-1}\ .$$
\end{no}

\begin{lemma} We have
$$\Fix(\sigma)=\left\{ \begin{pmatrix}  s & \beta\bar{t} \\ t & \bar{s}\end{pmatrix}\mid s,t\in\EE\right\}\cong \HH\ .$$
\end{lemma}

\begin{bew} Given $X:=(\begin{smallmatrix} s & u \\ t & v\end{smallmatrix})\in M_2(\EE)$, we have
\begin{align*}
X\in \Fix(\sigma)\ \Leftrightarrow\ \begin{pmatrix}  \bar{v} & \beta\bar{t} \\ \beta^{-1}\bar{u} & \bar{s}\end{pmatrix}=\begin{pmatrix}  s & u \\ t & v\end{pmatrix} \ \Leftrightarrow\ v=\bar{s}\ \wedge\ u=\beta\bar{t}\ \Leftrightarrow\ X\in \HH\ .
\end{align*}
\qed
\end{bew}

\begin{kor} If we extend $\sigma$ to an involution 
$$\sigma:M_6(\EE)\to M_6(\EE),\ X\mapsto \begin{pmatrix} B & & \\ & B & \\ & & B\end{pmatrix} \bar{X} \begin{pmatrix} B^{-1} & & \\ & B^{-1}& \\ & & B^{-1}\end{pmatrix}\ ,$$
we have
$$\bar{U}_+:=\Fix(\sigma)\cap U_+=\bigg\{\begin{pmatrix} I_2 & s& t \\ & I_2 & u \\ & & I_2\end{pmatrix} \mid s,t,u\in \HH\bigg\}\ .$$
\end{kor}

\begin{bem}\ 
\begin{enumerate}[label=(\alph*)]
\item  We translate the map $\sigma$ to root groups:\medskip

\noindent On the group $\langle U_{1,2},U_{2,3},U_{3,4},U_{2,1},U_{3,2},U_{4,3}\rangle$, the map $\sigma$ is given by
\begin{align*}
x_{1,2}(t)&\mapsto x_{2,1}(\beta^{-1}\bar{t})\ ,& x_{1,3}(t)&\mapsto x_{2,4}(\bar{t})\ ,\\
x_{2,3}(t)&\mapsto x_{1,4}(\beta\bar{t})\ , & x_{2,4}(t)&\mapsto x_{1,3}(\bar{t})\ , \\
x_{3,4}(t)&\mapsto x_{4,3}(\beta^{-1}\bar{t})\ , & x_{1,4}(t)&\mapsto x_{2,3}(\beta^{-1}\bar{t})\ .
\end{align*}
Moreover, the fixed point set in $\langle U_{1,2},U_{2,3},U_{3,4}\rangle$ is 
$$\{y^a(s,t):=x_{1,3}(s)x_{1,4}(\beta\bar{t})x_{2,3}(t)x_{2,4}(\bar{s})\mid (s,t)\in \HH\}\ .$$
On the group $\langle U_{3,4},U_{4,5},U_{5,6},U_{4,3},U_{5,4},U_{6,5}\rangle$, the map $\sigma$ is given by
\begin{align*}
x_{3,4}(t)&\mapsto x_{4,3}(\beta^{-1}\bar{t})\ , &x_{3,5}(t)&\mapsto x_{4,6}(\bar{t})\ , \\
x_{4,5}(t)&\mapsto x_{3,6}(\beta\bar{t})\ , & x_{4,6}(t)&\mapsto x_{3,5}(\bar{t})\ , \\
x_{5,6}(t)&\mapsto x_{6,5}(\beta^{-1}\bar{t})\ , & x_{3,6}(t)&\mapsto x_{4,5}(\beta^{-1}\bar{t})\ .
\end{align*}
Moreover, the fixed point set in $\langle U_{3,4},U_{4,5},U_{5,6}\rangle$ is
$$\{y^b(s,t):=x_{3,5}(s)x_{3,6}(\beta\bar{t})x_{4,5}(t)x_{4,6}(\bar{s})\mid (s,t)\in \HH\}\ .$$
If we set
$$\{y^{ab}(s,t):=x_{1,5}(s)x_{1,6}(\beta\bar{t})x_{2,5}(t)x_{2,6}(\bar{s})\mid (s,t)\in \HH\}\subseteq \bar{U}_+\ ,$$
we have
$$[y^a(s_1,t_1),y^b(s_2,t_2)]=y^{ab}\big((s_1,t_1)\cdot (s_2,t_2)\big)$$
for all $(s_1,t_1),(s_2,t_2)\in \HH$, hence
$$\big(\bar{U}_+, y^a(\HH),y^{ab}(\HH),y^b(\HH)\big)\cong \TTT(\HH)\ .$$
\item The corresponding Tits index is
\begin{center}
\begin{tikzpicture}[scale=0.35,>=stealth,thick]
\coordinate (1) at (-6,0);                    
\coordinate (2) at (-3,0);
\coordinate (3) at (0,0);
\coordinate (4) at (3,0);
\coordinate (5) at (6,0);
\draw (1)--(2)--(3)--(4)--(5);
\foreach \i in {1,2,3,4,5} {\fill (\i) circle (3pt);}
\foreach \i in {2,4} {\draw (\i) circle (7pt);}
\node () at (0,-1.5) {$\EE/\KK$};
\node () at (8,0) {.};
\node () at (-8,0) {};
\end{tikzpicture}\end{center}
\end{enumerate}
\end{bem}

\newglossaryentry{pn}{type=foundations,name={\ensuremath{\PPP_n^+(\HH)}},description=positive foundation with respect to the quaternion division algebra $\HH$, sort=simply laced}

\begin{satz}\label{40}
Let $\GGG$ be a complete graph over $I:=\{1,\ldots,n\}$. Then there exists an integrable positive foundation $\FFF:=\gls{pn}$ over $\HH$ such that $\GGG_F=\GGG$.

More precisely: For each $i\in I$, let $\BBB_i$ be a a copy of $\tilde{\BBB}$, and for all $1\leq i<j\leq n$ let $\BBB_{\{i,j\}}$ be a copy of $\BBB$. We identify $\BBB_i$ with $\BBB_{\{i,j\}}$ as in the first case and $\BBB_j$ with $\BBB_{\{i,j\}}$ as in the second case. Then the fixed point foundation $\FFF:=\PPP_n^+(\HH)$ is a positive foundation with defining field $\HH$ (cf. definition \ref{501}) and $\GGG_F=\GGG$. Moreover, we have $\gamma_{(i,j,k)}\in\{\sigma_s,\id^o\}$ for each $(i,j,k)\in G(F)$.
\end{satz}

\begin{bew}
Let 
\begin{align*}
\BBB:=\big(A_5(\EE)=(U_+,U_{1,2},U_{2,3},U_{3,4},U_{4,5},U_{5,6}),\sigma\big)\ , &&\tilde{\BBB}:=\big(A_3(\EE)=(\tilde{U}_+,\tilde{U}_{1,2},\tilde{U}_{2,3},\tilde{U}_{3,4}),\tilde{\sigma}\big)\end{align*} with the usual parametrizations.
\begin{enumerate}[label=(\roman*)]
\item We identify $\tilde{U}_+$ and $U_{1,2}U_{2,3}U_{3,4}$ via
\begin{align*}
\tilde{x}_{1,2}(t)&\mapsto x_{1,2}(t)\ , & \tilde{x}_{2,3}(t)&\mapsto x_{2,3}(t)\ , & \tilde{x}_{3,4}(t)&\mapsto x_{3,4}(t)\ , \\
\tilde{x}_{1,3}(t)&\mapsto x_{1,3}(t)\ , & \tilde{x}_{2,4}(t)&\mapsto x_{2,4}(t)\ , & \tilde{x}_{1,4}(t)&\mapsto x_{1,4}(t)\ .
\end{align*}\item We identify $\tilde{U}_+$ and $U_{3,4}U_{4,5}U_{5,6}$ via
\begin{align*}
\tilde{x}_{1,2}(t)&\mapsto x_{5,6}(t)\ , & \tilde{x}_{2,3}(t)&\mapsto x_{4,5}(-t)\ , & \tilde{x}_{3,4}(t)&\mapsto x_{3,4}(t)\ , \\
\tilde{x}_{1,3}(t)&\mapsto x_{4,6}(t)\ , & \tilde{x}_{2,4}(t)&\mapsto x_{3,5}(t)\ , & \tilde{x}_{1,4}(t)&\mapsto x_{3,6}(-t)\ .
\end{align*}
\end{enumerate}
We have
\begin{align*}
\forall\ t\in \EE:\qquad \alpha\big(\tilde{x}_{i,j}(t)\big)^\sigma=\alpha\big(\tilde{x}_{i,j}(t)^{\tilde{\sigma}}\big)
\end{align*}
for all $1\leq i<j\leq 3$ and therefore
$$\alpha\big(\Fix(\bar{\sigma})\big)=\Fix_{U_aU_bU_c}(\sigma)\ .$$
Thus we may identify $\Fix(\bar{\sigma})$ and $\Fix_{U_aU_bU_c}(\sigma)$ via
\begin{enumerate}[label=(\roman*)]
\item $\tilde{y}(s,t):=\tilde{x}_{1,3}(s)\tilde{x}_{1,4}(\beta\bar{t})\tilde{x}_{2,3}(t)\tilde{x}_{2,4}(\bar{s})\mapsto x_{1,3}(s)x_{1,4}(\beta\bar{t})x_{2,3}(t)x_{2,4}(\bar{s})=y^a(s,t)$.
\item $\tilde{y}(s,t):=\tilde{x}_{1,3}(s)\tilde{x}_{1,4}(\beta\bar{t})\tilde{x}_{2,3}(t)\tilde{x}_{2,4}(\bar{s})\mapsto x_{4,6}(s)x_{3,6}(-\beta\bar{t})x_{4,5}(-t)x_{3,5}(\bar{s})=y^b(\bar{s},-t)$.
\end{enumerate}
So if we have two copies $\BBB_1,\BBB_2$ of $\BBB$ and if we identify $\tilde{\BBB}$ with subgroups of $\BBB_1,\BBB_2$ as above, we have:
\begin{itemize}
\item $y_1^a(s,t)=y^a_2(s,t)$ if both the identifications are of type (i). Moreover, we have 
\begin{align*} \big(y^b_{2}(\HH^o),y_{2}^{ab}(\HH^o),y^a_{2}(\HH^o)\big)\cong\TTT(\HH^{o})\ , &&\big(y_{1}^a(\HH),y_{1}^{ab}(\HH),y_{1}^b(\HH)\big)\cong\TTT(\HH) \end{align*}
and therefore $\gamma_{(2,1)}=\id^o:\HH^o\to \HH$.
\item $y_1^b(s,t)=y_2^b(s,t)$ if both the identifications are of type (ii). Moreover, we have 
\begin{align*} \big(y_1^a(\HH),y_1^{ab}(\HH),y_1^b(\HH)\big)\cong\TTT(\HH)\ ,&&\big(y_2^b(\HH^o),y_2^{ab}(\HH^o),y_2^a(\HH^o)\big)\cong\TTT(\HH^{o})\end{align*}
and therefore $\gamma_{(1,2)}=\id^o:\HH\to \HH^o$.
\item $y_1^a(s,t)=y_2^b(\bar{s},-t)$ if the identifications are of different type. Moreover, we have \begin{align*} \big(y_2^a(\HH),y_2^{ab}(\HH),y_2^b(\HH)\big)\cong\TTT(\HH)\cong\big(y_1^a(\HH),y_1^{ab}(\HH),y_1^b(\HH)\big)\end{align*}
and therefore $\gamma_{(2,1)}=\sigma_s:\HH\to \HH$.
\end{itemize}\qed
\end{bew}

\begin{bem}
We obtain integrable mixed foundations over quaternions by glueing integrable positive foundations in a suitable way.
\end{bem}

\newglossaryentry{mf}{type=foundations,name={\ensuremath{\FFF(\GGG^P,\HH)}},description=mixed foundation with respect to the graph $\GGG^P$ and the quaternion division algebra $\HH$, sort=simply laced}

\begin{satz}\label{62}
Let $\GGG$ be a graph over $I$ and let $P$ be a collection of finite complete subgraphs such that
\begin{align*}
\bigcup_{\PPP\in P} \PPP=\GGG\ , && \forall\ \PPP_1\neq\PPP_2\in P: |V(\PPP_1)\cap V(\PPP_2)|\leq 1\ , && \forall\ i\in I: |\{ \PPP\in P\mid i\in V(\PPP)\}|\leq 2
\end{align*}
and such that the graph $\GGG^P$ with vertex set $P$ and
$$ \{ \PPP_1,\PPP_2\}\in E(\GGG^P)\ :\Leftrightarrow\ |V(\PPP_1)\cap V(\PPP_2)|=1\ $$
is a tree. Then there exists an integrable mixed foundation $\FFF:=\gls{mf}$ over $\HH$ such that:
\begin{itemize}
\item $\GGG_F=\GGG$,
\item The $V(\PPP)$-residue is isomorphic to $\PPP_{|V(\PPP)|}^+(\HH)$ for each $\PPP\in P$.
\item Given $\PPP_1,\PPP_2\in P$ such that $\PPP_1\cap\PPP_2=\{j\}$, we have
$$\forall\ (i,k)\in \big(V(\PPP_1)\sm\{j\}\big)\times \big(V(\PPP_2)\sm\{j\}\big):\qquad\gamma_{(i,j,k)}\in\{ \id,\sigma_s^o \}\ ,$$
i.e.,  these glueings are negative.
\end{itemize}
\end{satz}

\begin{bew}
Given a graph $\PPP\in P$, there is an integrable positive foundation $$\PPP_{|V(\PPP)|}^+(\HH)$$ which can be realized as fix foundation of a foundation $\BBB_\PPP$. We likewise want to realize the desired foundation $\FFF$ as a fixed point foundation. For this purpose, we will connect the foundations $\BBB_\PPP$ in a suitable way. Since $\GGG^P$ is a tree and $|\{ \PPP\in P\mid i\in V(\PPP)\}|\leq 2$ for each vertex $i\in I$, it suffices to prove the assertion for $|P|=2$.

Let $P=\{\PPP_1,\PPP_2\}$, let $V(\PPP_1)\cap V(\PPP_2)=\{j\}$, for $\lambda=1,2$ let $n_\lambda:=|V(\PPP_\lambda)|$, let $j_\lambda$ be the copy of $j$ in $\PPP_\lambda$ and let $\tilde{\BBB}_{j_\lambda}$ be the corresponding copy of $\tilde{\BBB}$ in the construction of $\BBB_{\PPP_\lambda}$ in theorem \ref{40}. We identify $\tilde{\BBB}_{j_1}$ and $\tilde{\BBB}_{j_2}$ via
\begin{align*}
\tilde{x}_{1,2}^{j_1}(t)\mapsto \tilde{x}_{3,4}^{j_2}(t)\ , && \tilde{x}_{2,3}^{j_1}(t)\mapsto \tilde{x}_{2,3}^{j_2}(-t)\ , && \tilde{x}_{3,4}^{j_1}(t)\mapsto\tilde{x}_{1,2}^{j_2}(t)\ .
\end{align*}
If $\BBB_{\{i,j_1\}}$ and $\BBB_{\{j_2,k\}}$ are copies of $\BBB$ in the construction of $\BBB_{\PPP_1}$ and $\BBB_{\PPP_2}$, respectively, we have the following possibilities:
\begin{itemize}
\item $y^{a}_{j_1}(s,t)=y_{j_2}^{a}(\bar{s},-t)$ if both the identifications are of type (i). Moreover, we have 
\begin{align*} \big(y_{j_1}^b(\HH^o),y^{ab}_{j_1}(\HH^o),y^{a}_{j_1}(\HH^o)\big)\cong\TTT(\HH^{o})\ , &&\big(y_{j_2}^a(\HH),y_{j_2}^{ab}(\HH),y_{j_2}^b(\HH)\big)\cong\TTT(\HH) \end{align*}
and therefore $\gamma_{(i,j,k)}=\sigma_s^o:\HH^o\to \HH$.
\item $y_{j_1}^b(s,t)=y_{j_2}^b(\bar{s},-t)$ if both the identifications are of type (ii). Moreover, we have 
\begin{align*} \big(y_{j_1}^a(\HH),y_{j_1}^{ab}(\HH),y_{j_1}^b(\HH)\big)\cong\TTT(\HH)\ , &&\big(y_{j_2}^b(\HH^o),y_{j_2}^{ab}(\HH^o),y_{j_2}^a(\HH^o)\big)\cong\TTT(\HH^o) \end{align*}
and therefore $\gamma_{(i,j,k)}=\sigma_s^o:\HH\to \HH^o$.
\item $y_{j_1}^b(s,t)=y_{j_2}^a(s,t)$ if both the identifications are of different type. Moreover, we have 
\begin{align*} \big(y_{j_1}^a(\HH),y_{j_1}^{ab}(\HH),y_{j_1}^b(\HH)\big)\cong\TTT(\HH)\cong\big(y_{j_2}^a(\HH),y_{j_2}^{ab}(\HH),y_{j_2}^b(\HH)\big)
\end{align*}
and therefore $\gamma_{(i,j,k)}=\id:\HH\to \HH$.
\end{itemize}
By construction, we have
\begin{align*} \GGG_F=\PPP_1\cup \PPP_2=\GGG\ .\end{align*}
\qed
\end{bew}

\newpage

\chapter{Triangle Foundations}
The smallest building bricks of foundations involving glueings are $A_3$- and $\tilde{A}_2$-residues.  Since there are $\tilde{A}_2$-foundations which do not have an integrable cover, it is not enough to consider $A_3$-residues. We will show that the only $\tilde{A}_2$-residues without an integrable cover are those constructed in chapter \ref{150}.

\begin{no} Throughout this chapter, $\BBB$ is of type $\tilde{A}_2$  with foundation
$$\FFF=\{\TTT(\AA_{(1,2)}),\TTT(\AA_{(2,3)}),\TTT(\AA_{(3,1)}),\gamma_2:=\gamma_{(1,2,3)},\gamma_3:=\gamma_{(2,3,1)},\gamma_1:=\gamma_{(3,1,2)}\}\ .$$
We denote its building at infinity by $\BBB^\infty$. 
\end{no}

\begin{bem}
Since $\BBB$ is a Moufang twin building, it is a Bruhat-Tits building by \cite{VV}, i.e., $\BBB^\infty$ is a Moufang triangle $\TTT(\AA^\infty)$ for some alternative division ring $\AA^\infty$. As a consequence, we may apply the results of \cite{W}.
\end{bem}

\section{The Defining Field}
First of all we show that each $A_2$-residue is isomorphic to the same Moufang triangle $\TTT(\AA)$, up to opposition. Then we prove that there is an embedding $\AA\hookrightarrow\AA^\infty$.

\begin{prop}\label{23} There is an automorphism $\alpha\in\Aut(\BBB)$ inducing $1\mapsto 2\mapsto 3\mapsto 1$ on $\GGG_F$.
\end{prop}

\begin{bew}
By corollary (18.15) of \cite{W}, $\Aut(\BBB)$ acts transitively on the set of gems, hence on the vertices and thus on $I$. As $|I|=3$ we are done.\qed
\end{bew}

\begin{kor}\label{65}
There is a $(1\ 2\ 3)$-automorphism of $\FFF$.
\end{kor}

\begin{bew}
Let $\UUU:=\UUU(\BBB,F,\Sigma,c)$ be a root group system such that $\FFF=\FFF(\UUU,\Lambda)$ for some parameter system $\Lambda$ for $\UUU$ and let $\alpha\in\Aut(\BBB)$ be as in proposition \ref{23}. Then $\alpha$ induces a $(1\ 2\ 3)$-isomorphism from $\UUU$ to $\UUU\big(\BBB,F,\alpha(\Sigma),\alpha(c)\big)$ which is specially isomorphic to $\UUU$ by theorem \ref{64}. By proposition \ref{71}, there is a $(1\ 2\ 3)$-automorphism of $\UUU$ which induces a $(1\ 2\ 3)$-automorphism of $\FFF$.\qed
\end{bew}

\begin{kor}\label{63}
We have $\AA:=\AA_{(1,2)}\cong \AA_{(2,3)}\cong \AA_{(3,1)}$.
\end{kor}

\begin{bew}
By corollary \ref{65}, we have
$$\TTT(\AA_{(1,2)})\cong\TTT(\AA_{(2,3)})\cong \TTT(\AA_{(3,1)})\ ,$$
thus the claim results from (35.6) of \cite{TW}.\qed
\end{bew}

\begin{satz}\label{500}
Let $\tilde{\FFF}$ be an integrable foundation. Then there is an alternative division ring $\AA$ such that
\begin{align*}
\forall\ (i,j)\in A(\tilde{F}):\qquad \tilde{\AA}_{(i,j)}\cong \AA\ \vee\ \tilde{\AA}_{(i,j)}\cong \AA^o\ .
\end{align*}
\end{satz}

\begin{bew}
By corollary \ref{34}  and corollary \ref{63}, this is true for each irr. rank 3 residue. The assertion results from an easy induction, starting with $\AA:=\tilde{\AA}_{(a,b)}$ for an arbitrary edge $(a,b)\in A(\tilde{F})$.
\qed
\end{bew}

\newpage
\begin{de}\label{501} Given an integrable foundation $\tilde{\FFF}$, we call the alternative division ring $\tilde{\AA}$ of theorem \ref{500} the \textit{defining field for $\mathit{\tilde{\FFF}}$}\index{defining field}, and we say that $\mathit{\tilde{\FFF}}$ is defined over $\tilde{\AA}$. By (35.6) of \cite{TW}, it is unique up to (anti-)isomorphism.
\end{de}

\begin{no}
Throughout the rest of this part, given an arbitrary integrable foundation $\tilde{\FFF}$, the alternative division ring $\tilde{\AA}$ always denotes its defining field.
\end{no}

\begin{satz}\label{41}
There is an embedding $\sigma:\AA\hookrightarrow \AA^\infty$.
\end{satz}

\begin{bew}
Let $$\RRR_0\cong \TTT(\AA)=\{x_1(\AA),x_2(\AA),x_3(\AA)\}$$ be the gem as in (18.1) of \cite{W}. As $\RRR_0$ is a special residue, it is a subbuilding of $$\BBB^\infty\cong \TTT(\AA^\infty)=\{\tilde{x}_1(\AA^\infty),\tilde{x}_2(\AA^\infty),\tilde{x}_3(\AA^\infty)\}$$ by theorem (6.3) of \cite{MV}. Given $t\in\AA$ and $i\in\{1,2,3\}$, there is a unique element $\tilde{x}_i\big(\sigma_i(t)\big)$ inducing $x_i(t)$, cf. proposition (29.61) of \cite{W}. By lemma \ref{362}, we may reparametrize $\TTT(\AA^\infty)$ such that
\begin{align*}
\sigma_1(1_\AA)=1_{\AA^\infty}\ , && \sigma_3(1_\AA)=1_{\AA^\infty}\ .
\end{align*}
Given $s,t\in\AA$ and $i\in \{1,2,3\}$, we have\footnote{We replace $\tilde{x}_i(\sigma_i(t))_{|\RRR_0}=x_i(t)$ by $\tilde{x}_i(\sigma_i(t))=x_i(t)$.}
\begin{align}\tilde{x}_i\big(\sigma_i(s+t)\big)=x_i(s+t)=x_i(s)x_i(t)=\tilde{x}_i\big(\sigma_i(s)\big)\tilde{x}_i\big(\sigma_i(t)\big)=\tilde{x}_i\big(\sigma_i(s)+\sigma_i(t)\big)\label{73}\end{align}
and
\begin{align}
\tilde{x}_2\big(\sigma_2(st)\big)&=x_2(st)=[x_1(s),x_3(t)]=\big[\tilde{x}_1\big(\sigma_1(s)\big),\tilde{x}_3\big(\sigma_3(t)\big)\big]=\tilde{x}_2\big(\sigma_1(s)\sigma_3(t)\big)\ .\label{72}
\end{align}
Putting $s:=1_\AA$ and $t:=1_\AA$ in equation \eqref{72} shows that
$$\sigma:=\sigma_3=\sigma_2=\sigma_1\ ,$$
hence
\begin{align*}
\forall\ s,t\in \AA: &&\sigma(s+t)\overset{\eqref{73}}{=}\sigma(s)+\sigma(t)\ , && \sigma(st)\overset{\eqref{72}}{=}\sigma(s)\sigma(t)\ .
\end{align*}
\qed
\end{bew}

\begin{bem}\label{74}\ 
\begin{enumerate}[label=(\alph*)]
\item The root group valuation $\phi:=\phi_{\RRR_0}$ with respect to $\RRR_0$ as in definition (13.8) of \cite{W} induces a discrete valuation $\nu$ of $\AA^\infty$. As a consequence, we have
$$\sigma(\AA)\subseteq \{ k\in \AA^\infty\mid \nu(k)=0 \}\ .$$
Moreover, $\sigma(\AA)$ is a set of representatives for the residue field $\bar{\AA}^\infty:=\OOO/m$, where $\OOO$ is the valuation ring of $\AA^\infty$ and $m$ is its unique maximal ideal.

In particular, $\tilde{x}_i\big(\sigma(\AA)\big)\subseteq U_{i,0}$ is a set of representatives for $\bar{U}_{i,0}$ and, given a uniformizer $\pi \in \AA^\infty$, $\tilde{x}_i(\sigma(\AA)\pi)\subseteq U_{i,1}$ and $\tilde{x}_i\big(\pi \sigma(\AA)\big)\subseteq U_{i,1}$ are sets of representatives for $\bar{U}_{i,1}$, where $U_{i,k}$ is defined as in definition (3.21) of \cite{W} and $\bar{U}_{i,k}$ is defined as in (18.21) of \cite{W}.
\item By the results of \cite{W}, we may suppose $\AA^\infty$ to be complete with respect to the valuation $\nu$ and $\nu(\AA^\infty)=\ZZ$. By lemma \ref{149}, $Z(\AA^\infty)$ is complete with respect to the induced valuation. As a consequence, we may apply the results of \cite{P} if $\AA^\infty$ is an octonion division algebra.
\item For brevity, we will write $\tilde{x}_i(t)$ instead of $\tilde{x}_i\big(\sigma(t)\big)$ in the following.
\end{enumerate}
\end{bem}

\newpage

\section{Triangles over Octonions}\label{153}

At this point, as well as in §\ref{151} and §\ref{152}, we heavily make use of the theory of affine buildings developed in \cite{W}. The main point is the fact that a parametrization for the building at infinity induces a parametrization for a given root group system. As uniformizers play a central role, we need some results of \cite{P}.

\begin{no} Throughout this paragraph, $\OO:=\AA$ is an octonion division algebra. As a consequence, $\OO^\infty$ is an octonion division algebra. \end{no}

\begin{prop}
We have $\nu\big(Z(\OO^\infty)\big)=\ZZ$. As a consequence, there is a uniformizer $\pi\in Z(\OO^\infty)$.
\end{prop}

\begin{bew}
As we have $\OO=\bar{\OO}^\infty$ by remark \ref{74}, this is a consequence of proposition 2 in \cite{P}.\qed
\end{bew}

\begin{satz}\label{94} The foundation $\FFF$ is specially  isomorphic to the foundation $\tilde{\AAA}_2(\OO)$.
\end{satz}

\begin{bew}
Let $( x_1,x_2,x_3)$ be a parametrization for $\TTT(\OO^\infty)$. Then the maps
\begin{align*} x_{(i,j)}^i&:\OO\to x_1(\OO^\infty)\ , && t\mapsto x_1(t)\in {U}_{1,0}\ ,\\
x_{(i,j)}^{ij}&:\OO\to x_2(\OO^\infty)\ , && t\mapsto x_2(t)\in {U}_{2,0}\ ,\\
x_{(i,j)}^{j}&:\OO\to x_3(\OO^\infty)\ , && t\mapsto x_3(t)\in {U}_{3,0}\end{align*}
yield a parametrization for the gem $\RRR_0=\RRR_{\{i,j\}}$, cf. remark \ref{74}. Moreover, the maps
\begin{align*} 
x_{(j,k)}^j&:\OO\to x_3(\OO^\infty)\ , &&t\mapsto x_3(t)\in {U}_{3,0}\ ,\\
x_{(j,k)}^{jk}&:\OO \to x_{-1}(\OO^\infty)\ , &&t\mapsto x_{-1}( t\pi)\in {U}_{-1,1}\ , \\
x_{(j,k)}^k&:\OO \to x_{-2}(\OO^\infty)\ , &&t\mapsto x_{-2}( t\pi)\in {U}_{-2,1}
\end{align*}
and
\begin{align*} 
x_{(k,i)}^k&:\OO\to x_{-2}(\OO^\infty)\ , &&t\mapsto x_{-2}(\pi t)\in {U}_{-2,1}\ ,\\
x_{(k,i)}^{ki}&:\OO \to x_{-3}(\OO^\infty)\ , &&t\mapsto x_{-3}(\pi t)\in {U}_{-3,1}\ , \\
x_{(k,i)}^i&:\OO \to x_{1}(\OO^\infty)\ , &&t\mapsto x_{1}( t)\in {U}_{1,0}
\end{align*}
yield parametrizations for two gems $\RRR_{\{j,k\}}$ and $\RRR_{(k,i)}$ at distance one to each other and to $\RRR_{\{i,j\}}$. By construction, this parameter system $\Lambda$ parametrizes a root group system $\UUU:=\UUU(\BBB, F,\Sigma,c)$. As 
\begin{align*} x_{(i,j)}^j(t)&=x_3(t)=x_{(j,k)}^j(t)\ , \\
x_{(j,k)}^k(t)&=x_{-2}(t\pi)=x_{-2}(\pi t)=x_{(k,i)}^k(t)\ , \\
x_{(k,i)}^i(t)&=x_1(t)=x_{(i,j)}^i(t)
\end{align*}
it follows that 
$$\gamma_{(i,j,k)}=\gamma_{(j,k,i)}=\gamma_{(k,i,j)}=\id_\OO\ .$$
Therefore, the resulting foundation
$$ \tilde{\FFF}:=\FFF(\UUU,\Lambda) $$
is the canonical foundation $\tilde{\AAA}_2(\OO)$. Finally we have $\FFF\cong \tilde{\FFF}=\tilde{\AAA}_2(\OO)$ by theorem \ref{64}.\qed
\end{bew}

\newpage
\section{Non-Existence of Tetrahedrons over Octonions}
Before we continue the examination of triangle foundations, we use the results of §\ref{153} to prove that there are no further integrable foundations over octonions. To apply them, we have to investigate the structure of $\Aut_J(\OO)$. In particular, we construct a certain subgroup $\Gamma\leq \Aut_J(\OO)$.

\begin{no}\ 
\begin{itemize}
\item Given a quaternion subalgebra $\HH$, $e\in\HH^\bot$ and $w,p\in \HH$, we set
\begin{align*} 
\psi_{(\HH,e,w)}&:\OO\to \OO,\ x+e\cdot y\mapsto x+e\cdot w^{-1}yw\ , \\
\phi_{(\HH,e,w,p)}&:\OO\to \OO,\ x+e\cdot y\mapsto w^{-1}xw+e\cdot w^{-1}ywp\ . \end{align*}
\item We set
\begin{align*} \Psi&:=\{ \psi_{(\HH,e,w)} \mid \HH\ \textrm{a quaternion subalgebra},\ e\in \HH^\bot,\ w\in\HH\}\ , \\
               \Phi&:=\{ \phi_{(\HH,e,w,p)} \mid \HH\ \textrm{a quaternion subalgebra},\ e\in \HH^\bot,\ w,p\in\HH,\ N(p)=1_{\OO}\}\ .
\end{align*}
\item Given $\psi_{(\HH,e,w)}\in \Psi$, we set
 $$\psi_{(\HH,e,w)}^o:\OO^o\to\OO^o,\ x+e\circ y\mapsto x+e\circ w^{-1}yw\ .$$
\item We set
\begin{align*} \Gamma:=\{ \psi\phi\mid \psi\in\Psi,\ \phi\in\Aut(\OO)\}\ , && \Gamma^o:=\{ \psi^o\phi\mid \psi\in \Psi,\ \phi\in \Aut(\OO^o) \}\ .\end{align*}
\end{itemize}
\end{no}

\begin{lemma}
Let $\psi:=\psi_{(\HH,e,w)}\in\Psi$. Then the following holds:
\begin{enumerate}[label=(\alph*)]
\item Given $s,t\in\OO$, we have
\begin{align*}
\psi(st)=\big(\psi(s)\cdot \psi(t)w\big)w^{-1}\ .
\end{align*}
\item $\psi\in \Aut_J(\OO)$.
\end{enumerate}
\end{lemma}

\begin{bewzwei}\ 
\begin{enumerate}[label=(\alph*)]
\item This is (20.24) of \cite{TW}. 
\item By the Moufang identities and the inverse properties, We have
\begin{align*} \psi(sts)&\overset{\mathclap{\textrm{(a)}}}{=}\big(\psi(s)\cdot \psi(ts)w\big)w^{-1}\overset{\mathclap{\textrm{(a)}}}{=}\big(\psi(s)\cdot \big[ (\psi(t)\cdot \psi(s)w)w^{-1}\big]w\big)w^{-1} \\
&=\big[\psi(s) \big(\psi(t)\cdot \psi(s)w\big)\big]w^{-1}=\big[\big(\psi(s)\psi(t)\cdot \psi(s)\big)w\big]w^{-1}=\psi(s)\psi(t)\psi(s)\ .
\end{align*}
for all $s,t\in\OO$.
\end{enumerate}\qed
\end{bewzwei}

\begin{lemma}\label{92}
We have $\Phi\subseteq \Aut_{Z(\OO)}(\OO)$. Given an element $\tilde{\phi}\in \Aut_{Z(\OO)}(\OO)$ leaving a quaternion subalgebra $\HH$ invariant, there is an element $\phi\in \Phi$ such that $\tilde{\phi}=\phi$. 
\end{lemma}

\begin{bew}
This is section (2.1) of \cite{S}. Notice that our point of view is that of the opposite multiplication. \qed
\end{bew}

\begin{lemma}\label{97}
Given $\gamma=\psi_{(\HH,e,w)}\phi\in \Gamma$ and $v\in \OO^*$, let $\alpha_{v,\gamma}:\TTT(\OO)\to \TTT(\OO)$,
\begin{align*}
x^1(t)\mapsto x^1\big( v\cdot \gamma(t)\big)\ , && x^{12}(t)\mapsto x^{12}\big(v\cdot \gamma(t)w\cdot v\big)\ , && x^2(t)\mapsto x^2\big(\gamma(t)w\cdot v\big)
\end{align*}
and given $\gamma^o=\psi^o_{(\HH,e,w)}\phi\in \Gamma$ and $v\in \OO^*$, let $\alpha_{v,\gamma^o}^o:\TTT(\OO^o)\to \TTT(\OO^o)$,
\begin{align*}
x^2(t)\mapsto x^2\big( v\circ \gamma^o(t)\big)\ , && x^{12}(t)\mapsto x^{12}\big(v\circ \gamma^o(t)w\circ v\big)\ , && x^1(t)\mapsto x^1\big(\gamma^o(t)w\circ v\big)\ .
\end{align*}
Then we have
$$\{ \alpha_{v,\gamma} \mid v\in \OO^*,\ \gamma\in\Gamma\}=\Aut \big(\TTT(\OO)\big)=\Aut\big(\TTT(\OO^o)\big)=\{ \alpha_{v,\gamma^o}^o \mid v\in \OO^*,\ \gamma^o\in \Gamma^o\}\ .$$
\end{lemma}

\begin{bew}
This holds by the proof of (37.12) in \cite{TW}.\qed
\end{bew}

\begin{lemma}\label{87}
We have $\Gamma=\Gamma^o$.
\end{lemma}

\begin{bew}
Given $\psi:=\psi_{(\HH,e,w)}\in \Psi$, we have
\begin{align*}
\psi^o(x+e\cdot {y})&=\psi^o(x+\bar{y}\cdot e)=\psi^o(x+e\circ \bar{y})\\
&=x+e\circ w^{-1}\bar{y}w=x+ w\bar{y}w^{-1}\cdot e=x+e\cdot \overline{w\bar{y}w^{-1}}=x+e\cdot wyw^{-1}
\end{align*}
for all $x,y\in\HH$ and therefore
$$\psi^o_{(\HH,e,w)}=\psi_{(\HH,e,w^{-1})}\ .$$\qed
\end{bew}

\begin{lemma}\label{119}
The set $\Gamma$ is a subgroup of $\Aut_J (\OO)$.
\end{lemma}

\begin{bew}
Given $\alpha_{v,\gamma}\in \Aut\big(\TTT(\OO)\big)$, we have
\begin{align} \alpha_{v,\gamma}^1(1_\OO)=v\cdot \gamma(1_\OO)=v\ .\label{86}\end{align}
\begin{itemize}
\item We have
$$\id_\OO=\psi_{(\HH,e,1_\OO)}\circ\id_\OO\in \Gamma\ .$$
\item Let $\gamma\in \Gamma$.  By equation \eqref{86} and lemma \ref{97}, there is an element $\rho\in \Gamma$ such that
$$\alpha_{1_\OO,\gamma}^{-1}=\alpha_{1_\OO,\rho}\ .$$
We obtain
$$\gamma^{-1}=(\alpha_{1_\OO,\gamma}^1)^{-1} =\alpha_{1_\OO,\rho}^1=\rho\in \Gamma \ .$$
\item Let $\gamma_1,\gamma_2\in \Gamma$. By equation \eqref{86} and lemma \ref{97}, there is an element $\rho\in \Gamma$ such that
$$\alpha_{1_\OO,\gamma_1}\circ \alpha_{1_\OO,\gamma_2}=\alpha_{1_\OO,\rho}\ .$$
We obtain
$$\gamma_1\gamma_2=\alpha_{1_\OO,\gamma_1}^1 \alpha_{1_\OO,\gamma_2}^1=\alpha_{1_\OO,\rho}^1=\rho\in \Gamma \ .$$
\end{itemize}\qed
\end{bew}

\begin{lemma}\label{466}
We have $\Psi\cup \Phi\subseteq \Gamma$.
\end{lemma}

\begin{bew}
Given a quaternion subalgebra $\HH$ and $e\in \HH^\bot$, we have
\begin{align*}
\id_\OO\in \Aut(\OO)\ , && \id_\OO=\psi_{(\HH,e,1_\OO)}\in \Psi\ .
\end{align*}\qed
\end{bew}

\begin{lemma}\label{96}
We have $\gamma_w\in \Gamma$ for each $w\in \OO$.
\end{lemma}

\begin{bew}
By lemma \ref{93}, there is a quaternion subalgebra $\HH$ containing $w$. Let $e\in \HH^\bot$. Then we have
$$\gamma_w(x+e\cdot y)=w^{-1}xw+w^{-1}(e\cdot y)w=w^{-1}xw+e\cdot \bar{w}^{-1}wy$$
for all $x,y\in\HH$. Notice that
$$N(\bar{w}^{-1}w)=N(\bar{w})^{-1}N(w)=N(w)^{-1}N(w)=1_{\OO}\ ,$$
thus we obtain
\begin{align*}
 \gamma_w=\phi_{(\HH,e,w,\bar{w}^{-1}w)}\circ \psi_{(\HH,e,w^{-2}\bar{w})}\in \Gamma\ .
\end{align*}\qed
\end{bew}

\begin{lemma}
We have $\sigma_s\notin \Gamma$.
\end{lemma}

\begin{bew}
Let $\psi=\psi_{(\HH,e,w)}\in \Psi$ and $\phi\in \Aut(\OO)$ such that $\sigma_s=\psi\phi$. Then
$$\id_{\HH}=\psi_{|\HH}=\sigma_s\phi^{-1}_{|\HH}$$
is both negative and positive$\qquad\lightning$.\qed
\end{bew}

\begin{lemma}\label{98} Let
$$\FFF=\{\TTT(\AA_{(1,2)}):=\TTT(\AA_{(2,3)}):=\TTT(\AA_{(3,1)}):=\TTT(\OO),\ \gamma_2:=\gamma_{(1,2,3)},\gamma_3:=\gamma_{(2,3,1)},\gamma_1:=\gamma_{(3,1,2)}\}$$
be an integrable triangle foundation over $\OO$. Then we have $\gamma_i\in \Gamma$ for $i=1,2,3$.
\end{lemma}

\begin{bew}
We have $\FFF\cong \tilde{\AAA}_2(\OO)$ by theorem \ref{94}, thus theorem \ref{141} yields a reparametrization $\alpha$ for $\tilde{\AAA}_2(\OO)$ such that
$$\tilde{\AAA}_2(\OO)_\alpha=\FFF\ .$$
From definition \ref{142} it follows that
$$\gamma_{(i,j,k)}=(\alpha_{(j,k)}^j)^{-1}\circ \id_\OO\circ \alpha_{(i,j)}^j\ .$$
By lemma \ref{97} and lemma \ref{87}, there are $\gamma_{(i,j)},\gamma_{(j,k)}\in \Gamma$ and $w\in \OO$ such that
\begin{align*}
\alpha_{(i,j)}^j(t)=w\cdot \gamma_{(i,j)}(t)\ , && \alpha_{(j,k)}^j(t)=\gamma_{(j,k)}(t)\cdot w 
\end{align*}
and thus
\begin{align*}
(\alpha_{(j,k)}^j)^{-1}(t)= \gamma_{(j,k)}^{-1}(t\cdot w^{-1})\ , && (\alpha_{(j,k)}^j)^{-1}\circ \alpha_{(i,j)}^j(t)=(\gamma_{(j,k)}^j)^{-1}\big(w\cdot \gamma_{(i,j)}(t)\cdot w^{-1}\big) 
\end{align*}
for each $t\in \OO$. We finally obtain
$$\gamma_{(i,j,k)}=(\gamma_{(j,k)}^j)^{-1}\circ \gamma_{w^{-1}}\circ \gamma_{(i,j)}\in \Gamma$$
by lemma \ref{96}.
\qed
\end{bew}

\newpage
\begin{prop}\label{99}
Let $\tilde{\FFF}$ be a foundation over $\OO$ such that $\GGG_{\tilde{F}}$ is a tetrahedron. Then $\tilde{\FFF}$ is not integrable.
\end{prop}

\begin{bew}
Assume that $\tilde{\FFF}$ is integrable. Then each rank 3 residue is integrable, hence specially isomorphic to $\tilde{\AAA}_2(\OO)$ by theorem \ref{94}. By remark \ref{143}, we may assume that the $\{1,2,3\}$-residue is $\tilde{\AAA}_2(\OO)$. Moreover, we may assume 
$$\TTT(\tilde{\AA}_{(3,4)})=\TTT({\tilde{\AA}_{(4,1)}})=\TTT(\tilde{\AA}_{(2,4)})=\TTT(\OO)\ .$$
Now we are in the situation of lemma \ref{98}, thus we have $\tilde{\gamma}_{(1,2,4)}\in \Gamma$ and
$$\tilde{\gamma}_{(4,2,3)}=\id^o\circ \id^o\circ \tilde{\gamma}_{(1,2,4)}^{-1}\circ \id^o=\tilde{\gamma}_{(1,2,4)}^{-1}\circ \id^o\ .$$
We extend the reparametrization 
$$ \alpha_{(4,2)}:=\{ \OO,\sigma^o_s,\sigma^o_s,\sigma^o_s\}$$
to a reparametrization $\alpha$ for $\tilde{\FFF}$ and obtain an integrable foundation
$$\FFF:=\tilde{\FFF}_\alpha$$
satisfying
$$\gamma_{(4,2,3)}=\id_\OO\circ \tilde{\gamma}_{(4,2,3)}\circ \sigma^o_s=\tilde{\gamma}_{(1,2,4)}^{-1}\circ \sigma_s\ .$$
As we have
$$\TTT(\AA_{(4,2)})=\TTT(\AA_{(2,3)})=\TTT(\AA_{(3,4)})=\TTT(\OO)\ ,$$
it follows that
\begin{align*}
\tilde{\gamma}_{(1,2,4)}^{-1}\circ \sigma_s=\gamma_{(4,2,3)}\in \Gamma\ , && \sigma_s\in \Gamma\qquad \lightning\ .
\end{align*}\qed
\end{bew}

\begin{satz}\label{444}
A foundation $\FFF$ over an octonion division algebra $\OO$ is integrable if and only if we have $\FFF\cong \AAA_2(\OO)$ or $\FFF\cong \tilde{\AAA}_2(\OO)$.
\end{satz}

\begin{bew}
As $\GGG_F$ is complete and as each residue is integrable, proposition \ref{99} implies $|V(F)|\leq 3$. If $|V(F)|=2$, we have $\FFF\cong \TTT(\OO)$, and if $|V(F)|=3$, we have $\FFF\cong \tilde{\AAA}_2(\OO)$ by theorem \ref{94}. Finally, $\tilde{\AAA}_2(\OO)$ is integrable by theorem \ref{445}.\qed
\end{bew}

\section{Triangles over Skew-Fields}\label{151}

As we are done with the octonion case, we next deal with positive foundations over skew-fields. The first step is to show that a triangle foundation $\FFF$ over a skew-field $\AA$ is negative if $\AA^\infty$ is a skew-field. Thus $\AA^\infty$ is necessarily an octonion division algebra if $\FFF$ is positive. Then we prove that $\AA$ is a quaternion division algebra, using the fact that $\AA$ embeds into $\AA^\infty$. At last we obtain a parametrization for a root group sequence via the building at infinity.

\begin{no} Throughout this paragraph, $\AA$ is a skew-field. \end{no}

\begin{lemma} \label{66}
Let $\TTT(\tilde{\AA})$ and $\TTT(\hat{\AA})$ be isomorphic Moufang triangles over skew-fields and let $\alpha=(\alpha_1,\alpha_2,\alpha_3):\TTT(\tilde{\AA})\to \TTT(\hat{\AA})$ be an isomorphism. Then there are $a,b\in \hat{\AA}$ and an isomorphism $\phi:\tilde{\AA}\to\hat{\AA}$ of skew-fields such that
$$\alpha=( \lambda_a\phi, \lambda_a\rho_b\phi, \rho_b \phi)\ .$$
\end{lemma}

\begin{bew}
Let 
\begin{align*}
a:=\alpha_1(1_{\tilde{\AA}})\ , && b:=\alpha_3(1_{\tilde{\AA}})\ .
\end{align*}
Since $\hat{\AA}$ is associative, the map
$$\tilde{\alpha}:=\big(\lambda_{a^{-1}}\alpha_1,\lambda_{a^{-1}}\rho_{b^{-1}}\alpha_2,\rho_{b^{-1}}\big):\TTT(\tilde{\AA})\to \TTT(\hat{\AA})$$
is an isomorphism. Moreover, we have
\begin{align*}
\lambda_{a^{-1}}\alpha_1(1_{\tilde{\AA}})=a^{-1}a=1_{\hat{\AA}}\ , && \rho_{b^{-1}}\alpha_3(1_{\tilde{\AA}})=bb^{-1}=1_{\hat{\AA}}\ .
\end{align*}
By (35.23) of \cite{TW} therefore, there is an isomorphism $\phi:\tilde{\AA}\to\hat{\AA}$ of skew-fields such that
$$\tilde{\alpha}=(\phi,\phi,\phi)\ .$$
Thus we have
$$\alpha=(\lambda_a\phi,\lambda_a\rho_b\phi,\rho_b\phi)\ .$$
\qed
\end{bew}

\begin{lemma}\label{67}
Let $\tilde{\FFF},\hat{\FFF}$ be isomorphic foundations over a skew-field, let $\alpha:\tilde{\FFF}\to \hat{\FFF}$ be an isomorphism and let $(i,j,k)\in G(F)$. If $\tilde{\gamma}_{(i,j,k)}$ is negative, then $\hat{\gamma}:=\hat{\gamma}_{(\pi(1),\pi(2),\pi(3))}$ is negative.
\end{lemma}

\begin{bew}
By lemma \ref{66}, there are isomorphisms $\phi_i:\tilde{\AA}_{(i,j)}\to \hat{\AA}_{(\pi(i),\pi(j))}$ and $\phi_k:\tilde{\AA}_{(j,k)}\to \hat{\AA}_{(\pi(j),\pi(k))}$ of skew-fields and elements $a_i\in \hat{\AA}_{(\pi(i),\pi(j))}$ and $a_k\in \hat{\AA}_{(\pi(j),\pi(k)k)}$ such that
\begin{align*}
\alpha_{(i,j)}^j=\rho_{a_i}\circ\phi_i\ , && \alpha_{(j,k)}^j=\lambda_{a_k}\circ \phi_k\ .
\end{align*}
By the definition of an isomorphism, we have
\begin{align*}
\hat{\gamma}=\alpha_{(j,k)}^j\circ \tilde{\gamma}_{(i,j,k)}\circ (\alpha_{(i,j)}^j)^{-1}\ .
\end{align*}
Combining these two facts implies that there are elements $a,b\in \hat{\AA}_{(\pi(j),\pi(k))}$ such that
\begin{align*}
\hat{\gamma}=\lambda_a\circ\rho_b\circ \phi_k\circ \tilde{\gamma}_{(i,j,k)}\circ \phi_i^{-1}\ .
\end{align*}
Moreover, we have
$$ab=\hat{\gamma}(1)=1$$
and therefore $b=a^{-1}$. It finally follows that
$$\hat{\gamma}=\gamma_b\circ \phi_k\circ \tilde{\gamma}_{(i,j,k)}\circ \phi_i^{-1}$$ is negative.
\qed
\end{bew}

\begin{kor} The following holds:
\begin{enumerate}[label=(\alph*)]
\item \label{68} The foundation $\FFF$ is either negative or positive.
\item \label{465} If a foundation isomorphic to $\FFF$ is negative, then $\FFF$ is negative.
\end{enumerate}
\end{kor}

\begin{bewzwei}\ 
\begin{enumerate}[label=(\alph*)]
\item This results from lemma \ref{67} and corollary \ref{65}.
\item This is an immediate consequence of lemma \ref{67}.
\end{enumerate}\qed
\end{bewzwei}

\begin{satz}\label{70}
Let $\AA^\infty$ be a skew-field. Then $\FFF$ is negative.
\end{satz}

\newpage

\begin{bew}
Let $( x_1,x_2,x_3)$ be a parametrization for $\TTT(\AA^\infty)$ and let $\pi\in\AA^\infty$ be a uniformizer. Then the maps
\begin{align*} x_{(i,j)}^i&:\AA\to x_1(\AA^\infty)\ , && t\mapsto x_1(t)\in {U}_{1,0}\ ,\\
x_{(i,j)}^{ij}&:\AA\to x_2(\AA^\infty)\ , && t\mapsto x_2(t)\in {U}_{2,0}\ ,\\
x_{(i,j)}^{j}&:\AA\to x_3(\AA^\infty)\ , && t\mapsto x_3(t)\in {U}_{3,0}\end{align*}
yield a parametrization for the gem $\RRR_0=\RRR_{\{i,j\}}$, cf. remark \ref{74}.
Moreover, the maps
\begin{align*} 
x_{(j,k)}^j&:\AA\to x_3(\AA^\infty)\ , &&t\mapsto x_3(t)\in {U}_{3,0}\ ,\\
x_{(j,k)}^{jk}&:\AA \to x_{-1}(\AA^\infty)\ , &&t\mapsto x_{-1}( t\pi)\in {U}_{-1,1}\ , \\
x_{(j,k)}^k&:\AA \to x_{-2}(\AA^\infty)\ , &&t\mapsto x_{-2}( t\pi)\in {U}_{-2,1}
\end{align*}
yield a parametrization for a gem $\RRR_{\{j,k\}}$ at distance one to $\RRR_{\{i,j\}}$. By construction, these two parametrizations are part of a parameter system $\Lambda$ for a root group system $\UUU:=\UUU(\BBB, F,\Sigma,c)$. As we have
$$x_{(i,j)}^j(t)=x_3(t)=x_{(j,k)}^j(t)\ ,$$
it follows that $\gamma_{(i,j,k)}=\id_\AA$ is negative. Therefore, the resulting foundation
$$ \tilde{\FFF}:=\FFF(\UUU,\Lambda) $$
is negative, cf. corollary \ref{68} with $\tilde{\FFF}$ in place of $\FFF$. As we have $\FFF\cong\tilde{\FFF}$ by theorem \ref{64}, corollary \ref{465} finally shows that $\FFF$ is negative.\qed
\end{bew}

\section{Skew-Fields inside Octonions}\label{147}

As we want to show that $\AA$ is a quaternion division algebra if $\FFF$ is positive, we prove that each non-commutative division subring of an octonion division algebra is a quaternion division algebra. Then the claim results from the fact that $\AA$ embeds into $\AA^\infty$.

In §\ref{152}, we will need once again the results of \cite{P} to find a suitable uniformizer for a parametrization. In fact, we can choose a uniformizer $\pi\in \AA^\bot$, thus conjugating elements of $\AA$ by $\pi$ is equal to applying the standard involution.

\begin{lemma}\label{44}
Let $\OO$ be an octonion algebra and let $\DD$ be a non-commutative division subring. Then we have $$Z(\DD)\subseteq Z(\OO)=:\KK\ .$$
\end{lemma}

\begin{bew}
The octonion division algebra $\OO$ quadratic over $\KK$. Then $\tilde{\DD}:=\langle \DD\rangle_{\KK}$ is non-commutative and quadratic over $\KK\subseteq Z(\tilde{\DD})$, thus we have $\KK=Z(\tilde{\DD})$ by proposition \ref{121} and therefore $$Z(\DD)\subseteq Z(\tilde{\DD})=\KK\ .$$
\qed
\end{bew}

\begin{lemma}\label{46}
Let $\DD$ be a skew-field such that $\dim_{Z(\DD)} \DD<\infty$. Then there is an $n\in\NN^*$ such that $$\dim_{Z(\DD)} \DD=n^2\ .$$
\end{lemma}

\begin{bew} This results from the fact that $\DD$ is central simple over $Z(\DD)$.
\qed
\end{bew}

\begin{satz}\label{43}
Let $\OO$ be an octonion division algebra, let $\KK:=Z(\OO)$ and let $\DD$ be a non-commutative alternative division subring. Then the following holds:
\begin{enumerate}[label=(\alph*)]
\item $\DD\otimes_{Z(\DD)} \KK$ is isomorphic to a division subalgebra of $\OO$.
\item $\DD$ is a quaternion division algebra or an octonion division algebra.
\end{enumerate}
\end{satz}

\begin{bew}
The non-commutative division subalgebra $\tilde{\DD}:=\langle \DD\rangle_{\KK}\subseteq \OO$ is quadratic over $\KK\subseteq Z(\tilde{\DD})$ and therefore a quaternion division algebra or an octonion division algebra by proposition \ref{121}. 
\begin{enumerate}[label=(\alph*)]
\item By lemma \ref{44}, we have $Z(\DD)\subseteq \KK$, hence
$$\DD\otimes_{Z(\DD)} \KK$$
is a central simple algebra over $\KK$. By the universal property, there is an epimorphism
$$\pi:\DD\otimes_{Z(\DD)} \KK\to \tilde{\DD}$$
which is injective due to simplicity. 
\item Part (a) yields
$$\dim_{Z(\DD)} \DD=\dim_{\KK} (\DD\otimes_{Z(\DD)} \KK)=\dim_{\KK} \tilde{\DD}\in\{4,8\}\ .$$
In particular, each $Z(\DD)$-basis of $\DD$ is a $\KK$-basis of $\tilde{\DD}$.
\begin{itemize}
\item If $\DD$ is non-associative, then $\DD$ is an octonion division algebra and thus $\dim_{Z(\DD)} \DD=8$.
\item If $\DD$ is associative, then we have $$\dim_{Z(\DD)} \DD\in\{4,8\}\cap\{n^2\mid n\in\NN^*\}=4$$
by lemma \ref{46}. Let $x\in \DD\sm Z(\DD)$ and let $\{1_\OO,x,y,z\}$ be a $Z(\DD)$-basis of $\DD$. Then there are $\lambda_1,\ldots,\lambda_4\in Z(\DD)\subseteq \KK$ such that
$$x^{-1}=\lambda_1\cdot 1_\OO+\lambda_2\cdot x+\lambda_3\cdot y+\lambda_4\cdot z\ .$$
Since $\{1_\OO,x,y,z\}$ is a $\KK$-basis of $\tilde{\DD}$, it follows that
\begin{align*}
\lambda_1=N(x)^{-1}T(x)\ , && \lambda_2=-N(x)^{-1}\ , && \lambda_3=0_\OO=\lambda_4
\end{align*}
and therefore
\begin{align*}
N(x)=-\lambda_2^{-1}\in Z(\DD)\ , && T(x)=\lambda_1\cdot N(x)\in Z(\DD)\ .
\end{align*}
As a consequence, $\DD$ is quadratic over $Z(\DD)$ and therefore a quaternion division algebra by proposition \ref{121}.
\end{itemize}
\end{enumerate}
\qed
\end{bew}

\section{Positive Triangles over Skew-Fields}\label{152}

\begin{no} Throughout this paragraph, $\AA$ is a skew-field and $\FFF$ is positive. \end{no}

\begin{bem} By theorem \ref{70}, $\AA^\infty$ is an octonion division algebra. \end{bem}

\begin{lemma}\label{37}
The defining field $\AA$ is a quaternion division algebra.
\end{lemma}

\begin{bew}
We have $\AA\subseteq \AA^\infty$ by theorem \ref{41}, thus  the claim results from theorem \ref{43}.\qed
\end{bew}

\begin{no} We set $\HH:=\AA$, $\OO:=\AA^\infty$ and $\KK:=Z(\OO)$.
\end{no}

\begin{prop}\label{54}
There is a uniformizer $\pi\in \HH^\bot$.
\end{prop}

\newpage

\begin{bew}
As we have $\HH=\bar{\OO}$ by remark \ref{74}, proposition 2 of \cite{P} implies $\nu(\KK)=2\ZZ$. As a consequence, $\OO$ is ramified and the quaternion division algebra $\HH':=\langle \HH\rangle_{\KK}$
is an unramified composition algebra such that $\bar{\HH}'=\HH$. Theorem 2 of \cite{P} shows that
$$\OO\cong (\HH',\pi')$$
for some $\nu_{|\KK}$-uniformizer $\pi'$ of $\KK$, i.e., we have  $\nu(\pi')=2$. As a consequence, there is an element $\pi\in \HH'^\bot\subseteq \HH^\bot$ satisfying $N(\pi)=-\pi'$ and, therefore,
$$\nu(\pi)=\frac{\nu\big(N(\pi)\big)}{2}=\frac{\nu(-\pi')}{2}=\frac{\nu(\pi')}{2}=1\ .$$\qed
\end{bew}

\newglossaryentry{ph}{type=foundations,name={\ensuremath{\PPP_3^+(\HH)}},description=positive triangle with respect to the quaternion division algebra $\HH$, sort=simply laced}

\begin{satz}\label{24} The foundation $\FFF$ is specially isomorphic to the standard positive foundation
$$\gls{ph}:=\{\TTT(\AA_{(1,2)}^+):=\TTT(\AA_{(2,3)}^+):=\TTT(\AA_{(3,1)}^+):=\TTT(\HH),\gamma_1^+:=\gamma_2^+:=\gamma_3^+:=\sigma_s\}\ .$$
\end{satz}

\begin{bew}
By proposition \ref{54}, there is a uniformizer $\pi\in \HH^\bot$. Let $( x_1,x_2,x_3)$ be a parametrization for $\TTT(\OO)$. Then the maps
\begin{align*} x_{(i,j)}^i&:\HH\to x_1(\OO)\ , && t\mapsto x_1(t)\in {U}_{1,0}\ ,\\
x_{(i,j)}^{ij}&:\HH\to x_2(\OO)\ , && t\mapsto x_2(t)\in {U}_{2,0}\ ,\\
x_{(i,j)}^{j}&:\HH\to x_3(\OO)\ , && t\mapsto x_3(t)\in {U}_{3,0} \end{align*}
yield a parametrization for the gem $\RRR_0=\RRR_{\{i,j\}}$, cf. remark \ref{74}. Moreover, the maps
\begin{align*} 
x_{(j,k)}^j&:\HH\to x_3(\OO)\ , &&t\mapsto x_3(\bar{t})\in {U}_{3,0}\ ,\\
x_{(j,k)}^{jk}&:\HH \to x_{-1}(\OO)\ , &&t\mapsto x_{-1}(\bar{t}\pi)\in {U}_{-1,1}\ , \\
x_{(j,k)}^k&:\HH\to x_{-2}(\OO)\ , &&t\mapsto x_{-2}(\bar{t}\pi)\in {U}_{-2,1}
\end{align*}
and
\begin{align*} 
x_{(k,i)}^k&:\HH\to x_{-2}(\OO)\ , &&t\mapsto x_{-2}(\pi\bar{t})\in {U}_{-2,1}\ ,\\
x_{(k,i)}^{ki}&:\HH \to x_{-3}(\OO)\ , &&t\mapsto x_{-3}(\pi\bar{t})\in {U}_{-3,1}\ , \\
x_{(k,i)}^i&:\HH \to x_{1}(\OO)\ , &&t\mapsto x_{1}(\bar{t})\in {U}_{1,0}
\end{align*}
yield parametrizations for two gems $\RRR_{\{j,k\}}$ and $\RRR_{(k,i)}$ at distance one to each other and to $\RRR_{\{i,j\}}$, e.g., we have
$$[x_{(j,k)}^j(s),x_{(j,k)}^k(t)]=[x_3(\bar{s}),x_{-2}(\bar{t}\pi)]=x_{-1}\big(\bar{s}(\bar{t}\pi)\big)=x_{-1}\big((\bar{t}\bar{s})\pi\big)=x_{-1}\big( (\overline{st})\pi\big)=x_{(j,k)}^{jk}(st)$$
for all $s,t\in\HH$. By construction, this parameter system $\Lambda$ parametrizes a root group system $\UUU:=\UUU(\BBB, F,\Sigma,c)$. As we have
\begin{align*} x_{(i,j)}^j(t)&=x_3(t)=x_{(j,k)}^j(\bar{t})\ , \\
x_{(j,k)}^k(t)&=x_{-2}(\bar{t}\pi)=x_{-2}\big(\pi (\pi^{-1}\bar{t}\pi)\big)=x_{-2}(\pi t)=x_{(k,i)}^k(\bar{t})\ , \\
x_{(k,i)}^i(t)&=x_1(\bar{t})=x_{(i,j)}^i(\bar{t})\ ,
\end{align*}
it follows that 
$$\gamma_{(i,j,k)}=\gamma_{(j,k,i)}=\gamma_{(k,i,j)}=\sigma_s\ .$$
Therefore, the resulting foundation
$$ \tilde{\FFF}:=\FFF(\UUU,\Lambda) $$
is $\PPP_3^+(\HH)$. Finally $\FFF\cong \tilde{\FFF}=\PPP_3^+(\HH)$ by theorem \ref{64}.\qed
\end{bew}

\chapter{Positive Foundations over Skew-Fields}

Since each residue of an integrable foundation is itself integrable, a positive foundation is built up of positive triangle foundations. Their uniqueness enables us to show the uniqueness of positive foundations for a given quaternion division algebra $\HH$ and a given value of $|I|$.

\begin{no} Throughout this chapter, $\FFF$ is a positive foundation over $I:=\{1,\ldots,n\}$. Since each rank 3 residue is a positive triangle, $\HH:=\AA$ is a quaternion division algebra by lemma \ref{37}. Moreover, $\GGG_F$ is complete by corollary $\ref{144}$.
\end{no}

\begin{lemma}\label{75}
Let $\tilde{\FFF},\hat{\FFF}$ be isomorphic positive foundations over a skew-field, let $\alpha:\tilde{\FFF}\to \hat{\FFF}$ be an isomorphism and let $(i,j,k)\in G(F)$. By lemma \ref{66}, there are isomorphisms 
\begin{align*}
\phi_i:\tilde{\AA}_{(i,j)}\to \hat{\AA}_{(\pi(i),\pi(j))}\ , &&\phi_k:\tilde{\AA}_{(j,k)}\to \hat{\AA}_{(\pi(j),\pi(k))}\end{align*}
of skew-fields and elements $a_i\in \hat{\AA}_{(\pi(i),\pi(j))}$ and $a_k\in \hat{\AA}_{(\pi(j),\pi(k)k)}$ such that
\begin{align*}
\alpha_{(i,j)}^j=\rho_{a_i}\circ\phi_i\ , && \alpha_{(j,k)}^j=\lambda_{a_k}\circ \phi_k\ .
\end{align*}
Then we have
$$\hat{\gamma}:=\hat{\gamma}_{(\pi(i),\pi(j),\pi(k))}=\phi_k\circ \tilde{\gamma}_{(i,j,k)}\circ \phi_i^{-1}\ .$$
\end{lemma}

\begin{bew}
By the definition of an isomorphism, we have
\begin{align*}
\hat{\gamma}=\alpha_{(j,k)}^j\circ \tilde{\gamma}_{(i,j,k)}\circ (\alpha_{(i,j)}^j)^{-1}=\lambda_{a_k}\circ \phi_k\circ \tilde{\gamma}_{(i,j,k)}\circ \circ \phi_i^{-1}\circ \rho_{a_i}^{-1}\ .
\end{align*}
As $\tilde{\gamma}_{(i,j,k)}$ is positive, it follows that there are elements $a,b\in \hat{\AA}_{(\pi(j),\pi(k))}$ such that
\begin{align*}
\hat{\gamma}=\lambda_a\circ \lambda_b\circ \phi_k\circ \tilde{\gamma}_{(i,j,k)}\circ \phi_i^{-1}\ .
\end{align*}
Moreover, we have
$$ab=\hat{\gamma}(1)=1$$
and therefore $a=b^{-1}$. It finally follows that
$$\hat{\gamma}=\phi_k\circ \tilde{\gamma}_{(i,j,k)}\circ \phi_i^{-1}\ .$$
\qed
\end{bew}

\begin{bem}\label{76}
The standard involution $\sigma_s$ satisfies
$$\sigma_s\circ \phi=\phi\circ \sigma_s$$
for all $\phi\in \Aut(\HH)$, cf corollary \ref{109}.
\end{bem}

\begin{lemma}\label{55}
Let 
$$\tilde{\FFF}:=(\TTT(\tilde{\AA}_{(1,2)}):=\TTT(\tilde{\AA}_{(2,3)}):=\TTT(\tilde{\AA}_{(3,1)}):=\TTT(\HH),\tilde{\gamma}_2:=\tilde{\gamma}_{(1,2,3)},\tilde{\gamma}_3:=\tilde{\gamma}_{(2,3,1)},\tilde{\gamma}_1:=\tilde{\gamma}_{(3,1,2)})$$
be a positive triangle foundation over $\HH$. Then we have $$\tilde{\gamma}_3\tilde{\gamma}_2\tilde{\gamma}_1=\sigma_s\ .$$
\end{lemma}

\begin{bew}
The reparametrizations 
\begin{align*}\tilde{\alpha}_{(1,2)}:=(\HH,\tilde{\gamma}_2^{-1}\tilde{\gamma}_3^{-1},\tilde{\gamma}_2^{-1}\tilde{\gamma}_3^{-1},\tilde{\gamma}_2^{-1}\tilde{\gamma}_3^{-1})\ , && \tilde{\alpha}_{(2,3)}:=(\HH, \tilde{\gamma}_3^{-1}\sigma_s,\tilde{\gamma}_3^{-1}\sigma_s,\tilde{\gamma}_3^{-1}\sigma_s)
\end{align*}
show that $\tilde{\FFF}$ is isomorphic to the foundation
$$\hat{\FFF}:=\{\TTT(\hat{\AA}_{(1,2)}):=\TTT(\hat{\AA}_{(2,3)}):=\TTT(\hat{\AA}_{(3,1)}):=\TTT(\HH),\hat{\gamma}_2:=\hat{\gamma}_3:=\sigma_s,\hat{\gamma}_1:=\tilde{\gamma}_3\tilde{\gamma}_2\tilde{\gamma}_1\}\ .$$
\noindent Since we have $\hat{\FFF}\cong \PPP_3^+(\HH)$ by theorem \ref{24}, there are isomorphisms
\begin{align*}
\alpha_{(i,j)}=( \lambda_{a_i}\phi_{(i,j)}, \lambda_{a_i}\rho_{b_j}\phi_{(i,j)}, \rho_{b_j}\phi_{(i,j)} ): \TTT(\tilde{\AA}_{(i,j)}) \to \TTT(\hat{\AA}_{(i,j)})
\end{align*}
satisfying
\begin{align*}
\hat{\gamma_1}=\phi_{(1,2)}\sigma_s\phi_{(3,1)}^{-1}\ ,&& \sigma_s=\phi_{(2,3)}\sigma_s\phi_{(1,2)}^{-1}\ , && \sigma_s=\phi_{(3,1)}\sigma_s\phi_{(2,3)}^{-1}
\end{align*}
by the definition of an isomorphism and lemma \ref{75}. Now remark \ref{76} implies
\begin{align*}
\phi_{(1,2)}=\phi_{(2,3)}=\phi_{(3,1)}\ , && \tilde{\gamma}_3\tilde{\gamma}_2\tilde{\gamma}_1=\hat{\gamma}_1=\sigma_s\ .
\end{align*}\qed
\end{bew}

\begin{satz}\label{39} The foundation $\FFF$ is specially isomorphic to $\tilde{\FFF}:=\PPP_n^+(\HH)$.
\end{satz}

\begin{bew} Induction on $n$:
\begin{itemize}
\item $n=3$: This is theorem \ref{24}.
\item $n\to n+1$: By induction assumption, we may assume $ \tilde{\FFF}_{[1,n]}=\FFF_{[1,n]}$, where $\tilde{\FFF}_{[1,n]}$ and $\FFF_{[1,n]}$ are the $[1,n]$-residues of $\tilde{\FFF}$ and $\FFF$, respectively. If we reparametrize the rank 2 residues of $\FFF\setminus \FFF_{[1,n]}$ and $\tilde{\FFF}\sm \tilde{\FFF}_{[1,n]}$ in such a way that
\begin{align*}
\tilde{\gamma}_{(j,j+1,n+1)}=\gamma_{(j,j+1,n+1)}=\id^o=\gamma_{(n+1,1,2)}=\tilde{\gamma}_{(n+1,1,2)}\ ,\qquad j=1,\ldots,n-1\ ,
\end{align*}
the remaining glueings are uniquely determined by lemma \ref{55}, thus corresponding glueings are equal.
\end{itemize}

\begin{center}
\begin{tikzpicture}[scale=0.7,>=stealth,thick]
\begin{scope}[xshift=-9cm]
\coordinate (1) at (0,0);                    
\coordinate (2) at (-30:3);
\coordinate (3) at (90:3);
\coordinate (4) at (210:3);
\draw[blue] (1)--(2)--(4)--cycle;
\draw (2)--(3)--(1);
\draw (4)--(3);
\foreach \i in {1,2,3,4} {\fill (\i) circle (3pt);}
\draw[-,red] (2)+(120:1.5)arc[radius=1.5,start angle=120, end angle=150];
\draw[-] (0,1.5)arc[radius=1.5,start angle=270, end angle=300];
\draw[-] (0,1.5)arc[radius=1.5,start angle=270, end angle=240];
\draw[-,red] (4)+(60:1.5)arc[radius=1.5,start angle=60, end angle=30];
\draw[-,red] (0,1)arc[radius=1,start angle=90, end angle=-30];
\draw[-] (0,1)arc[radius=1,start angle=90, end angle=210];
\draw[-,blue] (210:1)arc[radius=1,start angle=210, end angle=330];
\draw[-,blue] (2)+(150:1.5)arc[radius=1.5,start angle=150, end angle=180];
\draw[-,blue] (4)+(30:1.5)arc[radius=1.5,start angle=30, end angle=0];
\node () at (0,-2) {\ssi$\TTT(\HH)$};
\node () at (-30:3.5) {\ssi{$1$}};
\node () at (-150:3.5) {\ssi{$3$}};
\node () at (90:3.5) {\ssi{$4$}};
\node[blue] () at (217.5:2) {$+$};
\node[red] () at (202.5:2) {$+$};
\node[blue] () at (-37.5:2) {$+$};
\node[red] () at (-22.5:2) {$+$};
\node () at (82.5:2) {$+$};
\node () at (97.5:2) {$+$};
\node[red] () at (30:0.5) {$+$};
\node () at (150:0.5) {$+$};
\node[blue] () at (270:0.5) {$+$};
\node () at (-3,3.35) {(i)};
\end{scope}

\begin{scope}
\coordinate (5) at (0,0);                    
\coordinate (6) at (-30:3);
\coordinate (7) at (90:3);
\coordinate (8) at (210:3);
\draw[-,red] (6)+(120:1.5)arc[radius=1.5,start angle=120, end angle=150];
\draw[-,red] (8)+(60:1.5)arc[radius=1.5,start angle=60, end angle=30];
\draws{0.53}[blue] (6)--(5);
\draws{0.53}[blue] (5)--(8);
\draws{0.53}[blue]  (8)--(6);
\draw (6)--(7)--(5);
\draw (8)--(7);
\foreach \i in {5,6,7,8} {\fill (\i) circle (3pt);}
\draw[-] (0,1.5)arc[radius=1.5,start angle=270, end angle=300];
\draw[-] (0,1.5)arc[radius=1.5,start angle=270, end angle=240];
\draw[-,red] (0,1)arc[radius=1,start angle=90, end angle=-30];
\draw[-] (0,1)arc[radius=1,start angle=90, end angle=210];
\draw[<-,blue] (210:1)arc[radius=1,start angle=210, end angle=330];
\draw[<-,blue] (6)+(150:1.5)arc[radius=1.5,start angle=150, end angle=180];
\draw[->,blue] (8)+(30:1.5)arc[radius=1.5,start angle=30, end angle=0];
\node () at (0,-2) {\ssi$\TTT(\HH)$};
\node () at (-30:3.5) {\ssi{$1$}};
\node () at (-150:3.5) {\ssi{$3$}};
\node () at (90:3.5) {\ssi{$4$}};
\node[blue] () at (219:2) {\ssi{$\sigma_s$}};
\node[red] () at (202.5:2) {$+$};
\node[blue] () at (-39:2) {\ssi{$\sigma_s$}};
\node[red] () at (-22.5:2) {$+$};
\node () at (82.5:2) {$+$};
\node () at (97.5:2) {$+$};
\node[red] () at (30:0.5) {$+$};
\node () at (150:0.5) {$+$};
\node[blue] () at (270:0.6) {\ssi{$\sigma_s$}};
\node () at (-3,3.35) {(ii)};
\end{scope}

\begin{scope}[xshift=-9cm,yshift=-6.5cm]
\coordinate (9) at (0,0);                    
\coordinate (10) at (-30:3);
\coordinate (11) at (90:3);
\coordinate (12) at (210:3);
\draw[->,red] (10)+(120:1.5)arc[radius=1.5,start angle=120, end angle=150];
\draw[<-,red] (12)+(60:1.5)arc[radius=1.5,start angle=60, end angle=30];
\draws{0.53}[blue] (10)--(9);
\draws{0.53}[blue] (9)--(12);
\draws{0.53}[blue]  (12)--(10);
\draw (10)--(11)--(9);
\draw (12)--(11);
\foreach \i in {9,10,11,12} {\fill (\i) circle (3pt);}
\draw[-] (0,1.5)arc[radius=1.5,start angle=270, end angle=300];
\draw[-] (0,1.5)arc[radius=1.5,start angle=270, end angle=240];
\draw[<-,red] (0,1)arc[radius=1,start angle=90, end angle=-30];
\draw[-] (0,1)arc[radius=1,start angle=90, end angle=210];
\draw[<-,blue] (210:1)arc[radius=1,start angle=210, end angle=330];
\draw[<-,blue] (10)+(150:1.5)arc[radius=1.5,start angle=150, end angle=180];
\draw[->,blue] (12)+(30:1.5)arc[radius=1.5,start angle=30, end angle=0];
\node () at (0,-2) {\ssi{$\TTT(\HH)$}};
\node () at (-30:3.5) {\ssi{$1$}};
\node () at (-150:3.5) {\ssi{$3$}};
\node () at (90:3.5) {\ssi{$4$}};
\node[blue] () at (219:2) {\ssi{$\sigma_s$}};
\node[red] () at (204:2) {\ssi{$\id^o$}};
\node[blue] () at (-39:2) {\ssi{$\sigma_s$}};
\node[red] () at (-22.5:2) {\ssi{$\id^o$}};
\node () at (82.5:2) {$+$};
\node () at (97.5:2) {$+$};
\node[red] () at (30:0.5) {\ssi{$\id^o$}};
\node () at (150:0.5) {$+$};
\node[blue] () at (270:0.6) {\ssi{$\sigma_s$}};
\node () at (-3,3.35) {(iii)};
\end{scope}

\begin{scope}[yshift=-6.5cm]
\coordinate (9) at (0,0);                    
\coordinate (10) at (-30:3);
\coordinate (11) at (90:3);
\coordinate (12) at (210:3);
\draw[->,red] (10)+(120:1.5)arc[radius=1.5,start angle=120, end angle=150];
\draw[<-,red] (12)+(60:1.5)arc[radius=1.5,start angle=60, end angle=30];
\draws{0.53}[blue] (10)--(9);
\draws{0.53}[blue] (9)--(12);
\draws{0.53}[blue]  (12)--(10);
\draw (10)--(11)--(9);
\draw (12)--(11);
\foreach \i in {9,10,11,12} {\fill (\i) circle (3pt);}
\draw[->] (0,1.5)arc[radius=1.5,start angle=270, end angle=300];
\draw[<-] (0,1.5)arc[radius=1.5,start angle=270, end angle=240];
\draw[<-,red] (0,1)arc[radius=1,start angle=90, end angle=-30];
\draw[->] (0,1)arc[radius=1,start angle=90, end angle=210];
\draw[<-,blue] (210:1)arc[radius=1,start angle=210, end angle=330];
\draw[<-,blue] (10)+(150:1.5)arc[radius=1.5,start angle=150, end angle=180];
\draw[->,blue] (12)+(30:1.5)arc[radius=1.5,start angle=30, end angle=0];
\node () at (0,-2) {\ssi{$\TTT(\HH)$}};
\node () at (-30:3.5) {\ssi{$1$}};
\node () at (-150:3.5) {\ssi{$3$}};
\node () at (90:3.5) {\ssi{$4$}};
\node[blue] () at (219:2) {\ssi{$\sigma_s$}};
\node[red] () at (204:2) {\ssi{$\id^o$}};
\node[blue] () at (-39:2) {\ssi{$\sigma_s$}};
\node[red] () at (-22.5:2) {\ssi{$\id^o$}};
\node () at (81.5:1.9) {\ssi{$\sigma_s$}};
\node () at (98.5:1.9) {\ssi{$\id^o$}};
\node[red] () at (30:0.5) {\ssi{$\id^o$}};
\node () at (150:0.5) {\ssi{$\sigma_s$}};
\node[blue] () at (270:0.6) {\ssi{$\sigma_s$}};
\node () at (-3,3.35) {(iv)};
\end{scope}
\end{tikzpicture}\end{center}
\qed
\end{bew}

\newpage

\chapter{Mixed Foundations over Skew-Fields}\label{358}

In order to determine the integrable mixed foundations, we attach a graph to each of them. If this graph is a tree, it turns out that the corresponding foundation is isomorphic to one of the foundations constructed in §\ref{154}. If it is not a tree, then the corresponding foundation is covered by a foundation of §\ref{154}.

\begin{no} Throughout this chapter, $\FFF$ is a mixed foundation over a skew-field, i.e., there are positive and negative glueings.\end{no}

\begin{bem} Since $\FFF$ has positive residues, $\HH:=\AA$ is a quaternion division algebra by lemma \ref{37}. 
\end{bem}

\begin{de} Let $\FFF$ be a foundation over $I=V(F)$. Given $i\in I$, the \textit{set of neighbours of $\mathit{i}$}\index{set of neighbours} is
$$B_1(i):=\{ j\in I \mid \{i,j\}\in E(F) \}\ .$$
\end{de}

\begin{lemma}\label{47}
Let $\tilde{\FFF}$ be a foundation over $I=\{1,2,3,4\}$ such that its defining field is a skew-field and such that $|B_1(1)|=3$. Let
$$\tilde{\Gamma}:=\{\tilde{\gamma}_{(2,1,3)},\tilde{\gamma}_{(3,1,4)},\tilde{\gamma}_{(2,1,4)}\}\ .$$
Then we have
$$n:=|\{\tilde{\gamma}\in \tilde{\Gamma} \mid \tilde{\gamma}\ \textrm{positive}\}|\in\{1,3\}\ .$$
\end{lemma}

\begin{bew} Notice that a glueing is either negative or positive.
\begin{enumerate}[label=(\roman*)]
\item Suppose that $\tilde{\gamma}_{(2,1,3)}$ and $\tilde{\gamma}_{(3,1,4)}$ are negative. Then
$$\tilde{\gamma}_{(2,1,4)}=\tilde{\gamma}_{(3,1,4)}\circ\id^o\circ\tilde{\gamma}_{(2,1,3)}$$
is positive, thus $n\geq 1$.
\item Suppose that $\tilde{\gamma}_{(2,1,3)}$ is positive and that $\tilde{\gamma}_{(3,1,3)}$ is negative. Then
$$\tilde{\gamma}_{(2,1,4)}=\tilde{\gamma}_{(3,1,4)}\circ\id^o\circ\tilde{\gamma}_{(2,1,3)}$$
is negative, thus $n\neq 2$.
\end{enumerate}

\begin{center}\begin{tikzpicture}[scale=0.35,>=stealth,thick]
\begin{scope}
\coordinate (1) at (0,0);                    
\coordinate (2) at (-30:3);
\coordinate (3) at (90:3);
\coordinate (4) at (210:3);
\draw (3)--(1)--(2);
\draw (4)--(1);
\foreach \i in {1,2,3,4} {\fill (\i) circle (3pt);}
\draw[blue] (0,1)arc[radius=1,start angle=90, end angle=-30];
\draw[blue] (0,1)arc[radius=1,start angle=90, end angle=210];
\draw (210:1)arc[radius=1,start angle=210, end angle=330];
\node () at (-30:3.5) {\ssi{$2$}};
\node () at (-150:3.5) {\ssi{$3$}};
\node () at (90:3.5) {\ssi{$4$}};
\node[blue] () at (30:0.5) {\ssi{$-$}};
\node[blue] () at (150:0.5) {\ssi{$-$}};
\node[blue] () at (270:0.55) {};
\node () at (-3,3.35) {(i)};
\draw[->,double] (3,1)--(6,1);
\end{scope}

\begin{scope}[xshift=9cm]
\coordinate (1) at (0,0);                    
\coordinate (2) at (-30:3);
\coordinate (3) at (90:3);
\coordinate (4) at (210:3);
\draw (3)--(1)--(2);
\draw (4)--(1);
\foreach \i in {1,2,3,4} {\fill (\i) circle (3pt);}
\draw[blue] (0,1)arc[radius=1,start angle=90, end angle=-30];
\draw[blue] (0,1)arc[radius=1,start angle=90, end angle=210];
\draw[red] (210:1)arc[radius=1,start angle=210, end angle=330];
\node () at (-30:3.5) {\ssi{$2$}};
\node () at (-150:3.5) {\ssi{$3$}};
\node () at (90:3.5) {\ssi{$4$}};
\node[blue] () at (30:0.5) {\ssi{$-$}};
\node[blue] () at (150:0.5) {\ssi{$-$}};
\node[red] () at (270:0.55) {\ssi{$+$}};
\end{scope}

\begin{scope}[xshift=18cm]
\coordinate (1) at (0,0);                    
\coordinate (2) at (-30:3);
\coordinate (3) at (90:3);
\coordinate (4) at (210:3);
\draw (3)--(1)--(2);
\draw (4)--(1);
\foreach \i in {1,2,3,4} {\fill (\i) circle (3pt);}
\draw[blue] (0,1)arc[radius=1,start angle=90, end angle=-30];
\draw[red] (0,1)arc[radius=1,start angle=90, end angle=210];
\draw (210:1)arc[radius=1,start angle=210, end angle=330];
\node () at (-30:3.5) {\ssi{$2$}};
\node () at (-150:3.5) {\ssi{$3$}};
\node () at (90:3.5) {\ssi{$4$}};
\node[blue] () at (30:0.5) {\ssi{$-$}};
\node[red] () at (150:0.5) {\ssi{$+$}};
\node[blue] () at (270:0.55) {};
\node () at (-3,3.35) {(ii)};
\draw[->,double] (3,1)--(6,1);
\end{scope}

\begin{scope}[xshift=27cm]
\coordinate (1) at (0,0);                    
\coordinate (2) at (-30:3);
\coordinate (3) at (90:3);
\coordinate (4) at (210:3);
\draw (3)--(1)--(2);
\draw (4)--(1);
\foreach \i in {1,2,3,4} {\fill (\i) circle (3pt);}
\draw[blue] (0,1)arc[radius=1,start angle=90, end angle=-30];
\draw[red] (0,1)arc[radius=1,start angle=90, end angle=210];
\draw[blue] (210:1)arc[radius=1,start angle=210, end angle=330];
\node () at (-30:3.5) {\ssi{$2$}};
\node () at (-150:3.5) {\ssi{$3$}};
\node () at (90:3.5) {\ssi{$4$}};
\node[blue] () at (30:0.5) {\ssi{$-$}};
\node[red] () at (150:0.5) {\ssi{$+$}};
\node[blue] () at (270:0.55) {\ssi{$-$}};
\end{scope}
\end{tikzpicture}\end{center}
Notice that there could be more edges than the drawn ones.\qed
\end{bew}

\begin{no} Given a mixed foundation $\tilde{\FFF}$ over $\HH$, we denote the collection of all maximal positive residues of $\tilde{\FFF}$ by $\MMM(\tilde{\FFF})$. We set
$$P(\tilde{\FFF}):=\{ \GGG_{\tilde{P}} \mid \tilde{\PPP}\in \MMM(\tilde{\FFF})\}\ .$$
\end{no}

\begin{bem}\label{69}
Rank 2 two residues are considered to be positive, thus 
$$\bigcup_{\mathclap{ {\PPP}\in \MMM( {\FFF})}}\ \GGG_{ {P}}= \GGG_{ {F}}\ .$$
\end{bem}

\newpage

\begin{satz}\label{79}
The set $\MMM(\FFF)$ satisfies the following conditions:
\begin{enumerate}[label=(\roman*)]
\item Each  $ {\PPP}\in \MMM( {\FFF})$ is integrable, in particular, $\GGG_{ {P}}$ is complete.$\vphantom{\GGG_{ {F}}^{P( {\FFF})}}$
\item We have $$\bigcup_{\mathclap{ {\PPP}\in \MMM( {\FFF})}}\ \GGG_{ {P}}= \GGG_{ {F}}\ .$$
\item Given $ {\PPP}_1\neq {\PPP}_2\in \MMM( {\FFF})$, we have $|V( {P}_1)\cap V( {P}_2)|\leq 1$. $\vphantom{\GGG_{ {F}}^{P( {\FFF})}}$
\item Given $i\in I$, we have $|\{  {\PPP}\in \MMM( {\FFF})\mid i\in V( {P})\}|\leq 2$. $\vphantom{\GGG_{ {F}}^{P( {\FFF})}}$
\end{enumerate}
\end{satz}

\begin{bewzwei}\ 
\begin{enumerate}[label=(\roman*)]
\item Each $\PPP\in \MMM(\FFF)$ is integrable because $\FFF$ itself is integrable.
\item This holds by remark \ref{69}.
\item Since the elements of $\MMM(\FFF)$ are maximal and $G_P$ is complete for each $\PPP\in \MMM(\FFF)$, lemma \ref{47} implies that
$$\forall\ \PPP_1\neq\PPP_2\in \MMM(\FFF):\qquad  |V(P_1)\cap V(P_2)|\leq 1\ .$$
Otherwise $\PPP_1\cup\PPP_2$ would be positive with $\PPP_1\subsetneq \PPP_1\cup \PPP_2$.
\item By step (ii), the glueings connecting two elements $\PPP_1\neq \PPP_2\in \MMM(\FFF)$ are necessarily negative. From lemma \ref{47} again it follows that
$$\forall\ i\in I:\qquad |\{ \PPP\in \MMM(\FFF)\mid i\in V(P)\}|\leq 2\ .$$
\end{enumerate}\qed
\end{bewzwei}

\begin{prop}\label{80} Let $\tilde{\FFF}$ be a mixed foundation over $\HH$ such that $\MMM(\tilde{\FFF})$ satisfies the conditions (i)-(iv) of theorem \ref{79}.
If the graph $\GGG:=\GGG_{\tilde{F}}^{P(\tilde{\FFF})}$ is a tree, we have $\tilde{\FFF}\cong \FFF(\GGG,\HH)=:\hat{\FFF}$. In particular, $\tilde{\FFF}$ is integrable.
\end{prop}

\begin{bew}
First of all we observe that $\hat{\FFF}$ is well-defined the conditions (i)-(iv). By construction, we have $P(\tilde{\FFF})=P(\hat{\FFF})$, and by theorem \ref{39}, we have
$$\tilde{\PPP}\cong \PPP_{|V(\tilde{P})|}^+(\HH)=:\hat{\PPP}\in \MMM(\hat{\FFF})$$ for each
$\tilde{P}\in \MMM(\tilde{\FFF})$. Since the isomorphisms are special, we may assume $\tilde{\PPP}=\hat{\PPP}$ for each $\tilde{P}\in \MMM(\tilde{\FFF})$. It remains to adjust the glueings connecting two elements of $\MMM(\tilde{\FFF})$. Since $\GGG$ is a tree and $$|\{ \tilde{\PPP}\in \MMM(\tilde{\FFF})\mid i\in V(\tilde{P})\}|\leq 2$$ for each vertex $i\in I$, it suffices to show the following:\medskip

\noindent Given $\tilde{\PPP}_1,\tilde{\PPP}_2\in \MMM(\tilde{F})$ with $\{\GGG_{\tilde{P}_1},\GGG_{\tilde{P}_2}\}\in E(\GGG)$, there is an isomorphism $\alpha:\tilde{\PPP}_1\cup\tilde{\PPP}_2\to \hat{\PPP}_1\cup\hat{\PPP}_2$ fixing $\tilde{\PPP}_1$.\medskip

\noindent Let $V(\tilde{P}_1)\cap V(\tilde{P}_2)=\{b\}$, let $a\in V(\tilde{P}_1)\sm \{b\}$, let $c\in V(\tilde{P}_2)\sm \{b\}$ and let $\gamma:=\hat{\gamma}_{(a,b,c)}\circ \tilde{\gamma}_{(a,b,c)}^{-1}$. Then
$$\alpha:=\{ \id_F, \alpha_{(i,j)} \mid (i,j)\in A(\tilde{P}_1\cup \tilde{P}_2)\}$$
with
\begin{align*}
\forall\ (i,j)\in A(\tilde{P}_1):\alpha_{(i,j)}:=(\id,\id,\id )\ , && \forall\ (i,j)\in A(\tilde{P}_2):\alpha_{(i,j)}=(\gamma,\gamma,\gamma)
\end{align*}
satisfies the required condition.\medskip

\noindent Finally $\tilde{\FFF}$ is integrable by corollary \ref{82} and theorem \ref{62}.\qed
\end{bew}

\newpage

\begin{bsp}\ 
\begin{center}\begin{tikzpicture}[scale=0.35,>=stealth,thick]

\begin{scope}
\begin{scope}[xshift=-6cm]  
\node () at (-2,3) {(i)};        
\coordinate (1) at (0,0);
\coordinate (2) at (3,0);
\coordinate (3) at (60:3);
\draw[red] (3)--(1)--(2)--cycle;
\foreach \i in {1,2,3} {\fill[red] (\i) circle (3pt);}
\node[red] () at (1.5,1) {$+$};
\draw[red,dashed] (1)--(-2,0);
\node[blue] () at (3,-0.4) {$-$};
\node[blue] () at (0,-0.4) {$-$};
\draw[blue] (2,0)arc[radius=1,start angle=180, end angle=360];
\draw[blue] (1,0)arc[radius=1,start angle=0, end angle=-180];
\end{scope}
\begin{scope}[xshift=-3cm]          
\coordinate (1) at (0,0);
\coordinate (2) at (3,0);
\coordinate (3) at (60:3);
\draw[red] (3)--(1)--(2)--cycle;
\foreach \i in {1,2,3} {\fill[red] (\i) circle (3pt);}
\node[red] () at (1.5,1) {$+$};
\node[blue] () at (3,-0.4) {$-$};
\draw[blue] (2,0)arc[radius=1,start angle=180, end angle=360];
\end{scope}
\begin{scope}          
\coordinate (1) at (0,0);
\coordinate (2) at (3,0);
\coordinate (3) at (60:3);
\draw[red] (3)--(1)--(2)--cycle;
\foreach \i in {1,2,3} {\fill[red] (\i) circle (3pt);}
\node[red] () at (1.5,1) {$+$};
\node[blue] () at (3,-0.4) {$-$};
\draw[blue] (2,0)arc[radius=1,start angle=180, end angle=360];
\end{scope}
\begin{scope}[xshift=3cm]          
\coordinate (1) at (0,0);
\coordinate (2) at (3,0);
\coordinate (3) at (60:3);
\draw[red] (3)--(1)--(2)--cycle;
\foreach \i in {1,2,3} {\fill[red] (\i) circle (3pt);}
\node[red] () at (1.5,1) {$+$};
\node[blue] () at (3,-0.4) {$-$};
\draw[blue] (2,0)arc[radius=1,start angle=180, end angle=360];
\end{scope}
\begin{scope}[xshift=6cm]          
\coordinate (1) at (0,0);
\coordinate (2) at (3,0);
\coordinate (3) at (60:3);
\draw[red] (3)--(1)--(2)--cycle;
\foreach \i in {1,2,3} {\fill[red] (\i) circle (3pt);}
\node[red] () at (1.5,1) {$+$};
\draw[red,dashed] (2)--(5,0);
\node[blue] () at (3,-0.4) {$-$};
\draw[blue] (2,0)arc[radius=1,start angle=180, end angle=360];
\end{scope}
\end{scope}

\begin{scope}[xshift=21cm]
\begin{scope}[xshift=-6cm]    
\node () at (-2,3) {(ii)};       
\coordinate (1) at (0,0);
\coordinate (2) at (3,0);
\coordinate (3) at (60:3);
\draw[red] (3)--(1)--(2)--cycle;
\foreach \i in {1,2,3} {\fill[red] (\i) circle (3pt);}
\node[red] () at (1.5,1) {$+$};
\draw[red,dashed] (1)--(-2,0);
\node[blue] () at (3,-0.4) {$-$};
\node[blue] () at (0,-0.4) {$-$};
\draw[blue] (2,0)arc[radius=1,start angle=180, end angle=360];
\draw[blue] (1,0)arc[radius=1,start angle=0, end angle=-180];
\end{scope}
\begin{scope}[xshift=-3cm]          
\coordinate (1) at (0,0);
\coordinate (2) at (3,0);
\coordinate (3) at (60:3);
\draw[red] (3)--(1)--(2)--cycle;
\foreach \i in {1,2,3} {\fill[red] (\i) circle (3pt);}
\node[red] () at (1.5,1) {$+$};
\node[black] () at (3,-0.4) {$-$};
\draw[black] (2,0)arc[radius=1,start angle=180, end angle=360];
\end{scope}
\begin{scope}          
\coordinate (1) at (0,0);
\coordinate (2) at (3,0);
\coordinate (3) at (60:3);
\draw[red] (3)--(1)--(2)--cycle;
\foreach \i in {1,2,3} {\fill[red] (\i) circle (3pt);}
\node[red] () at (1.5,1) {$+$};
\node[black] () at (3,-0.4) {$-$};
\draw[black] (2,0)arc[radius=1,start angle=180, end angle=360];
\end{scope}
\begin{scope}[xshift=3cm]          
\coordinate (1) at (0,0);
\coordinate (2) at (3,0);
\coordinate (3) at (60:3);
\draw[red] (3)--(1)--(2)--cycle;
\foreach \i in {1,2,3} {\fill[red] (\i) circle (3pt);}
\node[red] () at (1.5,1) {$+$};
\node[blue] () at (3,-0.4) {$-$};
\draw[blue] (2,0)arc[radius=1,start angle=180, end angle=360];
\end{scope}
\begin{scope}[xshift=6cm]          
\coordinate (1) at (0,0);
\coordinate (2) at (3,0);
\coordinate (3) at (60:3);
\draw[red] (3)--(1)--(2)--cycle;
\foreach \i in {1,2,3} {\fill[red] (\i) circle (3pt);}
\node[red] () at (1.5,1) {$+$};
\draw[red,dashed] (2)--(5,0);
\node[blue] () at (3,-0.4) {$-$};
\draw[blue] (2,0)arc[radius=1,start angle=180, end angle=360];
\end{scope}
\end{scope}

\begin{scope}[yshift=-6cm]
\begin{scope}[xshift=-6cm]   
\node () at (-2,3) {(iii)};        
\coordinate (1) at (0,0);
\coordinate (2) at (3,0);
\coordinate (3) at (60:3);
\draw[red] (3)--(1)--(2)--cycle;
\foreach \i in {1,2,3} {\fill[red] (\i) circle (3pt);}
\node[red] () at (1.5,1) {$+$};
\draw[red,dashed] (1)--(-2,0);
\node[black] () at (3,-0.4) {$-$};
\node[blue] () at (0,-0.4) {$-$};
\draw[black] (2,0)arc[radius=1,start angle=180, end angle=360];
\draw[blue] (1,0)arc[radius=1,start angle=0, end angle=-180];
\end{scope}
\begin{scope}[xshift=-3cm]          
\coordinate (1) at (0,0);
\coordinate (2) at (3,0);
\coordinate (3) at (60:3);
\draw[red] (3)--(1)--(2)--cycle;
\foreach \i in {1,2,3} {\fill[red] (\i) circle (3pt);}
\node[red] () at (1.5,1) {$+$};
\node[black] () at (3,-0.4) {$-$};
\draw[black] (2,0)arc[radius=1,start angle=180, end angle=360];
\end{scope}
\begin{scope}          
\coordinate (1) at (0,0);
\coordinate (2) at (3,0);
\coordinate (3) at (60:3);
\draw[red] (3)--(1)--(2)--cycle;
\foreach \i in {1,2,3} {\fill[red] (\i) circle (3pt);}
\node[red] () at (1.5,1) {$+$};
\node[black] () at (3,-0.4) {$-$};
\draw[black] (2,0)arc[radius=1,start angle=180, end angle=360];
\end{scope}
\begin{scope}[xshift=3cm]          
\coordinate (1) at (0,0);
\coordinate (2) at (3,0);
\coordinate (3) at (60:3);
\draw[red] (3)--(1)--(2)--cycle;
\foreach \i in {1,2,3} {\fill[red] (\i) circle (3pt);}
\node[red] () at (1.5,1) {$+$};
\node[black] () at (3,-0.4) {$-$};
\draw[black] (2,0)arc[radius=1,start angle=180, end angle=360];
\end{scope}
\begin{scope}[xshift=6cm]          
\coordinate (1) at (0,0);
\coordinate (2) at (3,0);
\coordinate (3) at (60:3);
\draw[red] (3)--(1)--(2)--cycle;
\foreach \i in {1,2,3} {\fill[red] (\i) circle (3pt);}
\node[red] () at (1.5,1) {$+$};
\draw[red,dashed] (2)--(5,0);
\node[blue] () at (3,-0.4) {$-$};
\draw[blue] (2,0)arc[radius=1,start angle=180, end angle=360];
\end{scope}
\end{scope}

\begin{scope}[xshift=21cm,yshift=-6cm]
\begin{scope}[xshift=-6cm]  
\node () at (-2,3) {(iv)};         
\coordinate (1) at (0,0);
\coordinate (2) at (3,0);
\coordinate (3) at (60:3);
\draw[red] (3)--(1)--(2)--cycle;
\foreach \i in {1,2,3} {\fill[red] (\i) circle (3pt);}
\node[red] () at (1.5,1) {$+$};
\draw[red,dashed] (1)--(-2,0);
\node[black] () at (3,-0.4) {$-$};
\node[black] () at (0,-0.4) {$-$};
\draw[black] (2,0)arc[radius=1,start angle=180, end angle=360];
\draw[black] (1,0)arc[radius=1,start angle=0, end angle=-180];
\end{scope}
\begin{scope}[xshift=-3cm]          
\coordinate (1) at (0,0);
\coordinate (2) at (3,0);
\coordinate (3) at (60:3);
\draw[red] (3)--(1)--(2)--cycle;
\foreach \i in {1,2,3} {\fill[red] (\i) circle (3pt);}
\node[red] () at (1.5,1) {$+$};
\node[black] () at (3,-0.4) {$-$};
\draw[black] (2,0)arc[radius=1,start angle=180, end angle=360];
\end{scope}
\begin{scope}          
\coordinate (1) at (0,0);
\coordinate (2) at (3,0);
\coordinate (3) at (60:3);
\draw[red] (3)--(1)--(2)--cycle;
\foreach \i in {1,2,3} {\fill[red] (\i) circle (3pt);}
\node[red] () at (1.5,1) {$+$};
\node[black] () at (3,-0.4) {$-$};
\draw[black] (2,0)arc[radius=1,start angle=180, end angle=360];
\end{scope}
\begin{scope}[xshift=3cm]          
\coordinate (1) at (0,0);
\coordinate (2) at (3,0);
\coordinate (3) at (60:3);
\draw[red] (3)--(1)--(2)--cycle;
\foreach \i in {1,2,3} {\fill[red] (\i) circle (3pt);}
\node[red] () at (1.5,1) {$+$};
\node[black] () at (3,-0.4) {$-$};
\draw[black] (2,0)arc[radius=1,start angle=180, end angle=360];
\end{scope}
\begin{scope}[xshift=6cm]          
\coordinate (1) at (0,0);
\coordinate (2) at (3,0);
\coordinate (3) at (60:3);
\draw[red] (3)--(1)--(2)--cycle;
\foreach \i in {1,2,3} {\fill[red] (\i) circle (3pt);}
\node[red] () at (1.5,1) {$+$};
\draw[red,dashed] (2)--(5,0);
\node[black] () at (3,-0.4) {$-$};
\draw[black] (2,0)arc[radius=1,start angle=180, end angle=360];
\end{scope}
\end{scope}
\end{tikzpicture}\end{center}
\noindent The red triangles represent the elements of $\MMM(\tilde{\FFF})$, the blue glueings represent the corresponding glueings of $\tilde{\FFF}$, while the black ones represent those of $\hat{\FFF}$.
\end{bsp}

\begin{kor}
If $\GGG:=\GGG_F^{P(\FFF)}$ is a tree, then we have
$$\FFF\cong \FFF(\GGG,\HH)\ .$$
\end{kor}

\begin{bew}
This results immediately from theorem \ref{79} and proposition \ref{80}.\qed
\end{bew}

\begin{satz}\label{443}
Let $\tilde{\FFF}$ be a mixed foundation over $\HH$ such that $\MMM(\tilde{\FFF})$ satisfies the conditions (i)-(iv) of theorem \ref{79}. Then $\tilde{\FFF}$ is integrable.
\end{satz}

\begin{bew}
Let $\UUU$ be the the universal cover of $\GGG_{\tilde{F}}^{P(\tilde{\FFF})}$ and let $\hat{\FFF}$ be the cover corresponding to the induced cover of $\GGG_{\tilde{F}}$. By theorem \ref{145}, $\tilde{\FFF}$ is integrable if $\hat{\FFF}$ is integrable. But $\MMM(\hat{\FFF})$ equally satisfies the conditions (i)-(iii) of proposition \ref{80} and, moreover,
$$\GGG_{\hat{F}}^{P(\hat{\FFF})}\cong \UUU$$
is a tree. Therefore, $\hat{\FFF}$ is integrable by proposition \ref{80}.\qed
\end{bew}

\begin{bsp}\ 
\begin{center}\begin{tikzpicture}[scale=0.35,>=stealth,thick]

\coordinate (1) at (0,0);
\coordinate (2) at (3,0);
\coordinate (3) at (6,0);
\coordinate (4) at (9,0);
\coordinate (5) at (60:3);
\coordinate (6) at (60:6);
\coordinate (7) at (60:9);
\coordinate (8) at ($(4)+(120:3)$);
\coordinate (9) at ($(4)+(120:6)$);
\draw[red] (1)--(2)--(5)--cycle;
\draw[red] (3)--(4)--++(120:3)--cycle;
\draw[red] (6)--(7)--++(-60:3)--cycle;
\draw (8)--(9);
\draw (2)--(3);
\draw (5)--(6);
\draw[blue] (2,0)arc[radius=1,start angle=180, end angle=360];
\draw[blue] (5,0)arc[radius=1,start angle=180, end angle=360];
\draw[blue] (4)+(120:5)arc[radius=1,start angle=-60, end angle=120];
\draw[blue] (4)+(120:2)arc[radius=1,start angle=-60, end angle=120];
\draw[blue] (60:7)arc[radius=1,start angle=60, end angle=240];
\draw[blue] (60:4)arc[radius=1,start angle=60, end angle=240];
\node[blue] () at (3,-0.5) {$-$};
\node[blue] () at (6,-0.5) {$-$};
\node[blue] () at ($(8)+(30:0.5)$) {$-$};
\node[blue] () at ($(9)+(30:0.5)$) {$-$};
\node[blue] () at ($(5)+(150:0.5)$) {$-$};
\node[blue] () at ($(6)+(150:0.5)$) {$-$};
\node[red] () at (1.5,1) {$+$};
\node[red] () at (7.5,1) {$+$};
\node[red] () at ($(7.5,1)+(120:6)$) {$+$};
\node () at (1.5,-1) {\ssi$\PPP_3$};
\node () at (7.5,-1) {\ssi$\PPP_1$};
\node () at ($(7.5,-1)+(120:6)$) {\ssi$\PPP_2$};
\node () at (0,8) {$\tilde{\FFF}$};
\foreach \i in {1,...,9} {\fill (\i) circle (3pt);}

\begin{scope}[xshift=3cm]
\coordinate (10) at (3,-6);
\coordinate (11) at (6,-6);
\coordinate (12) at ($(10)+(60:3)$);
\draw[red] (10)--(11)--(12)--cycle;
\draw (11)--(9,-6);
\node[red] () at (4.5,-5) {$+$};
\node[blue] () at (6,-6.5) {$-$};
\node[blue] () at (3,-6.5) {$-$};
\node () at (4.5,-7) {\ssi$\PPP_1$};
\draw[blue] (5,-6)arc[radius=1,start angle=180, end angle=360];
\draw[blue] (2,-6)arc[radius=1,start angle=180, end angle=360];
\foreach \i in {10,...,12} {\fill (\i) circle (3pt);}
\end{scope}
\begin{scope}[xshift=9cm]
\coordinate (10) at (3,-6);
\coordinate (11) at (6,-6);
\coordinate (12) at ($(10)+(60:3)$);
\draw[red] (10)--(11)--(12)--cycle;
\draw (11)--(9,-6);
\node[red] () at (4.5,-5) {$+$};
\node[blue] () at (6,-6.5) {$-$};
\node[blue] () at (3,-6.5) {$-$};
\node () at (4.5,-7) {\ssi$\PPP_2$};
\draw[blue] (5,-6)arc[radius=1,start angle=180, end angle=360];
\draw[blue] (2,-6)arc[radius=1,start angle=180, end angle=360];
\foreach \i in {10,...,12} {\fill (\i) circle (3pt);}\end{scope}
\begin{scope}[xshift=15cm]
\coordinate (10) at (3,-6);
\coordinate (11) at (6,-6);
\coordinate (12) at ($(10)+(60:3)$);
\draw[red] (10)--(11)--(12)--cycle;
\draw[dashed] (11)--(8,-6);
\node[red] () at (4.5,-5) {$+$};
\node[blue] () at (6,-6.5) {$-$};
\node[blue] () at (3,-6.5) {$-$};
\node () at (4.5,-7) {\ssi$\PPP_3$};
\draw[blue] (5,-6)arc[radius=1,start angle=180, end angle=360];
\draw[blue] (2,-6)arc[radius=1,start angle=180, end angle=360];
\foreach \i in {10,...,12} {\fill (\i) circle (3pt);}\end{scope}
\begin{scope}[xshift=-3cm]
\coordinate (10) at (3,-6);
\coordinate (11) at (6,-6);
\coordinate (12) at ($(10)+(60:3)$);
\draw[red] (10)--(11)--(12)--cycle;
\draw (11)--(9,-6);
\node[red] () at (4.5,-5) {$+$};
\node[blue] () at (6,-6.5) {$-$};
\node[blue] () at (3,-6.5) {$-$};
\node () at (4.5,-7) {\ssi$\PPP_3$};
\draw[blue] (5,-6)arc[radius=1,start angle=180, end angle=360];
\draw[blue] (2,-6)arc[radius=1,start angle=180, end angle=360];
\foreach \i in {10,...,12} {\fill (\i) circle (3pt);}\end{scope}
\begin{scope}[xshift=-9cm]
\coordinate (10) at (3,-6);
\coordinate (11) at (6,-6);
\coordinate (12) at ($(10)+(60:3)$);
\draw[red] (10)--(11)--(12)--cycle;
\draw (11)--(9,-6);
\node[red] () at (4.5,-5) {$+$};
\node[blue] () at (6,-6.5) {$-$};
\node[blue] () at (3,-6.5) {$-$};
\node () at (4.5,-7) {\ssi{$\PPP_2$}};
\node () at (12,0) {};
\draw[blue] (5,-6)arc[radius=1,start angle=180, end angle=360];
\draw[blue] (2,-6)arc[radius=1,start angle=180, end angle=360];
\draw[dashed] (1,-6)--(10);
\node () at (1,-3) {$\hat{\FFF}$};
\foreach \i in {10,...,12} {\fill (\i) circle (3pt);}
\draw[->,double]  (-5,-5)--(-2,-5);
\end{scope}
\end{tikzpicture}\end{center}
\end{bsp}

\newpage

\chapter{Negative Foundations}
Negative foundations are quite easy to handle as we may apply theorems \ref{145} and \ref{146}. If $\AA$ is a field, there are no restrictions; each foundation is integrable. If $\AA$ is non-commutative, then lemma \ref{26} is very restrictive; the structure of an integrable foundation is very simple.

\begin{no} Throughout this chapter, $\FFF$ is a negative foundation over a skew-field.\end{no}

\begin{lemma}\label{26} If ${\FFF}$ has a residue ${\RRR}$ of type $D_4$, then $\AA$ is a field. In particular, ${\FFF}$ is negative.
\end{lemma}

\begin{bew} We label the vertices of ${\RRR}$ such that $\{1,2\},\{1,3\},\{1,4\}\in E({R})$. By lemma \ref{22}, each glueing is negative, and by lemma \ref{47}, at least one is positive.
\vspace*{0.3cm}
\begin{center}\begin{tikzpicture}[scale=0.55,>=stealth,thick]
\begin{scope}[xshift=-9cm]
\coordinate (1) at (0,0);                    
\coordinate (2) at (-30:3);
\coordinate (3) at (90:3);
\coordinate (4) at (210:3);
\draw (3)--(1)--(2);
\draw (4)--(1);
\foreach \i in {1,2,3,4} {\fill (\i) circle (3pt);}
\draw[blue] (0,1)arc[radius=1,start angle=90, end angle=-30];
\draw[blue] (0,1)arc[radius=1,start angle=90, end angle=210];
\draw[blue] (210:1)arc[radius=1,start angle=210, end angle=330];
\node () at (-30:3.5) {\ssi{$2$}};
\node () at (-150:3.5) {\ssi{$3$}};
\node () at (90:3.5) {\ssi{$4$}};
\node[blue] () at (30:0.5) {\ssi{$-$}};
\node[blue] () at (150:0.5) {\ssi{$-$}};
\node[blue] () at (0,-0.6) {\ssi{$-$}};
\node[blue] () at (270:0.55) {};
\node () at (-3,3.35) {\ssi{\ref{22}}};
\node () at (4.5,1) {\&};
\draw[->,double] (2,-4)--(5,-5.5);
\end{scope}

\begin{scope}
\coordinate (1) at (0,0);                    
\coordinate (2) at (-30:3);
\coordinate (3) at (90:3);
\coordinate (4) at (210:3);
\draw (3)--(1)--(2);
\draw (4)--(1);
\foreach \i in {1,2,3,4} {\fill (\i) circle (3pt);}
\draw[blue] (0,1)arc[radius=1,start angle=90, end angle=-30];
\draw[blue] (0,1)arc[radius=1,start angle=90, end angle=210];
\draw (210:1)arc[radius=1,start angle=210, end angle=330];
\node () at (-30:3.5) {\ssi{$2$}};
\node () at (-150:3.5) {\ssi{$3$}};
\node () at (90:3.5) {\ssi{$4$}};
\node[blue] () at (30:0.5) {\ssi{$-$}};
\node[blue] () at (150:0.5) {\ssi{$-$}};
\node[blue] () at (270:0.55) {};
\node () at (4.25,1.7) {{\ssi{\ref{47}}}};
\draw[->,double] (3,1)--(6,1);
\end{scope}

\begin{scope}[xshift=9cm]
\coordinate (1) at (0,0);                    
\coordinate (2) at (-30:3);
\coordinate (3) at (90:3);
\coordinate (4) at (210:3);
\draw (3)--(1)--(2);
\draw (4)--(1);
\foreach \i in {1,2,3,4} {\fill (\i) circle (3pt);}
\draw[blue] (0,1)arc[radius=1,start angle=90, end angle=-30];
\draw[blue] (0,1)arc[radius=1,start angle=90, end angle=210];
\draw[red] (210:1)arc[radius=1,start angle=210, end angle=330];
\node () at (-30:3.5) {\ssi{$2$}};
\node () at (-150:3.5) {\ssi{$3$}};
\node () at (90:3.5) {\ssi{$4$}};
\node[blue] () at (30:0.5) {\ssi{$-$}};
\node[blue] () at (150:0.5) {\ssi{$-$}};
\node[red] () at (0,-0.55) {\ssi{$+$}};
\node[blue] () at (270:0.55) {};
\draw[->,double] (-2,-4)--(-5,-5.5);
\end{scope}

\begin{scope}[yshift=-7cm]
\coordinate (1) at (0,0);                    
\coordinate (2) at (-30:3);
\coordinate (3) at (90:3);
\coordinate (4) at (210:3);
\draw (3)--(1)--(2);
\draw (4)--(1);
\foreach \i in {1,2,3,4} {\fill (\i) circle (3pt);}
\draw[blue] (0,1)arc[radius=1,start angle=90, end angle=-30];
\draw[blue] (0,1)arc[radius=1,start angle=90, end angle=210];
\draw[violet] (210:1)arc[radius=1,start angle=210, end angle=330];
\node () at (-30:3.5) {\ssi{$2$}};
\node () at (-150:3.5) {\ssi{$3$}};
\node () at (90:3.5) {\ssi{$4$}};
\node[blue] () at (30:0.5) {\ssi{$-$}};
\node[blue] () at (150:0.5) {\ssi{$-$}};
\node[violet] () at (0,-0.55) {\ssi{$\pm$}};
\node[blue] () at (270:0.6) {};
\end{scope}
\end{tikzpicture}\end{center}
\qed
\end{bew}

\begin{satz}\label{439} If $\AA$ is a non-commutative skew-field, then $\GGG_F$ is a string, a ray, a chain or a circle.
\end{satz}

\begin{bew}
By lemma \ref{26}, $\GGG_F$ has no branches.\qed
\end{bew}

\begin{satz}\label{440} Let $\tilde{\FFF}$ be a negative foundation. Then the following holds:
\begin{enumerate}[label=(\alph*)]
\item If $\tilde{\AA}$ is a field, then $\tilde{\FFF}$ is integrable.
\item If $\tilde{\AA}$ is a non-commutative skew-field and $\GGG_{\tilde{F}}$ is a string, a ray, a chain or a circle, then $\tilde{\FFF}$ is integrable.
\end{enumerate}
\end{satz}

\begin{bew} By theorem \ref{145}, $\tilde{\FFF}$ is integrable if its universal cover $\UUU$ is integrable. But since $\GGG_U$ is a tree, $\UUU$ is isomorphic to the corresponding canonical foundation by lemma \ref{446}, which is integrable by theorem \ref{146}. \qed
\end{bew}

\newpage

\chapter{Conclusion}

\newglossaryentry{CSimply}{type=results,name={{Classification of Simply Laced Twin Buildings}},description={},sort=res}

\begin{satz}[\textbf{\gls{CSimply}}]
Let $\FFF$ be an irr. simply laced foundation. Then $\FFF$ is integrable if and only if one of the following holds:
\begin{itemize}
\item The defining field is an octonion division algebra $\OO$, and $\FFF$ is isomorphic to one of the following foundations:
\begin{center}\begin{tikzpicture}[scale=0.7,>=stealth,thick]
\begin{scope}
\coordinate (1) at (0,0);                    
\coordinate (2) at (3,0);
\draws{0.55} (1)--(2);
\node () at (1.5,-0.5) {\ssi$\TTT(\OO)$};
\node () at (-1,2.5) {$\AAA_2(\OO)$:};
\fill (1) circle (3pt);
\fill (2) circle (3pt);
\end{scope}

\begin{scope}[xshift=6cm]
\coordinate (1) at (0,0);                    
\coordinate (2) at (3,0);
\coordinate (3) at (60:3);
\draws{0.55} (1)--(2);
\draws{0.55} (2)--(3);
\draws{0.55} (3)--(1);
\node()at (0.55,0.25) {\ssi{$\id_\OO$}};
\node()at (2.45,0.25) {\ssi$\id_\OO$};
\node()at (1.5,1.9) {\ssi$\id_\OO$};
\node()at (1.5,-0.5) {\ssi{$\TTT(\OO)$}};
\node()at (3.3,1.5) {\ssi{$\TTT(\OO)$}};
\node()at (-0.3,1.5) {\ssi{$\TTT(\OO)$}};
\node () at (-1,2.5) {$\tilde{\AAA}_2(\OO)$:};
\node()at (5.0,1.5) {\normalsize{.}};
\node()at (-2.0,1.5) {};
\draw[<-](60:2)arc[radius=1,start angle=240, end angle=300];
\draw[<-] (1,0) arc[radius=1,start angle=0, end angle=60];
\draw[->] (2,0) arc[radius=1,start angle=180, end angle=120];
\fill (1) circle (3pt);
\fill (2) circle (3pt);
\fill (3) circle (3pt);
\end{scope}
\end{tikzpicture}\end{center}
\item The defining field is a quaternion division algebra $\HH$, and $\MMM(\FFF)$ satisfies the conditions (i)-(iv) of theorem \ref{79}.

\newglossaryentry{An}{type=foundations,name={\ensuremath{\AAA_n(\DD)}},description=foundation of type $\AAA_n$ with respect to the skew-field $\DD$, sort=simply laced}
\newglossaryentry{Al}{type=foundations,name={\ensuremath{\AAA_\infty^l(\DD)}},description=ray foundation with respect to the skew-field $\DD$, sort=simply laced}
\newglossaryentry{Ar}{type=foundations,name={\ensuremath{\AAA_\infty^r(\DD)}},description=ray foundation with respect to the skew-field $\DD$, sort=simply laced}
\newglossaryentry{Ac}{type=foundations,name={\ensuremath{\AAA_\infty(\DD)}},description=chain foundation with respect to the skew-field $\DD$, sort=simply laced}
\newglossaryentry{Ant}{type=foundations,name={\ensuremath{\tilde{\AAA}_n(\DD,\gamma_i)}},description=foundation of type $\tilde{\AAA}_n$ with respect to the skew-field $\DD$ and the automorphisms ${\gamma_1,\ldots,\gamma_n\in\Aut(\DD)}$, sort=simply laced}

\item The defining field is a non-commutative skew-field $\DD$ different from a quaternion algebra, and $\FFF$ is isomorphic to one of the following foundations (where $\gamma_1,\ldots,\gamma_{n+1}\in \Aut(\DD)$):
\begin{center}\begin{tikzpicture}[scale=0.7,>=stealth,thick]
\begin{scope}
\coordinate (1) at (0,0);                    
\coordinate (2) at (3,0);
\coordinate (3) at (6,0);
\coordinate (4) at (9,0);
\draw[->] (2,0) arc[radius=1,start angle=180, end angle=0];
\draw[->] (5,0) arc[radius=1,start angle=180, end angle=0];
\draw[dashed] (2)--(3);
\draws{0.55} (1)--(2);
\draws{0.55} (3)--(4);
\node () at (7.5,-0.5) {\ssi$\TTT(\DD)$};
\node () at (1.5,-0.5) {\ssi$\TTT(\DD)$};
\node () at (-4,1) {$\gls{An}$:};
\node () at (0,-0.5) {\ssi$1$};
\node () at (3,-0.5) {\ssi$2$};
\node () at (6,-0.5) {\ssi$n-1$};
\node () at (9,-0.5) {\ssi$n$};
\node () at (3,0.5) {\ssi$\id_\DD$};
\node () at (6,0.5) {\ssi$\id_\DD$};
\foreach \i in {1,2,3,4} {\fill (\i) circle (3pt);}
\end{scope}

\begin{scope}[yshift=-3cm]
\coordinate (1) at (0,0);                    
\coordinate (2) at (3,0);
\coordinate (3) at (6,0);
\coordinate (4) at (9,0);
\coordinate (5) at (-2,0);
\draw[->] (2,0) arc[radius=1,start angle=180, end angle=0];
\draw[->] (5,0) arc[radius=1,start angle=180, end angle=0];
\draw[->] (-1,0) arc[radius=1,start angle=180, end angle=0];
\draw[dashed] (5)--(1);
\draws{0.55} (1)--(2);
\draws{0.55} (2)--(3);
\draws{0.55} (3)--(4);
\node () at (7.5,-0.5) {\ssi$\TTT(\DD)$};
\node () at (1.5,-0.5) {\ssi$\TTT(\DD)$};
\node () at (-4,1) {$\gls{Al}$:};
\node () at (0,-0.5) {\ssi$-3$};
\node () at (3,-0.5) {\ssi$-2$};
\node () at (6,-0.5) {\ssi$-1$};
\node () at (9,-0.5) {\ssi$0$};
\node () at (0,0.5) {\ssi$\id_\DD$};
\node () at (3,0.5) {\ssi$\id_\DD$};
\node () at (6,0.5) {\ssi$\id_\DD$};
\foreach \i in {1,2,3,4} {\fill (\i) circle (3pt);}
\end{scope}

\begin{scope}[yshift=-6cm]
\coordinate (1) at (0,0);                    
\coordinate (2) at (3,0);
\coordinate (3) at (6,0);
\coordinate (4) at (9,0);
\coordinate (5) at (11,0);
\draw[->] (2,0) arc[radius=1,start angle=180, end angle=0];
\draw[->] (5,0) arc[radius=1,start angle=180, end angle=0];
\draw[->] (8,0) arc[radius=1,start angle=180, end angle=0];
\draw[dashed] (4)--(5);
\draws{0.55} (1)--(2);
\draws{0.55} (2)--(3);
\draws{0.55} (3)--(4);
\node () at (7.5,-0.5) {\ssi$\TTT(\DD)$};
\node () at (1.5,-0.5) {\ssi$\TTT(\DD)$};
\node () at (-4,1) {$\gls{Ar}$:};
\node () at (0,-0.5) {\ssi$0$};
\node () at (3,-0.5) {\ssi$1$};
\node () at (6,-0.5) {\ssi$2$};
\node () at (9,-0.5) {\ssi$3$};
\node () at (9,0.5) {\ssi$\id_\DD$};
\node () at (3,0.5) {\ssi$\id_\DD$};
\node () at (6,0.5) {\ssi$\id_\DD$};
\foreach \i in {1,2,3,4} {\fill (\i) circle (3pt);}
\end{scope}

\begin{scope}[yshift=-9cm]
\coordinate (1) at (0,0);                    
\coordinate (2) at (3,0);
\coordinate (3) at (6,0);
\coordinate (4) at (9,0);
\coordinate (5) at (11,0);
\coordinate (6) at (-2,0);
\draw[->] (-1,0) arc[radius=1,start angle=180, end angle=0];
\draw[->] (2,0) arc[radius=1,start angle=180, end angle=0];
\draw[->] (5,0) arc[radius=1,start angle=180, end angle=0];
\draw[->] (8,0) arc[radius=1,start angle=180, end angle=0];
\draw[dashed] (4)--(5);
\draw[dashed] (6)--(1);
\draws{0.55} (1)--(2);
\draws{0.55} (2)--(3);
\draws{0.55} (3)--(4);
\node () at (7.5,-0.5) {\ssi$\TTT(\DD)$};
\node () at (1.5,-0.5) {\ssi$\TTT(\DD)$};
\node () at (-4,1) {$\gls{Ac}$:};
\node () at (0,-0.5) {\ssi$-1$};
\node () at (3,-0.5) {\ssi$0$};
\node () at (6,-0.5) {\ssi$1$};
\node () at (9,-0.5) {\ssi$2$};
\node () at (9,0.5) {\ssi$\id_\DD$};
\node () at (3,0.5) {\ssi$\id_\DD$};
\node () at (6,0.5) {\ssi$\id_\DD$};
\node () at (0,0.5) {\ssi$\id_\DD$};
\foreach \i in {1,2,3,4} {\fill (\i) circle (3pt);}
\end{scope}

\begin{scope}[yshift=-12cm]
\coordinate (1) at (3,0);                    
\coordinate (2) at (6,0);
\coordinate (3) at ($(2)+(-60:3)$);
\coordinate (4) at ($(3)+(-120:3)$);
\coordinate (5) at ($(4)+(-3,0)$);
\coordinate (6) at ($(5)+(120:3)$);
\draw[<-] (4,0) arc[radius=1,start angle=0, end angle=-120];
\draw[->] (5,0) arc[radius=1,start angle=180, end angle=300];
\draw[->] ($(3)+(120:1)$) arc[radius=1,start angle=120, end angle=240];
\draw[<-] ($(4)+(-1,0)$) arc[radius=1,start angle=180, end angle=60];
\draw[->] ($(5)+(1,0)$) arc[radius=1,start angle=0, end angle=120];
\draw[->] ($(6)+(-60:1)$) arc[radius=1,start angle=-60, end angle=60];
\draw[dashed] (4)--(5);
\draws{0.55} (6)--(1);
\draws{0.55} (1)--(2);
\draws{0.55} (2)--(3);
\draws{0.55} (3)--(4);
\draws{0.55} (5)--(6);
\node () at (-4,1) {$\gls{Ant}$:};
\node () at (3,0.5) {\ssi$n+1$};
\node () at (6,0.5) {\ssi$1$};
\node () at ($(3)+(0.5,0)$) {\ssi$2$};
\node () at ($(4)+(0,-0.5)$) {\ssi$3$};
\node () at ($(5)+(0,-0.5)$) {\ssi$n-1$};
\node () at ($(6)+(-0.5,0)$) {\ssi$n$};
\node () at (4.5,0.5) {\ssi$\TTT(\DD)$};
\node () at ($(3)+(0.3,1.5)$) {\ssi$\TTT(\DD)$};
\node () at ($(3)+(0.3,-1.5)$) {\ssi$\TTT(\DD)$};
\node () at ($(6)+(-0.3,1.5)$) {\ssi$\TTT(\DD)$};
\node () at ($(6)+(-0.3,-1.5)$) {\ssi$\TTT(\DD)$};
\node () at ($(2)+(-120:0.6)$) {\ssi$\gamma_1$};
\node () at ($(3)+(-180:0.5)$) {\ssi$\gamma_2$};
\node () at ($(4)+(120:0.5)$) {\ssi$\gamma_3$};
\node () at ($(5)+(60:0.5)$) {\ssi$\gamma_{n-1}$};
\node () at ($(6)+(0:0.5)$) {\ssi$\gamma_{n}$};
\node () at ($(1)+(-60:0.6)$) {\ssi$\gamma_{n+1}$};
\foreach \i in {1,...,6} {\fill (\i) circle (3pt);}
\end{scope}

\end{tikzpicture}\end{center}
\item The defining field is a field, and there are no further restrictions on the foundation $\FFF$.
\end{itemize}
\end{satz}

\newpage

\begin{bew}
This results from theorems \ref{444}, \ref{79}, \ref{443}, \ref{439} and \ref{440}.\qed
\end{bew}

\begin{bem}\ 
\begin{enumerate}[label=(\alph*)]
\item Given a non-commutative skew-field $\DD$, we have $$\AAA_\infty^l(\DD)\cong \AAA_\infty^r(\DD)\ \Leftrightarrow\ \DD\cong \DD^o\ .$$
\item The theorem can be stated in a more precise way by using classifying invariants. For example, given a skew-field $\DD$, the foundation $\tilde{\AAA}_n(\DD,\gamma_1,\ldots,\gamma_{n+1})$ only depends on the coset
$$\left(\prod_{i=1}^{n+1}\gamma_i\right)\cdot \mathrm{Inn}(\DD)\in \Aut(\DD)/\mathrm{Inn}(\DD)\ .$$
\item Each integrable mixed foundation over a quaternion division algebra $\HH$ arises in the following way: Start with some integrable positive foundations, then add some negative chains or strings of arbitrary length with the rule that each vertex of a positive foundation is part of at most one negative chain or string, e.g.,
\begin{center}\begin{tikzpicture}[scale=0.7,>=stealth,thick]

\begin{scope}
\coordinate (1) at (0,0);                    
\coordinate (2) at (0,-3);
\coordinate (3) at (-30:3);
\coordinate (4) at (0,-3.5);
\coordinate (5) at (0,-6.5);
\coordinate (6) at ($(4)+(-30:3)$);
\coordinate (7) at ($(3)+(6,-0.25)$);
\coordinate (8) at ($(6)+(6,0.25)$);
\coordinate (9) at ($(7)+(3,0)$);
\coordinate (10) at ($(8)+(3,0)$);
\draw[red] (1)--(2)--(3)--cycle;
\draw[red] (4)--(5)--(6)--cycle;
\draw[red] (7)--(8)--(9)--(10)--cycle;
\draw[red] (7)--(9);
\draw[red] (8)--(10);
\node[red] () at (0.9,-1.5) {$+$};
\node[red] () at (0.9,-5) {$+$};
\node[red] () at ($(8)+(45:2.12)$) {$+$};
\node () at (-2,0) {(i)};
\foreach \i in {1,...,10} {\fill (\i) circle (3pt);}
\end{scope}

\begin{scope}[yshift=-9cm]
\coordinate (1) at (0,0);                    
\coordinate (2) at (0,-3);
\coordinate (3) at (-30:3);
\coordinate (4) at (0,-3);
\coordinate (5) at (0,-6);
\coordinate (6) at ($(4)+(-30:3)$);
\coordinate (7) at ($(3)+(6,0)$);
\coordinate (8) at ($(6)+(6,0)$);
\coordinate (9) at ($(7)+(3,0)$);
\coordinate (10) at ($(8)+(3,0)$);
\coordinate (11) at ($(3)+(3,0)$);
\coordinate (12) at ($(6)+(-60:3)$);
\coordinate (13) at ($(8)+(-120:3)$);
\coordinate (14) at (3,0);
\coordinate (15) at (6,0);
\coordinate (16) at ($(9)+(150:3)$);
\coordinate (17) at (0,-9);
\coordinate (19) at (0,-11);
\coordinate (18) at ($(10)+(0,-3)$);
\draw[blue] (1,0) arc[radius=1,start angle=0, end angle=-30];
\draw[blue] ($(9)+(150:1)$) arc[radius=1,start angle=150, end angle=180];
\draw[blue] ($(3)+(1,0)$) arc[radius=1,start angle=0, end angle=150];
\draw[blue] ($(7)+(-1,0)$) arc[radius=1,start angle=180, end angle=0];
\draw[blue] ($(2)+(30:1)$) arc[radius=1,start angle=30, end angle=-30];
\draw[blue] ($(6)+(-60:1)$) arc[radius=1,start angle=-60, end angle=-150];
\draw[blue] ($(8)+(1,0)$) arc[radius=1,start angle=0, end angle=-120];
\draw[blue] ($(10)+(-1,0)$) arc[radius=1,start angle=180, end angle=270];
\draw[blue] ($(5)+(0,-1)$) arc[radius=1,start angle=-90, end angle=30];
\draw[red] (1)--(2)--(3)--cycle;
\draw[red] (4)--(5)--(6)--cycle;
\draw[red] (7)--(8)--(9)--(10)--cycle;
\draw[red] (7)--(9);
\draw[red] (8)--(10);
\draw[blue] (3)--(11)--(7);
\draw[blue] (6)--(12)--(13)--(8);
\draw[blue] (1)--(14)--(15)--(16)--(9);
\draw[blue] (5)--(17);
\draw[blue] (10)--(18);
\draw[dashed] (17)--(19);
\node[red] () at (0.8,-1.5) {$+$};
\node[red] () at (0.8,-4.5) {$+$};
\node[red] () at ($(8)+(45:2.12)$) {$+$};
\node () at (-2,0) {(ii)};
\foreach \i in {1,...,18} {\fill (\i) circle (3pt);}
\end{scope}

\end{tikzpicture}\end{center}
\end{enumerate}
\end{bem}

\addtocontents{toc}{\protect\newpage}
\addtocontents{toc}{\noindent\protect\mbox{}\protect\hrulefill\par}
\part{Jordan Automorphisms of Alternative Division Rings}
\addtocontents{toc}{\noindent\protect\mbox{}\protect\hrulefill\par}

\noindent In this part, we determine the structure of the group $\Aut_J(\AA)$ of Jordan automorphisms for an alternative division ring $\AA$. If $\AA$ is a skew-field, we know the answer by Hua's theorem. Since a non-associative alternative division ring is an octonion division algebra by the Bruck-Kleinfeld theorem \ref{490}, it remains to consider the octonion case. It turns out that we just have to add a certain class of Jordan automorphisms, each of them fixing a quaternion subalgebra.

One basic tool is the possibility to extend isomorphisms between subalgebras to an automorphism of the whole algebra, see \cite{S} for a detailed reference. Moreover, the crucial thing is the fact that nothing can go wrong with Jordan homomorphisms on fields, i.e., they are just monomorphisms of rings. This is not true for skew-fields, cf. §132.\\

\chapter{Composition Algebras and Norm Similarities}\label{134}
The theorems of this chapter provide the existence of automorphisms which we will need in §\ref{131}.

\begin{de}
A \textit{composition algebra over a field $\mathit{\KK}$}\index{composition algebra} is a unital algebra $\AA$ over $\KK$ together with a non-defective quadratic form $N:\AA\to \KK$ which \textit{permits composition}, i.e., we have
$$\forall\ x,y\in \AA:\qquad N(xy)=N(x)N(y)\ .$$
\end{de}

\begin{lemma} An octonion division algebra $\OO$ is a composition algebra over $\KK:=Z(\OO)$.
\end{lemma}

\begin{bew}
The norm $N$ is non-defective by corollary \ref{126} and  multiplicative by (9.9)(iii) of \cite{TW}.\qed
\end{bew}

\begin{de}
For $i=1,2$, let $V_i$ be a vector space over $\KK_i$ with non-defective quadratic form $N_i$, and let $\sigma:\KK_1\to \KK_2$ be an isomorphism of fields.
\begin{itemize}
\item A \textit{$\mathit{\sigma}$-similarity}\index{$\sigma$-similarity} is an isomorphism $(\varphi,\sigma):(V_1,\KK_1)\to (V_2,\KK_2)$ of vector spaces such that
$$\forall\ v\in V_1:\qquad N_2\big(\p(v)\big)=\rho_\p\cdot \sigma\big(N_1(v)\big)$$
for some element $\rho_\p\in \KK_2^*$, which is the \textit{multiplier of $\mathit{\p}$}\index{multiplier}.
\item A \textit{similarity}\index{similarity} is a $\sigma$-similarity such that $\KK_2=\KK_1$ and $\sigma=\id_{\KK_1}$.
\item A \textit{$\mathit{\sigma}$-isometry}\index{$\sigma$-isometry} is a $\sigma$-similarity such that $\rho_\p=1_{\KK_2}$.
\item An \textit{isometry}\index{isometry} is a $\sigma$-isometry such that $\KK_2=\KK_1$ and $\sigma=\id_{\KK_1}$.
\end{itemize}
\end{de}

\begin{lemma} \label{110}
For $i=1,2$, let $V_i$ be a vector space over $\KK_i$ with non-defective quadratic form $N_i$ and associated bilinear form $\langle \cdot,\cdot \rangle_i$, and let $\p:V_1\to V_2$ be a $\sigma$-similarity. Then $\p$ satisfies
$$\forall\ x,y\in V_1:\qquad \langle \p(x),\p(y)\rangle_2=\rho_\p\cdot \sigma_\p(\langle x,y\rangle_1)\ .$$
In particular, we have $\p(M^\bot)=\p(M)^\bot$ for each subset $M\subseteq V_1$.
\end{lemma}

\begin{bew}
Given $x\in M$ and $y\in M^\bot$, we have
\begin{align*}
\langle \p(x),\p(y)\rangle_2&= N_2\big(\p(x)+\p(y)\big)-N_2\big(\p(x)\big)-N_2\big(\p(y)\big)\\
&=\rho_\p\cdot \sigma_\p\big(N_1(x+y)-N_1(x)-N_1(y)\big)=\rho_\p\cdot \sigma_\p(\langle x,y\rangle_1)=0_{\KK_2}\ .
\end{align*}
\qed
\end{bew}

\newglossaryentry{NS}{type=symbols,name={\ensuremath{\GL_N(V,\KK)}},description=group of norm similarities with respect to the quadratic form $N$, sort=quadratic space}

\begin{no}
Let $V$ be a vector space over $\KK$ with non-defective quadratic form $N$. We denote the \textit{group of $\mathit{\sigma}$-isometries of $\mathit{V}$} by
$$\gls{NS}:=\left\{ (\varphi,\sigma)\in \GL(V) \mid \forall\ v\in V:\ N\big(\varphi(v)\big)=\sigma\big(N(v)\big)\right\}\ .$$
\end{no}

\begin{satz}\label{117}
Let $\AA$ be a composition algebra over $\KK$ and let $\sigma\in \Aut(\KK)$. Then there exists a $\sigma$-automorphism $\phi\in \Aut(\AA,\KK)$ if and only if there exists a $\sigma$-isometry $\p\in \GL_N(\AA,\KK)$.
\end{satz}

\begin{bew}
This is corollary (1.7.2) of [S].\qed
\end{bew}

\begin{satz}\label{118}
Let $\AA$ be a composition algebra over $\KK$, let $\sigma\in \Aut(\KK)$ and let $\mathbb{B}_1,\mathbb{B}_2$ be subalgebras of the same dimension. If there exists a $\sigma$-isometry $\p\in \GL_N(\AA,\KK)$, then each $\sigma$-isomorphism $\phi:\mathbb{B}_1\to \mathbb{B}_2$ of algebras can be extended to a $\sigma$-automorphism $\tilde{\phi}\in \Aut(\AA,\KK)$. In particular, each linear isomorphism $\psi:\mathbb{B}_1\to \mathbb{B}_2$ of algebras can be extended to a linear automorphism $\tilde{\psi}\in \Aut_\KK(\AA,\KK)$.
\end{satz}

\begin{bew}
This is corollary (1.7.3) of \cite{S}.\qed
\end{bew}

\chapter{Jordan Homomorphisms}
At this point we recall the definition of Jordan homomorphisms which play a central role in the classification of simply laced twin buildings. As the classification involves some detailed calculations concerning Jordan homomorphisms, it is natural to use this knowledge for a generalization of Hua's theorem to octonions.

\section{Basic Definitions and Basic Properties}

\begin{de}\label{356}\ 
\begin{itemize}  
\item Let $\AA,\tilde{\AA}$ be alternative division rings. A \textit{Jordan homomorphism} is an additive monomorphism $\gamma:\AA\to \tilde{\AA}$ such that
\begin{align*}
\gamma(1_{\AA})=1_{\tilde{\AA}}\ , && \forall\ x,y\in \AA:\ \gamma(xyx)=\gamma(x)\gamma(y)\gamma(x)\ .\end{align*}
\item Given an alternative division ring $\AA$, we denote the \textit{group of Jordan automorphisms of $\AA$} by
$$\gls{jora}:=\{ \gamma:\AA\to \AA \mid \gamma\ \textrm{Jordan automorphism}\}\ .$$
\end{itemize}
\end{de}

\begin{lemma}\label{100}
Let $\AA,\tilde{\AA}$ be alternative division rings. A Jordan homomorphism $\gamma:\AA\to\tilde{\AA}$ satisfies
\begin{align*}
\forall\ x,y\in \AA:\qquad \gamma(xy)+\gamma(yx)=\gamma(x)\gamma(y)+\gamma(y)\gamma(x)\ .
\end{align*}
\end{lemma}
 
\begin{bew}
The assertion is clearly true for $x=0_\AA$ or $y=0_\AA$, so assume $x\neq0_\AA\neq y$. As we have
$$\gamma(z^2)=\gamma(z\cdot 1_\AA\cdot z)=\gamma(z)\gamma(1_{{\AA}})\gamma(z)=\gamma(z)\cdot 1_{\tilde{\AA}}\cdot\gamma(z)=\gamma(z)^2$$
for each $z\in \AA$, it follows that
\begin{align*}
\gamma\big((x+y)^2\big)&=\gamma(x^2+xy+yx+y^2)=\gamma(x)^2+\gamma(xy)+\gamma(yx)+\gamma(y)^2\ ,\\
\gamma\big((x+y)^2\big)&=\gamma(x+y)^2=\big(\gamma(x)+\gamma(y)\big)^2=\gamma(x)^2+\gamma(x)\gamma(y)+\gamma(y)\gamma(x)+\gamma(y)^2
\end{align*}
and thus $$\gamma(xy)+\gamma(yx)=\gamma(x)\gamma(y)+\gamma(y)\gamma(x)\ .$$\qed
\end{bew}

\begin{lemma}\label{127} Let $\AA,\tilde{\AA}$ be alternative division rings. A Jordan homomorphism $\gamma:\AA\to\tilde{\AA}$ satisfies
$$\forall\ x\in \AA^*:\qquad \gamma(x^{-1})=\gamma(x)^{-1}\ .$$
\end{lemma}

\begin{bew}
Given $x\in \AA^*$, we have
\begin{align*}
\gamma(x^{-1})=\gamma(x^{-1}xx^{-1})=\gamma(x^{-1})\gamma(x)\gamma(x^{-1})
\end{align*}
and thus
$\gamma(x^{-1})=\gamma(x)^{-1}$ by lemma \ref{124}.\qed
\end{bew}

\section{Jordan Homomorphisms on Fields}

A Jordan homomorphism from a field $\FF$ to an alternative division ring $\AA$ is a monomorphism of rings. In particular, the image of $\AA$ is a subfield of $\AA$.

\begin{lemma}\label{102}
Let $\FF$ be a field, let $\AA$ be an alternative division ring and let $\gamma:\FF\to \AA$ be a Jordan homomorphism. Then $\gamma:\FF\to \gamma(\FF)$ is an isomorphism of fields.
\end{lemma}

\begin{bew}
By Hua's theorem, we may assume that $\AA$ is an octonion division algebra. By lemma \ref{100}, we have
 \begin{align} 2\gamma(xy)=\gamma(x)\gamma(y)+\gamma(y)\gamma(x) \label{101} \end{align}
for all $x,y\in\FF$.
\begin{itemize}
\item $\Char \FF\neq 2$: Given $x,y\in\FF^*$, we have
$$\gamma(xy)\gamma(x^{-1}y^{-1})\overset{\textrm{\ref{124}}}{=}\gamma(xy)\gamma\big((xy)^{-1}\big)\overset{\textrm{\ref{127}}}{=}\gamma(xy)\gamma(xy)^{-1}=1_\AA$$ 
and thus by equation \eqref{101} and lemma \ref{111}
\begin{align*}
4&=2\gamma(xy)2\gamma(x^{-1}y^{-1})=[\gamma(x)\gamma(y)+\gamma(y)\gamma(x)][\gamma(x^{-1})\gamma(y^{-1})+\gamma(y^{-1})\gamma(x^{-1})] \\
&=\gamma(x)\gamma(y)\gamma(x)^{-1}\gamma(y)^{-1}+1_\AA+1_\AA+\gamma(y)\gamma(x)\gamma(y)^{-1}\gamma(x)^{-1}\ .
\end{align*}
We set $z:=\gamma(x)\gamma(y)\gamma(x)^{-1}\gamma(y)^{-1}$ and obtain
\begin{align*}
2=z+z^{-1}\ , &&  (z-1_\AA)^2=z^2-2z+1_\AA=0_\AA
\end{align*}
and therefore
\begin{align*}  \gamma(x)\gamma(y)\gamma(x)^{-1}\gamma(y)^{-1}=z=1_\AA\ , && \gamma(x)\gamma(y)=\gamma(y)\gamma(x)\ .\end{align*}
From equation \eqref{101} it follows that $\gamma(xy)=\gamma(x)\gamma(y)\ .$
\item $\Char \FF=2$: In this case, equation \eqref{101} implies 
$$\gamma(x)\gamma(y)=\gamma(y)\gamma(x)$$
for all $x,y\in \FF$. As a consequence,
$$\tilde{\FF}:=\langle \gamma(\FF)\rangle_{Z(\AA)}$$
is a commutative subalgebra of $\AA$, hence a field. Since $\gamma:\FF\to \tilde{\FF}$ is a Jordan homomorphism, Hua's theorem implies that $\gamma:\FF\to \gamma(\FF)$ is an iso- or anti-isomorphism of skew-fields, and thus, in fact, an isomorphism of fields.
\end{itemize}\qed
\end{bew}

\newpage

\chapter{Jordan Automorphisms}
In this chapter, we determine the structure of $\Aut_J(\AA)$ for an alternative division ring $\AA$. If $\AA$ is a skew-field, then Hua's theorem gives the answer:

\begin{satz}
Let $\AA$ be skew-field. Then we have
$$\Aut_J(\AA)=\Aut(\AA)\cup \Aut^o(\AA)\ .$$
\end{satz}

\noindent Thus it remains to consider the Jordan automorphisms of an octonion division algebra $\OO$.

\begin{no} Throughout the rest of this chapter, $\OO$ is an octonion division algebra and $\KK:=Z(\OO)$ is its center.\end{no}

\section{Jordan Automorphisms on Subfields}\label{132}

As a consequence of the last paragraph, a Jordan automorphism restricted to a subfield is a monomorphism of rings. In particular, the image of a subfield is again a subfield.

This is not true for skew-subfields. Otherwise, a Jordan automorphism would be, in fact, an auto- or anti-automorphism by the proof of Hua's theorem. But there are indeed Jordan automorphism neither of the first nor of the second kind, cf. lemma \ref{352} and remark \ref{351}.

\begin{lemma}\label{349}
Let $\gamma\in \Aut_J(\OO)$ and let $\FF$ be a subfield of $\OO$. Then $\gamma_{|\FF}:\FF\to \gamma(\FF)$ is an isomorphism of fields.
\end{lemma}

\begin{bew}
This results from Lemma \ref{102}.\qed
\end{bew}

\begin{kor}\label{103}
An element $\gamma\in \Aut_J(\OO)$ satisfies
$$\gamma(\lambda)\cdot \gamma(x)=\gamma(x)\cdot \gamma(\lambda)$$
for all $\lambda\in \KK$, $x\in\OO$. 
\end{kor}

\begin{bew}
Since each element $x\in \OO$ is contained in a subfield $\FF$ of $\OO$ with $\KK\subseteq \FF$ by (20.9) of \cite{TW}, we may apply lemma \ref{349} to obtain
$$\gamma(\lambda)\cdot \gamma(x)=\gamma(\lambda\cdot x)=\gamma(x\cdot \lambda)=\gamma(x)\cdot \gamma(\lambda)$$
for all $\lambda\in \KK$, $x\in\OO$.\qed
\end{bew}

\section{Jordan Automorphisms and Norm Similarities}\label{130}

The results of §\ref{132} enable us to show that Jordan automorphisms are norm similarities. As a consequence, we may apply the results of chapter \ref{134}. Moreover, the results of \cite{J} show the reverse inclusion for $\Char \OO\neq 2$.

\begin{prop}\label{398}
Let $\gamma\in \Aut_J(\OO)$ and $\sigma:=\gamma_{|\KK}$. Then the following holds:
\begin{enumerate}[label=(\alph*)]
\item \label{104} $\sigma\in \Aut(\KK)$.
\item $(\gamma,\sigma)\in \GL(\OO,\KK)$.
\item \label{106} $(\gamma,\sigma)\in \GL_N(\OO,\KK)$.
\end{enumerate}
\end{prop}

\newpage

\begin{bewzwei}\ 
\begin{enumerate}[label=(\alph*)]
\item By corollary \ref{103} and lemma \ref{123}, we have $\gamma(\KK)=\KK$ and therefore $\sigma\in \Aut(\KK)$ by lemma \ref{102}.
\item This is a consequence of (a) and corollary \ref{103}.
\item By (a), the assertion is true for $x\in \KK$. Given $x\in\OO\sm \KK$, we have
\begin{align*} x^2-T(x)x+N(x)=0_\OO\ , &&\gamma(x)^2-T\big(\gamma(x)\big)\gamma(x)+N\big(\gamma(x)\big)=0_\OO
\end{align*} and hence
\begin{align*}
0_\OO&=\gamma^{-1}\big(\gamma(x)^2-T(\gamma(x))\gamma(x)+N(\gamma(x))\big)=\gamma^{-1}\gamma\big(x^2-\gamma^{-1}\big(T(\gamma(x))\big)x+\gamma^{-1}(N(\gamma(x)))\big) \\
&=x^2-\gamma^{-1}\big(T(\gamma(x))\big)x+\gamma^{-1}\big(N(\gamma(x))\big)\ .
\end{align*}
As the maps $T$ and $N$ are uniquely determined by the minimum equation, we obtain
\begin{align*}
\gamma^{-1}\big(N(\gamma(x))\big)=N(x)\ , && N\big(\gamma(x)\big)=\gamma\big(N(x)\big)=\sigma\big(N(x)\big)
\end{align*}
for all $x\in \OO\sm \KK$.
\end{enumerate}\qed
\end{bewzwei}

\begin{kor}\label{109}
We have
$$\sigma_s\in Z\big(\Aut_J(\OO)\big)\ .$$
\end{kor}

\begin{bew}
Given $\gamma\in \Aut_J(\OO)$ and $x\in \OO$, we have
\begin{align*}
\gamma\sigma_s(x)=\gamma(\bar{x})=\gamma\big(N(x)\cdot x^{-1}\big)=\sigma\big(N(x)\big)\cdot \gamma(x^{-1})=N\big(\gamma(x)\big)\cdot \gamma(x)^{-1}=\overline{\gamma(x)}=\sigma_s\gamma(x)\ .
\end{align*}\qed
\end{bew}

\begin{bem}
The results of \cite{J} are valid for Jordan algebras with characteristic different from 2, cf. definition (1.3) of \cite{J}.
\end{bem}

\begin{lemma}\label{105}  If $\Char \OO\neq 2$, the Jordan algebra $\OO$ is separable
\end{lemma}

\begin{bew}
By corollary \ref{126}, the bilinear form $\langle\cdot,\cdot \rangle$ and thus the trace form
$$\bar{T}:\OO\times \OO\to \KK,\ (x,y)\mapsto T(xy)=xy+\bar{y}\bar{x}=\langle x,\bar{y}\rangle$$ is non-degenerate. Now the assertion results from theorem (6.5) of \cite{J}.
\qed
\end{bew}

\begin{prop}\label{120}
If $\Char \OO\neq 2$, we have $$O(\OO,\KK)=GL_N(\OO,\KK)\subseteq \Aut_J(\OO)\ .$$
\end{prop}

\begin{bew}
Let $\p=(\p,\id_\KK)\in GL_N(\OO,\KK)$. As we have $\p(1_\OO)=\p(1_\OO\cdot 1_\OO)=1_\OO\cdot \p(1_\OO)$, the isometry $\p$ satisfies $\p(1_\OO)=1_\OO$, and because of $|\KK|>2=\deg \OO$, the assertion results from lemma \ref{105} and theorem (6.7)(a) of \cite{J}.\qed
\end{bew}

\begin{satz}
If $\Char \OO\neq 2$, we have
$$\GL_N(\OO,\KK)= \Aut_J(\OO)\ .$$
\end{satz}

\newpage

\begin{bewzwei}\ 
\begin{itemize}
\item[``$\subseteq $''] Let $(\p,\sigma)\in \GL_N(\OO,\KK)$. By theorem \ref{117}, there is a $\sigma$-automorphism $(\phi,\sigma)\in \Aut(\OO,\KK)$. By proposition \ref{120}, we have
$$(\phi^{-1},\sigma^{-1})\cdot (\p,\sigma)\in GL_N(\OO,\KK)\subseteq \Aut_J(\OO)$$
and thus
$$(\p,\sigma)\in (\phi,\sigma)\cdot \Aut_J(\OO)=\Aut_J(\OO)\ .$$
\item[``$\supseteq$''] This is proposition \ref{106}.
\end{itemize}
\qed
\end{bewzwei}

\section[The Structure of \texorpdfstring{$\Aut_J(\OO)$}{the Group of Jordan Automorphisms}]{The Structure of \texorpdfstring{$\boldsymbol{\Aut_J(\OO)}$}{the Group of Jordan Automorphisms}}\label{131}
Now we are ready to tackle the main problem which we split up into three steps. First of all we construct a subgroup $\Gamma\leq \Aut_J(\OO)$. The second step shows that we may suppose $\gamma\in \Aut_J(\OO)$ to fix a quaternion subalgebra pointwise. Then we finally prove that $\gamma$ is the product of elements of the given group $\Gamma$.

\begin{no}\label{412}\ 
\begin{itemize}
\item Given a quaternion subalgebra $\HH$, $e\in\HH^\bot$ and $w,p\in \HH$, we set
\begin{align*} 
\psi_{(\HH,e,w)}&:\OO\to \OO,\ x+e\cdot y\mapsto x+e\cdot w^{-1}yw\ , \\
\phi_{(\HH,e,w,p)}&:\OO\to \OO,\ x+e\cdot y\mapsto w^{-1}xw+e\cdot w^{-1}ywp\ . \end{align*}
\item We set
\begin{align*} \Psi&:=\{ \psi_{(\HH,e,w)} \mid \HH\ \textrm{a quaternion subalgebra},\ e\in \HH^\bot,\ w\in\HH\}\ , \\
               \Phi&:=\{ \phi_{(\HH,e,w,p)} \mid \HH\ \textrm{a quaternion subalgebra},\ e\in \HH^\bot,\ w,p\in\HH,\ N(p)=1_\KK\}\ .
\end{align*}
\item We set
\begin{align*} \Gamma:=\{ \psi\phi\mid \psi\in\Psi,\ \phi\in \Aut(\OO)\}\ .\end{align*}
\end{itemize}
\end{no}

\begin{bem}\label{351}
By remark (20.29) of \cite{TW}, a map $\psi\in\Psi$ is neither an auto- nor an anti-automorphism.
\end{bem}

\begin{lemma}\label{352} We have
$$\Psi\cup\Phi\subseteq \Gamma\leq \Aut_J(\OO)\ .$$
\end{lemma}
\begin{bew}
This results from lemma \ref{466} and lemma \ref{119}.\qed
\end{bew}

\begin{lemma}\label{112}
Let $\AA$ be a subalgebra of $\OO$ such that $\AA^\bot \nsubseteq \AA$ and let $e\in \AA^\bot\sm \AA$. Then we have
$$(e\cdot x)(e\cdot y)(e\cdot x)=-N(e)e\cdot x\bar{y}x$$
for all $x,y\in \AA$.
\end{lemma}

\begin{bew}
This results from lemma \ref{250}.\qed
\end{bew}

\begin{bem}
The following lemma is helpful since we know that Jordan homomorphisms on subfields are in fact isomorphisms between subfields and that Jordan automorphisms of skew-fields are iso- or anti-isomorphisms.
\end{bem}

\begin{lemma}\label{113}
Let $\AA$ be a subalgebra of $\OO$ such that $\AA^\bot \nsubseteq \AA$, let $e\in \AA^\bot\sm \AA$ and let $\gamma\in \Aut_J(\OO)\cap GL_N(\OO,\KK)$ such that $\gamma(e)=e$. Then the map $\tilde{\gamma}:\AA\to \OO$ defined by
$$\gamma(e\cdot x)=e\cdot \tilde{\gamma}(x)$$
is a linear Jordan homomorphism.
\end{lemma}

\begin{bew} Notice that we have $$\forall\ x\in \AA:\qquad \big(e,\tilde{\gamma}(a)\big)=\big(e,-N(e)^{-1}e\cdot \gamma(e\cdot a)\big)=-\big(1_\KK,\gamma(e\cdot a)\big)=0_\KK$$ by lemma \ref{116} and lemma \ref{110}, hence $e\in \tilde{\gamma}(\AA)^\bot\sm \tilde{\gamma}(\AA)$.
\begin{itemize}
\item We have
$$e\cdot \tilde{\gamma}(1_\OO)=\gamma(e\cdot 1_\OO)= \gamma(e)=e\cdot 1_\OO\ .$$
\item Given $\lambda\in \KK,\ x\in \AA$, we have
\begin{align*}
e\cdot \tilde{\gamma}(\lambda x)=\gamma(e\cdot \lambda x)=\lambda\gamma(e\cdot x)=\lambda\big(e\cdot \tilde{\gamma}(x)\big)=e\cdot\lambda\tilde{\gamma}(x)\ .
\end{align*}
\item Given $x,y\in \AA$, we have
\begin{align*}
e\cdot \tilde{\gamma}(x+y)=\gamma\big(e\cdot (x+y)\big)=\gamma(e\cdot x)+\gamma(e\cdot y)=e\cdot \big(\tilde{\gamma}(x)+\tilde{\gamma}(y)\big)\ .
\end{align*}
\item Given $x,y\in \AA$, we have
\begin{align*}
e\cdot \tilde{\gamma}(xyx)&=\gamma(e\cdot xyx)=-N(e)^{-1}\cdot \gamma\big( (e\cdot x)(e\cdot \bar{y})(e\cdot y)\big)\\
&=-N(e)^{-1}\cdot \big(e\cdot \tilde{\gamma}(x)\big)\big(e\cdot \tilde{\gamma}(\bar{y})\big)\big(e\cdot \tilde{\gamma}(x)\big)=e\cdot \tilde{\gamma}(x)\overline{\tilde{\gamma}(\bar{y})}\tilde{\gamma}(x)
\end{align*}
by lemma \ref{112}. We set $x:=1_\OO$ to obtain
\begin{align*}
\tilde{\gamma}(y)=\overline{\tilde{\gamma}(\bar{y})}
\end{align*}
for each $y\in \AA$ and thus $$\tilde{\gamma}(xyx)=\tilde{\gamma}(x)\tilde{\gamma}(y)\tilde{\gamma}(x)$$ for all $x,y\in \AA$.
\end{itemize}\qed
\end{bew}

\begin{de}
Let $e_1,e_2\in \OO^*$ such that $\lambda:=N(e_1)\neq 0_\KK, \mu:=N(e_2)\neq 0_\KK$. Then $(e_1,e_2)$ is a \textit{special $\mathit{(\lambda,\mu)}$-pair}\index{special $(\lambda,\mu)$-pair} if
\begin{itemize}
\item we have 
$$\langle e_1,1_\OO\rangle=\langle e_2,1_\OO\rangle=\langle e_1,e_2\rangle=0_\KK$$ for $\Char \OO\neq 2$;
\item we have \begin{align*} \langle e_1,1_\OO\rangle=1_\KK\ , && \langle e_2,1_\OO\rangle=\langle e_1,e_2\rangle=0_\KK\end{align*}
for $\Char \OO=2$.
\end{itemize}
\end{de}

\begin{bem}
We will need special $(\lambda,\mu)$-pairs to extend isomorphisms between quaternion subalgebras to the whole octonion division algebra.
\end{bem}

\newpage

\begin{lemma}\label{128} Let $(e_1,e_2)$ be a special $(\lambda,\mu)$-pair and $\EE:=\langle 1_\OO,e_1\rangle_\KK$. Then we have
\begin{align*} \bar{e}_1\neq e_1\ , && e_2\in \EE^\bot\sm \EE\ .\end{align*}
\end{lemma}

\begin{bewzwei}\ 
\begin{itemize}
\item $\Char \OO\neq2$: We have
\begin{align*}
e_1+\bar{e}_1=\langle e_1,1_\OO\rangle=0_\KK\ , && \bar{e}_1=-e_1\neq e_1\ .
\end{align*}
Let $x=s+te_1\in \EE^*$. If $s\neq 0_\KK$, we have
$$\langle x,1_\OO\rangle= \langle s+te_1,1_\OO\rangle=s\cdot \langle e_1,e_1\rangle=2s\cdot N(e_1)\neq 0_\KK$$
and thus $x\notin \EE^\bot$. If $s=0_\KK$ (and thus $t\neq 0_\KK$), we have
$$\langle x,e_1\rangle=\langle te_1,e_1\rangle=t\cdot \langle e_1,e_1\rangle=2t\cdot N(e_1)\neq 0_\KK$$
and thus $x\notin \EE^\bot$.
\item $\Char \OO= 2$: We have
\begin{align*} 
e_1+\bar{e}_1=\langle e_1,1_\OO\rangle=1_\KK\ , && \bar{e}_1=e_1+1_\KK\neq e_1\ .
\end{align*}
Let $x=s+te_1\in \EE^*$. If $t\neq 0_\KK$, we have
$$\langle x, 1_\OO\rangle =\langle s+te_1,1_\OO\rangle=s\cdot \langle 1_\OO,1_\OO\rangle+t\cdot \langle e_1,1_\OO\rangle=2s\cdot N(1_\OO)+t=t\neq 0_\KK$$
and thus $x\notin \EE^\bot$. If $t=0_\KK$ (and thus $s\neq 0_\KK$), we have
$$\langle e_1,x\rangle=\langle e_1,s\rangle=s\cdot \langle e_1,1_\OO\rangle=s\neq 0_\KK$$
and thus $x\notin \EE^\bot$.
\end{itemize}
\qed
\end{bewzwei}

\begin{bem}
The following lemma allows us to choose a suitable $\KK$-basis for a quaternion division algebra containing two elements $x,y\in \OO$.
\end{bem}

\begin{lemma}\label{135}
Let $x\in \OO\sm \KK$, let $\EE:=\langle 1_\OO, x\rangle_\KK$ and let $y\in \OO\sm \EE$. Then $$\{1_\OO,x,y,xy\}$$ is linearly independent over $\KK$.
\end{lemma}

\begin{bew}
Notice that $\EE$ is a field. Let $a,b,c,d\in \KK$ such that
$$a\cdot 1_\OO+b\cdot x+c\cdot y+d\cdot xy=0_\OO\ .$$
Then we have
$$(c\cdot 1_\OO+d\cdot x)\cdot y=-a\cdot 1_\OO-b\cdot x\ .$$
\begin{itemize}
\item $c\cdot 1_\OO+d\cdot x\neq 0_\OO$: In this case, we have
$$y=(c\cdot 1_\OO+d\cdot x)^{-1}\cdot (-a\cdot 1_\OO-b\cdot x)\in \EE \qquad \lightning\ .$$
\item $c\cdot 1_\OO+d\cdot x= 0_\OO$: In this case, we have
\begin{align*}
c\cdot 1_\OO+d\cdot x=0_\OO=a\cdot 1_\OO+b\cdot x\ , && c=d=0_\KK=a=b\ .
\end{align*}
\end{itemize}
\qed
\end{bew}

\begin{prop}\label{107} Given an element $\gamma\in \Aut_J(\OO)$, there is an element $\tilde{\gamma}\in\langle \sigma_s,\Gamma\rangle$ such that $\tilde{\gamma}\gamma$ fixes a quaternion subalgebra $\HH$ pointwise.
\end{prop}

\newpage

\begin{bewzwei}\ 
\begin{enumerate}[label=(\roman*)]
\item Let $\sigma:=\gamma_{|\KK}\in \Aut(\KK)$. Then $\gamma$ is a $\sigma$-isometry by proposition \ref{106}, thus $\gamma^{-1}$ is a $\sigma^{-1}$-isometry. By theorem \ref{117}, there is a $\sigma^{-1}$-automorphism $\bar{\gamma}\in \Aut(\OO,\KK)$. As a consequence, we have $\bar{\gamma}\gamma\in \Aut_J(\OO)\cap GL_N(\OO,\KK)$.
\item By lemma \ref{114}, there is an element $e_1\in \OO$ such that $\bar{e}_1\neq e_1$. Let $\EE:=\langle 1,e_1\rangle_\KK$ and $e_2\in \EE^\bot\sm \EE$. By lemma \ref{115}, $e_1$ and $e_2$ are contained in a quaternion subalgebra $\tilde{\HH}$ which contains a $(\lambda,\mu)$-special pair $(\tilde{e}_1,\tilde{e}_2)$ by definition (1.7.4) of \cite{S}. Thus we may assume that $(e_1,e_2)$ is $(\lambda,\mu)$-special by lemma \ref{128}.
\item As we may suppose $\gamma$ to be an isometry by (i), the pairs $(e_1,e_2)$ and $\big(\gamma(e_1),\gamma(e_2)\big)$ are $(\lambda,\mu)$-special. By corollary (1.7.5) of \cite{S}, there is a linear automorphism $\bar{\gamma}\in \Aut_\KK(\OO,\KK)$ extending
\begin{align*} \gamma(e_1)\mapsto e_1\ , && \gamma(e_2)\mapsto e_2\ .\end{align*}
As a consequence, $\gamma:=\bar{\gamma}\gamma$ fixes $\langle \EE,e_2\rangle_\KK$ pointwise.
\item By lemma \ref{113}, the map $\tilde{\gamma}:\EE\to \tilde{\gamma}(\EE)$ defined by
$$\gamma(e_2\cdot x)=e_2\cdot \tilde{\gamma}(x)$$
is a linear Jordan isomorphism and thus an isomorphism of fields by lemma \ref{102}.
\item If $\tilde{\gamma}(\EE)=\EE$, we have
$$\gamma(\EE+e_2\cdot \EE)=\EE+e_2\cdot \EE=:\HH\ .$$
In this case, $\gamma_{|\HH}$ is an auto- or anti-automorphism by Hua's theorem.
\begin{itemize}
\item If $\gamma_{|\HH}$ is an automorphism, we have
\begin{align*} \gamma(e_2\cdot e_1)=\gamma(e_2)\cdot \gamma(e_1)=e_2\cdot e_1\ , &&\gamma_{|\HH}=\id_\HH\ .\end{align*}
\item If $\gamma_{|\HH}$ is an anti-automorphism, then $\phi:=\sigma_s\gamma_{|\HH}$ is a linear automorphism. By theorem \ref{118}, we may extend $\phi$ to a linear automorphism $\tilde{\phi}\in \Aut_\KK(\OO,\KK)$. Then
$$\tilde{\phi}^{-1}\sigma_s\gamma$$
fixes $\HH$ pointwise.
\end{itemize}
Thus we may assume $\tilde{\gamma}(e_1)\notin \EE$.
\item By lemma \ref{115}, $e_1$ and $\tilde{\gamma}(e_1)$ are contained in a quaternion subalgebra $\tilde{\HH}$, and by (v) and lemma \ref{135}, we have 
$$\tilde{\HH}=\langle 1_\OO, e_1,\tilde{\gamma}(e_1),e_1\tilde{\gamma}(e_1)\rangle_\KK\ ,$$
thus we may extend $\tilde{\gamma}:\EE\to \tilde{\gamma}(\EE)$ to a linear automorphism $\phi\in \Aut_\KK(\tilde{\HH},\KK)$ by theorem \ref{118}. By the Skolem-Noether theorem, there is an element $w\in \tilde{\HH}$ such that
$$\phi=\gamma_w\ .$$
\item 
We have $e_2\cdot \EE\subseteq \EE^\bot$ and therefore
\begin{align*}
\langle e_2,1_\OO\rangle&=\langle e_2,e_1\rangle=0_\KK\ .
\end{align*} Moreover, we have
$$e_2\cdot \phi(\EE)=\gamma(e_2\cdot \EE)\subseteq \gamma(\EE^\bot)=\gamma(\EE)^\bot=\EE^\bot$$
by corollary \ref{110} and therefore
\begin{align*} \langle e_2, e_1\phi(e_1)\rangle&=\langle e_2\phi(e_1)^{-1}\phi(e_1),e_1\phi(e_1)\rangle= \langle e_2\phi(e_1)^{-1},e_1\rangle\cdot N\big(\phi(e_1)\big)=0_\KK\ , \\
\langle e_2,\phi(e_1)\rangle&=\langle e_2, -N(e_2)^{-1}e_2\cdot e_2\phi(e_1)\rangle=-N(e_2)^{-1}N(e_2)\cdot\langle 1_\OO, e_2\phi(e_1)\rangle=0_\KK\ ,
\end{align*}
hence $e_2\in \tilde{\HH}^\bot$, and $\gamma_{(\tilde{\HH},e_2,w^{-1})}\gamma$ fixes $\HH:=\EE+e_2\cdot \EE$ pointwise, cf. notation \ref{412}.
\end{enumerate}\qed
\end{bewzwei}

\newpage

\begin{prop}\label{108}
Let $\gamma\in \Aut_J(\OO)$ be a Jordan automorphism fixing a quaternion subalgebra $\HH$ pointwise. Then we have $\gamma\in \langle \sigma_s,\Gamma\rangle$.
\end{prop}

\begin{bew}
Let $e\in \HH^\bot\sm \HH$.
\begin{enumerate}[label=(\roman*)]
\item  As we have
$$\gamma(e\cdot \HH)=\gamma(\HH^\bot)=\gamma(\HH)^\bot=\HH^\bot =e\cdot \HH\ ,$$
by corollary \ref{110}, we may define a map $\tilde{\gamma}:\HH\to \HH$ via 
$$\gamma(e\cdot x)=e\cdot \tilde{\gamma}(x)\ .$$
\item We have
$$ N(e)=N\big(\gamma(e\cdot 1_\OO)\big)=N\big(e\cdot \tilde{\gamma}(1_\OO)\big)=N(e)\cdot N\big(\tilde{\gamma}(1_\OO)\big)$$
and thus $N\big(\tilde{\gamma}(1_\OO)\big)=1_\KK$. Therefore, the map
$$ \phi_{(\HH,e,1_\OO,\tilde{\gamma}(1_\OO)^{-1})}\gamma$$
fixes $\langle \HH,e\rangle_\KK$ pointwise.
\item By lemma \ref{113}, the map $\tilde{\gamma}$ is a linear Jordan automorphism and thus an auto- or an anti-automorphism by Hua's theorem. 
\item \begin{itemize}
\item $\tilde{\gamma}$ is an automorphism: As $\tilde{\gamma}$ is linear, the Skolem-Noether theorem yields an element $w\in \HH$ such that $\tilde{\gamma}=\gamma_w$, hence
$$\gamma=\gamma_{(\HH,e,w)}\in\Gamma\ .$$
\item $\tilde{\gamma}$ is an anti-automorphism: Then there is an element $w\in \HH$ such that
$$\gamma(x+e\cdot y)=x+e\cdot w^{-1}\bar{y}w$$
for all $x,y\in \HH$, hence
$$\gamma_e\gamma(x+e\cdot y)=e^{-1}xe+w^{-1}\bar{y}w\cdot e=\bar{x}+e\cdot w^{-1}yw$$
for all $x,y\in\HH$ and therefore
\begin{align*} \sigma_s\phi_{(\HH,e,1,-1)}\gamma_e\gamma=\gamma_{(\HH,e,w)}\in\Gamma\ , \gamma\in \langle \sigma_s,\Gamma\rangle\ .
\end{align*}\qed
\end{itemize}\end{enumerate}
\end{bew}

\section{Conclusion}

\newglossaryentry{JAD}{type=results,name={{Jordan Automorphisms of Octonion Division Algebras}},description={},sort=res}

\begin{satz}[\textbf{\gls{JAD}}]\label{399}
Given an octonion division algebra $\OO$, we have
$$\Aut_J(\OO)=\langle \sigma_s,\Gamma\rangle\cong \langle \sigma_s\rangle\times \Gamma\ .$$
\end{satz}

\begin{bew}
The first equality results from proposition \ref{107} and proposition \ref{108}, the second assertion from corollary \ref{109}.\qed
\end{bew}

\newpage

\addtocontents{toc}{\noindent\protect\mbox{}\protect\hrulefill\par}
\part{Moufang Sets}
\addtocontents{toc}{\noindent\protect\mbox{}\protect\hrulefill\par}

\noindent Now we turn to the general description of the root groups of Moufang Polygons. In fact, all of them are parametrized by Moufang sets, more precisely, by their associated groups. As the glueings of integrable foundations turn out to be Jordan isomorphisms, the classification of twin buildings is closely related to the solution of the isomorphism problem for Moufang sets.

We list the examples of Moufang sets which will appear in the sequel, then we give a complete overview of the Jordan isomorphisms between these Moufang sets before we give the missing proofs.\\

\chapter{Basic Definitions}\label{457}

\begin{de}
A \textit{Moufang set}\index{Moufang set} is a pair $\MM=\big( X, \{U_x\}_{x\in X}\big)$ consisting of a set $X$ with $|X|\geq 3$ and a set of \textit{root groups}\index{root group} $\{U_x\}_{x\in X}$ satisfying the following conditions:
\begin{enumerate}[label=(M\arabic*),leftmargin=30pt]
\item For each $x\in X$, the group $U_x\leq \Sym(X)$ fixes $x$ and acts regularly on $X\sm \{x\}$.
\item For each $x\in X$ and for each $\p\in \langle U_y \mid y\in X\rangle$, we have $U_x^\p=U_{\p(x)}$. 
\end{enumerate}
\end{de}

\begin{bem}
Let $(U,+)$ be a not necessarily commutative group, let $X:=U\cup \{\infty\}$ be the disjoint union of $U$ and $\{\infty\}$ and let $\tau\in \Sym(X)$ be a permutation interchanging $0$ and $\infty$, which means that we have $\tau_{|U^*}\in \Sym(U^*)$. By theorem 2 of \cite{DW}, the pair $(U,\tau)$ gives rise to a Moufang set $\MM(U,\tau)$ if and only if we have $h_a\in \Aut(U)$ for each $a\in U^*$, where $h_a$ is the \textit{Hua map with respect to $\mathit{a}$} as in definition 2 of \cite{DW}, more precisely, we consider $h_a$ to be the restriction to $U$ of that map given there. Conversely, each Moufang set $\MM$ arises in such a way, cf page 5 of \cite{DW}, or lemma 1.3.4 of \cite{DS} for a more precise statement,.

As both the descriptions are equivalent, we consider a Moufang set to be a pair $\MM=(U,\tau)$ consisting of a not necessarily commutative group $U=(U,+)$ and an element $\tau\in \Sym(U^*)$ such that $h_a\in \Aut(U)$ for each $a\in U^*$.
\end{bem}

\begin{de}\ 
\begin{itemize}
\item A Moufang set $\MM=(U,\tau)$ is \textit{commutative}\index{Moufang set!commutative} if the group $U$ is commutative.
\item A Moufang set $\MM=(U,\tau)$ is \textit{unital}\index{Moufang set!unital}\index{unital!Moufang set} if there is an element $1_\MM\in U^*$ such that $h_{1_\MM}=\id_U$.
\end{itemize}
\end{de}

\begin{bem}
If a Moufang set is unital, the element $1_\MM$ is not necessarily uniquely determined by the defining property. However, we will distinguish a canonical element for the examples we mainly deal with, cf. the next paragraph.
\end{bem}

\begin{de}
Let $\MM=(U,\tau)$, $\tilde{\MM}=(\tilde{U},\tilde{\tau})$ be Moufang sets. 
\begin{itemize} 
\item An \textit{isomorphism}\index{isomorphism!of Moufang sets} $\p:\MM\to\tilde{\MM}$ is an isomorphism of groups $\p:U\to \tilde{U}$ such that
$$\forall\ x\in U:\qquad \p\big(\tau(x)\big)=\tilde{\tau}\big(\p(x)\big)\ .$$
\item An \textit{automorphism of $\mathit{\MM}$}\index{automorphism!of a Moufang set} is an isomorphism $\p:\MM\to \MM$.
\item A \textit{Jordan isomorphism}\index{Jordan isomorphism!of Moufang sets} $\gamma:\MM\to\tilde{\MM}$ is an isomorphism of groups $\gamma:U\to \tilde{U}$ such that
$$\forall\ x\in U,\ a\in U^*:\qquad \gamma\big(h_a(x)\big)=\tilde{h}_{\gamma(a)}\big(\gamma(x)\big)\ ,$$
and, moreover, such that $\gamma(1_\MM)=1_{\tilde{\MM}}$ if $\MM$ and $\tilde{\MM}$ both are unital.
\item A \textit{Jordan automorphism of $\mathit{\MM}$} is a Jordan isomorphism $\gamma:\MM\to \MM$.
\end{itemize}
\end{de}

\begin{bem}
The list in the following chapter is not complete, we only list those Moufang sets appearing in triangles and quadrangles since we only classify foundations involving polygons of this type. Moreover, we exclude the non-commutative Moufang sets appearing in quadrangles of type $E_n$.
\end{bem}

\chapter{Examples}\label{458}

\begin{bsp}[\textbf{Moufang Sets of Linear Type}]
Given an alternative division ring $\AA$, the corresponding \textit{Moufang set of linear type}\index{Moufang set!of linear type} is
\begin{align*}
\MM(\AA):=(\AA,\tau)\ , && \tau:\AA^*\to \AA^*,\ x\mapsto -x^{-1}\ .
\end{align*}
Given $a\in \AA^*$, the Hua map with respect to $a$ is
$$h_a:\AA\to\AA,\ x\mapsto axa\ .$$
As a consequence, $\MM(\AA)$ is unital with $1_\MM=1_\AA$.
\end{bsp}

\begin{bsp}[\textbf{Moufang Sets of Involutory Type}]
Given an involutory set $(\KK,\KK_0,\sigma)$, the corresponding \textit{Moufang set of involutory type}\index{Moufang set!of involutory type} is
\begin{align*}
\MM(\KK,\KK_0,\sigma):=(\KK_0,\tau)\ , && \tau:\KK_0^*\to \KK_0^*,\ x\mapsto -x^{-1}\ .
\end{align*}
Given $a\in \KK_0^*$, the Hua map with respect to $a$ is
$$h_a:\KK_0\to\KK_0,\ x\mapsto axa\ .$$
As a consequence, $\MM(\KK,\KK_0,\sigma)$ is unital with $1_\MM=1_\KK\in \KK_0$.
\end{bsp}

\begin{bsp}[\textbf{Moufang Sets of Indifferent Type}]
Given an indifferent set $(\KK,\KK_0,\LL_0)$, the corresponding \textit{Moufang set of indifferent type}\index{Moufang set!of indifferent type} is
\begin{align*}
\MM(\KK,\KK_0,\LL_0):=(\KK_0,\tau)\ , && \tau:\KK_0^*\to \KK_0^*,\ x\mapsto -x^{-1}\ .
\end{align*}
Given $a\in \KK_0^*$, the Hua map with respect to $a$ is
$$h_a:\KK_0\to\KK_0,\ x\mapsto axa\ .$$
As a consequence, $\MM(\KK,\KK_0,\LL_0)$ is unital with $1_\MM=1_\KK\in \KK_0$.
\end{bsp}

\begin{bsp}[\textbf{Moufang Sets of Quadratic Form Type}]
Given a quadratic space $(L_0,\KK,q)$ with basepoint $\epsilon$, the corresponding \textit{Moufang set of quadratic form type with basepoint $\mathit{\epsilon}$}\index{Moufang set!of quadratic form type} is
\begin{align*}
\MM(L_0,\KK,q):=(L_0,\tau)\ , && \tau: L_0^*\to L_0^*,\ a\mapsto -a^\sigma\cdot q(a)^{-1}\ .
\end{align*}
Given $a\in L_0^*$, the Hua map with respect to $a$ is
$$h_a:L_0\to L_0,\ v\mapsto \pi_a\pi_\epsilon(v)\cdot q(a)\ .$$
As a consequence, $\MM(L_0,\KK,q)$ is unital with $1_\MM=\epsilon$.
\end{bsp}

\begin{bsp}[\textbf{Moufang Sets of Pseudo-Quadratic Form Type}]\ \\
Let $\Xi=(\KK,\KK_0,\sigma,L_0,q)$ be a pseudo-quadratic space and let $T=T(\Xi)$ be the group as in definition \ref{190} and corollary \ref{321}. The corresponding \textit{Moufang set of pseudo-quadratic form type}\index{Moufang set!of pseudo-quadratic form type} is
\begin{align*} \MM(\Xi):=(T,\tau)\ , && \tau:T^*\to T^*,\ (a,t)\mapsto (at^{-1},-t^{-1})\ .
\end{align*}
Given $(a,t)\in T^*$, the Hua map with respect to $(a,t)$ is
$$h_{(a,t)}: T\to T,\ (b,v)\mapsto \big(b t^\sigma-at^{-1}f(a,b)t^\sigma,t vt^\sigma\big)\ .$$
As a consequence, $\MM(\Xi)$ is unital with $1_\MM=(0,1_\KK)$.
\end{bsp}

\chapter{The Isomorphism Problem for Moufang Sets}\label{416}

Since the glueings appearing in a foundation are in fact Jordan isomorphisms, it is natural to solve the isomorphism problem for the appearing Moufang sets before we tackle the classification of foundations.

As mentioned, we consider Jordan isomorphisms, not isomorphisms of Moufang sets in the proper sense, which are a subset of the Jordan isomorphisms. It would be interesting in which cases both the definitions coincide, which should be correct in almost all the cases.

First of all we give an overview of the results, some of them already proved in the previous parts, then we give the missing proofs.

\section{Results}

\begin{de}\label{323} Let $\OO$ be an octonion division algebra.
\begin{itemize}
\item Given a quaternion subalgebra $\HH$, $e\in\HH^\bot$ and $w,p\in \HH$, we set
\begin{align*} 
\psi_{(\HH,e,w)}&:\OO\to \OO,\ x+e\cdot y\mapsto x+e\cdot w^{-1}yw\ .
\end{align*}
\item We set 
\begin{align*} \Psi&:=\{ \psi_{(\HH,e,w)} \mid \HH\ \textrm{a quaternion subalgebra},\ e\in \HH^\bot,\ w\in\HH\}\ . 
\end{align*}
and  $\Gamma:=\{ \psi\phi\mid \psi\in\Psi,\ \phi\in \Aut(\OO)\}$, which is a subgroup of $\Aut_J(\OO)$ by lemma \ref{119}.
\end{itemize}
\end{de}

\newglossaryentry{Linear}{type=results,name={{Moufang Sets of Linear Type}},sort=isomorphism problem,description={}}

\begin{satz}[\textbf{\gls{Linear}}] \label{433}
Let $\MM:=\MM(\AA)$, $\tilde{\MM}:=\MM(\tilde{\AA})$ be Moufang sets of linear type. A map $\gamma:\MM\to \tilde{\MM}$ is a Jordan isomorphism such that $\gamma(1_\MM)=\gamma(1_{\tilde{\MM}})$ if and only if one of the following holds:
\begin{enumerate}[label=(\roman*)]
\item The alternative division rings $\AA$ and $\tilde{\AA}$ are skew-fields and $\gamma$ is an iso- or anti-isomorphism of skew-fields. 
\item The alternative division rings $\AA$ and $\tilde{\AA}$ are isomorphic octonion division algebras and we have $$\phi^{-1}\gamma\in \langle \sigma_s,\Gamma\rangle\cong \langle \sigma_s\rangle\times \Gamma\ ,$$
where $\phi:\AA\to \tilde{\AA}$ is an isomorphism of alternative rings, $\sigma_s$ is the standard involution of $\AA$ and $\Gamma$ is the group defined in \ref{323}.
\end{enumerate}
\end{satz}

\begin{bewzwei}\ 
\begin{itemize}
\item[``$\R$''] If $\AA$ is a skew-field, then (i) holds by Hua's theorem. If $\AA$ (and thus $\tilde{\AA})$ is an octonion division algebra, we may adapt proposition \ref{398} to obtain that $\gamma:\AA\to \tilde{\AA}$ is a $\sigma$-isometry. By theorem (1.7.1) of \cite{S}, there is a $\sigma$-isomorphism $\phi:\AA\to \tilde{\AA}$, hence $\phi^{-1}\gamma\in \Aut_J(\AA)$, and we may apply theorem \ref{399}.
\item[``$\Leftarrow$''] $\checkmark$
\end{itemize}
 Notice that definition \ref{467} and definition \ref{356} include the condition $\gamma(1_\AA)=1_{\tilde{\AA}}$.\qed
\end{bewzwei}

\newglossaryentry{Involutory}{type=results,name={{Moufang Sets of Involutory Type}},sort=isomorphism problem,description={}}

\begin{satz}[\textbf{\gls{Involutory}}]
Let $(\KK,\KK_0,\sigma)$ be a proper involutory set, let $(\tilde{\KK},\tilde{\KK}_0,\tilde{\sigma})$ be an involutory set and let $\MM:=\MM(\KK,\KK_0,\sigma)$, $\tilde{\MM}:=\MM(\tilde{\KK},\tilde{\KK}_0,\tilde{\sigma})$ be the corresponding Moufang sets of involutory type. A map $\gamma:\MM\to \tilde{\MM}$ is a Jordan isomorphism such that $\gamma(1_\MM)=1_{\tilde{\MM}}$ if and only if there is an isomorphism $\phi:(\KK,\KK_0,\sigma)\to (\tilde{\KK},\tilde{\KK}_0,\tilde{\sigma})$ of involutory sets such that $\gamma=\phi_{|\KK_0}$.
\end{satz}

\begin{bew}
This is theorem \ref{175}. Notice that the proof is not complete yet.\qed
\end{bew}

\newpage

\begin{de}\label{324}
Let $(\AA,\FF,\sigma)$ be quadratic of type (iv) and suppose that $\dim_\AA L_0=2$. By \cite{DM}, there are exactly three pseudo-quadratic spaces 
\begin{align*}
(\AA,\FF,\sigma,L_0,q)=(\AA_1,\FF,\sigma,L_0,q_1)=\Xi_1\ ,&& (\AA_2,\FF,\sigma,L_0,q_2)=\Xi_2\ , &&(\AA_3,\FF,\sigma,L_0,q_3)=\Xi_3
\end{align*}
with pairwise non-isomorphic quaternion division algebras $\AA_1,\AA_2,\AA_3$ which define the group $T$. When we switch between the parametrizing pseudo-quadratic spaces, we indicate this by the map
$$\id_T^i:T\to T,\ (a,t)\mapsto (a,t)\ ,$$
i.e., after applying $\id_T^i$, we consider $T$ to be defined by $\Xi_i$.
\end{de}

\newglossaryentry{PQ}{type=results,name={Moufang Sets of Pseudo-Quadratic Form Type},sort=isomorphism problem,description={}}

\begin{satz}[\textbf{\gls{PQ}}]
Let $\Xi$ and $\tilde{\Xi}$ be proper pseudo-quadratic spaces and let $\MM:=\MM(\Xi)$, $\tilde{\MM}:=\MM(\tilde{\Xi})$ be the corresponding Moufang sets of pseudo-quadratic form type. A map $\gamma:\MM\to\tilde{\MM}$ is a Jordan isomorphism such that $\gamma(1_\MM)=1_{\tilde{\MM}}$ if and only if one of the following holds:
\begin{enumerate}[label=(\roman*)]
\item There is an isomorphism $\Phi:\Xi\to \tilde{\Xi}$ of pseudo-quadratic spaces that induces $\gamma$.
\item The involutory sets $(\KK,\KK_0,\sigma)$ and $(\tilde{\KK},\tilde{\KK}_0,\tilde{\sigma})$ both are quadratic of type (iv), we have
\begin{align*}
\KK\not\cong \tilde{\KK}\ , && \dim_\KK L_0=2=\dim_{\tilde{\KK}} \tilde{L}_0
\end{align*}
and there are an $i\in \{2,3\}$ and an isomorphism  $\Phi:\Xi \to \tilde{\Xi}_i$ of pseudo-quadratic spaces such that $\gamma$ is induced by $(\id_{\tilde{T}}^i)^{-1}\circ \Phi$, where $\id_{\tilde{T}}^i$ and $\tilde{\Xi}=:\tilde{\Xi}_1,\tilde{\Xi}_2,\tilde{\Xi}_3$ are as in definition \ref{324}.
\item The involutory sets $(\KK,\KK_0,\sigma)$ and $(\tilde{\KK},\tilde{\KK}_0,\tilde{\sigma})$ are quadratic of type (iv) and (iii), respectively, we have $\dim_\KK L_0=1,\ \dim_{\tilde{\KK}} \tilde{L}_0=2$ and $\gamma$ can be described by
$$\forall\ x=s+et\in \KK,\ u\in \KK_0:\qquad \gamma( ax,x^\sigma q(a)x+u)= \big(\tilde{a}\phi(s)+\tilde{b}\phi(t)^{\tilde{\sigma}}, \phi\big(N(x)q(a)+u\big)\big)\ ,$$
where $a\in L_0^*$ is arbitrary, $\phi:\EE_a\to \tilde{\KK}$ is an isomorphism of fields, $e\in \EE_a^\bot$, $\tilde{a}\in \tilde{L}_0$ and $\tilde{b}\in \tilde{a}^\bot$.
\item We have $\KK\cong \FF_4\cong \tilde{\KK}$, $\dim_\KK L=1$ and there are an isomorphism $\Phi:\Xi\to \tilde{\Xi}$ of pseudo-quadratic spaces and a non-trivial inner automorphism $\tilde{\gamma}\in \Aut(\tilde{T})$ such that $\gamma$ is induced by $\tilde{\gamma}\circ \Phi$.
\end{enumerate}
\end{satz}

\begin{bew}
This is theorem \ref{353}. Notice that definition \ref{354} includes the condition $\gamma(0,1_\KK)=(0,1_{\tilde{\KK}})$.\qed
\end{bew}

\newglossaryentry{Quadratic}{type=results,name={Moufang Sets of Quadratic Form Type},sort=isomorphism problem,description={}}

\begin{satz}[\textbf{\gls{Quadratic}}]\label{400}
Let $\Xi$ be a quadratic space with basepoint $\epsilon$, let $\tilde{\Xi}$ be a proper quadratic space with basepoint $\tilde{\epsilon}$ and let $\MM:=\MM(\Xi)$, $\tilde{\MM}:=\MM(\tilde{\Xi})$ be the corresponding Moufang sets of quadratic form type. A map $\gamma:\MM\to \tilde{\MM}$ is a Jordan isomorphism such that $\gamma(1_\MM)=1_{\tilde{\MM}}$ if and only if one of the following holds:
\begin{enumerate}[label=(\roman*)]
\item We have $\dim_{\tilde{\KK}} \tilde{L}_0\geq 3$ and there is an isomorphism $\phi:\KK\to \tilde{\KK}$ of fields such that the map
$$(\gamma,\phi):(L_0,\KK,q)\to (\tilde{L}_0,\tilde{\KK},\tilde{q})$$
is an isomorphism of quadratic spaces. In particular, we have $\dim_{{\KK}}{L}_0=\dim_{\tilde{\KK}}\tilde{L}_0$.
\item We have $\dim_{\tilde{\KK}} \tilde{L}_0\leq 2$, the map $$\phi:\KK\to \hat{\KK}:=\gamma(\langle \epsilon\rangle_\KK)\subseteq \hat{\FF}:=\FF(\tilde{L}_0,\tilde{\KK},\tilde{q}),\ s\mapsto \gamma(\epsilon\cdot s)$$
is an isomorphism of fields, the field $\hat{\FF}$ is quadratic over $\hat{\KK}$, and the map 
$$(\gamma,\phi):(L_0,\KK,q)\to (\hat{\FF},\hat{\KK},N^{\hat{\FF}}_{\hat{\KK}})$$
is an isomorphism of quadratic spaces. This is true even if $\tilde{\Xi}$ is non-proper.
\end{enumerate}
\end{satz}

\begin{bew}
This is theorem \ref{434}.\qed
\end{bew}

\newglossaryentry{Qub}{type=results,name={Moufang Sets of Quadratic Form and Linear Type}, sort=isomorphism problem,description={}}

\begin{satz}[\textbf{\gls{Qub}}]\label{379}  Let $\MM:=\MM(L_0,\KK,q)$ be a Moufang set of quadratic form type with basepoint $\epsilon$ and let $\tilde{\MM}:=\MM(\tilde{\AA})$ be a Moufang set of linear type. A map $\gamma:\MM\to \tilde{\MM}$ is a Jordan isomorphism such that $\gamma(1_\MM)=1_{\tilde{\MM}}$ if and only if the map 
$$\phi:\KK\to \tilde{\KK}:=\gamma(\langle \epsilon\rangle_\KK)\subseteq \tilde{\AA},\ s\mapsto \gamma(\epsilon\cdot s)$$
is an isomorphism of fields, $\tilde{\AA}$ is quadratic over $\tilde{\KK}$ and the map $(\gamma,\phi): (L_0,\KK,q)\to (\tilde{\AA},\tilde{\KK},N^{\tilde{\AA}}_{\tilde{\KK}})$ is an isomorphism of quadratic spaces.
\end{satz}

\begin{bew}
This is theorem \ref{350}.\qed
\end{bew}

\newglossaryentry{IL}{type=results,name={{Moufang Sets of Indifferent and Linear Type}},sort=isomorphism problem,description={}}

\begin{satz}[\textbf{\gls{IL}}] Let $\MM:=\MM(\KK,\KK_0,\LL_0)$ be a Moufang set of indifferent type and let $\tilde{\MM}:=\MM(\tilde{\AA})$ be a Moufang set of linear type. A map $\gamma:\tilde{\MM}\to {\MM}$ is a Jordan isomorphism such that $\gamma(1_{\tilde{\MM}})=1_{{\MM}}$ if and only if $\KK_0=\KK$, $\tilde{\AA}$ is a field and the map $\gamma:\tilde{\AA}\to \KK_0$ is an isomorphism of fields. In particular, the indifferent set $(\KK,\KK_0,\LL_0)$ is non-proper if we have $\MM\cong \tilde{\MM}$.
\end{satz}

\begin{bew}
This is theorem \ref{388}.\qed
\end{bew}

\newglossaryentry{InL}{type=results,name={{Moufang Sets of Involutory and Linear Type}},sort=isomorphism problem,description={}}

\begin{satz}[\textbf{\gls{InL}}] Let $\tilde{\MM}:=\MM(\tilde{\AA})$ be a Moufang set of linear type and let $\MM:=\MM(\KK,\KK_0,\sigma)$ be a Moufang set of involutory type. A map $\gamma:\tilde{\MM}\to {\MM}$ is a Jordan isomorphism such that $\gamma(1_{\tilde{\MM}})=1_{{\MM}}$ if and only if $(\KK,\KK_0,\sigma)$ is of quadratic type, $\tilde{\AA}$ and $\KK_0$ are fields and the map $\gamma:\tilde{\AA}\to \KK_0$ is an isomorphism of fields. In particular, the involutory set $(\KK,\KK_0,\sigma)$ is non-proper if we have $\MM\cong \tilde{\MM}$.
\end{satz}

\begin{bew}
This is theorem \ref{389}.\qed
\end{bew}

\begin{bem}
Most of the following proofs or different versions can also be found in \cite{K}. Notice, however, that some of them only give the idea for the proof so that we had to work out some details, especially in the following paragraph.
\end{bem}

\section[\texorpdfstring{$\MM(L_0,\KK,q)\cong \MM(\tilde{\AA})$}{Quadratic Form Type and Linear Type}]{\texorpdfstring{$\boldsymbol{\MM(L_0,\KK,q)\cong \MM(\tilde{\AA})}$}{Quadratic Form Type and Linear Type}}

If a Moufang set $\MM(L_0,\KK,q)$ of quadratic form type with basepoint $\epsilon$ is isomorphic to a Moufang set $\MM(\tilde{\AA})$ of linear type, then $\tilde{\AA}$ is quadratic over the field $\tilde{\KK}:=\gamma(\langle \epsilon_\KK\rangle)$ and $\gamma$ is induced by an isomorphism $(\gamma,\phi):(L_0,\KK,q)\to (\tilde{\AA},\tilde{\KK},N^{\tilde{\AA}}_{\tilde{\KK}})$ of quadratic spaces.

\begin{lemma}\label{472}
Given a quadratic space $(L_0,\KK,q)$ such that $\dim_\KK L_0\leq 2$, we have $$\MM:=\MM(L_0,\KK,q)=\MM\big(\FF(L_0,\KK,q)\big)=:\tilde{\MM}\ .$$ In particular, the corresponding Hua maps coincide.
\end{lemma}

\begin{bew}
We have
$$\forall\ x\in L_0:\qquad \tau(x)=-x^\sigma\cdot q(x)^{-1}=-x^\sigma\ast \epsilon\cdot q(x)^{-1}=-x^\sigma\ast N^{\tilde{\FF}}_{\tilde{\KK}}(x)^{-1}=-x^{-1}=\tilde{\tau}(x)\ .$$\qed
\end{bew}

\newpage

\begin{no}
Until proposition \ref{340}, $\MM:=\MM(L_0,\KK,q)$ is a Moufang set of quadratic form type with basepoint $\epsilon$, $\tilde{\MM}:=\MM(\tilde{\AA})$ is a Moufang set of linear type and $\gamma:\MM\to\tilde{\MM}$ is a Jordan isomorphism such that $\gamma(1_\MM)=1_{\tilde{\MM}}$.
\end{no}

\begin{lemma}\label{332}
The map
$$\phi:\KK\to \tilde{\KK}:=\gamma(\langle \epsilon\rangle_\KK)\subseteq \tilde{\AA},\ s\mapsto \gamma(\epsilon\cdot s)$$
is an isomorphism of fields.
\end{lemma}

\begin{bew}
By lemma \ref{328}, we have
$$\forall\ s,t\in \KK:\qquad  \phi(sts)=\gamma(\epsilon\cdot sts)=\gamma\big(h_{\epsilon \cdot s}(\epsilon\cdot t)\big)=\tilde{h}_{\gamma(\epsilon\cdot s)}\big(\gamma(\epsilon\cdot t)\big)=\phi(s)\phi(t)\phi(s)\ .$$
As a consequence, $\phi:\KK\to \tilde{\AA}$ is a Jordan homomorphism, hence $\phi:\KK\to \phi(\KK)=\gamma(\langle \epsilon\rangle_\KK)$ is an isomorphism of fields by lemma \ref{102}.\qed
\end{bew}

\begin{no}
Given $x\in L_0\sm \langle \epsilon\rangle_\KK$, we set $R_x:=\langle \epsilon, x\rangle_\KK$.
\end{no}

\begin{bem}
Given $x\in L_0\sm \langle \epsilon\rangle_\KK$, the triple $(R_x,\KK,q)$ is a quadratic space such that $\dim_\KK R_x=2$.
\end{bem}

\begin{no}
Given $x\in L_0\sm \langle \epsilon\rangle_\KK$, we set $\FF_x:=\FF(R_x,\KK,q)$.
\end{no}

\begin{lemma}
Given $x\in L_0\sm \langle \epsilon\rangle_\KK$, the map
$$\gamma_{|\FF_x}:\FF_x\to \gamma(\FF_x)\subseteq \tilde{\AA},\ y\mapsto \gamma(y)$$
is an isomorphism of fields.
\end{lemma}

\begin{bew}
By lemma \ref{472}, the Hua maps of $\MM(R_x,\KK,q)$ and $\MM(\FF_x)$ coincide, hence $\gamma_{|\FF_x}:\FF_x\to \tilde{\AA}$ is a Jordan homomorphism so that we may apply lemma \ref{102}.\qed
\end{bew}

\begin{prop}
Let $x\in L_0$. Then the following holds:
\begin{enumerate}[label=(\alph*)]
\item \label{475} We have
$$\forall\ s\in \KK:\qquad \gamma(x\cdot s)=\gamma(x)\cdot \phi(s)\ .$$
In particular, the map $(\gamma,\phi):(L_0,\KK)\to (\tilde{\AA},\tilde{\KK})$ is an isomorphism of vector spaces.
\item \label{474} We have
$$\forall\ s\in \KK:\qquad \gamma(x)\cdot \phi(s)=\phi(s)\cdot \gamma(x)\ .$$
In particular, we have $\tilde{\KK}\subseteq Z(\tilde{\AA})$.
\end{enumerate}
\end{prop}

\begin{bew} If $x\in \langle \epsilon\rangle_\KK$, the assertions result from lemma \ref{332}, so assume $x\in L_0\sm \langle \epsilon\rangle_\KK$.
\begin{enumerate}[label=(\alph*)]
\item Given $s\in \KK$, we have
$$\gamma(x\cdot s)=\gamma\big(x\ast(\epsilon\cdot s)\big)=\gamma(x)\cdot \gamma(\epsilon\cdot s)=\gamma(x)\cdot \phi(s)\ .$$
\item Given $s\in \KK$, we have
$$\gamma(x)\cdot \phi(s)=\gamma(x\cdot s)=\gamma\big((\epsilon\cdot s)\ast x\big)=\gamma(\epsilon\cdot s)\cdot \gamma(x)=\phi(s)\cdot \gamma(x)\ .$$
\end{enumerate}\qed
\end{bew}

\newpage

\begin{lemma}\label{335}
Given $x\in L_0$, we have
\begin{align*}
h_x(\epsilon)-x\cdot T(x)+\epsilon\cdot q(x)=0_{L_0}\ .
\end{align*}
\end{lemma}

\begin{bew}
Let $x\in L_0$. By remark \ref{325}, we have
\begin{align*}
h_x(\epsilon)-x\cdot T(x)+\epsilon\cdot q(x)&=x\cdot f_q(x,\bar{\epsilon})-\bar{\epsilon}\cdot q(x)-x\cdot T(x)+\epsilon\cdot q(x) \\
&=x\cdot T(x)-x\cdot T(x)-\epsilon\cdot q(x)+\epsilon\cdot q(x)=0_{L_0}\ .
\end{align*}\qed
\end{bew}

\begin{prop}\label{340}
Given $x\in L_0^*$, we have
\begin{align*} \gamma(x)^2-\gamma(x)\cdot \phi\big(T(x)\big)+\phi\big(q(x)\big)=0_{\tilde{\AA}}\ .
\end{align*}
In particular, the alternative division ring $\tilde{\AA}$ is quadratic over $\tilde{\KK}$ with norm $\tilde{N}:=N^{\tilde{\AA}}_{\tilde{\KK}}$, satisfying
\begin{align*}
\tilde{N}\circ \gamma=\phi\circ q\ .
\end{align*}
\end{prop}

\begin{bew}
Given $x\in L_0^*$, we have
\begin{align*}
0_{\tilde{\AA}}&=\gamma(0_{L_0})=\gamma\big(h_x(\epsilon)-x\cdot T(x)+\epsilon\cdot q(x)\big)=\gamma\big(h_x(\epsilon)\big)-\gamma\big(x\cdot T(x)\big)+\gamma\big(\epsilon\cdot q(x)\big)\\
&=\tilde{h}_{\gamma(x)}\big(\gamma(\epsilon)\big)-\gamma(x)\cdot \phi\big(T(x)\big)+\phi\big(q(x)\big)=\gamma(x)^2-\gamma(x)\cdot \phi\big(T(x)\big)+\phi\big(q(x)\big)
\end{align*}
by lemma \ref{335} and proposition \ref{475}, which shows that $\tilde{\AA}$ is quadratic over $\tilde{\KK}$. Given $x\in L_0\sm \langle \epsilon\rangle_\KK$, we have $\gamma(x)\in \tilde{\AA}\sm \tilde{\KK}$ and thus
\begin{align*}
\tilde{N}\big(\gamma(x)\big)=\phi\big( q(x)\big)
\end{align*}
since the minimum equation is unique. Finally, given $s\in \KK$, we have
\begin{align*}
\tilde{N}\big(\gamma(\epsilon\cdot s)\big)&=\tilde{N}\big(\phi(s)\big)=\phi(s)^2=\phi(s^2)=\phi\big(q(\epsilon\cdot s)\big)
\end{align*}
by lemma \ref{332}.
\qed
\end{bew}

\begin{satz}\label{350}  Let $\MM:=\MM(L_0,\KK,q)$ be a Moufang set of quadratic form type with basepoint $\epsilon$ and let $\tilde{\MM}:=\MM(\tilde{\AA})$ be a Moufang set of linear type. A map $\gamma:\MM\to \tilde{\MM}$ is a Jordan isomorphism such that $\gamma(1_\MM)=1_{\tilde{\MM}}$ if and only if the map 
$$\phi:\KK\to \tilde{\KK}:=\gamma(\langle \epsilon\rangle_\KK)\subseteq \tilde{\AA},\ s\mapsto \gamma(\epsilon\cdot s)$$
is an isomorphism of fields, $\tilde{\AA}$ is quadratic over $\tilde{\KK}$ with norm $\tilde{N}:=N^{\tilde{\AA}}_{\tilde{\KK}}$ and the map
$$(\gamma,\phi): (L_0,\KK,q)\to (\tilde{\AA},\tilde{\KK},\tilde{N})$$
is an isomorphism of quadratic spaces.
\end{satz}

\begin{bewzwei}\ 
\begin{itemize}
\item[``$\R$''] By lemma \ref{332}, the map $\phi:\KK\to \tilde{\KK}$ is an isomorphism of fields. By proposition \ref{340} and proposition \ref{475}, $\tilde{\AA}$ is quadratic over $\tilde{\KK}$, and the map $(\gamma,\phi):(L_0,\KK,q)\to (\tilde{\AA},\tilde{\KK},\tilde{N})$ is an isomorphism of quadratic spaces.

\item[``$\Leftarrow$''] We have $\gamma(1_\MM)=\gamma(\epsilon)=\phi(\epsilon)=1_{\tilde{\AA}}=1_{\tilde{\MM}}$. By lemma \ref{476}, we have
\begin{align*}
\gamma\big(h_a(x)\big)&=\gamma\big(a\cdot f_q(a,x^\sigma)-x^\sigma\cdot q(a)\big)=\gamma(a)\cdot \phi\big(f_q(a,x^\sigma)\big)-\gamma(x^\sigma)\cdot \phi\big(q(a)\big) \\
&=\gamma(a)\cdot \phi\big(q(a+x^\sigma)-q(a)-q(x^\sigma)\big)-\gamma(x)^{\tilde{\sigma}}\cdot \tilde{N}\big(\gamma(a)\big) \\
&=\gamma(a)\cdot \big(\tilde{N}\big(\gamma(a)+\gamma(x^\sigma)\big)-\tilde{N}\big(\gamma(a)\big)-\tilde{N}\big(\gamma(x^\sigma)\big)\big)-\gamma(x)^{\tilde{\sigma}}\cdot \tilde{N}\big(\gamma(a)\big) \\
&=\gamma(a)\cdot \big(\gamma(a)^{\tilde{\sigma}}\gamma(x^\sigma)+\gamma(x^\sigma)^{\tilde{\sigma}}\gamma(a)\big)-\gamma(x)^{\tilde{\sigma}}\cdot \tilde{N}\big(\gamma(a)\big) \\
&=\tilde{N}\big(\gamma(a)\big)\cdot \gamma(x)^{\tilde{\sigma}}+\gamma(a)(\gamma(x)^{\tilde{\sigma}})^{\tilde{\sigma}}\gamma(a)-\gamma(x)^{\tilde{\sigma}}\cdot \tilde{N}\big(\gamma(a)\big) \\
&=\gamma(a)\gamma(x)\gamma(a)=\tilde{h}_{\gamma(a)}\big(\gamma(x)\big)
\end{align*}
for all $a\in L_0^*$, $x\in L_0$.
\end{itemize}\qed
\end{bewzwei}

\begin{kor}
Let $(\AA,\KK,N^\AA_\KK)$ be a quadratic space of type (m). Then we have
$$\MM(\AA,\KK,N^\AA_\KK)\cong \MM(\AA)\ .$$
\end{kor}

\begin{bew}
Let $\gamma:=\id_{\AA}$. Then the map $\phi=\id_\KK: \KK\to \KK=\gamma(\langle 1_\AA\rangle_\KK)$ is an isomorphism of fields, and the map
$$(\gamma,\phi)=(\id_{\AA},\id_{\KK}): (\AA,\KK,N^{\AA}_{\KK})\to (\AA,\KK,N^{\AA}_{\KK})$$
is an isomorphism of quadratic spaces so that we may apply theorem \ref{350}.\qed
\end{bew}

\begin{lemma}\label{456}
More precisely: Let $(\AA,\KK,N^\AA_\KK)$ be a quadratic space of type (m). Then we have
$$\MM:=\MM\big(\AA,\KK,N^\AA_\KK\big)=\MM(\AA)=:\tilde{\MM}\ .$$
\end{lemma}

\begin{bew}
We have
$$\forall\ x\in \AA:\qquad \tau(x)=-x^\sigma\cdot N^{\AA}_{\KK}(x)^{-1}=-x^{-1}=\tilde{\tau}(x)\ .$$
\qed
\end{bew}

\section[\texorpdfstring{$\MM(\KK,\KK_0,\LL_0)\cong \MM(\tilde{\AA})$}{Indifferent Type and Linear Type}]{\texorpdfstring{$\boldsymbol{\MM(\KK,\KK_0,\LL_0)\cong \MM(\tilde{\AA})}$}{Indifferent Type and Linear Type}}

\begin{satz}\label{388} Let $\MM:=\MM(\KK,\KK_0,\LL_0)$ be a Moufang set of indifferent type and let $\tilde{\MM}:=\MM(\tilde{\AA})$ be a Moufang set of linear type. A map $\gamma:\tilde{\MM}\to {\MM}$ is a Jordan isomorphism such that $\gamma(1_{\tilde{\MM}})=1_{{\MM}}$ if and only if $\KK_0=\KK$, $\tilde{\AA}$ is a field and the map $\gamma:\tilde{\AA}\to \KK_0$ is an isomorphism of fields. In particular, the indifferent set $(\KK,\KK_0,\LL_0)$ is non-proper if we have $\MM\cong \tilde{\MM}$.
\end{satz}

\begin{bewzwei}\ 
\begin{itemize}
\item[``$\R$''] The map $\gamma:\tilde{\AA}\to \KK_0\subseteq \KK$ is a Jordan homomorphism. Since $\KK$ is associative, Hua's theorem implies that $\KK_0=\gamma(\tilde{\AA})$ is a skew-field, which is thus, in fact, a subfield of $\KK$. By Hua's theorem again, the map $\gamma$ is an isomorphism of fields. In particular, $\tilde{\AA}=\gamma^{-1}(\KK_0)$ is a field. Moreover, we have 
$$\KK_0=\langle \KK_0\rangle=\KK\ ,$$
hence $(\KK,\KK_0,\LL_0)=(\KK,\KK,\LL_0)$ is non-proper.
\item[``$\Leftarrow$''] We have $\gamma(1_{\tilde{\MM}})=\gamma(1_{\tilde{\AA}})=1_{\tilde{\KK}_0}=1_\KK=1_{\MM}$. Given $a\in \tilde{\AA}^*$ and $x\in \tilde{\AA}$, we have
$$\gamma\big(\tilde{h}_a(x)\big)=\gamma(axa)=\gamma(a)\gamma(x)\gamma(a)=h_{\gamma(a)}\big(\gamma(x)\big)\ .$$
\end{itemize}\qed
\end{bewzwei}

\newpage

\section[\texorpdfstring{$\MM(\KK,\KK_0,\sigma)\cong \MM(\tilde{\AA})$}{Involutory Type and Linear Type}]{\texorpdfstring{$\boldsymbol{\MM(\KK,\KK_0,\sigma)\cong \MM(\tilde{\AA})}$}{Involutory Type and Linear Type}}

\begin{lemma}\label{506} Let $(\KK,\KK_0, \sigma)$ be a non-proper involutory set with the additional assumption that $\KK_0$ is a field if $\sigma=\id_\KK$ and $\Char \KK=2$. Then $(\KK,\KK_0,\sigma)$ is of quadratic type.
\end{lemma}

\begin{bew} Assume that $(\KK,\KK_0,\sigma)$ is not quadratic of type (v), i.e., $\KK$ is a skew-field.
\begin{itemize}
\item $\sigma=\id_\KK$, $\Char \KK\neq 2$: Then $\KK$ is a field and $\KK_\sigma=\KK_0=\Fix(\sigma)=\KK$, hence $(\KK,\KK_0,\sigma)$ is quadratic of type (ii).
\item $\sigma=\id_\KK$, $\Char \KK=2$: Then $\KK$ is a field and $\KK^2\subseteq \KK_0\subseteq \KK$, hence $(\KK,\KK_0,\sigma)$ is quadratic of type (m)$\in\{\textrm{(i), (ii)}\}$.
\item $\sigma\neq \id_\KK$: By (23.23) of \cite{TW}, $(\KK,\KK_0,\sigma)$ is quadratic of type (m)$\in\{\textrm{(iii), (iv)}\}$.
\end{itemize}\qed
\end{bew}

\begin{satz}\label{389} Let $\tilde{\MM}:=\MM(\tilde{\AA})$ be a Moufang set of linear type and let $\MM:=\MM(\KK,\KK_0,\sigma)$ be a Moufang set of involutory type. A map $\gamma:\tilde{\MM}\to {\MM}$ is a Jordan isomorphism such that $\gamma(1_{\tilde{\MM}})=1_{{\MM}}$ if and only if $(\KK,\KK_0,\sigma)$ is of quadratic type, $\tilde{\AA}$ and $\KK_0$ are fields and the map $\gamma:\tilde{\AA}\to \KK_0$ is an isomorphism of fields. In particular, the involutory set $(\KK,\KK_0,\sigma)$ is non-proper if we have $\MM\cong \tilde{\MM}$.
\end{satz}

\begin{bewzwei}\ 
\begin{itemize}
\item[``$\R$''] If $(\KK,\KK_0,\sigma)$ is quadratic of type (v), $\KK_0$ is a field and the map $\gamma:\MM(\tilde{\AA})\to \MM(\KK_0)$ is a Jordan isomorphism, hence an isomorphism of fields by theorem \ref{433} since $\KK_0$ is associative and commutative. In particual, $\tilde{\AA}$ is a field. In the sequel we suppose $\KK$ to be associative.

The map $\gamma:\tilde{\AA}\to \KK_0\subseteq \KK$ is a Jordan homomorphism. Since $\KK$ is associative, Hua's theorem implies that $\KK_0=\gamma(\tilde{\AA})$ is a skew-subfield, which is thus, in fact, a subfield of $\KK$ since we have $\KK_0\subseteq \Fix(\sigma)$ and, therefore,
$$\forall\ x,y\in \KK_0:\qquad xy=(xy)^\sigma=y^\sigma x^\sigma=yx\ .$$ By Hua's theorem again, $\gamma$ is an isomorphism of fields. In particular, $\tilde{\AA}=\gamma^{-1}(\KK_0)$ is a field. Moreover, $(\KK,\KK_0,\sigma)$ is non-proper by lemma \ref{157} and thus quadratic of type (m)$\in\{\textrm{(i),\ldots,(iv)}\}$ by lemma \ref{506}.

\item[``$\Leftarrow$''] We have $$\gamma(1_{\tilde{\MM}})=\gamma(1_{\tilde{\AA}})=1_{\tilde{\KK}_0}=1_\KK=1_{\MM}\ .$$
Given $a\in \tilde{\AA}^*$ and $x\in \tilde{\AA}$, we have
$$\gamma\big(\tilde{h}_a(x)\big)=\gamma(axa)=\gamma(a)\gamma(xa)=\gamma(a)\gamma(x)\gamma(a)=h_{\gamma(a)}\big(\gamma(x)\big)\ .$$
\end{itemize}\qed
\end{bewzwei}

\section[\texorpdfstring{$\MM(L_0,\KK,q)\cong \MM(\tilde{L}_0,\tilde{\KK},\tilde{q})$}{Quadratic Form Type}]{\texorpdfstring{$\boldsymbol{\MM(L_0,\KK,q)\cong \MM(\tilde{L}_0,\tilde{\KK},\tilde{q})}$}{Quadratic Form Type}} \label{409}

\begin{bem}
Cf. chapter 4.6 in \cite{K} for another proof of the main result of this paragraph.
\end{bem}

\begin{no}
Until theorem \ref{434}, $\MM:=\MM(L_0,\KK,q)$ is a Moufang set of quadratic form type with basepoint $\epsilon$, $\tilde{\MM}:=\MM(\tilde{L}_0,\tilde{\KK},\tilde{q})$ is a Moufang set of proper quadratic form type with basepoint $\tilde{\epsilon}$ and $\dim_{\tilde{\KK}} \tilde{L}_0\geq 3$, and $\gamma:\MM\to\tilde{\MM}$ is a Jordan isomorphism such that $\gamma(1_\MM)=1_{\tilde{\MM}}$.
\end{no}

\begin{lemma}\label{505}
If we have 
$$\forall\ x\in L_0^*:\qquad \gamma( \langle x\rangle_\KK)\subseteq  \langle\gamma(x)\rangle_{\tilde{\KK}}\ ,$$
we have
$$\forall\ x\in L_0^*:\qquad \gamma(\langle x\rangle_\KK)=\langle \gamma(x) \rangle_{\tilde{\KK}}\ .$$
\end{lemma}

\begin{bew} Notice that we have $\dim_{\tilde{\KK}} \tilde{L}_0\geq 3$. First of all, assume $\dim_\KK L\leq 2$. Then $\gamma(L_0)$ is contained in a two-dimensial $\tilde{\KK}$-subspace of $\tilde{L}_0$ and hence in a proper subspace of $\tilde{L}_0$, which contradicts the fact that $\gamma:L_0\to \tilde{L}_0$ is a bijection. We obtain $\dim_\KK \geq 3$, ans thus, by symmetry,
$$\forall\ x\in L_0^*:\qquad \gamma^{-1}(\langle \gamma(x) \rangle_{\tilde{\KK}})\subseteq \langle \gamma^{-1}\big(\gamma(x)\big) \rangle_\KK=\langle x\rangle_\KK\ ,$$
hence
\begin{align*}
\forall\ x\in L_0^*: && \gamma( \langle x\rangle_\KK)\supseteq  \langle\gamma(x)\rangle_{\tilde{\KK}}\ .
\end{align*}\qed
\end{bew}

\begin{no}
Given $x\in L_0^*$, $y\in L_0$ and $s\in \KK$, we set
\begin{align*}
\psi_1(x,y,s)&:=\gamma(x\cdot s)\cdot f_{\tilde{q}}\big(\gamma(x\cdot s),\gamma(y)^{\tilde{\sigma}}\big)-\gamma(y)^{\tilde{\sigma}}\cdot \tilde{q}\big(\gamma(x\cdot s )\big)\ , \\
\psi_2(x,y,s)&:=\gamma(x)\cdot f_{\tilde{q}}\big(\gamma(x),\gamma(y\cdot s^2)^{\tilde{\sigma}}\big)-\gamma(y\cdot s^2)^{\tilde{\sigma}}\cdot \tilde{q}\big(\gamma(x)\big)\ , \\
\psi_3(x,y,s)&:=\gamma(x)\cdot f_{\tilde{q}}\big(\gamma(x\cdot s),\gamma(y)^{\tilde{\sigma}}\big)+\gamma(x\cdot s)\cdot f_{\tilde{q}}\big(\gamma(x),\gamma(y)^{\tilde{\sigma}}\big)-\gamma(y)^{\tilde{\sigma}}\cdot f_{\tilde{q}}\big(\gamma(x\cdot s),\gamma(x)\big)\ , \\
\psi_4(x,y,s)&:=\gamma(x)\cdot f_{\tilde{q}}\big(\gamma(x),2\gamma(y\cdot s)^{\tilde{\sigma}}\big)-2\gamma(y\cdot s)^{\tilde{\sigma}}\cdot \tilde{q}\big(\gamma(x)\big)\ .
\end{align*}
\end{no}

\begin{lemma}
Let $x\in L_0^*$, $y\in L_0$ and $s\in \KK$. Then  the following holds:
\begin{enumerate}[label=(\alph*)]
\item We have
$$\psi_1(x,y,s+1_\KK)=\psi_1(x,y,s)+\psi_1(x,y,1_\KK)+\psi_3(x,y,s)\ .$$
\item We have 
$$\psi_2(x,y,s+1_\KK)=\psi_2(x,y,s)+\psi_2(x,y,1_\KK)+\psi_4(x,y,s)\ .$$
\item We have $\psi_1(x,y,s)=\psi_2(x,y,s)$.
\item \label{423} We have $\psi_3(x,y,s)=\psi_4(x,y,s)$.
\end{enumerate}
\end{lemma}

\begin{bewzwei}\ 
\begin{enumerate}[label=(\alph*)]
\item This is a direct calculation using the facts that $\gamma$ is additive and that we have
$$ \forall\ x\in L_0,\ s\in \KK:\qquad \tilde{q}\big(\gamma(x\cdot s)+\gamma(x)\big)=\tilde{q}\big(\gamma(x\cdot s)\big)+\tilde{q}\big(\gamma(x)\big)+f_{\tilde{q}}\big(\gamma(x\cdot s),\gamma(x)\big)\ .$$
\item This a direct calculation using the fact that $\gamma $ is additive.
\item By lemma \ref{329} and lemma \ref{328}, we have
\begin{align*}
\psi_1(x,y,s)&=\gamma(x\cdot s)\cdot f_{\tilde{q}}\big(\gamma(x\cdot s),\gamma(y)^{\tilde{\sigma}}\big)-\gamma(y)^{\tilde{\sigma}}\cdot \tilde{q}\big(\gamma(x\cdot s )\big) \\
& =\tilde{h}_{\gamma(x\cdot s)}\big(\gamma(y)\big)=\gamma\big(h_{x\cdot s}(y)\big)=\gamma\big(h_x(y\cdot s^2)\big)=\tilde{h}_{\gamma(x)}\big(\gamma(y\cdot s^2)\big) \\
&=\gamma(x)\cdot f_{\tilde{q}}\big(\gamma(x),\gamma(y\cdot s^2)^{\tilde{\sigma}}\big)-\gamma(y\cdot s^2)^{\tilde{\sigma}}\cdot \tilde{q}\big(\gamma(x)\big)=\psi_2(x,y,s)\ .
\end{align*}
\item By (a), (b) and (c), we have 
\begin{align*}
\gamma_1(x,y,s)+\gamma_1(x,y,1_\KK)+\gamma_3(x,y,s)&=\gamma_1(x,y,s+1_\KK)=\gamma_2(x,y,s+1_\KK)\\
&=\gamma_2(x,y,s)+\gamma_2(x,y,1_\KK)+\gamma_4(x,y,s)\ ,\end{align*}
hence $\gamma_3(x,y,s)=\gamma_4(x,y,s)$ by (c) again.
\end{enumerate}\qed
\end{bewzwei}

\newpage
\begin{lemma} \label{425}
Assume $\Char \KK=2$. Given $x\in L_0$ such that $\gamma(x)\notin \tilde{L}_0^\bot$, we have $$\forall\ s\in \KK:\qquad \gamma(x\cdot s)\in \langle \gamma(x)\rangle_{\tilde{\KK}}\ .$$ 
\end{lemma}

\begin{bew}
Lemma \ref{423} simplifies to
\begin{align} \gamma(x)\cdot f_{\tilde{q}}\big(\gamma(x\cdot s),\gamma(y)^{\tilde{\sigma}}\big)+\gamma(x\cdot s)\cdot f_{\tilde{q}}\big(\gamma(x),\gamma(y)^{\tilde{\sigma}}\big)+\gamma(y)^{\tilde{\sigma}}\cdot f_{\tilde{q}}\big(\gamma(x\cdot s),\gamma(x)\big)=0_{\tilde{L}_0}  \label{424}\end{align}
for all $x\in L_0^*$, $y\in L_0$, $s\in \KK$. Since we have $\dim_{\tilde{\KK}} \tilde{L}_0\geq 3$, there is an element $y\in L_0$ such that $\gamma(y)^{\tilde{\sigma}}\in \gamma(x)^\bot\cap \gamma(x\cdot s)^\bot$. Therefore, we have $$\forall\ s\in \KK:\qquad f_{\tilde{q}}\big(\gamma(x\cdot s),\gamma(x)\big)=0_{\tilde{\KK}}\ ,$$
and equation \eqref{424} simplifies to
$$\gamma(x)\cdot f_{\tilde{q}}\big(\gamma(x\cdot s),\gamma(y)^{\tilde{\sigma}}\big)=\gamma(x\cdot s)\cdot f_{\tilde{q}}\big(\gamma(x),\gamma(y)^{\tilde{\sigma}}\big)$$
for all $x\in L_0^\bot$, $y\in L_0$, $s\in \KK$.
By assumption, there is an element $y\in L_0$ such that $f_{\tilde{q}}\big(\gamma(x),\gamma(y)^{\tilde{\sigma}}\big)\neq 0_{\tilde{\KK}}$, hence
$$\gamma(x\cdot s)=\gamma(x)\cdot f_{\tilde{q}}\big(\gamma(x\cdot s),\gamma(y)^{\tilde{\sigma}}\big)f_{\tilde{q}}\big(\gamma(x),\gamma(y)^{\tilde{\sigma}}\big)^{-1}\in \langle \gamma(x)\rangle_{\tilde{\KK}}\ .$$\qed
\end{bew}

\begin{no}\ 
\begin{itemize}
\item Until proposition \ref{429}, we assume $\Char \KK=2$.
\item Given $x\in L_0$ such that $\gamma(x)\notin L_0^\bot$, let $\phi_x:\KK\to \tilde{\KK}$ defined by
$$\forall\ s\in \KK:\qquad \gamma(x\cdot s)=\gamma(x)\cdot \phi_x(s)\ .$$
\end{itemize}
\end{no}

\begin{lemma}\label{426}
Let $x\in L_0$ be such that $\gamma(x)\notin \tilde{L}_0^\bot$ and let $y\in L_0$ be such that $\gamma(y)\in \tilde{L}_0^\bot$. Then we have 
\begin{align*}
\gamma(x+y)\notin \tilde{L}_0^\bot\ .
\end{align*}
\end{lemma}

\begin{bew}
Given $z\in L_0$, we have
$$f_{\tilde{q}}\big(\gamma(x+y),\gamma(z)\big)=f_{\tilde{q}}\big(\gamma(x),\gamma(z)\big)\ ,$$
hence $\gamma(x+y)\notin L_0^\bot$.\qed
\end{bew}

\begin{kor}\label{427}
Let $x\in L_0$ be such that $\gamma(x)\notin \tilde{L}_0^\bot$ and let $y\in L_0$ be such that $\gamma(y)\in \tilde{L}_0^\bot$. Then we have 
\begin{align*}
\forall\ t\in \KK^*:\qquad \gamma\big( x\cdot t +y\big)\notin \tilde{L}_0^\bot\ .
\end{align*}
\end{kor}

\begin{bew}
By lemma \ref{425}, we have 
$$\forall\ t\in \KK^*:\qquad \gamma(x\cdot t)=\gamma(x)\cdot \phi_x(t)\notin \tilde{L}_0^\bot$$
so that we may apply lemma \ref{426}.\qed
\end{bew}

\begin{lemma}\label{428}
Let $x,y\in L_0$ be such that $\gamma(y)\notin \langle \gamma(x)\rangle_{\tilde{\KK}}$ and $\gamma(x),\gamma(y)\notin \tilde{L}_0^\bot$. Then we have $$\phi_x=\phi_y\ .$$
\end{lemma}

\begin{bew}
The assertion is clearly true for $\KK=\FF_2$ so that we may assume $|\KK|\geq 4$.
\begin{enumerate}[label=(\roman*)]
\item $\gamma(y)\notin \langle \gamma(x)\rangle_{\tilde{\KK}}$, $\gamma(x+y)\notin L_0^\bot$: Given $s\in \KK$, we have
\begin{align*}
\gamma(x)\cdot \phi_{x+y}(s)+\gamma(y)\cdot \phi_{x+y}(s)&=\gamma\big((x+y)\cdot s\big)=\gamma(x\cdot s)+\gamma(y\cdot s)\\
&=\gamma(x)\cdot \phi_x(s)+\gamma(y)\cdot \phi_y(s)\ ,\end{align*}
hence
$$\phi_x(s)=\phi_{x+y}(s)=\phi_y(s)\ .$$
The condition in this step is always fulfilled if we have $\tilde{L}_0^\bot=\{0\}$ so that we may assume $\tilde{L}_0^\bot\neq \{0\}$ in the following.
\item $y\in \langle x\rangle_\KK$, which means that we allow $\gamma(y)\in \langle \gamma(x)_{\rangle{\tilde{\KK}}}$ in this case, cf. lemma \ref{425}: Since we have $|\KK|\geq 4$, there is an element $t\in \KK^*$ such that
$$x\cdot t\notin \{ x,y \}\ .$$
Let $z\in L_0^*$ be such that $\gamma(z)\in \tilde{L}_0^\bot$, which implies $\gamma(z)\notin \langle \gamma(x)\rangle_{\tilde{\KK}}$. Then we have
\begin{align*}
\gamma(x\cdot t+z)\notin \tilde{L}_0^\bot\ , && \gamma(x+x\cdot t+z)\notin \tilde{L}_0^\bot\ , &&  \gamma(y+x\cdot t+z)\notin \tilde{L}_0^\bot
\end{align*}
by corollary \ref{427} and thus $\phi_x=\phi_{x\cdot t+z}=\phi_y$ by (i).
\item $\gamma(y)\notin \langle\gamma(x)\rangle_{\tilde{\KK}}$, $\gamma(x+y)\in L_0^\bot$: Let $t\in \KK\sm\{0_\KK,1_\KK\}$. By lemma \ref{425} and corollary \ref{427}, we have
\begin{align*}
\gamma(x\cdot t)\notin \tilde{L}_0^\bot\ , && \gamma(x\cdot t+y)=\gamma\big(x\cdot (t+1_\KK)+(x+y)\big)\notin \tilde{L}_0^\bot\ , 
\end{align*}
hence $\phi_x\overset{\mathclap{\textrm{(i)}}}{=}\phi_{x\cdot t}\overset{\mathclap{\textrm{(ii)}}}{=}\phi_y$ by the previous steps.
\end{enumerate}\qed
\end{bew}

\begin{prop}\label{429}
If we have $\Char \KK=2$, the map $\gamma:L_0\to \tilde{L}_0$ is an isomorphism of vector spaces.
\end{prop}

\begin{bew}
We show that we have
$$\forall\ x\in L_0^*,\ s\in \KK:\qquad \gamma(x\cdot s)\in \langle \gamma(x)\rangle_{\tilde{\KK}}$$
so that we may apply the fundamental theorem of projective geometry by lemma \ref{505}. 

Let $x\in L_0^*$. By lemma \ref{425}, we may assume $\gamma(x)\in L_0^\bot$. Since $(\tilde{L}_0,\tilde{\KK},\tilde{q})$ is proper, there is an element $y\in L_0$ such that $\gamma(y)\notin L_0^\bot$. Moreover, we have $\gamma(y)\notin \langle \gamma(x+y)\rangle_{\tilde{\KK}}$, and, by lemma \ref{426}, $\gamma(x+y)\notin \tilde{L}_0$. By lemma \ref{428} therefore, we have
$$\phi_y=\phi_{x+y}$$
and thus
\begin{align*}
\gamma(x\cdot s)+\gamma(y)\cdot \phi_y(s)&=\gamma(x\cdot s)+\gamma(y\cdot s)=\gamma\big((x+y)\cdot s\big)\\ &=\gamma(x+y)\cdot \phi_{x+y}(s)=\gamma(x)\cdot \phi_y(s)+\gamma(y)\cdot \phi_y(s)
\end{align*}
for each $s\in \KK$, hence
$$\forall\ s\in \KK:\qquad \gamma(x\cdot s)=\gamma(x)\cdot \phi_y(s)\in \langle \gamma(x)\rangle_{\tilde{\KK}}\ .$$\qed
\end{bew}

\begin{bem}
We drop the condition $\Char \KK=2$.
\end{bem}

\begin{lemma}\label{431}
Assume $\Char \KK\neq 2$. Let $y\in L_0^*$ and $s\in \KK$. Given $x\in L_0^*$ such that 
$$f_{\tilde{q}}\big(\gamma(x),\gamma(y)^{\tilde{\sigma}}\big)=0_{\tilde{\KK}}=f_{\tilde{q}}\big( \gamma(x), \gamma(y\cdot s)\big)\ ,$$ we have 
$$f_{\tilde{q}}\big(\gamma(x\cdot s),\gamma(y)^{\tilde{\sigma}}\big)=0_{\tilde{\KK}}\ .$$
\end{lemma}

\begin{bew}
By lemma \ref{423}, we have
$$\gamma(x)\cdot f_{\tilde{q}}\big(\gamma(x\cdot s),\gamma(y)^{\tilde{\sigma}}\big)=\gamma(y)^{\tilde{\sigma}}\cdot f_{\tilde{q}}\big(\gamma(x\cdot s),\gamma(x)\big)-2\gamma(y\cdot s)^{\tilde{\sigma}}\cdot \tilde{q}\big(\gamma(x)\big)\ .$$
Assume $f_{\tilde{q}}\big(\gamma(x\cdot s),\gamma(y)^{\tilde{\sigma}}\big)\neq 0_{\tilde{\KK}}$. Then we have
$$\gamma(x)\in \langle \gamma(y)^{\tilde{\sigma}}, \gamma(y\cdot s)^{\tilde{\sigma}}\rangle_{\tilde{\KK}}\subseteq \gamma(x)^\bot\qquad \lightning\ .$$
\qed
\end{bew}

\begin{prop}\label{430}
If we have $\Char \KK\neq 2$, the map $\gamma:L_0\to \tilde{L}_0$ is an isomorphism of vector spaces.
\end{prop}

\begin{bew}
We show that we have
$$\forall\ y\in L_0^*,\ s\in \KK:\qquad \gamma(y\cdot s)\in \langle \gamma(y)\rangle_{\tilde{\KK}}$$
so that we may apply the fundamental theorem of projective geometry by lemma \ref{505}.

Let $y\in L_0^*$. By assumption, there is an element $x\in L_0^*$ such that
$$f_{\tilde{q}}\big(\gamma(x),\gamma(y)^{\tilde{\sigma}}\big)=0_{\tilde{\KK}}=f_{\tilde{q}}\big( \gamma(x), \gamma(y\cdot s)\big)\ .$$
By lemma \ref{423} and lemma \ref{431}, we have
$$\gamma(y)^{\tilde{\sigma}}\cdot f_{\tilde{q}}\big(\gamma(x\cdot s),\gamma(x)\big)=2\gamma(y\cdot s)^{\tilde{\sigma}}\cdot \tilde{q}\big(\gamma(x)\big)\ ,$$
hence
$$\gamma(y\cdot s)=\frac{1}{2}\gamma(y)\cdot f_{\tilde{q}}\big(\gamma(x\cdot s),\gamma(x)\big)\tilde{q}\big(\gamma(x)\big)^{-1}\in \langle \gamma(y)\rangle_{\tilde{\KK}}\ .$$
\qed
\end{bew}

\begin{satz}\label{434}
Let $\Xi$ be a quadratic space with basepoint $\epsilon$, let $\tilde{\Xi}$ be a proper quadratic space with basepoint $\tilde{\epsilon}$ and let $\MM:=\MM(\Xi)$, $\tilde{\MM}:=\MM(\tilde{\Xi})$ be the corresponding Moufang sets of quadratic form type. A map $\gamma:\MM\to \tilde{\MM}$ is a Jordan isomorphism such that $\gamma(1_\MM)=1_{\tilde{\MM}}$ if and only if one of the following holds:
\begin{enumerate}[label=(\roman*)]
\item We have $\dim_{\tilde{\KK}} \tilde{L}_0\geq 3$ and there is an isomorphism $\phi:\KK\to \tilde{\KK}$ of fields such that the map
$$(\gamma,\phi):(L_0,\KK,q)\to (\tilde{L}_0,\tilde{\KK},\tilde{q})$$
is an isomorphism of quadratic spaces. In particular, we have $\dim_{{\KK}} {L}_0=\dim_{\tilde{\KK}}\tilde{L}_0$.
\item We have $\dim_{\tilde{\KK}} \tilde{L}_0\leq 2$, the map $$\phi:\KK\to \hat{\KK}:=\gamma(\langle \epsilon\rangle_\KK)\subseteq \hat{\FF}:=\FF(\tilde{L}_0,\tilde{\KK},\tilde{q}),\ s\mapsto \gamma(\epsilon\cdot s)$$
is an isomorphism of fields, the field $\hat{\FF}$ is quadratic over $\hat{\KK}$, and the map 
$$(\gamma,\phi):(L_0,\KK,q)\to (\hat{\FF},\hat{\KK},N^{\hat{\FF}}_{\hat{\KK}})$$
is an isomorphism of quadratic spaces. This is true even if $\tilde{\Xi}$ is non-proper.
\end{enumerate}
\end{satz}

\begin{bewzwei}\ 
\begin{enumerate}
\item[``$\Rightarrow$'']
\begin{itemize}
\item $\dim_{\tilde{\KK}} \tilde{L}_0\geq 3$: By proposition \ref{429} and proposition \ref{430}, there is an isomorphism $\phi:\KK\to \tilde{\KK}$ of fields such that the map
$(\gamma,\phi):(L_0,\KK)\to (\tilde{L}_0,\tilde{\KK})$
is an isomorphism of vector spaces. By lemma \ref{335}, we have
\begin{align*}
\tilde{h}_{\gamma(x)}(\tilde{\epsilon})-\gamma(x)\cdot \phi\big(T(x)\big)+\tilde{\epsilon}\cdot \phi\big(q(x)\big)=\gamma\big(h_x(\epsilon)-x\cdot T(x)+\epsilon\cdot q(x)\big)=0_{\tilde{L}_0}
\end{align*}
and
\begin{align*}
\tilde{h}_{\gamma(x)}(\tilde{\epsilon})-\gamma(x)\cdot \tilde{T}\big(\gamma(x)\big)+\tilde{\epsilon}\cdot \tilde{q}\big(\gamma(x)\big)=0_{\tilde{L}_0}
\end{align*}
for each $x\in L^*_0$, hence
\begin{align*}
\forall\ x\in L^*_0:\qquad -\gamma(x)\cdot \phi\big(T(x)\big)+\tilde{\epsilon}\cdot \phi\big(q(x)\big)=-\gamma(x)\cdot \tilde{T}\big(\gamma(y)\big)+\tilde{\epsilon}\cdot \tilde{q}\big(\gamma(x)\big)\ ,
\end{align*}
which implies
$$\forall\ x\in L_0\sm \langle \epsilon\rangle_\KK:\qquad \tilde{q}\big(\gamma(x)\big)=\phi\big(q(x)\big)\ .$$
Given $s\in \KK$, we have
\begin{align*}
\tilde{q}\big(\gamma(\epsilon\cdot s)\big)=\tilde{q}\big(\tilde{\epsilon}\cdot \phi(s)\big)=\phi(s)^2=\phi(s^2)=\phi\big(q(\epsilon\cdot s)\big)\ .
\end{align*}

\item $\dim_{\tilde{\KK}} \tilde{L}_0\leq 2$: By lemma \ref{472}, we have
$$\MM(L_0,\KK,q)\cong \MM(\tilde{L}_0,\tilde{\KK},\tilde{q})=\MM\big(\FF(\tilde{L}_0,\tilde{\KK},\tilde{q})\big)\ .$$
Now we apply theorem \ref{379}.
\end{itemize}
\item[``$\Leftarrow$''] Let $(\gamma,\phi)$ be an isomorphism of quadratic spaces \big(independent of the target space $(\tilde{L}_0,\tilde{\KK},\tilde{q})$ or $(\hat{\FF},\hat{\KK},N^{\hat{\FF}}_{\hat{\KK}})$\big). By lemma \ref{476}, we have
\begin{align*}
\gamma\big(h_x(y)\big)&=\gamma\big(x\cdot f_q(x,y^\sigma)-y^\sigma\cdot q(x)\big)=\gamma(x)\cdot \phi\big(f_q(x,y^\sigma)\big)-\gamma(y^\sigma)\cdot \phi\big(q(x)\big) \\
&=\gamma(x)\cdot f_{\tilde{q}}\big( \gamma(x),\gamma(y)^{\tilde{\sigma}}\big)-\gamma(y)^{\tilde{\sigma}}\cdot \tilde{q}\big(\gamma(y)\big)=\tilde{h}_{\gamma(x)}\big(\gamma(y)\big)
\end{align*}
for all $x\in L_0^*$, $y\in L_0$ which completes the case $\dim_{\tilde{\KK}}\tilde{L}_0\geq 3$. In the case $\dim_{\tilde{\KK}}\tilde{L}_0\leq 2$ moreover, we have
$$\MM(\tilde{L}_0,\tilde{\KK},\tilde{q})=\MM(\hat{\FF})=\MM(\hat{\KK},\hat{\FF},N^{\hat{\FF}}_{\hat{\KK}})$$
by lemma \ref{472} and lemma \ref{456} so that
$$\gamma:\MM(L_0,\KK,q)\to \MM(\hat{\KK},\hat{\FF},N^{\hat{\FF}}_{\hat{\KK}})=\MM(\tilde{L}_0,\tilde{\KK},\tilde{q})$$
is a Jordan isomorphism.
\end{enumerate}\qed
\end{bewzwei}

\newpage

\addtocontents{toc}{\noindent\protect\mbox{}\protect\hrulefill\par}
\part{An Inventory of Moufang Polygons}
\addtocontents{toc}{\noindent\protect\mbox{}\protect\hrulefill\par}

\noindent As in the simply laced case, the parametrization for the appearing Moufang polygons in root group sequences of a given twin building is an important technical tool for the classification of twin buildings: We can make use of the knowledge about the parametrizing Moufang sets and Jordan isomorphisms between them, i.e., we may apply the results of chapter \ref{416}.\\

\chapter{Parametrized Moufang Polygons}

Before we give the appearing examples, we have to establish the concept of standard and opposite parametrized Moufang polygons as we can read a root group sequence in two directions.

\begin{de}\ 
\begin{itemize}
\item The symbol $\TTT$ is \textit{standard of type $\mathit{3}$}.
\item The symbols $\QQQ_I,\QQQ_Q,\QQQ_D,\QQQ_P,\QQQ_E$ and $\QQQ_F$ are \textit{standard of type $\mathit{4}$}.
\item The symbol $\HHH$ is \textit{standard of type $\mathit{6}$}.
\item The symbol $\OOO$ is \textit{standard of type $\mathit{8}$}.
\end{itemize} \index{symbol!standard}
\end{de}

\begin{de}\ 
\begin{itemize}
\item A \textit{parameter system of type $\mathit{\TTT}$} is an alternative division ring $\AA$.
\item A \textit{parameter system of type $\mathit{\QQQ_I}$} is a proper involutory set $(\KK,\KK_0,\sigma)$.
\item A \textit{parameter system of type $\mathit{\QQQ_Q}$} is a quadratic space $(L_0,\KK,q)$ with basepoint $\epsilon$.
\item A \textit{parameter system of type $\mathit{\QQQ_D}$} is a proper indifferent set $(\KK,\KK_0,\LL_0)$.
\item A \textit{parameter system of type $\mathit{\QQQ_P}$} is a proper right pseudo-quadratic space $(\KK,\KK_0,\sigma,L_0,q)$.
\item A \textit{parameter system of type $\mathit{\QQQ_E}$} is a quadratic space $(L_0,\KK,q)$ of type $E_6,E_7,E_8$.
\item A \textit{parameter system of type $\mathit{\QQQ_F}$} is a quadratic space $(L_0,\KK,q)$ of type $F_4$.
\item A \textit{parameter system of type $\mathit{\HHH}$} is a hexagonal system $(J,\FF,\#)$.
\item A \textit{parameter system of type $\mathit{\OOO}$} is an octagonal set $(\KK,\sigma)$.\index{parameter system}
\end{itemize}
\end{de}

\begin{de} \newglossaryentry{MP}{type=symbols,name={\ensuremath{\XXX(\Xi)}},description={parametrized Moufang polygon},sort={Moufang polygon}}
A \textit{parametrized standard Moufang $n$-gon}\index{parametrized Moufang polygon!standard} is a standard root group sequence
$$\gls{MP}=\big(U_{[1,n]}, x_1(\MM_1),\ldots, x_n(\MM_n)\big)\ ,$$
where $\XXX$ is a standard symbol of type $n\in \{3,4,6,8\}$, $\Xi$ is a parameter system of type $\XXX$ and $\MM_1,\ldots,\MM_n$ are the parameter groups with respect to the parametrizations $x_1,\ldots,x_n$ and the corresponding commutator relations, cf. chapter 16 of \cite{TW}.
\end{de}

\begin{bem}\ 
\begin{enumerate}[label=(\alph*)]
\item Given a parametrized standard Moufang $n$-gon $\XXX(\Xi)$, the parameter groups $\MM_1,\ldots,\MM_n$ are Moufang sets.
\item For reasons of brevity, we will write $$\XXX(\Xi)=\big( x_1(\MM_1),\ldots, x_n(\MM_n)\big)$$ instead of 
$\XXX(\Xi)=\big( U_{[1,n]}, x_1(\MM_1),\ldots, x_n(\MM_n)\big)$.
\end{enumerate}
\end{bem}

\begin{de}
Let $\XXX$ be a standard symbol of type $n\in\{3,4,6,8\}$.
\begin{itemize}
\item The symbol $\XXX^o$ is the corresponding \textit{opposite symbol of type $\mathit{n}$}\index{symbol!opposite}\index{opposite symbol}. We set $(\XXX^o)^o:=\XXX$.
\item A \textit{parameter system of type $\mathit{\XXX^o}$} is just a parameter system of type $\XXX$, except for the following standard symbols:
\begin{enumerate}[label=$\circ$]
\item A \textit{parameter system of type $\QQQ_P^o$} is a proper left pseudo-quadratic space.
\item Hexagonal systems and octagonal sets are not taken into account yet.
\end{enumerate}
\end{itemize}
\end{de}

\begin{no}\label{363}\ 
\begin{itemize}
\item In the following, a \textit{symbol $\mathit{\XXX}$}\index{symbol} denotes either a standard or an opposite symbol of type $n$ for some $n\in \{3,4,6,8\}$.
\item In the following, a \textit{parameter system $\mathit{\Xi}$}\index{parameter system} denotes a parameter system of type $\XXX$ for some symbol $\XXX$.
\end{itemize}
\end{no}

\begin{bem}
The following list is not complete since we restrict to parameter systems for $n$-gons with $n\in \{3,4\}$. The list can be extended easily.
\end{bem}

\begin{de}
Given a parameter system $\Xi$, there is a natural way to define an \textit{opposite parameter system $\mathit{\Xi^o}$}\index{parameter system!opposite}\index{opposite!parameter system}:
\begin{itemize}
\item Given an alternative division ring $\AA$, the corresponding opposite parameter system $\AA^o$ is just the opposite alternative division ring $\AA^o$ as in definition \ref{355}.
\item Given an involutory set $(\KK,\KK_0,\sigma)$, the corresponding opposite parameter system is $$(\KK,\KK_0,\sigma)^o:=(\KK^o,\KK_0,\sigma)\ ,$$ which is an involutory set.
\item Given an indifferent set $(\KK,\KK_0,\LL_0)$, the corresponding opposite parameter system is $$(\KK,\KK_0,\LL_0)^o:=(\KK,\KK_0,\LL_0)\ ,$$
which is the indifferent set itself.
\item Given a right (resp. left) pseudo-quadratic space $(\KK,\KK_0,\sigma,L_0,q)$, the corresponding opposite parameter system is $$(\KK,\KK_0,\sigma,L_0,q)^o:=(\KK^o,\KK_0,\sigma,L_0,q)\ ,$$ which is a left (resp. right) pseudo-quadratic space.
\item Given a quadratic space $(L_0,\KK,q)$, the corresponding opposite parameter system is $$(L_0,\KK,q)^o:=(L_0,\KK,q)\ ,$$ which is the quadratic space itself.
\end{itemize}
\end{de}

\begin{bem}
Let $\Xi$ be a parameter system. If $\Xi$ is of type $\XXX$, then $\Xi^o$ is of type $\XXX^o$.
\end{bem}

\begin{de}
An \textit{parametrized opposite Moufang $n$-gon}\index{parametrized Moufang polygon!opposite} is a root group sequence
$$\XXX^o(\Xi)=\big(x_1(\MM_1),\ldots, x_n(\MM_n)\big)$$
such that 
$$\big( x_n(\MM_n^o),\ldots,x_1(\MM_1^o)\big)=\XXX(\Xi^o)\ ,$$
where $\XXX$ is a standard symbol of type $n\in \{3,4,6,8\}$ and $\Xi$ is a parameter system of type $\XXX^o$.
\end{de}

\begin{bem}\ 
\begin{enumerate}[label=(\alph*)]
\item Parametrized opposite Moufang $n$-gons can be obtained as follows: Take the opposite root group sequence $\big(x_n(\MM_n),\ldots,x_1(\MM_1)\big)$ of some parametrized standard Moufang $n$-gon $\XXX(\Xi)=\big(x_1(\MM_1),\ldots,x_n(\MM_n)\big)$, calculate the commutator relations and interpret them in the opposite parameter system $\Xi^o$.
\item By definition, each parametrized opposite Moufang $n$-gon arises in this way.
\item We will make this more explicit in the next chapter.
\end{enumerate}
\end{bem}

\begin{de}
A \textit{parametrized Moufang polygon of type $\mathit{\XXX}$}\index{parametrized Moufang polygon} is a parametrized standard or opposite Moufang $n$-gon $\XXX(\Xi)$ for some symbol $\XXX$ of type $n$ and some parameter system of type $\XXX$.
\end{de}

\begin{bem}
Two isomorphic Moufang polygons are necessarily $n$-gons for the same value $n\in \{3,4,6,8\}$, cf. p. 419 of \cite{TW}. However, there are six families of Moufang quadrangles, and there are indeed quadrangles belonging to different families, cf. chapter 38 of \cite{TW}. But the six families of parametrized Moufang quadrangles are disjoint if we use the above list of parameter systems, i.e., two isomorphic parametrized quadrangles are necessarily of the same type $\XXX$, cf. (38.9) of \cite{TW}.
\end{bem}

\begin{de} Let $\XXX(\Xi)$, ${\XXX}(\tilde{\Xi})$ be parametrized Moufang polygons of the same type $\XXX$.
\begin{itemize}
\item An \textit{isomorphism}\index{isomorphism!of parametrized Moufang polygons} $\alpha:\XXX(\Xi)\to {\XXX}(\tilde{\Xi})$ is an ordered set $(\alpha_1,\ldots,\alpha_n)$ such that $\alpha_i:\MM_i\to \tilde{\MM}_i$ is an isomorphism of groups for each $i\in \{1,\ldots,n\}$ and such that
\begin{align*} \big(x_1\alpha_1(\MM_1),\ldots,x_n\alpha_n(\MM_n)\big)={\XXX}(\Xi)\ .\end{align*}
\item A \textit{reparametrization for $\mathit{\XXX(\Xi)}$}\index{reparametrization} is an ordered set $\alpha=(\tilde{\Xi},\alpha_1,\ldots,\alpha_n)$ such that $\tilde{\Xi}$ is a parameter system of type $\XXX$ and $\alpha_i:\tilde{\MM}_i\to {\MM}_i$ is an isomorphism of groups for each $i\in \{1,\ldots,n\}$ and such that
\begin{align*} \big(x_1\alpha_1(\tilde{\MM}_1),\ldots,x_n\alpha_n(\tilde{\MM}_n)\big)={\XXX}(\tilde{\Xi})\ .\end{align*}
\end{itemize}
\end{de}

\begin{bem}
The following results enable us to define parametrizations for root group sequences such that we have $\gamma(1_\MM)=1_\MM$ for each appearing glueing $\gamma$.
\end{bem}

\begin{lemma}\label{421}
Given a parametrized Moufang $n$-gon $\XXX(\Xi)$, we have 
\begin{align*}
1_{\MM_1}\in Z(\MM_1)\ , && 1_{\MM_n}\in Z(\MM_n)\ .
\end{align*}
\end{lemma}

\begin{bew}
This results from the fact that $\MM_1$ and $\MM_n$ are commutative except for the symbols $\QQQ_P$, $\QQQ_E$, $\OOO$. In these cases, the assertion results from definition (10.15) and (38.10) of \cite{TW}, cf. Fig. 5 on page 354 of \cite{TW}.\qed
\end{bew}

\begin{lemma}\label{365}
Let $\XXX(\Xi)$ be a parametrized Moufang $n$-gon and let $a\in Z(\MM_1)$, $b\in Z(\MM_n)$. Then the following holds:
\begin{enumerate}[label=(\alph*)]
\item If $\XXX$ is of type $\TTT$, $\QQQ_I$, $\QQQ_D$, $\QQQ_P$, $\HHH$, $\OOO$, there is a reparametrization $\alpha=(\tilde{\Xi},\alpha_1,\ldots,\alpha_n)$ s.t.
\begin{align*}
x_1\big(\alpha_1(1_{\tilde{\MM}_1})\big)=x_1(a)\ , && x_n\big(\alpha_n(1_{\tilde{\MM}_n})\big)=x_n(b)\ .
\end{align*}
\item \label{422} If $\XXX$ is of type $\QQQ_Q$, there is a reparametrization $\alpha=(\tilde{\Xi},\alpha_1,\ldots,\alpha_{n-1},\id_{\MM_n})$ such that
\begin{align*}
 x_1\big(\alpha_1(1_{\tilde{\MM}_1})\big)=x_1(a)\ , && \tilde{q}\big(b\big)=1_{\tilde{\KK}}\ .
\end{align*}
In particular, we have $1_{\tilde{\MM}_n}=b$.
\end{enumerate}
\end{lemma}

\newpage

\begin{bewzwei}\ 
\begin{enumerate}[label=(\alph*)]
\item \begin{itemize}
\item[$\TTT$:] This results from lemma \ref{362}.
\item[$\QQQ_I$:] This results from (35.16) and (22.39) in \cite{TW}.
\item[$\QQQ_D$:] This results from (35.18) in \cite{TW}.
\item[$\QQQ_P$:] This results from (35.19) and (25.20) in \cite{TW}.
\item[$\HHH$:] This results from (29.40) and (29.42) in \cite{TW}.
\item[$\OOO$:] This results from (31.35) in \cite{TW}.
\end{itemize}
\item \begin{itemize}
\item[$\QQQ_Q$:] This results from (35.17) and (23.25) in \cite{TW}.
\end{itemize}\end{enumerate}\qed
\end{bewzwei}

\begin{lemma}\label{420}
Given a parametrized quadrangle $\QQQ_E(L_0,\KK,q)$ and $a\in Z(\MM_1)$, $b\in Z(\MM_n)$, there is a reparametrization $\alpha=(\tilde{\Xi},\alpha_1,\ldots,\alpha_{n-1},\id_{\MM_n})$ such that
\begin{align*}
 x_1\big(\alpha_1(1_{\tilde{\MM}_1})\big)=x_1(a)\ , && \tilde{q}\big(b\big)=1_{\tilde{\KK}}\ .
\end{align*}
In particular, we have $1_{\tilde{\MM}_n}=b$.
\end{lemma}

\begin{bew}
By remark (21.17) of \cite{TW}, $\QQQ_E(L_0,\KK,q)$ is an extension of $\QQQ_Q(L_0,\KK,q)$, and by proposition (21.4) and (38.10) of \cite{TW}, the first root group is $Y_1=Z(U_1)$. We apply lemma \ref{422} to the quadrangle $\QQQ_Q(L_0,\KK,q)$ and extend the reparametrization for $\QQQ_Q(L_0,\KK,q)$ to a parametrization for $\QQQ_E(L_0,\KK,q)$, which is possible by theorem (21.12) of \cite{TW}, more precisely, by its proof.\qed
\end{bew}

\begin{lemma}\label{419}
Let $\QQQ_F(L_0,\KK,q)$ be a parametrized quadrangle, let $(\FF,\hat{L}_0,\hat{q})$ be as in definition (14.12) of \cite{TW} and let $a\in \mathrm{Def}(\hat{q})^*$, $b\in \mathrm{Def(q)}^*$. Then there is a reparametrization $\alpha=(\tilde{\Xi},\alpha_1,\ldots,\alpha_{n-1},\id_{\MM_n})$ such that
\begin{align*}
 x_1\big(\alpha_1(1_{\tilde{M}_1})\big)=x_1(a)\ , && \tilde{q}\big(b\big)=1_{\tilde{\KK}}\ .
\end{align*}
In particular, we have $1_{\tilde{\MM}_n}=b$.
\end{lemma}

\begin{bew}
By remark (21.18) of \cite{TW}, $\QQQ_E(L_0,\KK,q)$ is an extension of $\QQQ_Q(L_0,\KK,q)$, and by proposition (21.4), remark (21.18) and the proof of (14.13) of \cite{TW}, the first root group is \begin{align} Y_1=\{ x_1(0,t) \mid t\in \KK\}=x_1\big(\mathrm{Def}(\hat{q})\big)\ .\label{415}\end{align}
We apply lemma \ref{422} to the quadrangle $\QQQ_Q(L_0,\KK,q)$ and extend the reparametrization for $\QQQ_Q(L_0,\KK,q)$ to a parametrization for $\QQQ_E(L_0,\KK,q)$ which is possible by theorem (21.12) of \cite{TW}, more precisely, by its proof.\qed
\end{bew}

\begin{no} \newglossaryentry{paraa}{type=symbols,name={\ensuremath{\MM_{(i,j)}^i}},description={parameter group of the first root group},sort={Moufang polygon}}\newglossaryentry{parab}{type=symbols,name={\ensuremath{\MM_{(i,j)}^j}},description={parameter group of the last root group},sort={Moufang polygon}}
Let $I$ be an index set, let $i\neq j\in I$ and let
$$\BBB_{(i,j)}=\XXX(\Xi)=( U_{[1,n]},x_1(\MM_1),\ldots, x_n(\MM_n)\big)$$
be a parametrized Moufang $n$-gon. Then $\gls{paraa}:=\MM_1$ denotes the parameter group of the first root group and $\gls{parab}:=\MM_n$ denotes the parameter group of the last root group with corresponding parametrizations $x_{(i,j)}^i:=x_1$ and $x_{(i,j)}^j:=x_n$. Moreover, we set
$$\BBB_{(i,j)}^o:=\XXX^o(\Xi^o)\ .$$
\end{no}

\begin{bem}\label{417}
Let $\BBB_{(1,2)}:=\QQQ_F(L_0,\KK,q)$ and $\BBB_{(2,3)}:=\QQQ_F(\tilde{L}_0,\tilde{\KK},\tilde{q})$ be quadrangles of type $F_4$, and let $\gamma:\MM_{(1,2)}^2\to \MM_{(2,3)}^2$ be a Jordan isomorphism. Then we have $\gamma(0,\FF)=(0,\tilde{\KK})$, cf. theorem \ref{400} and equation \eqref{415}.
\end{bem}

\newpage

\chapter{Parametrized Quadrangles}\label{414}

\section{Quadrangles of Involutory Type}

\newglossaryentry{inv}{type=foundations,name={\ensuremath{\QQQ_I(\KK,\KK_0,\sigma)}},description={parametrized standard Moufang quadrangle with respect to the involutory set ${(\KK,\KK_0,\sigma)}$}, sort=Moufang polygon}
\newglossaryentry{invo}{type=foundations,name={\ensuremath{\QQQ_I^o(\KK,\KK_0,\sigma)}},description={parametrized opposite Moufang quadrangle with respect to the involutory set ${(\KK,\KK_0,\sigma)}$}, sort=Moufang polygon}

\begin{de} Let $(\KK,\KK_0,\sigma)$ be a (proper) involutory set.
\begin{itemize}
\item The root group sequence 
$$\gls{inv}:=\big(x_1(\KK_0), x_2(\KK),x_3(\KK_0),x_4(\KK)\big)$$
with commutator relations
\begin{align*}
\forall\ s,t\in \KK: &&[x_2(s),x_4(t)^{-1}]&:=x_3(s^\sigma t+t^\sigma s)\ , \\
\forall\ s\in \KK,\ u\in \KK_0:&&[x_1(u),x_4(s)^{-1}]&:=x_2(us)x_3(s^\sigma u s)
\end{align*}
is the \textit{parametrized standard quadrangle of involutory type with respect to $\mathit{(\KK,\KK_0,\sigma)}$}\index{parametrized quadrangle!of involutory type}.
\item The root group sequence 
$$\gls{invo}:=\big(x_1(\KK), x_2(\KK_0),x_3(\KK),x_4(\KK_0)\big)$$
with commutator relations
\begin{align*}
\forall\ s,t\in \KK: &&[x_1(s)^{-1},x_3(t)]&:=x_2(-st^\sigma-ts^\sigma)\ , \\
\forall\ s\in \KK,\ u\in \KK_0:&&[x_1(s)^{-1},x_4(u)]&:=x_2(-sus^\sigma)x_3(-su)
\end{align*}
is the \textit{parametrized opposite quadrangle of involutory type with respect to $\mathit{(\KK,\KK_0,\sigma)}$}.
\end{itemize}
\end{de}

\begin{lemma}\label{161}
Let $(\KK,\KK_0,\sigma)$ be a (proper) involutory set and let $$\QQQ_I^o(\KK,\KK_0,\sigma)=\big(x_1(\KK), x_2(\KK_0), x_3(\KK), x_4(\KK_0)\big)$$
be the corresponding opposite quadrangle. Then the action of the \textit{Hua automorphism}\index{Hua automorphism}
$$h_1(s):=\mu\big(x_1(1_\KK)\big)^{-1}\mu\big(x_1(s)\big)$$
on $x_1(\KK)\times x_4(\KK_0)$ corresponds to the map
$$(t,u)\mapsto (sts, s^{-1}us^{-\sigma})\ ,$$
and the action of the \textit{Hua automorphism}
$$h_4(s):=\mu\big(x_4(1_\KK)\big)^{-1}\mu\big(x_4(s)\big)$$
on $x_1(\KK)\times x_4(\KK_0)$ corresponds to the map
$$(t,u)\mapsto (ts^{-1},  s^\sigma us )\ .$$
\end{lemma}

\begin{bew}
If we consider the quadrangle $\QQQ_I(\KK^o,\KK_0^o,\sigma)$, then the action of $h_1(s)$ on $x_4(\KK_0)\times x_1(\KK)$ corresponds to the map
$$(u,t)\mapsto (  s^{-\sigma}\circ u\circ s^{-1}  ,s\circ t \circ s)=(s^{-1}us^{-\sigma}, sts)\ ,$$
and the action of $h_4(s)$ on $x_4(\KK_0)\times x_1(\KK)$ corresponds to the map
$$(u,t)\mapsto ( s\circ u\circ s^\sigma, s^{-1}\circ t)=(s^\sigma u s, ts^{-1})\ ,$$
cf. (33.13) of \cite{TW}.\qed
\end{bew}

\section{Quadrangles of Pseudo-Quadratic Form Type}

\newglossaryentry{pqf}{type=foundations,name={\ensuremath{\QQQ_P(\KK,\KK_0,\sigma,L_0,q)}},description={parametrized standard Moufang quadrangle with respect to the pseudo-quadratic space ${(\KK,\KK_0,\sigma,L_0,q)}$}, sort=Moufang polygon}
\newglossaryentry{pqfo}{type=foundations,name={\ensuremath{\QQQ_P^o(\KK,\KK_0,\sigma,L_0,q)}},description={parametrized opposite Moufang quadrangle with respect to the pseudo-quadratic space ${(\KK,\KK_0,\sigma,L_0,q)}$}, sort=Moufang polygon}

\begin{de}\ 
\begin{itemize}
\item Let $(\KK,\KK_0,\sigma,L_0,q)$ be a (proper) right pseudo-quadratic space. Then the root group sequence
$$\gls{pqf}:=\big(x_1(T),x_2(\KK), x_3(T), x_4(\KK)\big)$$
with commutator relations
\begin{align*}
[x_1(a,t),x_4(b,u)^{-1}]&:=x_2\big(f(a,b)\big)\ , \\
[x_2(v),x_4(w)^{-1}]&:=x_3(0,v^\sigma w+w^\sigma v )\ , \\
[x_1(a,t),x_4(v)^{-1}]&:=x_2(tv)x_3(av,v^\sigma tv)
\end{align*}
for all $v,w\in \KK,\ (a,t),(b,u)\in T$ is the \textit{parametrized standard quadrangle of pseudo-quadratic form type with respect to $(\KK,\KK_0,\sigma,L_0,q)$}.\index{parametrized quadrangle!of pseudo-quadratic form type}
\item Let $(\KK,\KK_0,\sigma,L_0,q)$ be a (proper) left pseudo-quadratic space. Then the root group sequence
$$\gls{pqfo}:=\big(x_1(\KK),x_2(T), x_3(\KK), x_4(T)\big)$$
with commutator relations
\begin{align*}
[x_2(b,u)^{-1},x_4(a,t)]&:=x_3\big(-f(a,b)\big)\ , \\
[x_1(w)^{-1},x_3(v)]&:=x_2(0,-wv^\sigma -vw^\sigma )\ , \\
[x_1(v)^{-1},x_4(a,t)]&:=x_2(-va,-v^\sigma t^\sigma v)x_3(-vt)
\end{align*}
for all $v,w\in \KK,\ (a,t),(b,u)\in T$ is the \textit{parametrized opposite quadrangle of pseudo-quadratic form type with respect to $(\KK,\KK_0,\sigma,L_0,q)$}.
\end{itemize}
\end{de}

\begin{lemma}\label{322}
Let $(\KK,\KK_0,\sigma,L_0,q)$ be a (proper) left pseudo-quadratic space and let $$\QQQ_P^o(\KK,\KK_0,\sigma,L_0,q)=\big(x_1(\KK), x_2(T),x_3(\KK),x_4(T)\big)$$
be the corresponding opposite quadrangle. Then the action of the \textit{Hua automorphism}\index{Hua automorphism}
$$h_1(s):=\mu\big(x_1(1_\KK)\big)^{-1}\mu\big(x_1(s)\big)$$
on $x_1(\KK)\times x_4(T)$ corresponds to the map
$$\big(u,(b,v) \big)\mapsto \big(sus, (s^{-1}b, s^{-1}vs^{-\sigma})\big)\ ,$$
and the action of the \textit{Hua automorphism}
$$h_4(a,t):=\mu\big(x_4(0,1_\KK)\big)^{-1}\mu\big(x_4(a,t)\big)$$
on $x_1(\KK)\times x_4(T)$ corresponds to the map
$$\big(u,(b,v)\big)\mapsto \big( ut^{-\sigma},( t^\sigma b- t^\sigma f(a,b)t^{-1}a, t^\sigma v t)\big)\ .$$
\end{lemma}

\begin{bew}
If we consider the quadrangle $\QQQ_P(\KK^o,\KK_0^o,\sigma,L_0,q)$, then the action of $h_1(s)$ on $x_4(T)\times x_1(\KK)$ corresponds to the map
$$\big((b,v),u\big)\mapsto \big(  (b\circ s^{-1}, s^{-\sigma}\circ  v\circ s^{-1}) ,s\circ u \circ s\big)=\big((s^{-1} b,s^{-1}vs^{-\sigma}), sus\big)\ ,$$
and the action of $h_4(a,t)$ on $x_4(T)\times x_1(\KK)$ corresponds to the map
$$\big((b,v),u\big)\mapsto \big( (b\circ t^\sigma-a\circ t^{-1}\circ f(a,b)\circ t^\sigma, t\circ v\circ t^{\sigma}), t^{-\sigma}\circ u\big)=\big(( t^\sigma b-t^\sigma f(a,b)t^{-1}a, t^\sigma v t), ut^{-\sigma}\big)\ ,$$
cf. (33.13) of \cite{TW}.\qed
\end{bew}

\section{Quadrangles of Quadratic Form Type}

\newglossaryentry{qf}{type=foundations,name={\ensuremath{\QQQ_Q(L_0,\KK,q)}},description={parametrized standard Moufang quadrangle with respect to the quadratic space ${(L_0,\KK,q)}$}, sort=Moufang polygon}
\newglossaryentry{qfo}{type=foundations,name={\ensuremath{\QQQ_Q^o(L_0,\KK,q)}},description={parametrized opposite Moufang quadrangle with respect to the quadratic space ${(L_0,\KK,q)}$}, sort=Moufang polygon}

\begin{de} Let $(L_0,\KK,q)$ be a quadratic space with basepoint $\epsilon$.
\begin{itemize}
\item The root group sequence
$$\gls{qf}:=\big(x_1(\KK),x_2(L_0), x_3(\KK), x_4(L_0)\big)$$
with commutator relations
\begin{align*}
[x_2(a),x_4(b)^{-1}]&:=x_3\big(f(a,b)\big)\ , \\
[x_1(t),x_4(a)^{-1}]&:=x_2(at)x_3\big(tq(a)\big)
\end{align*}
for all $a,b\in L_0,\ t\in \KK$ is the \textit{parametrized standard quadrangle of quadratic form type with respect to $(L_0,\KK,q)$}.\index{parametrized quadrangle!of quadratic form type}
\item The root group sequence
$$\gls{qfo}:=\big(x_1(L_0),x_2(\KK), x_3(L_0), x_4(\KK)\big)$$
with commutator relations
\begin{align*}
[x_1(b)^{-1},x_1(a)]&:=x_2\big(-f(a,b)\big)\ , \\
[x_4(a)^{-1},x_1(t)]&:=x_2\big(-tq(a)\big)x_3(-at) 
\end{align*}
for all $a,b\in L_0,\ t\in \KK$ is the \textit{parametrized opposite quadrangle of quadratic form type with respect to $(L_0,\KK,q)$}.\end{itemize}
\end{de}

\begin{lemma}
Let $(L_0,\KK,q)$ be a quadratic space with basepoint $\epsilon$ and let $$\QQQ_Q^o(L_0,\KK,q)=\big(x_1(L_0), x_2(\KK),x_3(L_0), x_4(\KK)\big)$$
be the corresponding opposite quadrangle. Then the action of the \textit{Hua automorphism}\index{Hua automorphism}
$$h_1(a):=\mu\big(x_1(\epsilon)\big)^{-1}\mu\big(x_1(a)\big)$$
on $x_1(L_0)\times x_4(\KK)$ corresponds to the map
$$(v,u)\mapsto \big(\pi_a\pi_\epsilon(v)\cdot q(a),  u/q(a)\big)\ ,$$
and the action of the \textit{Hua automorphism}
$$h_4(s):=\mu\big(x_4(1_\KK)\big)^{-1}\mu\big(x_4(s)\big)$$
on $x_1(L_0)\times x_4(\KK)$ corresponds to the map
$$(b,u)\mapsto ( b\cdot t^{-1},t^2u)\ .$$
\end{lemma}

\begin{bew}
If we consider the quadrangle $\QQQ_P(L_0,\KK,q)$, then the action of $h_1(a)$ on $x_4(\KK)\times x_1(L_0)$ corresponds to the map
$$(u,v)\mapsto \big(  u/q(a), \pi_a\pi_\epsilon(v)\cdot q(a)\big)\ ,$$
and the action of $h_4(s)$ on $x_4(\KK)\times x_1(L_0)$ corresponds to the map
$$(u,b)\mapsto (t^2u,b\cdot t^{-1})\ ,$$
cf. (33.11) of \cite{TW}.\qed
\end{bew}

\section{Quadrangles of Indifferent Type}

\newglossaryentry{ind}{type=foundations,name={\ensuremath{\QQQ_D(\KK,\KK_0,\LL_0)}},description={parametrized standard Moufang quadrangle with respect to the indifferent set ${(\KK,\KK_0,\LL_0)}$}, sort=Moufang polygon}
\newglossaryentry{indo}{type=foundations,name={\ensuremath{\QQQ_D^o(\KK,\KK_0,\LL_0)}},description={parametrized opposite Moufang quadrangle with respect to the indifferent set ${(\KK,\KK_0,\LL_0)}$}, sort=Moufang polygon}

\begin{de} Let $(\KK,\KK_0,\LL_0)$ be a (proper) indifferent set.
\begin{itemize}
\item The root group sequence
$$\gls{ind}:=\big(x_1(\KK_0),x_2(\LL_0), x_3(\KK_0), x_4(\LL_0)\big)$$
with commutator relations
\begin{align*}
\forall\ t\in \KK_0,a\in \LL_0:\qquad [x_1(t),x_4(a)]=x_2(t^2a)x_3(ta)
\end{align*}
is the \textit{parametrized standard quadrangle of indifferent type with respect to $(\KK,\KK_0,\LL_0)$}.\index{parametrized quadrangle!of indifferent type}
\item The root group sequence
$$\gls{indo}:=\big(x_1(\LL_0),x_2(\KK_0), x_3(\LL_0), x_4(\KK_0)\big)$$
with commutator relations
\begin{align*}
\forall\ t\in \KK_0,a\in \LL_0:\qquad [x_1(a),x_4(t)]=x_2(-ta)x_3(-t^2a)=x_2(ta)x_3(t^2a)\end{align*}
is the \textit{parametrized opposite quadrangle of indifferent type with respect to $(\KK,\KK_0,\LL_0)$}.\end{itemize}
\end{de}

\begin{lemma}
Let $(\KK,\KK_0,\LL_0)$ be a (proper) indifferent set and let $$\QQQ_D^o(\KK,\KK_0,\LL_0)=\big(x_1(\LL_0), x_2(\KK_0),x_3(\LL_0), x_4(\KK_0)\big)$$
be the corresponding opposite quadrangle. Then the action of the \textit{Hua automorphism}\index{Hua automorphism}
$$h_1(a):=\mu\big(x_1(1_\KK)\big)^{-1}\mu\big(x_1(a)\big)$$
on $x_1(\LL_0)\times x_4(\KK_0)$ corresponds to the map
$$(b,u)\mapsto (ba^2,ua^{-1})\ ,$$
and the action of the \textit{Hua automorphism}
$$h_4(t):=\mu\big(x_4(1_\KK)\big)^{-1}\mu\big(x_4(t)\big)$$
on $x_1(\LL_0)\times x_4(\KK_0)$ corresponds to the map
$$(b,u)\mapsto ( bt^{-2},ut^2)\ .$$
\end{lemma}

\begin{bew}
If we consider the quadrangle $\QQQ_D(\KK,\KK_0,\LL_0)$, then the action of $h_1(a)$ on $x_4(\KK_0)\times x_1(\LL_0)$ corresponds to the map
$$(u,b)\mapsto (ua^{-1},ba^2)\ ,$$
and the action of $h_4(s)$ on $x_4(\KK_0)\times x_1(\LL_0)$ corresponds to the map
$$(u,b)\mapsto (ut^2,bt^{-2})\ ,$$
cf. (33.12) of \cite{TW}.\qed
\end{bew}

\section[Quadrangles of Type \texorpdfstring{${E_n,F_4}$}{En,F4}]{Quadrangles of Type \texorpdfstring{$\boldsymbol{E_n,F_4}$}{En,F4}}

For an overview of those quadrangles, we refer to (16.6), (33.14), (16.7) and (33.15) of \cite{TW}.

\newpage

\chapter{The Moufang Sets of Moufang Polygons}

We give an overview of the Moufang sets appearing as root groups of Moufang triangles and quadrangles.

\begin{bem}
Let $\BBB_{(1,2)}:=\XXX(\Xi)$ be a parametrized standard Moufang polygon. Then one of the following holds:
\begin{enumerate}[label=(\roman*)]
\item We have 
$$\MM_{(1,2)}^1=\MM(\AA)=\MM_{(1,2)}^2$$
if $\BBB_{(1,2)}=\TTT(\AA)$ for some alternative division ring $\AA$.
\item We have
\begin{align*}
\MM_{(1,2)}^1=\MM(\KK,\KK_0,\sigma)\ , && \MM_{(1,2)}^2=\MM(\KK)
\end{align*}
if $\BBB_{(1,2)}=\QQQ_I(\KK,\KK_0,\sigma)$ for some (proper) involutory set $(\KK,\KK_0,\sigma)$.
\item We have
\begin{align*}
\MM_{(1,2)}^1=\MM(\KK,\KK_0,\sigma,L_0,q)\ , && \MM_{(1,2)}^2=\MM(\KK)
\end{align*}
if $\BBB_{(1,2)}=\QQQ_P(\KK,\KK_0,\sigma,L_0,q)$ for some (proper) pseudo-quadratic space $(\KK,\KK_0,\sigma,L_0,q)$.
\item We have
\begin{align*}
\MM_{(1,2)}^1=\MM(\KK)\ , && \MM_{(1,2)}^2=\MM(L_0,\KK,q)
\end{align*}
if $\BBB_{(1,2)}=\QQQ_Q(L_0,\KK,q)$ for some quadratic space $(L_0,\KK,q)$ with basepoint $\epsilon$.
\item We have
\begin{align*}
\MM_{(1,2)}^1=\MM(\KK,\KK_0,\LL_0)\ , && \MM_{(1,2)}^2=\MM(\LL,\LL_0,\KK_0^2)
\end{align*}
if $\BBB_{(1,2)}=\QQQ_I(\KK,\KK_0,\LL_0)$ for some (proper) indifferent set $(\KK,\KK_0,\LL_0)$.
\item \label{473} We have
\begin{align*}
\MM_{(1,2)}^1=\MM(S)\ , && \MM_{(1,2)}^2=\MM(L_0,\KK,q)
\end{align*}
if $\BBB_{(1,2)}=\QQQ_E(L_0,\KK,q)$ for some quadratic space $(L_0,\KK,q)$ of type $E_n$.
\item We have
\begin{align*}
\MM_{(1,2)}^1=\MM(\FF,\hat{L}_0,\hat{q})\ , && \MM_{(1,2)}^2=\MM(L_0,\KK,q)
\end{align*}
if $\BBB_{(1,2)}=\QQQ_4(L_0,\KK,q)$ for some quadratic space $(L_0,\KK,q)$ of type $F_4$.
\item We have
\begin{align*}
\MM_{(1,2)}^1=\MM(J,\FF,\#)\ , && \MM_{(1,2)}^2=\MM(\FF)
\end{align*}
if $\BBB_{(1,2)}=\HHH(J,\FF,\#)$ for some hexagonal system $(J,\FF,\#)$.
\item We have
\begin{align*}
\MM_{(1,2)}^1=\MM(\KK)\ , && \MM_{(1,2)}^2=\MM(\KK,\sigma)
\end{align*}
if $\BBB_{(1,2)}=\OOO(\KK,\sigma)$ for some octagonal system $(\KK,\sigma)$.
\end{enumerate}
\end{bem}

\begin{bem}
The Hua automorphisms of lemma \ref{361} and chapter \ref{414}, which can be defined for each polygon, cf. chapter 33 of \cite{TW}, induce the Hua maps on the corresponding Moufang sets. Each Hua automorphism of a polygon which is part of an integrable foundation is induced by an automorphism of the whole building, cf. theorem \ref{452}.
\end{bem}

\newpage

\addtocontents{toc}{\noindent\protect\mbox{}\protect\hrulefill\par}
\part{Foundations}\label{492}
\addtocontents{toc}{\noindent\protect\mbox{}\protect\hrulefill\par}

\noindent We generalize the definitions and results of chapter \ref{491}, i.e., we give the definition of a foundation involving arbitrary Moufang polygons and show that we can attach a foundation to each twin building. Once again, this foundation turns out to be a classifying invariant of the corresponding twin building, and the integrability criterions of chapter \ref{491} hold as well.\\

\chapter{Definition}

\begin{de}\ 
\begin{itemize}
\item Let $M$ be a Coxeter matrix. A \textit{foundation of type $\mathit{M}$}\index{foundation} is a set
$$\FFF:=\{\BBB_{(i,j)} , \gamma_{(i,j,k)} \mid (i,j)\in A(M),(i,j,k)\in G(M) \}$$
such that:
\begin{enumerate}[label=(F\arabic*)]
\item Given $(i,j)\in A(M)$, then $\BBB_{(i,j)}=\XXX_{(i,j)}(\Xi_{(i,j)})$ for some symbol $\XXX_{(i,j)}$ of type $m_{ij}$ as in notation \ref{363} and some parameter system $\Xi_{(i,j)}$ of type $\XXX_{(i,j)}$.
\item Given $(i,j)\in A(M)$, we have $\BBB_{(i,j)}=\BBB_{(j,i)}^o$.
\item Given $(i,j,k)\in G(M)$, then $\gamma_{(i,j,k)}:\MM^j_{(i,j)}\to \MM^j_{(j,k)}$ is an isomorphism of groups satisfying
\begin{align*} \gamma_{(i,j,k)}(1_{\MM})=1_{\MM}\ , && \gamma_{(i,j,k)}=\id^o\circ \gamma_{(k,j,i)}^{-1}\circ \id^o\ .\end{align*}
\item Given $(i,j,k),(i,j,l),(l,j,k)\in G(M)$, we have
$$\gamma_{(i,j,k)}=\gamma_{(l,j,k)}\circ \id^o\circ \gamma_{(i,j,l)}\ .$$
\end{enumerate}
\item Given a foundation $\FFF$, we denote the corresponding Coxeter Matrix by $F$.
\item A foundation $\FFF$ is a \textit{Moufang foundation}\index{foundation!Moufang}\index{Moufang foundation} if each \textit{glueing} $\gamma:=\gamma_{(i,j,k)}$ is a Jordan isomorphism, i.e., we have
\begin{align*}
\forall\ a\in \MM_{(i,j)}^*,\ x\in\MM_{(i,j)}:\qquad\gamma\big(h_a(x)\big)=h_{\gamma(a)}\big(\gamma(x)\big)\ .
\end{align*}
\end{itemize}
\end{de}

\begin{de}
Let $\FFF$ be a foundation over $I=V(F)$ and let $J\subseteq I$. The \textit{$\mathit{J}$-residue of $\mathit{\FFF}$}\index{$J$-residue!of a foundation}\index{residue!of a foundation} is the foundation
$$\gls{Fres}:=\{ \BBB_{(i,j)},\gamma_{(i,j,k)} \mid (i,j)\in J^2\cap A(F), (i,j,k)\in J^3\cap G(F)\}\ .$$
\end{de}

\begin{bem}
Since a foundation is, in fact, an amalgam of Moufang polygons, an isomorphism of foundations is a system of isomorphism of Moufang polygons preserving the glueings.
\end{bem}

\begin{de}
Let $\FFF,\tilde{\FFF}$ be foundations.
\begin{itemize}
\item An \textit{isomorphism}\index{isomorphism!of foundations} $\alpha:\FFF\to \tilde{\FFF}$ is a
system $\alpha=\{\pi, \alpha_{(i,j)} \mid (i,j)\in A(F)\}$ of isomorphisms
\begin{align*}
\pi:F\to \tilde{F}\ , && \alpha_{(i,j)}=(\alpha_{(i,j)}^i,\ldots,\alpha_{(i,j)}^j): \BBB_{(i,j)}\to \tilde{\BBB}_{(\pi(i),\pi(j))}
\end{align*}
such that
\begin{align*}
\forall\ (i,j,k)\in G(F): \qquad \tilde{\gamma}_{(\pi(i),\pi(j),\pi(k))}\circ \alpha_{(i,j)}^j=\alpha_{(j,k)}^j\circ \gamma_{(i,j,k)}
\end{align*}
and $$\forall\ (i,j)\in A(F):\qquad \alpha_{(i,j)}=\alpha_{(j,i)}^o\ .$$
\item An isomorphism $\alpha:\FFF\to\tilde{\FFF}$ is \textit{special}\index{isomorphism!special} if $F=\tilde{F}$ and $\pi=\id_F$.
\item An \textit{automorphism of $\mathit{\FFF}$}\index{automorphism!of a foundation} is an isomorphism $\alpha:\FFF\to\FFF$.
\end{itemize}
\end{de}

\newpage

\chapter{Root Group Systems}

The fact that a root group systems is a classifying invariant of the corresponding twin building is a fundamental result in twin building theory.

\begin{de}
Let $\BBB$ be a twin building of type $M$, let $\Sigma$ be a twin apartment of $\BBB$ and let $c\in \OOO_\Sigma$.  
\begin{itemize}
\item Given $(i,j)\in A(M)$, let $\alpha_i,\alpha_j$ be the simple roots with respect to $(\Sigma,c)$ and let $\Theta_{(i,j)}$ be as in theorem \ref{499}. Then
$$\gls{rgs}:=( U_{[i,j]}, U_{(i,j)}^i, \ldots, U_{(i,j)}^j ):=\Theta_{(i,j)}$$
denotes the root group sequence of $\BBB$ from $\alpha_i$ to $\alpha_j$, which is isomorphic to the root group sequence of $\BBB_{ij}$ from $\alpha_i \cap \BBB_{ij}$ to $\alpha_j\cap \BBB_{ij}$.
\item The resulting set
$$\gls{rgsys}:=\{ U_{(i,j)} \mid (i,j)\in A(M)\}$$
is the \textit{root group system of $\mathit{\BBB}$ based at $\mathit{(\Sigma,c)}$}\index{root group system}.
\end{itemize}
\end{de}

\begin{lemma}\label{364}
Given $(i,j,k)\in G(F)$, we have $U_{(i,j)}^j=U_{(j,k)}^j$.
\end{lemma}

\begin{bew}
This holds by definition.\qed
\end{bew}

\begin{de}
Let $\UUU:=\UUU(\BBB,M,\Sigma,c)$ and $\tilde{\UUU}:=\UUU(\tilde{\BBB},\tilde{M},\tilde{\Sigma},\tilde{c})$ be root group systems.\begin{itemize}
\item An \textit{isomorphism}\index{isomorphism!of root group systems} $\alpha:\UUU\to\tilde{\UUU}$ is a system
$$\alpha=\{ \pi,\alpha_{(i,j)} \mid (i,j)\in A(M)\} $$
of isomorphisms
\begin{align*}
\pi:M\to \tilde{M}\ , && \alpha_{(i,j)}:U_{(i,j)}\to \tilde{U}_{(\pi(i),\pi(j))}
\end{align*}
of root group sequences such that
\begin{align*}
\forall\ (i,j,k)\in G(M):\ {\alpha_{(i,j)}}_{|U_{(i,j)}^j}={\alpha_{(j,k)}}_{|U_{(j,k)}^j}\ , && \forall\ (i,j)\in A(M):\ \alpha_{(i,j)}=\alpha_{(j,i)}^o\ .
\end{align*}
\item An isomorphism $\alpha:\UUU\to\tilde{\UUU}$ is \textit{special}\index{isomorphism!special} if $M=\tilde{M}$ and $\pi=\id_M$.
\item An \textit{automorphism of $\mathit{\UUU}$}\index{automorphism!of a root group system} is an isomorphism $\alpha:\UUU\to\UUU$.
\end{itemize}
\end{de}

\begin{satz}
Two root group systems $\UUU(\BBB,M,\Sigma,c)$ and $\UUU(\BBB,M,\tilde{\Sigma},\tilde{c})$ of a twin building $\BBB$ are specially isomorphic.
\end{satz}

\begin{bew} This is a consequence of theorem \ref{497}. \qed
\end{bew}

\begin{satz}\label{368}
Let $\UUU:=\UUU(\BBB,M,\Sigma,c)$ be a root group system of a twin building $\BBB$ which satisfies condition (CO). Then the isomorphism class of $\UUU$ is a classifying invariant of the isomorphism class of $\BBB$.
\end{satz}

\begin{bew}
This is a consequence of the extension theorem \ref{498}.\qed
\end{bew}

\newpage

\chapter{Foundations and Root Group Systems}

Given a root group system, there is a natural way to attach a foundation to it.

\begin{de}\label{366}
Let $\UUU(\BBB,M,\Sigma,c)$ be a root group system.
\begin{itemize}
\item Given $(i,j)\in A(M)$, there is a symbol $\XXX_{(i,j)}$ and a parameter system $\Xi_{(i,j)}$ such that $U_{(i,j)}\cong \XXX_{(i,j)}(\Xi_{(i,j)})$. In particular, there is a system of parametrizations 
\begin{align*}
x_{(i,j)}^*:\MM_{(i,j)}^*\to U_{(i,j)}^*\ ,\ t\mapsto x_{(i,j)}^*(t)\ , && *\in\{i,j\}
\end{align*}
extending to the defining relations for $\XXX_{(i,j)}(\Xi_{(i,j)})$. Such a parametrization yields an opposite system of parametrizations
\begin{align*}
x_{(j,i)}^*:(\MM_{(i,j)}^*)^{o}\to U_{(j,i)}^*\ ,\ t\mapsto x_{(i,j)}^*\big(\id^o(t)\big)\ ,  && *\in\{j,i\}\ .
\end{align*}
The resulting set
$$\Lambda:=\{ \XXX_{(i,j)}(\Xi_{(i,j)}) \mid (i,j)\in A(M)\}$$
is a \textit{parameter system for $\mathit{\UUU}$}\index{parameter system}.
\item Given $(i,j,k)\in G(M)$ and parametrizations $\XXX_{(i,j)}(\Xi_{(i,j)})$ and $\XXX_{(j,k)}(\Xi_{(j,k)})$, we define the glueing  $\gamma_{(i,j,k)}:\MM_{(i,j)}^j\to \MM_{(j,k)}^j$ by
$$x_{(i,j)}^j(t)=x_{(j,k)}^j\big(\gamma_{(i,j,k)}(t)\big)$$
which is justified by lemma \ref{364}. Then $\gamma_{(i,j,k)}$ is an isomorphism of groups satisfying $\gamma_{(i,j,k)}=\id^o\circ\gamma_{(k,j,i)}^{-1}\circ\id^o$. By lemma \ref{421}, \ref{365}, \ref{420} and \ref{419}, we may adjust all the parametrizations such that $$\forall\ (i,j,k)\in G(F):\qquad \gamma_{(i,j,k)}(1_{\MM})=1_{\MM}\ .$$
In the first instance, we have to adjust the glueings connecting quadrangles of type $F_4$ (for which we need remark \ref{417}) since we have the least flexibility in this case: The element $1_\MM$ is an element in the corresponding defect.
\end{itemize}
\end{de}

\begin{lemma} Given a root group system $\UUU:=\UUU(\BBB,M,\Sigma,c)$, a parameter system $\Lambda$ as in definition \ref{366} induces a foundation $$\gls{fl}=\{\XXX_{(i,j)}(\Xi_{(i,j)}), \gamma_{(i,j,k)} \mid (i,j)\in A(M), (i,j,k)\in G(M)\}\ .$$ 
\end{lemma}

\begin{bew}
We emphasize that the glueings in definition \ref{366} are identifications with respect to directed edges. Given $(i,j,k),(i,j,l),(l,j,k)\in G(M)$ and $t\in \MM_{(i,j)}^j$, we have
\begin{align*} 
x_{(j,k)}^j\big(\gamma_{(i,j,k)}(t)\big)&=x_{(i,j)}^j(t)=x_{(j,l)}^j\big(\gamma_{(i,j,l)}(t)\big)\\
&=x_{(l,j)}^j\big(\id^o\circ \gamma_{(i,j,l)}(t)\big)=x_{(j,k)}^j\big(\gamma_{(l,j,k)}\circ\id^o\circ \gamma_{(i,j,l)}(t)\big)\end{align*}
and thus $\gamma_{(i,j,k)}=\gamma_{(l,j,k)}\circ \id^o\circ \gamma_{(i,j,l)}$.\qed
\end{bew}

\begin{de}
A foundation $\FFF$ is \textit{integrable}\index{integrable}\index{foundation!integrable} if it is the foundation of a twin building $\BBB$, i.e., if there are a root group system $\UUU:=\UUU(\BBB,M,\Sigma,c)$ and a parameter system $\Lambda$ for $\UUU$ such that $\FFF=\FFF(\UUU,\Lambda)$.
\end{de}

\begin{satz}
Let $\FFF$ be an integrable foundation. Then $\FFF$ is a Moufang foundation.
\end{satz}

\begin{bew}
Let $(i,j,k)\in G(F)$ and $\gamma:=\gamma_{(i,j,k)}$. If we set
\begin{align*}
h(a)&:=\mu\big(x_{(i,j)}^j(1_{\MM})\big)^{-1}\mu\big(x_{(i,j)}^j(a)\big)\ , && a\in (M_{(i,j)}^j)^*\ , \\
\tilde{h}(a)&:=\mu\big(x_{(j,k)}^j(1_{\MM})\big)^{-1}\mu\big(x_{(j,k)}^j(a)\big)\ , && a\in (M_{(j,k)}^j)^*\ ,
\end{align*}
we have
$$\tilde{h}\big(\gamma(a)\big)=\mu\big(x_{(j,k)}^j(1_\MM)\big)^{-1}\mu\big(x_{(j,k)}^j(\gamma(a))\big)=\mu\big(x_{(i,j)}^j(1_\MM)\big)^{-1}\mu\big(x_{(i,j)}^j(a)\big)=h(a)$$
for each $a\in \MM_{(i,j)}^j$. Moreover, we have
\begin{align*}
x_{(i,j)}^j(x)^{h(a)}=x_{(i,j)}^j\big(h_a(x)\big)\ , && x_{(j,k)}^j(y)^{\tilde{h}(b)}=x_{(j,k)}^j\big(\tilde{h}_b(y)\big)
\end{align*}
for all $a\in (\MM_{(i,j)}^j)^*$, $x\in \MM_{(i,j)}^j$, $b\in (\MM_{(j,k)}^j)^*$, $y\in \MM_{(j,k)}^j$. Combining these two facts yields
\begin{align*}
x_{(j,k)}^j\big(\gamma(h_a(x))\big)&=x_{(i,j)}^j\big(h_a(x)\big)=x_{(i,j)}^j(x)^{h(a)}=x_{(j,k)}^j\big(\gamma(x)\big)^{\tilde{h}(\gamma(x))}
=x_{(j,k)}^j\big(\tilde{h}_{\gamma(a)}(\gamma(x))\big)
\end{align*}
and thus $\gamma\big(h_a(x)\big)=\tilde{h}_{\gamma(a)}\big(\gamma(x)\big)$ for all $a\in (\MM_{(i,j)}^j)^*$, $x\in \MM_{(i,j)}^j$.\qed
\end{bew}

\begin{bem}
The following result provides an integrability criterion.
\end{bem}

\begin{de} Let $\FFF$ be a foundation.
\begin{itemize}
\item Let $(\tilde{F},\p)$ be a cover of $F$. Then the foundation
$$\gls{fc}:=\{ \tilde{\BBB}_{(i,j)}, \tilde{\gamma}_{(i,j,k)} \mid (i,j)\in A(\tilde{F}), (i,j,k)\in G(\tilde{F})\}$$
with
\begin{align*}
\forall\ (i,j)\in A(\tilde{F}):\ \tilde{\BBB}_{(i,j)}=\BBB_{(\p(i),\p(j))}\ , && \forall\ (i,j,k)\in G(\tilde{F}):\ \tilde{\gamma}_{(i,j,k)}=\gamma_{(\p(i),\p(j),\p(k))}
\end{align*}
is the \textit{cover corresponding to $\mathit(\tilde{F},\p)$}\index{cover!of a foundation}.
\item A foundation $\tilde{\FFF}$ is a \textit{cover of $\mathit{\FFF}$} if there is a cover $(\tilde{F},\p)$ of $F$ such that 
$$\tilde{\FFF}\cong \FFF(\tilde{F},\p)\ .$$
\end{itemize}
\end{de}

\begin{satz}
Let $\FFF$ be a foundation and let $\tilde{\FFF}$ be a cover of $\FFF$. Then $\FFF$ is integrable if $\tilde{\FFF}$ is integrable.
\end{satz}

\begin{bew}
This is a consequence of theorem C in \cite{M}.\qed
\end{bew}

\begin{bem}
The next step is to show that the foundation attached to a root group system is unique up to isomorphism. Moreover, we want to prove that the building corresponding to an integrable foundation is unique up to isomorphism.
\end{bem}

\begin{prop}\label{61}
Let $\UUU:=\UUU(\BBB,M,\Sigma,c)$ and $\tilde{\UUU}:=\UUU(\tilde{\BBB},\tilde{M},\tilde{\Sigma},\tilde{c})$ be root group systems and let $\Lambda$ and $\tilde{\Lambda}$ be parameter systems for $\UUU$ and $\tilde{\UUU}$, respectively. Then the following holds:
\begin{enumerate}[label=(\alph*)]
\item An isomorphism $\tilde{\alpha}:\FFF(\UUU,\Lambda)\to \FFF(\tilde{\UUU},\tilde{\Lambda})$ induces an isomorphism $\alpha:\UUU\to \tilde{\UUU}$.
\item An isomorphism $\alpha:\UUU\to \tilde{\UUU}$ induces an isomorphism $\tilde{\alpha}: \FFF(\UUU,\Lambda)\to \FFF(\tilde{\UUU},\tilde{\Lambda})$.
\end{enumerate}
\end{prop}

\begin{bew}
Each isomorphism $$\alpha_{(i,j)}:U_{(i,j)}\to \tilde{U}_{(\pi(i),\pi(j))}$$ induces an isomorphism
$$\tilde{\alpha}_{(i,j)}:\XXX_{(i,j)}(\Xi_{(i,j)})\to \XXX_{(i,j)}(\tilde{\Xi}_{(\pi(i),\pi(j))})$$
and vice versa. Given $(i,j)\in A(M)$, we have
$$ \alpha_{(i,j)}=\alpha_{(j,i)}^o\ \Leftrightarrow\ \tilde{\alpha}_{(i,j)}=\tilde{\alpha}_{(j,i)}^o\ .$$
Now we may go on as in the proof of proposition \ref{71}.\qed
\end{bew}

\chapter{Reparametrizations and Isomorphisms}

The concept of reparametrizations is quite similar to that of isomorphisms. However, we deal with a single foundation and produce (in fact, all the) foundations which are isomorphic to a given one. Moreover, this concept allows us to complete the proof that a foundation is a classifying invariant of the corresponding twin building.

\begin{de}
Let $\FFF$ be a foundation.
\begin{itemize}
\item A system of reparametrizations
$$\alpha:=\{ \alpha_{(i,j)} \mid (i,j)\in A(F) \}$$
satisfying $\alpha_{(i,j)}=\alpha_{(j,i)}^o$ for each $(i,j)\in A(F)$ and
$$\gamma_{(i,j,k)}\circ \alpha_{(i,j)}^j(1)=\alpha_{(j,k)}^j(1)$$
for each $(i,j,k)\in G(F)$ is a \textit{reparametrization for $\mathit{\FFF}$}\index{reparametrization}.
\item Given a reparametrization $\alpha$ for $\FFF$, we set
$$\gls{fr}:=\{ \XXX_{(i,j)}(\tilde{\Xi}_{(i,j)}),\tilde{\gamma}_{(i,j,k)}\mid (i,j)\in A(F),(i,j,k)\in G(F)\}$$
with
$$\tilde{\gamma}_{(i,j,k)}:=(\alpha_{(j,k)}^j)^{-1}\circ \gamma_{(i,j,k)}\circ \alpha_{(i,j)}^j$$
for each $(i,j,k)\in G(F)$.
\end{itemize}
\end{de}

\begin{lemma}
Let $\UUU:=\UUU(\BBB,M,\Sigma,c)$ be a root group system, let $\FFF:=\FFF(\UUU,\Lambda)$ for some parameter system $\Lambda$ for $\UUU$, let 
$\alpha$ be a reparametrization for $\FFF$ and let $\tilde{\Lambda}$ be the parameter system induced by $\alpha$. Then we have $\tilde{\FFF}:=\FFF(\UUU,\tilde{\Lambda})=\FFF_\alpha$.
\end{lemma}

\begin{bew}
We have
\begin{align*}
\tilde{x}_{(j,k)}^j\big(\tilde{\gamma}_{(i,j,k)}(t)\big)&=\tilde{x}_{(i,j)}^j(t)=x_{(i,j)}^j\big(\alpha_{(i,j)}^j(t)\big) \\
&=x_{(j,k)}^j\big(\gamma_{(i,j,k)}\circ \alpha_{(i,j)}^j(t)\big)=\tilde{x}_{(j,k)}^j\big((\alpha_{(j,k)}^j)^{-1}\circ \gamma_{(i,j,k)}\circ \alpha_{(i,j)}^j(t)\big)
\end{align*}
for each $t\in \tilde{\MM}_{(i,j)}^j$.\qed
\end{bew}

\begin{kor}\label{367}
Let $\UUU:=\UUU(\BBB,M,\Sigma,c)$ be a root group system, let $\FFF:=\FFF(\UUU,\Lambda)$ for some parameter system $\Lambda$ for $\UUU$ and let 
$$\alpha=\{\pi,\alpha_{(i,j)}\mid (i,j)\in A(F)\}:\FFF\to \tilde{\FFF}$$
be an isomorphism. Then $\tilde{\FFF}$ is integrable.
\end{kor}

\begin{bew}
Take $\big(\tilde{\Xi}_{(i,j)}:=\tilde{\Xi}_{(\pi(i),\pi(j))}, (\alpha_{(i,j)}^i)^{-1},\ldots,(\alpha_{(i,j)}^j)^{-1}\big)$ as reparametrization for $\XXX_{(i,j)}(\Xi_{(i,j)})$, then replace $i\in I$ by $\pi(i)\in \tilde{I}$. The resulting parameter system $\tilde{\Lambda}$
satisfies
$$\FFF(\UUU,\tilde{\Lambda})=\FFF_\alpha=\tilde{\FFF}\ .$$\qed \end{bew}

\begin{satz}
The isomorphism class of an integrable foundations $\FFF=\FFF(\UUU,\Lambda)$ is a classifying invariant of the isomorphism class of the corresponding building.
\end{satz}

\begin{bew}
This results from corollary \ref{367}, proposition \ref{61} and theorem \ref{368}.\qed
\end{bew}

\begin{bem}
The following theorem shows that the concept of reparametrization is useful if we want to determine all the foundations isomorphic to a given foundation $\FFF$.
\end{bem}

\begin{satz}
Let $\FFF,\tilde{\FFF}$ be foundations with $F=\tilde{F}$. Then the following holds:
\begin{enumerate}[label=(\alph*)]
\item \label{371} Let $\tilde{\alpha}=\{ \tilde{\alpha}_{(i,j)}\mid (i,j)\in A(F)\}:\FFF\to \tilde{\FFF}$ be a special isomorphism. Then there is a reparametrization $\alpha$ of $\FFF$ such that $\FFF_{\alpha}=\tilde{\FFF}$.
\item Let $\alpha=\{ \alpha_{(i,j)} \mid (i,j)\in A(F)$ be a reparametrization for $\FFF$ such that $\FFF_\alpha=\tilde{\FFF}$. Then there is a special isomorphism $\tilde{\alpha}:\FFF\to \tilde{\FFF}$.
\end{enumerate}
\end{satz}

\begin{bewzwei}\ 
\begin{enumerate}[label=(\alph*)]
\item If we take $\alpha:=\{ \alpha_{(i,j)} \mid (i,j)\in A(F)\}$ with
$$\alpha_{(i,j)}:=\{ \tilde{\Xi}_{(i,j)},(\tilde{\alpha}_{(i,j)}^{i})^{-1},\ldots,(\tilde{\alpha}_{(i,j)}^j)^{-1} \}$$
as reparametrization for $\FFF$, then $\FFF_\alpha=\tilde{\FFF}$.
\item We have $$\alpha_{(i,j)}=(\tilde{\Xi}_{(i,j)}, \alpha_{(i,j)}^i,\ldots,\alpha_{(i,j)}^j)$$
for each $(i,j)\in A(F)$, thus $\tilde{\alpha}:=\{ \id_F,\tilde{\alpha}_{(i,j)} \mid (i,j)\in A(F)\}:\FFF\to \tilde{\FFF}$ with
$$\tilde{\alpha}_{(i,j)}:=\big((\alpha_{(i,j)}^i)^{-1},\ldots,(\alpha_{(i,j)}^j)^{-1}\big)$$
is an isomorphism.
\end{enumerate}\qed
\end{bewzwei}

\begin{bem}\label{369}\ 
\begin{enumerate}[label=(\alph*)]
\item Let $\FFF$ and $\tilde{\FFF}$ be foundations and let $$\alpha=\{\pi,\alpha_{(i,j)}\mid (i,j)\in A(F)\}:\FFF\to \tilde{\FFF}$$ be an isomorphism. As we may replace $i\in V(F)$ by $\pi(i)\in V(\tilde{F})$, we may consider $\alpha$ as special. Thus it suffices to determine all foundations which are specially isomorphic to $\FFF$. The remaining foundations isomorphic to $\FFF$ are obtained by relabelings of the vertex set.
\item \label{372} The theorem is useful if we want to show that two given foundations $\FFF$ and $\tilde{\FFF}$ with isomorphic residues $\RRR$ and $\tilde{\RRR}$ are isomorphic. In this case, we may replace $\RRR$ by $\tilde{\RRR}$, observing that there is a relabeling of the corresponding vertices involved.
\end{enumerate}
\end{bem}

\newpage

\addtocontents{toc}{\protect\newpage}
\addtocontents{toc}{\noindent\protect\mbox{}\protect\hrulefill\par}
\part{443-Foundations}\label{201}
\addtocontents{toc}{\noindent\protect\mbox{}\protect\hrulefill\par}

\noindent Now we are ready to turn to the classification of integrable 443-foundations, whose Moufang polygons are two quadrangles and one triangle and whose Coxeter diagrams are complete graphs.

The first step is to exclude quadrangles of type $E_n$, of type $F_4$ and of indifferent type. Then we turn to unitary quadrangles, i.e., quadrangles of pseudo-quadratic form or involutory type. As we restrict to proper parameter systems, there are not many possibilities to glue these polygons together.

The final class is that of quadrangles of quadratic form type, which is rich in integrable foundations. In order to avoid characteristic 2 trouble, there is one point where we restrict to proper quadratic spaces although a small gap is the consequence.\\

\chapter{Definition}

\begin{de} A \textit{443-foundation}\index{443-foundation} is a foundation
$$\FFF:=\{ \BBB_{(1,2)},\BBB_{(2,3)}, \BBB_{(3,1)}, \gamma_{(1,2,3)},\gamma_{(2,3,1)},\gamma_{(3,1,2)}\}$$
such that $\BBB_{(1,2)}$ and $\BBB_{(2,3)}$ are quadrangles and $\BBB_{(3,1)}$ is a triangle.
\end{de}

\begin{no}\ 
\begin{itemize} 
\item Given a 443-foundation $\FFF$, we set
\begin{align*}
\gamma_1:=\gamma_{(3,1,2)}\ , &&\gamma_2:=\gamma_{(1,2,3)}\ , &&\gamma_3:=\gamma_{(2,3,1)}\ .\end{align*}
\item Throughout the rest of this part, $\FFF$ is an integrable 443-foundation.
\end{itemize}
\end{no}

\chapter{The Quadrangles Are Not of Type \texorpdfstring{$\boldsymbol{E_n}$}{En}}\label{390}

\begin{lemma} \label{380}
Let $\BBB=\big(x_1(\MM_1),\ldots,x_4(\MM_4)\big)$ be a parametrized standard quadrangle of type $E_n$. Then the following holds:
\begin{enumerate}[label=(\alph*)]
\item The Moufang set $\MM_1$ is non-commutative.
\item We have $\MM_4=\MM(L_0,\KK,q)$ for some quadratic space $(L_0,\KK,q)$ of type $E_n$.
\end{enumerate}
\end{lemma}

\begin{bewzwei}\ 
\begin{enumerate}[label=(\alph*)]
\item This results from (38.10) of \cite{TW}.
\item This holds by definition, cf. remark \ref{473} or example (16.6) of \cite{TW}.
\end{enumerate}\qed
\end{bewzwei}

\begin{satz}\label{451} The quadrangles are not of type $E_n$.
\end{satz}

\begin{bew}
Suppose that \vphantom{$\MM_{(1,3)}^1$}$\BBB_{(1,2)}$ or $\BBB_{(2,1)}$ is a parametrized standard quadrangle of type $E_n$. The glueing $\gamma_1: \MM_{(3,1)}^1\to \MM_{(1,2)}^1$ is a Jordan isomorphism. We have $\MM_{(3,1)}^1=\MM({\AA})$ for some alternative division ring ${\AA}$, thus $\MM_{(3,1)}^1$ is commutative. Hence we have $\MM_{(1,2)}^1=\MM(\tilde{L}_0,\tilde{\KK},\tilde{q})$ for some quadratic space $(\tilde{L}_0,\tilde{\KK},\tilde{q})$ of type $E_n$ by lemma \ref{380}. \vphantom{$\MM_{(1,3)}^1$}

Notice that $(\tilde{L}_0,\tilde{\KK},\tilde{q})$ is proper by remark \ref{382}. By theorem \ref{379}, the alternative division ring ${\AA}$ is quadratic over a subfield ${\FF}$ of its center ${\KK}:=Z({\AA})$ with $N:=N^\AA_\FF$, and $(\tilde{L}_0,\tilde{\KK},\tilde{q})$ and $({\AA},{\FF},{N})$ are isomorphic as quadratic spaces. In particular, we have $\tilde{\KK}\cong {\FF}$ and
\begin{align*}
\dim_{\tilde{\KK}} \tilde{L}_0=\dim_{{\FF}} {\AA}\in \{ 1,2,4,8 \}\cap \{6,8,12\}=\{8\}\ .
\end{align*}
As a consequence, $({\OO},{\KK},{N})=({\AA},{\FF},{N})$ is of type (v) and $(\tilde{L}_0,\tilde{\KK},\tilde{q})$ is of type $E_7$. But by corollary \ref{381}, we have $(\tilde{L}_0,\tilde{\KK},\tilde{q})\not\cong ({\OO},{\KK},{N})$ as quadratic spaces$\qquad \lightning.$\qed
\end{bew}

\chapter{The Quadrangles Are Not of Type \texorpdfstring{$\boldsymbol{F_4}$}{F4}}

\begin{bem}\label{384}
Let $\QQQ_F(L_0,\KK,q)=\big(x_1(\MM_1),x_2(\MM_2),x_3(\MM_3),x_4(\MM_4)\big)$ be a quadrangle of type $F_4$, where $(L_0,\KK,q)$ is a quadratic space of type $F_4$. 
\begin{enumerate}[label=(\alph*)]
\item By remark (21.18) of \cite{TW}, the quadrangle $\QQQ_F(L_0,\KK,q)$ is an extension of the quadrangle $\QQQ_Q(L_0,\KK,q)$, thus we have $\MM_4=\MM(L_0,\KK,q)$.
\item By remark (21.18) of \cite{TW} again, the quadrangle $\QQQ_F^o(L_0,\KK,q)^o$ is an extension of the quadrangle $\QQQ_Q(\hat{L}_0,\FF,\hat{q})$, thus we have $\MM_1=\MM(\hat{L}_0,\FF,\hat{q})$, where $(\hat{L}_0,\FF,\hat{q})$ is the quadratic space of type $F_4$ as defined in (14.12) of \cite{TW}, cf. (14.13) of \cite{TW}.
\end{enumerate}
\end{bem}

\begin{satz}\label{450} The quadrangles are not of type $F_4$.
\end{satz}

\begin{bew}
Suppose that $\BBB_{(1,2)}$ or $\BBB_{(2,1)}$ is a parametrized standard quadrangle of type $F_4$. The glueing $\gamma_1: \MM_{(3,1)}^1\to \MM_{(1,2)}^1$ is a Jordan isomorphism. We have $\MM_{(3,1)}^1=\MM(\tilde{\AA})$ for some alternative division ring $\tilde{\AA}$. By remark \ref{384}, we have $$\MM_{(1,2)}^1=\MM(L_0,\KK,q)$$ for some quadratic space $(L_0,\KK,q)$ of type $F_4$. By theorem \ref{379}, the alternative division ring $\tilde{\AA}$ is quadratic over a subfield $\tilde{\FF}$ of its center $\tilde{\KK}:=Z(\tilde{\AA})$, and $(L_0,\KK,q)$ and $(\tilde{\AA},\tilde{\FF},N^{\tilde{\AA}}_{\tilde{\FF}})$ are isomorphic as quadratic spaces, which contradicts lemma \ref{385}.\qed
\end{bew}

\chapter{The Quadrangles Are Not of Indifferent Type}

\begin{bem}\label{387}
Let $\BBB:=\QQQ_D(\KK,\KK_0,\LL_0)=\big(x_1(\MM_1),x_2(\MM_2),x_3(\MM_3),x_4(\MM_4)\big)$ be a quadrangle of indifferent type, where $(\KK,\KK_0,\LL_0)$ is a proper indifferent set.
\begin{enumerate}[label=(\alph*)]
\item By definition, the Moufang set $\MM_1=\MM(\KK,\KK_0,\LL_0)$ is of indifferent type.
\item By remark (35.9) of \cite{TW}, we have $\BBB^o=\QQQ_D(\LL,\LL_0,\KK_0^2)$, thus we have $\MM_4=\MM(\LL,\LL_0,\KK_0^2)$, where $(\LL,\LL_0,\KK_0^2)$ is the opposite of  $(\KK,\KK_0,\LL_0)$, which is proper by lemma \ref{395}.
\end{enumerate}
\end{bem}

\begin{satz}\label{449} The quadrangles are not of indifferent type.
\end{satz}

\begin{bew}
Suppose that $\BBB_{(1,2)}$ or $\BBB_{(2,1)}$ is a parametrized standard quadrangle of indifferent type. By remark \ref{387}, we have $$\MM_{(1,2)}^1=\MM(\KK,\KK_0,\LL_0)$$ for some proper indifferent set $(\KK,\KK_0,\LL_0)$, we have $\MM_{(3,1)}^1=\MM(\AA)$ for some alternative division ring $\AA$, and the glueing $\gamma_1: \MM_{(3,1)}^1\to \MM_{(1,2)}^1$ is a Jordan isomorphism, which contradicts theorem \ref{388}.\vphantom{$\MM_{(1,3)}^1$}\qed
\end{bew}

\begin{bem}
Now we are done with the exclusion of certain families of quadrangles. Next we pass to unitary 443-foundations, which can not be obtained as fixed point foundations of covers. As in the $\tilde{A}_2$-case with positive glueings, the parametrizing structures are quaternion division algebras, and the existence can be shown via Tits indices.

Then we finally come to 443-foundations involving quadrangles of quadratic form type, which can be constructed as fixed point foundations of covers. 

In both cases however, we leave off the existence proofs which require different kinds of techniques. 
\end{bem}

\newpage

\chapter{Unitary Quadrangles}

As quadrangles of pseudo-quadratic form type are extensions of quadrangles of involutory type (which are not necessarily of purely involutory type) and as quadrangles of involutory type can be considered as quadrangles of pseudo-quadratic form type of a non-proper pseudo-quadratic space, it is natural to treat them in a common setup. As a consequence, we sometimes omit the assumption of a proper parameter system to obtain a general statement for both the families.

\section{Definitions}

\begin{de}\label{186}\ 
\begin{itemize}
\item  A foundation
$$\FFF:=\{ \BBB_{(1,2)}=\QQQ_I^o(\hat{\Xi}),\BBB_{(2,3)}=\QQQ_I({\Xi}), \BBB_{(3,1)}=\TTT(\tilde{\KK}), \gamma_{(1,2,3)},\gamma_{(2,3,1)},\gamma_{(3,1,2)}\}$$
for some proper involutory sets $\Xi$ and $\hat{\Xi}$ is a \textit{443-foundation of involutory type}\index{443-foundation!of involutory type}.
\item  A foundation
$$\FFF:=\{ \BBB_{(1,2)}=\QQQ_P^o(\hat{\Xi}),\BBB_{(2,3)}=\QQQ_P({\Xi}), \BBB_{(3,1)}=\TTT(\tilde{\KK}), \gamma_{(1,2,3)},\gamma_{(2,3,1)},\gamma_{(3,1,2)}\}$$
for some proper pseudo-quadratic spaces $\Xi$ and $\hat{\Xi}$ is a \textit{443-foundation of pseudo-quadratic form type}\index{443-foundation!of pseudo-quadratic form type}.
\item A 443-foundation is of \textit{unitary type}\index{443-foundation!of unitary type} if it is either of involutory type or of pseudo-quadratic form type.
\end{itemize}
\end{de}

\begin{no}
Given a 443-foundation $\FFF$, we set
\begin{align*}
\gamma_1:=\gamma_{(3,1,2)}\ , && \gamma_2:=\gamma_{(1,2,3)}\ , && \gamma_3:=\gamma_{(2,3,1)}\ .
\end{align*}
\end{no}

\begin{lemma}
Let $\FFF$ be an integrable 443-foundation of unitary type. Then $\tilde{\KK}$ is associative.
\end{lemma}

\begin{bew}
The glueings $\gamma_3=\gamma_{(1,3,2)}$ and $\gamma_1=\gamma_{(3,1,2)}$ are positive or negative by Hua's theorem. In particular, $\tilde{\KK}$ is associative.
\qed
\end{bew}

\begin{lemma}\label{158}
Let $\KK$ be a skew-field, let $M\subseteq \KK$ and let $a,b,c\in \KK$ such that
\begin{align*}
1_\KK\in M\ , && \forall\ x\in M:\qquad axb=cx\ .
\end{align*}
Then we have
\begin{align*} b\in C_\KK(M)\ , && M\subseteq C_\KK(b)\ .\end{align*}
\end{lemma}

\begin{bew}
We have
$$c=c\cdot 1_\KK=a\cdot 1_\KK\cdot b=ab$$
and therefore
\begin{align*}
\forall\ x\in M:\qquad  xb=a^{-1}(axb)=a^{-1}(abx)=bx\ .
\end{align*}\qed
\end{bew}

\begin{satz}\label{192}
Let $\FFF$ be an integrable 443-foundation of pseudo-quadratic form type with $\hat{\Xi}=\Xi^o$ and $\gamma_2=\id_T^o$ and let
$$\pi_\KK:T\to \KK,\ (a,t)\mapsto t\ .$$ If one of the glueings is negative, we have
$$\pi_\KK(T)\subseteq Z(\KK)\ .$$
In particular, we have $\KK_0\subseteq Z(\KK)$.
\end{satz}

\begin{bew}
Assume that $\gamma_3=\gamma_{(2,3,1)}$ is negative (otherwise, we consider the opposite buildings and glueings). By taking $(\KK,\gamma_3,\gamma_3,\gamma_3)$ as reparametrization for $\TTT(\tilde{\KK})$, we may assume
\begin{align*}
\tilde{\KK}=\KK\ , && \gamma_3=\id_{\KK}\ .
\end{align*}
If we set 
\begin{align*}
h(s)&:= \mu\big(x_{(3,1)}^1(1_\KK)\big)^{-1}\mu\big(x_{(3,1)}^1(s^{-1})\big)\ , \qquad s\in \KK^*\ , \\
\tilde{h}(s)&:= \mu\big(x_{(1,2)}^1(1_\KK)\big)^{-1}\mu\big(x_{(1,2)}^1(s^{-1})\big)\ , \qquad s\in \KK^*\ ,
\end{align*}
we have
$$h(s)=\mu\big(x_{(3,1)}^1(1_\KK)\big)^{-1}\mu\big(x_{(3,1)}^1(s^{-1})\big)=\mu\big(x_{(1,2)}^1(1_\KK)\big)^{-1}\mu\big(x_{(1,2)}^1(\gamma_1(s)^{-1})\big)=\tilde{h}\big(\gamma_1(s)\big)$$
for each $s\in \KK^*$,
\begin{align*}
x_{(2,3)}^3(t)^{h(s)}&=x_{(3,1)}^3(t)^{h(s)}=x_{(3,1)}^3(ts)=x_{(2,3)}^3(ts)
\end{align*}
for all $s\in \KK^*,\ t\in \KK$ by lemma \ref{361} and
\begin{align*}
x_{(2,3)}^2(a,t)^{h(s)}&=x_{(1,2)}^2(a,t)^{\tilde{h}(\gamma_1(s))}\\
&=x_{(1,2)}^2\big(\gamma_1(s)\circ a,\gamma_1(s)\circ t\circ \gamma_1(s)^{\sigma}\big)=x_{(2,3)}^2\big(a\cdot\gamma(s),\gamma_1(s)^\sigma \cdot t\cdot \gamma_1(s)\big)
\end{align*}
for all $s\in \KK^*,\ (a,t)\in T$ by lemma \ref{322}. Given $s\in \KK^*$, the Hua automorphism $h\big(\gamma_1^{-1}(s)\big)$ induces an automorphism $\alpha_s\in \Aut(\BBB_{(2,3)})$ which satisfies 
\begin{align*}
x_{(2,3)}^2(a,t)\mapsto x_{(2,3)}^2(a\cdot s,s^\sigma\cdot  t\cdot s)\ , && x_{(2,3)}^3(t)\mapsto x_{(2,3)}^3(t\cdot \gamma_1^{-1}(s))\ .
\end{align*}
If we set $\tilde{\alpha}:=\alpha_{(\gamma_1^{-1}(s),1_\KK,(\id_{L_0},\id_\KK))}$ as in (37.33) of \cite{TW}, then $\tilde{\alpha}\alpha_s$ satisfies
\begin{align*}
x_{(2,3)}^2(a,t)\mapsto x_{(2,3)}^2(a\cdot s , s^\sigma\cdot  t \cdot s)\ , && x_{(2,3)}^3(t)\mapsto x_{(2,3)}^3(t)\ .
\end{align*}
By (37.33) of \cite{TW}, there is an element $c\in \KK_0^*$ such that
$$\forall\ t\in \pi_\KK(T):\qquad s^\sigma ts=c t\ .$$
Lemma \ref{158} implies that we have
$$\pi_\KK(T)\subseteq C_\KK(s)\ .$$
Since $s\in \KK^*$ is arbitrary, it follows that
$$\pi_\KK(T)\subseteq Z(\KK)\ .$$
\qed
\end{bew}

\begin{bem}
Notice that we don't need the fact that $\Xi$ is proper by definition of our parameter systems, i.e., if we allow $\Xi$ to be non-proper, the theorem remains true. As a consequence, we get a similar result for 443-foundations of involutory type which can be considered as 443-foundations of pseudo-quadratic form type for some non-proper pseudo-quadratic spaces, cf. lemma \ref{176}.
\end{bem}

\section{Quadrangles of Pseudo-Quadratic Form Type}

\begin{no} Throughout this paragraph, $\FFF$ is an integrable $443$-foundation such that at least one quadrangle is of pseudo-quadratic form type.\end{no}

\begin{lemma}\label{392}
The foundation $\FFF$ is of pseudo-quadratic form type.
\end{lemma}

\begin{bew}
\vphantom{$\MM_{(3,1)}^1$} We may assume that $\BBB_{(1,2)}$ or $\BBB_{(2,1)}$ is a standard quadrangle $\QQQ_P(\hat{\KK},\hat{\KK}_0,\hat{\sigma},\hat{L}_0,\hat{q})$ of pseudo-quadratic form type and thus $\MM_{(1,2)}^1=\MM(\hat{\KK})$ or $\MM_{(1,2)}^1=\MM(\hat{\KK},\hat{\KK}_0,\hat{\sigma},\hat{L}_0,\hat{q})$. \vphantom{$\MM_{(3,1)}^1$} But since the map $\gamma_1:\MM(\tilde{\KK})=\MM_{(3,1)}^1\to \MM_{(1,2)}^1$ is a Jordan isomorphism and Moufang sets of linear type are commutative while Moufang sets of pseudo-quadratic form type are not, we obtain that \vphantom{$\MM_{(3,1)}^1$}
$$\BBB_{(1,2)}=\QQQ^o_P\big((\hat{\KK},\hat{\KK}_0,\hat{\sigma},\hat{L}_0,\hat{q})^o\big)\ .$$
Now $\BBB_{(2,3)}$ is a Moufang quadrangle such that $\MM_{(2,3)}^2$ is non-commutative, and since we excluded quadrangles of type $E_n$ in chapter \ref{390}, we have
$$\BBB_{(2,3)}=\QQQ_P(\Xi)$$
for some proper pseudo-quadratic space $\Xi$.\qed
\end{bew}

\begin{no} Until proposition \ref{193}, $\FFF$ is an integrable $443$-foundation of pseudo-quadratic form type such that $\KK$ is non-commutative.\end{no}

\begin{bem} \label{194}
By theorem \ref{283}, the Jordan isomorphism $\gamma_2=\gamma_{(1,2,3)}:\hat{T}\to T$
is induced by an isomorphism $\Phi:\hat{\Xi}^o\to {\Xi}$ of pseudo-quadratic spaces. By taking $(\Xi^o,{\phi}^{-1},{\Phi}^{-1},{\phi}^{-1},{\Phi}^{-1})$ as reparametrization for $\QQQ_P^o(\hat{\Xi})$, we may assume
\begin{align*}
\hat{\Xi}=\Xi^o\ , && \gamma_2=\id^o_{T}\ .
\end{align*}
\end{bem}

\begin{lemma}
Both the glueings $\gamma_1$ and $\gamma_3$ are positive.
\end{lemma}

\begin{bew}
If one of the glueings is negative, then we have
$$\KK_0\subseteq \pi_\KK(T)\subseteq Z(\KK)$$
by theorem \ref{192}, thus $(\KK,\KK_0,\sigma)$ is non-proper by lemma \ref{157}. By remark \ref{191}, we have $$(\KK,\KK_0,\sigma)=(\HH,Z(\HH),\sigma_s)$$ for some quaternion division algebra $\HH$. Given $a\in L_0$, we have
\begin{align*} q(a)\in \pi_\HH(T)\subseteq Z(\HH)\ , && a=0_{L_0}\ \end{align*}
and thus $L_0=\{0_{L_0}\}$. But then $\Xi$ is non-proper$\qquad \lightning$.\qed
\end{bew}

\begin{prop}\label{193}
The skew-field $\KK$ is a quaternion division algebra and $\FFF$ is isomorphic to the foundation
$$\FFF_{443}(\Xi):=\{ \tilde{\BBB}_{(1,2)}=\QQQ_P^o(\Xi^o),\tilde{\BBB}_{(2,3)}=\QQQ_{P}(\Xi), \tilde{\BBB}_{(3,1)}=\TTT(\KK^o),\tilde{\gamma}_1=\sigma_s,\tilde{\gamma}_2=\id_{T}^o,\tilde{\gamma}_3=\id_\KK^o\} \ .$$
\end{prop}

\begin{bew}
As $\gamma_3=\gamma_{(2,3,1)}$ is positive, we may take $(\KK^o,\gamma_3^o,\gamma^o_3,\gamma_3^o)$ as reparametrization for $\TTT(\tilde{\KK})$. Therefore, we may assume
\begin{align*}
\tilde{\KK}=\KK^o\ , && \gamma_3=\id_{\KK}^o\ .
\end{align*}
Let $s\in \KK^o$. As in the proof of proposition \ref{162}, we obtain an automorphism $\alpha_s\in \Aut(\BBB_{(2,3)})$ satisfying
\begin{align*}
x_{(2,3)}^2(a,t)\mapsto x_{(2,3)}^2\big(a\cdot \gamma_1^o(s)\cdot \id_\KK^o(s),\id_\KK^o(s)^{\sigma}\cdot \gamma_1^o(s)^\sigma\cdot t\cdot \gamma_1^o(s)\cdot \id_\KK^o(s)\big)\ , && x_{(2,3)}^3(t)\mapsto x_{(2,3)}^3(t)\ .
\end{align*}
By (37.33) of \cite{TW}, the map
$$\p_1:L_0\to L_0,\ a\mapsto a\cdot \gamma_1^o(s)\cdot \id_\KK^o(s)$$
is an isomorphism of vector spaces satisfying
\begin{align*}
\forall\ a\in L_0,\ t\in \KK:\qquad \p_1(a\cdot t)=\p_1(a)\cdot t .
\end{align*}
As $\Xi$ is proper, we have $L_0\neq \{0_{L_0}\}$ and thus
\begin{align*}
\forall\ t\in \KK:\qquad t\cdot \gamma_1^o(s)\cdot \id_\KK^o(s)=\gamma_1^o(s)\cdot \id_\KK^o(s)\cdot t\ , && \gamma_1^o(s)\cdot \id_\KK^o(s)\in Z(\KK)\ .
\end{align*}
Since $\KK$ is non-commutative by assumption, lemma \ref{160} shows that $\KK$ is a quaternion division algebra and that we have
\begin{align*} \gamma_1= \sigma_s:\KK^o\to \KK^o\ . \end{align*}
\qed
\end{bew}

\begin{bem}\label{391} Let $\FFF$ be an integrable 443-foundation of pseudo-quadratic form type such that $\KK$ is a field, but $\KK\not\cong \FF_4$ if $\dim_\KK L_0=1$. Then $(\KK,\KK_0,\sigma)$ is non-proper by lemma \ref{157} and thus quadratic of type (iii) by remark \ref{191}. Moreover, we may reparametrize as in the non-commutative case, but we cannot get more information concerning the glueing $\gamma_1$, i.e., we have
$$\FFF\cong \FFF_{443}(\Xi,\gamma):=\{ \tilde{\BBB}_{(1,2)}=\QQQ_P^o(\Xi^o),\tilde{\BBB}_{(2,3)}=\QQQ_{P}(\Xi), \tilde{\BBB}_{(3,1)}=\TTT(\KK),\tilde{\gamma}_1=\gamma,\tilde{\gamma}_2=\id_{T}^o,\tilde{\gamma}_3=\id_\KK\}$$
for some pseudo-quadratic space $\Xi$ such that $(\KK,\KK_0,\sigma)$ is quadratic of type (iii) and for some $\gamma\in \Aut(\KK)$.
\end{bem}

\begin{satz}\label{448}
Let $\FFF$ be an integrable 443-foundation such that at least one quadrangle is of pseudo-quadratic form type and such that $\tilde{\KK}\not\cong \FF_4$ if $\dim_\KK L_0=1$, where $\BBB_{(3,1)}=\TTT(\tilde{\KK})$. Then one of the following holds:
\begin{enumerate}[label=(\roman*)]
\item We have
$$\FFF\cong \FFF_{443}(\Xi)=\{ \tilde{\BBB}_{(2,1)}=\tilde{\BBB}_{(2,3)}=\QQQ_{P}(\Xi), \tilde{\BBB}_{(3,1)}=\TTT(\KK^o),\tilde{\gamma}_1=\sigma_s,\tilde{\gamma}_2=\id_{T}^o,\tilde{\gamma}_3=\id_\KK^o\}$$
for some proper pseudo-quadratic space $\Xi=(\KK,\KK_0,\sigma,L_0,q)$ such that $\KK$ is a quaternion division algebra and $\sigma_s$ is its standard involution.
\item We have 
$$\FFF\cong \FFF_{443}(\Xi,\gamma)=\{ \tilde{\BBB}_{(2,1)}=\tilde{\BBB}_{(2,3)}=\QQQ_{P}(\Xi), \tilde{\BBB}_{(3,1)}=\TTT(\KK),\tilde{\gamma}_1=\gamma,\tilde{\gamma}_2=\id_{T}^o,\tilde{\gamma}_3=\id_\KK\}$$
for some proper pseudo-quadratic space $\Xi=(\KK,\KK_0,\sigma,L_0,q)$ such that $(\KK,\KK_0,\sigma)$ is quadratic of type (iii) and for some $\gamma\in \Aut(\KK)$.
\end{enumerate}
\end{satz}

\begin{bew}
This results from lemma \ref{392}, proposition \ref{193} and remark \ref{391}.\qed
\end{bew}

\section{Quadrangles of Involutory Type}

\begin{no} Throughout this paragraph, $\FFF$ is an integrable $443$-foundation such that at least one quadrangle is of involutory type.\end{no}

\begin{lemma}\label{393}
The foundation $\FFF$ is of involutory type.
\end{lemma}

\begin{bew}
We may assume that $\BBB_{(1,2)}$ or $\BBB_{(2,1)}$ is a standard quadrangle $\QQQ_I(\hat{\KK},\hat{\KK}_0,\hat{\sigma})$ of involutory type and thus $\MM_{(1,2)}^1=\MM(\hat{\KK})$ or $\MM_{(1,2)}^1=\MM(\hat{\KK},\hat{\KK}_0,\hat{\sigma})$. Since $(\hat{\KK},\hat{\KK}_0,\hat{\sigma})$ is proper, we have 
$$\MM(\tilde{\KK})=\MM_{(3,1)}^1\cong\MM_{(1,2)}^1\neq \MM( \hat{\KK},\hat{\KK}_0,\hat{\sigma})$$
by theorem \ref{389} and thus
$$\BBB_{(1,2)}=\QQQ^o_I\big((\hat{\KK},\hat{\KK}_0,\hat{\sigma})^o\big)\ .$$
Since we excluded quadrangles of type $E_n,F_4$, of indifferent type and of pseudo-quadratic form type in the previous paragraphs, the quadrangle $\BBB_{(2,3)}$ is either of quadratic form type or of involutory type. 

Assume that $\BBB_{(2,3)}$ is of quadratic form type. Then we have $\MM_{(2,3)}^2=\MM(\KK)$ or $\MM_{(2,3)}^3=\MM(\KK)$ for some field $\KK$. But we have
\begin{align*}
\MM_{(2,3)}^3\cong \MM_{(3,1)}^3=\MM(\tilde{\KK})=\MM_{(3,1)}^1\cong \MM_{(1,2)}^1=\MM(\hat{\KK})\ ,
\end{align*}
where $\hat{\KK}$ is a non-commutative skew-field since $(\hat{\KK},\hat{\KK}_0,\hat{q})$ is proper, and thus $\MM_{(2,3)}^3\not\cong \MM(\KK)$. Moreover, we have
$$\MM_{(2,3)}^2\cong \MM_{(1,2)}^2=\MM(\hat{\KK},\hat{\KK}_0,\hat{\sigma})$$
and thus $\MM_{(2,3)}^2\neq \MM(\KK)$ by theorem \ref{389}$\qquad \lightning$. 

Therefore, the quadrangle $\BBB_{(2,3)}=$ is of involutory type, i.e., $\BBB_{(2,3)}$ or $\BBB_{(3,2)}$ is a standard quadrangle $\QQQ_I(\KK,\KK_0,\sigma)$. Since we have
$$\MM_{(2,3)}^3\cong \MM_{(3,1)}^3=\MM(\tilde{\KK})$$
and $(\KK,\KK_0,\sigma)$ is proper, we have $\MM_{(2,3)}^3\not\cong \MM(\KK,\KK_0,\sigma)$ by theorem \ref{389} again and thus
$$\BBB_{(2,3)}=\QQQ_I(\KK,\KK_0,\sigma)\ .$$\qed
\end{bew}

\begin{no} Throughout the rest of this paragraph, $\FFF$ is an integrable $443$-foundation of involutory type.\end{no}

\begin{bem}\label{185} \ 
\item \label{159} By theorem \ref{175}, the Jordan isomorphism
$$\gamma_2=\gamma_{(1,2,3)}:\hat{\KK}_0\to \KK_0$$
is induced by an isomorphism $$\phi:(\hat{\KK},\hat{\KK}_0,\hat{\sigma})\to (\KK,\KK_0,\sigma)=:\Xi$$
of involutory sets. As a consequence, the map $$\tilde{\phi}:=\phi\circ \sigma^o: (\hat{\KK}^o,\hat{\KK}^o_0,\hat{\sigma})\to (\KK,\KK_0,\sigma)$$
is an isomorphism of involutory sets that induces $\gamma_2$ as well. By taking $(\Xi^o,\tilde{\phi}^{-1},\tilde{\phi}^{-1},\tilde{\phi}^{-1},\tilde{\phi}^{-1})$ as reparametrization for $\QQQ_I^o(\hat{\KK},\hat{\KK}_0,\hat{\sigma})$, we may assume
\begin{align*}
(\hat{\KK},\hat{\KK}_0,\hat{\sigma})=(\KK^o,\KK_0^o,\sigma)\ , && \gamma_2=\id^o_{\KK_0}\ .
\end{align*}
\end{bem}

\begin{lemma}\label{176}
Both the glueings $\gamma_1$ and $\gamma_3$ are positive.
\end{lemma}

\begin{bew}
Notice that we have
$$\QQQ_I(\KK,\KK_0,\sigma)=\QQQ_P(\KK,\KK_0,\sigma,L_0,q)$$
for $L_0:=\{0\}$, $q:=0$. If one of the glueings is negative, we have
$$\KK_0=\pi_\KK(T)\subseteq Z(\KK)$$
by theorem \ref{192}, where we did not use the fact that $(\KK,\KK_0,\sigma,L_0,q)$ is proper in the given setup. But then $(\KK,\KK_0,\sigma)$ is non-proper by lemma \ref{157}$\qquad \lightning$.\qed
\end{bew}

\begin{prop}\label{162}
The skew-field $\KK$ is a quaternion division algebra and $\FFF$ is isomorphic to the foundation
$$\FFF_{443}(\Xi):=\{ \tilde{\BBB}_{(1,2)}=\QQQ_I^o(\Xi^o),\tilde{\BBB}_{(2,3)}=\QQQ_{I}(\Xi), \tilde{\BBB}_{(3,1)}=\TTT(\KK^o),\tilde{\gamma}_1=\sigma_s,\tilde{\gamma}_2=\id_{\KK_0}^o,\tilde{\gamma}_3=\id_\KK^o\} \ .$$
\end{prop}

\begin{bew}
As $\gamma_3=\gamma_{(2,3,1)}$ is positive, we may take $(\KK^o,\gamma_3^o,\gamma^o_3,\gamma_3^o)$ as reparametrization for $\TTT(\tilde{\KK})$. Therefore, we may assume
\begin{align*}
\tilde{\KK}=\KK^o\ , && \gamma_3=\id_{\KK}^o\ .
\end{align*}
\begin{itemize}
\item If we set
\begin{align*}
h_1(s)&:= \mu\big(x_{(3,1)}^1(1_{\KK^o})\big)^{-1}\mu\big(x_{(3,1)}^1(s^{-1})\big)\ ,  && s\in \KK^o\ , \\
\tilde{h}_1(s)&:= \mu\big(x_{(1,2)}^1(1_\KK^o)\big)^{-1}\mu\big(x_{(1,2)}^1(s^{-1})\big)\ ,  && s\in \KK^o\ ,
\end{align*}
we have
\begin{align*} h_1(s)&=\mu\big(x_{(3,1)}^1(1_{\KK^o})\big)^{-1}\mu\big(x_{(3,1)}^1(s^{-1})\big)\\
&=\mu\big(x_{(1,2)}^1(1_\KK^o)\big)^{-1}\mu\big(x_{(1,2)}^1(\gamma_1(s)^{-1})\big)=\tilde{h}_1\big(\gamma_1(s)\big)\end{align*}
for each $s\in \KK^o$,
\begin{align*}
x_{(2,3)}^3(t)^{h_1(s)}&=x_{(3,1)}^3\big(\id_{\KK}^o(t)\big)^{h_1(s)}=x_{(3,1)}^3\big(\id_{\KK}^o(t)\circ s\big)=x_{(2,3)}^3\big(\id_\KK^o(s)\cdot t\big)
\end{align*}
for all $s\in \KK^o,\ t\in \KK$ by lemma \ref{361} and
\begin{align*}
x_{(2,3)}^2(t)^{h_1(s)}&=x_{(1,2)}^2\big(\id_{\KK_0}^o(t)\big)^{\tilde{h}_1(\gamma_1(s))}\\
&=x_{(1,2)}^2\big(\gamma_1(s)\circ \id_{\KK_0}^o(t)\circ \gamma_1(s)^{\sigma}\big)=x_{(2,3)}^2\big(\gamma_1^o(s)^\sigma\cdot  t\cdot  \gamma_1^o(s)\big)
\end{align*}
for all $s\in \KK^o,\ t\in \KK_0$ by lemma \ref{161}.
\item If we set
\begin{align*}
h_3(s)&:= \mu\big(x_{(3,1)}^3(1_{\KK^o})\big)^{-1}\mu\big(x_{(3,1)}^3(s)\big)\ , && s\in \KK^o\ , \\
\tilde{h}_3(s)&:= \mu\big(x_{(2,3)}^3(1_\KK)\big)^{-1}\mu\big(x_{(2,3)}^3(s)\big)\ , && s\in \KK\ ,
\end{align*}
we have
$$h_3(s)=\mu\big(x_{(3,1)}^1(1_{\KK^o})\big)^{-1}\mu\big(x_{(3,1)}^1(s^{-1})\big)=\mu\big(x_{(1,2)}^1(1_\KK)\big)^{-1}\mu\big(x_{(1,2)}^1(\id_\KK^o(s^{-1}))\big)=\tilde{h}_3\big(\id_\KK^o(s)\big)$$
for each $s\in \KK^o$,
\begin{align*}
x_{(2,3)}^3(t)^{h_3(s)}&=x_{(3,1)}^3\big(\id_\KK^o(t)\big)^{h_3(s)}\\
&=x_{(3,1)}^3\big(s^{-1}\circ \id_{\KK}^o(t)\circ s^{-1}\big)=x_{(2,3)}^3\big(\id_\KK^o(s)^{-1}\cdot t\cdot\id_\KK^o(s)^{-1}\big)
\end{align*}
for all $s\in \KK^o,\ t\in \KK$ by lemma \ref{361} and
\begin{align*}
x_{(2,3)}^2(t)^{h_3(s)}&=x_{(2,3)}^2(t)^{\tilde{h}_3(\id_\KK^o(s))}=x_{(2,3)}^2\big(\id_\KK^o(s)^{\sigma}\cdot t\cdot \id_\KK^o(s)\big)
\end{align*}
for all $s\in \KK^o,\ t\in \KK_0$ by (33.13) of \cite{TW}.
\end{itemize}

\newpage
 
\noindent 
Therefore, given $s\in \KK^o$, $h_3(s)h_1(s)$ induces an automorphism $\alpha_s\in \Aut(\BBB_{(2,3)})$ which satisfies 
\begin{align*}
x_{(2,3)}^2(t)\mapsto x_{(2,3)}^2\big(\id_\KK^o(s)^{\sigma}\cdot \gamma_1^o(s)^\sigma\cdot t\cdot \gamma_1^o(s)\cdot \id_\KK^o(s)\big)\ , && x_{(2,3)}^3(t)\mapsto x_{(2,3)}^3\big(t\cdot \id_\KK^o(s)^{-1}\big)\ .
\end{align*}
If we set $\tilde{\alpha}:=\alpha_{(\id_\KK^o(s),1_\KK,(\id_{L_0},\id_\KK))}$ as in (37.33) of \cite{TW}, then $\tilde{\alpha}\alpha_s\in \Aut(\BBB_{(2,3)})$ satisfies
\begin{align*}
x_{(2,3)}^2(t)\mapsto x_{(2,3)}^2\big(\id_\KK^o(s)^{\sigma}\cdot \gamma_1^o(s)^\sigma\cdot t\cdot \gamma_1^o(s)\cdot \id_\KK^o(s)\big)\ , && x_{(2,3)}^3(t)\mapsto x_{(2,3)}^3(t)\ .
\end{align*}
By (37.33) of \cite{TW}, there is an element $c\in \KK_0^*$ such that
$$\forall\ t\in \KK_0:\qquad \id_\KK^o(s)^{\sigma}\cdot \gamma_1^o(s)^\sigma\cdot t\cdot \gamma_1^o(s)\cdot \id_\KK^o(s)=c \cdot t\ .$$
Lemma \ref{158} implies that we have
$$\gamma_1^o(s) \cdot \id_\KK^o(s)\in C_\KK(\KK_0)=C_\KK(\langle \KK_0\rangle)=Z(\KK)\ .$$
Since $\KK$ is non-commutative by lemma \ref{157}, lemma \ref{160} shows that $\KK$ is a quaternion division algebra and that we have
\begin{align*} \gamma_1= \sigma_s:\KK^o\to \KK^o\ . \end{align*}
\qed
\end{bew}

\begin{satz}\label{165}
Let $\FFF$ be an integrable 443-foundation such that at least one quadrangle is of involutory type. Then we have
$$\FFF\cong \FFF_{443}(\Xi)=\{ \tilde{\BBB}_{(1,2)}=\QQQ_I^o(\Xi^o),\tilde{\BBB}_{(2,3)}=\QQQ_{I}(\Xi), \tilde{\BBB}_{(3,1)}=\TTT(\KK^o),\tilde{\gamma}_1=\sigma_s,\tilde{\gamma}_2=\id_{\KK_0}^o,\tilde{\gamma}_3=\id_\KK^o\}$$
for some proper involutory set $\Xi=(\KK,\KK_0,\sigma)$ such that $\KK$ is a quaternion division algebra.
\end{satz}

\begin{bew} 
This results from lemma \ref{393} and proposition \ref{162}.\qed
\end{bew}

\begin{bem}\ 
\begin{enumerate}[label=(\alph*)]
\item By remark (11.2) of \cite{TW}, the pair $(\KK,\sigma)$ uniquely determines $\KK_0$ if we have $\Char \KK\neq 2$. 
\item By lemma (3.1.6) of \cite{K} and (35.7) of \cite{TW}, we may assume $\sigma=\gamma_s$ if we have $\Char \KK=2$ and $\sigma$ is an involution of the first kind.
\end{enumerate}
\end{bem}

\chapter{Quadrangles of Quadratic Form Type}

By the previous chapters, there is only one case left: Both the quadrangles are of quadratic form type. 

\begin{no}\ \begin{itemize}
\item Throughout this chapter, $\FFF$ is an integrable $443$-foundation such that both the quadrangles are of quadratic form type. 
\item Given a foundation, if neither $\gamma_{(i,j,k)}$ nor $\gamma_{(k,j,i)}$ is specified, these glueings are supposed to be the identity map.
\end{itemize}
\end{no}

\begin{prop}\label{404}
If we have
$$\FFF=\{ \BBB_{(2,1)}=\QQQ_Q(\tilde{L}_0,\tilde{\KK},\tilde{q}),\BBB_{(2,3)}=\QQQ_Q(\hat{L}_0,\hat{\KK},\hat{q}),\BBB_{(3,1)}=\TTT(\AA),\gamma_1,\gamma_2,\gamma_3\}\ ,$$
then $\AA$ is quadratic over subfields $\FF_1,\FF_2$ of its center, and $\FFF$ is isomorphic to the foundation
$$\FFF_{443}\big(\AA,(\FF_1,\FF_2),\gamma_2\big):=\{ \BBB_{(2,1)}=\QQQ_Q(\AA,\FF_1,N^\AA_{\FF_1}),\BBB_{(2,3)}=\QQQ_Q(\AA,\FF_2,N^\AA_{\FF_2}),\BBB_{(3,1)}=\TTT(\AA),\gamma_2\}$$
for some isomorphism $\gamma_2:\FF_1\to \FF_2$ of fields.
\end{prop}

\begin{bew}
We have
$$\MM(\tilde{L}_0,\tilde{\KK},\tilde{q})=\MM_{(2,1)}^1\cong \MM_{(3,1)}^1=\MM(\AA)=\MM_{(3,1)}^3\cong \MM_{(2,3)}^3=\MM(\hat{L}_0,\hat{\KK},\hat{q})\ .$$
By theorem \ref{379}, the alternative division ring $\AA$ is quadratic over subfields $\FF_1,\FF_2$ of its center, and the maps
\begin{align*}
\gamma_1:\AA\to \tilde{L}_0\ , && \gamma_3^{-1}:\AA\to \hat{L}_0
\end{align*}
are induced by isomorphisms
\begin{align*}
(\gamma_1,\phi_1):(\AA,\FF_1,N^{\AA}_{\FF_1})\to (\tilde{L}_0,\tilde{\KK},\tilde{q})\ , && (\gamma_3^{-1},\phi_3^{-1}):(\AA,\FF_2,N^{\AA}_{\FF_2})\to (\hat{L}_0,\hat{\KK},\hat{q})
\end{align*}
of quadratic spaces. Now we take
\begin{align*}
\alpha_{(2,1)}:=(\Xi_1,\phi_1,\gamma_1,\phi_1,\gamma_1)\ , && \alpha_{(2,3)}:=(\Xi_2,\phi_3^{-1},\gamma_3^{-1},\phi_3^{-1},\gamma_3^{-1})
\end{align*}
as reparametrizations for $\BBB_{(2,1)}$ and $\BBB_{(2,3)}$, respectively, where $\Xi_i:=(\AA,\FF_i,N^\AA_{\FF_i})$.\qed
\end{bew}

\begin{prop}\label{405}
If we have
$$\FFF=\{ \BBB_{(1,2)}=\QQQ_Q(\tilde{L}_0,\tilde{\KK},\tilde{q}),\BBB_{(2,3)}=\QQQ_Q(\hat{L}_0,\hat{\KK},\hat{q}),\BBB_{(3,1)}=\TTT(\KK),\gamma_1,\gamma_2,\gamma_3\}\ ,$$
then $\KK$ is  a field which is quadratic over some subfield $\FF$ of its center, $\FF$ is quadratic over some subfield $\EE$ of its center, and we have
$$\FFF\cong \FFF_{443}(\KK,\FF,\EE,\gamma_1):=\{ \BBB_{(1,2)}=\QQQ_Q(\FF,\EE,N^\FF_\EE),\BBB_{(2,3)}=\QQQ_Q(\KK,\FF,N^\KK_\FF),\BBB_{(3,1)}=\TTT(\KK), \gamma_1\}$$
for some isomorphism $\gamma_1:\KK\to \EE$ of fields.
\end{prop}

\begin{bew}
Since we have
$$\MM(\KK)=\MM_{(3,1)}^1\cong \MM_{(1,2)}^1=\MM(\tilde{\KK})\ ,$$
the alternative division ring $\KK\cong \tilde{\KK}$ is a field by Hua's theorem. We have
$$\MM(\hat{L}_0,\hat{\KK},\hat{q})=\MM_{(2,3)}^3\cong\MM_{(3,1)}^3=\MM(\KK)\ ,$$
therefore, $\KK$ is quadratic over some subfield $\FF$ of its center by theorem \ref{379}, and the map $\gamma_3^{-1}:\KK\to \hat{L}_0$ is induced by an isomorphism $(\gamma_3^{-1},\phi_3^{-1}):(\KK,\FF,N^{\KK}_{\FF})\to (\hat{L}_0,\hat{\KK},\hat{q})$ of quadratic spaces. By taking
$$\alpha_{(2,3)}:=\big((\KK,\FF,N^\KK_\FF),\phi_3^{-1},\gamma_3^{-1},\phi_3^{-1},\gamma_3^{-1}\big)$$
as reparametrization for $\BBB_{(2,3)}$, we may assume $\BBB_{(2,3)}=\QQQ_Q(\KK,\FF,N^\KK_\FF)$ and $\gamma_3=\id_\KK$. Moreover, we have
$$\MM({\FF})=\MM_{(2,3)}^2\cong \MM_{(1,2)}^2=\MM(\tilde{L}_0,\tilde{\KK},\tilde{q})\ ,$$
therefore, $\FF$ is quadratic over some subfield $\EE$ of its center by theorem \ref{379} again, and the map $\gamma_2^{-1}:\FF\to \tilde{L}_0$ is induced by an isomorphism $(\gamma_2^{-1},\phi_2^{-1}):(\FF,\EE,N^{\FF}_{\EE})\to (\tilde{L}_0,\tilde{\KK},\hat{q})$ of quadratic spaces. By taking
$$\alpha_{(1,2)}:=\big((\FF,\EE,N^\FF_\EE),\phi_2^{-1},\gamma_2^{-1},\phi_2^{-1},\gamma_2^{-1}\big)$$
as reparametrization for $\BBB_{(1,2)}$, we may assume $\BBB_{(1,2)}=\QQQ_Q(\FF,\EE,N^\FF_\EE)$ and $\gamma_2=\id_\FF$. Finally, the map $\gamma_1:\KK\to \EE$ is an isomorphism of fields by Hua's theorem.\qed
\end{bew}

\begin{prop}\label{406}
If we have
$$\FFF=\{ \BBB_{(1,2)}=\QQQ_Q(\tilde{L}_0,\tilde{\KK},\tilde{q}),\BBB_{(3,2)}=\QQQ_Q(\hat{L}_0,\hat{\KK},\hat{q}),\BBB_{(3,1)}=\TTT(\AA),\gamma_1,\gamma_2,\gamma_3\}\ ,$$
such that $(\hat{L}_0,\hat{\KK},\hat{q})$ is proper with $\dim_{\hat{\KK}} \hat{L}_0\geq 3$, we have
$$\FFF\cong \FFF_{443}\big((\tilde{L}_0,\tilde{\KK},\tilde{q}),\gamma_3\big):=\{ \BBB_{(1,2)}=\BBB_{(3,2)}=\QQQ_Q(\tilde{L}_0,\tilde{\KK},\tilde{q}),\BBB_{(3,1)}=\TTT(\tilde{\KK}),\gamma_3\}$$
for some $\gamma_3\in \Aut(\tilde{\KK})$.
\end{prop}

\begin{bew}
By Hua's theorem, the map $\gamma_1^{-1}:\tilde{\KK}\to \AA$ is an isomorphism of fields. By taking
$$\alpha_{(3,1)}:=\big(\tilde{\KK},\gamma_1^{-1},\gamma_1^{-1},\gamma_1^{-1},\big)$$
as reparametrization for $\BBB_{(3,1)}$, we may assume $\BBB_{(3,1)}=\TTT(\tilde{\KK})$ and $\gamma_3=\id_{\tilde{\KK}}$. Moreover, by theorem \ref{400}, the map $\gamma_2:\tilde{L}_0\to \hat{L}_0$ is induced by an isomorphism $(\gamma_2,\phi_2):(\tilde{L}_0,\tilde{\KK},\tilde{q})\to (\hat{L}_0,\hat{\KK},\hat{q})$ of quadratic spaces. By taking
$$\alpha_{(3,2)}:=\big((\tilde{L}_0,\tilde{\KK},\tilde{q}), \phi_2,\gamma_2,\phi_2,\gamma_2\big)$$
as reparametrization for $\BBB_{(3,2)}$, we may assume $\BBB_{(3,2)}=\QQQ_Q(\tilde{L}_0,\tilde{\KK},\tilde{q})$ and $\gamma_2=\id_{\tilde{L}_0}$. Now the map $\gamma_3:\tilde{\KK}\to\tilde{\KK}$ is an isomorphism of fields by Hua's theorem.\qed
\end{bew}

\begin{prop}\label{407}
If we have
$$\FFF=\{ \BBB_{(1,2)}=\QQQ_Q(\tilde{L}_0,\tilde{\KK},\tilde{q}),\BBB_{(3,2)}=\QQQ_Q(\hat{L}_0,\hat{\KK},\hat{q}),\BBB_{(3,1)}=\TTT(\AA),\gamma_1,\gamma_2,\gamma_3\}\ ,$$
such that $\dim_{\hat{\KK}}\hat{L}_0\leq 2$, we have
$$\FFF\cong \FFF_{443}(\tilde{\Xi},\hat{\Xi},\gamma_2):=\{ \BBB_{(1,2)}=\QQQ_Q({\KK},\hat{L}_0,{q}),\BBB_{(3,2)}=\QQQ_Q(\hat{L}_0,\hat{\KK},\hat{q}),\BBB_{(3,1)}=\TTT({\KK}), \gamma_3\}\ ,$$
for the quadratic space $\hat{\Xi}:=(\hat{L}_0,\hat{\KK},\hat{q})$ with $\dim_{\hat{\KK}} \hat{L}_0\leq 2$ (and thus of type (m)$\in\{$(ii),(iii)$\}$), some quadratic space $\tilde{\Xi}:=(\hat{L}_0,\KK,q)$ of type (m)$\in\{$(i),(ii),(iii)$\}$ and some isomorphism $\gamma_3:\hat{\KK}\to {\KK}$ of fields.
\end{prop}

\begin{bew}
By theorem \ref{400}, there is a quadratic space $(\hat{L}_0,\KK,q)$ of type (m)$\in\{$(i),(ii),(iii)$\}$ such that the map $\gamma_2:\tilde{L}_0\to \hat{L}_0$ is induced by an isomorphism $(\gamma_2,\phi_2):(\tilde{L}_0,\tilde{\KK},\tilde{q})\to (\hat{L}_0,\KK,q)$ of quadratic spaces. Notice that we don't need $(\hat{L}_0,\hat{\KK},\hat{q})$ to be proper to establish case (ii) of theorem \ref{400}. By taking
$$\alpha_{(1,2)}:=\big(( \hat{L}_0,\KK,{q}), \phi_2^{-1},\gamma_2^{-1},\phi_2^{-1},\gamma_2^{-1}\big)$$
as reparametrization for $\BBB_{(1,2)}$, we may assume $\BBB_{(1,2)}=\QQQ_Q(\hat{L}_0,\KK,{q})$ and $\gamma_2=\id_{\hat{L}_0}$. By Hua's theorem, the map $\gamma_1^{-1}:{\KK}\to \AA$ is an isomorphism of fields. By taking
$$\alpha_{(3,1)}:=\big({\KK},\gamma_1^{-1},\gamma_1^{-1},\gamma_1^{-1},\big)$$
as reparametrization for $\BBB_{(3,1)}$, we may assume $\BBB_{(3,1)}=\TTT({\KK})$ and $\gamma_1=\id_\KK$. Now the map $\gamma_3:\hat{\KK}\to\KK$ is an isomorphism of fields by Hua's theorem.\qed
\end{bew}

\begin{satz}\label{447} Let $\FFF$ be an integrable $443$-foundation such that both the quadrangles are of quadratic form type and such that $(\hat{L}_0,\hat{\KK},\hat{q})$ is proper. Then $\FFF$ is isomorphic to one of the following foundations:
\begin{enumerate}[label=(\roman*)]
\item $\FFF_{443}\big(\AA,(\FF_1,\FF_2),\gamma_2\big)$ as in proposition \ref{404}
\item $\FFF_{443}(\KK,\FF,\EE,\gamma_1)$ as in proposition \ref{405}
\item $\FFF_{443}\big((\tilde{L}_0,\tilde{\KK},\tilde{q}),\gamma_3\big)$ as in proposition \ref{406}
\item $\FFF_{443}(\tilde{\Xi},\hat{\Xi},\gamma_3)$ as in proposition \ref{407}
\end{enumerate}
\end{satz}

\begin{bew}
This results from propositions \ref{404}, \ref{405}, \ref{406}, and \ref{407}.\qed\end{bew}

\begin{bem}\label{477}
Notice that we supposed $(\hat{L}_0,\hat{\KK},\hat{q})$ to be proper only in proposition \ref{406}. The remaining results are valid even if both the parametrizing quadratic spaces are non-proper. In particular, case (ii) of theorem \ref{400} doesn't require $(\hat{L}_0,\hat{\KK},\hat{q})$ to be proper.
\end{bem}

\newpage

\chapter{Conclusion}

We summarize the previous results to have a complete list of integrable 443-Foundations. By remark \ref{477}, the theorem can be extended to non-proper quadratic spaces by adding quadratic spaces of type (i) except for case (vi). However, we don't give the existence proofs.

\newglossaryentry{CFFT}{type=results,name={{Classification of 443 Twin Buildings}},description={},sort=res}
\newglossaryentry{443pqa}{type=foundations,name={\ensuremath{\FFF_{443}(\Xi)}},description={443 foundation w.r.t. the pseudo-quadratic space $\Xi$ over a quaternion division algebra},sort=443}
\newglossaryentry{443pqb}{type=foundations,name={\ensuremath{\FFF_{443}(\Xi,\gamma)}},description={443 foundation w.r.t. the pseudo-quadratic space $\Xi$ over a separable quadratic extension $\KK$ and $\gamma\in\Aut(\KK)$},sort=443}

\begin{satz}[\textbf{\gls{CFFT}}]\label{459}
An integrable 443-foundation $\FFF$ with proper parameter systems is isomorphic to one of the following foundations:
\begin{enumerate}[label=(\roman*)]
\item $\gls{443pqa}$ for some proper pseudo-quadratic space $\Xi=(\HH,\HH_0,\sigma,L_0,q)$ such that $\HH$ is a quaternion division algebra:

\begin{center}\begin{tikzpicture}[scale=0.8,>=stealth,thick]
\begin{scope}
\coordinate (1) at (0,0);                    
\coordinate (2) at (3,0);
\coordinate (3) at (60:3);
\draws{0.55} (2)--(1);
\drawd{0.55} (3)--(2);
\drawd{0.55} (3)--(1);
\node()at (0.55,0.25) {\ssi{$\id_\HH^o$}};
\node()at (2.45,0.25) {\ssi$\sigma_s$};
\node()at (1.5,1.9) {\ssi$\id_{T}^o$};
\node()at (0.0,-0.5) {\ssi$3$};
\node()at (3,-0.5) {\ssi$1$};
\node()at (1.5,3) {\ssi$2$};
\node()at (1.5,-0.5) {\ssi{$\TTT(\HH)$}};
\node()at (3.3,1.5) {\ssi{$\QQQ_P(\Xi)$}};
\node()at (-0.3,1.5) {\ssi{$\QQQ_P(\Xi)$}};
\draw[->] (1,0) arc[radius=1,start angle=0, end angle=60];
\draw[<-] (2,0) arc[radius=1,start angle=180, end angle=120];
\draw[<-] ($(3,0)+(120:2)$) arc[radius=1,start angle=-60, end angle=-120];
\foreach \i in {1,...,3} {\fill (\i) circle (3pt);}
\end{scope}
\end{tikzpicture}\end{center}

\item $\gls{443pqb}$ for some proper pseudo-quadratic space $\Xi=(\KK,\KK_0,\sigma,L_0,q)$ such that $(\KK,\KK_0,\sigma)$ is quadratic of type (iii) and for some automorphism $\gamma\in \Aut(\KK)$:

\begin{center}\begin{tikzpicture}[scale=0.8,>=stealth,thick]
\begin{scope}
\coordinate (1) at (0,0);                    
\coordinate (2) at (3,0);
\coordinate (3) at (60:3);
\draws{0.55} (2)--(1);
\drawd{0.55} (3)--(2);
\drawd{0.55} (3)--(1);
\node()at (0.55,0.25) {\ssi{$\id_\KK$}};
\node()at (2.45,0.25) {\ssi$\gamma$};
\node()at (1.5,1.9) {\ssi$\id_{T}^o$};
\node()at (0.0,-0.5) {\ssi$3$};
\node()at (3,-0.5) {\ssi$1$};
\node()at (1.5,3) {\ssi$2$};
\node()at (1.5,-0.5) {\ssi{$\TTT(\KK)$}};
\node()at (3.3,1.5) {\ssi{$\QQQ_P(\Xi)$}};
\node()at (-0.3,1.5) {\ssi{$\QQQ_P(\Xi)$}};
\draw[->] (1,0) arc[radius=1,start angle=0, end angle=60];
\draw[<-] (2,0) arc[radius=1,start angle=180, end angle=120];
\draw[<-] ($(3,0)+(120:2)$) arc[radius=1,start angle=-60, end angle=-120];
\foreach \i in {1,...,3} {\fill (\i) circle (3pt);}
\end{scope}
\end{tikzpicture}\end{center}

\newglossaryentry{443i}{type=foundations,name={\ensuremath{\FFF_{443}(\Xi)}},description={443 foundation with respect to the involutory set $\Xi$},sort=443}

\item $\gls{443i}$ for some proper involutory set $\Xi=(\HH,\HH_0,\sigma)$ such that $\HH$ is a quaternion division algebra:

\begin{center}\begin{tikzpicture}[scale=0.8,>=stealth,thick]
\begin{scope}
\coordinate (1) at (0,0);                    
\coordinate (2) at (3,0);
\coordinate (3) at (60:3);
\draws{0.55} (2)--(1);
\drawd{0.55} (3)--(2);
\drawd{0.55} (3)--(1);
\node()at (0.55,0.25) {\ssi{$\id_\HH^o$}};
\node()at (2.45,0.25) {\ssi$\sigma_s$};
\node()at (1.575,1.9) {\ssi$\id_{\HH_0}^o$};
\node()at (0.0,-0.5) {\ssi$3$};
\node()at (3,-0.5) {\ssi$1$};
\node()at (1.5,3) {\ssi$2$};
\node()at (1.5,-0.5) {\ssi{$\TTT(\HH)$}};
\node()at (3.3,1.5) {\ssi{$\QQQ_I(\Xi)$}};
\node()at (-0.3,1.5) {\ssi{$\QQQ_I(\Xi)$}};
\draw[->] (1,0) arc[radius=1,start angle=0, end angle=60];
\draw[<-] (2,0) arc[radius=1,start angle=180, end angle=120];
\draw[<-] ($(3,0)+(120:2)$) arc[radius=1,start angle=-60, end angle=-120];
\foreach \i in {1,...,3} {\fill (\i) circle (3pt);}
\end{scope}
\end{tikzpicture}\end{center}

\newglossaryentry{443qa}{type=foundations,name={\ensuremath{\FFF_{443}\big( \AA, (\FF_1,\FF_2),\gamma\big)}},description={443 foundation w.r.t. the quadratic spaces $\Xi_i:=(\AA,\FF_i,N^\AA_{\FF_i})$ of type (i)-(v) and the isomorphism $\gamma: \FF_1\to \FF_2$},sort=443}

\item $\gls{443qa}$ for some proper quadratic spaces $\Xi_i:=(\AA,\FF_i,N^\AA_{\FF_i})$, $i=1,2$, of type (ii)-(v) and some isomorphism $\gamma: \FF_1\to \FF_2$ of fields:

\begin{center}\begin{tikzpicture}[scale=0.8,>=stealth,thick]
\begin{scope}
\coordinate (1) at (0,0);                    
\coordinate (2) at (3,0);
\coordinate (3) at (60:3);
\draws{0.55} (1)--(2);
\drawd{0.55} (3)--(2);
\drawd{0.55} (3)--(1);
\node()at (0.55,0.25) {\ssi{$\id_\AA$}};
\node()at (2.45,0.25) {\ssi$\id_\AA$};
\node()at (1.5,1.9) {\ssi$\gamma$};
\node()at (0.0,-0.5) {\ssi$3$};
\node()at (3,-0.5) {\ssi$1$};
\node()at (1.5,3) {\ssi$2$};
\node()at (1.5,-0.5) {\ssi{$\TTT(\AA)$}};
\node()at (3.3,1.5) {\ssi{$\QQQ_Q(\Xi_1)$}};
\node()at (-0.3,1.5) {\ssi{$\QQQ_Q(\Xi_2)$}};
\draw[<-] (1,0) arc[radius=1,start angle=0, end angle=60];
\draw[->] (2,0) arc[radius=1,start angle=180, end angle=120];
\draw[->] ($(3,0)+(120:2)$) arc[radius=1,start angle=-60, end angle=-120];
\foreach \i in {1,...,3} {\fill (\i) circle (3pt);}
\end{scope}
\end{tikzpicture}\end{center}

\newglossaryentry{443qb}{type=foundations,name={\ensuremath{\FFF_{443}(\KK,\FF,\EE,\gamma)}},description={443 foundation w.r.t. the quadratic spaces $(\KK,\FF,N^\KK_\FF)$, $(\FF,\EE,N^\FF_\EE)$ of type (i)-(iii) and the isomorphism $\gamma: \KK\to\EE$},sort=443}

\item $\gls{443qb}$ for some proper quadratic spaces $(\KK,\FF,N^\KK_\FF)$, $(\FF,\EE,N^\FF_\EE)$ of type (ii)-(iii) and some isomorphism $\gamma:\KK\to \EE$ of fields:

\begin{center}\begin{tikzpicture}[scale=0.8,>=stealth,thick]
\begin{scope}
\coordinate (1) at (0,0);                    
\coordinate (2) at (3,0);
\coordinate (3) at (60:3);
\draws{0.55} (1)--(2);
\drawd{0.55} (2)--(3);
\drawd{0.55} (3)--(1);
\node()at (0.55,0.25) {\ssi{$\id_\KK$}};
\node()at (2.45,0.25) {\ssi$\gamma$};
\node()at (1.5,1.9) {\ssi$\id_\FF$};
\node()at (0.0,-0.5) {\ssi$3$};
\node()at (3,-0.5) {\ssi$1$};
\node()at (1.5,3) {\ssi$2$};
\node()at (1.5,-0.5) {\ssi{$\TTT(\KK)$}};
\node()at (3.5,1.5) {\ssi{$\QQQ_Q(\FF,\EE,N^\FF_\EE)$}};
\node()at (-0.5,1.5) {\ssi{$\QQQ_Q(\KK,\FF,N^\KK_\FF)$}};
\draw[<-] (1,0) arc[radius=1,start angle=0, end angle=60];
\draw[->] (2,0) arc[radius=1,start angle=180, end angle=120];
\draw[->] ($(3,0)+(120:2)$) arc[radius=1,start angle=-60, end angle=-120];
\foreach \i in {1,...,3} {\fill (\i) circle (3pt);}
\end{scope}
\end{tikzpicture}\end{center}

\newglossaryentry{443qc}{type=foundations,name={\ensuremath{\FFF_{443}(\Xi,\gamma)}},description={443 foundation w.r.t. the quadratic space $\Xi=(L_0,\KK,q)$ with $\dim_\KK L_0\geq 3$ and $\gamma\in \Aut(\KK)$},sort=443}

\item $\gls{443qc}$ for some proper quadratic space $\Xi=(L_0,\KK,q)$ such that $\dim_\KK L_0\geq 3$ and some automorphism $\gamma\in \Aut(\KK)$:

\begin{center}\begin{tikzpicture}[scale=0.8,>=stealth,thick]
\begin{scope}
\coordinate (1) at (0,0);                    
\coordinate (2) at (3,0);
\coordinate (3) at (60:3);
\draws{0.55} (1)--(2);
\drawd{0.55} (2)--(3);
\drawd{0.55} (1)--(3);
\node()at (0.55,0.25) {\ssi{$\gamma$}};
\node()at (2.45,0.25) {\ssi$\id_\KK$};
\node()at (1.575,1.875) {\ssi$\id_{L_0}$};
\node()at (0.0,-0.5) {\ssi$3$};
\node()at (3,-0.5) {\ssi$1$};
\node()at (1.5,3) {\ssi$2$};
\node()at (1.5,-0.5) {\ssi{$\TTT(\KK)$}};
\node()at (3.3,1.5) {\ssi{$\QQQ_Q(\Xi)$}};
\node()at (-0.3,1.5) {\ssi{$\QQQ_Q(\Xi)$}};
\draw[<-] (1,0) arc[radius=1,start angle=0, end angle=60];
\draw[->] (2,0) arc[radius=1,start angle=180, end angle=120];
\draw[->] ($(3,0)+(120:2)$) arc[radius=1,start angle=-60, end angle=-120];
\foreach \i in {1,...,3} {\fill (\i) circle (3pt);}
\end{scope}
\end{tikzpicture}\end{center}

\newglossaryentry{443qd}{type=foundations,name={\ensuremath{\FFF_{443}(\Xi,\tilde{\Xi},\gamma)}},description={443 foundation w.r.t. the quadratic space $\Xi=(L_0,\KK,q)$ with $\dim_\KK L_0\leq 2$, the quadratic space $\tilde{\Xi}=(L_0,\tilde{\KK},\tilde{q})$ of type (i)-(iii) and the isomorphism $\gamma:\tilde{\KK}\to \KK$},sort=443}

\item $\gls{443qd}$ for some proper quadratic space $\Xi=(L_0,\KK,q)$ such that $\dim_\KK L_0\leq 2$, some proper quadratic space $\tilde{\Xi}=(L_0,\tilde{\KK},\tilde{q})$ of type (ii)-(iii) and some isomorphism $\gamma:\tilde{\KK}\to \KK$ of fields:

\begin{center}\begin{tikzpicture}[scale=0.8,>=stealth,thick]
\begin{scope}
\coordinate (1) at (0,0);                    
\coordinate (2) at (3,0);
\coordinate (3) at (60:3);
\draws{0.55} (1)--(2);
\drawd{0.55} (2)--(3);
\drawd{0.55} (1)--(3);
\node()at (0.55,0.25) {\ssi{$\gamma$}};
\node()at (2.45,0.25) {\ssi$\id_\KK$};
\node()at (1.575,1.875) {\ssi$\id_{L_0}$};
\node()at (0.0,-0.5) {\ssi$3$};
\node()at (3,-0.5) {\ssi$1$};
\node()at (1.5,3) {\ssi$2$};
\node()at (1.5,-0.5) {\ssi{$\TTT(\KK)$}};
\node()at (3.3,1.5) {\ssi{$\QQQ_Q(\Xi)$}};
\node()at (-0.3,1.5) {\ssi{$\QQQ_Q(\tilde{\Xi})$}};
\draw[<-] (1,0) arc[radius=1,start angle=0, end angle=60];
\draw[->] (2,0) arc[radius=1,start angle=180, end angle=120];
\draw[->] ($(3,0)+(120:2)$) arc[radius=1,start angle=-60, end angle=-120];
\foreach \i in {1,...,3} {\fill (\i) circle (3pt);}
\end{scope}
\end{tikzpicture}\end{center}

\end{enumerate}
\end{satz}

\begin{bew}
This holds by theorems  \ref{451}, \ref{450}, \ref{449}, \ref{448}, \ref{165} and \ref{447}.\qed
\end{bew}

\begin{bem}\ 
\begin{enumerate}[label=(\alph*)]
\item Notice the restriction in theorem \ref{448} which we did not mention in the above formulation.
\item Let $(\AA,\KK,N^\AA_\KK)$ be a quadratic space of type (m). Then we have
$$\MM:=\MM\big(\AA,\KK,N^\AA_\KK\big)=\MM(\AA)=:\tilde{\MM}$$
by lemma \ref{456}.
\end{enumerate}
\end{bem}

\begin{kor}
If we have $\dim_{\FF_1} \AA\geq 3$ in case (iv) of theorem \ref{459}, we have $$\Xi_2=\Xi_1\ .$$
\end{kor}

\begin{bew}
This results from theorem \ref{400} as we have $\MM(\Xi_1)=M(\AA)=\MM(\Xi_2)$.\qed
\end{bew}

\addtocontents{toc}{\noindent\protect\mbox{}\protect\hrulefill\par}
\part*{Appendix}\addtocounter{chapter}{1}\chaptermark{${\tilde{A}_2}$-Buildings Revisited}\markboth{Appendix}{}\addtocounter{chapter}{-1}
\addtocontents{toc}{\noindent\protect\mbox{}\protect\hrulefill\par}

\chapter{\texorpdfstring{$\boldsymbol{\tilde{A}_2}$}{Triangle}-Buildings Revisited}

We prove theorem \ref{24} without using the building at infinity.

\begin{prop}\label{396}
Let 
$$\FFF:=\{ \BBB_{(1,2)}=\TTT(\hat{\AA}),\BBB_{(2,3)}=\TTT(\AA), \BBB_{(3,1)}=\TTT(\tilde{\AA}), \gamma_{(1,2,3)},\gamma_{(2,3,1)},\gamma_{(3,1,2)}\}$$
 be an integrable foundation of type $\tilde{A}_2$ such that the defining field is a non-commutative skew-field and such that at least one glueing is positive. Then each glueing is positive.
\end{prop}

\begin{bew}
Assume that $\gamma_2=\gamma_{(1,2,3)}$ is positive. Without loss of generality we may assume that $\gamma_3=\gamma_{(2,3,1)}$ is negative (otherwise, we consider the opposite buildings and glueings). By taking $(\AA,\gamma_3,\gamma_3,\gamma_3)$ as reparametrization for $\TTT(\tilde{\AA})$, we may assume
\begin{align*}
\tilde{\AA}=\AA\ , && \gamma_3=\id_{\AA}\ ,
\end{align*}
and by taking $(\AA^o,\gamma_2^{-1},\gamma_2^{-1},\gamma_2^{-1})$ as reparametrization for $\TTT(\hat{\AA})$, we may assume
\begin{align*}
\hat{\AA}=\AA^o\ , && \gamma_2=\id_{\AA}^o\ .
\end{align*}
If we set 
\begin{align*}
h(s)&:= \mu\big(x_{(3,1)}^1(1_\KK)\big)^{-1}\mu\big(x_{(3,1)}^1(s^{-1})\big)\ , \qquad s\in \AA^*\ , \\
\tilde{h}(s)&:= \mu\big(x_{(1,2)}^1(1_\KK)\big)^{-1}\mu\big(x_{(1,2)}^1(s^{-1})\big)\ , \qquad s\in \AA^*\ ,
\end{align*}
we have
$$h(s)=\mu\big(x_{(3,1)}^1(1_\AA)\big)^{-1}\mu\big(x_{(3,1)}^1(s^{-1})\big)=\mu\big(x_{(1,2)}^1(1_\AA)\big)^{-1}\mu\big(x_{(1,2)}^1(\gamma_1(s)^{-1})\big)=\tilde{h}\big(\gamma_1(s)\big)$$
for each $s\in \AA^*$,
\begin{align*}
x_{(2,3)}^3(t)^{h(s)}&=x_{(3,1)}^3(t)^{h(s)}=x_{(3,1)}^3(t)=x_{(2,3)}^3(t\cdot s)
\end{align*}
and
\begin{align*}
x_{(2,3)}^2(t)^{h(s)}&=x_{(1,2)}^2(t)^{\tilde{h}(\gamma_1(s))}\\
&=x_{(1,2)}^2\big(\gamma_1(s)\circ t\circ \gamma_1(s)\big)=x_{(2,3)}^2\big(\gamma_1(s)\cdot t\cdot \gamma_1(s)\big)
\end{align*}
for all $s\in \AA^*,\ t\in \AA$ by lemma \ref{361}. Given $s\in \KK$, the Hua automorphism $h(\gamma_1^{-1}\big(s)\big)$ induces an automorphism $\alpha_s\in \Aut(\BBB_{(2,3)})$ which satisfies 
\begin{align*}
x_{(2,3)}^2(t)\mapsto x_{(2,3)}^2(s\cdot t\cdot s)\ , && x_{(2,3)}^3(t)\mapsto x_{(2,3)}^3\big(t\cdot \gamma_1^{-1}(s)\big)\ .
\end{align*}
If we set $\tilde{\alpha}:=(\id_\AA,\rho_{\gamma_1(s)},\rho_{\gamma_1(s)})$, then $\tilde{\alpha}\alpha_s$ satisfies
\begin{align*}
x_{(2,3)}^2(t)\mapsto x_{(2,3)}^2( s\cdot t\cdot s)\ , && x_{(2,3)}^3(t)\mapsto x_{(2,3)}^3(t)\ .
\end{align*}
By lemma \ref{66}, there is an element $c\in \AA^*$ such that
$$\forall\ t\in \AA:\qquad s ts=c t\ .$$
Lemma \ref{158} implies that we have
$$s\in Z(\AA)\ .$$
Since $s\in \AA$ is arbitrary, it follows that
$$\AA\subseteq Z(\AA)\qquad \lightning\ .$$
\qed
\end{bew}

\begin{prop}\label{397}
Let 
$$\FFF:=\{ \BBB_{(1,2)}=\TTT(\hat{\AA}),\BBB_{(2,3)}=\TTT(\AA), \BBB_{(3,1)}=\TTT(\tilde{\AA}), \gamma_{(1,2,3)},\gamma_{(2,3,1)},\gamma_{(3,1,2)}\}$$
 be an integrable foundation of type $\tilde{A}_2$ such that the defining field is a non-commutative skew-field and such that each glueing is positive. Then the skew-field $\AA$ is a quaternion division algebra and $\FFF$ is isomorphic to the foundation
$$\PPP_3^+(\AA):=\{ \tilde{\BBB}_{(1,2)}=\TTT(\AA^o),\tilde{\BBB}_{(2,3)}=\TTT(\AA), \tilde{\BBB}_{(3,1)}=\TTT(\AA^o),\tilde{\gamma}_1=\sigma_s,\tilde{\gamma}_2=\id_{\AA}^o,\tilde{\gamma}_3=\id_\AA^o\} \ .$$
\end{prop}

\begin{bew}
As $\gamma_2=\gamma_{(2,3,1)}$ is positive, we may take $(\AA^o,\gamma_2^{-1},\gamma_2^{-1},\gamma_2^{-1})$ as reparametrization for $\TTT(\hat{\AA})$. Therefore, we may assume
\begin{align*}
\hat{\AA}=\AA^o\ , && \gamma_2=\id^o_{\AA_0}\ .\end{align*}
As $\gamma_3=\gamma_{(2,3,1)}$ is positive, we may take $(\AA^o,\gamma_3^o,\gamma^o_3,\gamma_3^o)$ as reparametrization for $\TTT(\tilde{\AA})$. Therefore, we may assume
\begin{align*}
\tilde{\AA}=\AA^o\ , && \gamma_3=\id_{\AA}^o\ .
\end{align*}
\begin{itemize}
\item If we set
\begin{align*}
h_1(s)&:= \mu\big(x_{(3,1)}^1(1_{\KK})\big)^{-1}\mu\big(x_{(3,1)}^1(s^{-1})\big)\ ,  && s\in \AA^o\ , \\
\tilde{h}_1(s)&:= \mu\big(x_{(1,2)}^1(1_\KK)\big)^{-1}\mu\big(x_{(1,2)}^1(s^{-1})\big)\ ,  && s\in \AA^o\ ,
\end{align*}
we have
\begin{align*} h_1(s)&=\mu\big(x_{(3,1)}^1(1_{\AA^o})\big)^{-1}\mu\big(x_{(3,1)}^1(s^{-1})\big)\\
&=\mu\big(x_{(1,2)}^1(1_\AA)\big)^{-1}\mu\big(x_{(1,2)}^1(\gamma_1(s)^{-1})\big)=\tilde{h}_1\big(\gamma_1(s)\big)\end{align*}
for each $s\in \AA^o$,
\begin{align*}
x_{(2,3)}^3(t)^{h_1(s)}&=x_{(3,1)}^3\big(\id_{\AA}^o(t)\big)^{h_1(s)}=x_{(3,1)}^3\big(\id_{\AA}^o(t)\circ s\big)=x_{(2,3)}^3\big(\id_\AA^o(s)\cdot t\big)
\end{align*}
for all $s\in \AA^o,\ t\in \AA$ by lemma \ref{361} and
\begin{align*}
x_{(2,3)}^2(t)^{h_1(s)}&=x_{(1,2)}^2\big(\id_{\AA}^o(t)\big)^{\tilde{h}_1(\gamma_1(s))}\\
&=x_{(1,2)}^2\big(\gamma_1(s)\circ \id_{\AA}^o(t)\circ \gamma_1(s)\big)=x_{(2,3)}^2\big(\gamma_1^o(s)\cdot  t \cdot \gamma_1^o(s)\big)
\end{align*}
for all $s\in \AA^o,\ t\in \AA$ by lemma \ref{361}.
\item If we set
\begin{align*}
h_3(s)&:= \mu\big(x_{(3,1)}^3(1_{\AA})\big)^{-1}\mu\big(x_{(3,1)}^3(s)\big)\ , && s\in \AA^o\ , \\
\tilde{h}_3(s)&:= \mu\big(x_{(2,3)}^3(1_\AA)\big)^{-1}\mu\big(x_{(2,3)}^3(s)\big)\ , && s\in \AA\ ,
\end{align*}
we have
$$h_3(s)=\mu\big(x_{(3,1)}^1(1_{\AA^o})\big)^{-1}\mu\big(x_{(3,1)}^1(s^{-1})\big)=\mu\big(x_{(1,2)}^1(1_\AA)\big)^{-1}\mu\big(x_{(1,2)}^1(\id_\AA^o(s^{-1}))\big)=\tilde{h}_3\big(\id_\AA^o(s)\big)$$
for each $s\in \AA^o$,
\begin{align*}
x_{(2,3)}^3(t)^{h_3(s)}&=x_{(3,1)}^3\big(\id_\AA^o(t)\big)^{h_3(s)}\\
&=x_{(3,1)}^3\big(s^{-1}\circ \id_{\AA}^o(t)\circ s^{-1}\big)=x_{(2,3)}^3\big(\id_\AA^o(s)^{-1}\cdot t\cdot\id_\AA^o(s)^{-1}\big)
\end{align*}
for all $s\in \AA^o,\ t\in \AA$ by lemma \ref{361} and
\begin{align*}
x_{(2,3)}^2(t)^{h_3(s)}&=x_{(2,3)}^2(t)^{\tilde{h}_3(\id_\AA^o(s))}=x_{(2,3)}^2\big(\id_\AA^o(s)^{\sigma}\cdot t\cdot \id_\AA^o(s)\big)
\end{align*}
for all $s\in \AA^o,\ t\in \AA_0$ by lemma \ref{361}.
\end{itemize} 
Given $s\in \AA^o$, then $h_3(s)h_1(s)$ induces an automorphism $\alpha_s\in \Aut(\BBB_{(2,3)})$ which satisfies 
\begin{align*}
x_{(2,3)}^2(t)\mapsto x_{(2,3)}^2(\id_\AA^o(s)\cdot \gamma_1^o(s)\cdot t\cdot \gamma_1^o(s)\cdot \id_\AA^o(s))\ , && x_{(2,3)}^3(t)\mapsto x_{(2,3)}^3(t\cdot \id_\AA^o(s)^{-1})\ .
\end{align*}
If we set $\tilde{\alpha}:=(\id_{\AA^o},\rho_{\id_{\AA^o}(s)},\rho_{\id_{\AA^o}(s)})$, then $\tilde{\alpha}\alpha_s\in \Aut(\BBB_{(2,3)})$ satisfies
\begin{align*}
x_{(2,3)}^2(t)\mapsto x_{(2,3)}^2(\id_\AA^o(s)\cdot \gamma_1^o(s)\cdot t\cdot \gamma_1^o(s)\cdot \id_\AA^o(s))\ , && x_{(2,3)}^3(t)\mapsto x_{(2,3)}^3(t)\ .
\end{align*}
By lemma \ref{66}, there is an element $c\in \AA^*$ such that
$$\forall\ t\in \AA:\qquad \id_\AA^o(s)\cdot \gamma_1^o(s)\cdot t\cdot \gamma_1^o(s)\cdot \id_\AA^o(s)=c \cdot t\ .$$
Lemma \ref{158} implies that we have
$$\gamma_1^o(s) \cdot \id_\AA^o(s)\in C_\AA(\AA)=Z(\AA)\ .$$
Since $\AA$ is non-commutative, lemma \ref{160} shows that $\AA$ is a quaternion division algebra and that we have
\begin{align*} \gamma_1= \sigma_s:\AA^o\to \AA^o\ . \end{align*}
\qed
\end{bew}

\begin{satz} 
Let $\FFF$ be an integrable foundation of type $\tilde{A}_2$ such that the defining field is a non-commutative skew-field and such that at least one glueing is positive.  Then there is is a quaternion $\HH$ such that
$$\FFF\cong\PPP_3^+(\HH)=\{ \tilde{\BBB}_{(1,2)}=\TTT(\HH^o),\tilde{\BBB}_{(2,3)}=\TTT(\HH), \tilde{\BBB}_{(3,1)}=\TTT(\HH^o),\tilde{\gamma}_1=\sigma_s,\tilde{\gamma}_2=\id_{\HH}^o,\tilde{\gamma}_3=\id_\HH^o\} \ .$$
\end{satz}

\begin{bew}
This results from proposition \ref{396} and proposition \ref{397}.\qed
\end{bew}

\begin{bem}
As we did not use the fact that the residues embed into the building at infinity, which is a result in twin building theory, the above theorem is true for any foundation of an arbitrary affine building of type $\tilde{A}_2$ in which each Hua automorphism is induced by an automorphism of the whole building.
\end{bem}

\chapter{Jordan Automorphisms of Octonion Division Algebras Revisited}

We give a direct proof for proposition \ref{120}, independent of the characteristic.

\begin{satz}\label{502} Let $\OO$ be an octonion division algebra with center $\KK:=Z(\OO)$ and norm $N:=N^\OO_\KK$. Then we have
$$\Aut_J(\OO)=\Gamma L_N(\OO,\KK)\ .$$
\end{satz}

\begin{bew}
Let $(\p,\sigma)\in \Gamma L_N(\OO,\KK)$. Then we have $$\p(1_\OO)=\p(1_\OO\cdot 1_\OO)=\sigma(1_\OO)\cdot \p(1_\OO)=1_\OO\cdot \p(1_\OO)$$
and thus $\p(1_\OO)=1_\OO$. By definition, lemma \ref{456} and theorem \ref{434}, we have
\begin{align*}
\Aut_J(\OO)&=\Aut_J\big(\MM(\OO)\big)=\Aut_J\big( \MM(\OO,\KK,N)\big)\\
&=\{ (\p,\sigma)\in \Gamma L_N(\OO,\KK) \mid \p(1_\OO)=1_\OO\}=\Gamma L_N(\OO,\KK)\ .\end{align*}\qed
\end{bew}

\chapter{The Defining Field Revisited}

With theorem \ref{502}, we can prove theorem \ref{500} without using proposition \ref{23}.

\begin{lemma}\label{503}
Let $\FFF$ be a foundation such that there exists an edge $(a,b)\in A(F)$ with $\AA:=\AA_{(a,b)}$ an octonion division algebra. Then we have
\begin{align*} \forall\ (i,j)\in A(F):&&  \AA_{(i,j)}\cong \AA\ \vee\ \AA_{(i,j)}\cong \AA^o\ , && \TTT(\AA_{(i,j)})\cong \TTT(\AA)\ \vee\ \TTT(\AA_{(i,j)})\cong \TTT(\AA^o)\ .
\end{align*}
\end{lemma}

\begin{bew}
By theorem \ref{502}, each glueing $\gamma_{(i,j,k)}$ is a norm similarity. Therefore, we have $$\forall\ (i,j,k)\in G(F):\qquad \AA_{(i,j)}\cong \AA_{(j,k)}$$ by theorem (1.7.1) of \cite{S}.\qed
\end{bew}

\begin{satz}
Let $\tilde{\FFF}$ be an integrable foundation. Then there is an alternative division ring $\AA$ such that
\begin{align*}
\forall\ (i,j)\in A(\tilde{F}):\qquad \tilde{\AA}_{(i,j)}\cong \AA\ \vee\ \tilde{\AA}_{(i,j)}\cong \AA^o\ .
\end{align*}
\end{satz}

\begin{bew}
This is an immediate consequence of corollary \ref{34} and lemma \ref{503}.\qed
\end{bew}

\noindent \rule[0pt]{\textwidth}{1pt}
\vfill
\noindent \rule[10pt]{\textwidth}{1pt}

\chapter{Introduction in German}\selectlanguage{ngerman}

\section*{Historischer und theoretischer Hintergrund}

\noindent Bei der folgenden Darstellung orientieren wir uns stark an den Ausführungen in \cite{M} und \cite{AB}.

\subsection*{Zwillingsgebäude}

Gebäude wurden von J. Tits eingeführt, um halbeinfache algebraische Gruppen von einem geometrischen Standpunkt aus zu untersuchen. Eines der wichtigsten Resultate in der Gebäude-Theorie ist die Klassifikation der irreduziblen sphärischen Gebäuden vom Rang mindestens 3 in \cite{Ti74}. Mittlerweile gibt es einen vereinfachten Beweis in \cite{TW}, der auf der Klassifikation der Moufang-Polygone beruht. 

Vor über 25 Jahren definierten M. Ronan und J. Tits die Klasse der Zwillingsgebäude, eine natürliche Verallgemeinerung der sphärischen Gebäude. Motiviert wurde diese Definition durch die Theorie der Kac-Moody-Gruppen. Wir verweisen an diesem Punkt auf \cite{Ti} für allgemeine, weitergehende Informationen über Zwillingsgebäude.

Zwillingsgebäude verallgemeinern sphärische Gebäude in folgender Hinsicht: Bei sphärischen Gebäuden gibt es eine natürliche \textit{Oppositions-Relation} auf der Menge der Kammern, die die Struktur des Gebäudes wesentlich einschränkt. Die oben erwähnte Klassifikation der irreduziblen sphärischen Gebäuden vom Rang mindestens 3 basiert letzten Endes genau auf dieser Oppositions-Relation. Ein Zwillingsgebäude besteht nun aus zwei verschieden Gebäuden gleichen Typs, auf deren Kammern eine symmetrische Relation eingeführt wird, die ähnliche Eigenschaften wie die Oppositions-Relation von sphärischen Gebäuden besitzt. Ein Zwillingsgebäude ist also ein Tripel bestehend aus zwei Gebäuden gleichen Typs und einer Oppositions-Relation auf den zwei "`Hälften"' des Zwillingsgebäudes.

\newpage

\subsection*{Das Klassifikations-Programm für 2-sphärische Zwillingsgebäude}

Im Hinblick auf die Klassifikation der sphärischen Gebäude ergibt sich ganz natürlich die Frage, ob es auch möglich ist, Zwillingsgebäude höheren Ranges zu klassifizieren. Ein großer Teil von  \cite{Ti} beschäftigt sich mit genau dieser Problemstellung. Zunächst stellt sich heraus, dass eine solche Klassifikation nur unter der zusätzlichen Annahme möglich ist, dass die Einträge der zugehörigen Coxeter-Matrizen endlich sind. Wir nennen diese Gebäude \textit{2-sphärisch}. Das in \cite{Ti} beschriebene Klassifikations-Programm basiert auf der Vermutung, dass es für jede 2-sphärische Coxeter-Matrix vom Typ $M$ eine Bijektion zwischen Zwillingsgebäuden vom Typ $M$ und bestimmten Moufang-Fundamenten vom Typ $M$ gibt.

Fundamente wurden von M. Ronan und J. Tits in \cite{RT} eingeführt, um Kammer-Systeme zu beschreiben, die Kandidaten für die lokale Struktur eines Gebäudes sind. Grob gesagt sind Fundamente Amalgame von Rang-2-Gebäuden, die entlang bestimmter Rang-1-Residuen verklebt sind. Ist $c$ eine Kammer eines Gebäudes $\BBB$ vom Typ $M$, so bildet die Vereinigung $E_2(c)$ der Rang-2-Residuen, die diese Kammer $c$ enthalten, ein Fundament vom Typ $M$, das \textit{Fundament von $\mathit{\BBB}$ in $\mathit{c}$}. Der Ausdruck "`lokale Struktur"' ist also als eine Art 2-Umgebung einer Kammer $c$ des Gebäudes zu verstehen. 

Es ist eine (nicht triviale) Tatsache, dass die Fundamente zweier Kammern in derselben Hälfte eines Zwillingsgebäudes isomorph sind. Darüber hinaus ist der Isomorphie-Typ des Fundamentes der einen Hälfte durch den Isomorphie-Typ des Fundaments der anderen Hälfte eindeutig bestimmt. Umgekehrt besagt eine Verallgemeinerung von Tits' Erweiterungs-Satz durch B. Mühlherr und M. Ronan in \cite{MR}, dass ein Zwillingsgebäude in fast allen Fällen durch das Fundament einer seiner Hälften eindeutig bestimmt ist, siehe (5.10), (*5.11), (*9.11) und (*9.12) in \cite{AB} für eine Zusammenfassung. Somit ist das Fundament in einer Kammer eine klassifizierende Invariante des zugehörigen Zwillingsgebäudes, falls die folgende Bedingung erfüllt ist:
\begin{enumerate}[leftmargin=30pt]
\item[(CO)]
Kein Rang-2-Residuum ist isomorph zu einem Gebäude, das zu einer der Gruppen $B_2(2)$, $G_2(2)$, $G_2(3)$ und $^{2}{F_4}(2)$ gehört.
\end{enumerate}
Diese Bedingung garantiert, dass die Menge $c^o$ der einer Kammer $c\in \BBB_\epsilon$ ($\epsilon\in \{\pm\}$) gegenüber liegenden Kammern stets eine Galerie-zusammenhängende Teilmenge von $\BBB_{-\epsilon}$ ist.

In Anbetracht der bisherigen Überlegungen reduziert sich die Klassifikation der 2-sphärischen Zwillingsgebäude also auf die Klassifikation aller Fundamente, die als lokale Struktur eines Zwil-lingsgebäudes realisiert werden können. Wir nennen ein solches Fundament \textit{integrierbar}. Bei der Bestimmung der integrierbaren Fundamente verfährt man in zwei Schritten.

\subsection*{Schritt 1: Herausfiltern der nicht integrierbaren Fundamente}

In \cite{Ti} wird bewiesen, dass ein integrierbares Fundament Moufang ist, die Rang-2-Gebäude also Moufang-Polygone sind, deren Verklebungen mit den induzierten Moufang-Mengen auf den Rang-1-Residuen kompatibel sind. Eine erste notwendige Bedingung für die Integrierbarkeit eines Fundaments ist also die Moufang-Eigenschaft.

Als Folge dessen sind die Klassifikation der Moufang-Polygone in \cite{TW} und die Lösung des Isomorphie-Problems für Moufang-Mengen grundlegend bei der Untersuchung, welche Moufang-Polygone zu einem Fundament zusammengefügt werden können. Zudem kann man die Liste der möglicherweise integrierbaren Fundamente weiter einschränken, indem man bestimmte Automorphismen des Zwillingsgebäudes betrachtet, die sogenannten  \textit{Hua-Automorphismen}, die in einem engen Zusammenhang mit den Doppel-$\mu$-Maps der auftauchenden Moufang-Mengen stehen.

\subsection*{Schritt 2: Existenz- / Integrierbarkeits-Beweis}

Schließlich muss man beweisen, dass jeder der verbleibenden Kandidaten auch wirklich integrierbar ist. Das zugehörige Zwillingsgebäude ist dann bis auf Isomorphie eindeutig. In \cite{M} und seiner Habilitationsschrift \cite{MHab} entwickelte B. Mühlherr Techniken, die bestimmte Zwillingsgebäude als Fixpunktmenge in zu Kac-Moody-Gruppen gehörenden Zwillingsgebäuden realisieren. Er, H. Petersson und R. Weiss bereiten momentan ein Buch vor, das weitergehende, substanzielle Hintergründe liefert.

\newpage

\section*{Ziele und Ergebnisse}

Diese Arbeit widmet sich der Erstellung vollständiger Listen integrierbarer Fundamente für bestimmte Diagramm-Typen. Wir folgen hierbei dem Ansatz für die Klassifikation sphärischer Gebäude in \cite{TW}, wobei wir jedoch die dort verwendeten Techniken verfeinern müssen, da Fundamente im Allgemeinen nicht nur vom zugehörigen Diagramm und dem definierenden \protect{(Alternativ-)} Körper abhängen. Zum Beispiel gibt es für einen festen Schiefkörper $\AA$ in der Regel mehrere nicht isomorphe Fundamente vom Typ $\tilde{A}_n$: Automorphismen von $\AA$ spielen ebenfalls eine Rolle, was der Tatsache Rechnung trägt, dass es mehrere Möglichkeiten gibt, Moufang Polygone entlang eines Rang-1-Residuums zu verkleben.

Ein wesentlicher Aspekt ist die passende Parametrisierung von Sequenzen von Moufang-Polygonen bzw. deren Wurzelgruppen-Sequenzen mit den zugehörigen Kommutator-Relationen, um die verschiedenen Verklebungen sichtbar zu machen. Die entscheidende Subtilität ist die folgende: Jedes Moufang-Polygon wird zweimal parametrisiert, einmal für jede Richtung, in der die zugehörige Wurzelgruppen-Sequenz gelesen werden kann. Folglich erhalten wir Verklebungen zwischen gerichteten Moufang-Polygonen, und es macht einen Unterschied, ob wir $\id_\AA:\AA\to \AA$ oder $\id_\AA^o:\AA\to \AA^o$ betrachten, wobei $\AA^o$ der zu $\AA$ entgegengesetzte Schiefkörper ist: Ersteres ist ein Isomorphismus, während Letzteres ein Anti-Isomorphismus von Schiefkörpern ist.

Wie bereits erwähnt, ist das Herausfiltern nicht integrierbarer Fundamente eng verknüpft mit der Betrachtung von Moufang-Mengen und deren Isomorphismen. Deshalb beschäftigt sich ein großer Teil dieser Arbeit mit der Einführung der zugrundeliegenden Parameter-Systeme und der Lösung des Isomorphie-Problems für Moufang-Mengen. Viele Probleme wurden bereits gelöst, siehe \cite{K}, aber wir müssen die existierenden Ergebnisse für unsere Anforderungen sowohl verfeinern als auch erweitern und übersetzen deren Beweise in unser Setup. 

\subsection*{Fundamente mit einfachen Kanten (Simply Laced Foundations)}

Das Hauptergebnis dieser Arbeit ist die vollständige Klassifikation der Zwillingsgebäude mit einfachen Kanten via ihrer Fundamente. Natürlich ist die Hauptvoraussetzung für ein integrierbares Fundament die Moufang-Eigenschaft: Die Verklebungen sind Jordan-Isomorphismen, d.h., sie sind mit dem Jordan-Produkt $xyx$ verträglich. 

Ein mächtiges Werkzeug ist der Satz von Hua, siehe \cite{H} für eine Referenz, der das Isomorphie-Problem für Moufang-Mengen von Schiefkörpern löst: Jeder Jordan-Isomorphismus ist letztlich ein Iso- oder Anti-Isomorphismus von Schiefkörpern. Leider beinhaltet die Klasse von Parameter-Systemen für Moufang-Dreiecke zusätzlich Oktaven-Divisionsalgebren, was wegen der fehlenden Assoziativität zu einem gewissen Mehraufwand führt. Ein Nebenprodukt ist die Existenz von Jordan-Isomorphismen, die weder Iso- noch Anti-Isomorphismen von alternativen Ringen sind. Der aufwändigste Teil handelt von den Ausnahmefällen, in denen Oktaven auftauchen.

Wir geben an dieser Stelle einen Überblick über den Klassifikations-Prozess und streichen die Hauptideen heraus. Die folgenden Beobachtungen liefern die erste Einschränkung an Möglichkeiten: 
\begin{enumerate}[label=(\arabic*)]
\item Jedes Moufang-Dreieck ist über demselben alternativen Divisionsring $\AA$ definiert.
\item Ein integrierbares Fundament vom Typ $A_3$ ist notwendigerweise über einem Schiefkörper definiert, und die Verklebung ist notwendigerweise ein Isomorphismus von Schiefkörpern. 
\end{enumerate}

\noindent Folglich ist der entscheidende Schritt die Klassifikation der integrierbaren Fundamente vom Typ $\tilde{A}_2$, da dies die kleinsten sind, bei denen "`Nicht-Standard"'-Phänomene auftreten können. Die Theorie affiner Gebäude, von Bruhat-Tits-Gebäuden sowie die Theorie von Kompositions-Algebren, die bzgl. einer diskreten Bewertung komplett sind, ermöglichen uns weitere Einschränkungen:
\begin{enumerate}
\item[(3)] Zu einer Oktaven-Divisionsalgebra $\OO$ gibt es genau ein Zwillingsgebäude vom Typ $\tilde{A}_2$.
\item[(4)] Ein integrierbares Fundament vom Typ $\tilde{A}_2$, dessen Verklebungen Anti-Isomorphismen sind, ist notwendigerweise über einer Quaternionen-Divisionsalgebra definiert, und zu einem Quaternionen-Schiefkörper $\HH$ gibt es genau ein solches "`positives"' Zwillingsgebäude vom Typ $\tilde{A}_2$.
\end{enumerate}

\noindent Eine genauere Betrachtung der Gruppe der Jordan-Automorphismen einer \protect{Oktaven-Divisions-} \nohyphens{algebra} hilft dabei, die Klassifikation der über Oktaven definierten Zwillingsgebäude abzuschließen:

\begin{enumerate}
\item[(5)] Es gibt keine integrierbaren Fundamente über Oktaven, deren zugehöriger Graph ein Tetraeder ist. Insbesondere sind $\AAA_2(\OO)=\TTT(\OO)$ und $\tilde{\AAA}_2(\OO)$ die einzigen integrierbaren Fundamente über einer Oktaven-Divisionsalgebra $\OO$. 
\end{enumerate}

\noindent Schließlich liefert die folgende Beobachtung in Verbindung mit (4) eine wesentliche Einschränkung der Liste integrierbarer Fundamente über echten Schiefkörpern, die keine Quaternionen-Schiefkör- per sind:
\begin{enumerate}
\item[(6)] Ein integrierbares Fundament vom Typ $D_4$ ist notwendigerweise über einem Körper definiert.
\end{enumerate}
Als Konsequenz sind über echten Schiefkörpern, die keine Quaternionen-Schiefkörper sind, nur Fundamente integrierbar, deren Diagramm ein Kreis, ein String, ein Strahl oder eine Kette ist und deren Verklebungen Isomorphismen von Schiefkörpern sind.

Schließlich stellt die Kac-Moody-Theorie die Integrierbarkeit sicher, solange das zugehörige Coxeter-Diagramm ein Baum ist. Die restlichen Integrierbarkeits-Beweise basieren auf Techniken, die von B. Mühlherr entwickelt wurden. 

\subsection*{Jordan-Automorphismen von alternativen Divisionsringen}

In Hinblick auf den Satz von Hua, nach dem
$$\Aut_J(\DD)=\Aut(\DD)\cup \Aut^o(\DD)$$
für jeden Schiefkörper $\DD$, seine Gruppe $\Aut_J(\DD)$ von Jordan-Automorphismen, seine Untergruppe $\Aut(\DD)$ von Automorphismen und seine Menge $\Aut^o(\DD)$ von Anti-Automorphismen gilt, stellt sich die Frage nach einem ähnlichen Ergebnis für Oktaven-Divisionsalgebren. 

Im Beweis, dass integrierbare Tetraeder-Fundamente über Oktaven nicht existieren, definieren wir eine Teilmenge $\Gamma\subseteq \Aut_J(\OO)$, für die sich herausstellt, dass sie nicht die Standard-Involution $\sigma_s$ enthält. Die Elemente von $\Gamma$ sind Automorphismen von $\OO$, die zusätzlich mit einem der in \cite{TW} definierten "`Ausnahme"'-Jordan-Automorphismen multipliziert werden. Diese fixieren eine Quaternionen-Unteralgebra $\HH$ punktweise und wirken auf dem orthogonalen Komplement von $\HH$ als Konjugation. 

Dass $\Gamma$ eine Untergruppe von $\Aut_J(\OO)$ ist, kann man aus der Kenntnis der Automorphismen-Gruppe des zugehörigen Moufang-Dreiecks $\TTT(\OO)$ ableiten. Diese Untergruppe $\Gamma$ entspricht der Untergruppe $\Aut(\DD)$ im Satz von Hua, d.h., wir erhalten
$$\Aut_J(\OO)=\langle \sigma_s,\Gamma\rangle=\Gamma\cup \sigma_s\Gamma\ .$$
Die Strategie beim Beweis ist wie folgt:
\begin{enumerate}[label=(\arabic*)]
\item Jordan-Automorphismen eingeschränkt auf Unterkörper sind Ring-Monomorphismen, d.h., das Bild eines Unterkörpers ist wieder ein Unterkörper.
\item Als unmittelbare Konsequenz ergibt sich, dass Jordan-Automorphismen von Oktaven Norm-Ähnlichkeiten sind.
\item Die Ergebnisse aus \cite{S} erlauben uns, die Problemstellung auf Isometrien zurückzuführen, die eine Quaternionen-Unteralgebra punktweise fixieren.
\item Der Satz von Hua und das Skolem-Noether-Theorem erlauben uns zu zeigen, dass jeder Jordan-Automorphismus wirklich ein Produkt in $\langle \sigma_s,\Gamma\rangle$ ist.
\end{enumerate}

\subsection*{443-Fundamente}

Das zweite Ergebnis im Zusammenhang mit der Klassifikation der Zwillingsgebäude ist die Durchführung von Schritt 1 für 443-Fundamente. Dies sind Fundamente, deren Diagramm ein Dreieck ist und deren Moufang-Polygone zwei Vierecke und ein Dreieck sind. Obwohl wir uns in diesem Fall mit nur einem einzigen Diagramm beschäftigen, gibt es dennoch eine Vielzahl verschiedener integrierbarer 443-Fundamente, da es sechs Familien von Moufang-Vierecken gibt, die in dieser Konfiguration zudem auch noch oft zusammenpassen. Allerdings treten keine Vierecke vom Typ $E_n$, vom Typ $F_4$ und vom indifferenten Typ auf, weil ihre Moufang-Mengen nicht vom linearen Typ, also keine projektiven Geraden sind. 

Das Gleiche gilt zwar für Moufang-Mengen vom pseudo-quadratischen und involutorischen Typ, allerdings ist das zweite Panel des zugehörigen \textit{unitären} Vierecks vom linearen Typ, sodass es genau eine Möglichkeit für die Orientierung des Vierecks gibt. Die Lösung des Isomorphie-Problems für die auftauchenden Moufang-Mengen und die Kenntnis der Automorphismen-Gruppe eines unitären Vierecks führen zu folgendem Ergebnis: 
\begin{enumerate}[label=(\arabic*)]
\item Die auftauchenden pseudo-quadratischen Räume sind über einem Quaternionen-Schiefkörper $\HH$ oder einer separablen quadratischen Erweiterung $\EE$ definiert. Im ersten Fall gibt es genau ein integrierbares 443-Fundament über einem solchen pseudo-quadratischen Raum $\Xi$; im zweiten Fall hängt die Isomorphie-Klasse eines integrierbaren Fundaments zusätzlich von einem Automorphismus $\gamma\in \Aut(\EE)$ ab.
\item Die auftauchenden involutorischen Mengen sind über einem Quaternionen-Schiefkörper $\HH$ definiert, und es gibt genau ein integrierbares 443-Fundament über einer solchen involutorischen Menge $\Xi$.
\end{enumerate}
Schließlich bilden Vierecke vom quadratischen Typ die flexibelste Familie, da es Moufang-Mengen gibt, die sowohl vom quadratischen als auch vom linearen Typ sind, sodass diese Vierecke in jeder Orientierung verklebbar sind. Des Weiteren gibt es eine Stelle, an der wir uns auf echte quadratische Räume als parametrisierende Strukturen beschränken müssen, um für eine zufriedenstellende Darstellung Charakteristik-2-Phänomene ausschließen zu können. 

Im Gegensatz zur Klassifikation der integrierbaren Fundamente mit einfachen Kanten verzichten wir im Rahmen dieser Arbeit jedoch auf Schritt 2 des Klassifikations-Programms, da die Integrierbarkeits-Beweise andersartige, von B. Mühlherr, H. Petersson und R. Weiss eingeführte Techniken verwenden. Wie zuvor gibt es zwei Möglichkeiten, die Integrierbarkeit eines gegebenen Fundamentes zu beweisen: Entweder ist die universelle Überlagerung isomorph zu einem kanonischen Fundament, also einem Fundament, dessen Verklebungen alle die Identität sind und das integrierbar ist, falls das zugehörige Diagramm ein Baum ist, oder das Fundament kann als Fixpunkt-Struktur eines Tits-Index realisiert werden. Die erste Methode funktioniert bei 443-Fundamenten mit Vierecken vom quadratischen Typ, während die zweite bei 443-Fundamenten mit unitären Vierecken angewendet wird. 

\subsection*{Jordan-Isomorphismen von pseudo-quadratischen Räumen}

Genauso, wie der Satz von Hua bei der Klassifikation der integrierbaren Fundamente mit einfachen Kanten entscheidend ist, ist die Lösung des Isomorphie-Problems für die auftauchenden Moufang-Mengen ein wesentlicher Bestandteil bei der Klassifikation der integrierbaren 443-Fundamente. Wie bereits erwähnt, wurden viele Fälle von R. Knop in seiner Dissertation \cite{K} abgehandelt. Allerdings beschäftigt er sich nur mit kommutativen Moufang-Mengen, sodass wir die entsprechenden Resultate für Moufang-Mengen vom pseudo-quadratischen Typ ergänzen müssen.

Wir erhalten, dass Jordan-Isomorphismen zwischen zwei Moufang-Mengen vom pseudo-quadratischen Typ in der Regel von Isomorphismen zwischen den zugehörigen pseudo-quadrati- schen Räumen induziert werden. Genauer ist dies immer der Fall, wenn die Dimension mindestens 3 ist oder die beteiligte involutorische Menge echt ist. Folglich tauchen Ausnahmen nur bei pseudo-quadratischen Räumen kleiner Dimension auf, die über einem Quaternionen-Schiefkörper oder einer separablen quadratischen Erweiterung definiert sind. Glücklicherweise tauchen diese Ausnahme-Fälle nicht bei der Klassifikation integrierbarer 443-Fundamente auf, sodass beide Vierecke über demselben pseudo-quadratischen Raum $\Xi$ definiert sind. 

\newpage

\section*{Ausblick und offene Probleme}

\subsection*{Jordan-Isomorphismen}

In der Theorie der Moufang-Mengen spielen die $\mu$-Maps und die Hua-Maps eine entscheidende Rolle, da sie viele Informationen enthalten. Folglich stehen Jordan-Isomorphismen -- dies sind additive Isomorphismen, die mit den Hua-Maps verträglich sind -- in engem Zusammenhang mit Isomorphismen von Moufang-Mengen. Da die Hua-Maps in Summen und der Permutation $\tau$ dargestellt werden können, ist jeder Isomorphismus von Moufang-Mengen letztlich auch ein Jordan-Isomorphismus.

In diesem Zusammenhang stellt sich ganz natürlich die Frage, ob jeder Jordan-Isomorphismus auch ein Isomorphismus von Moufang-Mengen ist. Die Hua-Maps von scharf 2-fach transitiven Moufang-Mengen sind trivial, sodass die Frage in diesem Fall natürlich negativ beantwortet werden muss. Experten auf diesem Gebiet wie R. Weiss und T. De Medts gehen aber davon aus, dass beide Definitionen äquivalent sind, solange man sich auf "`echte"' Moufang-Mengen beschränkt.

\subsection*{Das Klassifikations-Programm}

Die Hauptvermutung im Zusammenhang mit dem Klassifikations-Programm ist die folgende, siehe Seite 5 in \cite{MHab}:
\begin{itemize}
\item[] Ein Moufang-Fundament vom 2-sphärischen Typ ist genau dann integrierbar, falls jedes Rang-3-Residuum integrierbar ist.
\end{itemize}
In seiner Habilitationsschrift \cite{MHab} deutet B. Mühlherr an, dass die Vermutung unter der zusätzlichen Annahme beweisbar wäre, dass alle Rang-3-Residuen sphärisch sind, was natürlich eine starke Einschränkung ist. Allerdings gibt es bislang noch keinen veröffentlichten Beweis.

Sobald die Vermutung bewiesen ist, reduziert sich das Klassifikations-Programm auf die Klassifikation der integrierbaren Moufang-Fundamente vom Rang 3. Die meisten Fälle können mit den in \cite{MHab} und \cite{M} eingeführten Methoden abgehandelt werden. Allerdings gibt es ein paar Ausnahmen, von denen die kompliziertesten vom Typ $\tilde{C}_2$, $\tilde{A}_2$ und $443$-Fundamente sind. Der $\tilde{A}_2$- und der $443$-Fall werden in dieser Arbeit gelöst, während es für den $\tilde{C}_2$-Fall (unveröffentlichte) Teilergebnisse von T. De Medts, B. Mühlherr, H. Van Maldeghem und R. Weiss gibt.

\subsection*{Die Klassifikation der Zwillingsgebäude mit einfachen Kanten}

Obwohl die Klassifikation der integrierbaren Fundamente mit einfachen Kanten abgeschlossen ist, machen wir keine Aussage darüber, ob zwei Fundamente unserer Liste isomorph sind. Verwendet man klassifizierende Invarianten und führt passende Parameter ein, kann man eine Liste paarweise nicht isomorpher Fundamente erstellen.

Ist das zugehörige Coxeter-Diagramm $\GGG_F$ ein Baum, so hängt das Fundament $\FFF$ nur vom definierenden Körper ab. Kreise im Diagramm sorgen für eine zusätzliche Abhängigkeit von "`Twists"', also von Automorphismen des definierenden Körpers $\AA$. Genauer:
\begin{itemize}
\item Ist $\AA$ ein Körper, so ist ein integrierbares Fundament $\FFF$ durch $\GGG_F$ und einen Homomorphismus $\p: \Pi_1(\GGG_F)\to \Aut(\AA)/\mathrm{Inn}(\AA)\cong \Aut(\AA)$ eindeutig bestimmt, wobei $\Pi_1(\GGG_F)$ die Fundamentalgruppe von $\GGG_F$ ist.
\item Ist $\AA$ ein Schiefkörper, aber kein Quaternionen-Schiefkörper, und $\FFF$ ein integrierbares Fundament vom Typ $\tilde{A}_n$, so ist das Fundament durch $n$ und ein Element von $\Aut(\AA)/\mathrm{Inn}(\AA)$ eindeutig bestimmt.
\item Ist der definierende Körper ein Quaternionen-Schiefkörper, so gilt ein ähnliches Resultat wie für einen Körper.
\end{itemize}

\noindent Zudem könnte man die Integrierbarkeits-Beweise in einigen Punkten überarbeiten, sobald die angewandte Theorie von B. Mühlherr, H. Petersson und R. Weiss vollständig entwickelt wurde.

\subsection*{Endliche Moufang-Fundamente}

Die eingeführte Terminologie und die Methoden von \cite{MHab} können verwendet werden, um zu zeigen, dass jedes lokal endliche Zwillingsgebäude vom 2-sphärischen Typ das Fixpunkt-Gebäude einer Galois-Wirkung im Sinne von B. Rémy ist, was bedeutet, dass es algebraischen Ursprungs ist.

\section*{Danksagungen}

Zunächst möchte ich meinem Betreuer Bernhard Mühlherr meinen ganz besonderen Dank dafür aussprechen, dass er meine Aufmerksamkeit auf das interessante Gebiet der Moufang-Fundamente und deren Moufang-Mengen gelenkt hat. Es war mir eine Freude, meinen Teil zum Klassifikations-Programm für Zwillingsgebäude beisteuern zu können. Viele ergiebige Diskussionen zeigten mir die richtigen Ansätze, also die möglicherweise zu beweisenden Aussagen und die zugehörigen Beweisideen. Seine Intuition ist beeindruckend.

Des Weiteren möchte ich Richard Weiss danken, der die Frage nach der Verallgemeinerung des Satz von Hua aufwarf und der mit seiner wunderbaren und detaillierten Arbeit über Moufang-Polygone, sphärische und affine Gebäude sowie deren Klassifikation die Grundlage für diese Arbeit gelegt hat. Viele kleine Fragen konnten mit der Hilfe von \cite{TW} beantwortet werden, und falls nicht, hatte er immer eine Idee, wo man nachschlagen könnte. 

Ralf Köhl erweiterte die Arbeitsgruppe um viele nette Leute und gab mir die Möglichkeit, an einem Forschungsprojekt über kompakte Untergruppen von Kac-Moody-Gruppen teilzunehmen. Seine Begeisterung und seine Hingabe sind bewundernswert und bilden die Grundlage für unsere gedeihende Arbeitsgruppe. 

Dank gebührt auch den folgenden Personen: Tom De Medts, der mir einen angenehmen Aufenthalt in Ghent bereitete, während dessen wir am Isomorphie-Problem für Moufang-Mengen vom pseudo-quadratischen Typ arbeiteten. Neben meiner Familie gibt es zuletzt noch zahlreiche Freunde, die mein Leben mit gemeinsam verbrachter Zeit, Unterhaltungen und Unternehmungen bereicherten, die einen angenehmen und wichtigen Gegenpol zu dem Mathematik genannten "`abstrakten Nonsense"' bildeten, mit dem ich mich die letzten 10 Jahre beschäftigt habe.

Der wichtigste Punkt, die die letzten 5 Jahre zu einer wunderbaren Zeit gemacht haben, ist der folgende: Bernhard und Ralf erlaubten mir die größtmögliche Flexibilität bei der Gestaltung meiner akademischen Arbeit. Dies ermöglichte mir insbesondere, zusammen mit meinen langjährigen Freunden Steffen Presse und Joram Gornowitz unser ambitioniertes Kino-Projekt zu realisieren.\\
 


\selectlanguage{english}

\newpage
\rhead[]{Bibliography}
\addtocontents{toc}{\protect\mbox{}\hfill\par}

\makeatletter
\renewcommand{\gls@doclearpage}{%
\ifthenelse{\equal{\@@glossarysec}{chapter}}{%
\@ifundefined{cleardoublepage}{\clearpage}{\clearpage}}{}%
}
\makeatother

\setlength{\glsdescwidth}{0.7\linewidth}
\rhead[]{List of Symbols}
\printglossary[type=symbols, style=long3col]
\setlength{\glsdescwidth}{0.65\linewidth}
\printglossary[type=foundations, style=long3col]
\printglossary[type=results, style=super4col]

\rhead[]{Index}
\lhead[Index]{}
\addcontentsline{toc}{chapter}{Index}
\printindex

\newpage

\thispagestyle{empty}

\vspace*{4cm}
\noindent \rule[10pt]{\textwidth}{1pt}
\section*{Selbständigkeitserklärung}
\addtocontents{toc}{\protect\mbox{}\hfill\par}
\addcontentsline{toc}{chapter}{Selbständigkeitserklärung}
Ich erkläre: Ich habe die vorgelegte Dissertation selbständig und ohne unerlaubte fremde Hilfe und nur mit den Hilfen angefertigt, die ich in der Dissertation angegeben habe.
 
Alle Textstellen, die wörtlich oder sinngemäß aus veröffentlichten Schriften entnommen sind, und alle Angaben, die auf mündlichen Auskünften beruhen, sind als solche kenntlich gemacht.

Bei den von mir durchgeführten und in der Dissertation erwähnten Untersuchungen habe ich die Grundsätze guter wissenschaftlicher Praxis, wie sie in der "`Satzung der Justus-Liebig-Universität Gießen zur Sicherung guter wissenschaftlicher Praxis"' niedergelegt sind, eingehalten.\\
\\
\\
\\
\underline{\qquad\qquad\qquad\qquad\qquad\qquad\qquad\qquad\qquad}\\
\\
Nidderau, im September 2013\\
\\
\\
\noindent \rule[0pt]{\textwidth}{1pt}


\end{document}